\documentclass[10pt]{article}
\usepackage[T1]{fontenc}
\usepackage[utf8]{inputenc}
\usepackage{lmodern}
\usepackage[margin=1in]{geometry}
\usepackage{amsmath,amssymb,amsthm,mathtools,mathrsfs}
\usepackage{microtype}
\newcommand{\gaussian}[2]{\genfrac{[}{]}{0pt}{}{#1}{#2}_{p}}
\usepackage{booktabs,longtable,array,enumitem,needspace}
\usepackage{tocloft}
\usepackage{xcolor}
\usepackage{hyperref}
\usepackage[nameinlink,noabbrev]{cleveref}

\hypersetup{
	pdftitle={Structural Packaging: Foundations, Faithful Amalgamation, and Rigidity},
	pdfauthor={Rongchao Ma},
	colorlinks=true,
	linkcolor=blue,
	citecolor=blue,
	urlcolor=blue,
}

\usepackage{xurl}
\numberwithin{equation}{section}

\newtheorem{theorem}{Theorem}[section]
\newtheorem{proposition}[theorem]{Proposition}
\newtheorem{lemma}[theorem]{Lemma}
\newtheorem{corollary}[theorem]{Corollary}
\theoremstyle{definition}
\newtheorem{definition}[theorem]{Definition}
\newtheorem{example}[theorem]{Example}
\newtheorem{problem}[theorem]{Problem}
\crefname{lemma}{lemma}{lemmas}
\Crefname{lemma}{Lemma}{Lemmas}
\crefname{problem}{Problem}{Problems}
\Crefname{problem}{Problem}{Problems}
\theoremstyle{remark}
\newtheorem{remark}[theorem]{Remark}

\AddToHook{env/proposition/begin}{\crefalias{theorem}{proposition}}
\AddToHook{env/lemma/begin}{\crefalias{theorem}{lemma}}
\AddToHook{env/corollary/begin}{\crefalias{theorem}{corollary}}
\AddToHook{env/definition/begin}{\crefalias{theorem}{definition}}
\AddToHook{env/example/begin}{\crefalias{theorem}{example}}
\AddToHook{env/problem/begin}{\crefalias{theorem}{problem}}
\AddToHook{env/remark/begin}{\crefalias{theorem}{remark}}

\newcommand{\Pkg}{\operatorname{Pkg}}

\newcommand{\Match}{\operatorname{Match}}

\newcommand{\LL}{\mathsf L}
\newcommand{\RR}{\mathsf R}
\newcommand{\MM}{\mathsf M}
\newcommand{\WW}{\mathsf W}

\newcommand{\NN}{\mathbb N}

\newcommand{\id}{\mathrm{id}}

\newcommand{\im}{\operatorname{im}}

\newcommand{\hr}{\operatorname{hr}}

\newcommand{\Ind}{\mathbf 1}
\newcommand{\Esep}{\ensuremath{\mathsf E_{\mathrm{sep}}}}

\newcommand{\Con}{\operatorname{Con}}

\newcommand{\ff}{\mathbb F}

\providecommand{\NP}{\mathsf{NP}}

\title{Structural Packaging: Foundations, Faithful Amalgamation, and Rigidity}

\author{
	Rongchao Ma \\
	\textit{Department of Physics, University of Alberta} \\
	\textit{Edmonton, Alberta, Canada} \\
}

\date{\today}

\begin{document}
	\maketitle

	\begin{abstract}
		A packaging system records packages, their incident units, candidate spaces, and selected families.
		We introduce structural packaging to specify how realized packages form assemblies:
		interfaces compare observations, formation relations select the admitted combinations of states, and native realization laws specify how formations produce	packaged endpoints.
		We study when assemblies preserve local states,	functional dependencies and nondependencies, and prescribed common information.
		
		We describe the composition of formation witnesses and construct, for finite presentations, the coarsest observation-preserving	congruence quotient that lifts every tree witness with an arbitrary representative prescribed at one variable.
		For finite equational join-tree atlases with jointly separating homomorphic readouts and equal separator images, a faithful native assembly exists in the specified presentation context exactly when the matching algebra has the determination closure generated by the local charts.
		An ordered-pair test either supplies the required witnesses or identifies a nondependency that no compatible source can preserve.
		
		For finite Mal'cev join-tree atlases, chart coverage and exact common information at every full tree cut force the source onto the matching	algebra, giving the sharp minimum source cardinality and its equality classification.
		For two copies of a finite group quotient with kernel of order $m\ge2$, preservation of chart states, whole-chart determination, and exact common information has sharp minimum $m^2$ states per quotient value for group sources and $2m-1$ for relational sources.
		The ratio is unbounded.
		
		Further results characterize substitution of finite group fragments, determine source minima under prescribed cuts, and give witness and repair criteria for cyclic formations.
	\end{abstract}

	\noindent\textbf{Keywords.} Structural packaging, faithful amalgamation,
	functional dependence, common information, subdirect products, source rigidity.
	
	\clearpage
	\tableofcontents
	\clearpage

	\section{Introduction}
	\label{sec:intro}

	The \textbf{packaging systems} introduced in \cite{Ma2026} record packages, their incident units, candidate spaces, and selected families.
	A realized system also records a formation datum and the judgment relating it to a packaged endpoint in a specified mathematical context.
	These data describe the realized packages.
	To describe their assembly, one must additionally specify how their states are compared, which combinations are admitted, and what information the resulting system must preserve.
	
	We introduce \textbf{structural packaging} to supply these assembly data.
	Shared observation sites and interface maps specify how states are compared.
	Formation relations select the matching combinations of input and output states that are admitted together.
	A witness assigns states to the package occurrences in a formation so that every local formation relation is satisfied.
	Composite formations require these assignments to agree at shared occurrences.
	
	A two-state example shows why the formation relations are additional
	data. Take two realized packaging systems, each with one package
	whose selected states are $0$ and $1$, and a prospective output with
	the same selected states. Let each package expose its state through
	the identity map at a common interface. The matching
	input--input--output triples are then $(0,0,0)$ and $(1,1,1)$.
	One formation relation admits only the first triple; another admits
	both. Both respect the interface comparisons and the hierarchy in
	which the inputs have rank zero and the output has rank one
	(\cref{ex:two-base-operational-choice}). Thus the packages, interfaces,
	and hierarchy do not determine which matching formations are admitted.
	
	An operational witness specifies mutually admissible states.
	A native realization additionally requires a formation datum and
	an output endpoint satisfying the prescribed realization law.
	The central problem is to determine when such an assembly preserves
	the required local states, recovery laws, and common information.
	A functional dependency expresses that one collection of observations
	determines another observation. Preserving its failure requires two
	states that agree on the first collection and differ on the latter.
	Faithful assembly must retain these distinctions as well as the
	dependencies themselves.
	
	The choice of native realization law can change the minimum size
	of a source even when the observed charts and preservation
	requirements are fixed. Let $\pi:G\twoheadrightarrow Q$ be a
	surjective homomorphism of finite groups with
	$m=|\ker\pi|\ge2$, and take two copies of this quotient map.
	Require both chart images to be $G$, the dependencies among the two
	whole-chart readouts and their common quotient readout to equal
	those of $G\times_QG$, and the common information of the two charts
	to be exactly the quotient value. The sharp minimum source sizes are
	\[
	N_{\mathrm{grp}}=|Q|m^2,
	\qquad
	N_{\mathrm{rel}}=|Q|(2m-1).
	\]
	A minimum relational source uses a spanning tree of pairs over each
	quotient value. Every admissible group source with homomorphic chart
	maps must cover the full matching group
	(\cref{thm:unbounded-native-resource-separation}). The ratio of the
	two minima is unbounded.
	
	These questions organize the paper: how formation witnesses compose
	and descend through reductions; when prescribed packages admit a
	faithful native assembly; and which realization requirements force
	a particular source or a sharp minimum size.

	\subsection{From matching to faithful assembly}
	\label{subsec:intro-motivation}

	Two small examples explain why both directions of recovery matter. First
	consider the three scalar observations $(a,b,a+b)$ over $\mathbb F_p$
	with constant interface. On this source, any two coordinates determine
	the third. Recombining the first coordinate freely with the last two
	produces all of $\mathbb F_p^3$ and destroys that determination, even
	though the three individual chart images stay the same.
	
	Conversely, take charts $(A,B)$ and $(B,C)$ with $B\in\{0,1,2,3\}$.
	On the first chart $A$ takes the respective values $0,0,1,1$, and on the
	second $C$ takes those same values. The local functional dependencies
	include $B\to A$ and $B\to C$. They do not entail $A\to C$. Yet matching
	the two fixed tables on $B$ forces $A=C$. Every compatible subsource
	inherits the extra dependency (\cref{ex:correlated-join}).
	
	For a source $H$ with descriptor set $D$, let $f_T$ collect the
	observations named in $T\subseteq D$. Its recovery closure is
	\[
	c_H(T)=\{d\in D:\ f_T(x)=f_T(y)\Longrightarrow f_d(x)=f_d(y)
	\text{ for all }x,y\in H\}.
	\]
	Write $c_0$ for the closure generated by the local recovery laws.
	Compatibility gives $c_0(T)\subseteq c_H(T)$. The reverse inclusion
	requires, for each $d\notin c_0(T)$, a pair of states that agree on $T$
	and differ on $d$. We write $\mathsf F$ for this equality of global
	recovery closures. A full faithful amalgamation must also cover the local
	states ($\mathsf S$), realize precisely the prescribed common information
	of the two sides of each protected cut ($\mathsf E$), and arise by the
	specified native formation law ($\mathsf N$). Chartwise equality of
	recovery closures is denoted by $\mathsf C$
	(\cref{def:preservation-contract}). The choice of cuts,
	readouts, and formation law is part of the data.

	\subsection{Principal results}
	\label{subsec:principal-results}

	\paragraph{Formation and reduction.}
	To compose two formations, identify their shared package occurrences and
	project away the internal states (\cref{prop:substitution}). The marked
	relational presentation retains both these witnesses and the judgments
	on the realized bases (\cref{prop:rel-representation}). Identifying local
	states can create a spurious composite witness. For a finite presentation
	with specified readouts, operations, and formation relations,
	\cref{thm:canonical-reduction-canonical} constructs the coarsest quotient
	through which every witness on a tree lifts with any one local representative
	prescribed. Without a prescribed representative, even depth-two formations can realize every graph coloring partition as a central quotient (\cref{thm:foundational-chromatic-reductions}). The native interpretation of the canonical quotient uses the additional quotient law
	\textup{(NQ)}. In the equational case, \cref{prop:equational-levelii-bridge}
	constructs the package formation from compatible chart maps and recovers
	those maps from a sound one-step formation under its stated output-coverage
	and readout hypotheses.

	\paragraph{Existence and rigidity.}
	Let $W$ be the full matching algebra of a finite equational join-tree
	atlas with jointly separating homomorphic readouts and equal actual
	images at each separator. Then a faithful native source exists precisely
	when $c_W=c_0$. If the equality fails, restricting to a subsource cannot
	restore a missing witness. If it holds, the prescribed equational
	presentation realizes $W$
	(\cref{prop:equational-levelii-bridge,thm:amalgamation-existence-criterion}).
	For a finite Mal'cev atlas, any source $H$ that covers the charts and has
	exact common information at every full tree cut maps onto $W$. Thus
	$H/\ker\psi\cong W$ for the joint chart map $\psi$, and $|H|=|W|$
	exactly when that map is injective (\cref{thm:native-source-rigidity,thm:faithful-matching-minimum}).
	The uniform permutability thresholds for lifting through complete separator
	signatures and for full-cut rigidity are respectively three and two
	(\cref{prop:finite-goursat-boundary,cor:finite-join-tree-rigidity-boundary}).

	\paragraph{What an interface must retain.}
	Two marked fragments are equivalent for a specified collection of
	tests if substituting one for the other gives the same outcomes under
	every allowed continuation. For finite group fragments the recovery
	record consists of boundary images of agreement kernels and
	determinations within the boundary kernel. For numerical range tests
	one keeps the boundary images and the kernel fiber sizes. Equality of
	the appropriate records characterizes equivalence, including on finite
	matching graphs with cycles. Each distinction between records is
	detected by an explicit attachment
	(\cref{thm:complete-compositional-invariants}).
	When continuation must preserve state coverage, the attainable record
	sets classify the families of realizations for the specified tests.
	When hidden kernels must be coupled by a specified subgroup, the
	classification uses marked extension classes together with readout
	kernels (\cref{thm:contract-filtered-spectrum,thm:coupled-extension-spectrum,thm:affine-coupled-spectrum}).
	Normalizer complements detect the existence of native lifts, and an
	abelian normalizer kernel permits a cohomological enumeration subject
	to the contract (\cref{thm:chief-normalizer-complement,thm:cohomological-contract-enumerator,cor:abelian-complement-contracts}).

	\paragraph{Recombination and the size of a source.}
	For a finite family of charts over one boundary, the minimal premises of
	whole-chart determination conditional on that boundary define a hypergraph. Its components give the
	finest partition across which one can freely recombine states and
	preserve the original recovery closure, even when the state sets are
	infinite (\cref{thm:determinants}). For a product boundary of finite
	groups, shared primes among central chief factors govern which ports must
	be observed together to guarantee faithfulness for the block-local
	quotient vocabulary (\cref{thm:central-prime-blocks}).
	
	For central elementary abelian extensions of a finite group $Q$, let
	$\Pi$ be the partition into the atoms of the protected cuts. A
	representable polymatroid $f$ determined by the charts gives
	\[
	\min_H |H|=|Q|p^{\sum_{B\in\Pi}f(B)} .
	\]
	The minimum ranges over chart-surjective subgroups that permit
	recombination at those cuts. Minima for all families of singleton cuts and the chart kernel
	dimensions recover $f$. Conditional range minima for several requests
	need not have a common minimizer. We give exact simultaneous-attainment
	tests for whole charts and for fixed chart quotient observations, given the common boundary
	(\cref{thm:min,cor:rank,thm:conditional,thm:incomparable-test,thm:partial-chart-reads}).
	With a prescribed quotient image and kernel the answer changes sharply.
	For every odd prime, a fixed group chart and interface type on cycles
	realize every finite downward-closed cut family containing the empty set
	(\cref{thm:complex}).

	\paragraph{Globalization and repair.}
	On cycles, a witness must satisfy all return constraints at once.
	Fixed-point and linear tests describe that obstruction, and simultaneous
	repair requires compatible choices of witnesses for all missing
	nonimplications (\cref{thm:native-witness-holonomy,thm:witness-globalization,thm:joint-faithful-repair}).
	For a finite nonabelian simple group $S$, we count exactly the subdirect
	subgroups of $S^n\times S^n$ that can be distinguished by a native
	probe language naming every such subgroup
	(\cref{thm:probe-complete-noncollapse}).

	\subsection{Related work}
	\label{sec:related}

	The closest results concern three different choices of what must be
	preserved under assembly. Join theory asks whether tables fit together.
	Continuation equivalence tests how an attachment behaves after
	substitution. Common-information and source-minimization results fix
	observations or resource costs.

	\paragraph{Consistency and dependence.}
	Vorob'ev studied extensions of consistent families of measures
	\cite{Vorobev1962}. The join-tree theorems give a local-to-global
	criterion for relational tables on acyclic schemas
	\cite{BeeriEtAl1981,Beeri1983,Yannakakis1981}. Related formulations use
	presheaves and valuations \cite{Stacks,ShenoyShafer,AbramskyBrandenburger}.
	Atserias and Kolaitis extended the acyclicity characterization to
	positive semirings \cite{AtseriasKolaitisSemirings}. Recent work treats consistency over
	monoids and functional dependence on locally consistent relations
	\cite{AtseriasKolaitis2025,BarlagEtAl2026}. For annotated relations,
	even the usual join of two consistent relations may fail to witness
	consistency. This motivates the study of consistency witness functions
	\cite{AtseriasKolaitisWitnesses,AtseriasKolaitis2026}.

	Armstrong's dependency calculus and Armstrong relations characterize
	implication by agreement of tuples \cite{Armstrong1974,BeeriDFS1984,Fagin1982}.
	Dependencies in joins and integrated views have also been studied directly
	\cite{SagivShmueli1989,ComignaniEtAl2022}. The extension of local closure operators is studied in
	\cite{Ma2026,Kim2026ClosureAtlases}.
	Our existence problem fixes the chart instances and tests all
	nonimplications in their matching source. Minimal premises in our
	decomposition form the canonical direct unit basis
	\cite{BertetMonjardet2010}. The $D$-base and $D$-relation are studied
	in \cite{ANV2025}. Earlier work on Armstrong relations and decompositions
	of closure operators gives related structural results
	\cite{DemetrovicsFurediKatona1985,Libkin1993}. Here the partition also
	determines a least recombinant enlargement in the marked chart product.

	\paragraph{Interfaces, quotients, and invariants.}
	The Myhill--Nerode principle identifies states by their possible
	continuations \cite{Nerode1958,ComonTATA}. For higher-dimensional
	automata, regularity is equivalent to finiteness of the prefix quotient
	\cite{FahrenbergZiemianski2024}. Graphical realizations provide
	notions of observability, controllability and local reduction
	\cite{Forney2001,Forney2014,ForneyGluesingLuerssen2013}.
	Our finite group records compute the outcomes of two specific
	continuation languages. Their construction uses the structure of
	subdirect products \cite{Goursat1889,BauerSenZvengrowski2015,LaiMengZhang2025},
	Scott's diagonal decomposition for nonabelian simple factors
	\cite{Scott1980}, and normal subgroups of direct products
	\cite{DeMedtsTarnauceanu2011}. A related commutator criterion for
	equational domains occurs in \cite{AichingerBehrischRossi2025}.
	Direct images of congruences in Goursat varieties provide the algebraic
	input for the interface saturation criterion
	\cite{HagemannMitschke1973,CarboniEtAl1993,Gran2020Regular}.

	\paragraph{Common information, minima, and obstructions.}
	For two finite observations, their common functions are the functions on
	the connected components of the bipartite support graph
	\cite{GacsKorner1973,Nair2013}. Rectangular recombination is the
	lossless-join property associated with multivalued dependencies
	\cite{Fagin1977,Hannula2024}. Fleischer's description of subdirect
	products in congruence-permutable varieties underlies the full-cut
	rigidity argument \cite{Malcev1954,Fleischer1955,MayrRuskuc2019,BurrisSankappanavar1981}.
	Tree realizations of a prescribed code minimize state or constraint
	complexity on the tree \cite{Forney2001,Kashyap2009}. Here the charts
	and separator maps are prescribed, and the quantity minimized is the
	order of a global source.

	Matroid parallel connections, Horn representations of dependence, and
	gain graphs provide precedents for the tree representations and cycle
	transport tests below
	\cite{Brylawski1975,BakerLorscheidWalshZhang2025,BercziMatroid2024,Zhu2026,Zaslavsky1989Gain}.
	Query-preserving minimization and database repairs optimize related
	changes to relational instances
	\cite{LivshitsKimelfeldRoy2018,LivshitsKimelfeldWijsen2021,CarmeliEtAl2024,HuSintos}.
	Our repair problem also protects common information and witnesses to
	failed dependencies. Product homomorphisms and generating sets of
	powers of finite algebras give related complexity problems
	\cite{TenCateDalmau2015,ZhukPowers}.

	\subsection{Organization}
	\label{subsec:organization}

	\Cref{sec:foundations} passes from realized packages to formation rules,
	their composition, and quotients preserving witnesses. Within it,
	\cref{subsec:canonical-reduction-subsection} gives the lifting criterion
	for trees. \Cref{sec:native-amalgamation} constructs the native formation
	for finite equational chart data and gives the exact test for faithfulness.
	
	\Cref{sec:compositional-invariants,sec:recovery,sec:minima,sec:limits}
	derive the group records, cut decompositions, and sharp source minima.
	\Cref{subsubsec:nonsplit-chief-monodromy,subsec:contract-filtered-invariants}
	treat extension lifting and composition with constraints on the kernel.
	\Cref{sec:faithful-tree,sec:terminal-representations} determine what
	forces a minimum source on a tree.
	\Cref{sec:native-witness-descent,sec:simultaneous-faithful-repair,sec:fixed-law-noncollapse}
	treat cyclic witnesses, repair, and comparison of native laws.
	\Cref{sec:discussion} states the remaining classification questions. Native quotient constructions accompany the finite equational presentation law in
	\cref{subsec:canonical-reduction-native-verification}; quotient counterexamples
	follow the canonical reduction theorem in \cref{app:canonical-reduction-limits}.
	The linear and group tests for cycle constraints are developed in
	\cref{sec:native-witness-descent}.

	\section{Foundation of structural packaging}
	\label{sec:foundations}

	A realized system in the packaging framework of \cite[Section~2]{Ma2026} records a context, a formation datum, and the judgment relating that datum to a packaged endpoint.
	This section adds observation sites, interfaces, and relations for assembling states.
	Every realized base admits such an enrichment, and the three coherence and hierarchy conditions are independent (\cref{prop:conservative,prop:independence}).
	The substitution and reduction theorems	(\cref{prop:substitution,thm:canonical-reduction-canonical}) describe how admitted tuples compose and which identifications preserve their tree witnesses.
	For the weaker requirement of unpinned exactness, \cref{thm:foundational-chromatic-reductions} realizes the coloring partitions of any finite graph as exact reductions of depth-two formations.

	The marked relational encoding recovers the whole presentation (\cref{prop:rel-representation}).
	When only current readout values matter, their joint image is the coarsest quotient (\cref{thm:universal-observable-quotient}).
	When a summary must also support future formations, observations after common continuations determine which histories it can identify (\cref{thm:canonical-contextual-summary}).

	\Needspace{12\baselineskip}
	\paragraph{Foundational notation.}
	\begin{center}
		\small
		\begin{tabular}{@{}p{0.32\textwidth}p{0.61\textwidth}@{}}
			\toprule
			Symbol & Role \\
			\midrule
			$\mathbb P=(\Xi,P,\mathfrak P)$ & Native context, formation datum, and realized endpoint. \\
			$\mathfrak L$ & Level-II structural presentation over realized bases. \\
			$\LL(K),\RR(K)$ & Candidate space and selected family at a package. \\
			$D,\rho$ & Observation sites and their restriction maps. \\
			$\alpha,\MM_\alpha,\WW_T$ & Primitive formation, admitted tuples, and tree witnesses. \\
			$q:X\twoheadrightarrow Y$ & Typed reduction or quotient with declared retained data. \\
			$\Phi_{I_0},\operatorname{Sol}_{\Phi}(I_0)$ & Preservation contract and successful candidate outputs. \\
			\bottomrule
		\end{tabular}
	\end{center}

	\subsection{Inherited Level-I systems and source recovery}
	\label{subsec:standalone-native-data}

	Source recovery asks whether the endpoint, together with its retained markings, determines a specified distinction between formation data.
	We recall the Level-I typing conditions from \cite[Definitions~2.5, 2.9, and~2.14]{Ma2026}.

	\begin{definition}[Native contexts and realized packaging systems]
		\label{def:standalone-realized-system}
		Fix Grothendieck universes $\mathfrak U\in\mathfrak V$, and reserve $\mathcal U$ for a set of edits.
		A \emph{structural domain} $\mathscr C$ is a $\mathfrak U$-small set of possible carriers.		
		For a structural domain $\mathscr C$, let $\mathbf{PkgDat}_{\mathscr C}$ denote the set of five-tuples
		\[
		\mathfrak P=(\Pkg_{\mathfrak P},E_{\mathfrak P},
		\triangleleft_{\mathfrak P},\LL_{\mathfrak P},\RR_{\mathfrak P})
		\]
		with domain type $\mathscr C$, satisfying the following typing conditions:
		\[
		\Pkg_{\mathfrak P}\subseteq\mathscr C,\qquad
		E_{\mathfrak P}\in\mathfrak U,\qquad
		\triangleleft_{\mathfrak P}\subseteq
		E_{\mathfrak P}\times\Pkg_{\mathfrak P}.
		\]
		Here $\Pkg_{\mathfrak P}$ is the collection of packages, $E_{\mathfrak P}$ is the collection of units, and $\triangleleft_{\mathfrak P}$ is their incidence relation.
		For each $K\in\Pkg_{\mathfrak P}$, write
		\[
		I_K=\{e\in E_{\mathfrak P}:e\triangleleft_{\mathfrak P}K\}
		\]
		for its incident units.
		For the carrier $K$ with incident units $I_K$, the endpoint specifies a $\mathfrak U$-small candidate space $\LL_{\mathfrak P}(K)$,	and the selected family satisfies
		$\RR_{\mathfrak P}(K)\subseteq\LL_{\mathfrak P}(K)$.
		No additional structure on the candidate set is required here.		
		\begin{enumerate}[label=\textup{(\roman*)},leftmargin=*]
			\item \textbf{Native context.}
			A \textbf{native context} is a quadruple
			\[
			\Xi=(\mathscr C_\Xi,\operatorname{Car}_\Xi,
			\operatorname{Form}_\Xi,\mathsf J_\Xi),
			\qquad
			\operatorname{Car}_\Xi\subseteq\mathscr C_\Xi,
			\qquad
			\operatorname{Form}_\Xi\in\mathfrak V,
			\]
			where $\mathscr C_\Xi$ is a structural domain, 	$\operatorname{Car}_\Xi$ is the set of admissible carriers, and	$\operatorname{Form}_\Xi$ is the space of formation data.
			Define
			\[
			\mathbf{PkgDat}_\Xi
			=\{\mathfrak P\in\mathbf{PkgDat}_{\mathscr C_\Xi}:
			\Pkg_{\mathfrak P}\subseteq\operatorname{Car}_\Xi\}.
			\]
			The \emph{native realization law} is the specified relation
			\[
			\mathsf J_\Xi\subseteq
			\operatorname{Form}_\Xi\times\mathbf{PkgDat}_\Xi.
			\]
			This relation need not be total or single-valued.
			
			\item \textbf{Realized packaging system.}
			A \textbf{primitive realized packaging system} over a native context $\Xi$ is a triple
			$
			\mathbb P=(\Xi,P,\mathfrak P),
			$
			where $P\in\operatorname{Form}_\Xi$ is the formation datum, $\mathfrak P\in\mathbf{PkgDat}_\Xi$ is the packaged endpoint, and
			$
			(P,\mathfrak P)\in\mathsf J_\Xi.
			$
			This membership is the \textbf{realization judgment}, written
			$
			\Xi\vdash P\rightsquigarrow\mathfrak P.
			$
			In particular, a realized system retains its native context and
			formation datum together with its realized endpoint.
		\end{enumerate}
	\end{definition}

	\begin{definition}[Recoverability of formation data from packaged endpoints]
		\label{def:standalone-source-recovery}
		Fix a native context $\Xi$, a binary relation $\mathcal S_\Xi$ on
		$\operatorname{Form}_\Xi$, and a binary relation $\mathcal E_\Xi$ on
		$\mathbf{PkgDat}_\Xi$.
		The endpoint comparison $\mathcal E_\Xi$ includes the prescribed comparison of any retained markings.
		We say that \textbf{formation data are recoverable relative to
			$(\mathcal S_\Xi,\mathcal E_\Xi)$} if, for every
		$(P,\mathfrak P),(P',\mathfrak P')\in\mathsf J_\Xi$,
		\[
		\mathfrak P\mathrel{\mathcal E_\Xi}\mathfrak P'
		\quad\Longrightarrow\quad
		P\mathrel{\mathcal S_\Xi}P'.
		\]
		For literal equality on both sides this says that each realized endpoint
		has exactly one formation datum.
		A single datum may still realize several endpoints.
	\end{definition}

	This is the endpoint-to-source direction of the comparison in
	\cite[Definitions~2.31--2.32]{Ma2026}. For a specified quantity
	$a:\operatorname{Form}_\Xi\to A$, take $P\mathrel{\mathcal S_\Xi}P'$
	to mean $a(P)=a(P')$. Suppose $\mathcal E_\Xi$ is reflexive on the set
	$Y_\Xi=\{\mathfrak P:\exists P,\ (P,\mathfrak P)\in\mathsf J_\Xi\}$
	of realized endpoints. Then recovery is equivalent to the existence of a
	unique map $d:Y_\Xi\to A$ satisfying
	\[
	d(\mathfrak P)=a(P)\quad ((P,\mathfrak P)\in\mathsf J_\Xi),
	\qquad
	\mathfrak P\mathrel{\mathcal E_\Xi}\mathfrak P'
	\ \Longrightarrow\ d(\mathfrak P)=d(\mathfrak P')
	\quad(\mathfrak P,\mathfrak P'\in Y_\Xi).
	\]
	Indeed, reflexivity and recovery give the same value $a(P)$ for every
	datum realizing a fixed endpoint. This value defines $d$ uniquely, and
	recovery makes it equal on related endpoints. Conversely, any such $d$
	gives the defining implication. For example, with generator names fixed,
	a marked algebra determines the congruence generated by its presentation
	through the kernel of term evaluation \eqref{eq:equational-source-decoder}.
	Different equation lists can generate that same congruence.

	\subsection{Assembly data, explicit axioms, and coherence}
	\label{subsec:level-two-backbone}

	The assembly data specify which observations can be compared and which formations are admitted.
	We first give their types, then state the coherence and rank conditions.

	Let $D$ be a small poset of observation sites, with distinguished objects for package types and additional objects for interfaces.
	The order $U\leq V$ represents restriction from $V$ to $U$.
	Give each site a candidate set $\LL(U)$ and each comparable pair a map $\rho_{V,U}:\LL(V)\to\LL(U)$. Give every formation type $K$ a selected family $\RR(K)\subseteq\LL(K)$.
	At inherited package sites both sets agree with the realized base.
	At a prospective type the candidate set supplies the ambient space and the selected family supplies the operational states.

	A primitive formation signature has the form
	\[
	\alpha:(K_1,\ldots,K_m)\longrightarrow K,
	\qquad m\geq1.
	\]
	The input and output positions are typed separately even when their underlying candidate sets happen to coincide.
	For each signature, specify a family of interface comparisons.
	A comparison between positions of types $V,W$ has a site $U\leq V,W$ and requires equality of their $U$-observations.
	The resulting subset of the product of selected states is denoted
	\[
	\Match_\alpha\subseteq
	\RR(K_1)\times\cdots\times\RR(K_m)\times\RR(K).
	\]
	A comparison needs only the specified common lower site $U$.
	The poset need not have binary meets.

	\begin{definition}[Raw structural packaging presentation]
		\label{def:level-ii-presentation}
		A \textbf{raw structural packaging presentation}, also called a \textbf{raw Level-II packaging system}, consists of the data in \textup{(A0)}--\textup{(A3)} below.
		The numbering (A0)--(A3) distinguishes these additional data clauses from (S1)--(S5), (P1)--(P3), and (Y1)--(Y4) in \cite{Ma2026}.
		\begin{enumerate}[label=\textup{(A\arabic*)},start=0]
			\item \emph{Inherited realized bases:}
			a nonempty indexed family of realized base systems $\mathbb	P_i=(\Xi_i,P_i,\mathfrak P_i)$, each with its own context, realization judgment, and package-side data as above;
			
			\item \emph{Typed observation sites:}
			an observation-site poset $D$, candidate sets $\LL(U)$, restriction maps $\rho_{V,U}$, and selected families at formation types as above.
			Package types are occurrence-labeled when necessary:
			for an original type $(i,K)$, the candidate and selected sets agree with those of $K$ in $\mathfrak P_i$.
			A prospective output type carries declared operational candidates only;
			
			\item \emph{Typed signatures and interfaces:}
			a collection $\mathcal A$ of typed primitive formation signatures, together with
			the chosen interface comparisons for each signature; and
			
			\item \emph{Operational admission:}
			for each $\alpha:(K_1,\ldots,K_m)\to K$, an operational relation
			\[
			\MM_\alpha\subseteq
			\RR(K_1)\times\cdots\times\RR(K_m)\times\RR(K).
			\]
		\end{enumerate}
		We write the resulting presentation as $\mathfrak L=((\mathbb P_i)_{i \in I};D,\LL,\RR,\rho;\mathcal A,\mathcal I,\MM)$, where $\mathcal I$ records the chosen comparisons.
		For one base we suppress the occurrence index when it is unambiguous.
		An admission relation may be empty or multivalued.
	\end{definition}

	The relation $\Match_\alpha$ tests agreement of observations, $\MM_\alpha$ records operational admission, and $\mathsf J_\Xi$ records realization in the native context.

	\subsubsection{Coherence and optional hierarchy}

	Coherent restrictions make each observation independent of the chosen chain of restrictions.
	Interface soundness checks the observations of each admitted tuple.
	The rank condition separately bounds the length of a formation branch.

	\begin{definition}[Coherent structural presentations and hierarchical enrichments]
		\label{def:enrichment}
		A Level-II presentation is \textbf{coherent} when it satisfies \textup{(H1)} and \textup{(H2)} below.
		It is a \textbf{coherent hierarchical enrichment} when it also satisfies \textup{(H3)}.
		The relations $\MM_\alpha$ are those already specified in Definition~\ref{def:level-ii-presentation}.
		\begin{description}[style=nextline,leftmargin=2.8em]
			\item[(H1) Coherent restriction.]
			For $U\leq V\leq W$,
			\[
			\rho_{U,U}=\id,\qquad
			\rho_{W,U}=\rho_{V,U}\circ\rho_{W,V}.
			\]
			
			\item[(H2) Interface soundness.]
			Every admitted formation respects its declared comparisons:
			$\MM_\alpha\subseteq\Match_\alpha$.
			
			\item[(H3) Finite-valued formation rank.]
			There is a map $\ell$ from package types to $\NN$ such that	$\ell(K_i)<\ell(K)$ for every input of every primitive signature $\alpha:(K_1,\ldots,K_m)\to K$.
		\end{description}
		Formation relations may be empty, partial, or multivalued.
		Composite formations are evaluated by the witness semantics of \cref{sec:composition}.
	\end{definition}

	\paragraph{The three levels of assembly axioms.}
	Write $\mathsf{Ax}_{\mathrm{Pkg}}$ for the packaging-system theory of \cite{Ma2026}.
	The three assembly theories satisfy
	\begin{align*}
		\mathsf{Ax}_{\mathrm{II}}^{\mathrm{raw}}
		&=\mathsf{Ax}_{\mathrm{Pkg}}^{\mathrm{indexed\ bases}}+(\mathrm{A1})+(\mathrm{A2})+(\mathrm{A3}),\\
		\mathsf{Ax}_{\mathrm{II}}^{\mathrm{coh}}
		&=\mathsf{Ax}_{\mathrm{II}}^{\mathrm{raw}}+(\mathrm{H1})+(\mathrm{H2}),\\
		\mathsf{Ax}_{\mathrm{II}}^{\mathrm{hier}}
		&=\mathsf{Ax}_{\mathrm{II}}^{\mathrm{coh}}+(\mathrm{H3}).
	\end{align*}
	Condition (A0) supplies the indexed family of Level-I models.
	In the last case, a branch ending at type $K$ has at most $\ell(K)$ primitive steps.
	A finite type-dependency graph admits such a rank exactly when it is acyclic.
	An infinite well-founded type order may require ranks beyond $\NN$.

	Condition (H3) concerns formation edges.
	The overlap graph may contain cycles.
	Interface agreement in (H2) is tested in $\LL(U)$.
	When $U$ and $V$ have selected families, stability under restriction amounts to $\rho_{V,U}(\RR(V))\subseteq\RR(U)$.

	\begin{example}[Two realized bases and a prospective assembly]
		\label{ex:two-base-operational-choice}
		Take two distinct realized Level-I systems $\mathbb P_1,\mathbb P_2$,
		each with one package type $K_i$ whose selected states are $\{0,1\}$.
		Their native contexts can each have a singleton formation-data space
		and exactly the displayed realization judgment. Add a prospective
		output type $K_*$ with provisional selected set $\{0,1\}$, and one
		common observation site $U$ below $K_1,K_2,K_*$, with all candidate
		spaces $\{0,1\}$ and identity restrictions. For
		$\alpha:(K_1,K_2)\to K_*$ the matching relation is
		\[
		\Match_\alpha=\{(0,0,0),(1,1,1)\}.
		\]
		Both choices $\MM_\alpha^{(0)}=\{(0,0,0)\}$ and
		$\MM_\alpha^{(1)}=\Match_\alpha$ satisfy (H1)--(H2); assigning
		ranks zero to $K_1,K_2$ and one to $K_*$ gives (H3).
		The two assemblies have \emph{identical realized input systems},
		sites, and interfaces. They admit different output witnesses.
		Their marked structural presentations are not isomorphic because
		the operational relations have different cardinalities.
		The same rules can also act within one realized base. Take three carriers
		with empty unit set and incidence, binary candidate and selected sets,
		and a native law consisting of the single judgment realizing this endpoint.
		Use the two input carriers and the output carrier in the three positions above.
	\end{example}

	\begin{proposition}[Conservativity over primitive packaging]\label{prop:conservative}
		Every primitive realized packaging system admits an enrichment satisfying
		\emph{(H1)--(H3)}. Forgetting the enrichment recovers exactly the original
		system. Thus every sentence in the original signature that holds in all
		expanded models also holds in all original models.
	\end{proposition}

	\begin{proof}
		Take $D$ to be the discrete poset of its package carriers, retain their
		candidate sets and selected families, supply only identity restrictions, and
		use no primitive formation signatures. Set $\ell=0$ on every type.
		All three requirements hold, and none changes $\Xi$, $P$, $\mathfrak P$, or
		its realization judgment. Empty package families cause no difficulty.
		Every original model therefore expands to a model of the extended theory.
	\end{proof}

	\begin{proposition}[Independence of the structural hypotheses]
		\label{prop:independence}
		Within the raw typed data preceding Definition~\ref{def:enrichment}, none of
		\emph{(H1)}, \emph{(H2)}, or \emph{(H3)} follows from the other two.
	\end{proposition}

	\begin{proof}\leavevmode
		\begin{itemize}
			\item \emph{Failure of (H1).}
			Use $a<b<c$, put $\LL(a)=\LL(b)=\LL(c)=\{0,1\}$, let the two adjacent restrictions be identities, and let $\rho_{c,a}$ exchange $0$ and $1$.
			Use no formations.
			Then (H2) is vacuous and (H3) holds.
			
			\item \emph{Failure of (H2).}
			Take $u<a$ and $u<b$, all three candidate sets equal to	$\{0,1\}$, and both restrictions to $u$ equal to the identity.
			Select all states.
			Give the signature $a\to b$ the relation $\{(0,1)\}$ and declare
			comparison at $u$.
			Rank $a$ below $b$.
			The restrictions are coherent and the hierarchy is well-founded, but the admitted pair disagrees at $u$.
			
			\item \emph{Failure of (H3).}
			Use two discrete package sites with singleton selected sets	and signatures $a\to b$ and $b\to a$, each with its unique state pair.
			There are no interface comparisons.
			Conditions (H1) and (H2) hold.
			A rank would have to increase strictly around the cycle and return to its initial value.
		\end{itemize}
	\end{proof}

	Condition (H1) is precisely functoriality of the candidate sets as a presheaf $D^{\mathrm{op}}\to\mathbf{Set}$ on the observation poset \cite{Stacks}.
	This need not define a presheaf of selected states:
	give $u<v$ the identity restriction on $\{0,1\}$ and select $\RR(v)=\{1\}$,
	$\RR(u)=\{0\}$.
	With no formations, (H1)--(H3) hold.
	On a poset of open sets with the usual covers and all overlap comparisons, a sheaf gluing law would identify the section set of a covered open with its matching families.
	Neither this existence-and-uniqueness condition nor preservation of recovery laws follows from (H1)--(H2).
	Operational admission and native realization remain additional data.

	\subsection{Formation, operational witnesses, and substitution}
	\label{sec:composition}

	Two assignments may agree on all leaves and on the root while using different intermediate states.
	The substitution law (\cref{prop:substitution}) describes the boundary relation of a grafted tree by quantifying over the states shared at each graft.

	Let $T$ be a finite rooted tree. Each vertex $v$ has a package type $K_v$.
	At each nonleaf vertex $v$ with ordered children $v_1,\ldots,v_m$, choose a
	signature $\alpha_v:(K_{v_1},\ldots,K_{v_m})\to K_v$.
	Occurrence labels distinguish vertices, even when their types agree.
	Under (H3), each edge respects the type rank.

	\begin{definition}[Witnesses, boundary relations, and operational realizability]
		\label{def:witness}
		Let $T$ be a finite rooted formation tree with root $r$ and leaf set $L(T)$.
		\begin{itemize}
			\item A \textbf{witness on $T$} is an assignment $x=(x_v)_{v\in V(T)}$ of a selected state $x_v\in\RR(K_v)$ to each vertex such that the prescribed formation rule at every nonleaf vertex admits its ordered input states and its output state.
			Thus, writing $v_1,\ldots,v_{m_v}$ for the ordered children of $v$, the ambient state space $X_T$ and the set $\WW_T$ of all witnesses are
			\begin{align}
				X_T &= \prod_{v\in V(T)}\RR(K_v),\\
				\WW_T &= \bigl\{x\in X_T \mid
				(x_{v_1},\ldots,x_{v_{m_v}},x_v)\in\MM_{\alpha_v}
				\text{ for every nonleaf vertex }v\bigr\}.
			\end{align}
			In particular, the state assigned to an intermediate vertex is both the output of its own formation and an input to its parent's formation.
			
			\item The \textbf{boundary} of $T$ is the vertex set $\partial T=L(T)\cup\{r\}$.
			The \textbf{boundary relation} of $T$ is
			\[
			\{\,x|_{\partial T}\mid x\in\WW_T\,\}.
			\]
			A boundary assignment belongs to this relation exactly when it extends to a witness on the whole tree.
			For prescribed leaf states $\ell$, the possible root states are $\{\,x_r\mid x\in\WW_T,\ x|_{L(T)}=\ell\,\}$.
			If $T$ has only one vertex, the root is also the sole leaf, the boundary has one coordinate, and every selected root state is a witness.
			
			\item The tree is \textbf{operationally realizable} if $\WW_T\ne\varnothing$.
			A native realization of a witness additionally requires a formation	datum $P$ interpreting that witness in the specified context $\Xi$ and a packaged endpoint $\mathfrak P$ satisfying $(P,\mathfrak P)\in\mathsf J_\Xi$.
		\end{itemize}
	\end{definition}

	\begin{proposition}[Boundary relations under substitution]
		\label{prop:substitution}
		\leavevmode
		\begin{enumerate}[label=\textup{(\roman*)},leftmargin=*]
			\item \emph{Boundary relation of a graft.}
			Grafting typed formation trees and then taking their boundary relation gives
			the same relation as composing their boundary relations by existentially
			quantifying the shared states.
			\item \emph{Associativity and units.}
			Repeated substitution is associative.
			The identity relation on a selected state space is a unit when formal identity
			operations are adjoined.
		\end{enumerate}
	\end{proposition}
	
	\begin{proof}\leavevmode
		\begin{enumerate}[label=\textup{(\roman*)},leftmargin=*]
			\item \emph{Boundary relation of a graft.}
			A witness on a grafted tree consists exactly of a witness on the outer tree
			and witnesses on the inserted trees with equal states at identified boundary
			vertices. This is precisely the existential definition of relational
			composition.
			\item \emph{Associativity and units.}
			All bracketings express the existence of the same finite
			collection of internal states satisfying the same conjunction of primitive
			relations. Finite existential quantifiers can be regrouped, proving
			associativity. Equality at an inserted identity vertex forces its two states
			to coincide and changes no witness condition.
		\end{enumerate}
	\end{proof}
	
	The rank condition (H3) applies to primitive formations. Associativity compares bracketings of a fixed substitution. A different formation tree may define a different boundary relation.
	
	A rank can order successive formations even when their required intermediate states are incompatible.	
	
	\begin{example}[Identifying intermediate states can create a boundary pair]
		\label{ex:chain}
		Consider two successive formations $K_0\to K_1\to K_2$ with
		selected state sets
		\[
		\RR(K_0)=\{a\},\qquad
		\RR(K_1)=\{0,1\},\qquad
		\RR(K_2)=\{c\}.
		\]
		Use only identity restrictions, impose no interface comparisons, and
		give $K_i$ rank $i$. The formation relations
		\[
		\MM_{01}=\{(a,0)\},\qquad
		\MM_{12}=\{(1,c)\}
		\]
		satisfy \textup{(H1)--(H3)}. Each step has a witness. A witness for
		the chain, however, must assign a single state to the intermediate
		vertex. The first step requires that state to be $0$, and the second
		requires it to be $1$. Thus the boundary relation of the chain is empty.
		
		Now identify $0$ and $1$ with a single state $*$, keeping the endpoint
		states $a$ and $c$ fixed. The images of the two formation relations
		are $\{(a,*)\}$ and $\{(*,c)\}$. Their composite contains $(a,c)$,
		witnessed by $(a,*,c)$. This quotient witness has no lift to the
		original chain.
	\end{example}
	
	\subsection{Observations, comparisons, and preservation contracts}
	\label{subsec:typed-native-observations}
	
	An assembly can cover its inputs and still change what their observations determine. The following definitions specify the observation maps and the conditions an admissible comparison must preserve.
	
	\begin{definition}[Observation presentation over realized inputs]
		\label{def:associated-observation-presentation}
		Fix a family of realized packaging systems $(\mathbb P_i)_{i\in I}$.
		An \textbf{observation presentation} over these inputs consists of the following data:
		\begin{enumerate}
			\item \emph{Sites and restrictions.} An observation-site poset $D^{\mathrm{obs}}$, an ambient set $O(U)$ at each site, and a restriction map		$\rho^{\mathrm{obs}}_{V,U}:O(V)\to O(U)$ for each $U\leq V$.
			Each shared site has a specified readout interpretation.
			
			\item \emph{Selected observations and formations.} A selected set $S(K)\subseteq O(K)$ at each formation type $K$ and, for each signature $\alpha:(K_1,\ldots,K_m)\to K$,	a relation of admitted observation tuples
			\[
			N_\alpha\subseteq
			S(K_1)\times\cdots\times S(K_m)\times S(K).
			\]
			
			\item \emph{Native interpretation.} At each input package site $K$, a specified bijection
			\[
			\iota_K:O(K)\longrightarrow
			\operatorname{Obs}(\mathbb P_{i(K)},K;\mu_K)
			\]
			identifying $S(K)$ with the selected native observations.
			Here $\operatorname{Obs}$ is a specified construction from the native source or endpoint with retained markings $\mu_K$.
			A record type may use a disjoint union of such observation sets indexed by the named bases, with the base tag retained.
			Each prospective output type has a specified observation type and native context.
		\end{enumerate}
		The restriction maps and admission relations satisfy the diagram equations and formation conditions declared for the presentation.
		Restrictions act on the ambient sets. Preservation of selected sets under restriction is an additional condition when required.
	\end{definition}
	
	\paragraph{Native realization of observation data.}
	The witness set and boundary relation of \cref{def:witness} are formed using $S(K)$ and $N_\alpha$. Let $Z$ be the specified set of witnesses or boundary tuples to be realized.
	A native realization consists of a judgment $\Xi'\vdash P'\rightsquigarrow\mathfrak P'$ in the prescribed
	output context, a selected family $S'$ of observations of that output, and a map $h:S'\to Z$. For each retained readout, its output and target maps must satisfy
	\[
	r_{\mathrm{out}}=r_Z\circ h
	\]
	after the value identifications specified by the contract. The contract also specifies any injectivity or surjectivity required of $h$ (\cref{subsec:framework-contract}).
	
	Taking $O(U)=\LL(U)$, $S(K)=\RR(K)$, $\rho^{\mathrm{obs}}=\rho$, and $N_\alpha=\MM_\alpha$ recovers the selected-state semantics of Definitions~\ref{def:level-ii-presentation} and~\ref{def:witness}.
	Other choices require the specified native interpretation.
	In the finite spectral and algebraic models below all observation states are selected, so $S(K)=O(K)$.
	
	\subsubsection{Structural comparisons and retained observations}
	\label{subsec:structural-comparisons}
	
	A structural isomorphism preserves the specified data over
	fixed realized bases.
	Comparing different bases requires a specified invertible map that preserves the native law.
	The recovery relations $\mathcal S_\Xi$ and $\mathcal E_\Xi$ of
	Definition~\ref{def:standalone-source-recovery} define quotients when they are equivalence relations compatible with the relevant structure.
	
	\begin{definition}[Structural isomorphism over native bases]
		\label{def:structural-isomorphism}
		Let $\mathfrak L$ and $\mathfrak L'$ be structural presentations with base index sets $I,I'$, observation-site posets $D,D'$, and signature
		sets $\mathcal A,\mathcal A'$.
		Fix the permitted invertible transports	of native bases.
		A \textbf{structural isomorphism over the native bases}
		is a tuple
		\[
		\Phi=\bigl(\beta,(\tau_i)_{i\in I},\delta,(h_U)_{U\in D},\sigma\bigr)
		\colon \mathfrak L \xrightarrow{\ \sim\ } \mathfrak L'
		\]
		satisfying the following conditions:		
		\begin{enumerate}[label=\textup{(\roman*)}]
			\item \emph{Realized bases.}
			The map $\beta:I\to I'$ is a bijection, and
			$\tau_i:\mathbb P_i\to\mathbb P'_{\beta(i)}$ is a permitted invertible transport of the native context, formation datum, and marked endpoint.
			Write $\kappa_i$, $u_i$, and $t_{i,K}$ for its bijections on package carriers, units, and the candidate states of package $K$, respectively.
			They satisfy
			\[
			e\mathrel{\triangleleft_i}K
			\iff
			u_i(e)\mathrel{\triangleleft'_{\beta(i)}}\kappa_i(K),
			\qquad
			t_{i,K}\bigl(\RR_i(K)\bigr)
			=\RR'_{\beta(i)}\bigl(\kappa_i(K)\bigr).
			\]
			Over pointwise fixed indexed bases, $\beta$ and all $\tau_i$ are identities.
			
			\item \emph{Sites and restrictions.}
			The map $\delta:D\to D'$ is an order isomorphism carrying formation types onto formation types.
			Each $h_U:\LL(U)\to\LL'(\delta U)$ is a bijection.
			At an inherited package	site $(i,K)$ these maps agree with the native transport:
			\[
			\delta(i,K)=(\beta(i),\kappa_i(K)),
			\qquad h_{(i,K)}=t_{i,K}.
			\]
			For every formation type $K$ and every pair of sites $U\le V$, respectively,
			\[
			h_K\bigl(\RR(K)\bigr)=\RR'(\delta K),
			\qquad
			h_U\circ\rho_{V,U}
			=\rho'_{\delta V,\delta U}\circ h_V .
			\]
			
			\item \emph{Formation rules.}
			The map $\sigma:\mathcal A\to\mathcal A'$ is a bijection.
			For each $\alpha:(K_1,\ldots,K_m)\to K$, its image has signature
			\[
			\sigma(\alpha):(\delta K_1,\ldots,\delta K_m)\to\delta K.
			\]
			Let $\mathcal I_\alpha$ denote the declared interface comparisons of $\alpha$.
			Define $\delta_*\mathcal I_\alpha$ by replacing each
			comparison $(p,q,U)$ with $(p,q,\delta U)$, where $p,q$ are its two argument positions, including the output position when applicable.
			The comparisons and admitted tuples satisfy
			\[
			\mathcal I'_{\sigma(\alpha)}=\delta_*\mathcal I_\alpha,
			\qquad
			\MM'_{\sigma(\alpha)}
			=(h_{K_1}\times\cdots\times h_{K_m}\times h_K)(\MM_\alpha).
			\]
		\end{enumerate}		
		If the contract retains a readout $r:X\to Z$, with corresponding readout $r':X'\to Z'$, the induced bijection $\Phi_X:X\to X'$ must satisfy $r'\circ\Phi_X=\lambda\circ r$, where $\lambda:Z\to Z'$ is the prescribed relabeling and is the identity for fixed value labels.
		A retained cost function satisfies $c'(\Phi_X(x))=c(x)$ on its domain.		
		We write $\mathfrak L\cong\mathfrak L'$ over the specified native bases if such a tuple $\Phi$ exists.
	\end{definition}
	
	The native transports used here have concrete meanings. In
	$\Xi_{\mathcal V}^{\mathrm{pres}}$ (the presentation context defined in \cref{def:equational-native-context}), a marked-algebra isomorphism
	transports states, units, operations, and generator markings.
	Compare presentations on the same named generators by equality of
	their generated congruences. A specified relabeling of generators first
	identifies the free algebras. The evaluation kernels in
	\eqref{eq:equational-source-decoder} verify preservation of the native
	judgments. In $\Xi_{\mathcal H}^{\mathrm{com}}$ (the access--commutant context of \eqref{eq:standalone-commutant-law}), a unitary
	$U:\mathcal H\to\mathcal H'$ transports a source $M$, its endpoint
	carrier $M'$, and its retained readout algebras by conjugation:
	\[
	(UMU^*)'=UM'U^*,\qquad (UM'U^*)'=UMU^*.
	\]
	Conjugation gives the unit bijection and transports refinements of
	the underlying carrier set. Thus it preserves the original candidate
	and selected families as well as the native law. The decoded source is
	transported with its representation; comparison across ambient spaces
	uses that declared source comparison. Identity transports,
	their inverses and their composites are again transports of these
	types.

	\begin{proposition}[Invariance under structural isomorphism]
		\label{prop:structural-isomorphism-invariance}
		\leavevmode
		\begin{enumerate}[label=\textup{(\roman*)},leftmargin=*]
			\item \emph{Structural equivalence and axioms.}
			Structural equivalence is an equivalence relation in each of the declared native comparison classes.			
			It preserves \emph{(H1)--(H3)} when imposed.
			
			\item \emph{Tree witnesses and boundaries.}
			Every finite typed formation tree is carried to one whose witness set and boundary relation are bijective to the originals.
			
			\item \emph{Readouts and native judgments.}
			Retained readout values, their determination and non-determination judgments, and the specified native realization judgments are preserved.
		\end{enumerate}
	\end{proposition}

	\begin{proof}
		\leavevmode
		\begin{enumerate}[label=\textup{(\roman*)},leftmargin=*]
			\item \emph{Structural equivalence and axioms.}
			Identity maps, inverse maps, and composites satisfy all displayed
			conditions.
			The commuting restriction equations transport (H1);
			the comparison sites and bijections transport $\Match_\alpha$ and
			therefore (H2).
			A rank for (H3) is transported along the type bijection.
			
			\item \emph{Tree witnesses and boundaries.}
			For a tree, send a state assignment $(x_v)_v$ to
			$(h_{K_v}(x_v))_v$.
			The selected-family equalities give a bijection	on the ambient products of selected states, and exact transport of each formation relation restricts it to the witness sets.
			It commutes with every boundary	projection.
			
			\item \emph{Readouts and native judgments.}
			A pair has equal specified readouts before transport exactly when its transported pair does, so it has the same recovery and non-recovery judgments.
			Native preservation follows from the chosen base transports, independently of the witness calculation.
		\end{enumerate}
	\end{proof}

	\paragraph{Comparisons over a fixed atlas.}
	For fixed chart state spaces $X_i$ and separator readouts, compare sources $(H,(\pi_i))$ and $(H',(\pi_i'))$ by maps $h:H\to H'$ satisfying $\pi_i'h=\pi_i$ for every $i$.
	These are set maps or homomorphisms in the specified algebraic class. 
	An invertible comparison is a \emph{source isomorphism over the atlas}.
	A quotient comparison additionally requires surjectivity.
	Native formation data and realization judgments remain separately specified.
	The contract specifies any presentation lists or generator markings that the comparison must retain.

	\paragraph{Observation reductions over retained bases.}
	For the associated observation presentations defined in
	Definition~\ref{def:associated-observation-presentation}, write $O(U)$ for the state sets, $\rho^{\rm obs}_{V,U}$ for restrictions, and $N_\alpha$ for relations.
	A local image reduction on a fixed site	and signature diagram consists of surjections	$q_U:O(U)\twoheadrightarrow\overline O(U)$ satisfying
	\[
	q_U\rho^{\rm obs}_{V,U}
	=\overline\rho^{\rm obs}_{V,U}q_V,
	\qquad
	\overline N_\alpha=(\textstyle\prod_{p\in\alpha}q_p)(N_\alpha),
	\]
	where the product runs over all input and output positions, and
	$q_K(S(K))=\overline S(K)$ at every formation type. Here $S(K)$ and
	$\overline S(K)$ are the selected observation families. The relation
	images use the restrictions of the $q_K$ to those selected families.
	Native bases, their formation data and their selected families remain
	recorded. At an input package site, the reduced interpretation is the
	quotient of the original observation construction induced by
	$q_K\circ\iota_K^{-1}$. This coding map is part of the reduced
	observation data. The reduced
	sets are observation codes; replacing an original
	native selected family by these codes would additionally require an
	output native context and a realization judgment.

	For a finite formation tree $T$, let $\mathsf W_T^{\rm obs}$ and
	$\overline{\mathsf W}_T^{\rm obs}$ be the two witness sets, and put
	$q_T=\prod_{v\in V(T)}q_{K_v}$ on the product of selected observation
	families. Local image construction always gives
	\[
	q_T(\mathsf W_T^{\rm obs})
	\subseteq\overline{\mathsf W}_T^{\rm obs}.
	\]
	Call the reduction \emph{witness-exact on $T$} when equality holds.
	Its boundary image then also equals the reduced boundary relation, by commuting the coordinate projections with $q_T$.
	Example~\ref{ex:chain} exhibits the obstruction at a shared variable when the local lifts select different representatives.

	\begin{definition}[Retained-observation reduction and non-collapse]
		\label{def:retained-observation-reduction}
		Let $X$ be a specified state, witness, or source set, and let $f_j:X\to D_j$ ($j\in J$) be the declared retained observations.
		\begin{itemize}
			\item A surjection $q:X\twoheadrightarrow\overline X$ is a \textbf{retained-observation reduction} when there are maps
			$\overline f_j:\overline X\to D_j$ with
			$f_j=\overline f_jq$.
			For this family, non-collapse means that $q$ identifies states only when their retained observations agree.
			The contract may also require the map to preserve restrictions, substitutions, or realized outputs.
			
			\item \textbf{Non-collapse} of a whole presentation means preservation of every clause in its contract.
			Unobserved data may be identified when the contract permits it.
		\end{itemize}
	\end{definition}

	\begin{proposition}[Readout-preserving reductions]
		\label{prop:observation-noncollapse}
		Let $(f_j:X\to O_j)_{j\in J}$ be a family of readouts, and let
		$q:X\twoheadrightarrow\overline X$ be a surjection.
		\begin{enumerate}[label=\textup{(\roman*)}]
			\item \emph{Separation of states.}
			The map $q$ keeps every pair of states distinguished by a readout
			distinct if and only if each $f_j$ descends to a unique map
			$\overline f_j:\overline X\to O_j$ with	$f_j=\overline f_j\circ q$.
			The exact criterion is
			\[
			\ker q\subseteq\bigcap_{j\in J}\ker f_j.
			\]
			In particular, if the readouts separate points of $X$, then every
			readout-preserving reduction is bijective.
			
			\item \emph{Functional determination.}
			Suppose the readouts descend along $q$.
			For every $S\subseteq J$ and $j\in J$, the family $(f_s)_{s\in S}$ determines $f_j$ on $X$ if and only if $(\overline f_s)_{s\in S}$ determines $\overline f_j$ on $\overline X$.
		\end{enumerate}
	\end{proposition}

	\begin{proof}
		\begin{enumerate}[label=\textup{(\roman*)}]
			\item \emph{Separation of states.} Keeping readout-distinguished states distinct means that
			$q(x)=q(y)$ implies $f_j(x)=f_j(y)$ for every $j\in J$.
			This is exactly the kernel inclusion.
			Under this condition, $\overline f_j(q(x))=f_j(x)$ is well defined, and surjectivity of $q$ makes $\overline f_j$ unique.
			Conversely, the factorization $f_j=\overline f_j\circ q$ makes $f_j$ constant on each fiber of $q$.
			If the readouts separate points, each fiber is a singleton.
			
			\item \emph{Functional determination.} A pair witnesses failure of determination precisely when it
			agrees on every readout in $S$ and differs on readout $j$.
			The factorizations in \textup{(i)} carry each such pair to a
			witness in $\overline X$.
			Conversely, surjectivity provides preimages of both states of any witness in $\overline X$, and those preimages have the same respective readout values.
			Thus a witness exists on one space exactly when it exists on the other.
		\end{enumerate}
	\end{proof}

	The joint image gives the usual image factorization in sets.
	For homomorphic readouts of one algebraic signature, the first isomorphism
	theorem identifies $X/\ker F$ with $F(X)$, where $F=(f_j)_{j\in J}$
	\cite{BurrisSankappanavar1981}.
	The following proposition gives the	factorization for the declared readouts.

	\begin{proposition}[Canonical observable quotient]
		\label{thm:universal-observable-quotient}
		Let $X$ be a set with a declared family of readouts	$f_j:X\to D_j$, indexed by $j\in J$, and put
		\[
		F=(f_j)_{j\in J}:X\longrightarrow
		\prod_{j\in J}D_j,
		\qquad X_{\mathrm{obs}}=F(X).
		\]
		\begin{enumerate}[label=\textup{(\roman*)},leftmargin=*]
			\item \emph{Universal property.}
			$F:X\twoheadrightarrow X_{\mathrm{obs}}$ is the coarsest surjective reduction preserving all the declared readouts.		
			More precisely, every readout-preserving surjection
			$q:X\twoheadrightarrow Y$ admits a unique surjection
			$p:Y\twoheadrightarrow X_{\mathrm{obs}}$ satisfying
			$
			F=pq.
			$
			
			\item \emph{Determination and non-determination.}
			The observable quotient preserves every determination and non-determination judgment among the $f_j$.
			
			\item \emph{Irreducibility and minimum cardinality.}
			This quotient admits no further proper readout-preserving surjective reduction, and every irreducible such quotient is uniquely isomorphic to $X_{\mathrm{obs}}$ over $X$.			
			For finite $X$, it therefore has minimum cardinality among these quotients, with equality precisely for quotients isomorphic to $X_{\mathrm{obs}}$ over $X$.
		\end{enumerate}
	\end{proposition}

	\begin{proof}
		Write $\pi_j:X_{\mathrm{obs}}\to D_j$ for the restriction of the $j$th coordinate projection.
		Since $f_j=\pi_j\circ F$ for every $j$,	the surjection $F:X\twoheadrightarrow X_{\mathrm{obs}}$ preserves all the declared readouts.		
		\begin{enumerate}[label=\textup{(\roman*)},leftmargin=*]
			\item \emph{Universal property.}
			Let $q:X\twoheadrightarrow Y$ be any readout-preserving surjection,
			and write $g_j:Y\to D_j$ for its descended readouts, so that
			$f_j=g_j\circ q$.
			Define
			\[
			p(y)=(g_j(y))_{j\in J}.
			\]
			For each $y\in Y$, choose $x\in X$ with $q(x)=y$.
			Then $p(y)=(f_j(x))_{j\in J}=F(x)$, so $p(y)\in X_{\mathrm{obs}}$.
			Thus $p\circ q=F$, and $p$ is surjective because $F$ is surjective.
			This identity determines $p$ uniquely because $q$ is surjective.
			Moreover,
			\[
			\pi_j\circ p=g_j \qquad (j\in J),
			\]
			so $p$ itself preserves all the readouts.
			Hence every readout-preserving quotient reduces further to $X_{\mathrm{obs}}$, which proves coarseness.
			
			\item \emph{Determination and non-determination.}
			Fix $S\subseteq J$ and $j\in J$.
			A pair $x,x'\in X$ witnesses failure of determination of $f_j$ by $(f_s)_{s\in S}$ exactly when $F(x),F(x')$ agree in every coordinate in $S$ and differ in coordinate $j$.
			Every pair of points in $X_{\mathrm{obs}}$ has preimages in $X$.
			Thus a witness exists on one space exactly when	it exists on the other.
			
			\item \emph{Irreducibility and minimum cardinality.}
			Let $r:X_{\mathrm{obs}}\twoheadrightarrow Z$ preserve the coordinate
			readouts. If $r(u)=r(v)$, then $\pi_j(u)=\pi_j(v)$ for every $j$.
			Since $u$ and $v$ are tuples, $u=v$.
			Thus $r$ is bijective, proving that $X_{\mathrm{obs}}$ admits no proper readout-preserving reduction.
			
			Now suppose $q:X\twoheadrightarrow Y$ is an irreducible quotient.
			Its comparison map $p:Y\twoheadrightarrow X_{\mathrm{obs}}$ preserves
			the readouts, so irreducibility makes $p$ bijective.
			The identities $\pi_j\circ p=g_j$ show that its inverse also preserves the readouts.
			It is the unique isomorphism over $X$ because $p\circ q=F$.
			
			Finally, if $X$ is finite, the surjection $p$ gives
			$|Y|\geq |X_{\mathrm{obs}}|$, with equality exactly when $p$ is
			bijective. The minimum is attained by $F$ itself.
		\end{enumerate}
	\end{proof}

	\begin{example}[Redundant hidden states and indispensable readouts]
		\label{ex:observable-quotient}
		Let $X=\{0,1\}^2\times\{r,s\}$ and retain
		$f_1(a,b,t)=a$ and $f_2(a,b,t)=b$. The observable quotient
		$F(a,b,t)=(a,b)$ reduces eight states to four while retaining
		all dependencies and nondependencies among these readouts. It is
		irreducible under the stated set-level contract. If the third
		coordinate $t$ is also retained, the canonical quotient has eight
		states. A protected native operation using $t$ descends to the
		quotient that forgets $t$ only if the quotient equivalence is
		compatible with the operation.
	\end{example}
	
	\subsubsection{Structural problems and preservation contracts}
	\label{subsec:framework-contract}
	
	A \textbf{structural problem} specifies the following data:
	\begin{enumerate}[label=(\roman*)]
		\item \emph{Input class.} The supplied realized bases and structural,
		native, or observational data, with the hypotheses defining the class
		under study. When states are given by observations, the input
		includes the native source and the construction of each observation.
		
		\item \emph{Allowed constructions.} The witnesses, common sources,
		native formation data, allowed edits, source summaries, or interface
		maps that may be constructed, together with their types and the maps
		they must carry to the original data. The allowed changes to carriers,
		ambient representations, markings, and observations are stated here.
		
		\item \emph{Preservation contract.} The observations and properties
		required of the output, with their languages, target values and
		comparisons. The targets, such as the least recovery closure
		generated by the input's local rules, are fixed by the input.
		
		\item \emph{Resource objective, when present.} The admissible output
		class, the cost to be bounded or minimized, and the comparison under
		which any equality or uniqueness claim is made.
	\end{enumerate}
	
	For an input $I_0$, write $\operatorname{Cand}(I_0)$ for its well-typed candidate outputs and $\Phi_{I_0}$ for its preservation predicate.
	The solutions are
	\[
	\operatorname{Sol}_{\Phi}(I_0)
	=\{z\in\operatorname{Cand}(I_0):\Phi_{I_0}(z)\}.
	\]
	A native assembly problem fixes its realization law.
	A contract may allow the choice of ambient Hilbert space within the fixed
	access--commutant family of \eqref{eq:standalone-commutant-law}.
	
	\paragraph{Preservation conditions.}
	Each clause of $\Phi_{I_0}$ specifies its inputs and outputs, the objects being compared, the allowed maps, and the quantifiers.
	\begin{center}
		\small
		\begin{tabular}{@{}>{\raggedright\arraybackslash}p{0.27\textwidth}>{\raggedright\arraybackslash}p{0.66\textwidth}@{}}
			\toprule
			Clause & Required output or retained information \\
			\midrule
			Operational realization & A common selected-state witness for the
			specified formation instance. \\
			Witness-exact reduction & Every quotient witness on the stated diagrams lifts to a
			witness in the original diagram, as in $q_T(W_T)=\overline W_T$. \\
			Original-chart coverage & Prescribed images on each complete original
			chart; in a state-exact atlas $\pi_i(W)=X_i$. \\
			Assembled-side coverage & Prescribed images on each full side of a named cut, in addition to
			any prescribed chart images. \\
			Recovery faithfulness & The target descriptor determination and
			non-determination judgments, expressed by $c_W=c_0$ for target closure
			$c_0$. Protected cross-chart pairs are included as separate
			conditions when required. \\
			Native realization & A datum $P^z$ and endpoint $\mathfrak P^z$ with
			$(P^z,\mathfrak P^z)\in\mathsf J_{\Xi^z}$ and all named markings
			and comparison maps. \\
			Separator and cut exactness & Prescribed separator images or embeddings,
			represented intersections where required, and the specified full-cut
			information under the theorem's source category. \\
			Source-observation descent & Factorization of observation semantics
			through the source map, including any required continuation or
			substitution compatibility. \\
			Exact interface answers & For $q:X\twoheadrightarrow Y$ and a declared
			constraint language, obtain $q(\bigcap\Gamma)$ for every permitted
			finite constraint instance $\Gamma$. \\
			Truth preservation & A named observation $f:X\to O$ factors as
			$f=\bar f\circ q$; existential feasibility alone is weaker. \\
			\bottomrule
		\end{tabular}
	\end{center}
	
	\begin{definition}[Operational, contract, and faithful native completeness]
		\label{def:relative-completion}
		Fix a Level-II presentation $\mathfrak L$.		
		\begin{enumerate}[label=\textup{(\roman*)}]
			\item \textbf{Operational completeness at specified boundaries.}
			Fix a class $\mathscr T$ of finite typed formation trees.
			For each $T\in\mathscr T$, choose a boundary $B_T\subseteq V(T)$ and write $b_T(x)=x|_{B_T}$ for the boundary projection on $X_T$.
			Let $M_T$ be the set of matching assignments
			\[
			M_T=\bigl\{x\in X_T \mid
			(x_{v_1},\ldots,x_{v_{m_v}},x_v)\in\Match_{\alpha_v}
			\text{ at every internal vertex }v\bigr\},
			\]
			where $v_1,\ldots,v_{m_v}$ are the ordered children of $v$.			
			The presentation is \textbf{operationally complete} for these trees
			and boundaries if
			\[
			b_T(\WW_T)=b_T(M_T)
			\qquad\text{for every }T\in\mathscr T.
			\]
			Under \textup{(H2)}, $\WW_T\subseteq M_T$, so the condition is
			equivalent to the following witness requirement:
			\[
			\forall T\in\mathscr T,\quad
			\forall x\in M_T,\quad
			\exists w\in\WW_T
			\quad\text{such that}\quad w|_{B_T}=x|_{B_T}.
			\]
			The witness $w$ may differ from $x$ at interior vertices.
			
			\item \textbf{Contract completeness and faithful native completeness.}
			Fix a class $\mathscr D$ of structural problems.
			For each $I_0\in\mathscr D$, specify the candidate outputs, native contexts and realization laws, comparison maps, and preservation predicate $\Phi_{I_0}$.
			Let $\operatorname{Sol}_{\Phi}(I_0)$ be the set of candidates satisfying $\Phi_{I_0}$.			
			The presentation is \textbf{contract-complete} on $\mathscr D$ for $\Phi$ if
			\[
			\operatorname{Sol}_{\Phi}(I_0)\ne\varnothing
			\qquad\text{for every }I_0\in\mathscr D.
			\]
			When $\Phi_{I_0}$ requires the same candidate to have a realization
			in the prescribed native context and to satisfy all stipulated
			faithfulness conditions, we call this property
			\textbf{faithful native completeness}.
		\end{enumerate}
	\end{definition}

	For the two admissions in \cref{ex:two-base-operational-choice}, take the boundary to consist of the leaves and the root.
	The second rule admits every matching boundary tuple, while the first excludes $(1,1,1)$.

	\paragraph{Cut terminology.}
	For a finite matching graph a \textbf{complete-boundary cut} is a bipartition of its vertex set into two nonempty parts, with the joint readout of all crossing-edge observations, including their actual correlations.
	We use ``actual cut'' for such a cut of the specified diagram, rather than a hypothetical independent boundary.
	A \textbf{connected-side cut} additionally requires both induced subgraphs to be connected.
	A \textbf{protected} or \textbf{named} cut is one included in the contract.
	Full range requires the candidate's boundary image to equal the stipulated ambient boundary image.
	Represented common-information exactness and equality with the fiber product of actual side images agree under the appropriate permutability hypotheses, not for arbitrary relations.
	In a tree, a \textbf{full tree cut} means the bipartition from removing one edge, with its entire separator retained.

	Definition~\ref{def:preservation-contract} combines chart coverage, recovery faithfulness, exact cut information, and native realization in the full-cut native contract.

	\subsection{Relational representation and relative non-collapse}
	\label{subsec:relational-representation}

	The multisorted encoding records the complete marked presentation.
	Finite-tree witnesses then give primitive-positive queries
	\cite{ChandraMerlin1977}.

	\begin{definition}[Well-formed marked typed relational encoding]
		\label{def:rel-encoding}
		Fix a small occurrence-labeled skeleton $\Sigma$ with base index set
		$I$, observation-site poset $D$, formation types, ordered primitive
		signatures, and position-labeled interface comparisons.
		A \textbf{well-formed marked typed relational encoding over $\Sigma$} is a
		many-sorted relational structure
		\[
		\mathfrak R=
		\bigl((\mathcal B_i)_i,\ (L_U)_U,\ (S_K)_K,\
		(G_{V,U})_{U\leq V},\ (N_\alpha)_\alpha,\
		(H_{i,K})_{(i,K)\,\mathrm{inherited}}\bigr)
		\]
		satisfying the following conditions:		
		\begin{enumerate}[label=\textup{(\roman*)}]
			\item \textbf{Marked native bases.}
			At each $i\in I$, the native table $\mathcal B_i$ has a structural
			domain $\mathscr C_i$, an admissible-carrier predicate
			$\operatorname{Car}_i\subseteq\mathscr C_i$, a sort $F_i$ of
			formation data, and a sort $Q_i$ of endpoint records.
			For each $e\in Q_i$, its membership and incidence relations specify
			the five-tuple
			\[
			\mathfrak P(e)=
			(\Pkg_e,E_e,\triangleleft_e,\LL_e,\RR_e).
			\]
			The map $e\mapsto\mathfrak P(e)$ must be a bijection from $Q_i$
			onto
			\[
			\bigl\{\mathfrak P\in\mathbf{PkgDat}_{\mathscr C_i}
			\mid \Pkg_{\mathfrak P}\subseteq\operatorname{Car}_i\bigr\}.
			\]
			The table also contains a judgment relation $J_i\subseteq F_i\times Q_i$
			and distinguished elements $p_i\in F_i$, $e_i\in Q_i$ satisfying
			$(p_i,e_i)\in J_i$.
			Any additional native data named by the	prescribed transports are included with their relation symbols and function graphs.
			
			\item \textbf{Sites and formation rules.}
			Each site $U$ has a sort $L_U$, and each formation type $K$ has a
			selected-state predicate $S_K\subseteq L_K$.
			For $U\leq V$, the relation $G_{V,U}\subseteq L_V\times L_U$ satisfies
			\[
			\forall x\in L_V\quad \exists!\,y\in L_U
			\quad (x,y)\in G_{V,U}.
			\]
			For each signature $\alpha:(K_1,\ldots,K_m)\to K$, the admission
			relation has type
			\[
			N_\alpha\subseteq
			S_{K_1}\times\cdots\times S_{K_m}\times S_K.
			\]
			Its argument positions and their interface comparisons are those
			specified by $\Sigma$.
			
			\item \textbf{Identification of inherited sites.}
			For each inherited site $(i,K)$, the relation $H_{i,K}$ is the
			graph of a named bijection
			\[
			h_{i,K}:L_{(i,K)}\xrightarrow{\ \sim\ }\LL_{e_i}(K),
			\qquad
			h_{i,K}\bigl(S_{(i,K)}\bigr)=\RR_{e_i}(K).
			\]
		\end{enumerate}		
		Write $\mathsf{Rel}^{\mathrm{wf}}_\Sigma$ for the groupoid of these
		encodings.
		Its morphisms are sortwise bijections fixing $\Sigma$,
		carrying each named relation onto its counterpart, and preserving
		the distinguished elements.
		On each native base they must induce a permitted invertible transport.
		Over fixed native bases their native components are the identity.
	\end{definition}
	
	\begin{proposition}[Relational representation and witness semantics]
		\label{prop:rel-representation}
		Fix an occurrence-labeled skeleton $\Sigma$.
		Let $\mathsf{Pres}_\Sigma$ be the groupoid of Level-II presentations and structural isomorphisms over $\Sigma$.
		Use the same permitted native transports in	$\mathsf{Pres}_\Sigma$ and $\mathsf{Rel}^{\mathrm{wf}}_\Sigma$.		
		\begin{enumerate}[label=\textup{(\roman*)}]
			\item \emph{Representation of marked presentations.}
			The encoding of \cref{def:rel-encoding} defines an equivalence
			of groupoids
			\[
			\operatorname{Enc}_\Sigma:
			\mathsf{Pres}_\Sigma \longrightarrow
			\mathsf{Rel}^{\mathrm{wf}}_\Sigma.
			\]
			In particular, every mark-preserving isomorphism between two
			encodings is induced by a unique structural isomorphism between
			the corresponding presentations.
			
			\item \emph{Witnesses and boundary relations.}
			Let $T$ be a finite typed formation tree.
			Give each vertex $v$ a variable $x_v$ of sort $L_{K_v}$.
			For each internal vertex $v$, write $v_1,\ldots,v_{m_v}$ for its ordered children.
			In the encoding, the witness set $\WW_T$ is defined by the formula
			\[
			\varphi_T((x_v)_{v\in V(T)})
			:=
			\bigwedge_{v\in V(T)} S_{K_v}(x_v)
			\;\land\!
			\bigwedge_{v\text{ internal}}
			N_{\alpha_v}(x_{v_1},\ldots,x_{v_{m_v}},x_v).
			\]
			For any designated boundary $B\subseteq V(T)$, the boundary
			relation is defined by
			\[
			\varphi_{T,B}((x_v)_{v\in B})
			:=
			\exists (x_v)_{v\in V(T)\setminus B}\;
			\varphi_T((x_v)_{v\in V(T)}).
			\]
			Agreement of states $x\in L_V$ and $y\in L_W$ at a common
			site $U\leq V,W$ is defined by
			\[
			\exists z\;\bigl(G_{V,U}(x,z)\land G_{W,U}(y,z)\bigr),
			\qquad z\text{ of sort }L_U.
			\]
			All three formulas are primitive positive.
		\end{enumerate}
	\end{proposition}

	\begin{proof}
		\leavevmode
		\begin{enumerate}[label=\textup{(\roman*)}]
			\item \emph{Representation of marked presentations.} To decode a well-formed encoding, read each endpoint record
			as its packaged five-tuple. The marked datum and endpoint satisfy
			the recorded native judgment, so they give a realized base at
			each index. Each restriction graph defines a unique total map.
			The named bijection at each inherited site identifies its ambient
			and selected states with those of the marked endpoint. Finally,
			the skeleton specifies the ordered formation signatures and
			comparison sites, and each $N_\alpha$ supplies the corresponding
			admission relation. These data give a Level-II presentation.
			Encoding it again recovers the given structure up to the named
			site identifications and sort tags. Thus the encoding is
			essentially surjective.
			
			A structural isomorphism induces a relational isomorphism by its
			action on the native data and site sorts. Conversely, a
			mark-preserving relational isomorphism gives the permitted native
			transports and site bijections. Preservation of each restriction
			graph is exactly the commuting equation for its restriction map.
			Preservation of the selected-state predicates, inherited-site
			identifications, and admission relations gives the remaining
			conditions of \cref{def:structural-isomorphism}. These two
			assignments are inverse on isomorphisms. Hence the encoding is
			fully faithful and therefore an equivalence.
			
			\item \emph{Witnesses and boundary relations.} By \cref{def:witness}, an assignment is a witness precisely
			when every vertex has a selected state and every internal vertex
			has an admitted tuple. These are the two conjunctions in
			$\varphi_T$. For a one-vertex tree with root $r$, the formula
			reduces to $S_{K_r}(x_r)$. A boundary assignment extends to a
			witness precisely when values can be assigned to the variables
			outside $B$, which gives $\varphi_{T,B}$. Finally, the two
			restriction graphs share a value $z$ exactly when the two
			restrictions agree at $U$. Each formula uses only atomic
			predicates, conjunction, and existential quantification.
		\end{enumerate}
	\end{proof}
	
	\paragraph{Which tests survive forgetting a relation?}
	Fix the types, marks and permissible isomorphisms.
	After forgetting some relations, a terminal test can still be evaluated if all presentations with the same reduct give the same answer.
	The next result states this for one test.
	Applying it to a family gives the joint profile of all its answers.
	
	\begin{proposition}[Tests and relations determined by retained data]
		\label{prop:forget-criterion}
		Let $F(\mathfrak L)$ be the presentation obtained from $\mathfrak L$
		by forgetting a specified collection of relations.
		Let $Q(\mathfrak L)$ be the outcome of a fixed isomorphism-invariant test.		
		\begin{enumerate}[label=\textup{(\roman*)}]
			\item \emph{Determination of a test.} The retained data determine $Q$ if and only if, for every
			pair of presentations $\mathfrak L,\mathfrak L'$,
			\[
			F(\mathfrak L)\cong F(\mathfrak L')
			\quad\Longrightarrow\quad
			Q(\mathfrak L)=Q(\mathfrak L').
			\]
			Here ``determine'' means that there is a function $\bar Q$ on the
			isomorphism classes of the reducts in the image of $F$ such that
			\[
			Q(\mathfrak L)=\bar Q\bigl([F(\mathfrak L)]\bigr)
			\qquad\text{for every }\mathfrak L.
			\]
			This function is unique on its stated domain.
			
			\item \emph{Reconstruction of relations.} Fix a formation signature $\alpha$ and a common named set of well-typed candidate tuples.
			The retained data determine the	admission relation $\MM_\alpha$ exactly if and only if they	determine, for each candidate tuple $t$, the answer to $t\in\MM_\alpha$.
			The same criterion applies to the native judgment relation $\mathsf J_\Xi$ when the candidates are named formation-datum--endpoint pairs.
		\end{enumerate}
	\end{proposition}

	\begin{proof}\leavevmode
		\begin{enumerate}[label=\textup{(\roman*)},leftmargin=*]
			\item \emph{Determination of a test.}
			If $\bar Q$ exists, isomorphic reducts give the same
			argument to $\bar Q$ and hence the same outcome.
			Conversely, given the displayed implication, define $\bar Q([R])$ to be $Q(\mathfrak L)$ for any presentation with $F(\mathfrak L)\cong R$.
			The implication	makes this value independent of the choice of $\mathfrak L$.
			Every class in the stated domain has such a representative, so the
			formula also proves uniqueness.
			
			\item \emph{Reconstruction of relations.}
			Knowing $\MM_\alpha$ gives the answer to every membership test.
			Conversely, the tuples whose tests return true form exactly $\MM_\alpha$.
			Since the tuple names are fixed, these answers identify the relation on the specified candidate set.
			Replacing tuples by pairs $(P,\mathfrak P)$ gives the assertion for $\mathsf J_\Xi$.
		\end{enumerate}
	\end{proof}

	\begin{example}[Two distinctions lost by forgetting data]
		\label{ex:two-reduct-failures}
		Forgetting $\MM_\alpha$ identifies the two presentations in \cref{ex:two-base-operational-choice}, although one admits a tuple excluded by the other.		
		For the native judgment, take contexts $\Xi_0,\Xi_1$ with the same admissible carriers, formation data $\{p,q\}$, and packaged endpoint $\mathfrak P$.
		Give them the respective laws
		\[
		\mathsf J_0=\{(p,\mathfrak P)\},
		\qquad
		\mathsf J_1=\{(p,\mathfrak P),(q,\mathfrak P)\}.
		\]
		Both contexts realize $p$ at $\mathfrak P$, whereas only $\Xi_1$ realizes $q$.
		The realized triples $(\Xi_j,p,\mathfrak P)$ have identical data after forgetting the judgment relation, even if all operational relations are retained.
	\end{example}

	\begin{proposition}[One-step trees suffice to check interface agreement]
		\label{prop:h2-iff}
		For each finite typed formation tree $T$, let $\WW_T$ be its
		operational witness set and let $M_T\subseteq X_T$ be the set of
		assignments satisfying all declared interface comparisons at every
		internal vertex.
		The following are equivalent:
		\begin{enumerate}[label=\textup{(\roman*)}]
			\item Every admitted primitive formation satisfies its interface
			comparisons, that is, $\MM_\alpha\subseteq\Match_\alpha$ for
			every signature $\alpha$ \textup{(H2)}.
			
			\item Every operational witness satisfies the interface comparisons
			throughout its tree, that is, $\WW_T\subseteq M_T$ for every
			finite typed formation tree $T$.
		\end{enumerate}
		In \textup{(ii)}, it suffices to check the one-step tree associated
		with each primitive signature.
	\end{proposition}

	\begin{proof}
		\leavevmode
		\begin{itemize}
			\item \textup{(i)} $\implies$ \textup{(ii)}. Take $x\in\WW_T$. At each internal vertex
			$v$, the tuple of input and output states assigned by $x$ belongs
			to $\MM_{\alpha_v}$ by \cref{def:witness}. By \textup{(i)}, this
			tuple belongs to $\Match_{\alpha_v}$. Thus $x\in M_T$, proving
			\textup{(ii)}.
			
			\item \textup{(ii)} $\implies$ \textup{(i)}. Fix a primitive signature $\alpha$ and form the one-step
			tree $T_\alpha$ with $\alpha$ at the root and one leaf for each
			ordered input position. Writing the leaf states first and the root
			state last, we have
			\[
			\WW_{T_\alpha}=\MM_\alpha,
			\qquad
			M_{T_\alpha}=\Match_\alpha.
			\]
			Hence $\WW_{T_\alpha}\subseteq M_{T_\alpha}$ gives precisely
			$\MM_\alpha\subseteq\Match_\alpha$. Checking these one-step trees
			therefore proves \textup{(i)}.
		\end{itemize}
	\end{proof}

	For a native obstruction, interpret the admitted tuple in \cref{ex:two-base-operational-choice} by the datum $q$ of \cref{ex:two-reduct-failures}.
	In the context with law $\mathsf J_0$, this sound tuple has no native endpoint for that datum.

	\subsection{Contract-relative irreducibility and source universality}

	Two formation histories can share a summary if every permitted continuation gives them the same observation.
	We first specify what minimality means for a given contract, then construct the coarsest summary for histories under union.

	\paragraph{Minimal hypotheses and minimum costs.}
	\label{subsec:hypothesis-minimality}
	Fix a baseline class $\mathcal B$, a conclusion $C$, and a finite set of additional hypotheses $\Sigma$.
	The set $\Sigma$ is	\emph{inclusion-minimal sufficient for $C$ over $\mathcal B$} if it implies $C$ there and, for each $H\in\Sigma$, some member of $\mathcal B$ satisfies $\Sigma\setminus\{H\}$ but fails $C$.
	Logical independence requires a model satisfying $\Sigma\setminus\{H\}$ and violating $H$, for each $H\in\Sigma$.

	\subsubsection{Universal source summaries}	
	\label{subsec:contextual-foundation}

	\begin{definition}[Irreducibility under a preservation contract]
		\label{def:contract-irreducibility}
		Fix the permitted comparisons of structural objects and a preservation
		contract $\Phi$.
		For objects $X$ and $Y$, let $\mathcal R(X,Y)$ denote the allowed reductions $q:X\to Y$, with their specified maps on the structural data.
		Write $\Phi(q)$ when $q$ satisfies the contract.		
		\begin{enumerate}[label=\textup{(\roman*)},leftmargin=*]
			\item \emph{Invertible and proper reductions.}
			A reduction $q:X\to Y$ is \textbf{invertible} if there is a permitted
			comparison $r:Y\to X$ such that
			\[
			r\circ q=\mathrm{id}_X
			\qquad\text{and}\qquad
			q\circ r=\mathrm{id}_Y.
			\]
			A reduction with no such inverse is \textbf{proper}.
			
			\item \emph{Irreducible objects.}
			The object $X$ is \textbf{irreducible relative to $(\mathcal R,\Phi)$}
			if every contract-preserving reduction out of $X$ is invertible:
			\[
			\forall Y,\quad \forall q\in\mathcal R(X,Y),\qquad
			\Phi(q)\ \Longrightarrow\ q\text{ is invertible}.
			\]
		\end{enumerate}
	\end{definition}

	The allowed reductions determine the meaning of irreducibility.
	Quotient irreducibility, product indecomposability and generator minimality are
	different tests.
	Even in a finite quotient poset an irreducible object is maximal under admissible coarsening, not necessarily of minimum cost.
	A greatest quotient, when it exists, is stronger than maximality.
	The	graph construction of \cref{thm:foundational-chromatic-reductions} makes
	these distinctions sharp.
	Here we first study quotients of source summaries.

	\begin{definition}[Exact summaries of static histories and their reductions]
		\label{def:contextual-source-summary}
		Fix a set $\mathcal U$ of edits.
		A \textbf{static formation history} is a subset $E\subseteq\mathcal U$.
		Histories combine by union, with empty history $\varnothing$.
		Write $\mathscr H=\mathcal P(\mathcal U)$ and fix an observation map $\omega:\mathscr H\to O$.		
		\begin{enumerate}[label=\textup{(\roman*)}]
			\item \textbf{Exact source summary.}
			An exact source summary of $(\mathscr H,\omega)$ is a triple
			$(a,S,\omega_S)$ consisting of a join-semilattice $(S,\vee_S,0_S)$
			with zero, a surjection $a:\mathscr H\twoheadrightarrow S$, and
			a decoder $\omega_S:S\to O$ such that, for all $E,F\in\mathscr H$,
			\[
			a(\varnothing)=0_S,\qquad
			a(E\cup F)=a(E)\vee_S a(F),\qquad
			\omega_S(a(E))=\omega(E).
			\]
			The last identity determines the decoder uniquely because $a$ is surjective.
			
			\item \textbf{Reduction of summaries.}
			Let $(a,S,\omega_S)$ and $(b,T,\omega_T)$ be exact summaries of	the same $(\mathscr H,\omega)$.
			A reduction from the first to the second is a surjective zero-preserving join homomorphism $p:S\twoheadrightarrow T$ satisfying
			\[
			p(a(E))=b(E)\qquad\text{for every }E\in\mathscr H.
			\]
			The reduction is \textbf{proper} if $p$ identifies two distinct
			elements of $S$.
		\end{enumerate}
	\end{definition}

	Every reduction preserves the declared observation.
	Indeed, each $s\in S$ has the form $s=a(E)$, and hence
	\[
	\omega_T(p(s))=\omega_T(b(E))=\omega(E)=\omega_S(s).
	\]
	
	The coarsest such summary is determined by observations after all common
	extensions of a history.
	Define
	\begin{equation}\label{eq:contextual-source-equivalence}
		E\equiv_\omega F
		\quad\Longleftrightarrow\quad
		\omega(E\cup C)=\omega(F\cup C)
		\quad\text{for every }C\subseteq\mathcal U.
	\end{equation}
	This is the syntactic-congruence construction applied to the union
	algebra of formation histories \cite[Lemma~2.1]{ClarkDaveyFreeseJackson2004}.

	\begin{theorem}[Canonical contextual source summary]
		\label{thm:canonical-contextual-summary}
		In Definition~\ref{def:contextual-source-summary}, put
		$Q_\omega=\mathscr H/{\equiv_\omega}$ and $q(E)=[E]$.
		\begin{enumerate}[label=(\roman*)]
			\item \emph{Construction of the quotient.}
			The relation $\equiv_\omega$ is the largest join-semilattice
			congruence contained in the equality relation of $\omega$.
			The operations $[E]\vee[F]=[E\cup F]$ and $0=[\varnothing]$ are
			well-defined, and $\bar\omega([E])=\omega(E)$ makes $q$ an exact
			source summary.
			It retains exactly the distinctions detected by	the observation profile after some common extension.
			
			\item \emph{Universal property.}
			Every exact source summary $(a,S,\omega_S)$ has a unique surjective zero-preserv\-ing join homomorphism
			\[
			p:S\twoheadrightarrow Q_\omega,
			\qquad q=pa,\qquad \omega_S=\bar\omega p.
			\]
			Thus $q$ is coarsest among the exact summaries of the stated class.
			
			\item \emph{Irreducibility and minimum cardinality.}
			Such a summary $(a,S,\omega_S)$ is irreducible under the reductions in Definition~\ref{def:contextual-source-summary} if and only if
			\begin{equation}\label{eq:contextual-separation-criterion}
				s\ne t\quad\Longrightarrow\quad
				\text{some }u\in S\text{ satisfies }
				\omega_S(s\vee u)\ne\omega_S(t\vee u).
			\end{equation}
			This is equivalent to the comparison $p$ in \emph{(ii)} being an
			isomorphism.
			Consequently $Q_\omega$ is the unique irreducible exact source summary up to the unique isomorphism commuting with the history maps.
			It has minimum cardinality among these summaries.
			If $Q_\omega$ is finite (in particular if the edit set is finite), equality of cardinality holds exactly for the irreducible summaries.
			
			\item \emph{Monotonicity.}
			If $\omega=X:\mathscr H\to\mathcal P(\Omega)$ is monotone, its decoder
			on $Q_\omega$ is monotone for the quotient semilattice order.
		\end{enumerate}
	\end{theorem}

	\begin{proof}
		\leavevmode
		\begin{enumerate}[label=\textup{(\roman*)}]
			\item \emph{Construction of the quotient.}
			By \eqref{eq:contextual-source-equivalence}, $E\equiv_\omega F$
			means that $\omega(E\cup C)=\omega(F\cup C)$ for every $C\in\mathscr H$.
			This is an equivalence relation, since it is equality of the functions $C\mapsto\omega(E\cup C)$.
			If $E\equiv_\omega F$, then $E\cup D\equiv_\omega F\cup D$ for every $D$: evaluate the defining equality at $D\cup C$.
			Applying this observation in each argument shows that $\equiv_\omega$ is a congruence for union.
			Taking $C=\varnothing$ gives $\equiv_\omega\,\subseteq\ker\omega$.
			
			To prove maximality, let $\theta$ be any join-semilattice congruence contained in $\ker\omega$.
			If $E\mathrel\theta F$,	then $E\cup C\mathrel\theta F\cup C$ for every $C$.
			Hence $\omega(E\cup C)=\omega(F\cup C)$ for every $C$, proving
			$\theta\subseteq\equiv_\omega$.
			
			The congruence therefore gives the quotient join $[E]\vee[F]=[E\cup F]$ with zero $[\varnothing]$.
			Its inclusion in $\ker\omega$ makes $\bar\omega([E])=\omega(E)$ well defined.
			Thus $q$ is an exact summary, and $q(E)\ne q(F)$ precisely when	some common extension distinguishes $E$ and $F$.
			
			\item \emph{Universal property.}
			Let $(a,S,\omega_S)$ be an exact summary.
			Since $a$ preserves	joins, $\ker a$ is a congruence.
			Since $\omega=\omega_S\circ a$,	we have $\ker a\subseteq\ker\omega$.
			Part~\textup{(i)} now gives	$\ker a\subseteq\equiv_\omega$.
			Consequently,
			\[
			p(a(E))=[E]\qquad(E\in\mathscr H)
			\]
			defines a map $p:S\to Q_\omega$.
			It is onto because every class $[E]$ occurs.
			It sends $0_S=a(\varnothing)$ to $[\varnothing]$ and preserves joins because
			\[
			p(a(E)\vee a(F))
			=p(a(E\cup F))
			=[E\cup F]
			=[E]\vee[F].
			\]
			The defining formula gives $p\circ a=q$.
			It also gives $\bar\omega(p(a(E)))=\omega(E)=\omega_S(a(E))$ for every $E$.
			Surjectivity of $a$ therefore yields $\bar\omega\circ p=\omega_S$ and makes $p$ unique.
			
			\item \emph{Irreducibility and minimum cardinality.}
			Write $s=a(E)$ and $t=a(F)$.
			Every $u\in S$ has the form $a(C)$,	and $\omega_S(s\vee a(C))=\omega(E\cup C)$.
			Thus
			\[
			p(s)=p(t)
			\quad\Longleftrightarrow\quad
			\omega_S(s\vee u)=\omega_S(t\vee u)
			\quad\text{for every }u\in S.
			\]
			Hence \eqref{eq:contextual-separation-criterion} holds exactly when $p$ is injective.
			Since $p$ is a surjective homomorphism,	this is equivalent to $p$ being an isomorphism of summaries.
			
			If $S$ is irreducible, its reduction $p$ must be an isomorphism.
			Conversely, suppose $p$ is injective, and let
			$r:S\twoheadrightarrow T$ be any permitted reduction.
			Applying \textup{(ii)} to the summary with history map $r\circ a$
			gives a comparison $p_T:T\twoheadrightarrow Q_\omega$.
			Since $(p_T\circ r)\circ a=q$, uniqueness in \textup{(ii)} gives $p=p_T\circ r$.
			Injectivity of $p$ forces $r$ to be injective.
			Thus every permitted reduction of $S$ is invertible.
			
			The comparison for $Q_\omega$ itself is the identity, so $Q_\omega$ is irreducible.
			Every irreducible summary is therefore isomorphic to $Q_\omega$ by the unique map commuting with the history maps.
			Finally, $p:S\twoheadrightarrow Q_\omega$ gives	$|Q_\omega|\leq |S|$.
			If $Q_\omega$ is finite, equality holds	exactly when $p$ is bijective, equivalently when $S$ is irreducible.
			
			\item \emph{Monotonicity.}
			Suppose $\omega=X:\mathscr H\to\mathcal P(\Omega)$ is monotone.
			If $[E]\leq[F]$, then $[E\cup F]=[F]$.
			Hence
			\[
			\bar\omega([E])=X(E)
			\subseteq X(E\cup F)=X(F)=\bar\omega([F]),
			\]
			which proves that the decoder is monotone.
		\end{enumerate}
	\end{proof}

	\begin{corollary}[Recovery of the contextual summary from native endpoints]
		\label{cor:contextual-native-recovery}
		Suppose an exact source summary $a$ also has chosen native realizations
		$\Xi\vdash P(E)\rightsquigarrow e(E)$ and a decoder	$d(e(E))=a(E)$ invariant under the declared endpoint comparison.
		Then $p\circ d$ recovers $q(E)$ from $e(E)$ and is invariant under
		that comparison.
		In particular, endpoints identified by the declared comparison have equal observations after every common added history.
	\end{corollary}

	\begin{proof}
		Use $q(E)=p(a(E))=p(d(e(E)))$ and the invariance of $d$.
	\end{proof}

	\begin{example}[Equal current observations distinguished by an extension]
		\label{ex:contextual-two-edits}
		Let $\mathcal U=\{u,v\}$ and let $X(E)=\{x\}$ precisely when
		$\{u,v\}\subseteq E$, with $X(E)=\varnothing$ otherwise.
		The witness requires both edits.
		The histories initially have equal observations, $X(\varnothing)=X(\{u\})$.
		Adjoining $\{v\}$ distinguishes them.
		Indeed every two distinct histories are contextually distinguishable: for $E\setminus F\ne\varnothing$, adjoining	$C=\mathcal U\setminus E$ makes $x$ occur after extending $E$ and leaves the witness set empty after extending $F$.
		Thus $Q_X=\mathcal P(\mathcal U)$ has four elements, although $X$
		has only two values.
		A quotient preserving only the current observation does not provide the required join-preserving summary.
	\end{example}

	\subsection{Contract-relative equivalence and compositional invariants}
	\label{subsec:n37-invariant-foundation}

	A continuation can observe only some of a component's marked data.
	We specify its available tests before comparing components.
	This gives a criterion for the invariant records used in later sections.

	\begin{definition}[Contextual equivalence and complete compositional invariants]
		\label{def:compositional-invariant}
		Fix a native formation language closed under typed substitution and
		a family of terminal tests.
		For each marked boundary type $\tau$, let $\mathcal X_\tau$ be the class of components of type $\tau$.
		
		A \textbf{continuation} for $\tau$ is a permitted finite formation
		expression $C[-]$ with one open position of type $\tau$. The marked
		types determine which continuations may be used. For each terminal
		test $t$ at the output type of $C[-]$, write $o_{C,t}(X)$ for the
		declared outcome of $t$ on $C[X]$, using the prescribed comparison
		of outcome values. The contract specifies the outcomes for
		unsuccessful formations and includes a success indicator whenever
		success itself is observed.
		
		Components $X,X'\in\mathcal X_\tau$ are \textbf{contextually equivalent},
		written $X\equiv_\tau X'$, if
		$
		o_{C,t}(X)=o_{C,t}(X')
		$
		for every continuation $C[-]$ for $\tau$ and every terminal test $t$ at its output type.	
		\begin{enumerate}[label=\textup{(\roman*)}]
			\item \textbf{Invariant.}
			An invariant is a family of maps $I_\tau:\mathcal X_\tau\to R_\tau$, where $R_\tau$ is a set of records, such that for all $X,X'\in\mathcal X_\tau$,
			\[
			X\equiv_\tau X'
			\quad\Longrightarrow\quad
			I_\tau(X)=I_\tau(X').
			\]
			
			\item \textbf{Completeness.}
			The invariant is complete if the converse holds for every $\tau$:
			\[
			I_\tau(X)=I_\tau(X')
			\quad\Longrightarrow\quad
			X\equiv_\tau X'.
			\]
			
			\item \textbf{Compositionality.}
			For each allowed formation $\alpha:(\tau_1,\ldots,\tau_m)\to\tau$, fix its declared interface data.
			The invariant is compositional if there is an operation
			$\widehat\alpha$ on the corresponding input records such that
			\[
			I_\tau\bigl(\alpha(X_1,\ldots,X_m)\bigr)
			=
			\widehat\alpha\bigl(
			I_{\tau_1}(X_1),\ldots,I_{\tau_m}(X_m)
			\bigr)
			\]
			for every tuple for which the formation is defined.
			The same operation must work for every choice of components with those
			records.
			If the contract observes whether a partial formation is defined, the input records must also determine its definedness.
		\end{enumerate}
		An invariant is a \textbf{complete compositional invariant} if it
		is both complete and compositional.
	\end{definition}

	Structural isomorphisms fixing the types, formation language and
	terminal tests give contextual equivalence whenever those formations
	and tests are invariant under the specified transports.
	The converse need not hold.
	For a set of components of a common type, contextual equivalence is compatible with formation.
	A terminal test performed after replacing one input is a test in the composed one-hole continuation.
	Replacing inputs one at a time proves compatibility for every finite formation.
	The resulting quotient gives a formal classification.
	A computational classification also needs a way to construct its record from the component and to compose records.
	
	\paragraph{The information required by different comparisons.}
	
	For group-native matching, \cref{thm:complete-compositional-invariants} constructs complete recovery and numerical records, including joint boundary images of agreement	kernels; proper marginal images and recovery closure alone are insufficient (\cref{ex:marginal-failure}).
	These records decide the reachable-boundary	contract (\cref{cor:arbitrary-graph-faithfulness}); block-local observation
	restrictions give the distributed criteria of
	\cref{thm:distributed-enrichment,thm:central-prime-blocks}.
	For independently chosen admissible lifts, attainable joint spectra replace
	single records; neither ambient records nor separate attainable coordinates
	suffice (\cref{ex:contract-spectrum-obstructions}).
	With literal marked	kernel couplings, diagonal tests recover decorated extension records
	(\cref{thm:coupled-extension-spectrum,ex:coupled-affine-separation}).
	The coupled records retain private readout kernels.
	A contract that also fixes private value labels requires those labels as additional data.
	
	\paragraph{Intrinsic recovery versus compositional information.}
	For a finite nonempty state set $X$ with labeled readouts $f_d$, the
	agreement relations $\theta_T=\bigcap_{d\in T}\ker f_d$ transport under
	bijections of the state set and compatible relabelings of readout
	values. The cardinalities $r_X(T)=|f_T(X)|$ are invariant under such
	comparisons, where $f_T$ is the joint readout and the empty readout
	has one value. The recovery law is determined by either record:
	\[
	d\in c_X(T)\ \Longleftrightarrow\
	\theta_T\subseteq\ker f_d\ \Longleftrightarrow\
	r_X(T\cup\{d\})=r_X(T).
	\]
	The last equivalence holds because the coordinate projection from
	$f_{T\cup\{d\}}(X)$ onto $f_T(X)$ is surjective and is bijective exactly
	when $T$ determines $d$. These records concern the specified
	observations, not necessarily all native data; $r_X$ is a cardinality
	profile, not the recovery-rank convention of \cite{Ma2026}.
	In the finite group model with trivial boundary, the qualitative
	record of \cref{thm:complete-compositional-invariants} reduces to
	the recovery law; its numerical record consists of the agreement
	kernel orders, from which $r_G(T)=|G|/|K_G(T)|$ follows.
	With nontrivial exposed boundaries one must also retain joint images.
	\Cref{ex:marginal-failure} exhibits a failure of marginal
	and recovery data to determine composition.
	
	\subsection{Canonical quotients and exact tree lifting}
	\label{subsec:canonical-reduction-subsection}
	
	Consider a quotient that retains the observation maps and descends the declared
	operations. A witness in the quotient specifies a class at each variable.
	When its incidence graph is a tree, one local test on the primitive
	relations decides whether such a witness lifts with any representative
	prescribed at any one variable.
	
	\paragraph{Finite source presentation.}
	Fix finitely many sorts $s\in\mathcal S$, finite nonempty sets $X_s$, a
	finite family $\mathcal F$ of typed total operations, and finitely many
	typed operational relations
	\[
	R_\alpha\subseteq
	X_{s_{\alpha,1}}\times\cdots\times X_{s_{\alpha,m_\alpha}},
	\qquad m_\alpha\geq 1.
	\]
	For a formation signature the last position may designate its output;
	the following lifting condition treats every position explicitly.
	Include each restriction map that must survive in $\mathcal F$ as a
	unary operation. Retain observation maps $o_s:X_s\to O_s$, allowing a
	tuple of readouts as one map.
	The $X_s$ are interpreted source or observation sets.
	At each formation sort all the displayed states are selected.
	A sort to be retained literally can be protected by including its identity map in $o_s$.
	
	A family of equivalence relations $\theta=(\theta_s)_s$ is an $\mathcal F$-congruence if every operation descends to the quotient sets $X_s/\theta_s$.
	For such a family write $q_s(x)=[x]_{\theta_s}$, and give every operational relation its coordinatewise image
	\[
	\overline R_\alpha
	=\{(q_{s_{\alpha,j}}(x_j))_{j=1}^{m_\alpha}:
	(x_j)_j\in R_\alpha\}.
	\]
	Coherence identities between restrictions descend because they are equations between operations.
	Interface soundness also descends when the stated comparisons use these restriction maps.
	No equality of global witness images follows merely from these primitive image
	equalities.
	
	\begin{definition}[Admissible quotient of a finite source presentation]
		\label{def:canonical-reduction-reduction}
		Let $\theta=(\theta_s)_{s\in\mathcal S}$ be a family of equivalence
		relations, one on each source set $X_s$, and let
		$q_s^\theta:X_s\twoheadrightarrow X_s/\theta_s$ be the quotient map.
		The quotient determined by $\theta$ is \textbf{admissible} if the
		following three conditions hold:		
		\begin{enumerate}[label=(\roman*)]
			\item \emph{Observations are preserved.}
			For every sort $s$ and all $x,y\in X_s$,
			\[
			x\mathrel{\theta_s}y \quad\Longrightarrow\quad o_s(x)=o_s(y).
			\]
			
			\item \emph{Operations descend.}
			For every operation
			$f:X_{s_1}\times\cdots\times X_{s_n}\to X_t$ in $\mathcal F$,
			the value
			\[
			\overline f([x_1]_{\theta_{s_1}},\ldots,[x_n]_{\theta_{s_n}})
			=[f(x_1,\ldots,x_n)]_{\theta_t}
			\]
			is independent of the chosen representatives.
			Equivalently, $\theta$ is an $\mathcal F$-congruence.
			
			\item \emph{Relation tuples lift with one prescribed coordinate.}
			For each $\alpha$, let
			$q_\alpha^\theta=\prod_{k=1}^{m_\alpha}q_{s_{\alpha,k}}^\theta$
			and $\overline R_\alpha=q_\alpha^\theta(R_\alpha)$.
			For every $\overline x\in\overline R_\alpha$, every position
			$j\in\{1,\ldots,m_\alpha\}$, and every representative
			$a\in X_{s_{\alpha,j}}$ of $\overline x_j$, there exists
			$x\in R_\alpha$ satisfying
			\begin{equation}\label{eq:canonical-reduction-back}
				q_\alpha^\theta(x)=\overline x
				\qquad\text{and}\qquad x_j=a.
			\end{equation}
		\end{enumerate}		
		The allowed reductions are the admissible quotient maps and their factor maps over the source.
		Explicitly, if $\theta$ and $\eta$ are admissible and $\theta_s\subseteq\eta_s$ for every $s$, the factor map at sort $s$ is
		\[
		p_s:X_s/\theta_s\longrightarrow X_s/\eta_s,
		\qquad p_s([x]_{\theta_s})=[x]_{\eta_s}.
		\]
		It is the unique map satisfying $p_s\circ q_s^\theta=q_s^\eta$.
		All reductions keep the sort, operation, and relation names fixed and preserve the observation values.
	\end{definition}

	For clarity, a tree of relation occurrences has a bipartite incidence graph:
	one kind of vertex is a typed variable occurrence, and the other is a relation occurrence joined to its argument positions.
	Require this graph to be a finite tree.
	Different positions in one relation are	different occurrences, even when they have the same sort.
	A witness assigns a state to every variable so that every relation holds.
	Ordinary finite formation trees are examples.
	Fixing a witness in the quotient fixes the quotient value at every variable.
	A separator consisting of several coordinates must be treated as one
	specified joint interface sort when applying this argument across that
	separator.
	The lift must realize all its coordinates at once.

	\paragraph{Native quotient hypothesis.}
	Fix the native class, its realization laws, its retained markings, and
	its source interpretation before forming any quotient.
	Assume a specified quotient construction with the following property:
	\begin{quote}
		\textup{(NQ)} Every $\mathcal F$-congruence quotient of the supplied source
		is assigned a valid realization in that same native class, whose decoded
		source is the stated quotient.
		The declared native comparisons retain the quotient maps and are compatible with identity maps and iterated quotients.
		The specified operations and markings descend by those maps.
	\end{quote}
	\Cref{subsec:canonical-reduction-native-verification} verifies (NQ) for marked
	finite algebra presentations when $\mathcal F$ includes the basic operations of each native algebra.
	Additional operations and interfaces descend by the same $\mathcal F$-congruence test.

	\begin{theorem}[Canonical reduction with exact tree lifting]
		\label{thm:canonical-reduction-canonical}
		Fix a finite source presentation as above, with observation maps $o_s$, operations $\mathcal F$, and relations $R_\alpha$.		
		\begin{enumerate}[label=(\roman*)]
			\item \emph{Greatest admissible congruence.}
			There is a greatest admissible congruence family
			$\theta^*=(\theta_s^*)_{s\in\mathcal S}$ under inclusion at	each sort.
			Its quotient
			\[
			q_s^*:X_s\twoheadrightarrow Q_s^*,
			\qquad Q_s^*=X_s/\theta_s^*,
			\]
			is the coarsest admissible quotient over $X$.
			Starting from the observation partitions $\ker o_s$, finite partition
			refinement computes $\theta^*$ from the operation and relation tables.
			
			\item \emph{Lifting a tree witness with one prescribed representative.}
			Let $q:X\twoheadrightarrow\overline X$ be an observation-preserving $\mathcal F$-congruence quotient.
			The local lifting condition \eqref{eq:canonical-reduction-back}	holds if and only if the following property holds for every finite tree $T$ of relation occurrences.
			Given a quotient witness $\overline w\in\overline W_T$, a variable $v$ of sort $s(v)$, and a state $a\in X_{s(v)}$ with $q_{s(v)}(a)=\overline w_v$, there is a witness $w\in W_T$ such that
			\[
			q_T(w)=\overline w
			\qquad\text{and}\qquad
			w_v=a,
			\]
			where $q_T$ applies the quotient map at each variable.
			Under these equivalent conditions,
			\begin{equation}\label{eq:canonical-reduction-tree-image}
				q_T(W_T)=\overline W_T.
			\end{equation}
			The corresponding equality holds after projection to any designated boundary.
			The joint relation of retained readout values is therefore unchanged, as are all determination and non-determination judgments among those readouts.
			These conclusions also hold after any typed substitution whose expanded incidence graph is a tree.
			
			\item \emph{Universal comparison and minimum size.}
			For every admissible quotient $a:X\twoheadrightarrow Y$,
			there is a unique admissible comparison
			$p:Y\twoheadrightarrow Q^*$ such that $q^*=p\circ a$.
			The quotient $Q^*$ is irreducible, and every irreducible admissible quotient is uniquely isomorphic to $Q^*$ over $X$.
			Moreover,
			$
			|Q_s^*|\leq |Y_s| \quad (s \in \mathcal S).
			$
			Equality holds at every sort if and only if $p$ is an isomorphism.
			
			\item \emph{Native realization.}
			Under \textup{(NQ)}, $Q^*$ has a realization in the specified native class, and the comparisons in \textup{(iii)} are permitted native quotient comparisons over the supplied realized sources.
			
			\item \emph{Static formation histories.}
			Let $\mathcal U$ be finite.
			In the one-sort presentation $X=\mathcal P(\mathcal U)$ with operations $\cup,\varnothing$, observation $\omega$, and no operational relations, the canonical congruence is contextual equivalence:
			\begin{equation}\label{eq:canonical-reduction-static}
				E\mathrel{\theta^*}F
				\quad\Longleftrightarrow\quad
				E\equiv_\omega F
				\quad\Longleftrightarrow\quad
				\omega(E\cup C)=\omega(F\cup C)
				\quad\text{for every }C\subseteq\mathcal U.
			\end{equation}
		\end{enumerate}
	\end{theorem}

	\begin{proof}
		\leavevmode
		\begin{enumerate}[label=(\roman*)]
			\item \emph{Greatest admissible congruence.}
			An elementary translation $t:X_s\to X_r$ is obtained by
			fixing all but one argument of a positive-arity operation
			in $\mathcal F$. A family of equivalence relations is an
			$\mathcal F$-congruence if and only if every elementary
			translation preserves it. Indeed, compatibility with an
			operation follows by changing its arguments one at a time.
			Nullary operations impose no further condition.
			
			For a partition family $\pi=(\pi_s)_s$, a relation position
			$(\alpha,j)$, and a state $x\in X_{s_{\alpha,j}}$, define
			\begin{equation}\label{eq:canonical-reduction-profile}
				\Lambda_{\alpha,j}^{\pi}(x)
				=
				\left\{
				\bigl([z_k]_{\pi_{s_{\alpha,k}}}\bigr)_{k\ne j}
				\;\middle|\;
				z\in R_\alpha,\ z_j=x
				\right\}.
			\end{equation}
			This set records which tuples of classes can occur in the
			remaining positions when position $j$ is fixed at $x$.
			The positions retain their order, including when two positions
			have the same sort. For a unary relation, the set contains
			the empty tuple exactly when $x$ belongs to the relation.
			
			Start with $\pi_s^0=\ker o_s$. For $x,y\in X_s$, put
			$x\mathrel{\pi_s^{n+1}}y$ exactly when all three conditions
			below hold:
			\begin{equation}\label{eq:canonical-reduction-refine}
				\begin{aligned}
					x&\mathrel{\pi_s^n}y,\\
					t(x)&\mathrel{\pi_r^n}t(y)
					&&\text{for every elementary translation }
					t:X_s\to X_r,\\
					\Lambda_{\alpha,j}^{\pi^n}(x)
					&=\Lambda_{\alpha,j}^{\pi^n}(y)
					&&\text{for every }(\alpha,j)
					\text{ with }s_{\alpha,j}=s.
				\end{aligned}
			\end{equation}
			Each condition defines an equivalence relation on $X_s$.
			Their intersection is therefore an equivalence relation
			refining $\pi_s^n$. Write $N=\sum_s|X_s|$ and
			$b_0=\sum_s|X_s/\pi_s^0|$. Each strict refinement increases
			the total number of blocks, so there are at most $N-b_0$
			strict rounds. Let $\theta^*$ be the stable family.
			
			We first prove that $\theta^*$ is admissible. It refines the
			observation partitions, and stability under elementary
			translations makes it an $\mathcal F$-congruence. To verify
			\eqref{eq:canonical-reduction-back}, choose a quotient tuple
			$\overline z\in\overline R_\alpha$ and a preimage
			$z\in R_\alpha$. If $a\mathrel{\theta^*_{s_{\alpha,j}}}z_j$,
			stability gives
			$\Lambda_{\alpha,j}^{\theta^*}(a)
			=\Lambda_{\alpha,j}^{\theta^*}(z_j)$.
			Thus some $z'\in R_\alpha$ has $z'_j=a$ and lies in the
			same $\theta^*$-class as $z$ at every other position.
			This is the required lift of $\overline z$.
			
			To prove greatestness, let $\phi$ be any admissible family.
			We show by induction that $\phi_s\subseteq\pi_s^n$ for
			every sort $s$. At $n=0$ this is observation preservation.
			Suppose the inclusion holds at $n$, and let $x\mathrel{\phi_s}y$.
			For every elementary translation $t:X_s\to X_r$, congruence
			compatibility gives $t(x)\mathrel{\phi_r}t(y)$, hence
			$t(x)\mathrel{\pi_r^n}t(y)$. Now take any $z\in R_\alpha$
			with $z_j=x$. Apply the lifting condition for $\phi$ to the
			image of $z$, prescribing $y$ at position $j$. The resulting
			tuple has the same $\phi$-classes, and hence the same
			$\pi^n$-classes, at all remaining positions. This proves
			$\Lambda_{\alpha,j}^{\pi^n}(x)
			\subseteq\Lambda_{\alpha,j}^{\pi^n}(y)$.
			Interchanging $x$ and $y$ gives equality. All three
			conditions in \eqref{eq:canonical-reduction-refine} hold,
			so $\phi_s\subseteq\pi_s^{n+1}$. At stabilization we obtain
			$\phi_s\subseteq\theta_s^*$ at every sort.
			
			\item \emph{Lifting a tree witness with one prescribed representative.}
			Assume \eqref{eq:canonical-reduction-back}. Fix a quotient
			witness $\overline w$ and root the incidence tree at the
			variable $v$ where the representative $a$ is prescribed.
			Set $w_v=a$ and process the relation occurrences in order
			of increasing distance from $v$. Each occurrence meets the
			assigned variables at exactly one position, its parent
			variable. Apply the local lifting condition to its quotient
			tuple, prescribing the value already assigned at that
			position. The resulting tuple assigns values to its child
			variables. Since the incidence graph is a tree, each child
			receives its value exactly once. This constructs a witness
			$w$ with $q_T(w)=\overline w$ and $w_v=a$. For a tree
			consisting of a single variable, the initial assignment
			already gives the lift.
			
			Conversely, take the tree with one occurrence of $R_\alpha$
			and a distinct variable at each argument position. Its
			witnesses are exactly the tuples of $R_\alpha$. Prescribing
			a representative at position $j$ gives precisely
			\eqref{eq:canonical-reduction-back}.
			
			Every source witness maps to a quotient witness because
			each quotient relation is the image of its source relation.
			The lifting argument gives the reverse inclusion in
			\eqref{eq:canonical-reduction-tree-image}. For a set $B$ of
			boundary variables, coordinate projection commutes with the
			quotient maps, so
			\[
			q_B\bigl(\operatorname{pr}_B W_T\bigr)
			=\operatorname{pr}_B\overline W_T.
			\]
			Each retained readout factors through its quotient map.
			The joint readout values on the two witness sets therefore
			coincide. A failure of determination is witnessed by two
			such value tuples agreeing on the determining readouts and
			differing on the target readout. Thus determination and
			non-determination are both preserved. The same argument
			applies to a substituted formation whenever its expanded
			incidence graph is a tree.
			
			\item \emph{Universal comparison and minimum size.}
			Let $a:X\twoheadrightarrow Y$ be admissible and put
			$\phi_s=\ker a_s$. By \textup{(i)},
			$\phi_s\subseteq\theta_s^*$, so
			\[
			p_s(a_s(x))=q_s^*(x)
			\]
			defines a surjection $p_s:Y_s\twoheadrightarrow Q_s^*$.
			Surjectivity of $a_s$ makes this the unique map satisfying
			$q_s^*=p_s\circ a_s$. The operations and readouts on both
			quotients descend from $X$, so $p$ preserves them. It also
			maps each relation on $Y$ onto the corresponding relation
			on $Q^*$, since both are images of the same source relation.
			
			To check lifting for $p$, take a relation tuple in $Q^*$ and
			prescribe a representative $y\in Y_{s_{\alpha,j}}$ at
			position $j$. Choose $x\in X_{s_{\alpha,j}}$ with $a(x)=y$.
			Lift the tuple through $q^*$ with $x$ prescribed at $j$.
			Applying $a$ to this lift gives a relation tuple in $Y$
			with the prescribed value $y$ and the required image under
			$p$. Hence $p$ is admissible.
			
			Admissible quotients compose. For the lifting condition,
			first lift a tuple through the second quotient with the
			image of the prescribed state, then lift through the first
			quotient with the prescribed state itself. Thus a proper
			admissible quotient of $Q^*$ would yield an admissible
			congruence on $X$ strictly larger than $\theta^*$ at some
			sort. Greatestness rules this out, so $Q^*$ is irreducible.
			If $Y$ is irreducible, its admissible comparison $p$ must
			be invertible, and the factorization over $X$ makes this
			isomorphism unique.
			
			Finally, each $p_s$ is surjective, giving
			$|Q_s^*|\leq|Y_s|$. Since the sets are finite, equality
			at every sort holds exactly when all $p_s$ are bijective.
			In that case their inverses preserve the descended
			operations, readouts, and relation tables, so $p$ is an
			isomorphism.
			
			\item \emph{Native realization.}
			Apply the construction in \textup{(NQ)} to the congruence
			$\theta^*$. It supplies a native realization of $Q^*$.
			Compatibility with successive quotients makes the
			factorizations $q^*=p\circ a$ in \textup{(iii)} permitted
			native comparisons over the original sources.
			
			\item \emph{Static formation histories.}
			In this specialization, admissibility means exactly that
			the equivalence relation is a join-semilattice congruence
			contained in $\ker\omega$. By
			\cref{thm:canonical-contextual-summary}, the greatest such
			congruence is $\equiv_\omega$. Its defining condition is
			equality of observations after every common extension $C$,
			which gives \eqref{eq:canonical-reduction-static}.\qedhere
		\end{enumerate}
	\end{proof}
	
	\begin{corollary}[Admissibility of the contextual quotient after adding operations and relations]
		\label{cor:canonical-reduction-promotion}
		Let $\mathcal U$ be finite, let $X=\mathcal P(\mathcal U)$, and fix
		an observation profile $\omega:X\to O$.
		Write
		\[
		Q_\omega=X/{\equiv_\omega},
		\qquad q_\omega(E)=[E]_{\equiv_\omega},
		\]
		where $\equiv_\omega$ is the contextual equivalence in
		\eqref{eq:canonical-reduction-static}.
		Enlarge the presentation $(X,\cup,\varnothing,\omega)$ by finitely many total operations and operational relations on $X$, keeping $\omega$ fixed.
		Let $\theta^*$ be the greatest admissible congruence of the enlarged presentation.		
		\begin{enumerate}[label=(\roman*)]
			\item \emph{Criterion for admissibility.}
			The map $q_\omega$ is an admissible quotient of the enlarged
			presentation if and only if both of the following hold.
			\begin{enumerate}[label=(\alph*)]
				\item Every added operation respects $\equiv_\omega$.
				
				\item For every added relation $R_\alpha$, every tuple
				$(E_1,\ldots,E_{m_\alpha})\in R_\alpha$, every position
				$j$, and every $A\equiv_\omega E_j$, there is a tuple
				$(F_1,\ldots,F_{m_\alpha})\in R_\alpha$ such that
				\[
				F_j=A
				\qquad\text{and}\qquad
				F_k\equiv_\omega E_k
				\quad\text{for every }k\ne j.
				\]
			\end{enumerate}
			When these conditions hold, $\theta^*=\equiv_\omega$, and every quotient witness on a finite tree lifts with any one prescribed representative.
			Under \textup{(NQ)}, $q_\omega$ also gives an admissible native reduction.
			
			\item \emph{Comparison when admissibility fails.}
			If either condition fails, then $\theta^*\subsetneq\equiv_\omega$.
			The canonical map
			\[
			p:X/\theta^*\twoheadrightarrow Q_\omega,
			\qquad p([E]_{\theta^*})=[E]_{\equiv_\omega},
			\]
			preserves joins, zero, and the decoded observations.
			Every pair of histories identified by $q_\omega$ but separated by $X/\theta^*$ is separated by every admissible quotient of the enlarged presentation.
		\end{enumerate}
	\end{corollary}
	
	\begin{proof}
		\leavevmode
		\begin{enumerate}[label=(\roman*)]
			\item By \cref{thm:canonical-contextual-summary}, $\equiv_\omega$ is the greatest join-semilattice congruence contained in $\ker\omega$.
			Thus $q_\omega$ already preserves the observations, and union and zero descend to $Q_\omega$.
			Condition~(a) is exactly what is required for the added operations to descend.
			Every tuple in a quotient relation is the image of a source tuple, so condition~(b) is exactly \eqref{eq:canonical-reduction-back} for $q_\omega$.
			This proves the admissibility criterion.
			
			Every admissible congruence of the enlarged presentation still respects union and is contained in $\ker\omega$.
			Hence every such congruence is contained in $\equiv_\omega$.
			When \textup{(a)} and \textup{(b)} hold, $\equiv_\omega$ is itself admissible and therefore equals $\theta^*$.
			The tree lifting assertion follows from
			\cref{thm:canonical-reduction-canonical}\textup{(ii)}, and the native assertion follows from part~\textup{(iv)} of that theorem.
			
			\item We have shown that $\theta^*\subseteq\equiv_\omega$.
			If either condition fails, $\equiv_\omega$ is not admissible, whereas $\theta^*$ is admissible.
			The inclusion is therefore strict.
			It also makes $p$ well-defined.
			Both quotient maps are surjective and carry union, zero, and $\omega$ to their respective quotient operations and decoders, so $p$ preserves all three.
			
			Finally, if $a:X\twoheadrightarrow Y$ is any admissible	quotient, greatestness gives $\ker a\subseteq\theta^*$.
			Thus $[E]_{\theta^*}\ne[F]_{\theta^*}$ implies $a(E)\ne a(F)$, including when $q_\omega(E)=q_\omega(F)$.
		\end{enumerate}
	\end{proof}
	
	\begin{remark}[Partition refinement and the back condition]
		\label{rem:canonical-reduction-scope}
		The proof uses syntactic congruences and stable partitions
		\cite{ClarkDaveyFreeseJackson2004,PaigeTarjan1987}.
		In the binary case, \eqref{eq:canonical-reduction-back} applies the
		bounded-morphism back condition at each argument position
		\cite[Section 2.2.1]{GorankoOtto2007}.
		The following examples separate exact
		tree lifting at one prescribed variable from unpinned and cyclic
		lifting.
	\end{remark}
	
	\subsubsection{Examples of observation and lifting obstructions}
	\label{app:canonical-reduction-limits}

	Three examples separate the roles of the native operations, prescribed
	local representatives, and the shape of the attachment diagram. Each
	changes one of the requirements used to construct or lift a quotient.

	\begin{example}[A native operation can obstruct the static quotient]
		\label{ex:canonical-reduction-operation}
		Let $X=\mathcal P(\{u,v\})$ with union, zero, and the operation
		$\tau$ interchanging $u$ and $v$. Retain
		$\omega(E)=\mathbf 1_{\{u\in E\}}$. The static contextual quotient
		has two classes, determined by membership of $u$. But
		\[
		\varnothing\equiv_\omega\{v\},\qquad
		\omega(\tau\varnothing)=0\ne1=\omega(\tau\{v\}).
		\]
		Thus $\tau$ cannot descend to that quotient. Any observation-preserving
		congruence for both union and $\tau$ must preserve membership of both
		$u$ and $v$, and hence must be the identity. The canonical admissible
		source has four states, although the static summary has two.
		This occurs within the equational class of join-semilattices with
		zero and an involutive join automorphism.
	\end{example}

	\begin{example}[Unpinned witness exactness need not have a coarsest quotient]
		\label{ex:canonical-reduction-no-coarsest}
		Take three sorts
		\[
		A=\{a,d\},\quad B=\{0,1,2,3\},\quad C=\{c,e\},
		\]
		with operational relations
		\[
		R=\{(a,0),(d,3)\}\subseteq A\times B,\qquad
		S=\{(2,c),(3,e)\}\subseteq B\times C.
		\]
		Keep $A,C$ pointwise observed and give $B$ a constant observation.
		There are no additional operations. For the directed signatures
		$A\to B$ and $B\to C$, the only nontrivial two-step formation has
		the nonempty witness set
		\[
		W=\{(d,3,e)\}.
		\]
		Consider the two partitions of $B$
		\[
		\phi_1=\{\{0,1\},\{2\},\{3\}\},\qquad
		\phi_2=\{\{0\},\{1,2\},\{3\}\}.
		\]
		Each quotient is exact on this witness set and on each primitive
		relation. Each is therefore exact on every ordinary directed
		formation tree in this two-signature presentation. Any common
		coarsening identifies $0$ and $2$, so its quotient chain has the
		additional witness $(a,[0],c)$. The fixed observations at $A,C$ prevent any preimage of this
		witness in $W$.
		The example applies to finite sets and to any fixed native
		interpretation that supplies these set quotients.
		
		Both $\phi_i$ fail the lifting test with a prescribed representative.
		For each partition one primitive relation has a quotient tuple
		in the merged block that cannot lift with representative $1$.
		Theorem~\ref{thm:canonical-reduction-canonical} requires that stronger test.
	\end{example}

	\begin{example}[Tree lifting does not imply cyclic lifting]
		\label{ex:canonical-reduction-cycle}
		Let $X=\{0,1\}$, retain a constant observation, and take
		$R=\{(0,1),(1,0)\}$. The quotient $X\to\{*\}$ satisfies
		\eqref{eq:canonical-reduction-back} at both argument positions. Every
		quotient tree witness lifts after either representative is prescribed
		at any one vertex. But the triangle
		\[
		R(x,y),\quad R(y,z),\quad R(z,x)
		\]
		has no concrete witness, whereas its one-state quotient has one.
		Even a single edge has no lift when both representatives are
		prescribed as zero.
	\end{example}

	For a classification and sharp minimum-size obstruction for unpinned
	exact reductions in coherent deterministic hierarchies, see
	\cref{thm:foundational-chromatic-reductions,cor:foundational-depth-boundary}.
	For cyclic diagrams, \cref{prop:deterministic-cyclic-quotient,thm:multivalued-linear-quotient}
	give exact criteria for isomorphism transport and linear relations.

	\subsection{Sharp non-collapse and the complexity of exact reductions}
	\label{subsec:foundational-chromatic-reductions}
	
	Here a quotient witness may lift through any representatives. A greatest exact quotient can already fail to exist on a chain (\cref{ex:canonical-reduction-no-coarsest}). Already at depth two,
	its admissible central partitions can be exactly the proper coloring
	partitions of any finite graph. All original formations have witnesses.
	
	\begin{definition}[Exact observation quotient for directed formation trees]
		\label{def:foundational-unpinned-quotient}
		Fix a finite coherent hierarchical structural presentation with its
		selected-state observation interpretation.
		Assume that every observation state at a formation type is selected.
		At each site $U$, let $X_U$ be the ambient observation set and fix a retained readout $o_U:X_U\to O_U$, which may be constant.		
		A local image reduction $q=(q_U:X_U\twoheadrightarrow\overline X_U)_U$ in the sense of \cref{subsec:structural-comparisons} is an
		\textbf{exact observation quotient} if it satisfies both conditions below:
		\begin{enumerate}[label=(\roman*)]
			\item \emph{Readout preservation.}
			For every site $U$, there is a map $\overline o_U:\overline X_U\to O_U$ such that
			$
			o_U=\overline o_U\circ q_U.
			$
			
			\item \emph{Exactness on directed formation trees.}
			For every ordinary finite directed formation tree $T$, let $W_T$ and $\overline W_T$ be its witness sets in the original and reduced presentations, respectively.
			If $K_v$ is the	formation type at vertex $v$, then
			\[
			q_T(W_T)=\overline W_T,
			\qquad
			q_T=\prod_{v\in V(T)}q_{K_v}.
			\]
			Equivalently, every $\overline w\in\overline W_T$ admits a witness $w\in W_T$ satisfying	$q_{K_v}(w_v)=\overline w_v$ at every vertex $v$.
		\end{enumerate}
	\end{definition}

	We identify exact observation quotients that differ only by relabeling the quotient states at each site.
	Each quotient $q$ is then represented by its kernel family $\theta=(\ker q_U)_U$.
	We order these families by
	\[
	\theta\leq\phi
	\quad\Longleftrightarrow\quad
	\theta_U\subseteq\phi_U\quad\text{for every site }U.
	\]
	Thus $\phi$ is a coarsening of $\theta$ when $\theta\leq\phi$.
	An exact quotient is maximal if it has no strictly coarser exact quotient.

	For two such quotient families $\theta\subseteq\phi$, the induced site
	maps $X_U/\theta_U\to X_U/\phi_U$ commute with restrictions and carry the
	image relations onto the coarser image relations. They are exact on every
	formation tree: if the corresponding product maps are $p_T$, then
	$p_T(q_T^\theta(W_T))=q_T^\phi(W_T)=W_T^\phi$. Conversely, composing two
	exact reductions is exact by the same calculation. Maximality is therefore
	irreducibility relative to these reductions in the sense of
	\cref{def:contract-irreducibility}.

	\begin{lemma}[Exactness of a two-step deterministic formation]
		\label{lem:foundational-two-step-quotient}
		Let $f:A\to B$ and $g:B\to C$ be total functions, with their graphs as the two primitive unary formation relations.
		Keep $A,C$ pointwise observed, and let $q:B\twoheadrightarrow\overline B$ be a quotient.
		\begin{enumerate}[label=\textup{(\roman*)},leftmargin=*]
			\item \emph{Witness exactness.}
			The two-step formation is witness-exact if and only if
			\begin{equation}\label{eq:foundational-two-step-quotient}
				q(b)=q(f(a))\quad\Longrightarrow\quad g(b)=g(f(a))
				\qquad(a\in A,\ b\in B).
			\end{equation}
			In particular, the condition is that $g$ be constant on each quotient
			class meeting $f(A)$.
			\item \emph{Boundary exactness.}
			The reduced leaf-to-root boundary relation equals the image of the original boundary relation if and only if \eqref{eq:foundational-two-step-quotient} holds.
		\end{enumerate}
	\end{lemma}
	
	\begin{proof}\leavevmode
		\begin{enumerate}[label=\textup{(\roman*)},leftmargin=*]
			\item \emph{Witness exactness.}
			Identify the endpoint quotients with $A,C$ by their readouts.
			The reduced	witness set consists of the triples $(a,q(f(a)),g(b))$ with
			$q(b)=q(f(a))$.
			The image of the original witness set consists precisely of $(a,q(f(a)),g(f(a)))$ for $a\in A$.
			Equality is thus equivalent to \eqref{eq:foundational-two-step-quotient}.
			\item \emph{Boundary exactness.}
			Projecting away the middle coordinate gives the same condition because that coordinate is uniquely determined by $a$.
			One-step and single-vertex witnesses are always exact by the definition of local image reduction.
		\end{enumerate}
	\end{proof}
	
	A coloring partition of a finite simple graph is a partition of its
	vertex set into nonempty independent sets. Its blocks are the color
	classes, without names for the colors. Write $\chi(\Gamma)$ for the
	chromatic number and $\mathcal P_{\mathrm{col}}(\Gamma)$ for the poset
	of coloring partitions ordered by coarsening (equivalently, inclusion of
	the corresponding equivalence relations).
	
	\begin{theorem}[Graph coloring and exact observation quotients]
		\label{thm:foundational-chromatic-reductions}
		Let $\Gamma=(V,E)$ be a finite simple graph with $V\ne\varnothing$,
		and put $m=|E|$.
		There is a finite Level-II structural presentation satisfying \textup{(H1)--(H3)} with a distinguished site $B$ whose state set is
		$
		X_B=V\sqcup\{z\}.
		$
		The interface value of $z$ distinguishes it from every vertex in $V$.
		Every primitive formation is unary, total, and deterministic.
		All declared comparisons use the same nonconstant binary interface.
		Every directed formation tree has depth at most two and admits a
		witness. The source sets and their site quotients have native
		realizations under the finite marked-set presentation
		law of \cref{def:equational-native-context} with empty algebraic
		signature.		
		Order exact observation quotients by inclusion of their kernel families, and order the coloring partitions
		$\mathcal P_{\mathrm{col}}(\Gamma)$ by coarsening.		
		\begin{enumerate}[label=\textup{(\roman*)},leftmargin=*]
			\item \emph{Classification by coloring partitions.}
			The map sending an exact observation quotient to its kernel
			partition on $V$ at site $B$ is an order-preserving surjection
			onto $\mathcal P_{\mathrm{col}}(\Gamma)$ and has an
			order-preserving section. Thus two distinct vertices can be
			identified in some exact quotient if and only if they are
			nonadjacent. Every exact quotient preserves the joint relation
			of retained readout values on every directed formation tree,
			and hence all determination and non-determination judgments
			among those readouts.
			
			\item \emph{Minimum size and computational complexity.}
			Among all exact observation quotients,
			\begin{equation}\label{eq:foundational-chromatic-minimum}
				\min_q|\overline X_B|=1+\chi(\Gamma).
			\end{equation}
			The minimum of the total site cost
			$\sum_U|\overline X_U|$ is $\chi(\Gamma)+3+8m$.
			On this explicit family, deciding whether an exact quotient
			satisfies a given bound on $|\overline X_B|$ is $\NP$-complete
			when the input consists of the graph, its associated
			presentation, and the bound.
			
			\item \emph{Greatest quotients and the effect of retaining copy
				bijections.}
			A greatest achievable partition on $V$ exists if and only if
			$\Gamma$ is complete multipartite, including the one-part
			edgeless case. A greatest quotient family over all sites exists
			if and only if $E=\varnothing$ or $|V|\leq 2$.
			
			Suppose in addition that each specified copy bijection
			$B\to B_e$, one for each edge $e$, must remain bijective
			after reduction. Then restriction to $V$ at $B$ is an order
			isomorphism from the full poset of exact quotients onto
			$\mathcal P_{\mathrm{col}}(\Gamma)$. Under this requirement,
			a greatest quotient exists if and only if $\Gamma$ is complete
			multipartite, and the minimum total site cost is
			\[
			(m+1)(1+\chi(\Gamma))+5m+2.
			\]
			
			\item \emph{Irreducibility and the gap from minimum size.}
			When the copy bijections are retained, an exact quotient is
			irreducible if and only if every two distinct color classes
			have an edge between them. Even among bipartite graphs with
			chromatic number two, the ratio of the number of central
			states in an irreducible quotient to the minimum number of
			central states is unbounded.
		\end{enumerate}
	\end{theorem}
	
	\begin{proof}
		\emph{Construction.}
		Choose $z\notin V$. Take a common interface site $L$ with
		$X_L=\{0,1\}$ and a central site $B$ with $X_B=V\sqcup\{z\}$.
		For each edge $e$, choose an ordering $(u_e,v_e)$ of its endpoints
		and introduce three distinct sites with state sets
		\[
		X_{A_e}=\{a_e,d_e\},\qquad
		X_{B_e}=V\sqcup\{z\},\qquad
		X_{C_e}=\{c_e^0,c_e^1,c_e^*\}.
		\]
		Each $X_{B_e}$ is a named copy of $X_B$. Put $L$ below every
		other site and $B_e$ below $B$ for each edge. These are all the
		strict order relations. Use the copy bijection as the restriction
		from $B$ to $B_e$. The restrictions to $L$ send $z$ (at $B$ and
		each $B_e$), $d_e$, and $c_e^*$ to $1$, and all other states to
		$0$. In particular, restriction from $B$ to $L$ agrees with
		restriction through any $B_e$, proving \textup{(H1)}.
		Retain the identity readouts at $L,A_e,C_e$ and the restrictions
		to $L$ as readouts at $B,B_e$.
		
		The only primitive formations are $A_e\to B_e$ and $B_e\to C_e$.
		Their admission relations are the graphs of the total functions
		\[
		\begin{aligned}
			f_e(a_e)&=u_e, & f_e(d_e)&=z,\\
			g_e(v_e)&=c_e^1, & g_e(z)&=c_e^*,\\
			g_e(w)&=c_e^0 &&(w\in V\setminus\{v_e\}).
		\end{aligned}
		\]
		Compare each input and output at $L$. Both functions preserve
		the interface value, which proves \textup{(H2)}. Give $A_e,B_e,C_e$
		formation ranks $0,1,2$, respectively, and give rank zero to any
		remaining formation types. This proves \textup{(H3)}. Select all
		states. The only trees of depth two are the chains
		$A_e\to B_e\to C_e$. All state sets are nonempty and both maps
		in each chain are total, so every formation tree admits a witness.
		
		For the native realization, use the empty algebraic signature
		with no identities. At each site take one generator name per
		state. The free algebra on these names is the set of names itself,
		so the empty equation list realizes the marked state set under
		\cref{eq:equational-native-law}. To realize a site quotient, add
		the equations identifying names in the same kernel class. The
		resulting marked set is precisely the quotient, with each original
		name marking its quotient class. The associated observation
		presentation retains the original bases and carries the restriction
		and formation tables.
		
		\begin{enumerate}[label=\textup{(\roman*)},leftmargin=*]
			\item \emph{Classification by coloring partitions.}
			The identity readouts force the quotients at $L,A_e,C_e$ to
			be bijective. At $B$ and each $B_e$, the binary readout forces
			$\{z\}$ to be a singleton class. It therefore suffices to
			specify the kernel equivalence $\theta$ on $V$ at $B$ and
			the kernel equivalence $\theta_e$ on its copy at each $B_e$.
			The copy restriction $B\to B_e$ descends exactly when
			$\theta\subseteq\theta_e$. All restrictions to $L$ then
			descend because their values are constant on each class.
			
			Single-vertex and one-step trees are exact by the local image
			construction. For the chain $A_e\to B_e\to C_e$,
			\cref{lem:foundational-two-step-quotient} says that exactness
			holds precisely when $g_e$ is constant on each quotient class
			meeting $f_e(A_e)=\{u_e,z\}$. The class of $z$ is a singleton.
			On $V$, the only vertex whose $g_e$-value differs from that
			of $u_e$ is $v_e$. Hence the quotient family is exact if and
			only if
			\begin{equation}\label{eq:foundational-chromatic-classification}
				\theta\subseteq\theta_e
				\qquad\text{and}\qquad
				(u_e,v_e)\notin\theta_e
				\qquad(e\in E).
			\end{equation}
			
			It follows that $\theta$ separates the endpoints of every
			edge, so its classes form a coloring partition. Conversely,
			given any coloring partition $\theta$, set $\theta_e=\theta$
			for every edge. This gives an exact family by
			\eqref{eq:foundational-chromatic-classification}. Restriction
			to the central site and this choice of section both preserve
			the order. For any nonadjacent $u,v$, the partition with
			$\{u,v\}$ as its only nonsingleton block is a coloring
			partition. Thus some exact quotient identifies them, whereas
			every exact quotient separates adjacent vertices.
			
			Finally, every retained readout factors through its site
			quotient. Since $q_T(W_T)=\overline W_T$ on every directed
			formation tree, the joint readout tuples on the two witness
			sets coincide. This preserves every determination judgment
			and every pair of readout tuples witnessing its failure.
			
			\item \emph{Minimum size and computational complexity.}
			If $\theta$ has $k$ color classes, then $B$ has $k+1$
			quotient states, including the class of $z$. Part~\textup{(i)}
			therefore gives \eqref{eq:foundational-chromatic-minimum}.
			At each $B_e$, the two endpoints of $e$ and the protected
			state $z$ require at least three classes. The fixed sites
			$L,A_e,C_e$ contribute $2+5m$ states in total. Thus
			\[
			\sum_U|\overline X_U|
			\geq (1+\chi(\Gamma))+3m+(2+5m)
			=\chi(\Gamma)+3+8m.
			\]
			To attain the bound, choose a coloring partition with
			$\chi(\Gamma)$ classes. At each $B_e$, take two classes on
			$V$: the central class containing $u_e$ and its complement.
			The complement contains $v_e$. These two classes coarsen
			$\theta$ and separate the endpoints, so
			\eqref{eq:foundational-chromatic-classification} holds
			simultaneously for all edges.
			
			The presentation has $O((m+1)(|V|+1))$ entries in its
			explicit tables. By \eqref{eq:foundational-chromatic-minimum},
			a central bound of $k+1$ is feasible exactly when $\Gamma$
			has a proper coloring with at most $k$ colors. This gives a
			polynomial reduction from graph coloring with a specified
			color bound \cite{Karp1972}. Conversely, a coloring is a
			polynomial certificate: one checks its number of colors and
			the colors at the endpoints of each edge. Part~\textup{(i)}
			then supplies an exact quotient. This proves membership in
			$\NP$ and hence $\NP$-completeness on the stated family.
			
			\item \emph{Greatest quotients and the effect of retaining copy
				bijections.}
			By part~\textup{(i)}, a greatest central partition must
			identify each nonadjacent pair and separate each adjacent
			pair. Its equivalence relation must therefore be
			\[
			u\sim v
			\quad\Longleftrightarrow\quad
			u=v\ \text{or}\ \{u,v\}\notin E.
			\]
			If this is an equivalence relation, each class is independent
			and every two vertices in distinct classes are adjacent.
			Thus $\Gamma$ is complete multipartite. Conversely, in a
			complete multipartite graph every independent set lies in one
			part, so the partition into parts is the greatest coloring
			partition.
			
			For a greatest family over all sites, suppose there is an edge
			$e=\{u,v\}$ and a third vertex $w$. Keep the central partition
			and all edge-copy partitions other than $\theta_e$ discrete.
			At $B_e$, one exact family merges only $u,w$, and another
			merges only $v,w$. A common upper bound would identify
			$u$ and $v$ at $B_e$, contrary to
			\eqref{eq:foundational-chromatic-classification}. Hence no
			greatest family exists. If $E=\varnothing$, the one-block
			partition on $V$ is greatest. If $E\ne\varnothing$ and
			$|V|\leq2$, the graph consists of a single edge and both
			partitions on $V$ are forced to be discrete. This proves
			the criterion for the original contract.
			
			Now require every copy restriction to remain bijective.
			The induced map from the quotient at $B$ to the quotient
			at $B_e$ is surjective because the original copy map is
			bijective. It is injective exactly when
			$\theta=\theta_e$. Thus the central coloring partition
			determines the entire quotient family, and the section in
			part~\textup{(i)} is the inverse of the central-partition
			map. The criterion for a greatest quotient follows. With
			$k$ color classes, the $m+1$ sites $B,B_e$ each contribute
			$k+1$ states, and the fixed sites contribute $5m+2$.
			Taking $k=\chi(\Gamma)$ gives the asserted minimum.
			
			\item \emph{Irreducibility and the gap from minimum size.}
			By part~\textup{(iii)}, an exact quotient with the copy
			bijections retained is irreducible exactly when its coloring
			partition has no proper coloring coarsening. Two classes can
			be merged exactly when there is no edge between them.
			Moreover, any proper coloring coarsening merges at least
			one such pair. This proves the stated criterion.
			
			For each $r\geq2$, take a disjoint union of $\binom r2$ edges,
			one for each unordered pair $\{i,j\}$ of distinct elements
			of $\{1,\ldots,r\}$. Color the endpoints of the edge indexed
			by $\{i,j\}$ with colors $i$ and $j$. This is a proper
			coloring with $r$ classes, and each pair of classes has an
			edge between them. Its quotient is therefore irreducible
			and has $r+1$ central states. The graph is a nonempty
			matching, so its chromatic number is two and the minimum
			central size is three. The ratio $(r+1)/3$ is unbounded.
			\qedhere
		\end{enumerate}
	\end{proof}
	
	\begin{corollary}[Formation-depth threshold for exact quotient minimization]
		\label{cor:foundational-depth-boundary}
		\begin{enumerate}[label=\textup{(\roman*)},leftmargin=*]
			\item \emph{Depth at most one.}
			If all ordinary directed formation trees of a presentation have depth
			at most one, every readout-preserving local image reduction is
			witness-exact.
			\item \emph{Depth two.}
			For presentations of depth at most two, minimizing the number of
			states at a distinguished site is already $\NP$-hard, even with the
			coherence, total determinism, nonempty original witnesses and fixed
			native marked-set law in \cref{thm:foundational-chromatic-reductions}.
		\end{enumerate}
	\end{corollary}
	
	\begin{proof}\leavevmode
		\begin{enumerate}[label=\textup{(\roman*)},leftmargin=*]
			\item \emph{Depth at most one.}
			A depth-zero tree has the selected state set as its witness set, and a
			depth-one tree has its one primitive admission relation as its witness
			set. Both images are exact by construction.
			\item \emph{Depth two.}
			The depth-two examples of
			\cref{thm:foundational-chromatic-reductions} give the stated hardness.
		\end{enumerate}
	\end{proof}
	
	\begin{remark}[Colorings and the quotient contract]
		\label{rem:foundational-chromatic-scope}
		The image of a total deterministic formation under a local quotient can
		be multivalued. This freedom permits the exact coloring construction.
		
		With the copy bijections protected, the maximal partitions are the
		classical complete colorings. Thus the largest irreducible central quotient
		has $1+\psi(\Gamma)$ states, where $\psi$ is the achromatic number
		\cite{HararyHedetniemi1970}. Classical minimization of incompletely
		specified finite-state machines also leads to hardness results
		\cite{Pfleeger1973,AbelReineke2015}.
	\end{remark}
	
	\section{Existence and faithfulness of amalgamations}
	\label{sec:native-amalgamation}
	
	When packages are assembled along a tree, each shared interface constrains
	which local states can occur together.
	A different separator may be used at each edge.
	Equal separator images then guarantee that the full matching source covers
	every chart.
	Faithfulness also requires a pair of matching states for each dependency excluded by the local rules.
	We establish this criterion under the spectral and equational formation laws, then ask which interfaces preserve the requisite pairs for every permitted
	attachment.
	The fixed-atlas existence test is \cref{thm:amalgamation-existence-criterion}. The universal attachment criterion is \cref{thm:congruence-interface-saturation}.
	
	The joint-boundary criteria in \cref{sec:compositional-invariants} extend matching-source tests to arbitrary finite group graphs, using complete reachable boundaries rather than independent separator values.
	
	\subsection{Spectral native interpretations and preservation contracts}
	\label{subsec:spectral-native-interpretation}
	
	In a finite spectral atlas, a descriptor is a function on a local state set.
	Its fibers determine both the functional dependencies among observations and the represented algebra generated by them. \Cref{prop:spectral-observation-interpretation} identifies these two descriptions.
	We use this identification to state precisely what an assembled source must preserve.
	
	\paragraph{The access--commutant realization law.}
	\label{subsec:standalone-commutant-context}
	The spectral constructions use the following realization law
	\cite[Section~2.15]{Ma2026}.
	For a nonzero finite-dimensional Hilbert space $\mathcal H$, let $\mathscr C_{\mathcal
		H}$ be the concrete unital $*$-subalgebras of $B(\mathcal H)$.
	For a set $I$, write $\mathbf{Ref}(I)=\{q_\theta:I\to I/\theta:\theta\text{ is an
		equivalence relation on }I\}$ for canonical refinement representatives.
	For $A\in\mathscr C_{\mathcal H}$, let $\operatorname{Und}(A)$ be its underlying set and put
	\[
	\mathfrak P_A=
	(\{A\},\operatorname{Und}(A),
	\operatorname{Und}(A)\times\{A\},\LL_A,\RR_A),
	\qquad
	\LL_A(A)=\RR_A(A)=\mathbf{Ref}(\operatorname{Und}(A)).
	\]
	Both the structural domain and the admissible carriers are $\mathscr C_{\mathcal H}$; the formation-data space consists of all
	$(\mathrm{com},M)$ with $M\in\mathscr C_{\mathcal H}$.
	The realization relation contains exactly the pairs given by
	\begin{equation}\label{eq:standalone-commutant-law}
		\Xi_{\mathcal H}^{\mathrm{com}}\vdash
		(\mathrm{com},M)\rightsquigarrow\mathfrak P_{M'},
		\qquad
		M'=\{T\in B(\mathcal H):Tm=mT\ (m\in M)\}.
	\end{equation}
	The represented carrier and the ambient space are retained.
	Thus $M''=M$ recovers the concrete source, and $(M\vee N)'=M'\cap N'$ is the native composition law.
	
	\paragraph{Spectral observations of commutant sources.}
	For the commutant law, the selected family at a carrier $K$ is $\mathbf{Ref}(\operatorname{Und}(K))$.
	The finite commuting constructions use characters of the recovered commutative source $K'$ as their observation states.
	
	\begin{proposition}[Finite spectral interpretation of the native law]
		\label{prop:spectral-observation-interpretation}
		Let $H$ be a finite nonempty set with readouts $f_d:H\to D_d$ for a finite descriptor set $V$, and suppose their joint map is injective.
		Write $f_S=(f_d)_{d\in S}$ and
		\[
		c_H(S)=\{d\in V:\ f_S(x)=f_S(x')\Longrightarrow
		f_d(x)=f_d(x')\text{ for all }x,x'\in H\}.
		\]
		Represent $\mathbb C^H$ by multiplication on $\mathcal H_H=\ell^2(H)$.
		For $S\subseteq V$ put
		\[
		A_d=\mathbb C[f_d],\qquad
		M_H(S)=\operatorname{alg}^*\{A_d:d\in S\},\qquad
		K_H(S)=M_H(S)'.
		\]
		Here $\mathbb C[f_d]=\{g\circ f_d:g:f_d(H)\to\mathbb C\}$ is represented by multiplication and $M_H(\varnothing)=\mathbb C1$.
		\begin{enumerate}[label=\textup{(\roman*)},leftmargin=*]
			\item \emph{Represented source and spectrum.}
			The datum $(\mathrm{com},M_H(S))$ has the native endpoint $\mathfrak P_{K_H(S)}$ in $\Xi_{\mathcal H_H}^{\mathrm{com}}$.
			\begin{equation}\label{eq:spectral-observation-source}
				K_H(S)'=M_H(S),\qquad \operatorname{Spec}(M_H(S))\cong f_S(H).
			\end{equation}
			In particular, $\operatorname{Spec}(M_H(V))\cong H$ with all descriptor readouts retained.
			\item \emph{Recovery and composition.}
			For $S,T\subseteq V$, one has
			\begin{equation}\label{eq:spectral-observation-recovery}
				d\in c_H(S)\ \Longleftrightarrow\ A_d\subseteq M_H(S),
				\qquad K_H(S\cup T)=K_H(S)\cap K_H(T).
			\end{equation}
			\item \emph{Pullback along a chart map.}
			If $\pi:H\to X$ is a map to a finite local spectrum, its source map is the unital homomorphism $\pi^*:\mathbb C^X\to\mathbb C^H$, and $\pi^*$ is injective exactly when $\pi$ is surjective.
		\end{enumerate}
	\end{proposition}
	
	\begin{proof}
		\leavevmode
		\begin{enumerate}[label=\textup{(\roman*)},leftmargin=*]
			\item \emph{Represented source and spectrum.}
			The algebra generated by the $A_d$ for $d\in S$ consists exactly
			of functions constant on the fibers of $f_S$. Indeed, products of
			the indicator functions of the individual readout fibers give the
			indicator of each nonempty joint fiber. Thus
			$M_H(S)\cong\mathbb C^{f_S(H)}$, whose characters are evaluations
			at those joint values. On the orthogonal decomposition of
			$\ell^2(H)$ into its $f_S$-fibers, the commutant is the direct sum
			of the full matrix algebras on the fibers; its commutant again
			consists of their scalar blocks. This proves
			\eqref{eq:spectral-observation-source} and the stated native
			judgment by \eqref{eq:standalone-commutant-law}.
			\item \emph{Recovery and composition.}
			The inclusion $A_d\subseteq M_H(S)$ says precisely that $f_d$
			is constant on every $f_S$-fiber, which is the definition of
			$d\in c_H(S)$. Consequently every non-recovery counterexample
			pair in $H$ remains the same pair of source observations.
			Taking the commutant of the algebra generated by two families
			gives the intersection of their commutants, proving the second
			identity in \eqref{eq:spectral-observation-recovery}.
			\item \emph{Pullback along a chart map.}
			$\pi^*$ preserves pointwise operations. It is injective when
			$\pi$ is onto; if $x\notin\pi(H)$, the nonzero indicator of
			$\{x\}$ lies in its kernel, proving the converse.\qedhere
		\end{enumerate}
	\end{proof}
	
	For two partial diagrams with interface maps $p_L,p_R$ to $Y$, the formation relation on their spectral observations is
	\[
	N_{L,R;Y}
	=\{(x_L,x_R,(x_L,x_R)):
	p_L(x_L)=p_R(x_R)\},
	\qquad X_{L\cup R}=X_L\times_Y X_R.
	\]
	The input and output types are the indicated partial diagrams, and restrictions are tuple projections.
	
	A finite spectral atlas specifies the local states and their readouts. State coverage on a join tree with equal separator images follows from the classical join-tree theorem \cite{Beeri1983}. The preservation contract below also asks which dependencies hold after assembly.
	
	\begin{definition}[Finite native spectral atlas]
		\label{def:spectral-atlas}
		Fix a finite descriptor set $V$, a value set $D_e$ for each $e\in V$, and a cover $(B_i)_{i\in I}$ of $V$ indexed by a finite nonempty set $I$.
		A \textbf{finite native spectral atlas} on this cover consists of one chart for each $i\in I$, equipped with the native input base determined by that chart.		
		\begin{itemize}
			\item \emph{Chart observables and overlap maps.}
			For each $i\in I$, specify a finite nonempty state set $X_i$ and observables $f_{i,e}:X_i\to D_e$ for $e\in B_i$.
			For $S\subseteq B_i$, write
			\[
			f_{i,S}:X_i\longrightarrow\prod_{e\in S}D_e,
			\qquad f_{i,S}(x)=(f_{i,e}(x))_{e\in S}.
			\]
			Require the full chart map $f_i=f_{i,B_i}$ to be injective.
			On each overlap $S_{ij}=B_i\cap B_j$, the maps $f_{i,S_{ij}}$ and $f_{j,S_{ij}}$ take values in the common set $\prod_{e\in S_{ij}}D_e$.
			
			\item \emph{Canonical native input base.}
			For each $i\in I$, let $A_i=\mathbb C^{X_i}$ act by multiplication on $\ell^2(X_i)$, and let $A_i'$ denote its commutant in this representation.
			The realized input base is
			\[
			\bigl(\Xi_{\ell^2(X_i)}^{\mathrm{com}},
			(\mathrm{com},A_i),\mathfrak P_{A_i'}\bigr).
			\]
			Under the spectral interpretation of
			\cref{prop:spectral-observation-interpretation}, each state
			$x\in X_i$ corresponds to the evaluation character $h\mapsto h(x)$ of $A_i$.
		\end{itemize}
	\end{definition}
	
	\paragraph{Joint spectrum and fixed-atlas subamalgams.}
	A tuple of local states is compatible if its two chart values agree	on every overlap.
	The \textbf{full finite spectrum} of the atlas is the natural join
	\begin{equation}\label{eq:native-joint-spectrum}
		X=\left\{
		(x_i)_{i\in I}\in\prod_{i\in I}X_i
		\;\middle|\;
		f_{i,S_{ij}}(x_i)=f_{j,S_{ij}}(x_j)
		\text{ for all }i,j\in I
		\right\}.
	\end{equation}
	Let $\pi_i:X\to X_i$ be the coordinate projection.
	When $X\ne\varnothing$, the common source is $A=\mathbb C^X$, and the chart maps induce the pullback homomorphisms
	\[
	\pi_i^*:A_i\longrightarrow A,
	\qquad \pi_i^*(h)=h\circ\pi_i.
	\]
	
	A \textbf{fixed-atlas subamalgam} is a nonempty subset $R\subseteq X$.
	It is \textbf{state-exact} if every local state occurs in a tuple of $R$,
	that is, if $\pi_i(R)=X_i$ for every $i\in I$.
	
	For finite maps $X_i\to Y$, give $\mathbb C^{X_i}$ its $\mathbb C^Y$-algebra structure by pullback.
	In the balanced tensor product
	$\mathbb C^{X_1}\otimes_{\mathbb C^Y}\mathbb C^{X_2}$, tensors of point
	indicators at $x_1,x_2$ vanish when their $Y$-values differ:
	the indicator of the first value acts as the identity on the first factor and as zero on the second.
	The remaining tensors map to the point indicators of $X_1\times_Y X_2$ and give an algebra isomorphism with $\mathbb C^{X_1\times_Y X_2}$.
	An empty fiber product gives the zero algebra.
	
	A \textbf{join tree} on the chart indices is a tree $J$ such that, for every descriptor $d$, the charts containing $d$ form a connected subtree.
	Equal realized separator images on an edge are an additional input condition.
	
	\paragraph{The finite-atlas preservation contract.}
	\label{subsec:preservation-contract}
	Fix a finite spectral atlas as in \cref{def:spectral-atlas}.
	Write $c_i$ for the determination closure of its actual local state space $X_i$, and $c_0$ for the least closure generated by all local implications.
	A compatible candidate source is a nonempty finite set $H$ with maps $\pi_i:H\to X_i$ whose descriptor readouts agree on overlaps.
	Thus $f_d=f_{i,d}\pi_i$ is independent of the chart containing $d$.
	Let $c_H$ denote their determination closure.
	For $T\subseteq V$, let $\theta_T=\bigcap_{e\in T}\ker f_e$ denote the equivalence relation of joint agreement (the empty intersection is $H\times H$).
	In the finite setting this gives the equivalent tests
	\[
	d\in c_H(T)\quad\Longleftrightarrow\quad
	\theta_T\subseteq\ker f_d\quad\Longleftrightarrow\quad
	\mathbb C[f_d]\subseteq\mathbb C[f_T],
	\]
	where the function algebras are pulled back to $H$. The last
	equivalence holds because functions on a finite image separate its
	points. On the joint observation image this is exactly the ordinary
	functional-dependency closure \cite{Armstrong1974}; a failed inclusion
	is witnessed by a pair in $\theta_T\setminus\ker f_d$. For group
	homomorphisms the relation inclusion becomes inclusion of the
	corresponding normal kernels. For scalar linear forms on a vector
	space it becomes membership of $f_d$ in $\operatorname{span}\{f_e:e\in T\}$,
	by annihilator duality, the represented matroid closure \cite{Oxley2011}.
	For a conditional closure the boundary readout is also included in the premise.
	The native comparison specifies whether states with equal joint readouts may remain distinct.
	
	Whenever an edge $e=ij$ is used in the separator contract, assume that $f_{i,S_{ij}}(X_i)=f_{j,S_{ij}}(X_j)$ and denote this common actual image by $Y_e$, where $S_{ij}=B_i\cap B_j$.
	If a join tree $J$ is specified, deleting an edge $e$ gives two	entire chart components $L_e,R_e$.
	For a set of charts $I'$, let $M_{I'}(H)$ be the algebra of all complex-valued functions of the joint readout $(\pi_i)_{i\in I'}$.
	Let $p_e:H\to Y_e$ be the induced separator readout and put $A_e(H)=p_e^*\mathbb C^{Y_e}$, as an image algebra in $\mathbb C^H$.
	If all local maps are onto, these separator pullbacks are faithful embeddings.
	
	\begin{definition}[Preservation contract for a finite spectral atlas]
		\label{def:preservation-contract}
		Fix a finite spectral atlas.
		A \textbf{preservation contract} specifies a set of required clauses from the list below.
		A compatible candidate source $H$ with chart maps $\pi_i:H\to X_i$ satisfies the contract if it satisfies every required clause:	
		\begin{description}
			\item[$\mathsf C$] \emph{Chart-closure exactness.}
			For every chart $i\in I$ and every $S\subseteq B_i$,
			\[
			c_H(S)\cap B_i=c_i(S).
			\]
			
			\item[$\mathsf S$] \emph{Local-state surjectivity.}
			Every local state occurs in the candidate source:
			\[
			\pi_i(H)=X_i \qquad (i\in I).
			\]
			The chart maps use the specified observable values and
			interface identifications. This condition is also called
			\emph{state-exact amalgamation}.
			
			\item[$\mathsf F$] \emph{Dependency faithfulness.}
			For every $S\subseteq V$,
			\[
			c_H(S)=c_0(S),
			\]
			where $c_0$ is the least closure generated by the local
			implications.
			
			\item[$\mathsf E$] \emph{Exact represented separators.}
			For every edge $e$ of the specified join tree, let $L_e$
			and $R_e$ be the two chart components obtained by deleting
			$e$. The algebras of the two full components satisfy
			\[
			M_{L_e}(H)\cap M_{R_e}(H)=A_e(H)
			\qquad\text{inside }\mathbb C^H.
			\]
			
			\item[$\mathsf N$] \emph{Native realization.}
			The candidate is equipped with a formation datum $P_H$,
			an endpoint $\mathfrak P_H$, and a realization judgment
			\[
			\Xi\vdash P_H\rightsquigarrow\mathfrak P_H
			\]
			in the specified native context. The realization supplies
			the required source representation, retained markings,
			structural maps, and map interpreting the output observations.
		\end{description}		
		The \emph{full faithful native contract} requires
		$\mathsf S$, $\mathsf F$, $\mathsf E$, and $\mathsf N$.
	\end{definition}
	
	The chartwise optimization problem requires $\mathsf C$ alone.
	
	\paragraph{Transfer through the joint chart image.}
	Fix the observed atlas, including the actual chart and separator values,
	and put $T=\psi(H)\subseteq W$, where $\psi$ is the joint chart map
	and $W$ is the full matching set. Every retained descriptor and complete-side
	observation factors through $H\twoheadrightarrow T$. By
	\cref{prop:observation-noncollapse,thm:universal-observable-quotient},
	$c_H=c_T$; chart images are also unchanged. Pullback identifies each
	function algebra of these observations on $T$ with its counterpart on $H$.
	Since pullback is injective, it preserves and reflects their intersections
	and equalities. Thus all four observation clauses
	$\mathsf C,\mathsf S,\mathsf F,\mathsf E$ have the same truth values on
	$H$ and $T$ (for $\mathsf E$ in the displayed finite-function formulation).
	Sources with the same joint chart image can differ in their native realization, representation multiplicities, and hidden kernels. Those requirements must be checked on the represented sources.
	
	\paragraph{Logical relations among the clauses.}
	Compatibility of the maps already gives $c_0\leq c_H$, because every local implication holds after pullback.
	Surjectivity supplies all local counterexample pairs, so
	$
	\mathsf S\Longrightarrow\mathsf C.
	$
	Whenever $c_0$ has the prescribed local traces, $c_0(S)\cap B_i=c_i(S)$ for $S\subseteq B_i$, one also has $\mathsf F\Longrightarrow\mathsf C$.
	\Cref{prop:closure-descent} proves the trace condition for join trees with equal overlap images.
	Chart-closure exactness and dependency faithfulness can both hold even when a chart state is missing.
	For one chart $X_1=\{00,01,10,11\}$ on two binary descriptors, the proper subset
	$H=\{00,01,10\}$ has the same free determination closure and misses one actual local state.
	
	The full join on a compatible join tree satisfies $\mathsf S$ and $\mathsf E$ before any pair-lifting hypothesis is imposed.
	The join in \cref{ex:correlated-join} also admits the constructed commutant realization, so it satisfies $\mathsf C,\mathsf S, \mathsf E,\mathsf N$ while failing $\mathsf F$.
	The source defined by the parity constraint in
	\cref{rem:native-group-adjacent-insufficient} is native, state-exact,
	and dependency-faithful. Its adjacent-chart tests all pass, yet its complete sides have extra common information at a full cut.
	
	For finite Mal'cev join-tree atlases with onto interfaces and
	$c_W=c_0$, \cref{thm:complete-contract-logic} realizes exactly the nine
	truth assignments allowed by $\mathsf S\Rightarrow\mathsf C$,
	$\mathsf F\Rightarrow\mathsf C$ and $\mathsf S\wedge\mathsf E\Rightarrow\mathsf F$
	on one fixed module atlas. Within these already native candidates,
	$\mathsf S$ and $\mathsf E$ together suffice, and neither alone suffices for
	the full contract. For scalar linear charts over a common finite field, jointly separating observations give
	$\mathsf F\Longleftrightarrow\mathsf S\wedge\mathsf E$
	(\cref{cor:scalar-contract-equivalence}).
	The dependencies among the clauses therefore depend on the native class.

	\paragraph{Preservation properties of the amalgamations.}
	A check mark indicates that the construction satisfies the clause under its stated hypotheses.
	The notation ``RPL'' marks relative pair lifting as a sufficient condition for faithfulness. The exact test for a fixed atlas is
	\cref{thm:amalgamation-existence-criterion}.
	``Com'' uses the access--commutant law on the constructed ambient space;
	``Pres'' uses the fixed marked presentation law.
	\begin{center}
		\small
		\begin{tabular}{@{}>{\raggedright\arraybackslash}p{0.47\textwidth}ccccc@{}}
			\toprule
			Construction & $\mathsf C$ & $\mathsf S$ & $\mathsf F$
			& $\mathsf E$ & $\mathsf N$ \\
			\midrule
			Spectral pair lifting, \cref{thm:local-pair-native}
			& $\checkmark$ & $\checkmark$ & $\checkmark$ & $\checkmark$ & Com \\
			Finite equational gluing, \cref{thm:equational-native-amalgamation}
			& $\checkmark$ & $\checkmark$ & RPL & $\checkmark$ & Pres \\
			Joint two-state equational gluing, \cref{cor:two-state-equational-descent}
			& $\checkmark$ & $\checkmark$ & $\checkmark$ & $\checkmark$ & Pres \\
			Complete three-permutable interfaces, \cref{thm:complete-port-native-descent}
			& $\checkmark$ & $\checkmark$ & $\checkmark$ & $\checkmark$ & Pres \\
			Visible scalar-linear interfaces, \cref{thm:linear-native-amalgamation}
			& $\checkmark$ & $\checkmark$ & $\checkmark$ & $\checkmark$ & Com/Pres \\
			Fixed terminal-tree atlas,
			\cref{thm:native-spectral-tree,thm:exact-native-separator}
			& $\checkmark$ & $\checkmark$ & $\checkmark$ & $\checkmark$ & Com \\
			\bottomrule
		\end{tabular}
	\end{center}
	The scalar-linear construction has both native realizations.
	\Cref{thm:fixed-ambient-criterion} gives the separate criterion for preserving preassigned Hilbert-space multiplicities.
	
	\subsection{Obstructions and the exact two-witness criterion}
	\label{subsec:early-obstructions}
	
	A join can cover every local state and still lose a pair that
	witnesses a failed dependency. The loss introduces an extra dependency.
	Deleting states or changing positive weights can leave it in force.
	
	\begin{example}[State consistency does not preserve negative recovery]
		\label{ex:correlated-join}
		Consider two charts $B_1=\{A,B\}$ and $B_2=\{B,C\}$, with
		\[
		R_{AB}=\{(0,0),(0,1),(1,2),(1,3)\},\qquad
		R_{BC}=\{(0,0),(1,0),(2,1),(3,1)\}.
		\]
		The first chart has exactly the nontrivial dependency $B\to A$ and the
		second exactly $B\to C$. The overlap values $B=0,1,2,3$ match
		bijectively, so the full join has four sectors and both chart
		projections are onto. Yet globally $A=C$ in every row, hence
		$A\to C$ and $C\to A$. These are extra dependencies beyond the two local rules. They
		hold on every subset of the join because deleting rows preserves
		valid dependencies. Thus the join covers each individual chart state
		and still lacks the required counterexample pairs.
	\end{example}
	
	\begin{proposition}[Deletion cannot restore a nonimplication]
		\label{prop:subamalgam-monotonicity}
		For nonempty relations $R\subseteq X$ on the same descriptors, $c_X(S)\subseteq c_R(S)$ for every $S$.
		Consequently, if the full fixed-atlas join has $c_X\not\subseteq c_0$, no nonempty
		subamalgam of that join has determination closure $c_0$.
	\end{proposition}
	
	\begin{proof}
		Every pair of rows in $R$ is a pair in $X$, so each dependency valid on $X$ remains valid on $R$.
		An unintended dependency of $X$ therefore survives every row deletion.
	\end{proof}
	
	\begin{proposition}[Positive annotations cannot repair a support obstruction]
		\label{prop:positive-weight-support-obstruction}
		Let $K=(K,+,0)$ be a nontrivial positive commutative monoid, where $a+b=0$ implies $a=b=0$.
		Let $(R_i)_{i\in I}$ be a nonempty finite family of finitely supported $K$-valued
		relations on descriptor sets $B_i$, with nonempty supports $X_i$, and put
		$V=\bigcup_iB_i$.
		Marginalization sums the annotation values of all tuples with the same restriction.
		Let $T$ be a finitely supported $K$-relation on $V$ with marginals $R_i$.
		\begin{enumerate}[label=\textup{(\roman*)},leftmargin=*]
			\item \emph{Marginal supports.}
			The support of $T$ projects onto every local support and lies in their full join:
			\begin{equation}\label{eq:weighted-support-inclusion}
				\pi_i(\operatorname{supp}T)=X_i\quad(i\in I),\qquad
				\operatorname{supp}T\subseteq W:=\Join_i X_i.
			\end{equation}
			\item \emph{Dependency obstruction.}
			For every $S\subseteq V$,
			\[
			c_W(S)\subseteq c_{\operatorname{supp}T}(S).
			\]
			Consequently, if the full support join has a dependency not in the closure $c_0$
			generated by the local support theories, no such
			weighted marginal witness is recovery-faithful for that support semantics.
			This holds even when $K$ has the transportation property.
		\end{enumerate}
	\end{proposition}
	
	\begin{proof}
		\leavevmode
		\begin{enumerate}[label=\textup{(\roman*)},leftmargin=*]
			\item \emph{Marginal supports.}
			A finite sum in $K$ is nonzero exactly when at least one summand
			is nonzero. Thus the support of a marginal is the projection of
			the support, giving the first equality. Every tuple in
			$\operatorname{supp}T$ therefore belongs to the natural join $W$.
			\item \emph{Dependency obstruction.}
			The closure inclusion follows from
			\cref{prop:subamalgam-monotonicity}.\qedhere
		\end{enumerate}
	\end{proof}
	
	The support identity is the standard positive-monoid observation
	\cite[Lemma~1]{AtseriasKolaitis2025}; the proposition records its elementary consequence for the present faithfulness contract.
	For \cref{ex:correlated-join}, give every displayed local tuple the	same nonzero annotation $k\in K$.
	The four matching global tuples, each annotated by $k$, have exactly these marginals.
	Every tuple of nonzero weight still satisfies $A=C$.
	Thus the extra dependencies persist under positive weighting.

	\subsubsection{Exact witness descent on a join tree}

	On a join tree $J$, running intersection allows matching local ordered pairs to be assembled into a pair of global states.
	First assume that the two charts at every edge $ij$ have equal images on $S_{ij}=B_i\cap B_j$.

	\begin{proposition}[State coverage and local recovery closures on a join tree]
		\label{prop:closure-descent}
		\leavevmode
		\begin{enumerate}[label=\textup{(\roman*)},leftmargin=*]
			\item \emph{State coverage.}
			Under the stated join-tree hypotheses, the full join $X$ is nonempty and projects surjectively onto each $X_i$.
			
			\item \emph{Closure and local traces.}
			The induced closure satisfies $c_0\le c_X$, with exact local traces
			\begin{equation}\label{eq:closure-chart-exact}
				c_0(S)\cap B_i=c_i(S)\qquad(S\subseteq B_i).
			\end{equation}
			The inclusion $c_0\le c_X$ can be strict.
		\end{enumerate}
	\end{proposition}

	\begin{proof}
		\leavevmode
		\begin{enumerate}[label=\textup{(\roman*)},leftmargin=*]
			\item \emph{State coverage.}
			Root the tree at $i$.
			Starting from any $x_i\in X_i$, extend the separator value of each parent to a child state.
			Equality of the two	separator images permits each extension.
			Running intersection implies that charts sharing a descriptor agree on that descriptor along their entire connecting path.
			The resulting assignment lies in $X$.
			Thus every chart projection is onto, and $X$ is nonempty.
			\item \emph{Closure and local traces.}
			Every local implication holds after pullback, so $c_0\le c_X$.
			For	$S\subseteq B_i$, a dependency in $X$ with premise $S$ and target in
			$B_i$ holds in $X_i$:
			lift any two chart states to global states by surjectivity.
			Consequently
			\[
			c_i(S)\subseteq c_0(S)\cap B_i\subseteq c_X(S)\cap B_i=c_i(S).
			\]
			The first inclusion uses the definition of $c_0$.
			The strict global inclusion occurs in \cref{ex:correlated-join}.\qedhere
		\end{enumerate}
	\end{proof}
	
	To test the missing inclusion $c_X\le c_0$, it is essential to synchronize
	two states.
	For $S\subseteq V$ and $e\in V$, define the local ordered pair constraint
	\begin{equation}\label{eq:pair-local}
		P_i(S,e)=
		\left\{(x,x')\in X_i^2:
		\begin{gathered}
			f_{i,s}(x)=f_{i,s}(x')\ (s\in S\cap B_i),\\
			f_{i,e}(x)\ne f_{i,e}(x')\ \text{if }e\in B_i
		\end{gathered}\right\}.
	\end{equation}
	A pair restricts to the ordered pair of actual separator value tuples,
	$\partial_{ij}(x,x')=(p_{i,S_{ij}}(x),p_{i,S_{ij}}(x'))\in Y_{ij}^2$,
	where $Y_{ij}$ is the common realized separator image.
	For the pair test alone the two chart images may differ.
	In that case take $Y_{ij}$ to be their union inside the common descriptor product and use each chart's map into this set.
	Root $J$ at $r$ and define bottom-up messages on oriented edges $i\to j$ (with $j$ the parent):
	\begin{equation}
		\label{eq:pair-message}
		\mathscr M_{i\to j}^{S,e}=
		\left\{\partial_{ij}(x,x'):(x,x')\in P_i(S,e),\quad
		\partial_{ik}(x,x')\in\mathscr M_{k\to i}^{S,e}
		\ (k\text{ a child of }i)\right\}.
	\end{equation}
	The root test accepts if a pair in $P_r(S,e)$ matches all its incoming messages.
	In particular, a message records two exact interface value tuples.
	An agreement/disagreement mask alone is generally insufficient.
	
	\begin{theorem}[Ordered-pair test for dependency faithfulness on a join tree]
		\label{thm:pair-descent}
		Let $J$ be a join tree for a finite atlas with descriptor set $V$ and full join $X$.
		Fix any root of $J$.
		For each edge $ij$, let	$Y_{ij}$ be the union of the two realized separator images in the common descriptor product.		
		\begin{enumerate}[label=\textup{(\roman*)},leftmargin=*]
			\item \emph{Witness for a single nondependency.}
			For every $S\subseteq V$ and $e\in V$, the root test for
			$(S,e)$ accepts if and only if there exist $x,x'\in X$
			such that
			\[
			x|_S=x'|_S
			\qquad\text{and}\qquad
			x_e\ne x'_e.
			\]
			In particular, when $X\ne\varnothing$, the test accepts	exactly when $e\notin c_X(S)$.
			
			\item \emph{Faithfulness of the full join.}
			If $X\ne\varnothing$, then
			\begin{equation}\label{eq:faithful-pair-criterion}
				c_X=c_0
				\quad\Longleftrightarrow\quad
				\text{the root test accepts every $(S,e)$
					with $e\notin c_0(S)$.}
			\end{equation}
			
			\item \emph{Message size and number of premises.}
			For each test $(S,e)$ and each oriented edge $i\to j$,
			\[
			\mathscr M_{i\to j}^{S,e}\subseteq Y_{ij}^2,
			\qquad
			|\mathscr M_{i\to j}^{S,e}|\leq |Y_{ij}|^2.
			\]
			The test uses only local chart pairs and separator messages.
			It requires neither enumeration of $X$ nor precomputation of coatoms.
			For a fixed target $e$, the criterion in \textup{(ii)} can contain $2^{|V|-1}$ distinct premises.
		\end{enumerate}
	\end{theorem}
	
	\begin{proof}
		\leavevmode
		\begin{enumerate}[label=\textup{(\roman*)},leftmargin=*]
			\item \emph{Witness for a single nondependency.}
			Fix $(S,e)$.
			We first show that the message $\mathscr M_{i\to j}^{S,e}$ consists exactly of the ordered separator values that extend to a compatible choice of local	pairs throughout the subtree rooted at $i$, with the pair
			at each chart $k$ belonging to $P_k(S,e)$.
			
			At a leaf $i$, this follows from \eqref{eq:pair-message}:
			the message is the separator projection of $P_i(S,e)$.
			For the induction step, fix a pair in $P_i(S,e)$.
			By the induction hypothesis, it extends through a child subtree
			exactly when its ordered separator value belongs to that
			child's message.
			These extensions can be chosen simultaneously.
			Indeed, if a descriptor occurs in two different child subtrees,
			running intersection places it also at $i$ and at every chart
			on the paths connecting those occurrences to $i$. Both
			extensions therefore agree with the pair at $i$ on that
			descriptor, in each coordinate.
			Conversely, any compatible extension through the whole subtree restricts to an extension through each child. Projecting to the parent separator gives exactly \eqref{eq:pair-message}.
			
			At the root, acceptance thus gives a compatible choice of
			local pairs on the entire tree. Their first coordinates form
			a state $x\in X$, and their second coordinates form a state
			$x'\in X$. The conditions defining $P_i(S,e)$ give agreement
			on every descriptor in $S$ and disagreement on $e$.
			Conversely, any such $x,x'$ restrict to compatible local pairs
			and hence give acceptance. This argument also covers the
			empty join, for which no test accepts. When $X\ne\varnothing$,
			the displayed pair is precisely a witness to $e\notin c_X(S)$.
			
			\item \emph{Faithfulness of the full join.}
			Every implication valid on a chart remains valid after
			pullback along the projection $X\to X_i$. Since $c_0$ is
			generated by the chart implications, $c_0(S)\subseteq c_X(S)$
			for every $S\subseteq V$. Equality therefore holds for every
			$S$ exactly when
			\[
			e\notin c_0(S)\quad\Longrightarrow\quad e\notin c_X(S)
			\qquad(S\subseteq V,\ e\in V).
			\]
			Part~\textup{(i)} identifies the conclusion of each implication
			with acceptance of the corresponding root test. This proves
			\eqref{eq:faithful-pair-criterion}.
			
			\item \emph{Message size and number of premises.}
			Each entry of \eqref{eq:pair-message} is an ordered pair of
			separator values in $Y_{ij}$, giving the stated cardinality
			bound. The recursion computes each message from $P_i(S,e)$
			and the incoming child messages, so it uses only local data.
			
			For the premise count, take one singleton chart $\{d\}$ for
			each $d\in V$, with both binary values realized at each chart.
			Then $c_0(S)=S$ for every $S\subseteq V$. For a fixed target
			$e$, every subset of $V\setminus\{e\}$ is therefore a premise
			in \eqref{eq:faithful-pair-criterion}, giving
			$2^{|V|-1}$ distinct premises.\qedhere
		\end{enumerate}
	\end{proof}
	
	The state-coverage argument of \cref{prop:closure-descent} recovers
	the classical consistency implication for a join tree. The root
	messages of \cref{thm:pair-descent} use the same acyclic-join mechanism
	on the doubled chart relations $X_i^2$
	\cite{Beeri1983,Yannakakis1981,ComignaniEtAl2022}.
	
	\begin{corollary}[Faithfulness from extension of exact agreement pairs]
		\label{cor:relative-pair-homogeneity}
		Suppose each local closed set $C_i$ is the exact agreement set
		of an ordered pair in $X_i^2$. Also require the following at every
		oriented join-tree edge $i\to j$. Take any local closed sets
		$C_i,C_j$ with the same trace on $S_{ij}$. Every pair with
		agreement set $C_i$ must extend to a pair in $X_j^2$ with agreement
		set $C_j$ and the same ordered separator values. Then every global
		$c_0$-closed set is the exact agreement set of two rows of $X$.
		In particular $c_X=c_0$.
	\end{corollary}
	
	\begin{proof}
		For a globally closed $C$, the traces $C\cap B_i$ are locally closed.
		Choose a pair at the root with that exact agreement set.
		Extend it outward, at each edge using the stated pair-extension property with the prescribed closed agreement set in the child chart.
		Running intersection ensures global pair compatibility.
		The assembled pair agrees precisely	on $C$.
		Thus whenever $e\notin c_0(S)$, the globally closed set	$C=c_0(S)$ is realized as an agreement set; that pair agrees on $S$
		and disagrees on $e$.
		Apply \cref{thm:pair-descent}.
	\end{proof}
	
	The three-sector charts in \cref{subsec:tree-spectral-gluing} glue faithfully even though some locally closed sets never occur as exact agreement sets of distinct sectors. \Cref{ex:binary-xor-lifting} distinguishes exact agreement from agreement on a prescribed set.
	
	\subsubsection{A sharp fixed-atlas repair boundary}
	\label{app:port-boundaries}
	
	To measure the repair cost in \cref{ex:correlated-join}, suppose $B$ is a key, $A$ and $C$ each have two possible values, and the local charts are fixed up to overlap-value permutations.
	If $n_{ac}$ counts the $B$-states of type $(A=a,C=c)$, then
	\begin{equation}\label{eq:contingency-conditions}
		\begin{aligned}
			A\not\to C&\quad\Longleftrightarrow\quad
			\text{some row of }N=(n_{ac})\text{ has two positive entries},\\
			C\not\to A&\quad\Longleftrightarrow\quad
			\text{some column has two positive entries},\\
			AC\not\to B&\quad\Longleftrightarrow\quad
			\text{some entry of }N\text{ is at least two}.
		\end{aligned}
	\end{equation}
	These tests assume the only intended nontrivial FDs are $B\to A,C$,
	and the separate local relations already have no extra FDs.

	\begin{proposition}[A five-state repair of a four-state obstruction]
		\label{prop:five-sector-repair}
		\leavevmode
		\begin{enumerate}[label=\textup{(\roman*)},leftmargin=*]
			\item \emph{Four-state obstruction.}
			If the four distinct $B$-states have two states of each $A$-value and
			two of each $C$-value, then no permutation of the $B$-interface
			identification gives a faithful representation of precisely
			$B\to A,C$.
			\item \emph{Five-state repair.}
			For every original four-state alignment, one new $B$-state suffices.
			Append its rows to the two local sources. This keeps the existing
			rows and the two local closure theories. Thus
			five is the sharp cost of an append-only faithful repair for this class.
		\end{enumerate}
	\end{proposition}
	
	\begin{proof}
		\leavevmode
		\begin{enumerate}[label=\textup{(\roman*)},leftmargin=*]
			\item \emph{Four-state obstruction.}
			The possible $2\times2$ contingency matrices with both row and column margins $(2,2)$ are, up to reindexing, $\begin{psmallmatrix}k&2-k\\2-k&k\end{psmallmatrix}$,
			$k=0,1,2$.
			At $k=0,2$, $A\leftrightarrow C$ is an extra dependency;
			at $k=1$, all four cells have one state and $AC\to B$ is an extra dependency.
			Thus no alignment of those four fixed sectors is faithful.
			\item \emph{Five-state repair.}
			For $k=0$ or $2$, append one state in an originally empty diagonal or off-diagonal cell, respectively;
			for $k=1$, append a state in any cell.
			In all three cases the enlarged table has a row with two positive cells, a column
			with two positive cells, and an entry at least two, so it passes all three tests in
			\eqref{eq:contingency-conditions}.
			For the original diagonal example, one may append the cell $(A,C)=(0,1)$, obtaining $\begin{psmallmatrix}2&1\\0&2\end{psmallmatrix}$.
			The appended state has a fresh $B$-value, so $B$ remains a key.
			Each old pair witnessing the failure of $A\to B$ or $C\to B$ is still
			present. Every generating local implication is therefore preserved,
			and the four-state lower bound proves optimality for this repair class.\qedhere
		\end{enumerate}
	\end{proof}
	
	\subsection{Local pair lifting and faithful spectral descent}
	\label{sec:constructive-native}
	
	The pair test examines the full matching source. Relative pair lifting can be checked one incidence at a time.
	
	\subsubsection{Local pair lifting and actual native formation}
	\label{subsec:local-pair-native}
	
	To extend a counterexample pair across a chart, relative pair lifting prescribes its separator values and agreement on a closed set. The lifted states may agree on additional descriptors.
	
	\begin{definition}[Relative pair lifting]\label{def:relative-pair-lifting}
		For a chart $i$ and an incident join-tree edge $ij$, put $S=B_i\cap B_j$ and write $p_{i,S}:X_i\twoheadrightarrow Y_S$ for its map to the common realized separator spectrum.
		For a locally closed set $C\subseteq B_i$, set
		\[
		E_i(C)=\{(x,x')\in X_i^2:f_{i,d}(x)=f_{i,d}(x')\text{ for all }d\in C\}.
		\]
		The incidence $(i,S)$ has \textbf{relative pair lifting} if, for every locally closed $C$,
		\begin{equation}\label{eq:relative-pair-lifting}
			(p_{i,S}\times p_{i,S})E_i(C)
			=\{(y,y')\in Y_S^2:y|_{C\cap S}=y'|_{C\cap S}\}.
		\end{equation}
		An empty separator has its singleton spectrum and satisfies this condition automatically.
		The condition is required at both ends of each edge when the entire atlas is said to have relative pair lifting.
	\end{definition}
	
	The inclusion from left to right in \eqref{eq:relative-pair-lifting} is automatic.
	Its reverse says that the actual local states can lift any specified pair of separator values compatible with the declared local agreement.
	The number of local closed sets may still be exponential in $|B_i|$.
	
	\begin{theorem}[Constructive faithful native amalgamation by local pair lifting]
		\label{thm:local-pair-native}
		Let a finite spectral atlas have a join tree, equal realized
		separator images on each edge, and relative pair lifting at every
		incidence. Let $c_0$ be the least closure generated by its local
		determination closures, and let $X$ be the full natural join
		in \eqref{eq:native-joint-spectrum}. Then the following hold.
		
		\begin{enumerate}[label=(\roman*)]
			\item \emph{State coverage and dependency faithfulness.}
			Every chart map $\pi_i:X\to X_i$ is onto, and
			\begin{equation}\label{eq:local-pair-faithful}
				c_X=c_0.
			\end{equation}
			Consequently, all local states, local recovery judgments,
			and global non-recovery judgments of $c_0$ are preserved.
			
			\item \emph{Exactness across a tree cut.}
			For an edge $t=ij$, let $X_L,X_R$ be the partial joins
			on the two components of the cut tree, and let $Y_t$
			be its separator spectrum. Both maps to $Y_t$ are onto, and
			\begin{equation}\label{eq:local-pair-native-cut}
				X=X_L\times_{Y_t}X_R,
				\qquad
				M_L\cap M_R=\mathbb C[Y_t],
				\qquad
				M_L\vee M_R=\mathbb C^X.
			\end{equation}
			Here all algebras are represented by multiplication on
			$\ell^2(X)$. The algebras $M_L,M_R$ consist of functions
			of the indicated partial states, and $\mathbb C[Y_t]$
			denotes the pulled-back separator algebra.
			
			\item \emph{Native realization.}
			The construction supplies an actual access--commutant
			formation presentation. On $\mathcal H_X=\ell^2(X)$,
			let $A_d=\mathbb C[f_d]$ be the algebra of all complex-valued
			functions of the finite readout $f_d$, and define
			\begin{equation}\label{eq:local-pair-native-endpoints}
				M(T)=\operatorname{alg}^*\{A_d:d\in T\},
				\qquad
				P_T=(\mathrm{com},M(T)),
				\qquad
				K(T)=M(T)'.
			\end{equation}
			In the access--commutant context on $\mathcal H_X$,
			$P_T$ realizes the represented endpoint
			$\mathfrak P_{K(T)}$, and
			\begin{equation}\label{eq:local-pair-native-composition}
				\begin{aligned}
					K(T)'&=M(T),\\
					K(T\cup U)&=K(T)\cap K(U),\\
					\{d\in V:A_d\subseteq M(T)\}&=c_0(T).
				\end{aligned}
			\end{equation}
			The local spectra and interface maps determine this
			construction of the common native source.
		\end{enumerate}
	\end{theorem}
	
	\begin{proof}
		\leavevmode
		\begin{enumerate}[label=(\roman*)]
			\item \emph{State coverage and dependency faithfulness.}
			The one-state extension argument in \cref{prop:closure-descent} proves that $X$ is nonempty and every chart projection $\pi_i:X\to X_i$ is surjective.
			Every local implication holds after pullback to $X$, so
			$c_0(S)\subseteq c_X(S)$ for every $S\subseteq V$.
			
			To prove the reverse inclusion, fix $S_0\subseteq V$ and
			$e\notin c_0(S_0)$, and put $C=c_0(S_0)$.
			For each chart,	set $C_i=C\cap B_i$.
			Since $C$ is closed under all local	implications,
			$
			c_i(C_i)=C_i.
			$
			Root the join tree at a chart $r$ containing $e$.
			Then $e\notin c_r(C_r)$, so there are two states in $X_r$ agreeing on $C_r$ and disagreeing on $e$.
			
			Suppose a pair has been chosen at a parent chart $i$, agreeing on $C_i$. Its ordered separator values belong to the common realized separator image and agree on $C\cap S_{ij}$.
			Since $C_j$ is locally closed, relative	pair lifting at the child incidence extends these	particular separator values to a pair in $X_j$ agreeing on $C_j$.
			Apply this step successively throughout the tree.
			
			The running-intersection property ensures that the first
			coordinates of the chosen pairs form a state $x\in X$,
			and their second coordinates form a state $x'\in X$.
			These states agree on $C$ and disagree on $e$. Since
			$S_0\subseteq C$, they witness $e\notin c_X(S_0)$.
			Thus $c_X=c_0$. Surjectivity supplies every local state
			and preserves each local determination judgment, while
			the equality of closures preserves every prescribed
			global non-recovery judgment.
			
			\item \emph{Exactness across a tree cut.}
			Fix an edge $t=ij$. By running intersection, the unions
			of descriptors on the two components of the cut tree
			intersect exactly in $S_{ij}$. Consequently, two partial
			states combine to a global state precisely when their
			separator values agree. This gives the canonical
			identification
			$
			X=X_L\times_{Y_t}X_R.
			$
			
			For any $y\in Y_t$, choose a state of $X_i$ with separator
			value $y$ and extend it through the left component.
			The one-state extension argument gives an element of
			$X_L$ mapping to $y$. Applying the same argument from
			$j$ gives an element of $X_R$ mapping to $y$. Thus both
			partial-join maps onto $Y_t$ are surjective, and each
			separator fiber is a nonempty rectangle
			$
			X_y=X_{L,y}\times X_{R,y}.
			$
			
			Let $h\in M_L\cap M_R$.
			There are functions $h_L$ and $h_R$ on the respective partial joins such that
			\[
			h(x_L,x_R)=h_L(x_L)=h_R(x_R)
			\qquad ((x_L,x_R)\in X).
			\]
			Fixing a right state in $X_{R,y}$ shows that $h_L$ is
			constant on $X_{L,y}$. Hence $h$ is constant on the whole
			rectangle $X_y$ and depends only on $y$. Conversely, every
			function of the separator value belongs to both $M_L$
			and $M_R$. Therefore
			$
			M_L\cap M_R=\mathbb C[Y_t].
			$
			For each compatible pair $(x_L,x_R)$, the product of the
			indicators of its left and right states is the indicator
			of that global state. These point indicators span
			$\mathbb C^X$, proving
			$
			M_L\vee M_R=\mathbb C^X.
			$
			
			\item \emph{Native realization.}
			The joint descriptor readout on $X$ is injective. Indeed,
			equality of all descriptor values gives equality of the
			joint readout on every chart, and each local joint readout
			is injective. The multiplication algebras $A_d$ therefore
			generate $\mathbb C^X$.
			
			For each $T\subseteq V$, form $M(T)$, $P_T$, and $K(T)$
			as in \eqref{eq:local-pair-native-endpoints}, with
			$M(\varnothing)=\mathbb C1$. The access--commutant law
			\eqref{eq:standalone-commutant-law} supplies the judgment
			\[
			\Xi_{\mathcal H_X}^{\mathrm{com}}
			\vdash P_T\rightsquigarrow\mathfrak P_{K(T)}.
			\]
			Each $M(T)$ is a finite-dimensional unital $*$-subalgebra
			of $\mathcal B(\mathcal H_X)$. Bicommutant recovery gives
			$K(T)'=M(T)$, and
			\[
			\begin{aligned}
				K(T\cup U)
				&=(M(T)\vee M(U))'\\
				&=M(T)'\cap M(U)'\\
				&=K(T)\cap K(U).
			\end{aligned}
			\]
			
			The algebra $M(T)$ consists exactly of the functions
			constant on each fiber of the joint $T$-readout.
			Consequently,
			\[
			A_d\subseteq M(T)
			\quad\Longleftrightarrow\quad
			d\in c_X(T)
			\quad\Longleftrightarrow\quad
			d\in c_0(T),
			\]
			where the final equivalence follows from \textup{(i)}.
			This proves \eqref{eq:local-pair-native-composition}.
			The join, its descriptor readouts, and the represented
			algebras are determined by the given chart data and
			interface maps.\qedhere
		\end{enumerate}
	\end{proof}
	
	\begin{remark}[Pullback multiplicities]
		\label{rem:local-pair-native-scope}
		Pullback embeds each local source algebra $\mathbb C^{X_i}$ on $\ell^2(X)$.
		The multiplicity of the point $x_i$ is the size of its projection fiber in $X$.
		The commutant on this common Hilbert space uses those multiplicities.
	\end{remark}
	
	\begin{proposition}[Exact size of the full join from local fibers]
		\label{prop:local-pair-source-cost}
		Let $J$ be a finite tree on a nonempty vertex set $I$, with a finite nonempty state set $X_i$ at each vertex and surjections
		\[
		p_{i,ij}:X_i\twoheadrightarrow Y_{ij},\qquad
		p_{j,ij}:X_j\twoheadrightarrow Y_{ij}
		\]
		at each edge $ij$. Let $X$ be the set of tuples satisfying $p_{i,ij}(x_i)=p_{j,ij}(x_j)$ at every edge, and write $N_J(i)$ for the neighbors of $i$.
		For a spectral atlas with equal realized separator images, take $p_{i,ij}=p_{i,S_{ij}}$; running intersection identifies $X$ with its full natural join.
		All empty products below equal one.		
		\begin{enumerate}[label=(\roman*)]
			\item \emph{Counting by separator messages.}
			For an oriented edge $i\to j$ and $y\in Y_{ij}$, let $N_{i\to j}(y)$ be the number of compatible partial states	on the component containing $i$ after deleting $ij$, with separator value $y$.
			These numbers satisfy
			\begin{equation}\label{eq:local-pair-count-message}
				N_{i\to j}(y)=
				\sum_{\substack{x_i\in X_i\\p_{i,ij}(x_i)=y}}
				\prod_{k\in N_J(i)\setminus\{j\}}
				N_{k\to i}\bigl(p_{i,ik}(x_i)\bigr).
			\end{equation}
			For any choice of root $r$, the size of the full join is
			\begin{equation}\label{eq:local-pair-source-size}
				|X|=
				\sum_{x_r\in X_r}
				\prod_{k\in N_J(r)}
				N_{k\to r}\bigl(p_{r,rk}(x_r)\bigr).
			\end{equation}
			Thus the messages computed from the leaves toward the root determine $|X|$.
			
			\item \emph{Product formula for constant fiber sizes.}
			Suppose there is a root $r$ such that, for each nonroot vertex $i$ with parent $j$, the fiber cardinality $|p_{i,ij}^{-1}(y)|$ is independent of $y\in Y_{ij}$.
			Then
			\begin{equation}\label{eq:local-pair-uniform-cost}
				|X|=
				\frac{\prod_{i\in I}|X_i|}
				{\prod_{ij\in E(J)}|Y_{ij}|}.
			\end{equation}
		\end{enumerate}
	\end{proposition}
	
	\begin{proof}
		\leavevmode
		\begin{enumerate}[label=(\roman*)]
			\item \emph{Counting by separator messages.}
			Fix an oriented edge $i\to j$ and a state $x_i\in X_i$.
			For each neighbor $k\ne j$, the number of compatible partial states on the $k$-side matching $x_i$ is
			$
			N_{k\to i}\bigl(p_{i,ik}(x_i)\bigr).
			$
			These partial states can be chosen independently.
			In a tree, distinct components have no matching constraints between them after $x_i$ is fixed. For a spectral atlas, running intersection also ensures that any descriptor shared by two components occurs at $i$, so their values agree.
			Each collection of choices therefore gives a unique	compatible partial state on the $i$-side.
			
			Multiplying the numbers of choices and summing over the states $x_i$ with separator value $y$ proves	\eqref{eq:local-pair-count-message}.
			At a leaf, the product is empty, so the message is simply the size of the
			corresponding separator fiber.
			The same argument at the root, using all its neighboring components and summing over every $x_r\in X_r$, proves \eqref{eq:local-pair-source-size}.
			
			\item \emph{Product formula for constant fiber sizes.}
			Choose a root $r$ satisfying the hypothesis.
			For each child $i$ with parent $j$, surjectivity gives the constant fiber cardinality
			\[
			\frac{|X_i|}{|Y_{ij}|}.
			\]
			We prove inductively that each message toward the root
			is constant on its separator spectrum. For a leaf, this
			follows directly from the constant fiber hypothesis.
			Suppose that the incoming messages at $i$ have constant
			values $n_{k\to i}$.
			Then \eqref{eq:local-pair-count-message} gives
			\[
			N_{i\to j}(y)
			=
			\frac{|X_i|}{|Y_{ij}|}
			\prod_{k\in N_J(i)\setminus\{j\}}n_{k\to i},
			\]
			which is independent of $y$.
			Denote this value by $n_{i\to j}$.
			
			The root formula now becomes
			\[
			|X|=|X_r|\prod_{k\in N_J(r)}n_{k\to r}.
			\]
			Substituting the recursive expressions for the constants
			contributes one factor $|X_i|$ for each chart and one
			divisor $|Y_{ij}|$ for each tree edge.
			This proves \eqref{eq:local-pair-uniform-cost}.
		\end{enumerate}
	\end{proof}
	
	\subsubsection{Symmetric interfaces and the universal binary boundary}
	
	Symmetry sometimes supplies all the required ordered pairs.
	A separator with at most two realized joint values supplies them automatically, even when several descriptors name the separator.
	
	\begin{corollary}[Pair lifting and source size at symmetric singleton separators]
		\label{cor:symmetric-point-separators}
		Let a finite spectral atlas have a join tree and equal actual images at each separator.
		Suppose every nonempty separator is a singleton $\{s\}$.
		For each incident chart, take a permutation group on $X_i$ that	preserves the kernel of every descriptor map.
		When at least two $s$-values occur, require the induced action on those values to be two-transitive.
		\begin{enumerate}[label=\textup{(\roman*)},leftmargin=*]
			\item \emph{Pair lifting.}
			Relative pair lifting holds, and all conclusions of
			\cref{thm:local-pair-native} follow.
			\item \emph{Source size.}
			The local separator fibers have	constant size, so \eqref{eq:local-pair-uniform-cost} applies.
		\end{enumerate}
	\end{corollary}
	
	\begin{proof}
		\leavevmode
		\begin{enumerate}[label=\textup{(\roman*)},leftmargin=*]
			\item \emph{Pair lifting.}
			At an empty separator the only separator pair is the diagonal and lifts by repeating a state. At a nonempty separator, the group preserves $E_i(C)$ for every $C$.
			If $s\in C$, its image on the separator is exactly the diagonal, since every value is realized by a repeated state pair.
			If $s\notin C$ and $C$ is closed, some pair in $E_i(C)$ has distinct $s$-values.
			Two-transitivity moves this ordered pair of values to any other distinct ordered pair while preserving agreement on $C$.
			Repeated pairs again give all diagonal values.
			This proves \eqref{eq:relative-pair-lifting}.
			\item \emph{Source size.}
			For a one-point separator image the unique fiber is the whole chart. Otherwise transitivity on the separator values gives bijections between their fibers.\qedhere
		\end{enumerate}
	\end{proof}
	
	A local symmetry may act on several coordinates at once, with correlated permutations.
	In particular, finite linear spectra with scalar observables supply many examples:
	translations by local states and common nonzero scalar multiplication induce the affine two-transitive group on every nonconstant scalar readout.
	A more general linear separator criterion is developed separately below.
	
	\begin{theorem}[Automatic relative pair lifting for two-state separators]
		\label{thm:two-state-joint-separators}
		Let a finite spectral atlas have a join tree $J$ and equal realized separator images at every edge.
		Suppose that each common joint separator spectrum satisfies
		\begin{equation}\label{eq:two-state-joint-separator}
			|Y_{ij}|\leq 2
			\qquad(ij\in E(J)).
		\end{equation}
		The separator descriptor sets $S_{ij}=B_i\cap B_j$ may have arbitrary finite cardinality.
		Let $X$ be the full natural join and $c_0$ the least closure generated by the local determination closures.		
		\begin{enumerate}[label=(\roman*)]
			\item \emph{Automatic relative pair lifting.}
			Every chart--separator incidence satisfies \eqref{eq:relative-pair-lifting}.
			
			\item \emph{Faithful native amalgamation.}
			The join $X$ is nonempty, every chart projection $\pi_i:X\to X_i$ is onto, and
			$
			c_X=c_0.
			$
			The full-cut algebra identities and the access--commutant realization of \cref{thm:local-pair-native}\textup{(ii)--(iii)} also hold.
		\end{enumerate}
	\end{theorem}
	
	\begin{proof}
		\leavevmode
		\begin{enumerate}[label=(\roman*)]
			\item \emph{Automatic relative pair lifting.}
			Fix an incidence $(i,S)$ and a locally closed set $C\subseteq B_i$, so that $c_i(C)=C$.
			Write $Y$ for the realized separator spectrum and
			$p_{i,S}:X_i\twoheadrightarrow Y$ for the separator map.
			
			Every pair in $E_i(C)$ has separator values agreeing on	$C\cap S$.
			To prove the reverse inclusion in \eqref{eq:relative-pair-lifting}, fix $y,y'\in Y$ satisfying
			$
			y|_{C\cap S}=y'|_{C\cap S}.
			$
			If $y=y'$, choose $x\in X_i$ with $p_{i,S}(x)=y$.
			The pair $(x,x)$ belongs to $E_i(C)$ and lifts $(y,y')$.
			
			Suppose $y\ne y'$.
			Since $|Y|\leq2$, we have $Y=\{y,y'\}$.
			These distinct joint separator values differ on some descriptor $s\in S$.
			Their agreement on $C\cap S$ implies $s\notin C$, and local closedness gives $s\notin c_i(C)$.
			Hence there are states $x,x'\in X_i$ agreeing on every descriptor in $C$ and differing on $s$.
			Their separator images are distinct, so they are $y,y'$	in some order.
			Exchanging $x$ and $x'$ if necessary gives
			\[
			p_{i,S}(x)=y,
			\qquad
			p_{i,S}(x')=y'.
			\]
			Thus $(x,x')\in E_i(C)$ is the required ordered lift.
			This proves relative pair lifting at every incidence.
			
			\item \emph{Faithful native amalgamation.}
			Part~\textup{(i)} verifies the relative pair lifting hypothesis of \cref{thm:local-pair-native}.
			The remaining hypotheses follow from the assumed join tree and equality
			of realized separator images.
			That theorem gives nonemptiness, surjectivity of every chart projection,
			and $c_X=c_0$, together with the full-cut algebra identities and the access--commutant realization.
			\qedhere
		\end{enumerate}
	\end{proof}
	
	\begin{theorem}[Sharp separator-size threshold for universal faithful descent]
		\label{thm:binary-separator-threshold}
		Fix an integer $q\ge1$.
		The following are equivalent.
		\begin{enumerate}[label=(\roman*)]
			\item $q\le2$.
			
			\item Every finite spectral atlas with a join tree,	single-descriptor edge separators, and equal realized	separator images of cardinality at most $q$ has a 	nonempty full join $X$ satisfying
			$
			c_X=c_0.
			$
		\end{enumerate}
		For every $q\ge3$, a counterexample to \textup{(ii)} exists with two charts on $\{a,s\}$ and $\{s,b\}$, exactly $q$	realized separator values, and bijective chart projections.
		It satisfies
		\[
		c_0(\{a,b\})=\{a,b\},
		\qquad
		c_X(\{a,b\})=\{a,s,b\}.
		\]
	\end{theorem}
	
	\begin{proof}
		\leavevmode
		\begin{itemize}
			\item \emph{\textup{(i)} $\implies$ \textup{(ii)}.}
			If $q\le2$, every joint separator spectrum has at most two states.
			The conclusion follows from \cref{thm:two-state-joint-separators}.
			
			\item \emph{\textup{(ii)} $\implies$ \textup{(i)}.}
			We prove the contrapositive by constructing, for every $q\ge3$,
			an atlas that violates \textup{(ii)}.
			Let $Y=\{0,1,\ldots,q-1\}$, and choose quotient maps
			\[
			a:Y\to D_a,\qquad b:Y\to D_b
			\]
			whose only nonsingleton fibers are $\{0,1\}$ and $\{1,2\}$, respectively.
			Use the local spectra
			\[
			X_1=\{(a(y),y):y\in Y\},
			\qquad
			X_2=\{(y,b(y)):y\in Y\}
			\]
			on descriptor sets $\{a,s\}$ and $\{s,b\}$, where $s(y)=y$.
			The two charts form a join tree with common separator $\{s\}$ and realized separator image $Y$.
			
			Both $a$ and $b$ are nonconstant and noninjective.
			Consequently, the local determination closures are generated by $s\to a$ and $s\to b$, respectively.
			Neither rule adds a descriptor to $\{a,b\}$, so
			$
			c_0(\{a,b\})=\{a,b\}.
			$
			
			The full join is canonically identified with
			$
			X=\{(a(y),y,b(y)):y\in Y\}.
			$
			It has $q$ states, and each chart projection is bijective, since the separator coordinate determines the entire tuple.
			Furthermore, the joint map
			\[
			Y\longrightarrow D_a\times D_b,
			\qquad y\longmapsto(a(y),b(y)),
			\]
			is injective.
			Indeed, the only distinct pair identified by $a$ is $\{0,1\}$, and the only distinct pair identified by $b$ is $\{1,2\}$.
			Thus equality of both values forces equality of $y$.
			
			It follows that $X$ satisfies $\{a,b\}\to s$, and hence
			$
			c_X(\{a,b\})=\{a,s,b\}.
			$
			This disproves \textup{(ii)} whenever $q\ge3$ and proves the equivalence.
			
			For $q=3$, suitable value labels give the explicit local spectra
			\[
			X_1=\{(0,0),(0,1),(1,2)\},
			\qquad
			X_2=\{(0,0),(1,1),(2,1)\}.
			\qedhere
			\]
		\end{itemize}
	\end{proof}
	
	\begin{example}[Lifting without exact realization of every closed set]
		\label{ex:binary-xor-lifting}
		Take $X_i=\mathbb F_2^2$ with observables $a(u,v)=u$, $b(u,v)=v$,
		and $s(u,v)=u+v$. Its local closure has the three all-but-one
		implications, and the empty set is closed. Every two distinct states
		agree on exactly one of $a,b,s$, so no pair has exact agreement set
		$\varnothing$. Thus the hypothesis of
		\cref{cor:relative-pair-homogeneity} fails. Nevertheless every
		singleton interface has two values, so relative pair lifting holds
		at all such interfaces. Arbitrarily large join-tree atlases made from
		these charts therefore admit the constructive faithful native
		amalgamation of \cref{thm:local-pair-native}. This example separates
		the present sufficient condition from the earlier demand for exact
		realization of every local closed set.
	\end{example}
	
	\begin{remark}[Attachment obstructions from one-sided lifting failure]
		\label{rem:point-lifting-boundary}
		If a closed set $C$ omitting $s$ admits no lift of two distinct interface values $a,b$, the elementary partner in \cref{thm:universal-singleton-port} creates the extra dependency $C\cup\{z\}\to s$ by identifying $a,b$ at $z$.
		The pair-lifting condition on one chart still allows a failure
		when a partner lacks that condition. For example, on descriptors $(s,a,e)$ let
		\[
		R=\{(0,0,0),(1,0,0),(2,0,1),
		(0,1,0),(1,1,0),(2,1,0)\}.
		\]
		Its nontrivial dependencies are generated by $as\to e$. The closed sets omitting
		$s$ are $\varnothing,\allowbreak\{a\},\allowbreak\{e\},\allowbreak\{a,e\}$, and the three states with
		$(a,e)=(1,0)$ lift all pairs of separator values for each of them.
		Hence $R$ has relative pair lifting at $s$. Attach
		$Q=\{(0,0),(1,0),(2,1)\}$ on $(s,z)$. The join satisfies $az\to e$,
		whereas the local rules $as\to e$ and $s\to z$ do not imply it.
		The companion $Q$ fails relative pair lifting.
	\end{remark}
	
	\subsubsection{Faithfulness under all attachments at a singleton interface}
	\label{subsec:universal-singleton-port}
	
	For a fixed chart, ask whether it attaches faithfully to every partner with the same single-descriptor interface.
	The target of a failed dependency can lie on either side.
	The test must therefore	lift agreement on the premise and, in certain cases, disagreement on the target.
	
	Fix a finite spectral chart $(X,B,(f_e)_{e\in B})$ as in	\cref{def:spectral-atlas}, a descriptor $s\in B$, and its actual interface
	alphabet $Y=f_s(X)$.
	Write $c$ for the determination closure of $X$ and
	$Y^2_{\ne}=\{(a,b)\in Y^2:a\ne b\}$.
	
	\begin{definition}[Universal singleton-interface condition]
		\label{def:universal-singleton-port}
		The chart satisfies the \textbf{universal singleton-interface condition}
		at $s$ if, for every $c$-closed $C\subseteq B$ and every
		$e\in c(C\cup\{s\})\setminus C$,
		\begin{equation}\label{eq:universal-port-pairs}
			\left\{(f_s(x),f_s(x')):\begin{array}{l}
				x,x'\in X,\quad f_u(x)=f_u(x')\quad(u\in C),\\
				f_e(x)\ne f_e(x')
			\end{array}\right\}=Y^2_{\ne}.
		\end{equation}
		The inclusion from left to right is automatic, because
		$e\in c(C\cup\{s\})$.
		The requirement is that every ordered pair of distinct interface values occurs in such a witness.
	\end{definition}
	
	Whenever $s\notin C$, choosing $e=s$ in	\eqref{eq:universal-port-pairs} gives full lifting of interface pairs by states agreeing on $C$.
	The additional targets record exactly the conditional distinctions that the interface can supply beyond $C$.
	There is no assumption that every closed set is the exact agreement
	set of a pair.
	
	\begin{theorem}[Local criterion for universal faithful attachment at a singleton interface]
		\label{thm:universal-singleton-port}
		Fix the finite spectral chart $(X,B,(f_e)_{e\in B})$, a descriptor
		$s\in B$, and its realized interface image $Y=f_s(X)$.
		A partner is a finite spectral chart $(Z,B',(g_e)_{e\in B'})$
		with $B\cap B'=\{s\}$ and $g_s(Z)=Y$.
		For each partner, form the join
		\begin{equation}\label{eq:universal-port-join}
			R=X\times_Y Z
			=\{(x,\zeta)\in X\times Z:f_s(x)=g_s(\zeta)\},
		\end{equation}
		and let $c_0$ be the least closure generated by the two local
		determination closures. The projections $R\to X$ and $R\to Z$
		are surjective, and the following conditions are equivalent.	
		\begin{enumerate}[label=(\roman*)]
			\item \emph{Local interface condition.}
			The universal singleton-interface condition holds at $s$.
			
			\item \emph{Faithfulness for all partners.}
			For every partner, the join satisfies $c_R=c_0$.
			
			\item \emph{Faithfulness for key partners.}
			The equality in \textup{(ii)} holds for every partner on
			$B'=\{s,z\}$ whose $s$-observable is a key.
			
			\item \emph{Faithfulness for elementary partners.}
			The equality in \textup{(ii)} holds for each of the $\binom{|Y|}{2}$ partners with states
			$
			Z_h=\{(y,h(y)):y\in Y\}
			$
			on $\{s,z\}$, where $z\notin B$ and $h$ identifies one prescribed unordered pair of distinct values of $Y$ and no other distinct values.
		\end{enumerate}	
		If \textup{(i)} fails, one of the elementary partners in
		\textup{(iv)} gives a counterexample to faithful attachment
		with exactly $|Y|$ partner states.
	\end{theorem}
	
	\begin{proof}
		Write $c$ for the closure of the fixed chart and $d$ for the
		closure of a partner. Equality of the realized interface images
		makes both projections $R\to X,Z$ surjective. Every local
		implication holds after pullback to $R$, so $c_0\leq c_R$.	
		\begin{itemize}
			\item \emph{\textup{(i)} $\implies$ \textup{(ii)}.}
			Assume the universal singleton-interface condition.
			Let $M\subseteq B\cup B'$ be $c_0$-closed, put
			\[
			C=M\cap B,\qquad D=M\cap B',
			\]
			and fix $e\notin M$. Then $c(C)=C$ and $d(D)=D$.
			We construct two states of $R$ agreeing on $M$ and
			disagreeing on $e$.
			
			First suppose $e\in B'$. Since $e\notin d(D)$, there
			are $\zeta,\zeta'\in Z$ agreeing on $D$ and disagreeing
			on $e$. Write
			\[
			a=g_s(\zeta),\qquad b=g_s(\zeta').
			\]
			If $a=b$, choose any $x\in X$ with $f_s(x)=a$.
			Then $(x,\zeta)$ and $(x,\zeta')$ are the required states.
			If $a\ne b$, agreement on $D$ implies $s\notin M$ and
			hence $s\notin C$. Apply
			\eqref{eq:universal-port-pairs} to $C$ with target $s$.
			It gives $x,x'\in X$ agreeing on $C$ with interface
			values $a,b$, respectively. The states
			$(x,\zeta)$ and $(x',\zeta')$ give the required pair.
			
			Now suppose $e\in B\setminus B'$. If
			$e\notin c(C\cup\{s\})$, choose $x,x'\in X$ agreeing
			on $C\cup\{s\}$ and disagreeing on $e$. A partner state
			with their common interface value extends them to two
			states of $R$ agreeing on $M$ and disagreeing on $e$.
			
			Otherwise,
			$e\in c(C\cup\{s\})\setminus C$. Closedness of $C$
			implies $s\notin C$, so $s\notin M$ and $s\notin D$.
			Since $d(D)=D$, choose $\zeta,\zeta'\in Z$ agreeing
			on $D$ and having distinct interface values $a,b$.
			Condition~\eqref{eq:universal-port-pairs} supplies
			$x,x'\in X$ agreeing on $C$, disagreeing on $e$, and
			having interface values $a,b$, respectively. Again,
			$(x,\zeta)$ and $(x',\zeta')$ give the required pair.
			
			For any premise $S$, apply this construction to
			$M=c_0(S)$. Every $e\notin c_0(S)$ then has a
			non-recovery witness in $R$, proving $c_R\leq c_0$.
			Together with the reverse inclusion, this gives $c_R=c_0$.
			
			\item \emph{\textup{(ii)} $\implies$ \textup{(iii)}.}
			Every two-descriptor key partner in \textup{(iii)} is
			a partner covered by \textup{(ii)}.
			
			\item \emph{\textup{(iii)} $\implies$ \textup{(iv)}.}
			Each elementary partner has exactly one state for each
			value $y\in Y$. Its $s$-observable is therefore a key,
			so \textup{(iii)} applies.
			
			\item \emph{\textup{(iv)} $\implies$ \textup{(i)}.}
			We prove the contrapositive. Suppose \textup{(i)} fails.
			Then there are a $c$-closed set $C$, a target
			$
			e\in c(C\cup\{s\})\setminus C,
			$
			and an ordered pair $(a,b)\in Y^2_{\ne}$ for which
			the required witness in \eqref{eq:universal-port-pairs}
			does not exist. Necessarily $s\notin C$.
			
			Let $h$ identify $a$ and $b$ and no other distinct
			values of $Y$. Choose a descriptor $z\notin B$ and form
			the elementary partner
			\begin{equation}\label{eq:universal-port-test-partner}
				Z=\{(y,h(y)):y\in Y\},
				\qquad B'=\{s,z\}.
			\end{equation}
			Its $s$-observable is a key. Since $h(a)=h(b)$ with
			$a\ne b$, the descriptor $z$ does not determine $s$,
			and hence $d(\{z\})=\{z\}$. Thus
			$
			M=C\cup\{z\}
			$
			has locally closed intersections with both charts.
			It is therefore $c_0$-closed and omits $e$.
			
			Nevertheless, $M$ determines $e$ in $R$. To see this,
			take two global states agreeing on $C$ and on $z$.
			If their $s$-values agree, their $e$-values agree because
			$e\in c(C\cup\{s\})$. If their $s$-values differ, those
			values must be $a,b$ in some order. The missing witness
			for $(a,b)$ forces agreement on $e$, and reversing
			a pair gives the same conclusion for $(b,a)$.
			Consequently,
			$
			e\in c_R(M)\setminus c_0(M),
			$
			contradicting \textup{(iv)}. The partner has exactly
			$|Y|$ states, which also proves the counterexample
			assertion.\qedhere
		\end{itemize}
	\end{proof}
	
	\begin{corollary}[Binary and symmetry criteria for universal faithful attachment]
		\label{cor:universal-port-symmetry}
		The fixed chart satisfies the universal singleton-interface
		condition at $s$ under either of the following hypotheses.
		Consequently, under either hypothesis it attaches faithfully
		to every finite spectral partner with the same realized
		singleton-interface image.	
		\begin{enumerate}[label=(\roman*)]
			\item \emph{At most two interface values.}
			The realized interface image satisfies $|Y|\leq2$.
			
			\item \emph{Double transitivity.}
			A group of permutations of $X$ preserves the equality
			relation of every observable $f_e$ and induces a
			doubly transitive action on $Y$.
		\end{enumerate}
	\end{corollary}
	
	\begin{proof}
		By \cref{thm:universal-singleton-port}, it suffices to verify
		the universal singleton-interface condition under each hypothesis.	
		\begin{enumerate}[label=(\roman*)]
			\item \emph{At most two interface values.}
			If $|Y|=1$, then $s\in c(\varnothing)$, so every closed
			set $C$ contains $s$. Hence
			$c(C\cup\{s\})\setminus C$ is empty, and the condition
			holds.
			
			Suppose $|Y|=2$. Fix a closed set $C$ and a target
			$e\in c(C\cup\{s\})\setminus C$. Since $e\notin c(C)$,
			there are $x,x'\in X$ agreeing on $C$ and disagreeing
			on $e$. Their interface values must be distinct, because
			$C\cup\{s\}$ determines $e$. This pair and its reversal
			realize the two ordered pairs in $Y^2_{\ne}$, proving
			\eqref{eq:universal-port-pairs}.
			
			\item \emph{Double transitivity.}
			Fix a closed set $C$ and a target
			$e\in c(C\cup\{s\})\setminus C$. Choose $x,x'\in X$
			agreeing on $C$ and disagreeing on $e$. Their interface
			values are distinct.
			
			Each group element preserves observable equality and
			inequality, since both it and its inverse preserve the
			observable equality relations. In particular, it maps
			$f_s$-fibers to $f_s$-fibers and induces a permutation
			of $Y$.
			
			Given any $(a,b)\in Y^2_{\ne}$, double transitivity
			provides a group element sending the ordered interface
			values of $x,x'$ to $a,b$. The transformed states still
			agree on $C$ and disagree on $e$. Thus every required
			ordered interface pair has a witness, proving
			\eqref{eq:universal-port-pairs}.\qedhere
		\end{enumerate}
	\end{proof}
	
	\begin{corollary}[Faithful descent with an oriented interface condition]
		\label{cor:oriented-universal-ports}
		Let a finite spectral atlas have a join tree rooted at $r$,
		with equal realized separator images and every nonempty edge
		separator a singleton. For each nonroot chart $i$, let $p(i)$
		be its parent and put
		$
		S_i=B_i\cap B_{p(i)}.
		$
		Whenever $S_i\ne\varnothing$, assume that chart $i$ satisfies
		the universal singleton-interface condition at $S_i$.
		Let $X$ be the full join and $c_0$ the least closure generated
		by the local determination closures. For an empty separator,
		take $Y_{i,p(i)}$ to be a one-point set.	
		\begin{enumerate}[label=(\roman*)]
			\item \emph{Faithful descent.}
			The join $X$ is nonempty, every chart projection
			$\pi_i:X\to X_i$ is onto, and $c_X=c_0$.
			
			\item \emph{Full-source cardinality.}
			If every child-to-parent separator map has constant
			fiber cardinality, then
			\begin{equation}\label{eq:oriented-port-source-size}
				|X|=|X_r|
				\prod_{i\ne r}
				\frac{|X_i|}{|Y_{i,p(i)}|}.
			\end{equation}
			
			\item \emph{Symmetry certificate.}
			The constant-fiber hypothesis in \textup{(ii)} holds
			if every nonempty child interface satisfies
			\cref{cor:universal-port-symmetry}\textup{(ii)}.
			Hence these symmetry certificates imply
			\eqref{eq:oriented-port-source-size}.
		\end{enumerate}
	\end{corollary}
	
	\begin{proof}
		\leavevmode
		\begin{enumerate}[label=(\roman*)]
			\item \emph{Faithful descent.}
			Add charts in a parent-before-child order, starting with
			$X_r$. Inductively, let $H$ be the partial join on the
			charts already added, and let $B_H$ be the union of their
			descriptor sets. Assume that $H$ projects onto every
			chart already added and that its closure is generated
			by their local closure theories.
			
			The joint descriptor readout on $H$ is injective, because
			the local chart readouts are injective. Thus $H$ is a
			finite spectral chart on $B_H$. For the next chart $i$,
			running intersection gives
			$
			B_i\cap B_H=S_i.
			$
			Moreover, the image of $H$ on $S_i$ is the prescribed
			separator image, since $H$ projects onto the parent chart.
			
			If $S_i$ is a singleton, apply
			\cref{thm:universal-singleton-port} with $X_i$ as the
			chart satisfying the universal condition and $H$ as
			its partner. The new partial join projects onto both
			sources, and its closure is generated by $c_H$ and $c_i$.
			By the induction hypothesis, this is precisely the
			closure generated by all chart theories added so far.
			
			If $S_i=\varnothing$, the new partial join is
			$H\times X_i$ on disjoint descriptor sets. For every
			$T\subseteq B_H\cup B_i$,
			\[
			c_{H\times X_i}(T)
			=
			c_H(T\cap B_H)\cup c_i(T\cap B_i).
			\]
			Indeed, local implications hold in the product, and
			a failed local determination extends to the product
			by holding the other coordinate fixed. This equality
			is the closure generated by the two local theories.
			The product also projects onto both nonempty factors.
			
			The induction therefore preserves state coverage and
			the asserted closure equality at every step. At the
			final step it gives \textup{(i)}.
			
			\item \emph{Full-source cardinality.}
			Apply \cref{prop:local-pair-source-cost}\textup{(ii)} at the specified root. Each edge occurs once as a child-to-parent edge, giving \eqref{eq:oriented-port-source-size}. An empty separator contributes a divisor of one.
			
			\item \emph{Symmetry certificate.}
			At a nonempty child interface, the symmetry certificate
			induces a transitive action on its realized image.
			Given two interface values $y,y'$, choose a group element
			sending $y$ to $y'$. Its action on $X_i$ maps the fiber
			over $y$ bijectively onto the fiber over $y'$. Thus all
			fibers have the same cardinality. A one-point interface
			has a single fiber, and the same conclusion holds for
			an empty separator under the stated convention.
			Part~\textup{(ii)} now gives the formula.\qedhere
		\end{enumerate}
	\end{proof}
	
	\subsubsection{Pair lifting at homomorphic interfaces}
	\label{subsec:homomorphic-ports}
	
	For homomorphic readouts, a kernel test on one chart decides faithful attachment to every finite spectral partner. It uses the classical subgroup description of functional dependencies \cite[Proposition~3]{Chan2007} \cite[Section~5.3]{GogaczTorunczyk2017}.
	
	Let $G$ be a finite group, let $B$ be a finite descriptor set, and let
	the jointly separating readouts $\varphi_d:G\to H_d$, $d\in B$, be
	group homomorphisms.
	Their joint image is the prescribed spectral chart, identified with $G$.
	Fix $s\in B$, put $Y=\varphi_s(G)$, and write $p:G\twoheadrightarrow Y$ for the actual interface map.
	For $A\subseteq B$, define
	\begin{equation}\label{eq:homomorphic-local-kernels}
		K_A=\bigcap_{d\in A}\ker\varphi_d,
		\qquad K_\varnothing=G,
		\qquad
		c_G(A)=\{e\in B:K_A\leq\ker\varphi_e\}.
	\end{equation}
	The last equality follows because $x,x'$ agree on $A$ exactly when
	$x^{-1}x'\in K_A$.
	In particular, $K_{c_G(A)}=K_A$, and $p(K_A)$ is a normal subgroup of the actual interface group $Y$.
	
	\begin{theorem}[Kernel criterion for faithful attachment at a homomorphic singleton interface]
		\label{thm:homomorphic-singleton-ports}
		For the fixed chart and interface above, the following are equivalent.
		\begin{enumerate}[label=(\roman*),nosep]
			\item \emph{Kernel images.} For every $A\subseteq B$,
			\begin{equation}\label{eq:homomorphic-port-dichotomy}
				p(K_A)\in\{\{1\},Y\}.
			\end{equation}
			It is equivalent to test only locally closed $A$.
			
			\item \emph{Relative pair lifting.} Relative pair lifting holds at the singleton interface $\{s\}$.
			
			\item \emph{Universal attachment.} The universal singleton-interface condition
			\eqref{eq:universal-port-pairs} holds at $s$.
			Equivalently, the full join with every compatible finite spectral partner sharing only $s$ is state-exact and has exactly the generated recovery closure.
			
			\item \emph{Quotient companions.} The full join is faithful for every normal-quotient companion
			\begin{equation}\label{eq:normal-quotient-port-companion}
				Z_L=\{(y,yL):y\in Y\}
				\quad\text{on descriptors }\{s,z\},
				\qquad \{1\}<L<Y,\quad L\trianglelefteq Y,
			\end{equation}
			where $z\notin B$.
		\end{enumerate}
		Consequently these conditions are also equivalent to faithful attachment
		to every compatible group-homomorphic partner.
		Every test in \textup{(i)} is internal to the given chart.
		If \textup{(i)} fails at $A$, the subgroup $L=p(K_A)$ supplies a companion $Z_L$ for which faithfulness fails. This companion has exactly $|Y|$ states and one new descriptor.
	\end{theorem}
	
	\begin{proof}
		\leavevmode
		\begin{itemize}
			\item \emph{\textup{(i)} $\Longleftrightarrow$ \textup{(ii)}.}
			Fix a locally closed set $C\subseteq B$.
			An ordered interface pair $(a,b)\in Y^2$ lifts to two states agreeing
			on $C$ if and only if
			$
			a^{-1}b\in p(K_C).
			$
			Indeed, any such lift $x,x'$ satisfies $x^{-1}x'\in K_C$ and
			$p(x^{-1}x')=a^{-1}b$.
			Conversely, choose $x\in G$	with $p(x)=a$ and $k\in K_C$ with
			$p(k)=a^{-1}b$.
			Then $x,xk$ agree on $C$ and have interface values $a,b$, respectively.
			
			If $s\in C$, then $p(K_C)=\{1\}$, and every required
			diagonal interface pair lifts by repeating a state.
			If $s\notin C$, local closedness gives $K_C\not\leq\ker p$, so $p(K_C)\ne\{1\}$.
			In this case relative pair lifting requires every ordered interface pair to lift, which is equivalent to $p(K_C)=Y$.
			
			Thus relative pair lifting is equivalent to the dichotomy in \textup{(i)} for locally closed sets.
			Since $K_{c_G(A)}=K_A$ for every $A\subseteq B$, these tests are also equivalent to the tests for all subsets $A$.
			
			\item \emph{\textup{(i)} $\implies$ \textup{(iii)}.}
			Assume \textup{(i)}.
			Fix a locally closed $C$ and
			$
			e\in c_G(C\cup\{s\})\setminus C.
			$
			Necessarily $s\notin C$.
			Put
			\[
			K=K_C,\qquad N=\ker p,\qquad
			K'=K\cap\ker\varphi_e=K_{C\cup\{e\}}.
			\]
			The closure conditions and \textup{(i)} give
			\[
			p(K)=Y,\qquad
			K\cap N\leq\ker\varphi_e,\qquad
			K\not\leq\ker\varphi_e.
			\]
			
			Condition~\textup{(i)} also applies to $K'$.
			Suppose	$p(K')=Y$.
			For every $k\in K$, choose $k'\in K'$ with $p(k')=p(k)$.
			Then
			$
			(k')^{-1}k\in K\cap N\leq\ker\varphi_e.
			$
			Since $k'\in\ker\varphi_e$, it follows that	$k\in\ker\varphi_e$.
			This contradicts $K\not\leq\ker\varphi_e$.
			Therefore $p(K')=\{1\}$.
			
			Now fix any distinct $a,b\in Y$.
			Choose $x\in G$ with $p(x)=a$ and $k\in K$ with $p(k)=a^{-1}b$.
			Since $p(k)\ne1$, we have $k\notin K'$, and hence $\varphi_e(k)\ne1$. Thus $x,xk$ agree on $C$, have interface values $a,b$, and disagree on $e$.
			
			This realizes every ordered pair in $Y^2_{\ne}$.
			Conversely, a pair agreeing on $C$ and disagreeing on $e$ must have distinct interface values, because $e\in c_G(C\cup\{s\})$.
			This proves	\eqref{eq:universal-port-pairs}.
			Its attachment interpretation follows from 			\cref{thm:universal-singleton-port}.
			
			\item \emph{\textup{(iii)} $\implies$ \textup{(iv)}.}
			Each $Z_L$ is a finite spectral partner with descriptor	overlap $\{s\}$ and realized interface image $Y$.
			Consequently, the universal attachment assertion in	\textup{(iii)} applies to every companion in \textup{(iv)}.
			
			\item \emph{\textup{(iv)} $\implies$ \textup{(i)}.}
			We prove the contrapositive.
			Suppose \textup{(i)} fails for some $A\subseteq B$.
			Then
			$
			L=p(K_A)
			$
			is a proper nontrivial normal subgroup of $Y$.
			Normality follows because $K_A\trianglelefteq G$ and $p:G\twoheadrightarrow Y$ is surjective.
			Put $C=c_G(A)$. Then $K_C=K_A$, and $s\notin C$	because $p(K_C)\ne\{1\}$.
			
			Consider the companion $Z_L$ in \eqref{eq:normal-quotient-port-companion}.
			Its observations are the identity map on $Y$ and the
			quotient map $Y\to Y/L$. Both are nonconstant, and
			the quotient map is noninjective. Its local closure
			is therefore generated by $s\to z$.
			
			View $C$ as a subset of $B\cup\{z\}$. Its intersection
			with the original chart is locally closed, and its
			intersection with $\{s,z\}$ is empty. The empty set
			is closed in the companion because both observations
			are nonconstant. Hence $C$ is closed under the generated
			closure $c_0$, and $z\notin c_0(C)$.
			
			Nevertheless, the join determines $z$ from $C$.
			Indeed, if two global states agree on $C$, their $G$-coordinates $x,x'$ satisfy
			\[
			x^{-1}x'\in K_C,
			\qquad
			p(x)^{-1}p(x')\in p(K_C)=L.
			\]
			Their interface values therefore have the same $L$-coset, so their $z$-values agree.
			Thus
			$
			z\in c_R(C)\setminus c_0(C),
			$
			which disproves \textup{(iv)}.
			
			The companion $Z_L$ has exactly $|Y|$ states and one
			new descriptor. It is group-homomorphic, since its
			readouts are the identity and quotient homomorphisms.
			Condition~\textup{(iii)} implies faithful attachment
			to every compatible group-homomorphic partner.
			Conversely, faithful attachment to every such partner
			implies \textup{(iv)}, because each $Z_L$ belongs to
			that class. This proves the additional equivalence
			and the counterexample assertion.\qedhere
		\end{itemize}
	\end{proof}
	
	For a fixed homomorphic chart, attachment is controlled by the images $p(K_A)$. Requiring the test for every chart gives a condition on the interface group $Y$ itself.
	
	\begin{corollary}[Faithful attachment at simple interfaces]
		\label{cor:simple-algebraic-ports}
		Let $Y$ be a fixed nontrivial finite group.
		The following are equivalent.
		\begin{enumerate}[label=(\roman*),nosep]
			\item \emph{Simplicity.} $Y$ is simple.
			
			\item \emph{All spectral partners.} Every finite group-homomorphic chart whose actual singleton interface group is $Y$ attaches faithfully to every compatible finite spectral partner across that interface.
			
			\item \emph{Homomorphic partners.} Every pair of compatible group-homomorphic charts with that
			actual singleton interface group has a faithful full join.
		\end{enumerate}
		The same equivalence holds for a fixed nonzero finite module $Y$ over a fixed finite ring, with ``simple group'' replaced by ``simple	module'' and all chart readouts taken to be module homomorphisms.
		In particular, every one-dimensional $\mathbb F_q$-valued linear interface has the universal attachment property, for every prime power $q$.
	\end{corollary}
	
	\begin{proof}
		\leavevmode
		\begin{itemize}
			\item \emph{\textup{(i)} $\implies$ \textup{(ii)}.}
			For a group chart, each $p(K_A)$ is normal in $Y$.
			Simplicity therefore forces \eqref{eq:homomorphic-port-dichotomy}, proving	\textup{(i)}$\Rightarrow$\textup{(ii)} by		\cref{thm:homomorphic-singleton-ports}.
			\item \emph{\textup{(ii)} $\implies$ \textup{(iii)}.}
			Compatible group-homomorphic partners are finite spectral partners.
			\item \emph{\textup{(iii)} $\implies$ \textup{(i)}.}
			If $Y$ is not simple, choose $\{1\}<N<Y$ with $N\trianglelefteq Y$.
			Take the left and right local source groups both to be $Y$, with
			readouts
			\[
			\text{left: }s(y)=y,\quad e(y)=yN;
			\qquad
			\text{right: }s(y)=y,\quad a(y)=yN.
			\]
			The local theories say $s\to e$ and $s\to a$.
			Each quotient value	comes from several values of $s$.
			In the full join the shared identity readout forces	the two source elements to coincide, so $e=a$ and the unintended dependencies $e\to a$ and $a\to e$ both hold.
			Thus \textup{(iii)} fails.
			\item \emph{Module specialization.}
			For modules the proof is identical, with kernels and their images
			submodules and with $0$ in place of $1$.
			Equivalently, apply the	kernel-test proof to the underlying additive groups, using the fact that all the images $p(K_A)$ are submodules.
			A nonzero proper submodule gives the same two quotient charts when simplicity fails.
			A one-dimensional vector space has no nonzero proper subspace.\qedhere
		\end{itemize}
	\end{proof}
	
	A fixed chart can satisfy \eqref{eq:homomorphic-port-dichotomy} even when its interface group is nonsimple. For example, in the single-readout chart $G=Y$ with $s=\operatorname{id}_Y$, the only possible kernel images are $Y$ and $\{1\}$.
	
	\begin{corollary}[Faithful native trees of homomorphic interfaces]
		\label{cor:homomorphic-port-trees}
		Consider a finite rooted join-tree atlas of finite groups with
		group-homomorphic readouts, equal actual separator images, and separators
		of size at most one.
		At each nonroot chart with a nonempty parent separator, require \eqref{eq:homomorphic-port-dichotomy} only for its
		interface to its parent.
		\begin{enumerate}[label=\textup{(\roman*)},leftmargin=*]
			\item \emph{State coverage and recovery.}
			The full join $G$ is a group under coordinatewise operations,
			all chart projections are surjective homomorphisms, and its recovery
			closure is exactly $c_0$.
			\item \emph{Source order.}
			Its order is
			\begin{equation}\label{eq:homomorphic-tree-source-order}
				|G|=|G_r|\prod_{i\ne r}
				\frac{|G_i|}{|Y_{i,p(i)}|},
			\end{equation}
			where an empty separator has the one-element group.
			\item \emph{Native realization.}
			In the fixed equational group context of		\cref{def:equational-native-context}, this construction supplies a
			native formation retaining the group operations and the marked local
			and separator quotient maps.
		\end{enumerate}
		All three conclusions also hold for finite modules over a fixed finite ring.
		Simple parent-interface groups or simple parent-interface modules automatically satisfy the required local tests.
	\end{corollary}
	
	\begin{proof}
		\leavevmode
		\begin{enumerate}[label=\textup{(\roman*)},leftmargin=*]
			\item \emph{State coverage and recovery.}
			By \cref{thm:homomorphic-singleton-ports}, every child interface has the
			universal singleton-interface property.
			Apply \cref{cor:oriented-universal-ports} to obtain state exactness and
			$c_G=c_0$.
			Equality of homomorphic	readouts is preserved by coordinatewise multiplication and inversion, so the full join is a group and its projections are homomorphisms.
			\item \emph{Source order.}
			Each homomorphic interface has constant fiber size $|G_i|/|Y_{i,p(i)}|$, giving \eqref{eq:homomorphic-tree-source-order}.
			\item \emph{Native realization.}
			The native formation and its exact separator congruences are supplied
			by the operation-preserving construction of
			\cref{thm:equational-native-amalgamation}; faithfulness here is
			provided by the oriented interface argument.
		\end{enumerate}
		The same reasoning applies to addition and scalar multiplication for modules. Simplicity forces the kernel-image dichotomy at the corresponding parent interfaces.
	\end{proof}
	
	The order formula counts the full matching source.
	Requiring a native group law can increase the minimum for the same observations (\cref{prop:native-group-versus-spectral-cost}).
	Terminal-tree indicator maps and their complements take two values, so the binary-interface theorem also applies to those charts.

	\subsubsection{Linear sources with visible separator information}
	\label{subsec:linear-native-descent}

	A linear separator can transmit a combination absent from the named descriptors. The following condition detects, directly from the local maps, when this limits the pairs that can be lifted.

	Fix a finite field $\mathbb F_q$ and a join tree $J$ on the charts $B_i$.
	Let $X_i$ be a finite-dimensional $\mathbb F_q$-vector space and let
	$f_{i,e}:X_i\to\mathbb F_q$, $e\in B_i$, be linear observables whose joint
	map is injective.
	Put $W_i=X_i^*$ and write $v_{i,e}=f_{i,e}\in W_i$;
	these forms span $W_i$.
	Zero forms and parallel descriptor labels are allowed.
	On every edge $a=ij$, assume equality of the actual separator
	images
	\[
	Y_a=f_{i,S_a}(X_i)=f_{j,S_a}(X_j),\qquad S_a=B_i\cap B_j.
	\]
	The dual space $U_a=Y_a^*$ embeds in both $W_i$ and $W_j$ by pullback,
	identifying the corresponding separator forms.
	Write
	\[
	r_i=\dim W_i,\qquad s_a=\dim U_a.
	\]
	As before, $c_i$ is the actual local determination closure, and $c_0$ is
	the least closure generated by all local implications. Here
	\[
	c_i(A)=\{e\in B_i:v_{i,e}\in
	\operatorname{span}\{v_{i,a}:a\in A\}\}.
	\]

	\begin{definition}[Visible linear separators]
		\label{def:linear-separator-visibility}
		An incident separator $U_a\subseteq W_i$ is \textbf{visible} if, for every
		$A\subseteq B_i$,
		\begin{equation}\label{eq:linear-visibility}
			\operatorname{span}(v_{i,A})\cap U_a
			=\operatorname{span}\{v_{i,s}:s\in c_i(A)\cap S_a\}.
		\end{equation}
		Equivalently, the equality holds with $A$ replaced by every local closed
		set $C$.
		Thus every linear separator functional determined by $A$ is generated by named separator descriptors already recoverable from $A$.
	\end{definition}

	The condition is local to one chart and one incident separator.
	It can be checked by finite-dimensional linear algebra without forming a global
	relation or enumerating global nonimplications.
	It is automatic when $s_a\leq1$.
	More generally, it holds if the separator descriptor forms contain a representative of every one-dimensional subspace of $U_a$.
	Indeed, each vector in the left side of \eqref{eq:linear-visibility} is
	a multiple of a separator form lying in the span of $A$, and that descriptor belongs to $c_i(A)$.

	If $S_a$ is a local flat, \eqref{eq:linear-visibility} is exactly the modular rank identity against every local flat.
	This is the	modular amalgamation mechanism for representable matroids
	\cite{Brylawski1975,BakerLorscheidWalshZhang2025}.

	\begin{lemma}[Linear visibility and relative pair lifting]
		\label{lem:linear-visible-pairs}
		\leavevmode
		\begin{enumerate}[label=\textup{(\roman*)},leftmargin=*]
			\item \emph{Liftable separator pairs.}
			For a local closed set $C\subseteq B_i$, the pairs of separator states
			that lift to a pair of local states agreeing on $C$ are exactly the
			pairs whose difference annihilates $\operatorname{span}(v_{i,C})\cap U_a$.
			
			\item \emph{Visibility criterion.}
			Consequently \eqref{eq:linear-visibility} is equivalent to the following condition:
			every ordered pair in $Y_a^2$ agreeing on $C\cap S_a$ lifts to an
			ordered pair in $X_i^2$ agreeing on $C$.
		\end{enumerate}
	\end{lemma}

	\begin{proof}
		\leavevmode
		\begin{enumerate}[label=\textup{(\roman*)},leftmargin=*]
			\item \emph{Liftable separator pairs.}
			Put $H=\operatorname{span}(v_{i,C})$.
			The differences of pairs agreeing on $C$ form the annihilator $H^\perp\subseteq X_i$.
			Restriction to $U_a$ maps this annihilator onto $(H\cap U_a)^\perp\subseteq U_a^*$:
			a functional on $U_a$ vanishing on $H\cap U_a$ extends first to
			$U_a+H$, by setting it zero on $H$, and then to $W_i$.
			Given a separator pair $(y,y')$, lift $y$ to any $x\in X_i$ and lift
			$y'-y$ by such a difference $z$; then $(x,x+z)$ is the required pair.
			
			\item \emph{Visibility criterion.}
			Equality of the annihilators is equivalent to equality of the
			subspaces in \eqref{eq:linear-visibility}.\qedhere
		\end{enumerate}
	\end{proof}

	\begin{theorem}[Faithful native amalgamation for linear join-tree atlases]
		\label{thm:linear-native-amalgamation}
		Assume that every incidence of a join-tree edge satisfies
		\eqref{eq:linear-visibility}.
		Then the actual full join $X$ has the following properties.
		\begin{enumerate}[label=(\roman*)]
			\item \emph{Linear source and dimension.}
			It is a vector space, every projection $X\to X_i$ is surjective, and
			\begin{equation}\label{eq:linear-source-size}
				R=\dim X=\sum_i r_i-\sum_{a\in E(J)}s_a,
				\qquad |X|=q^R.
			\end{equation}
			It is constructed by solving the specified local matching equations;
			the solutions are the global states.
			
			\item \emph{Global determination.}
			The global observables have determination closure exactly $c_0$.
			In particular, $c_0$ is a representable matroid closure, of rank $R$.
			All prescribed dependencies and nondependencies are preserved.
			
			\item \emph{Cut algebras.}
			For an edge cut $a$, let $X_L,X_R$ be the full partial joins on
			the two sides.
			Then
			\begin{equation}\label{eq:linear-native-cut}
				X=X_L\times_{Y_a}X_R,
				\qquad
				\mathbb C^{X_L}\cap\mathbb C^{X_R}=\mathbb C^{Y_a}
				\quad\text{inside }\mathbb C^X,
			\end{equation}
			where all three algebras are included by their actual pullbacks.
			
			\item \emph{Native spectral realization.}
			On $\mathcal H=\ell^2(X)$, regard each finite-field readout as
			its diagonal spectral projections and put
			\[
			M(T)=\operatorname{alg}^{*}\{\Ind[f_e=t]:e\in T,\ t\in\mathbb F_q\},\qquad
			K(T)=M(T)',\qquad T\subseteq V.
			\]
			Then $(\mathrm{com},M(T))$ is an access--commutant native formation
			datum with represented endpoint $\mathfrak P_{K(T)}$, and
			\begin{equation}\label{eq:linear-native-formation}
				K(T)'=M(T),\qquad
				K(T\cup U)=K(T)\cap K(U),\qquad M(V)=\mathbb C^X.
			\end{equation}
			The local sources embed by pullback and their joins construct this
			common native source.
		\end{enumerate}
	\end{theorem}

	\begin{proof}
		\leavevmode
		\begin{enumerate}[label=\textup{(\roman*)},leftmargin=*]
			\item \emph{Linear source and dimension.}
			The matching equations are homogeneous linear equations in the chart coordinates, so their solution set $X$ is a vector space.
			The surjectivity assertion follows by extending a local state outward
			along the join tree.
			Since every edge map is linear and onto $Y_a$, each extension across an edge to a new chart $i$ has $q^{r_i-s_a}$ choices.
			Multiplying these counts from any root proves			\eqref{eq:linear-source-size}; it also constructs the source from the
			local input.
			
			\item \emph{Global determination.}
			Let $W$ be the
			quotient of $\bigoplus_i W_i$ by the relations identifying the two
			copies of every $u\in U_a$ across every tree edge. Its dual is exactly
			$X$. Successively attaching a leaf space by a pushout of injective
			linear maps shows that every $W_i$ embeds in $W$ and gives the same
			dimension formula. Denote the common image of a descriptor form by
			$v_e\in W$.
			The global determination closure is
			$
			c_X(A)=\{e:v_e\in\operatorname{span}(v_A)\}.
			$
			
			Let $C$ be $c_0$-closed, and put $C_i=C\cap B_i$ and $H_i=\operatorname{span}(v_{i,C_i})$.
			Each $C_i$ is locally closed. Visibility gives
			\[
			H_i\cap U_a
			=\operatorname{span}(v_{C\cap S_a})
			=H_j\cap U_a\qquad(a=ij).
			\]
			Write this common subspace as $L_a$. The quotient spaces $W_i/H_i$
			therefore form a tree diagram with injective separator maps
			$U_a/L_a\to W_i/H_i$. Its pushout embeds every $W_i/H_i$ by the same
			leaf argument. This quotient pushout is $W/H$, where
			$H=\sum_i H_i=\operatorname{span}(v_C)$: both are the direct sum
			quotient obtained by imposing the edge identifications and setting
			each $H_i$ to zero. Hence $H\cap W_i=H_i$. If $e\notin C$, choose a
			chart containing $e$; local closedness gives $v_{i,e}\notin H_i$, and
			therefore $v_e\notin H$. Thus $c_X(C)=C$. All local implications hold
			in $X$, so $c_0\leq c_X$; applying the preceding conclusion to
			$C=c_0(A)$ gives the reverse inclusion and proves (ii).
			
			\item \emph{Cut algebras.}
			Deleting a tree edge separates all matching equations except for its
			separator, so the first identity in \eqref{eq:linear-native-cut} holds.
			Both partial joins project onto $Y_a$. Each separator fiber is
			therefore a full nonempty Cartesian rectangle. A function of the left
			partial state that is also a function of the right partial state must
			be constant on every such rectangle, which proves the algebra
			intersection identity.
			
			\item \emph{Native spectral realization.}
			The observables separate global states, so their full diagonal algebra is $\mathbb C^X$.
			Finite-dimensional bicommutant recovery and the commutant-of-a-join identity now prove \eqref{eq:linear-native-formation}.
			The realization judgment follows from \eqref{eq:standalone-commutant-law}.
			Surjectivity of each chart projection makes its pullback on functions injective, and the pulled-back algebras together distinguish all global states.\qedhere
		\end{enumerate}
	\end{proof}

	\Cref{cor:native-linear-source-rigidity} classifies the minimum over the actual atlas.
	Its rank lower bound also applies to any full scalar-linear representation of $c_0$.

	\begin{remark}[Local amplification in the common Hilbert space]
		\label{rem:linear-native-amplification}
		Every fiber of $X\to X_i$ has size $m_i=q^{R-r_i}$.
		After a splitting, the source representation on $\ell^2(X)$ is
		$\mathbb C^{X_i}\otimes I_{m_i}$ and its commutant is
		$\mathbb C^{X_i}\otimes M_{m_i}(\mathbb C)$.
	\end{remark}

	\begin{example}[Invisible separator information creates a new dependency]
		\label{ex:linear-invisible-separator}
		Let $X_1=X_2=\mathbb F_q^2$. The first chart has descriptors
		$\{a,u,v\}$ with $a=u+v$, and the second has
		$\{u,v,b\}$ with $b=u+v$. Match the actual $(u,v)$-values.
		The full join is $\mathbb F_q^2$, both chart projections are
		bijective, and the native separator intersection is exact.
		Nevertheless $a\to b$ holds globally, whereas
		$c_0(\{a\})=\{a\}$. The observation $a$ determines the separator direction
		$\operatorname{span}(u+v)$. The named separator observations in
		$c_1(\{a\})$ span a smaller space. Thus \eqref{eq:linear-visibility}
		fails. This shows why actual state and
		algebra gluing alone do not establish faithful recovery.
	\end{example}

	\begin{example}[The full linear optimum is not an unrestricted optimum]
		\label{ex:linear-source-cost-scope}
		Take $r$ independent binary one-descriptor charts with empty
		separators. The full linear join has $2^r$ states. The subset
		$\{0,e_1,\ldots,e_r\}\subseteq\mathbb F_2^r$ has $r+1$ states, is
		state-exact, and has the same free determination closure: if a
		descriptor $j$ is absent from a premise, the pair $(0,e_j)$ refutes
		its recovery. The separator algebra on any cut is also exact. For
		the last claim, the bipartite incidence graph of the two partial
		state projections is connected, since the zero state connects the
		zero projections and each other row differs on only one side.
		For a full scalar-linear representation with jointly separating
		forms, the dimension of its vector space equals the rank of its
		represented closure. Here that rank is $r$, giving $2^r$ states
		for every such binary linear realization.
	\end{example}

	\subsection{Native construction in finite equational classes}
	\label{subsec:equational-native-bridge}

	To retain the algebraic operations on local states through assembly, fix a native presentation law before choosing the atlas. Identities ensure closure under the algebraic constructions, and the presentation specifies the realized endpoint.
	
	Fix a finite finitary algebraic signature $\tau$ and a finite set
	$\mathcal E$ of identities, and let $\mathcal V$ be the variety of
	nonempty $\tau$-algebras satisfying $\mathcal E$. Bold letters denote
	algebras and unbolded letters their underlying sets. Groups, unital
	rings, lattices, and modules over a fixed finite ring are examples.
	The constructions use the classical closure of equational classes
	under products, subalgebras, and homomorphic images
	\cite{Birkhoff1935,BurrisSankappanavar1981}.

	\begin{definition}[Native presentation context for finite marked algebras]
		\label{def:equational-native-context}
		Fix a $\mathfrak U$-small set of labels $\Lambda$ closed under forming
		finite sets and tuples of its elements.
		Obtain it by closing a $\mathfrak U$-small pool of labels in countably many stages, including the supplied local states and generator names in the initial pool.
		The \textbf{native presentation context for finite marked algebras}
		is the quadruple
		\[
		\Xi_{\mathcal V}^{\mathrm{pres}}
		=\bigl(\mathscr C_{\mathcal V}^{\mathrm{pres}},
		\mathscr C_{\mathcal V}^{\mathrm{pres}},
		\operatorname{Form}_{\mathcal V}^{\mathrm{pres}},
		\mathsf J_{\mathcal V}^{\mathrm{pres}}\bigr),
		\]
		with formation data, carriers, and realization law defined below.
		The prescribed source and endpoint comparisons are given in \textup{(iii)}.
		\begin{enumerate}[label=\textup{(\roman*)},leftmargin=*]
			\item \emph{Finite presentations as formation data.}
			For a finite nonempty set $G\subseteq\Lambda$ of generator names, let
			$F_{\mathcal V}(G)$ be the free algebra in $\mathcal V$.
			The set $\operatorname{Form}_{\mathcal V}^{\mathrm{pres}}$ consists of all
			finite presentations $P=(G,\mathcal R)$, where $\mathcal R$ is a
			finite set of equations between terms on $G$.
			Write $\theta_P$ for the congruence generated by these equations in
			$F_{\mathcal V}(G)$, and set
			\[
			\mathbf A_P=F_{\mathcal V}(G)/\theta_P,
			\qquad \eta_P:G\longrightarrow A_P,\quad g\longmapsto[g]_{\theta_P}.
			\]
			For fixed $G$, comparisons keep generator names fixed unless an
			explicit renaming is specified.
			
			\item \emph{Marked carriers and endpoints.}
			The structural domain and the admissible-carrier set are both
			$\mathscr C_{\mathcal V}^{\mathrm{pres}}$, the $\mathfrak U$-small set
			of marked algebras $K=(G,\mathbf A,\eta)$ in which $G\subseteq\Lambda$
			is finite and nonempty, $\mathbf A\in\mathcal V$ has finite carrier
			$A\subseteq\Lambda$, and $\eta:G\to A$ has image generating $\mathbf A$.
			For each such $K$, define its packaged state endpoint by
			\[
			\mathfrak Q_K=
			\bigl(\{K\},A,\triangleleft_K,\mathsf L_K,\mathsf R_K\bigr),
			\qquad
			a\triangleleft_K K\ (a\in A),\qquad
			\mathsf L_K(K)=\mathsf R_K(K)=A.
			\]
			Thus its selected states are exactly the elements of the displayed
			algebra.
			
			\item \emph{Realization law and comparisons.}
			The realization relation $\mathsf J_{\mathcal V}^{\mathrm{pres}}$ has
			pairs of the form $(P,\mathfrak Q_K)$, with membership defined by
			\begin{equation}\label{eq:equational-native-law}
				\Xi_{\mathcal V}^{\mathrm{pres}}\vdash P
				\rightsquigarrow\mathfrak Q_K
				\quad\Longleftrightarrow\quad
				(\mathbf A,\eta)\cong(\mathbf A_P,\eta_P)
				\text{ as marked algebras}.
			\end{equation}
			A presentation with infinite quotient has no endpoint in this finite
			context.
			A source comparison uses equality of the generated congruences
			$\theta_P$ for the same named generators, or after a declared relabeling;
			an endpoint comparison is marked-algebra isomorphism.
		\end{enumerate}
	\end{definition}

	The marked algebra supplies the five packaged coordinates of \cite[Section~2]{Ma2026} and the selected states used in \cref{def:witness}.
	Its marking recovers the source congruence by evaluation:
	\begin{equation}\label{eq:equational-source-decoder}
		\theta_P
		=\ker\bigl(\operatorname{ev}_{\eta}:
		F_{\mathcal V}(G)\longrightarrow\mathbf A\bigr).
	\end{equation}
	The marking matters.
	For example, the presentations of $C_2$ on $x,y$ given by $x^2=1,y=1$ and
	$x=1,y^2=1$ give isomorphic unmarked groups.
	They generate different congruences on the same free group.

	\begin{lemma}[Finite tables provide effective native formation data]
		\label{lem:finite-table-presentation}
		For a finite algebra $\mathbf A\in\mathcal V$ on an admitted carrier (or transported to one by a specified bijection), introduce one allowed name
		$\dot a$ for each $a\in A$, and let $P[\mathbf A]$ contain the equations
		\begin{equation}\label{eq:algebra-table-presentation}
			\omega(\dot a_1,\ldots,\dot a_r)
			=\dot{\omega^{\mathbf A}(a_1,\ldots,a_r)}
			\qquad(\omega\in\tau\text{ of arity }r).
		\end{equation}
		Nullary operations are included with the corresponding equation
		$\omega=\dot{\omega^{\mathbf A}}$.
		\begin{enumerate}[label=\textup{(\roman*)},leftmargin=*]
			\item \emph{Realization by a table.}
			The table presentation $P[\mathbf A]$ realizes the marked state endpoint of $\mathbf A$ in \eqref{eq:equational-native-law}.
			
			\item \emph{Maps induced by homomorphisms.}
			For a second such finite algebra $\mathbf B$, a homomorphism $h:\mathbf A\to\mathbf B$ induces the explicit	presentation map $\dot a\mapsto\dot{h(a)}$.
			It respects the table equations and induces $h$ on the quotients.
		\end{enumerate}
	\end{lemma}

	\begin{proof}
		\leavevmode
		\begin{enumerate}[label=\textup{(\roman*)},leftmargin=*]
			\item \emph{Realization by a table.}
			Evaluation in $\mathbf A$ satisfies the equations and maps every
			$\dot a$ to $a$. Induction on terms using the table equations reduces
			every term to the name of its value in $\mathbf A$. Distinct names
			cannot become equal in the presented quotient, since evaluation would
			then identify distinct elements of $A$. Consequently evaluation induces
			a marked isomorphism of the presented quotient with $\mathbf A$.
			
			\item \emph{Maps induced by homomorphisms.}
			For the assertion about $h$, the image of each displayed equation is
			a valid table equation in $\mathbf B$, precisely because $h$ preserves
			$\omega$. This also proves compatibility with identities and with
			composition of homomorphisms.\qedhere
		\end{enumerate}
	\end{proof}

	Membership is decidable from the finite operation tables by testing
	the identities $\mathcal E$ and checking that the specified maps
	preserve every basic operation.

	\subsubsection{Native quotients of marked presentations}
	\label{app:canonical-native-proofs}
	\label{subsec:canonical-reduction-native-verification}
	Use the context of \cref{def:equational-native-context}, with finite
	equation lists as formation data. At each native algebra sort include
	its basic operations in $\mathcal F$. Then the component of an
	$\mathcal F$-congruence on that sort is an algebra congruence, so it
	can be used in a native algebra quotient. Let
	$P=(G,\mathcal R)$ realize the finite marked algebra
	$(\mathbf A,\eta)$ and write
	$\gamma=\ker(\operatorname{ev}_\eta)=\theta_P$. For every $a\in A$
	choose a term $t_a$ on $G$ with $t_a^{\mathbf A}(\eta)=a$; such terms
	exist because the retained marking generates $\mathbf A$. Given a
	congruence $\theta\in\Con(\mathbf A)$, form the finite presentation
	\begin{equation}\label{eq:canonical-reduction-finite-quotient-presentation}
		P_\theta=(G,\mathcal R_\theta),\qquad
		\mathcal R_\theta
		=\mathcal R\cup\{t_a=t_b:(a,b)\in\theta\}.
	\end{equation}
	Its generated congruence is exactly
	\[
	\theta_{P_\theta}
	=(\operatorname{ev}_\eta\times\operatorname{ev}_\eta)^{-1}(\theta).
	\]
	Indeed the displayed equations hold in $\mathbf A/\theta$, giving
	one inclusion. Conversely, every term $t$ is congruent modulo
	$\gamma$ to $t_{t^{\mathbf A}(\eta)}$. If the values of terms
	$t,u$ are $\theta$-equivalent, the additional equation between their
	chosen representatives therefore makes $t,u$ congruent modulo
	$\theta_{P_\theta}$. This proves the reverse inclusion and gives
	\[
	F_{\mathcal V}(G)/\theta_{P_\theta}
	\cong\mathbf A/\theta
	\]
	over the original named generators. The right side is finite and
	belongs to the same variety, so \eqref{eq:equational-native-law}
	supplies its native endpoint without changing the realization law.
	The retained marking is $q\eta$, and
	\eqref{eq:equational-source-decoder} recovers
	$\theta_{P_\theta}$ from that marked endpoint.
	
	The comparison is the homomorphism $[t]_\gamma\mapsto
	[t]_{\theta_{P_\theta}}$. Identity and successive quotient maps
	compose on these term classes. Different choices of representative
	terms yield the same generated congruence and hence the same source
	comparison. The algebraic quotient construction
	\cite{BurrisSankappanavar1981} keeps the named generators and
	verifies (NQ) for the presentation law.
	For several supplied sources the construction applies at each native
	sort, while the $\mathcal F$-congruence condition makes all declared
	cross-sort operations and interfaces descend.

	\subsubsection{Native realization of the matching algebra}

	\begin{theorem}[Native amalgamation for finite equational join-tree atlases]
		\label{thm:equational-native-amalgamation}
		Let $(B_i)_{i\in I}$ be a finite descriptor cover with a join tree $J$.
		Suppose that each local state algebra $\mathbf A_i$ and each descriptor
		algebra $\mathbf D_d$ is a finite member of $\mathcal V$, and that
		$
		f_{i,d}:\mathbf A_i\longrightarrow\mathbf D_d
		\quad(d\in B_i)
		$
		are homomorphisms whose joint map is injective.
		On each edge $ij$, assume equality of the actual separator images on $S_{ij}=B_i\cap B_j$, and denote the common image algebra by $\mathbf Y_{ij}$.
		Retain a specified realized local presentation
		$P_i=(G_i,\mathcal R_i)$ and its generating marking
		$\eta_i:G_i\to A_i$ at each chart, so that
		\[
		\Xi_{\mathcal V}^{\mathrm{pres}}\vdash P_i
		\rightsquigarrow\mathfrak Q_{(G_i,\mathbf A_i,\eta_i)}.
		\]
		The finite table presentation is one effective instance of these local data.
		For other supplied presentations, their displayed realization judgments are part of the local input.
		Construct
		\begin{equation}\label{eq:equational-witness-algebra}
			W=\Bigl\{(a_i)\in\prod_i A_i:
			f_{i,d}(a_i)=f_{j,d}(a_j)
			\text{ whenever }d\in B_i\cap B_j\Bigr\}.
		\end{equation}
		Then the following hold.
		\begin{enumerate}[label=(\roman*)]
			\item \emph{Matching algebra and native formation.}
			Coordinatewise operations make $W$ a nonempty algebra
			$\mathbf W\in\mathcal V$.
			Every local projection $\pi_i:\mathbf W\twoheadrightarrow\mathbf A_i$ is a surjective homomorphism.
			The formation datum $P[\mathbf W]$ and its marked state endpoint satisfy \eqref{eq:equational-native-law}.
			The	projections are induced by explicit presentation maps into the
			retained $P_i$, with their original local markings and realization
			judgments.
			This constructs a native common source from locally	checked algebra and interface data.
			
			\item \emph{Universal property.}
			The algebra $\mathbf W$ has the pullback universal property.
			Take homomorphisms from one algebra to the $\mathbf A_i$ that agree on the declared separators.
			They factor through a unique homomorphism into $\mathbf W$.
			
			\item \emph{Cut kernels.}
			For every cut edge, let $\mathbf W_L,\mathbf W_R$ be the partial witness algebras and $\mathbf Y$ the separator algebra.
			Then $\mathbf W\cong\mathbf W_L\times_{\mathbf Y}\mathbf W_R$.
			If $\alpha_L,\alpha_R$ are the kernels of the two projections and $\sigma$ is the kernel of the separator map, then
			\begin{equation}\label{eq:equational-exact-kernels}
				\alpha_L\cap\alpha_R=\Delta_W,
				\qquad
				\alpha_L\circ\alpha_R
				=\alpha_R\circ\alpha_L
				=\alpha_L\vee\alpha_R=\sigma.
			\end{equation}
			Here $\vee$ is join of congruences.
			In particular
			$\mathbf W/\alpha_L\cong\mathbf W_L$,
			$\mathbf W/\alpha_R\cong\mathbf W_R$, and
			$\mathbf W/\sigma\cong\mathbf Y$.
			
			\item \emph{Boundary algebras and substitution.}
			The image of $\mathbf W$ on any specified descriptor boundary
			is a finite algebra in the same variety and has its own native
			table presentation. For a grafting whose combined witness diagram
			satisfies the same compatibility hypotheses, the substitution
			identity is an identity of native boundary algebras under the
			canonical identifications of their coordinates.
			
			\item \emph{Separator information and recovery.}
			The represented separator function algebras have the exact
			intersections of \cref{thm:local-pair-native} on every full tree cut.
			If the underlying atlas also satisfies relative pair lifting at every
			incidence, its descriptor determination closure is the least closure
			generated by the local closures:
			$
			c_W=c_0.
			$
			
			\item \emph{Source size.}
			The construction is finite, with $|W|\leq\prod_i|A_i|$.
			For finite groups, or finite modules over a fixed ring, its size is
			\begin{equation}\label{eq:equational-group-source-count}
				|W|=\frac{\prod_i|A_i|}
				{\prod_{ij\in E(J)}|Y_{ij}|}.
			\end{equation}
		\end{enumerate}
	\end{theorem}

	\begin{proof}
		\leavevmode
		\begin{enumerate}[label=\textup{(\roman*)},leftmargin=*]
			\item \emph{Matching algebra and native formation.}
			Root $J$ at any chart. A state at that chart extends successively
			along the tree, because the two maps on each edge have the same
			image. Running intersection ensures that the extended states match
			on every repeated descriptor. This proves nonemptiness and, by
			rooting at each chart in turn, every local surjectivity assertion.
			
			For a basic $r$-ary operation $\omega$ and $r$ matching tuples,
			apply $\omega$ at every chart. The two separator images of the
			result agree, since the separator maps are homomorphisms and their
			images agreed in each input tuple. The same argument includes
			nullary operations. Thus $W$ is a subalgebra of the product.
			Every identity in $\mathcal E$ holds coordinatewise in the product
			and hence in this subalgebra. The projections are homomorphisms.
			The complete table of $\mathbf W$ gives an actual global formation
			datum by \cref{lem:finite-table-presentation}. To retain a given
			local presentation $P_i$, choose for each $a\in A_i$ a term
			$t_{i,a}$ on the original generators $G_i$ whose value under
			$\eta_i$ is $a$. These terms can be found by closing the marked
			generators under the finite operation tables. The assignment
			$
			\dot w\longmapsto t_{i,\pi_i(w)}
			$
			defines a map of free algebras. It carries every global table
			equation to a pair with equal values in $\mathbf A_i$, hence to
			the original congruence $\theta_{P_i}$. It therefore induces the
			specified projection onto the original marked quotient. Different
			choices of representative terms give the same induced map.
			The projection maps retain the original equation lists and markings.
			One can enumerate the extensions above, compute all global
			operation values coordinatewise, and output this presentation and
			its source maps.
			\item \emph{Universal property.}
			The coordinate tuple of any compatible family of
			homomorphisms is the unique factorization in (ii), and it is a
			homomorphism by the same coordinatewise calculation.
			\item \emph{Cut kernels.}
			Across a tree cut the only matching requirements between the two
			partial diagrams are their common separator values. Both partial
			maps onto $Y$ are surjective. Each separator fiber is therefore a
			full nonempty rectangle. Equality of both partial coordinates gives
			$\alpha_L\cap\alpha_R=\Delta_W$. If two global states have the same
			separator value, the tuple with the left coordinate of the first
			and the right coordinate of the second is another global state.
			It witnesses membership in $\alpha_L\circ\alpha_R$; the reversed
			intermediate tuple gives the other composition. Conversely either
			composition preserves the separator value. The common composition
			is thus the congruence $\sigma$, which is also the join of the two
			kernels. The quotient identifications follow from surjectivity and
			the homomorphism theorem.
			\item \emph{Boundary algebras and substitution.}
			A boundary readout is a homomorphism to a product of the
			$\mathbf D_d$. Its image inherits the operations, satisfies every
			identity in $\mathcal E$, and is isomorphic to the corresponding
			quotient of $\mathbf W$. It therefore has a native table
			presentation. For a combined diagram satisfying the theorem's
			hypotheses, nonemptiness has already been proved. Existential
			composition of its boundary relations has the same tuples as the
			boundary of the grafted witness diagram,
			by \cref{prop:substitution}. Both carry the inherited coordinatewise
			operations, so this is also equality of the boundary algebras.
			\item \emph{Separator information and recovery.}
			The rectangular-fiber argument gives the separator
			function-algebra identities independently of the algebraic
			presentation law. Under relative pair lifting, \cref{thm:local-pair-native} applied to the underlying atlas gives $c_W=c_0$.
			\item \emph{Source size.}
			Since $W$ is a subset of the product of the local states, $|W|\leq\prod_i|A_i|$. Each surjective group or module map
			onto $Y_{ij}$ has constant fiber size $|A_i|/|Y_{ij}|$.
			Multiplying the extension counts along the rooted tree proves
			\eqref{eq:equational-group-source-count}.\qedhere
		\end{enumerate}
	\end{proof}
	
	\begin{proposition}[Construction and recognition of finite native atlas formations]
		\label{prop:equational-levelii-bridge}
		Fix the finite equational join-tree atlas of
		\cref{thm:equational-native-amalgamation}, including its original
		realized, marked local charts and jointly separating readouts
		$f_{i,d}:\mathbf A_i\to\mathbf D_d$.
		For every finite $\mathbf H\in\mathcal V$ with an admitted
		state carrier in the fixed native context and with compatible
		homomorphisms $h_i:\mathbf H\to\mathbf A_i$, the construction in \textup{(i)}--\textup{(iii)} gives a coherent, hierarchical Level-II presentation. Part~\textup{(iv)} characterizes the resulting formation relations.
		\begin{enumerate}[label=\textup{(\roman*)},nosep]
			\item \emph{Sites and restrictions.} Its inherited input types $K_i$ have selected states $A_i$;
			its prospective output type $K_*$ has selected states $H$.
			For each descriptor $d$ there is one typed site $U_d$, with
			$\LL(U_d)=D_d$, lying below $K_*$ and each $K_i$ that names $d$.
			The input restrictions are $f_{i,d}$ and the output restriction is
			$g_d=f_{i,d}h_i$, independent of the chart containing $d$.
			\item \emph{Matching and witnesses.} A single signature $\alpha:(K_i)_{i\in I}\to K_*$
			compares each input with the output at all its named descriptors.
			Its operational relation is exactly the full matching relation,
			\begin{equation}\label{eq:levelii-native-graph-relation}
				\MM_\alpha=\Match_\alpha
				=\{((h_i(z))_{i\in I},z):z\in H\}.
			\end{equation}
			Consequently the one-step witness set is canonically $H$, and
			its projection to the input vertices is $\psi(H)\subseteq W$. The
			boundary including the output is the whole graph in
			\eqref{eq:levelii-native-graph-relation}. Here
			$\psi=(h_i)_i:\mathbf H\to\mathbf W$ is the native comparison map.
			\item \emph{Native realization and preservation.} The output has the native realization $P[\mathbf H]$ in the
			fixed presentation context; the inherited native
			realization judgments and markings are unchanged.
			All four tests---chart-state coverage, descriptor determination,
			exact full-cut common information, and native realization---have
			exactly the same truth values when evaluated on the Level-II
			witnesses and on $H$ with its prescribed maps $h_i$.
			\item \emph{Recognition from the formation relation.} Conversely, restrict attention to one-step presentations with precisely
			these input and descriptor types, with a finite native output algebra
			$\mathbf H$, homomorphic output restrictions
			$g_d:\mathbf H\to\mathbf D_d$, and the declared comparisons of
			\textup{(ii)}. Assume interface soundness $\MM_\alpha\subseteq\Match_\alpha$
			and that every selected output $z\in H$ occurs in an admitted
			formation. Then there are \emph{unique} compatible homomorphisms
			$h_i:\mathbf H\to\mathbf A_i$ with $g_d=f_{i,d}h_i$, and
			\eqref{eq:levelii-native-graph-relation} holds.
			Thus the compatible native source cones over the fixed atlas correspond
			to these one-step formations through their readouts and admitted tuples.
			Within this class, soundness and output coverage already force admission
			of every matching tuple.
		\end{enumerate}
	\end{proposition}
	\begin{proof}
		\leavevmode
		\begin{enumerate}[label=\textup{(\roman*)},leftmargin=*]
			\item \emph{Sites and restrictions.}
			Retain the original realized charts as bases. Use distinct site
			tags for the input occurrences, the output, and the descriptors,
			even when the underlying algebras coincide. The only strict incidences
			are $U_d<K_i$ for $d\in B_i$ and $U_d<K_*$. Give these incidences
			the restrictions $f_{i,d}$ and $g_d$, respectively, and use identities
			on reflexive incidences. Compatibility of the $h_i$ makes
			each $g_d$ independent of the chart in which $d$ occurs. This
			establishes restriction coherence \textup{(H1)}.
			\item \emph{Matching and witnesses.}
			A tuple $((a_i)_i,z)$ satisfies every input--output comparison
			precisely when
			$f_{i,d}(a_i)=f_{i,d}h_i(z)$ for all $d\in B_i$ and all $i$.
			Since the chart readouts are jointly injective, this holds if and
			only if $a_i=h_i(z)$ for every $i$. Hence the full matching
			relation is exactly the graph in
			\eqref{eq:levelii-native-graph-relation}. Taking it as the
			operational relation establishes \textup{(H2)}, and assigning rank
			zero to the inputs and rank one to the output establishes
			\textup{(H3)}.
			\item \emph{Native realization and preservation.}
			The native table presentation of $H$ exists by
			\cref{lem:finite-table-presentation}; this does not alter the
			original chart presentations.
			
			Projection to the output coordinate identifies the
			witness set with $H$. Under this identification, chart readouts
			are precisely the $h_i$, descriptor readouts are $g_d$, and
			complete-side and separator readouts are their native composites.
			Thus chart surjectivity, all positive and negative descriptor
			determinations, and every equality of finite function algebras
			at a full cut are literally the same tests on the two sides.
			The output judgment is supplied by $P[\mathbf H]$.
			
			If a finite algebra is initially given with an unadmitted literal
			carrier, one first transports it along a chosen bijection to an
			admitted finite set of labels; such relabeling is an additional
			identification, not a change of the native law.
			\item \emph{Recognition from the formation relation.}
			For the converse, every $z\in H$ has some admitted input tuple
			by output-totality. Joint injectivity of the readouts on $A_i$
			makes its $i$th input state unique; call it $h_i(z)$.
			The comparison equations give $g_d=f_{i,d}h_i$.
			Because every $g_d$ and $f_{i,d}$ is a homomorphism,
			for each basic $r$-ary operation $\omega$ we have
			\[
			f_{i,d}\bigl(h_i(\omega^H(z_1,\ldots,z_r))\bigr)
			=\omega^{D_d}(g_d(z_1),\ldots,g_d(z_r))
			=f_{i,d}\bigl(\omega^{A_i}(h_i(z_1),\ldots,h_i(z_r))\bigr).
			\]
			Joint injectivity proves that $h_i$ preserves $\omega$;
			nullary operations follow by the same calculation.
			Shared $g_d$ give compatibility. For a fixed $z$, any matching
			input has coordinate $h_i(z)$ at every $i$ by joint injectivity.
			Output coverage puts this unique tuple in $\MM_\alpha$. Thus
			$\Match_\alpha\subseteq\MM_\alpha$, and soundness gives the reverse inclusion.
			The two constructions recover one another on the specified
			readouts and witnessed formation data, up to the evident
			relabeling of finite-table generators.\qedhere
		\end{enumerate}
	\end{proof}
	
	\paragraph{Structural meaning of the native faithfulness test.}
	For the normalized graph formations of \cref{prop:equational-levelii-bridge}, the full matching algebra $\mathbf W$ supplies a candidate for the entire Level-II contract. Any compatible output maps to $\mathbf W$. Hence a counterexample pair missing from $\mathbf W$ is missing from every candidate, including those with extra states invisible to the charts.
	
	\begin{theorem}[Faithful native existence and the ordered-pair obstruction]
		\label{thm:amalgamation-existence-criterion}
		Assume the finite equational join-tree data, equal actual separator
		images, and realized marked local presentations of
		\cref{thm:equational-native-amalgamation}. Let $\mathbf W$ be the
		canonical matching algebra of \eqref{eq:equational-witness-algebra}
		and let $c_0$ be the determination closure generated by the actual
		local chart theories. The following are equivalent.
		\begin{enumerate}[label=(\roman*)]
			\item \emph{Existence of a native amalgamation.} There exists a finite native algebra $\mathbf H$ with compatible
			homomorphic maps $h_i:\mathbf H\to\mathbf A_i$ satisfying the full
			$\mathsf S,\mathsf F,\mathsf E,\mathsf N$ preservation contract in the
			fixed presentation context.
			\item \emph{Faithfulness of the matching algebra.} The canonical matching algebra is recovery-faithful:
			\begin{equation}\label{eq:existence-canonical-faithfulness}
				c_W=c_0.
			\end{equation}
			\item \emph{Ordered-pair witnesses.} For every $S\subseteq V$ and every $d\in V\setminus c_0(S)$,
			the ordered-pair root test of \cref{thm:pair-descent} accepts $(S,d)$;
			equivalently, two tuples of $W$ agree on every descriptor in $S$
			and disagree on $d$.
		\end{enumerate}
		
		\begin{itemize}
			\item \emph{Realization.} Whenever these conditions hold, $\mathbf W$ itself is a native
			amalgamation satisfying the full contract. A finite table
			presentation $P[\mathbf W]$ supplies its actual formation datum;
			the local projections are induced by maps to the \emph{retained}
			marked local presentations. The resulting matching algebra has the
			pullback universal property, the exact full-cut separator algebras,
			and the native boundary-substitution properties established in
			\cref{thm:equational-native-amalgamation}.
			
			\item \emph{Witnesses and obstructions.} The tests in \textup{(iii)} give an exact finite procedure for this
			fixed atlas. Successful root tests can be traced back to explicit
			matching witness pairs. If a test fails for $(S,d)$, then
			$d\notin c_0(S)$ but $d\in c_W(S)$; consequently \emph{no}
			compatible nonempty source over the prescribed charts---even one
			with no algebraic operations or state-surjectivity requirement---can
			have descriptor closure $c_0$.
		\end{itemize}
	\end{theorem}
	\begin{proof}
		\leavevmode
		\begin{itemize}
			\item \emph{\textup{(ii)} $\implies$ \textup{(i)} and native realization.}
			By \cref{thm:equational-native-amalgamation}, the finite matching
			algebra $\mathbf W$ is nonempty, all its chart projections are onto,
			and every local implication holds in $W$. Hence
			$c_0\subseteq c_W$. Moreover, the same theorem constructs
			$P[\mathbf W]$, maps into the original marked local presentations,
			and exact separator algebras at every full cut. Thus
			\textup{(ii)} implies \textup{(i)}, including $\mathsf N$ and
			$\mathsf E$.
			\item \emph{\textup{(i)} $\implies$ \textup{(ii)}.}
			Conversely, let $\mathbf H$ satisfy \textup{(i)}. Its compatible
			chart maps induce a homomorphism
			$q:\mathbf H\to\mathbf W$. Every dependency valid on $W$ holds
			on $H$ after pullback. Therefore
			\[
			c_0\subseteq c_W\subseteq c_H=c_0,
			\]
			which proves \textup{(ii)}. This particular implication does not
			use state surjectivity, cut exactness, or the algebraic operations
			of $H$: it holds for any compatible set with the prescribed readouts.
			\item \emph{\textup{(ii)} $\Longleftrightarrow$ \textup{(iii)}, witnesses and obstructions.}
			The doubled-relation join-tree test of \cref{thm:pair-descent}
			accepts $(S,d)$ exactly when $d\notin c_W(S)$. Since
			$c_0\subseteq c_W$, it follows that \textup{(ii)} and
			\textup{(iii)} are equivalent. The finite recursion stores
			extendible ordered separator pairs and permits traceback to two
			global matching tuples. If a root test rejects, the forced
			implication $S\to d$ on $W$ persists under every compatible
			pullback, proving the final obstruction assertion.\qedhere
		\end{itemize}
	\end{proof}
	
	The construction combines finite presentations by operation tables with the acyclic join algorithm on pairs of chart states \cite{Birkhoff1935,Beeri1983,Yannakakis1981}.
	
	\begin{corollary}[Faithful native descent over arbitrary two-state equational interfaces]
		\label{cor:two-state-equational-descent}
		In \cref{thm:equational-native-amalgamation}, suppose in addition
		that $|Y_{ij}|\leq2$ at every join-tree edge.
		\begin{enumerate}[label=\textup{(\roman*)},leftmargin=*]
			\item \emph{Faithfulness and state coverage.}
			Relative pair lifting holds at each incidence.
			The full matching algebra $\mathbf W$ has $c_W=c_0$ and covers
			every original state.
			\item \emph{Native formation and cut data.}
			It has exact information at the full tree cuts,
			native boundary images, and the stated substitution properties
			under the same presentation law. The result applies in every finite-signature equational variety with this law.
		\end{enumerate}
	\end{corollary}

	\begin{proof}
		\leavevmode
		\begin{enumerate}[label=\textup{(\roman*)},leftmargin=*]
			\item \emph{Faithfulness and state coverage.}
			Forget the operations temporarily and apply
			\cref{thm:two-state-joint-separators} to the underlying spectral
			atlas: joint two-state interfaces automatically satisfy relative
			pair lifting. Part~\textup{(v)} of
			\cref{thm:equational-native-amalgamation} then gives $c_W=c_0$.
			
			\item \emph{Native formation and cut data.}
			Its parts~\textup{(i)}--\textup{(iv)} provide the native algebra,
			original marked presentation maps, cut identities, boundary
			formation, and compatible substitution without any Mal'cev
			assumption.\qedhere
		\end{enumerate}
	\end{proof}

	\begin{corollary}[Native realization of algebraic witnesses]
		\label{cor:operational-native-algebra}
		In the formation-tree semantics of \cref{def:witness}, suppose that
		each selected state space is the underlying set of a finite native
		state algebra in $\Xi_{\mathcal V}^{\mathrm{pres}}$. Suppose each
		primitive formation relation is a nonempty subalgebra of the product
		of its input and output state algebras. At each state variable shared
		by two adjacent formation occurrences, assume the two relations have
		the same projection onto that variable.
		\begin{enumerate}[label=\textup{(\roman*)},leftmargin=*]
			\item \emph{Witness and boundary algebras.}
			The witness set and its boundary relation have native algebraic presentations in this same context.
			
			\item \emph{Local extension and substitution.}
			Every local	formation tuple extends to a full witness, and grafting and
			substitution agree as native boundary-algebra constructions.
			
			\item \emph{Dependency faithfulness.}
			If the formation-tuple charts also satisfy relative pair lifting,
			the full witness algebra preserves exactly their generated
			descriptor recovery theory.
		\end{enumerate}
	\end{corollary}

	\begin{proof}
		Index a chart by each nonleaf formation occurrence, with its parent
		and child state variables as descriptors and its formation tuples
		as local states. Two such charts are adjacent precisely when they
		share the state variable at the child occurrence. They form a join tree. At each edge the two projections have the
		same image by hypothesis, and coordinate readouts jointly separate
		formation tuples.
		\begin{enumerate}[label=\textup{(\roman*)},leftmargin=*]
			\item \emph{Witness and boundary algebras.}
			Parts~\textup{(i)} and \textup{(iv)} of \cref{thm:equational-native-amalgamation} give the native witness and boundary presentations.
			For a one-vertex tree the state algebra itself is both the witness and the boundary algebra.
			
			\item \emph{Local extension and substitution.}
			Surjectivity of the chart projections extends every local tuple to a full witness. Boundary substitution is part~\textup{(iv)} of the same theorem.
			
			\item \emph{Dependency faithfulness.}
			Under the pair-lifting hypothesis, part~\textup{(v)} identifies the determination closure with the one generated by the formation-tuple charts. In the one-vertex case this is the determination closure of the single state chart.\qedhere
		\end{enumerate}
	\end{proof}

	The equal-image condition is also necessary for every local formation
	tuple to extend: at either occurrence, every projected value in a local
	tuple must then occur at the other end of the edge. In the additive-group version of
	\cref{ex:linear-invisible-separator}, the readouts are homomorphisms
	and the source is native. The join still gives the extra implication
	$a\to b$.

	\begin{remark}[Closure under boundary images]
		\label{rem:equational-native-scope}
		The boundary image is a subalgebra of a product of descriptor algebras. It therefore satisfies every universal Horn axiom shared by those algebras. Extending this observation to arbitrary congruence quotients would require closure under those quotients.
	\end{remark}

	\subsection{Pair lifting across homomorphic interfaces}
	\label{subsec:multidescriptor-native-ports}

	Pair lifting at a homomorphic interface requires the image of each agreement congruence to be transitive and to equal the intersection of the named separator kernels that contain it. Three-permutability supplies transitivity. The remaining condition tests the actual observation signature.

	Fix an equational class $\mathcal V$ and its presentation
	context as in \cref{def:equational-native-context}. Write
	$\operatorname{Con}(\mathbf A)$ for the congruence lattice, ordered by
	inclusion, with bottom $\Delta_A$ and top $\nabla_A=A^2$.
	For a finite chart algebra $\mathbf A\in\mathcal V$ with jointly
	separating homomorphic readouts $f_d:\mathbf A\to\mathbf D_d$ to
	finite algebras in $\mathcal V$, $d\in B$, and $T\subseteq B$, put
	\begin{equation}\label{eq:malcev-local-agreement}
		\theta_T=\bigcap_{d\in T}\ker f_d,
		\qquad c_{\mathbf A}(T)=\{d:\theta_T\subseteq\ker f_d\}.
	\end{equation}
	The empty intersection is $\nabla_A$. For a separator $S\subseteq B$,
	let $p:\mathbf A\twoheadrightarrow\mathbf Y$ be its actual joint
	image, and let $\rho_s:\mathbf Y\to\mathbf D_s$, $s\in S$, be the
	jointly separating coordinate readouts. Set
	\[
	\kappa=\ker p,\qquad \sigma_s=\ker\rho_s,\qquad
	\beta_T=(p\times p)(\theta_T),\qquad
	\lambda_T=(\theta_T\vee\kappa)/\kappa.
	\]
	Here $\mathbf A/\kappa$ is identified with $\mathbf Y$.
	The least congruence containing $\beta_T$ is $\lambda_T$. The
	inclusion can be strict.
	For $\lambda\in\operatorname{Con}(\mathbf Y)$ define
	\begin{equation}\label{eq:congruence-interface-closure}
		\gamma_S(\lambda)=
		\bigcap_{\substack{s\in S\\\lambda\subseteq\sigma_s}}\sigma_s,
		\qquad \bigcap\varnothing=\nabla_Y.
	\end{equation}

	\begin{proposition}[The two exact obstructions to homomorphic pair lifting]
		\label{prop:two-interface-obstructions}
		\leavevmode
		\begin{enumerate}[label=\textup{(\roman*)},leftmargin=*]
			\item \emph{Local criterion.}
			In any fixed equational class, relative pair lifting at
			$(\mathbf A,S)$ holds if and only if, for every $T\subseteq B$,
			\begin{equation}\label{eq:two-interface-obstructions}
				\underbrace{\kappa\circ\theta_T\circ\kappa
					=\theta_T\vee\kappa}_{\text{congruence-image condition}}
				\quad\text{and}\quad
				\underbrace{\gamma_S(\lambda_T)=\lambda_T}_{\text{named-kernel condition}}.
			\end{equation}
			Equivalently, $\beta_T=\lambda_T=\gamma_S(\lambda_T)$ for all $T$.
			It suffices to test locally closed $T$.
			
			\item \emph{Obstruction to repair.}
			If the first condition fails at $T$, the failure persists under all
			additions of homomorphic observations of that separator. This
			statement fixes the chart, its interface map, and the old observations.
			
			\item \emph{Repair by quotient readouts.}
			Conversely, if the first condition holds
			for every old premise, adjoining the quotient readouts
			$\mathbf Y\to\mathbf Y/\lambda_T$ for those premises gives pair
			lifting for the entire enlarged descriptor vocabulary.
		\end{enumerate}
	\end{proposition}

	\begin{proof}
		\leavevmode
		\begin{enumerate}[label=\textup{(\roman*)},leftmargin=*]
			\item \emph{Local criterion.}
			The inverse image of $\beta_T$ along $p\times p$ is
			$\kappa\circ\theta_T\circ\kappa$. The least congruence of
			$\mathbf Y$ containing $\beta_T$ is $\lambda_T$: by quotient
			correspondence, a congruence above $\beta_T$ pulls back to a congruence
			containing both $\theta_T$ and $\kappa$, and conversely. Hence
			\[
			\beta_T\subseteq\lambda_T\subseteq\gamma_S(\lambda_T),\qquad
			s\in c_{\mathbf A}(T)\cap S
			\Longleftrightarrow \beta_T\subseteq\sigma_s
			\Longleftrightarrow \lambda_T\subseteq\sigma_s.
			\]
			The right side of the pair-lifting equation is therefore
			$\gamma_S(\lambda_T)$. Equality with $\beta_T$ holds exactly when
			both inclusions above are equalities. Surjectivity of $p$ identifies
			the first equality with the first condition in
			\eqref{eq:two-interface-obstructions}. Also
			$\theta_{c_{\mathbf A}(T)}=\theta_T$, which proves the closed-premise
			reduction.
			
			\item \emph{Obstruction to repair.}
			If $\beta_T\ne\lambda_T$, choose a pair in
			$\lambda_T\setminus\beta_T$. Every additional homomorphic readout
			locally determined by $T$ has a kernel containing $\beta_T$, and
			therefore $\lambda_T$. The pair passes every such test. Its set of lifts agreeing on $T$
			is still empty. This proves the
			non-repair assertion.
			
			\item \emph{Repair by quotient readouts.}
			For the converse, let $S^+$ be the enlarged
			separator signature. Its named kernels include every $\lambda_T$,
			so $\gamma_{S^+}(\lambda_T)=\lambda_T$. Any new premise consists
			of an old premise $T$ and some added separator readouts, whose
			kernels on $Y$ have intersection $\eta$. Its agreement image is
			\[
			(p\times p)\bigl(\theta_T\cap(p\times p)^{-1}(\eta)\bigr)
			=\beta_T\cap\eta=\lambda_T\cap\eta.
			\]
			Both $\lambda_T$ and $\eta$ are intersections of named kernels;
			so is their intersection. Hence that image is a congruence fixed by
			$\gamma_{S^+}$, and the pair-lifting equation holds for every new
			premise as well.\qedhere
		\end{enumerate}
	\end{proof}
	
	The first condition of \cref{prop:two-interface-obstructions} is the classical obstruction to preserving congruences under regular images \cite{CarboniEtAl1993,Gran2020Regular}. Once it holds, naming the finitely many realized quotients by $\lambda_T$ suffices. The enriched chart theories determine the new target closure. \Cref{cor:minimal-congruence-interface,thm:minimum-compositional-interface} give the minimum signatures for all charts and for a specified collection, respectively.

	\paragraph{Three-permutable interface algebras.}
	For the automatic interface theorem suppose that every finite algebra
	in $\mathcal V$ is \emph{three-permutable}, meaning
	\begin{equation}\label{eq:finite-three-permutability}
		\alpha\circ\kappa\circ\alpha
		=\kappa\circ\alpha\circ\kappa
		\qquad(\alpha,\kappa\in\operatorname{Con}(\mathbf A)).
	\end{equation}
	A standard sufficient equational condition is a pair of fixed ternary
	terms $u,v$ with
	\begin{equation}\label{eq:fixed-goursat-terms}
		u(x,y,y)=x,\qquad u(x,x,y)=v(x,y,y),\qquad v(x,x,y)=y.
	\end{equation}
	These are the classical Hagemann--Mitschke terms
	\cite{HagemannMitschke1973}.
	A variety satisfying them is also called a Goursat variety.
	For the finite attachment results, assume \eqref{eq:finite-three-permutability} on the finite algebras.
	A Mal'cev term
	\begin{equation}\label{eq:fixed-malcev-term}
		m(x,y,y)=x,\qquad m(x,x,y)=y
	\end{equation}
	is sufficient, by taking $u=m$ and $v(x,y,z)=z$.
	Groups, modules, and rings retain their earlier role as examples.
	The native signature, identities, and presentation law remain fixed before the atlas.
	
	\begin{lemma}[Congruence images and pair lifting under three-permutability]
		\label{lem:malcev-port-image}
		Assume \eqref{eq:finite-three-permutability}.
		\begin{enumerate}[label=\textup{(\roman*)},leftmargin=*]
			\item \emph{Congruence image.} Every $\beta_T$ is the congruence $\lambda_T$.
			
			\item \emph{Lifting test.}
			Relative pair lifting is equivalent to
			\begin{equation}\label{eq:multidescriptor-congruence-test}
				\qquad
				(p\times p)(\theta_T)
				=\bigcap_{s\in c_{\mathbf A}(T)\cap S}\sigma_s
				=\gamma_S(\beta_T)
				\quad(T\subseteq B).
				\qquad
			\end{equation}
			It is equivalent to test only locally closed $T$.
		\end{enumerate}
	\end{lemma}

	\begin{proof}
		\leavevmode
		\begin{enumerate}[label=\textup{(\roman*)},leftmargin=*]
			\item \emph{Congruence image.}
			For congruences $\alpha,\kappa$, the relation
			$R=\kappa\circ\alpha\circ\kappa$ is reflexive, symmetric, and
			compatible.
			Three-permutability gives
			\[
			R\circ R
			=\kappa\circ(\alpha\circ\kappa\circ\alpha)\circ\kappa
			=\kappa\circ(\kappa\circ\alpha\circ\kappa)\circ\kappa=R.
			\]
			Thus $R$ is a congruence and, since it contains $\alpha,\kappa$,
			it equals $\alpha\vee\kappa$.
			Taking $\alpha=\theta_T$ proves $\beta_T=\lambda_T$ by quotient correspondence.
			
			\item \emph{Lifting test.}
			Applying the named-kernel criterion of \cref{prop:two-interface-obstructions} now gives the displayed test and its reduction to locally closed premises.
			This is the algebraic instance of the classical stability of equivalence relations under regular images in Goursat categories \cite{CarboniEtAl1993};
			see also \cite[Proposition~4.3]{Gran2020Regular}.\qedhere
		\end{enumerate}
	\end{proof}

	\begin{definition}[Separator signatures, congruence completeness, and charts]
		\label{def:congruence-complete-separator}
		\leavevmode
		\begin{enumerate}[label=\textup{(\roman*)},leftmargin=*]
			\item \textbf{Separator signature.}
			A finite native separator signature consists of
			$\mathbf Y\in\mathcal V$ and jointly separating homomorphisms
			$\rho_s:\mathbf Y\to\mathbf D_s$, $s\in S$.
			It is \textbf{congruence-complete} if
			\begin{equation}\label{eq:congruence-complete-signature}
				\gamma_S(\beta)=\beta
				\qquad\bigl(\beta\in\operatorname{Con}(\mathbf Y)\bigr).
			\end{equation}
			
			\item \textbf{Charts over the signature.}
			A chart over this signature has an onto homomorphism
			$p:\mathbf A\to\mathbf Y$ and shared readouts $f_s=\rho_s p$.
			Its private readouts are homomorphisms which, together with the
			shared readouts, separate its actual states.
			All these algebras belong to the same fixed class $\mathcal V$.
		\end{enumerate}
	\end{definition}

	The condition is recognizable from finite separator tables:
	enumerate the equivalence relations compatible with the basic operations, and test the indicated intersections of named kernels.

	\begin{theorem}[Faithful attachment at congruence-complete interfaces]
		\label{thm:congruence-interface-saturation}		
		Suppose every finite algebra in the fixed class $\mathcal V$ is
		three-permutable.
		For a fixed finite separator signature, the following are equivalent.
		\begin{enumerate}[label=(\roman*),nosep]
			\item \emph{Complete signature.}
			The signature is congruence-complete.
			
			\item \emph{Lifting on every chart.}
			Every finite homomorphic chart in $\mathcal V$ over this signature has relative pair lifting at the separator.
			
			\item \emph{Faithfulness of every attachment.}
			Every pair of finite homomorphic charts in $\mathcal V$ sharing exactly this separator, with disjoint private descriptor names, has a state-exact full join whose recovery closure is the	least closure generated by their local closures.
		\end{enumerate}
		If \textup{(i)} fails at $\beta$, two charts with local state algebra
		$\mathbf Y$ and one private quotient readout to $\mathbf Y/\beta$
		each already refute \textup{(iii)}.
		No global source or list of global nonimplications is needed to check \textup{(i)} or the chart-specific criterion \eqref{eq:multidescriptor-congruence-test}.
	\end{theorem}

	\begin{proof}
		\leavevmode
		\begin{itemize}
			\item \emph{\textup{(i)} $\implies$ \textup{(ii)}.}
			Condition \textup{(i)} makes
			\eqref{eq:multidescriptor-congruence-test} automatic, proving
			\textup{(ii)}.
			\item \emph{\textup{(ii)} $\implies$ \textup{(iii)}.}
			Applying \cref{thm:local-pair-native} to the two
			underlying spectra gives \textup{(ii)}$\Rightarrow$\textup{(iii)}.
			\item \emph{\textup{(iii)} $\implies$ \textup{(i)}.}
			Conversely choose $\beta<\gamma_S(\beta)$, and use the charts
			\begin{equation}\label{eq:multidescriptor-quotient-counterexample}
				\begin{array}{c|c|c}
					\text{state algebra}&\text{shared readouts}&\text{private readout}\\ \hline
					\mathbf Y&(\rho_s)_{s\in S}&a(y)=y/\beta\\
					\mathbf Y&(\rho_s)_{s\in S}&b(y)=y/\beta.
				\end{array}
			\end{equation}
			The quotient remains in $\mathcal V$. Put
			$T=\{s:\beta\subseteq\sigma_s\}$. In the first chart
			$C=\{a\}\cup T$ is closed, with agreement congruence $\beta$.
			In the second chart the agreement congruence of $T$ is
			$\gamma_S(\beta)$. A shared descriptor recovered from $T$ already
			belongs to $T$, while $b$ is not recovered because
			$\gamma_S(\beta)\nsubseteq\beta$. Hence $C$ is closed under both
			local theories and $b\notin c_0(C)$. Since the shared readouts
			jointly separate $Y$, matching them forces identical local states.
			The full join is the diagonal copy of $Y$, where $a=b$, creating
			the unintended implication $C\to b$. This contradicts
			\textup{(iii)}, while retaining every local state and every actual
			separator value.\qedhere
		\end{itemize}
	\end{proof}
	
	Stability of congruences under regular images characterizes Goursat categories \cite{HagemannMitschke1973,CarboniEtAl1993,Gran2020Regular}. The counterexample above also uses the partition interpretation of functional dependencies \cite{Lee1983,CosmadakisKanellakisSpyratos1986}. The same quotient charts can test three-permutability in the finite algebras of a variety.
	
	\begin{proposition}[Three-permutability and automatic lifting at complete interfaces]
		\label{prop:finite-goursat-boundary}
		For any fixed equational class $\mathcal V$, the following are equivalent:
		\begin{enumerate}[label=(\roman*),nosep]
			\item \emph{Three-permutability.} Every finite algebra in $\mathcal V$ satisfies
			\eqref{eq:finite-three-permutability}.
			\item \emph{Automatic relative pair lifting.} Every finite homomorphic chart over every finite
			congruence-complete separator signature in $\mathcal V$ has
			relative pair lifting.
		\end{enumerate}
		It suffices in \textup{(ii)} to test signatures consisting of all
		congruence quotient maps of a finite algebra, and charts with two
		additional private readouts: one quotient map and the identity.
	\end{proposition}
	\begin{proof}
		\leavevmode
		\begin{itemize}
			\item \emph{\textup{(i)} $\implies$ \textup{(ii)}.}
			The forward implication is \cref{lem:malcev-port-image}.
			\item \emph{\textup{(ii)} $\implies$ \textup{(i)}.}
			Conversely fix a finite $\mathbf A$ and
			$\alpha,\kappa\in\operatorname{Con}(\mathbf A)$. Give
			$\mathbf Y=\mathbf A/\kappa$ its complete quotient signature.
			On $\mathbf A$ use the shared quotient readouts composed with
			$p:\mathbf A\to\mathbf Y$, the private readout
			$a:\mathbf A\to\mathbf A/\alpha$, and a separately named
			identity readout to ensure state separation. The agreement congruence
			of $\{a\}$ is $\alpha$. By \textup{(ii)} and
			\cref{prop:two-interface-obstructions},
			\[
			\kappa\circ\alpha\circ\kappa=\alpha\vee\kappa.
			\]
			Interchanging $\alpha$ and $\kappa$ proves
			\eqref{eq:finite-three-permutability}. All chart algebras and
			quotients are finite and remain in the same native class.\qedhere
		\end{itemize}
	\end{proof}
	
	\paragraph{Strict enlargement beyond Mal'cev classes.}
	The classical variety of implication algebras, with binary operation
	$\to$ and identities
	\[
	\begin{gathered}
		(x\to y)\to x=x,\qquad
		(x\to y)\to y=(y\to x)\to x,\\
		x\to(y\to z)=y\to(x\to z),
	\end{gathered}
	\]
	is three-permutable \cite{Mitschke1971}; it need not be
	congruence-permutable. For an explicit finite witness, take the
	subalgebra $\{a=(0,1),b=(1,0),1=(1,1)\}$ of the square of
	Boolean implication, with its coordinatewise operation. The kernels of the second and first coordinate projections, respectively, are
	$\alpha=\{\{a,1\},\{b\}\}$ and
	$\kappa=\{\{a\},\{b,1\}\}$.
	The path $a\mathrel\alpha1\mathrel\kappa b$ exists, whereas
	there is no path $a\mathrel\kappa z\mathrel\alpha b$.
	This finite algebra satisfies the three-permutability hypothesis and exhibits two nonpermuting congruences.
	
	\begin{example}[A complete lattice interface with an unrepairable image obstruction]
		\label{ex:complete-port-image-obstruction}
		In the fixed variety of lattices let
		$\mathbf A=\{0<1<2<3\}$ and $\mathbf Y=\{0<1<2\}$, with
		$p(0)=0$, $p(1)=p(2)=1$, $p(3)=2$.
		The convex partitions
		$\kappa=\{\{0\},\{1,2\},\{3\}\}$ and
		$\alpha=\{\{0,1\},\{2,3\}\}$ are congruences of $\mathbf A$,
		and $\kappa=\ker p$. Name on $\mathbf Y$ the two quotient maps
		with kernels
		\[
		\sigma_1=\{\{0,1\},\{2\}\},\qquad
		\sigma_2=\{\{0\},\{1,2\}\}.
		\]
		They form a congruence-complete signature: the congruence lattice of
		the three-element chain consists of $\Delta_Y,\sigma_1,\sigma_2,
		\nabla_Y$, and $\sigma_1\cap\sigma_2=\Delta_Y$.
		Use the private quotient $a:A\to A/\alpha$ and the identity readout.
		For the premise $T=\{a\}$,
		\[
		\beta_T=\Delta_Y\cup\{(0,1),(1,0),(1,2),(2,1)\}
		\ne\nabla_Y=\lambda_T.
		\]
		The pair $(0,2)$ agrees on all shared readouts determined by $T$,
		but has no lift agreeing on $T$.
		By \cref{prop:two-interface-obstructions}, no added homomorphic readouts
		of this same $\mathbf Y$ can correct the pair-lifting failure.
	\end{example}

	\begin{corollary}[The sharp minimum complete separator signature]
		\label{cor:minimal-congruence-interface}
		\leavevmode
		\begin{enumerate}[label=\textup{(\roman*)},leftmargin=*]
			\item \emph{Necessary and sufficient kernels.}
			In any fixed equational class, a finite native separator signature
			is congruence-complete exactly when its named kernels include every
			proper meet-irreducible element	of $\operatorname{Con}(\mathbf Y)$.
			
			\item \emph{Minimum size and equality cases.}
			Its minimum number of named homomorphic descriptors is therefore
			\begin{equation}\label{eq:minimum-congruence-signature}
				\mu(\mathbf Y)=
				\#\{\beta<\nabla_Y:\beta\text{ is meet-irreducible in }
				\operatorname{Con}(\mathbf Y)\}.
			\end{equation}
			The quotient maps $\mathbf Y\to\mathbf Y/\beta$ for these congruences attain the minimum.
			Every minimum signature has exactly	these kernels, up to descriptor relabeling.
			Its actual quotient images are fixed up to the isomorphisms induced from $\mathbf Y$.
		\end{enumerate}
	\end{corollary}

	\begin{proof}
		\leavevmode
		\begin{enumerate}[label=\textup{(\roman*)},leftmargin=*]
			\item \emph{Necessary and sufficient kernels.}
			A proper element $\beta$ is meet-irreducible when
			$\beta=\beta_1\cap\beta_2$ implies $\beta=\beta_1$ or $\beta=\beta_2$.
			Every element of a finite lattice is the meet of the proper meet-irreducibles above it:
			a proper nonirreducible element is a meet of two strictly larger elements,
			and iteration terminates.
			The top is the empty meet.
			Thus all the specified kernels together give \eqref{eq:congruence-complete-signature}.
			Conversely, expressing a proper meet-irreducible as a finite intersection of named kernels containing it forces one of those kernels to equal it.
			Each is	therefore indispensable.
			
			\item \emph{Minimum size and equality cases.}
			The quotient family attains the bound and is jointly separating because its kernels intersect in $\Delta_Y$.
			A surjective map with kernel $\beta$ has actual image isomorphic to
			$\mathbf Y/\beta$, proving the last assertion.
			For a one-element algebra the empty signature suffices and $\mu(\mathbf Y)=0$.\qedhere
		\end{enumerate}
	\end{proof}

	The smaller signature needed for a specified family of charts is computed in \cref{thm:minimum-compositional-interface}.

	\begin{corollary}[Kernel criteria for pair lifting and complete group or module interfaces]
		\label{cor:homomorphic-kernel-saturation}
		For a finite group chart with readouts $\varphi_d$ and actual interface
		$p:G\twoheadrightarrow Y$, put
		\[
		K_T=\bigcap_{d\in T}\ker\varphi_d,
		\qquad N_s=\ker\rho_s,\qquad
		\gamma_S^{\rm n}(N)=\bigcap_{\substack{s\in S\\N\leq N_s}}N_s.
		\]
		\begin{enumerate}[label=\textup{(\roman*)},leftmargin=*]
			\item \emph{Fixed chart.}
			Relative pair lifting is exactly
			\begin{equation}\label{eq:multidescriptor-kernel-test}
				\qquad
				p(K_T)=\bigcap_{s\in c_G(T)\cap S}N_s
				=\gamma_S^{\rm n}\bigl(p(K_T)\bigr)
				\quad(T\subseteq B).
				\qquad
			\end{equation}
			
			\item \emph{Universal signatures.}
			The signature is automatic for all finite group-homomorphic charts
			exactly when its named kernels contain all proper meet-irreducible
			normal subgroups of $Y$. The analogous statement for finite modules
			uses images of kernel submodules and meet-irreducibles of the
			submodule lattice.
		\end{enumerate}
	\end{corollary}
	\begin{proof}
		\leavevmode
		\begin{enumerate}[label=\textup{(\roman*)},leftmargin=*]
			\item \emph{Fixed chart.}
			Group congruences are normal subgroups through the identity class.
			The image of the agreement congruence of $T$ has identity class
			$p(K_T)$: a separator pair $(y,y')$ lifts to a pair agreeing on $T$
			exactly when $y^{-1}y'\in p(K_T)$. Apply \cref{lem:malcev-port-image} to obtain the fixed-chart test.
			\item \emph{Universal signatures.}
			For all charts at once, \cref{thm:congruence-interface-saturation,cor:minimal-congruence-interface} give the meet-irreducible criterion.
			For modules use additive differences and the zero classes of their
			congruences, which are precisely submodules.\qedhere
		\end{enumerate}
	\end{proof}
	
	For a singleton group separator, the only intersections of named kernels are $\{1\}$ and $Y$, recovering \cref{thm:homomorphic-singleton-ports}. With several descriptors, such intersections can realize intermediate normal subgroups.
	
	\begin{corollary}[Projective completeness and minimum size of scalar interface signatures]
		\label{cor:projectively-complete-separators}
		Let $Y=\mathbb F_q^k$ with $k\geq1$, and let the named separator
		readouts be scalar linear forms that jointly separate $Y$.
		\begin{enumerate}[label=\textup{(\roman*)},leftmargin=*]
			\item \emph{Completeness criterion.}
			The interface is automatic for all compatible linear charts exactly
			when its nonzero forms contain a representative of every
			one-dimensional subspace of $Y^*$.
			\item \emph{Minimum size.}
			Its minimum number of scalar
			descriptors is
			\begin{equation}\label{eq:projective-interface-minimum}
				\frac{q^k-1}{q-1}.
			\end{equation}
			Thus projective completeness is necessary for this automatic class,
			as well as sufficient for the visibility condition of
			\cref{def:linear-separator-visibility}.
		\end{enumerate}
	\end{corollary}
	\begin{proof}
		\leavevmode
		\begin{enumerate}[label=\textup{(\roman*)},leftmargin=*]
			\item \emph{Completeness criterion.}
			The proper meet-irreducible subspaces of $Y$ are precisely its
			hyperplanes. A hyperplane is maximal proper, hence meet-irreducible.
			Any subspace of codimension at least two is the intersection of two
			strictly larger subspaces, by choosing independent directions in
			the quotient, so it is not meet-irreducible. Hyperplanes are exactly
			the kernels of nonzero scalar forms, with two such forms having the
			same kernel exactly when they are nonzero scalar multiples. Apply
			the module version of \cref{thm:congruence-interface-saturation,cor:minimal-congruence-interface} to obtain the criterion.
			\item \emph{Minimum size.}
			There are $q^k-1$ nonzero elements in $Y^*$ and $q-1$ representatives on each line, giving \eqref{eq:projective-interface-minimum}.\qedhere
		\end{enumerate}
	\end{proof}
	
	For a concrete failure, use the separator $Y=\mathbb F_q^2$ with
	only the two coordinate readouts $u,v$. The missing quotient $a=u+v$ has kernel $N=\ker(u+v)$.
	Each named separator kernel misses a nonzero element of $N$,
	so $\gamma_S^{\rm n}(N)=Y$. Put this same quotient
	on both sides under the private names $a,b$. The full join preserves
	all states and the exact separator spectrum, but creates $a\to b$.
	This is precisely the obstruction in
	\cref{ex:linear-invisible-separator}. More generally, omitting any
	one hyperplane kernel from a scalar signature leaves the same
	quotient counterexample for that hyperplane, provided the remaining
	shared forms still separate $Y$.
	
	\begin{corollary}[Minimum complete interfaces for products of nonabelian simple groups]
		\label{cor:simple-product-separators}
		Let $Y=\prod_{j=1}^r S_j$, where each $S_j$ is a finite nonabelian
		simple group, and name the $r$ coordinate projections as the
		separator descriptors.
		\begin{enumerate}[label=\textup{(\roman*)},leftmargin=*]
			\item \emph{Complete signature of minimum size.}
			The signature is congruence-complete
			and has the minimum possible number $r$ of named homomorphic
			descriptors.
			\item \emph{Pair lifting.}
			Every finite group-homomorphic chart over this interface
			has relative pair lifting, with no restriction on its private
			descriptor homomorphisms beyond joint state separation.
		\end{enumerate}
	\end{corollary}
	\begin{proof}
		\leavevmode
		\begin{enumerate}[label=\textup{(\roman*)},leftmargin=*]
			\item \emph{Complete signature of minimum size.}
			Every normal subgroup $N$ of $Y$ is a product of a subset of the
			factors. Suppose an element of $N$ has a nonidentity $j$th coordinate.
			Since $S_j$ has trivial center, choose an element supported in
			$S_j$ that fails to commute with it. Their commutator is a
			nonidentity element of $N\cap S_j$. The intersection $N\cap S_j$ is normal in $S_j$, hence
			equals $S_j$. It follows that $N$ contains every factor on which
			it has nontrivial projection, proving the assertion. The normal
			lattice is therefore Boolean, and its proper meet-irreducibles
			are exactly the kernels of the $r$ coordinate projections.
			By \cref{cor:minimal-congruence-interface}, these kernels give a complete signature of minimum size $r$.
			\item \emph{Pair lifting.}
			Pair lifting now follows from \cref{thm:congruence-interface-saturation}.\qedhere
		\end{enumerate}
	\end{proof}
	
	In a finite algebra, the proper meet-irreducible congruences are
	the kernels of the nontrivial subdirectly irreducible quotients
	\cite[Chapter~II, Sections~5,~8, and~12]{BurrisSankappanavar1981} \cite{Birkhoff1944}.
	They give a complete interface signature and identify the quotient
	probes used in the necessity proof.
	
	\subsubsection{Minimum interfaces for a compositional collection of charts}
	\label{subsec:minimum-compositional-interfaces}
	
	Prescribing the permitted attachments can reduce the minimum
	interface cost. Local chart data determine the least cost of
	additional homomorphic readouts that guarantee faithfulness for
	assemblies of arbitrary size over the fixed interface.
	
	Fix a finite algebra $\mathbf Y$ in a three-permutable class with
	the finite presentation law of \cref{def:equational-native-context}.
	Fix a jointly separating initial signature $\Sigma_0$ on $Y$ and a
	finite library of onto homomorphisms
	$p_i:\mathbf A_i\twoheadrightarrow\mathbf Y$ with finite chart
	algebras $\mathbf A_i$. Each library chart has a finite private
	homomorphic readout family $B_i$ which, together with the shared
	signature, separates its states. For $T\subseteq B_i$
	put
	\[
	\beta_{i,T}=(p_i\times p_i)
	\left(\bigcap_{d\in T}\ker f_{i,d}\right),\qquad
	\mathcal B=\{\beta_{i,T}:i,\ T\subseteq B_i\}.
	\]
	These are congruences by \cref{lem:malcev-port-image}.
	The permitted attachment language consists of copies of the library
	charts, the quotient test charts with state algebra $Y$, interface $\id_Y$,
	and one private readout $Y\to Y/\beta$ for $\beta\in\mathcal B$,
	and repeated fiber products over $Y$. Private descriptor names of
	different occurrences are disjoint. A formed component may replace
	its assembly expression, retaining its actual native source, its
	interface, and any selected private readouts. Thus substitution and
	forgetting private readouts are permitted.
	
	Let $\mathcal Q$ be a finite set of candidate quotient readouts
	$Y\to Y/\eta_j$, with nonnegative costs $w_j$. Write
	$\gamma_0$ for the intersection operator of the initial signature.
	Define the finite set of requirements and the subset covered by
	candidate $j$ by
	\begin{align}
		\mathcal O&=\{(\beta,y,y'):
		\beta\in\mathcal B,\ (y,y')\in\gamma_0(\beta)\setminus\beta\},
		\label{eq:interface-cover-obligations}\\
		C_j&=\{(\beta,y,y')\in\mathcal O:
		\beta\subseteq\eta_j,\ (y,y')\notin\eta_j\}.
		\label{eq:interface-cover-items}
	\end{align}
	Every object here is obtained from the local operation and readout
	tables; no assembled source is an input.
	
	\begin{theorem}[Minimum compositional interface enrichment]
		\label{thm:minimum-compositional-interface}
		For a menu selection $J\subseteq\mathcal Q$, the following are equivalent.
		\begin{enumerate}[label=\textup{(\roman*)}]
			\item \emph{Covering the obligations.} The sets $(C_j)_{j\in J}$ cover $\mathcal O$.
			\item \emph{Lifting on components.} Every component in the permitted attachment language has
			relative pair lifting at the enriched interface.
			\item \emph{Faithfulness of assemblies.} Every finite assembly in that language is recovery-faithful
			to the least closure generated by its constituent chart theories
			in the chosen signature and retained descriptor names.
		\end{enumerate}
		\begin{itemize}
			\item \emph{Preservation under assembly.} In these conditions every full assembly also has chart-state
			surjectivity, exact represented information at every assembly-tree
			cut, and a realization in the same fixed presentation law.
			Replacing any subassembly by its formed component preserves these
			conclusions and its retained boundary-readout relation.
			
			\item \emph{Minimum cost and distinguishing attachments.} The exact minimum interface cost is
			\begin{equation}\label{eq:minimum-compositional-interface-cost}
				\tau_{\rm int}=
				\min\left\{\sum_{j\in J}w_j:
				\bigcup_{j\in J}C_j=\mathcal O\right\}.
			\end{equation}
			If no cover exists, the minimum is $+\infty$.
			Every minimum cover is an attaining native signature, and these are
			all the minimum selections. If a selection leaves an obligation
			uncovered, a single library chart and its corresponding quotient
			probe already witness failure of recovery faithfulness. In
			particular every descriptor in an inclusion-minimal feasible
			selection has a distinguishing two-chart attachment after its
			deletion.
		\end{itemize}
	\end{theorem}
	\begin{proof}
		\leavevmode
		\begin{itemize}
			\item \emph{\textup{(i)} $\implies$ \textup{(ii)}.}
			Let $\Sigma_J$ be the enriched signature and $\gamma_J$ its
			intersection operator. For $\beta\in\mathcal B$,
			\[
			\gamma_J(\beta)=
			\gamma_0(\beta)\cap
			\bigcap_{\substack{j\in J\\\beta\subseteq\eta_j}}\eta_j.
			\]
			Consequently the cover condition is exactly
			$\gamma_J(\beta)=\beta$ for every $\beta\in\mathcal B$.
			An agreement premise in a library chart consists of a private
			premise $T$ and some shared readouts, whose kernels on $Y$ have
			intersection $\eta$. Its image relation is exactly
			$\beta_{i,T}\cap\eta$. This follows by lifting each pair in
			$\beta_{i,T}$ and observing that the shared equalities depend only
			on the prescribed interface pair. Under the cover condition both
			$\beta_{i,T}$ and $\eta$ are intersections of named kernels, so
			their intersection is fixed by $\gamma_J$. The local pair-lifting
			criterion applies. The same argument applies to each quotient
			probe, whose private-premise images are $\nabla_Y$ and $\beta$.
			
			For two components over $Y$, an interface pair lifts to states agreeing
			on a prescribed premise in their fiber product exactly when it
			has such lifts in both components: choose the two lifts with the
			same prescribed interface values and combine them. Hence the agreement
			image of the product is the intersection of the two component
			agreement images. This proves by induction that every component
			has pair lifting. Forgetting private readouts merely restricts the
			premises tested, so it preserves this property. This proves
			\textup{(i)}$\Rightarrow$\textup{(ii)} and makes the induction
			valid after substitution.
			\item \emph{\textup{(ii)} $\implies$ \textup{(iii)}.}
			Let $C$ be closed
			under the local theories of the constituent charts in their retained
			language, and let $d\notin C$. Choose a chart carrying $d$. Its trace
			of $C$ is closed, so there is a pair in that chart agreeing on the trace
			and differing at $d$. The pair's values in $Y$ agree on the shared
			descriptors that lie in $C$. At each other chart, lift these same
			ordered interface values by relative pair lifting at the corresponding
			trace of $C$. The resulting two tuples in the common fiber product
			agree on $C$ and differ at $d$. Every local implication holds on the
			product by pullback. Taking $C=c_0(S)$ for each premise $S$ therefore
			proves exact recovery. This argument uses the actual chart states
			and the retained readout equalities, so it applies also when forgetting
			private readouts has left some chart states indistinguishable.
			
			\item \emph{\textup{(iii)} $\implies$ \textup{(i)}.}
			For necessity, suppose $\beta=\beta_{i,T}$ is not fixed by
			$\gamma_J$. Attach to $A_i$ the probe with private readout
			$b(y)=y/\beta$. On the full join the premise $T$ determines $b$.
			Let $C=c_{A_i}^{J}(T)$ and let $S_C=C\cap\Sigma_J$. Precisely the
			shared kernels containing $\beta$ are named in $S_C$, and their
			intersection is $\gamma_J(\beta)>\beta$. In the test chart, $b$ varies on a fiber of the observations in
			$S_C$. Every shared observation determined by $S_C$ already belongs
			to $S_C$, since $\gamma_J(\beta)\subseteq\ker\rho$ implies
			$\beta\subseteq\ker\rho$. Thus $C$ is closed under both local
			theories and excludes $b$. The full join has the extra implication
			$T\to b$, contradicting \textup{(iii)}. This supplies the promised
			finite distinguishing attachment.
			\item \emph{Preservation under assembly.}
			All constituent maps onto $Y$ are surjective, so the matching tuples
			cover every original chart state, and every nonempty cut splits each
			$Y$-fiber into a full rectangle of partial tuples. Those rectangles
			give exact common information. The matching tuples are a nonempty
			subalgebra, and the finite table construction gives the native
			judgment and the homomorphisms to the original charts. A formed
			component can be substituted by the canonical rebracketing of its
			fiber product over $Y$. This map preserves each constituent state
			coordinate, and hence the relation on every retained boundary
			readout. The same pair argument applies with the formed component
			as a chart if its constituent expression is no longer part of the
			retained language.
			\item \emph{Minimum cost and distinguishing attachments.}
			The equivalence gives the minimum formula and all its attaining
			selections. Removing an item from an inclusion-minimal cover
			uncovers some obligation, to which the preceding construction
			applies.\qedhere
		\end{itemize}
	\end{proof}
	
	\begin{corollary}[Sharp library-dependent savings for module interfaces]
		\label{cor:hyperplane-library-interface}
		Let $Y=\mathbb F_q^d$, $d\ge2$, with its whole-state readout as
		the initial signature. Fix $t$ distinct hyperplanes
		$H_1,\ldots,H_t$. Let the library consist of the charts with
		state $Y$, identity interface, and private quotient $Y\to Y/H_i$.
		Allow every linear quotient as an additional readout, at unit cost.
		\begin{enumerate}[label=\textup{(\roman*)},leftmargin=*]
			\item \emph{Specified library.}
			The minimum compositional enrichment has exactly $t$
			descriptors, and every minimum has precisely the kernels
			$H_1,\ldots,H_t$.
			\item \emph{All module attachments.}
			For all homomorphic module attachments, the minimum enrichment of the same initial signature has $(q^d-1)/(q-1)$ additional descriptors. Thus a fixed finite attachment
			library admits a sharp, arbitrarily large saving as $d$ grows.
		\end{enumerate}
	\end{corollary}
	\begin{proof}
		\leavevmode
		\begin{enumerate}[label=\textup{(\roman*)},leftmargin=*]
			\item \emph{Specified library.}
			The library images are $\nabla_Y$ and the $H_i$.
			A quotient kernel containing a hyperplane is either that
			hyperplane or $\nabla_Y$. Consequently fixing $H_i$ requires
			and is achieved by naming $H_i$ itself. The preceding theorem
			gives the minimum, its uniqueness at the kernel level, and
			stability under arbitrary permitted compositions.
			\item \emph{All module attachments.}
			The proper
			meet-irreducibles of the subspace lattice are exactly its
			hyperplanes. Each is the kernel of a nonzero linear form,
			and two such forms have the same kernel exactly when they
			differ by a nonzero scalar. Their number is therefore
			$(q^d-1)/(q-1)$; apply
			\cref{cor:minimal-congruence-interface}.\qedhere
		\end{enumerate}
	\end{proof}
	
	For repeated assembly over one full product boundary,
	\cref{thm:distributed-enrichment} adds an exact feasibility test
	when the quotient menu above is restricted to blocks of ports.
	The fixed-library and universal group-class decompositions are
	\cref{cor:finest-library-blocks,thm:central-prime-blocks}.
	
	\subsection{Faithful native descent through complete interfaces}
	\label{subsec:complete-port-native-descent}
	
	Given local operation tables and interface maps, \cref{thm:complete-port-native-descent} combines the two incidence tests with the matching-algebra construction. It gives an explicit native presentation, maps to the original marked charts, and exact recovery and separator information.
	
	\begin{theorem}[Faithful native descent under the local pair-lifting criterion]
		\label{thm:complete-port-native-descent}
		Fix any equational context of \cref{def:equational-native-context}.
		Let a finite join-tree atlas have local algebras $\mathbf A_i\in\mathcal V$,
		finite descriptor algebras $\mathbf D_d\in\mathcal V$, jointly
		separating homomorphic readouts, and equal actual separator images
		$\mathbf Y_{ij}$.
		Retain the original realized marked local presentations as in \cref{thm:equational-native-amalgamation}.
		Assume each incidence satisfies both conditions in
		\eqref{eq:two-interface-obstructions}.
		If every finite algebra in $\mathcal V$ is three-permutable, this is exactly
		\eqref{eq:multidescriptor-congruence-test}, and is automatic when every separator signature is congruence-complete.
		Let $c_0$ be the least closure generated by the local determination closures and construct
		\[
		W=\{(a_i)\in\prod_i A_i:
		f_{i,d}(a_i)=f_{j,d}(a_j)\text{ on every shared descriptor}\}.
		\]
		Then the following conclusions hold simultaneously.
		\begin{enumerate}[label=(\roman*)]
			\item \emph{Actual-state preservation.}
			Coordinatewise operations make $W$ a nonempty algebra $\mathbf W\in\mathcal V$.
			Every local projection $\mathbf W\to\mathbf A_i$ is an onto homomorphism, retaining the original states and descriptor values.
			
			\item \emph{Recovery and non-recovery preservation.}
			The actual global determination closure is exactly $c_W=c_0$.
			Every nonimplication of $c_0$ has a pair of states in $W$ witnessing it.
			
			\item \emph{Native formation and source recovery.}
			The finite table of $\mathbf W$ gives an actual presentation $P[\mathbf W]$ in the
			fixed context $\Xi_{\mathcal V}^{\mathrm{pres}}$.
			Projection maps are induced by explicit maps into the original local
			presentations, retaining their generators, equations, and marks.
			The marked global endpoint recovers its source congruence by term evaluation, as in \eqref{eq:equational-source-decoder}.
			
			\item \emph{Exact separators and substitution.}
			For every tree cut, let $\mathbf W_L,\mathbf W_R$ be the partial sources and $\mathbf Y$ the actual separator algebra.
			Write $\pi_L:W\to W_L$, $\pi_R:W\to W_R$ for the partial-state projections, and $p:W\to Y$ for the common separator readout.
			Then
			\[
			\mathbf W\cong\mathbf W_L\times_{\mathbf Y}\mathbf W_R,
			\qquad
			\pi_L^*(\mathbb C^{W_L})
			\cap\pi_R^*(\mathbb C^{W_R})
			=p^*(\mathbb C^Y)
			\quad\text{inside }\mathbb C^W.
			\]
			For the projection congruences $\alpha_L,\alpha_R$ and separator congruence $\sigma$,
			\[
			\alpha_L\cap\alpha_R=\Delta_W,\qquad
			\alpha_L\circ\alpha_R=\alpha_R\circ\alpha_L
			=\alpha_L\vee\alpha_R=\sigma.
			\]
			Boundary images have native presentations in the same context.
			Grafting and substitution preserve their operations under the canonical coordinate
			identifications whenever the combined diagram satisfies the same compatibility
			hypotheses.
			\item \emph{Source size.} The construction uses finite local tables and maps.
			The exact number of its states is computed by the local fiber recursion of
			\cref{prop:local-pair-source-cost}.
			For groups and modules this specializes to
			\begin{equation}\label{eq:complete-port-source-count}
				|W|=\frac{\prod_i|A_i|}
				{\prod_{ij\in E(J)}|Y_{ij}|}.
			\end{equation}
		\end{enumerate}
	\end{theorem}
	
	\begin{proof}
		By \cref{prop:two-interface-obstructions}, the two local incidence conditions are equivalent to relative pair lifting. The simplification in the three-permutable case follows from \cref{lem:malcev-port-image,thm:congruence-interface-saturation}.
		\begin{enumerate}[label=\textup{(\roman*)},leftmargin=*]
			\item \emph{Actual-state preservation.}
			The matching-algebra construction in \cref{thm:equational-native-amalgamation} gives the nonempty subalgebra and the surjective chart projections.
			\item \emph{Recovery and non-recovery preservation.}
			Relative pair lifting gives $c_W=c_0$ by \cref{thm:local-pair-native}. By the definition of $c_W$, each $d\notin c_0(T)$ has a witness pair in $W$ agreeing on $T$ and differing on $d$.
			\item \emph{Native formation and source recovery.}
			The finite table and the maps to the retained local presentations are given in \cref{thm:equational-native-amalgamation}. Term evaluation at the global marking recovers the kernel by \eqref{eq:equational-source-decoder}.
			\item \emph{Exact separators and substitution.}
			At a tree cut the matching fiber is the entire nonempty rectangle of its two partial fibers.
			This gives the displayed function-algebra intersection and projection-congruence identities.
			Native boundary images and substitutions are the operation-preserving conclusions
			of \cref{thm:equational-native-amalgamation,cor:operational-native-algebra}.
			\item \emph{Source size.}
			Apply \cref{prop:local-pair-source-cost} to the finite matching sets. For groups and modules all fibers of a surjective homomorphism have equal size, giving \eqref{eq:complete-port-source-count}.\qedhere
		\end{enumerate}
	\end{proof}
	
	The graph construction of \cref{prop:equational-levelii-bridge} realizes
	these tests on one-step formation witnesses. Under the additional Mal'cev
	hypothesis, \cref{thm:native-source-rigidity} determines the minimum for
	the full-cut contract.
	
	The same source $W$ also has the access--commutant realization on $\ell^2(W)$ from \cref{thm:local-pair-native}, with the descriptor algebras and recovery identities in \eqref{eq:local-pair-native-composition}. The presentation endpoint in part~\textup{(iii)} retains the algebraic operations.
	
	\section{Compositional invariants and interface decomposition}
	\label{sec:compositional-invariants}
	To replace a fragment of a package assembly, one must know what a future
	attachment can observe at the exposed interface. Its boundary values form
	a tuple whose coordinates can be correlated. \Cref{thm:arbitrary-graph-cut-information}
	identifies the common information of its two sides with the functions
	on the jointly reachable boundary. For finite group charts,
	\cref{thm:complete-compositional-invariants} determines exactly
	what an attachment can detect through recovery and counting tests.
	\Cref{thm:distributed-enrichment,thm:central-prime-blocks}
	then identify which ports of a common product boundary must be read
	together.
	
	\subsection{Exact common information on arbitrary matching graphs}
	\label{subsec:arbitrary-graph-cut-information}
	
	A partial assignment on either side of a cut extends across it precisely
	when the opposite side can produce the same tuple of crossing values.
	This characterization determines both the common function algebra and
	the assignments omitted by the full matching relation.
	
	\begin{theorem}[Boundary fiber products and exact cut information on finite matching graphs]
		\label{thm:arbitrary-graph-cut-information}
		Let $G=(I,E)$ be any finite graph. At each vertex $i$ let $X_i$
		be a finite nonempty state set. For each edge $e=ij$, specify a
		finite interface set $B_e$ and maps $p_{i,e}:X_i\to B_e$ and
		$p_{j,e}:X_j\to B_e$. Let $W\subseteq\prod_iX_i$ be the
		full relation satisfying all edge equalities, and assume $W\ne\varnothing$.
		For any prescribed bipartition $I=L\sqcup R$ with both sides nonempty,
		let $U_L,U_R$ be the partial matching spaces imposing all edges
		internal to their respective sides. Collect the \emph{complete}
		crossing-edge values into
		\[
		B_{\partial}=\prod_{e\in\delta(L)}B_e,
		\qquad b_L:U_L\to B_{\partial},\quad
		b_R:U_R\to B_{\partial},
		\]
		using the relevant incidence map on each side. Set
		\[
		Y_{\partial}=b_L(U_L)\cap b_R(U_R),\qquad
		U_L^Y=b_L^{-1}(Y_{\partial}),\quad
		U_R^Y=b_R^{-1}(Y_{\partial}).
		\]
		For $y\in Y_{\partial}$ write $U_{L,y}=b_L^{-1}(y)$ and
		$U_{R,y}=b_R^{-1}(y)$. On $W$, let $b_{\partial}$ denote the
		common crossing tuple. Then the following hold.
		\begin{enumerate}[label=(\roman*)]
			\item \emph{Boundary fibers.} The natural restrictions identify
			\begin{equation}\label{eq:arbitrary-graph-cut-fiber}
				W\cong U_L^Y\times_{Y_{\partial}}U_R^Y.
			\end{equation}
			Every $y\in Y_{\partial}$ has a full nonempty rectangular
			fiber $U_{L,y}\times U_{R,y}$.
			\item \emph{Common and generated algebras.} Inside the pointwise function algebra $\mathbb C^W$, let
			$M_L,M_R$ consist of all functions of the respective \emph{full}
			restricted side states. Then
			\begin{equation}\label{eq:arbitrary-graph-cut-functions}
				M_L\cap M_R=\mathbb C[Y_{\partial}],\qquad
				M_L\vee M_R=\mathbb C^W,
			\end{equation}
			where $\mathbb C[Y_{\partial}]$ denotes the pullback along
			$b_{\partial}$ of all complex functions on $Y_{\partial}$.
			\item \emph{Exactness of a named cut observation.} If a preassigned cut observation on $W$ is
			$h=\rho\circ b_{\partial}$ with
			$\rho:Y_{\partial}\to Z$, then
			\begin{equation}\label{eq:arbitrary-graph-named-cut}
				M_L\cap M_R=\mathbb C[h]
				\quad\Longleftrightarrow\quad
				\rho\text{ is injective on }Y_{\partial}.
			\end{equation}
			\item \emph{Coverage of the partial spaces.} The two projections onto the \emph{untrimmed} partial spaces
			$U_L$ and $U_R$ are simultaneously surjective exactly when
			\begin{equation}\label{eq:arbitrary-graph-side-coverage}
				b_L(U_L)=b_R(U_R).
			\end{equation}
		\end{enumerate}
	\end{theorem}
	\begin{proof}
		\leavevmode
		\begin{enumerate}[label=\textup{(\roman*)},leftmargin=*]
			\item \emph{Boundary fibers.} An assignment to all vertices satisfies every edge constraint
			if and only if its restrictions satisfy the internal constraints
			and their complete crossing-interface vectors agree. Conversely any
			pair of such partial assignments with the same boundary vector
			satisfies every edge equality. This gives the fiber product
			\eqref{eq:arbitrary-graph-cut-fiber} and its full rectangles.
			
			\item \emph{Common and generated algebras.} If a function on $W$ depends on both the left and the right
			side separately, it is constant on every rectangle, and hence
			is precisely a function of $y\in Y_{\partial}$. Conversely each
			function of $y$ factors through both sides. For $w\in W$, the product of the indicators of its two
			restricted side states is the indicator of $\{w\}$. These
			point indicators span $\mathbb C^W$.
			
			\item \emph{Exactness of a named cut observation.} The algebra of a finite observation is the entire function
			algebra on its realized image; hence its pullback equals
			$\mathbb C[Y_{\partial}]$ exactly when $\rho$ separates the
			reachable boundary values.
			
			\item \emph{Coverage of the partial spaces.} A left partial assignment
			extends if and only if its boundary vector belongs to
			$b_R(U_R)$, and symmetrically on the right. Both projections
			are onto their untrimmed spaces if and only if their boundary
			images are equal.\qedhere
		\end{enumerate}
	\end{proof}
	
	\begin{corollary}[Native finite-algebra form of the graph-cut identity]
		\label{cor:algebraic-arbitrary-graph-cut}
		In \cref{thm:arbitrary-graph-cut-information}, suppose in addition
		that the local state spaces and interface spaces are finite
		algebras in the fixed equational context of
		\cref{def:equational-native-context}, and that every incidence
		map is a homomorphism.
		\begin{enumerate}[label=\textup{(\roman*)},leftmargin=*]
			\item \emph{Native realization.} If $W$ is nonempty, it is a subalgebra
			of $\prod_i\mathbf X_i$, and its finite table presentation is an
			actual native realization. In particular, for finite groups the identity tuple guarantees $W\ne\varnothing$.
			
			\item \emph{Cut identities.} All the cut identities of
			\cref{thm:arbitrary-graph-cut-information} remain true for
			the underlying finite relation $W$.
		\end{enumerate}
	\end{corollary}
	\begin{proof}
		\leavevmode
		\begin{enumerate}[label=\textup{(\roman*)},leftmargin=*]
			\item \emph{Native realization.} Each matching equality is preserved by every basic operation,
			so the nonempty intersection of the matching subalgebras is
			itself a subalgebra of the product. Invoke
			\cref{lem:finite-table-presentation} for its native realization. In the group case all interface maps preserve identity.
			
			\item \emph{Cut identities.} Apply the preceding set-theoretic theorem to the underlying relation.\qedhere
		\end{enumerate}
	\end{proof}
	
	\Cref{app:cyclic-full-source-test}\textup{(iii)} applies the reachable-boundary identity to a group-native cycle.
	
	\subsection{Complete compositional invariants of group-native fragments}
	\label{subsec:n37-joint-invariants}
	
	Two fragments may have the same image at every individual port and
	interact differently with the same attachment. We retain the image of
	each agreement kernel in the entire boundary group. One additional
	entry per premise records either which private readouts are determined
	at fixed boundary or how many states have that boundary and premise.
	
	\begin{definition}[Marked group fragments: assemblies, continuations, and tests]
		\label{def:fragment}
		\leavevmode
		\begin{enumerate}[label=\textup{(\roman*)},leftmargin=*]
			\item \emph{Marked group fragment.} A \textbf{marked group fragment} consists of a finite group $G$, a finite labeled family of
			boundary homomorphisms $p_a:G\to Y_a$ ($a\in A$), and a finite labeled
			family of retained homomorphisms $f_d:G\to D_d$ ($d\in D$). The target
			groups and labels are fixed in a comparison. Put
			$P=\prod_{a\in A}Y_a$ and $p=(p_a)_{a\in A}$. An empty product is the
			trivial group. The joint boundary image can be a proper subgroup of $P$. Distinct
			states may also agree on all retained readouts, as happens when
			readouts are forgotten.
			
			\item \emph{Assembly.} An allowed assembly of finitely many fragments $G_i$ has a subgroup
			$C\leq\prod_i P_i$ of boundary constraints and source
			\[
			W=\{(g_i)_i\in\prod_i G_i:(p_i(g_i))_i\in C\}.
			\]
			It retains any specified subset of the occurrence-labeled component
			readouts and may introduce homomorphic readouts $h:C\to D_h$ evaluated
			on the joint boundary tuple. Its new boundary map is the composite of
			the joint boundary tuple with a specified homomorphism $\rho:C\to P'$,
			where $P'$ is a product of new interface groups. Disjoint copies of
			private descriptor labels are used before any required identifications.
			
			\item \emph{Continuation.} Allowed continuations are finite expressions in these constructions,
			with one distinguished fragment occurrence of the given type. The
			identity continuation is included. Equality of interface values on an
			arbitrary finite graph is one such subgroup constraint. A subgroup
			constraint can itself be implemented by a native chart with its
			coordinate homomorphisms and the ordinary matching equalities.
			
			\item \emph{Continuation tests and equivalence.} A \emph{recovery test} in a continuation asks whether a named readout
			is determined by a specified family of retained readouts. A
			\emph{numerical test} asks for the order of the assembled group or the
			number of possible joint values of a requested readout family at an
			attained value of an observed family. All tests use only equality and
			cardinalities, not the actual labels of private readout values. A
			continuation can also retain a homomorphic boundary readout as a
			private readout, in particular a boundary coordinate or the joint
			boundary tuple. Two fragments are recovery-equivalent (respectively
			numerically equivalent) when all the corresponding tests in all
			common allowed continuations give equal answers.
		\end{enumerate}
		All groups use the finite group presentation law of
		\cref{def:equational-native-context}. An attachment may use any subgroup of a boundary group as a probe,
		including subgroups whose attachment restricts the chart images.
	\end{definition}
	
	For $T\subseteq D$, write
	\begin{equation}\label{eq:invariants}
		\begin{aligned}
			K_G(T)&=\bigcap_{d\in T}\ker f_d,& K_G(\varnothing)&=G,\\
			B_G(T)&=p(K_G(T))\leq P,\\
			V_G(T)&=\{d\in D:K_G(T)\cap\ker p\leq\ker f_d\},\\
			m_G(T)&=|K_G(T)\cap\ker p|.
		\end{aligned}
	\end{equation}
	The subgroup $B_G(T)$ is normal in the actual joint boundary image
	$B_G(\varnothing)$, not necessarily in the ambient product $P$.
	Define the recovery and numerical records
	\[
	\mathcal I_{\rm rec}(G)=(B_G(T),V_G(T))_{T\subseteq D},\qquad
	\mathcal I_{\rm num}(G)=(B_G(T),m_G(T))_{T\subseteq D}.
	\]
	A native isomorphism fixing the marked maps preserves both records.
	A boundary relabeling transports each subgroup entry by the same
	relabeling. The numerical record determines the
	recovery record, since
	\begin{equation}\label{eq:vertical-cardinality}
		d\in V_G(T)\quad\Longleftrightarrow\quad
		m_G(T)=m_G(T\cup\{d\}).
	\end{equation}
	The assertion follows from inclusion of the corresponding vertical
	kernels and finiteness.
	
	\begin{theorem}[Composition and complete continuation invariants for finite group fragments]
		\label{thm:complete-compositional-invariants}
		Use the fragment and continuation class of \cref{def:fragment}.
		The following statements hold.
		\begin{enumerate}[label=\textup{(\roman*)},leftmargin=*]
			\item \emph{Composition of records.} Both records in \eqref{eq:invariants} are closed under every
			allowed assembly and forgetting of readouts. In an assembly $W$, let
			$T_i$ denote the component readouts occurring in a premise $T$, and
			let $T_*$ denote its new boundary-derived readouts $h$. Set
			\begin{equation}\label{eq:matching-kernel-image}
				M_T=C\cap\prod_i B_{G_i}(T_i)\cap\bigcap_{h\in T_*}\ker h,
				\qquad L_T=M_T\cap\ker\rho.
			\end{equation}
			Then the new boundary image and vertical multiplicity are
			\begin{equation}\label{eq:composition}
				B_W(T)=\rho(M_T),\qquad
				m_W(T)=|L_T|\prod_i m_{G_i}(T_i).
			\end{equation}
			For a retained component descriptor $d$ of component $i$,
			\begin{equation}\label{eq:vertical-composition}
				d\in V_W(T)\quad\Longleftrightarrow\quad
				d\in V_{G_i}(T_i)\ \text{and}\
				\operatorname{pr}_i L_T\leq B_{G_i}(T_i\cup\{d\}).
			\end{equation}
			For a new boundary-derived descriptor $h$, the corresponding test is
			$L_T\leq\ker h$. Hence the qualitative record can be composed without
			computing any multiplicities. Replacing a component by one with the
			same respective record leaves every subsequent test unchanged.
			
			\item \emph{Evaluation of tests.} Source order, determination and conditional range sizes can be
			computed directly. Namely,
			\begin{equation}\label{eq:kernel-size}
				|K_W(T)|=|M_T|\prod_i m_{G_i}(T_i),\qquad
				|W|=|K_W(\varnothing)|.
			\end{equation}
			For every attained observation on $T$, the number of possible joint
			values on $J$ is
			\begin{equation}\label{eq:range}
				R_W(J\mid T)=\frac{|K_W(T)|}{|K_W(T\cup J)|}.
			\end{equation}
			The determination of a component descriptor is tested by
			\eqref{eq:vertical-composition} with $M_T$ in place of $L_T$;
			for a new boundary-derived descriptor the test is $M_T\leq\ker h$.
			
			\item \emph{Completeness.} The records are complete for their respective equivalences:
			\begin{align}
				G\equiv_{\rm rec}G'&\quad\Longleftrightarrow\quad
				\mathcal I_{\rm rec}(G)=\mathcal I_{\rm rec}(G'),
				\label{eq:rec-complete}\\
				G\equiv_{\rm num}G'&\quad\Longleftrightarrow\quad
				\mathcal I_{\rm num}(G)=\mathcal I_{\rm num}(G').
				\label{eq:num-complete}
			\end{align}
			The necessity of the first equality is witnessed by elementary
			native probes only: the trivial subgroup of $P$ and the cyclic
			subgroups of $P$ of nontrivial prime-power order, each supplied with
			its quotient by its unique maximal subgroup. The numerical tests on
			the trivial probe additionally recover every $m_G(T)$.
			\item \emph{Universal summary and native realization.} In particular, the image of each record map is the coarsest exact
			summary for its declared tests. Every other exact summary, restricted to its attained codes, admits a
			unique map onto that image commuting with the summary maps. Every assembled group admits a native table presentation in the fixed group context.
		\end{enumerate}
	\end{theorem}
	
	\begin{proof}
		\leavevmode
		\begin{enumerate}[label=\textup{(\roman*)},leftmargin=*]
			\item \emph{Composition of records.} For any premise, two group elements agree on that premise exactly
			when their quotient $g^{-1}g'$ belongs to its agreement kernel.
			The restriction of the joint boundary map to $K_W(T)$ is an onto
			homomorphism to $M_T$: for each tuple in $M_T$, independently lift its
			$i$th coordinate to an element of $K_{G_i}(T_i)$. These lifts satisfy
			$C$ and every new boundary-readout kernel in $T_*$. Conversely every
			element of $K_W(T)$ gives such a tuple. Its kernel is exactly
			$\prod_i(K_{G_i}(T_i)\cap\ker p_i)$. Thus
			\[
			1\longrightarrow\prod_i(K_{G_i}(T_i)\cap\ker p_i)
			\longrightarrow K_W(T)\longrightarrow M_T\longrightarrow1
			\]
			is exact. Restricting the last surjection to the preimage of $L_T$
			and counting its fibers proves \eqref{eq:composition}; counting
			all fibers gives \eqref{eq:kernel-size}.
			
			For the qualitative test, first suppose
			$d\notin V_{G_i}(T_i)$. A local element in
			$K_{G_i}(T_i)\cap\ker p_i$ with nonidentity $d$-value, combined with
			identity elements at all other components, lies in the vertical
			kernel of $W$ and refutes the desired determination. Suppose instead
			$d\in V_{G_i}(T_i)$. Among elements of $K_{G_i}(T_i)$ with a fixed
			boundary value $b$, the value of $f_d$ is constant: any two differ
			by an element of $K_{G_i}(T_i)\cap\ker p_i$, which is killed by $f_d$.
			This constant is the identity exactly when $b\in
			B_{G_i}(T_i\cup\{d\})$. All tuples in $L_T$ lift, so every vertical
			lift kills $d$ exactly under the inclusion in
			\eqref{eq:vertical-composition}. For a new $h$ the value depends
			only on the tuple in $L_T$. Forgetting readouts simply
			restricts the record indices and the retained target labels. Repeated
			substitution follows by induction on finite formation expressions;
			no acyclicity assumption is made on a matching diagram inside an
			assembly.
			
			\item \emph{Evaluation of tests.} A fiber of the observed homomorphism on $T$ is a left coset of
			$K_W(T)$. Its image under the requested homomorphism on $J$ is a
			translate of the image of $K_W(T)$. The homomorphism theorem gives
			\eqref{eq:range}. The preceding qualitative argument, without
			fixing the new boundary, uses $M_T$ and proves the remaining tests. In particular, numerical equivalence implies recovery
			equivalence, since a determination is precisely range size one.
			
			\item \emph{Completeness.} Equality of the corresponding records gives equality of
			all tests by (i)--(ii), including the identity continuation. Conversely,
			attach the trivial subgroup via the boundary equality; the resulting
			source is $\ker p$. Determination of $d$ by $T$ on that group is exactly
			membership in $V_G(T)$. Now let $E\leq P$ be cyclic of order $\ell^r$
			for a prime $\ell$ and $r\geq1$. Attach $E$ via its inclusion in $P$,
			retaining its quotient readout $z:E\to E/E^\ell$ under a fresh label.
			The matching source is canonically $p^{-1}(E)$ and its agreement
			kernel for $T$ has image $B_G(T)\cap E$ in $E$. Therefore
			\begin{equation}\label{eq:cyclic-probe}
				T\not\to z\quad\Longleftrightarrow\quad
				B_G(T)\cap E\nleq E^\ell\quad\Longleftrightarrow\quad
				E\leq B_G(T).
			\end{equation}
			The last equivalence uses that every proper subgroup of a cyclic
			prime-power group is contained in its unique maximal subgroup.
			Every subgroup of a finite group is generated by its elements of
			prime-power order: each element is the product of the prime-power
			parts in the cyclic group it generates. Consequently the tests
			\eqref{eq:cyclic-probe} recover $B_G(T)$ for every $T$ and prove
			\eqref{eq:rec-complete}.
			
			If all numerical tests agree, the recovery tests agree as just noted.
			In the trivial-probe source $Z=\ker p$, the source order and the
			unconditional range of the readouts $T$ recover
			$m_G(T)=|Z|/R_Z(T\mid\varnothing)$. Thus the numerical record agrees,
			proving \eqref{eq:num-complete}.
			
			\item \emph{Universal summary and native realization.} An exact summary identifies
			only contextually equivalent fragments, so sending its code to the
			record of a representative is well-defined, surjective onto the
			record image, and unique. Each matching source is a subgroup of a
			finite product, hence has the finite table presentation of
			\cref{lem:finite-table-presentation} with the stated homomorphic maps.\qedhere
		\end{enumerate}
	\end{proof}
	
	\begin{corollary}[Minimal cyclic probes and storage bounds for continuation records]
		\label{cor:probe-minimality}
		Let $P$ be a fixed finite group and $\mathscr E$ its set of nontrivial
		cyclic prime-power subgroups.
		\begin{enumerate}[label=\textup{(\roman*)},leftmargin=*]
			\item \emph{Minimal probe family.} For the class of fragments with no
			private readouts and arbitrary boundary map into $P$, the tests
			\eqref{eq:cyclic-probe} with empty premise and $E\in\mathscr E$
			distinguish all possible boundary images. No proper subfamily of
			these tests does so. For $P=(C_\ell)^r$, the sharp number within this
			probe family is $(\ell^r-1)/(\ell-1)$.
			
			\item \emph{Storage bounds.} For fixed $P$ and a descriptor set of size $d$, the qualitative
			record of any fragment can be stored in at most
			$O(2^d((1+\log_2|P|)^2+d))$ bits with a fixed labeling of $P$ and the
			descriptors. This bound is independent of the order of the hidden
			source group. The numerical record additionally requires
			$O(2^d(1+\log_2|G|))$ bits for the multiplicities.
		\end{enumerate}
	\end{corollary}

	\begin{proof}
		\leavevmode
		\begin{enumerate}[label=\textup{(\roman*)},leftmargin=*]
			\item \emph{Minimal probe family.} Sufficiency is \eqref{eq:cyclic-probe}. For necessity, if $E$ is
			omitted, compare the two inclusion fragments $E\hookrightarrow P$
			and $E^\ell\hookrightarrow P$. Every proper subgroup of $E$ lies in
			$E^\ell$, so the only test of the specified family distinguishing
			them is the one for $E$. In an elementary abelian group the
			nontrivial cyclic subgroups are the lines; their count is as stated.
			
			\item \emph{Storage bounds.} Each subgroup of $P$ can be generated by at
			most $\lfloor\log_2|P|\rfloor$ elements: adjoining an element outside
			the group already generated at least doubles its order. Each generator
			needs at most $\lceil\log_2|P|\rceil$ bits, with a length marker or
			fixed padding, while a set of $d$ target labels needs $d$ bits.
			There are $2^d$ premise indices. A multiplicity is an integer between
			one and $|G|$.\qedhere
		\end{enumerate}
	\end{proof}

	\begin{example}[All proper marginals and recovery closure can agree]
		\label{ex:marginal-failure}
		For an odd prime $\ell$ and $r\geq2$, use additive groups and take
		$P=\mathbb F_\ell^r$.
		Let
		\[
		H=\{x:\textstyle\sum_{j=1}^r x_j=0\},\qquad
		H'=\{x:\textstyle\sum_{j<r}x_j+2x_r=0\}.
		\]
		Use the inclusions as joint boundary maps, with the individual
		coordinates also retained as named readouts. Both sources have order
		$\ell^{r-1}$, all proper coordinate projections are full, and both
		have the same whole-coordinate recovery closure: any $r-1$ coordinates
		determine the last, and fewer do not determine any omitted coordinate.
		The coefficient vectors are nonproportional. Attaching $H$ itself by
		joint boundary equality produces $H$ from the first source and
		$H\cap H'$ from the second, of orders $\ell^{r-1}$ and $\ell^{r-2}$.
		The partner chart is covered in the first assembly but not in the
		second. Thus even the entire family of proper marginal images,
		order, and recovery closure fails to determine the joint record.
		The attached chart projects onto each individual port.
	\end{example}

	\begin{remark}[Comparison of the two equivalences]
		\label{rem:scope}
		With trivial boundary and no retained readouts, all finite groups
		have the same recovery record, while their numerical record gives
		their order. In particular, $C_4$ and $C_2\times C_2$ have the same
		numerical record. Related constructions that compose fragments through
		interfaces appear in graphical realizations of codes
		\cite{Forney2001,Forney2014,ForneyGluesingLuerssen2013}.
	\end{remark}

	\begin{corollary}[Faithful-source existence on finite group matching graphs]
		\label{cor:arbitrary-graph-faithfulness}
		Take an arbitrary finite graph of finite group charts with homomorphic
		edge interfaces and retained homomorphic descriptor readouts, and use
		its full matching group $W$. Identify shared descriptor labels only
		when the wiring enforces their equality. Fix the target $c_0$ generated
		by the local recovery laws. At each protected cut take the complete
		joint crossing-edge boundary with its actual reachable image in $W$,
		as in \cref{thm:arbitrary-graph-cut-information}. Require coverage of
		all original chart states, the closure $c_0$, exact common information
		on those complete reachable boundaries, and native finite group
		realization.
		\begin{enumerate}[label=\textup{(\roman*)},leftmargin=*]
			\item \emph{Coverage and recovery tests.} There exists a compatible native group source for this
			contract if and only if
			\begin{equation}\label{eq:graph-obstruction}
				\operatorname{pr}_i M_{\varnothing}=B_{G_i}(\varnothing)
				\quad\text{for every }i,\qquad c_W=c_0.
			\end{equation}
			Here $C$ consists of the edge equalities; both tests are determined
			by the local recovery records of \cref{thm:complete-compositional-invariants},
			with the incidence maps retained. The latter equality tests every
			nonimplication of $c_0$ by the qualitative composition rule.
			
			\item \emph{Realization and obstruction.} In the
			positive case $W$ is a native solution, whose order is given by
			\eqref{eq:kernel-size}. If either condition fails, no compatible
			source can satisfy the coverage and recovery clauses, even if the
			source is not required to be a group.
		\end{enumerate}
	\end{corollary}
	\begin{proof}
		\leavevmode
		\begin{enumerate}[label=\textup{(\roman*)},leftmargin=*]
			\item \emph{Coverage and recovery tests.} A state $g_i\in G_i$ extends to $W$ exactly when its entire boundary
			tuple belongs to $\operatorname{pr}_i M_{\varnothing}$. Indeed choose
			an extending boundary tuple in $M_{\varnothing}$ and lift the other
			coordinates independently. This proves the coverage test. Every
			compatible source maps into $W$ by its chart tuple, hence cannot
			restore a chart state missing from $W$. Every local implication
			holds in $W$ and every dependency of $W$ holds on each compatible
			source by pullback, so $c_0\leq c_W\leq c_H$ for any compatible
			source $H$. This proves necessity of recovery equality.
			
			\item \emph{Realization and obstruction.} In the
			positive case $W$ is a finite group, and its exact complete-cut
			common information is \cref{thm:arbitrary-graph-cut-information}.
			Its finite table presentation supplies the native solution. If either test fails, the obstruction established in (i) applies to every compatible source.
			The tests and the count are those of the preceding theorem. A
			shared readout can be evaluated using any one of its equal
			occurrences; local and global closure tests use these same labels.\qedhere
		\end{enumerate}
	\end{proof}
	
	To check chart coverage after hiding a chart's internal structure, retain
	its projection as a named descriptor. The numerical record gives the
	size of the projection image, which equals the chart order exactly
	when the projection is onto.
	
	\subsection{Exact distributed interface obstructions and canonical blocks}
	\label{subsec:n37-distributed-interfaces}
	
	For a common boundary $Y=\prod_{a\in A}Y_a$, with $A$ finite and
	nonempty, suppose that each shared readout can access only one block
	of ports. An agreement image that couples two blocks can prevent
	global recovery implications from being derived locally.
	\Cref{thm:distributed-enrichment} detects this failure by a
	subgroup index.
	Fix a finite library of group charts $p_i:G_i\twoheadrightarrow Y$
	with private homomorphic readouts and a jointly separating initial
	shared signature $\Sigma_0$ of homomorphisms on $Y$. Let
	\[
	\mathcal B=\{B_{G_i}(T): i\text{ in the library},\ T\subseteq D_i\}
	\subseteq\{N:N\trianglelefteq Y\}.
	\]
	Here $D_i$ denotes the private descriptor labels. Permit repeated
	fiber products over the entire joint $Y$, forgetting private readouts,
	and the quotient probes on $Y$ with one private quotient $Y\to Y/N$
	for each $N\in\mathcal B$, as in
	\cref{thm:minimum-compositional-interface}. For a partition $\Pi$ of
	$A$ write $Y_B=\prod_{a\in B}Y_a$ for a block $B$.

	\begin{theorem}[Block-local interface enrichment: factorization, obstructions, and minimum cost]
		\label{thm:distributed-enrichment}
		Let $\Pi$ be a partition of the boundary port labels.
		Suppose the initial shared descriptors factor through individual blocks. Allow as new descriptors all homomorphic quotient readouts of the blocks,
		pulled back to $Y$.
		Then the following statements hold.
		\begin{enumerate}[label=\textup{(\roman*)},leftmargin=*]
			\item \emph{Equivalence.}
			The following conditions are equivalent.
			\begin{enumerate}[label=\textup{(\alph*)},leftmargin=*]
				\item \emph{Relative pair lifting.}
				Some such enrichment gives relative pair lifting for every component in the permitted common-joint-boundary language.
				
				\item \emph{Faithful recovery.}
				Some such enrichment gives recovery-faithfulness for every
				finite assembly in that language, relative to its enriched vocabulary.
				
				\item \emph{Block factorization.}
				Every $N\in\mathcal B$ has the block factorization
				\begin{equation}\label{eq:block-factorization}
					N=\prod_{B\in\Pi}\operatorname{pr}_B N.
				\end{equation}
			\end{enumerate}
			
			\item \emph{Index obstruction.}
			The obstructions are measured by the intrinsic positive integers
			\begin{equation}\label{eq:block-defect}
				\Delta_\Pi(N)=\left[\prod_{B\in\Pi}\operatorname{pr}_B N:N\right].
			\end{equation}
			If any such index exceeds one, no number of block-local homomorphic
			readouts can repair the corresponding defect; one library chart and
			its quotient probe already refute universal faithfulness.
			
			\item \emph{Minimum cost.}
			For a restricted finite menu of block-local quotient kernels
			$L_j\trianglelefteq Y$ with costs $w_j\geq0$, let $\gamma_0(N)$ be the intersection of the kernels of initial readouts containing $N$.
			The exact additional cost is the weighted covering minimum for the obligations
			\[
			\mathcal O=\{(N,y):N\in\mathcal B,\ y\in\gamma_0(N)\setminus N\},
			\]
			where candidate $j$ covers $(N,y)$ precisely when $N\leq L_j$ and
			$y\notin L_j$.
			If no cover exists the cost is $+\infty$.
			When the factorization condition holds and every block quotient is available, a cover always exists.
		\end{enumerate}
	\end{theorem}

	\begin{proof}
		\leavevmode
		\begin{enumerate}[label=\textup{(\roman*)},leftmargin=*]
			\item \emph{Equivalence.}
			For any signature $\Sigma$, put
			$\gamma_\Sigma(N)=\bigcap_{\rho\in\Sigma,\,N\leq\ker\rho}\ker\rho$.
			If every readout factors through a block, its kernel contains $N$
			only if it contains the product of all block projections of $N$.
			Indeed a kernel has the form $\operatorname{pr}_B^{-1}(L)$ for some
			$L\trianglelefteq Y_B$, and $N\leq\operatorname{pr}_B^{-1}(L)$
			means $\operatorname{pr}_B N\leq L$. The full menu attains equality:
			use the quotients of each $Y_B$ by $\operatorname{pr}_B N$, which
			are normal subgroups because $N\trianglelefteq Y$. Thus the intersection
			of every block-local kernel containing $N$ is exactly
			$\prod_B\operatorname{pr}_B N$. All $N\in\mathcal B$ can therefore
			be fixed by the enriched $\gamma$ exactly under
			\eqref{eq:block-factorization}.
			
			\smallskip
			\noindent\emph{\textup{(c)} $\implies$ \textup{(a)} and \textup{(b)}.}
			Assume \textup{(c)} and choose all block quotient readouts.
			The kernel calculation above gives $\gamma_\Sigma(N)=N$ for every
			$N\in\mathcal B$. The
			agreement image of a library premise, including any shared readouts,
			is $N\cap L$, where $L$ is the intersection of those shared kernels.
			It is again an intersection of named kernels. The quotient probes
			have agreement images $Y$ or $N$ for their private premises.
			Fiber-product composition over the entire $Y$ intersects the
			agreement images, by \cref{thm:complete-compositional-invariants};
			forgetting readouts only removes tested premises. The group version
			of relative pair lifting in \cref{cor:homomorphic-kernel-saturation}
			then holds throughout the language. The witness argument of
			\cref{thm:minimum-compositional-interface} supplies faithfulness of
			every finite assembly, as well as coverage, exact common-joint cut
			information, and native realization.
			
			\smallskip
			\noindent\emph{\textup{(a)} $\implies$ \textup{(c)}.}
			Suppose an enriched signature $\Sigma$ gives relative pair lifting.
			By \cref{cor:homomorphic-kernel-saturation}, every library agreement
			image $N\in\mathcal B$ satisfies $\gamma_\Sigma(N)=N$.
			The kernel calculation gives
			\[
			N\leq\prod_{B\in\Pi}\operatorname{pr}_B N
			\leq\gamma_\Sigma(N)=N,
			\]
			which proves \textup{(c)}.
			
			\smallskip
			\noindent\emph{\textup{(b)} $\implies$ \textup{(c)}.}
			Suppose \textup{(c)} fails, and fix any permitted enriched signature
			$\Sigma$. For some $N=B_{G_i}(T)\in\mathcal B$, the kernel calculation
			gives $N<\gamma_\Sigma(N)$. Attach the quotient
			probe $Y\to Y/N$ under a fresh private name $z$. The full join
			satisfies $T\to z$. Let $C_i$ be the closure of $T$ in the first
			chart. Its shared names are exactly the readouts whose kernels contain
			$N$. Their kernels intersect in $\gamma_\Sigma(N)$. In the probe these
			names determine neither $z$ nor any additional shared name: a
			shared kernel containing $\gamma_\Sigma(N)$ already contains $N$.
			Thus $C_i$ is closed under both local theories and omits $z$.
			The joined implication is extra. Hence this enrichment fails \textup{(b)}.
			Since $\Sigma$ was arbitrary, this proves the contrapositive.
			
			\item \emph{Index obstruction.}
			Since $N\leq\prod_{B\in\Pi}\operatorname{pr}_B N$ and the groups
			are finite, $\Delta_\Pi(N)$ is a positive integer. It equals one
			exactly when \eqref{eq:block-factorization} holds for $N$.
			If it exceeds one, the kernel calculation in \textup{(i)} shows
			that $N<\gamma_\Sigma(N)$ for every block-local enrichment,
			including the full quotient signature. The two-chart construction
			in \textup{(i)} then supplies the claimed obstruction using the
			library chart realizing $N$ and its quotient probe.
			
			\item \emph{Minimum cost.} For a finite menu the enriched intersection is
			$\gamma_0(N)\cap\bigcap_{j:N\leq L_j,\ j\text{ selected}}L_j$.
			Equality with $N$ is precisely coverage of every obligation.
			The kernel criterion established in \textup{(i)} yields the exact cost and all minimizing
			selections. This is the group-difference form of the covering
			criterion of \cref{thm:minimum-compositional-interface}, now with
			an exact structural test for feasibility under block-locality.\qedhere
		\end{enumerate}
	\end{proof}

	\begin{corollary}[Finest port grouping for a prescribed attachment library]
		\label{cor:finest-library-blocks}
		For any finite family $\mathcal B$ of normal subgroups of a direct
		product $Y$, the following hold.
		\begin{enumerate}[label=\textup{(\roman*)},leftmargin=*]
			\item \emph{Finest admissible partition.} There is a unique finest partition $\Pi_{\mathcal B}$ of
			its port labels at which all the indices in
			\eqref{eq:block-defect} equal one. A partition has this property
			if and only if it coarsens $\Pi_{\mathcal B}$.
			
			\item \emph{Construction by bipartitions.} The finest partition
			can be constructed as the common refinement of all bipartitions at
			which every $N\in\mathcal B$ factors; if there are no such nontrivial
			bipartitions it has one block. With an initial signature consisting
			of the individual port identities, this gives the finest grouping
			permitting universally faithful block-local enrichment for the
			specified library and quotient probes.
		\end{enumerate}
	\end{corollary}
	\begin{proof}
		\leavevmode
		\begin{enumerate}[label=\textup{(\roman*)},leftmargin=*]
			\item \emph{Finest admissible partition.} For a set of labels $B$, let $e_B:Y\to Y$ retain the coordinates
			in $B$ and replace the others by the identity. A subgroup $N$
			factors over $\Pi$ exactly when $e_B(N)\subseteq N$ for every block
			$B$: under these inclusions each combination of projected tuples
			is the product of elements of $N$ on disjoint supports. Conversely
			a factorization gives those inclusions. The maps satisfy
			$e_B e_{B'}=e_{B\cap B'}$. Hence the family of admissible partitions
			is closed under common refinement and under coarsening. It is
			nonempty since the one-block partition is admissible. Its common
			refinement is therefore the unique finest one.
			
			\item \emph{Construction by bipartitions.} Each admissible
			partition is the common refinement of its block-versus-complement
			bipartitions (with the one-block case understood as an empty
			refinement). This proves the construction and the final statement
			by \cref{thm:distributed-enrichment}.\qedhere
		\end{enumerate}
	\end{proof}
	
	\begin{theorem}[Universal interaction blocks from central chief factors]
		\label{thm:central-prime-blocks}
		For a finite group $Y_a$ define
		\[
		\mathsf p_c(Y_a)=\{\ell:\ell\text{ prime},\
		\ell\mid |Z(Y_a/L)|\text{ for some }L\trianglelefteq Y_a\}.
		\]
		Equivalently, these are the orders of its central chief factors.
		On the port labels join distinct $a,b$ by an edge when
		$\mathsf p_c(Y_a)\cap\mathsf p_c(Y_b)\ne\varnothing$.
		Let $\Pi_c$ be the connected-component partition of this graph.
		The following are equivalent for any partition $\Pi$:
		\begin{enumerate}[label=\textup{(\roman*)},leftmargin=*]
			\item \emph{Normal subgroups.} Every normal subgroup of $Y=\prod_aY_a$ factors over $\Pi$.
			\item \emph{Interaction blocks.} $\Pi$ coarsens $\Pi_c$.
			\item \emph{Universal enrichment.} Starting with the individual port identities, a finite
			signature of block-local homomorphic quotient readouts suffices for
			relative pair lifting of every finite group chart onto the joint
			boundary $Y$, and for recovery-faithfulness of all their finite
			common-joint-boundary assemblies.
		\end{enumerate}
		If a partition separates the endpoints of an edge corresponding to
		a prime $\ell$, there is a normal subgroup $N\trianglelefteq Y$ with
		$\Delta_\Pi(N)=\ell$. Two charts with state group $Y$, identity joint
		boundary, and private quotient readouts to $Y/N$ then refute universal
		faithfulness for every block-local signature. Thus $\Pi_c$ is the
		unique finest universally sufficient grouping. In particular,
		individual ports suffice exactly when their central-prime sets are
		pairwise disjoint; pairwise coprime group orders are sufficient but
		not necessary.
	\end{theorem}
	\begin{proof}
		By Cauchy's theorem, divisibility of $|Z(Y_a/L)|$ by $\ell$ gives a
		central subgroup of order $\ell$ in $Y_a/L$. Its preimage and $L$
		are consecutive normal subgroups and can occur in a chief series.
		Conversely a central chief factor is a minimal nontrivial normal
		subgroup in a quotient and is central there, hence is cyclic of
		prime order. This proves the stated equivalent description.
		\begin{itemize}
			\item \emph{\textup{(ii)} $\implies$ \textup{(i)}.} Fix $N\trianglelefteq Y$. Identify the factors with their canonical
			coordinate subgroups of $Y$, and set $N_a=N\cap Y_a$ and
			$Q_a=\operatorname{pr}_a N/N_a$. For $x\in N$ and $g\in Y_a$ the
			commutator of $x$ with the element supported at $a$ having value $g$
			belongs to $N$ and is supported only at $a$. Thus $Q_a$ embeds
			into $Z(Y_a/N_a)$. In particular all prime divisors of $|Q_a|$
			belong to $\mathsf p_c(Y_a)$. The quotient
			$N/(\prod_a N_a)$ is a subgroup of $\prod_a Q_a$, and $N$ is
			its full preimage in $\prod_a\operatorname{pr}_a N$.
			
			If $\Pi$ coarsens $\Pi_c$, no prime occurs in the central-prime
			sets of ports in different blocks. Consequently the abelian groups
			$\prod_{a\in B}Q_a$ have pairwise coprime orders as $B$ runs over
			$\Pi$. Every subgroup of a product of groups of pairwise coprime
			orders is the product of its block projections: for an element of
			the subgroup, take an integer power which is one modulo the order
			of the chosen block group and zero modulo all other block orders.
			It isolates that block coordinate. Such powers exist by the Chinese
			remainder theorem, and the isolated coordinates generate the
			projected product. Applying this to the displayed quotient and
			taking full preimages proves (i).
			
			\item \emph{\textup{(i)} $\implies$ \textup{(ii)}.} If (ii) fails, some edge has endpoints $a,b$ in different blocks.
			Choose normal subgroups $L_a,L_b$ and cyclic central subgroups
			$E_a\leq Y_a/L_a$, $E_b\leq Y_b/L_b$, both of order $\ell$ for
			a prime common to their sets. The graph of any isomorphism
			$E_a\to E_b$ is central in the product quotient. Let $N$ be its
			full preimage in $Y_a\times Y_b$, multiplied by the other coordinate
			groups. This is normal in $Y$. Its projected product over $\Pi$
			replaces the diagonal subgroup of order $\ell$ by the product of
			order $\ell^2$, with the same quotient kernels. The index is $\ell$.
			Thus (i) fails, proving the equivalence of (i) and (ii).
			
			\item \emph{\textup{(i)} $\Longleftrightarrow$ \textup{(iii)}.} Every normal $N\trianglelefteq Y$ occurs as an agreement image in
			the group chart $Y$ with identity joint boundary and private quotient
			$Y\to Y/N$. Thus if (iii) holds, the necessity argument of
			\cref{thm:distributed-enrichment} forces (i) for every $N$.
			Conversely if (i) holds, the finite signature of all quotient readouts
			on the blocks fixes every normal subgroup of $Y$. The group
			congruence-completeness criterion
			\cref{thm:congruence-interface-saturation} gives relative pair
			lifting, and the finite common-boundary composition argument gives
			faithfulness for all finite assemblies. In the negative case the
			two quotient-chart proof of \cref{thm:distributed-enrichment}
			applies to the normal subgroup of index defect $\ell$. This proves
			all claims.\qedhere
		\end{itemize}
	\end{proof}
	
	\begin{remark}[Classical group-theoretic input and the structural conclusion]
		The analysis of normal subgroups of direct products is classical
		Goursat theory; the role of common central chief factors also appears
		in criteria for distributive normal subgroup lattices
		\cite{BauerSenZvengrowski2015,DeMedtsTarnauceanu2011}.
		Nonabelian simple port groups have empty central-prime sets
		and thus contribute isolated vertices, even when their orders have
		common prime divisors. Two $S_3$ ports, on the other hand, share the prime
		$2$ through their sign quotients; the subgroup of pairs with equal
		sign gives index defect two. 
	\end{remark}
	
	\section{Recovery under cuts over a common boundary}
	\label{sec:recovery}
	
	A protected cut allows the two sides of a source to recombine freely
	at each fixed boundary value. This operation can destroy a
	determination between charts on opposite sides. \Cref{thm:determinants}
	constructs the least recombinant enlargement of a source and gives an
	exact test for preserving every whole-chart determination. The test
	uses the minimal premises for each chart. For a matroid closure it
	reduces to the circuit criterion in \cref{cor:circuit}.

	\subsection{Whole-chart determination and cut atoms}

	Let $I$ be a finite nonempty set of chart labels.
	For each $i\in I$ fix a surjection $p_i:X_i\to Q$ of nonempty sets.
	A chart state is an element of $X_i$ and its boundary value is its image in $Q$.
	The full matching relation and its boundary map are
	\[
	W=\{x\in\prod_{i\in I}X_i:p_i(x_i)=p_j(x_j)\text{ for every }i,j\},\qquad p(x)=p_i(x_i).
	\]
	A source is a nonempty relation $H\subseteq W$.
	Write $H_q=H\cap p^{-1}(q)$ and let $\operatorname{pr}_S$ retain the coordinates indexed by $S\subseteq I$.
	A source is \textbf{chart-surjective} if $\operatorname{pr}_{\{i\}}H=X_i$ for every $i$.
	Fix one such source $H_0$.

	A cut is a bipartition $L\mid R$ of $I$ into nonempty sides.
	A source is \textbf{exactly recombinant} at the cut if, in every attained boundary fiber,
	$
	H_q=\operatorname{pr}_{L}H_q\times\operatorname{pr}_{R}H_q
	$
	under the coordinate identification.
	For a family of protected cuts, its atoms are the classes of labels lying on the same side of every cut.
	Denote the resulting partition by $\Pi$.
	For the empty cut family, $\Pi=\{I\}$.

	For a source $H\subseteq W$, let $c_H(S)$ consist of the labels $i$ for
	which agreement on $Q$ and on all charts of $S$ forces agreement on
	chart $i$. This is the recovery closure conditional on the common
	boundary. Every nonempty set of whole charts already reads $Q$.
	For chart-surjective $H$, the set $c_H(\varnothing)$ consists
	of the $i$ for which $X_i\to Q$ is injective. Indeed, chart coverage
	puts every pair of possible $i$-values with the same boundary into
	that boundary fiber of $H$. Thus equality of these conditional closures is
	equivalent to equality of the ordinary whole-chart closures for
	chart-surjective sources with these same charts. Indeed ordinary
	empty-premise determination means the chart itself is a singleton.

	Call $P\subseteq I\setminus\{i\}$ a minimal premise for $i$ if
	$i\in c_{H_0}(P)$ and no proper subset of $P$ has this property.
	The map $c_H$ is a closure operator. Extensivity and monotonicity
	follow from agreement on the specified charts. Two states in one fiber
	that agree on $S$ agree on every element of $c_H(S)$, which gives
	idempotence. The implications $P\to i$ just defined form the canonical
	direct unit basis of the conditional closure \cite{BertetMonjardet2010}. Form a hypergraph on $I$ with an edge $P\cup\{i\}$ for each such pair.
	Include isolated vertices when taking its connected components.
	For an atom partition $\Pi$, define $R_\Pi$ fiberwise by
	\[
	(R_\Pi)_q=\prod_{B\in\Pi}\operatorname{pr}_B (H_0)_q,
	\qquad (H_0)_q=\{x\in H_0:x\text{ has boundary }q\}.
	\]

	\begin{theorem}[Minimal premises and simultaneous exact recombination]
		\label{thm:determinants}
		For a cut family with atoms $\Pi$, the following are equivalent.
		\begin{enumerate}[label=\textup{(\roman*)},leftmargin=*]
			\item \emph{Existence.}
			Some chart-surjective source is exactly recombinant at every protected cut and has the same whole-chart recovery closure as $H_0$.
			
			\item \emph{Recovery by the least enlargement.}
			$R_\Pi$ has that closure.
			
			\item \emph{Minimal premises.}
			Each connected component of the minimal-premise hypergraph lies in one atom of $\Pi$.
		\end{enumerate}
		For every $i\in B\in\Pi$ the closure of the least recombinant
		enlargement is given by
		\begin{equation}\label{eq:localpremise}
			i\in c_{R_\Pi}(S)\quad\Longleftrightarrow\quad
			i\in c_{H_0}(S\cap B).
		\end{equation}
		In particular, the components give the finest partition preserving
		recovery, and separate feasibility of each protected cut implies
		joint feasibility. These statements also hold for subgroup sources
		when the charts and boundary maps are groups and homomorphisms
		and $H_0$ is a subgroup.
		The chart groups and their kernels may be infinite.
	\end{theorem}

	\begin{proof}
		Within one boundary fiber, recombination permits splicing two tuples
		along either side of a cut. Starting with two tuples and splicing
		successively with the second tuple keeps the first tuple on the
		intersection of the chosen sides. This realizes splicing on a single
		atom. Repeating over atoms proves that an exact source has product
		fibers over the atoms, and that $R_\Pi$ is the least such enlargement
		of $H_0$. Its chart images and all atom images agree with those of
		$H_0$. In the group setting it is the intersection, in the matching
		group, of the inverse images of the subgroups
		$\operatorname{pr}_B H_0$ under the atom projections. Hence it is
		itself a subgroup.
		
		For any source with product fibers, determination of $i\in B$ by
		$S$ is equivalent, conditional on $Q$, to determination by $S\cap B$.
		Indeed two $B$-tuples in one fiber that agree on $S\cap B$ and differ
		on $i$ can be completed using identical tuples on every other atom.
		This produces a witness to failure at $S$. The converse follows by
		inclusion of premises. Since the $B$-image of $R_\Pi$ is the seed's
		$B$-image, this proves \eqref{eq:localpremise}.
		\begin{itemize}
			\item \emph{\textup{(i)} $\implies$ \textup{(iii)}.} Apply the product-fiber observation to a source as in (i).
			A minimal premise for $i$ must be contained in its atom, since its
			intersection with the atom still determines $i$. Thus no hyperedge
			meets two atoms, proving (iii).
			
			\item \emph{\textup{(iii)} $\implies$ \textup{(ii)}.} Assume (iii) and take
			$i\in c_{H_0}(S)$ with $i\notin S$. Finiteness supplies a minimal
			premise $P\subseteq S$ for $i$. Then $P\subseteq S\cap B$, so
			\eqref{eq:localpremise} gives $i\in c_{R_\Pi}(S)$.
			The reverse inclusion follows from $H_0\subseteq R_\Pi$.
			The case $i\in S$ follows from extensivity.
			
			\item \emph{\textup{(ii)} $\implies$ \textup{(i)}.} The source $R_\Pi$ witnesses (i).
		\end{itemize}
		Only the label set $I$ needed to be finite:
		the splicings, the partition, and the search for a minimal premise
		all take place on these labels. Finally, an edge
		lies in a joint atom precisely when no individual cut separates
		its vertices. This proves the separate-to-joint assertion.
	\end{proof}
	
	\subsection{Circuits and the information retained by whole charts}
	In a matroid, each element of a circuit is determined by the other elements of that circuit.
	\begin{corollary}[Circuit criterion for preservation of whole-chart recovery]
		\label{cor:circuit}
		Suppose $c_{H_0}$ is the closure of a matroid $M$.
		A chart-surjective source can preserve this closure and recombine
		exactly at all protected cuts if and only if no circuit of $M$
		crosses a protected cut.
	\end{corollary}
	\begin{proof}
		For a circuit $C$ and $i\in C$, the set $C\setminus\{i\}$ is a
		minimal premise for $i$. This includes a loop, whose only minimal
		premise is empty and whose singleton circuit crosses no cut.
		Conversely, a minimal premise $P$ for
		$i\notin P$ is independent, and $P\cup\{i\}$ is minimally dependent.
		Indeed, a dependent $P$ has an element whose removal preserves its
		closure, contrary to minimality. Every proper subset of $P$ is
		independent and fails to span $i$, so adjoining $i$ to it gives
		an independent set. The hyperedges of \cref{thm:determinants} are therefore exactly
		the circuits.
	\end{proof}
	
	At the level of closures the atom condition recovers the direct-product
	identity $c(S)=\bigcup_{B\in\Pi}(c(S\cap B)\cap B)$
	\cite{DemetrovicsFurediKatona1985,Libkin1993}. For a matroid $M$ it is
	the direct-sum criterion $M=\bigoplus_{B\in\Pi}M|B$
	\cite[Corollary~4.2.9 and identity~4.2.12]{Oxley2011}.
	In relational language, let $b\notin I$ name the boundary coordinate
	with values in $Q$. Recombination at $L\mid R$ is exactly the
	multivalued dependency $b\twoheadrightarrow L$: the augmented source
	is the natural join of its projections on $\{b\}\cup L$ and $\{b\}\cup R$
	\cite[Theorem~1]{Fagin1977}. Functional and multivalued dependencies have a joint inference calculus \cite{BeeriFaginHoward1977}.
	
	In the opening example the seed consists of the $p^2$ triples
	$(a,b,a+b)$ over $\mathbb F_p$ at a single boundary value. All $p^2$
	pairs occur on the last two charts, so isolating the first chart makes
	the least recombinant enlargement the full set of $p^3$ triples.
	Every chart-surjective source exact at this cut loses the determination
	of the first chart by the other two, since all $p$ first readings
	must occur with each attained pair.
	
	Vector-valued charts illustrate what whole-chart recovery can miss.
	Take $Q=1$ and two charts with private values $(a,x)$ and $(a,y)$
	over $k$, with $a,x,y$ independent coordinates in the source.
	Neither whole chart determines the other. Recombination preserves
	this recovery law while increasing the dimension of the source
	from three to four. It removes equality of the first scalar readings.
	Naming those readings therefore changes the preservation problem.
	
	\section{Absolute source minima over a common quotient}
	\label{sec:minima}
	
	A chart-surjective subgroup of a finite matching group can realize
	fewer states than the full group. For central extensions of a common
	quotient with elementary abelian kernels, \cref{thm:min} identifies
	the subspace forced by the extension and computes the least source
	order compatible with the cuts. \Cref{thm:conditional} minimizes
	conditional ranges at that order. The subsequent results determine
	when several observations share an optimizer, extend the calculation
	to quotient readouts, and realize representable rank functions as
	cut costs.
	
	\subsection{A sharp minimum over all sources}
	\label{sub:min}
	Let $k=\mathbb F_p$. For each $i$ in a finite nonempty set $I$, let
	\[
	1\longrightarrow A_i\longrightarrow E_i
	\xrightarrow{\pi_i} Q\longrightarrow1
	\]
	be a central extension of the same finite group $Q$, with $A_i$ an
	elementary abelian $p$-group of finite dimension over $k$. The extensions may be nonsplit. Put
	$A=\bigoplus_{i\in I}A_i$ and let
	$G=\{(x_i)\in\prod_i E_i:\pi_i(x_i)\text{ is independent of }i\}$.
	This is the native matching group on any connected chart graph that
	uses the displayed maps to $Q$ on every edge. Every nontrivial
	vertex cut has complete boundary image a diagonal copy of $Q$.
	A finite tuple of matching states generates a subgroup of $G$.
	Every subgroup occurs under this formation law because $G$ is finite. For $X\subseteq I$, identify $A_X=\bigoplus_{i\in X}A_i$ with
	its coordinate subspace in $A$ and write $\operatorname{pr}_X$
	for the corresponding projection.
	
	Define $P_p(G)=[G,G]G^p$, where $G^p$ is generated by all $p$th
	powers. The quotient $G/P_p(G)$ is the largest elementary abelian
	$p$-quotient of $G$. Set
	\begin{equation}\label{eq:support-atoms}
		W=A\cap P_p(G),\qquad
		W_X=\operatorname{pr}_X W,\qquad
		d_i=\dim_k A_i-\dim_k W_{\{i\}}.
	\end{equation}
	For a finite $p$-group $G$, $P_p(G)$ is its Frattini subgroup
	\cite[Theorem~5.3.2]{Robinson1996}.
	
	\begin{theorem}[Admissible kernels and minimum source order under prescribed cuts]
		\label{thm:min}
		\leavevmode
		\begin{enumerate}[label=\textup{(\roman*)},leftmargin=*]
			\item \emph{Admissible kernels.} A subspace $K\le A$ occurs as $H\cap A$ for a chart-surjective
			native source $H\le G$ exactly recombinant at all protected cuts
			if and only if
			\begin{equation}\label{eq:kernel}
				W\le K,\qquad
				K=\bigoplus_{B\in\Pi}\operatorname{pr}_B K,\qquad
				\operatorname{pr}_{\{i\}} K=A_i\quad(i\in I).
			\end{equation}
			
			\item \emph{Sources with a prescribed kernel.} For each such $K$, the sources with kernel $K$ are in bijection with
			the homomorphic sections of $G/K\to Q$. There are
			$|\operatorname{Hom}(Q,A/K)|$ of them.
			
			\item \emph{Minimum source order.} The minimum source order over
			all admissible kernels is
			\begin{equation}\label{eq:minimum}
				\boxed{\min_H |H|=|Q|\,p^{m(\Pi)},\qquad
					m(\Pi)=\sum_{B\in\Pi}
					\left(\dim_k W_B+\max_{i\in B}d_i\right).}
			\end{equation}
		\end{enumerate}
	\end{theorem}
	\begin{proof}
		\leavevmode
		\begin{enumerate}[label=\textup{(\roman*)},leftmargin=*]
			\item \emph{Admissible kernels.} A source mapping onto $Q$ satisfies $G=HA$. Since $A$ is central
			of exponent $p$, commutators and $p$th powers in $G$ can be computed
			using the $H$-factors alone. Hence $P_p(G)=P_p(H)\le H$, so
			$W\le H\cap A$. Conversely, suppose $W\le K\le A$.
			The map $A/K\to G/(KP_p(G))$ is injective, since
			$A\cap KP_p(G)=K(A\cap P_p(G))=K$. In the elementary abelian
			quotient choose a linear retraction onto this copy of $A/K$.
			Its composite with $G/K\to G/(KP_p(G))$ is a group retraction onto
			$A/K$. The kernel maps isomorphically onto $Q$ and hence supplies
			a section of $G/K\to Q$.
			
			Let $H$ be the inverse image in $G$ of the image of this section. Then $H\cap A=K$ and $H$ maps onto $Q$. For $q\in Q$ choose $h_q\in H$ over $q$. The fiber of $H$ over $q$ is
			$h_qK$. Consequently its $i$th chart covers the entire fiber of
			$E_i$ exactly when $\operatorname{pr}_{\{i\}} K=A_i$.
			For a cut $S\mid S^c$, taking arbitrary sides in this fiber yields
			the translate of
			$\operatorname{pr}_S K\oplus\operatorname{pr}_{S^c}K$ by $h_q$.
			Equality with the source fiber means that both coordinate projections,
			extended by zero, preserve $K$. Products of these projections over
			the protected cuts give the individual atom projections. Conversely,
			summing the appropriate atom projections gives each side projection.
			This proves \eqref{eq:kernel}.
			
			\item \emph{Sources with a prescribed kernel.} Every subgroup over $Q$ with intersection $K$ is the inverse image of the image of a unique section. Given one section, every other is obtained by multiplying it by a unique homomorphism $Q\to A/K$, since $A/K$ is central.
			
			\item \emph{Minimum source order.} It remains to minimize separately on each atom $B$. Any
			subspace $K_B\le A_B$ containing $W_B$ and surjecting onto
			each $A_i$ has $\dim K_B/W_B\ge d_i$ for every $i\in B$. Set
			$t=\max_{i\in B}d_i$. Choose surjections
			$k^t\to A_i/W_{\{i\}}$ for all $i\in B$, lift each linearly to
			$A_i$, and combine the lifts in a map $j:k^t\to A_B$. At an index
			attaining $t$, the map into the quotient is an isomorphism, so
			$j$ is injective modulo $W_B$. Then $K_B=W_B+j(k^t)$ attains the
			bound and projects onto every $A_i$. This also works at $t=0$.
			Taking the sum over atoms and using $|H|=|Q|\,|K|$ proves the formula.\qedhere
		\end{enumerate}
	\end{proof}

	The same forced subspace can be read from the extension class. The class
	$e\in H^2(Q,A)=H^2(Q,k)\otimes_k A$ of the joint extension maps
	to zero in $H^2(Q,A/K)$ exactly when the quotient extension splits
	\cite{BrownCohomology1982}. By part~\textup{(i)} of the proof, this
	is equivalent to $W\le K$. More explicitly, write $e=\sum_j e_j\otimes a_j$ in a basis
	$(e_j)$ of $H^2(Q,k)$. Its image in $H^2(Q,A/K)$ vanishes exactly
	when every $a_j$ belongs to $K$. The least such $K$ is therefore
	$\langle a_j:j\rangle$, the image obtained by contracting $e$ with
	the dual of $H^2(Q,k)$. Thus the contraction image is $W$.

	\begin{corollary}[Polymatroid ranks determined by cut minima]
		\label{cor:rank}
		\leavevmode
		\begin{enumerate}[label=\textup{(\roman*)},leftmargin=*]
			\item \emph{Rank representation.}
			There are subspaces $U_i$ of a single $k$-space such that
			\[
			f(S)=\dim_k\sum_{i\in S}U_i
			=\dim_k W_S+\max_{i\in S}d_i,\qquad
			m(\Pi)=\sum_{B\in\Pi} f(B),
			\]
			where $f(\varnothing)=0$ and the empty maximum is zero.
			
			\item \emph{Recovery from singleton cuts.}
			The chart dimensions and the table of minima for all families of
			singleton cuts recover $f(S)$ for every $S$.
			
			\item \emph{Equality of cut increments.}
			Put $\Delta_f(\Pi)=\sum_{B\in\Pi} f(B)-f(I)$.
			Two such rank functions	$f,g$ on the same labeled set have identical increases $\Delta$ for all families of singleton cuts if and only if there are weights $w_i$ such that $f(S)-g(S)=\sum_{i\in S}w_i$ for every $S\subseteq I$.
			In that case the increases agree for every partition.
		\end{enumerate}
	\end{corollary}

	\begin{proof}
		\leavevmode
		\begin{enumerate}[label=\textup{(\roman*)},leftmargin=*]
			\item \emph{Rank representation.} The restriction of the $i$th coordinate dual to $W$ has image
			$R_i\le W^*$, with $\dim\sum_{i\in S}R_i=\dim W_S$ by duality.
			Choose a flag $0=F_0\subset F_1\subset\cdots\subset F_{\max d_i}$
			with $\dim F_j=j$ in a separate space. Take $U_i=R_i\oplus F_{d_i}$.
			
			\item \emph{Recovery from singleton cuts.} For nonempty $S$, protecting the singleton cuts at all indices
			outside $S$ produces the atom $S$ and singleton outsiders. Thus its
			minimum exponent is $f(S)+\sum_{i\notin S}\dim A_i$. For a single chart the partition is obtained from the empty cut family.
			
			\item \emph{Equality of cut increments.} Let $\Pi_S$ be this partition for nonempty $S$, and let $\Pi_{\rm disc}$
			be the singleton partition. Subtraction gives
			\[
			\Delta_f(\Pi_{\rm disc})-\Delta_f(\Pi_S)
			=\sum_{i\in S} f(\{i\})-f(S).
			\]
			If the increases agree for $f$ and $g$, subtract their identities.
			This gives the stated weights $w_i=f(\{i\})-g(\{i\})$, also for the
			empty set since both ranks there are zero. Conversely any such weights
			sum to the same total over the blocks of a partition as over $I$,
			so their contribution to every increase vanishes.\qedhere
		\end{enumerate}
	\end{proof}

	For a single cut $S\mid S^c$ the increase in the minimum exponent
	is $f(S)+f(S^c)-f(I)$, the dimension of the intersection of the
	two subspace sums. For a binary hierarchy of cuts, the largest such
	individual increase is the width of its subspace branch decomposition
	under the intersection-dimension convention
	\cite{JeongKimOum2021}. Imposing every cut in that hierarchy has the
	joint cost \eqref{eq:minimum}. Linear code realizations on trees
	and lossless joins supply related cut constructions
	\cite{Kashyap2009,Fagin1977,Kobayashi1994}.

	\subsection{Least conditional range at minimum source order}

	For a chart-surjective subgroup as in \cref{thm:min}, fix sets
	$P,J\subseteq I$ of observed and requested whole charts,
	respectively. Fix an attained common quotient value $q$ and an attained
	tuple of values on $P$ above $q$. The number of possible tuples on
	$J$ will be called the conditional range size. If $P\ne\varnothing$,
	observing even one whole chart already determines $q$. If
	$P=\varnothing$ and $J\ne\varnothing$, omitting the conditioning on
	$q$ multiplies the range size by $|Q|$.

	\begin{theorem}[Conditional range minima and simultaneous attainment at minimum source order]
		\label{thm:conditional}
		Under the central-extension hypotheses of \cref{thm:min}, write
		$W_X=\operatorname{pr}_X W$, $r(X)=\dim_k W_X$, and
		$d_i=\dim_k A_i-r(\{i\})$. For each cut atom $B\in\Pi$ set
		\[
		S_B=P\cap B,\qquad U_B=(J\setminus P)\cap B,\qquad
		a_B=\dim_k A_{S_B}-r(S_B),\qquad
		e_B=\max_{i\in U_B}d_i.
		\]
		Take $A_\varnothing=0$, $r(\varnothing)=0$ and the empty maximum
		equal to zero.
		\begin{enumerate}[label=\textup{(\roman*)},leftmargin=*]
			\item \emph{A fixed observed set.} Among all chart-surjective subgroups exact at the
			protected cuts, the least conditional range size is $p^{b(P,J)}$, where
			\begin{equation}\label{eq:conditional}
				b(P,J)=\sum_{B\in\Pi}
				\left[r(S_B\cup U_B)-r(S_B)+\max\{0,e_B-a_B\}\right].
			\end{equation}
			For each fixed observed set $P$ there is a source of the absolute
			minimum order \eqref{eq:minimum} that attains $p^{b(P,J)}$ for
			every requested set $J\subseteq I$ and every attained conditioning
			value. In particular, conditional on $Q$, exact recovery of $J$ is attainable
			at that order precisely when, for every $B$,
			\[
			r(S_B\cup U_B)=r(S_B)\quad\text{and}\quad e_B\le a_B.
			\]
			
			\item \emph{All observed sets simultaneously.} One source attains these minima for \emph{every} pair $(P,J)$ if and
			only if the following test holds in each atom $B$. Put
			$t_B=\max_{i\in B}d_i$ and $L_B=\{i\in B:d_i<t_B\}$. The test is
			\begin{equation}\label{eq:all-observations}
				\dim_k A_{L_B}-r(L_B)=\max_{i\in L_B}d_i .
			\end{equation}
			Whenever it holds, every source optimal for all requests at
			$P=\varnothing$ is a simultaneous optimizer of minimum order.
		\end{enumerate}
	\end{theorem}

	\begin{proof}
		\leavevmode
		\begin{enumerate}[label=\textup{(\roman*)},leftmargin=*]
			\item \emph{A fixed observed set.} Write $K=H\cap A$. Above any attained boundary value the source
			fiber is a translate of $K$. Once the charts in $P$ are also fixed,
			the possible outputs are a translate of
			$\operatorname{pr}_J(K\cap\ker\operatorname{pr}_P)$.
			By \cref{thm:min}, $K$ is the sum of its atom subspaces, so
			the exponents of the range sizes add over atoms.
			
			Fix one atom $B$. Write $K_B=\operatorname{pr}_B K$, $S=S_B$,
			$U=U_B$, $a=a_B$, $e=e_B$ and $D_i=A_i/W_{\{i\}}$. For
			$X\subseteq B$ put $Y_X=\operatorname{pr}_X K_B/W_X$.
			By rank-nullity the exponent on this atom equals
			\[
			\dim\operatorname{pr}_{S\cup U}K_B-\dim\operatorname{pr}_S K_B
			=r(S\cup U)-r(S)+\dim Y_{S\cup U}-\dim Y_S.
			\]
			The last difference is nonnegative since restriction to $S$ is
			surjective. Coverage makes $Y_{S\cup U}$ surject onto every $D_i$
			with $i\in U$, whereas $Y_S\le A_S/W_S$. The difference is thus at
			least $\max\{0,e-a\}$.
			
			We construct one $K_B$ for all requests with the given observed set.
			Put $t=\max_{i\in B}d_i$, $h=\min(t,a)$, and $E=k^t$ with ordered
			coordinates $x_1,\ldots,x_t$. First choose an $h$-plane
			$Y\le A_S/W_S$ projecting onto each $D_i$ for $i\in S$.
			Indeed, set $d=\max_{i\in S}d_i$, using $d=0$ for empty $S$.
			Then $d\le h$, since $A_S/W_S$ projects onto every observed $D_i$.
			Choose surjections $k^d\to D_i$ and combine them. The image has
			dimension $d$, as one coordinate attains $d$ when $S$ is nonempty.
			Lift the image through the surjection
			$A_S/W_S\to\bigoplus_{i\in S}D_i$ and extend it to an $h$-plane.
			For $S=\varnothing$ take $Y=0$.
			
			Identify $k^h$ with a linear lift of $Y$ in $A_S$. Define
			$j:E\to A_B$ on $S$ by sending the first $h$ coordinates through
			this identification. For every $i\in B\setminus S$, choose a
			complement to $W_{\{i\}}$ in $A_i$, identify it with $k^{d_i}$,
			and set $j_i(x)=(x_1,\ldots,x_{d_i})$ in that complement. Put
			$K_B=W_B+j(E)$. It covers each $A_i$. At a coordinate attaining
			$d_i=t$, the residual map $E\to D_i$ is onto and hence injective,
			so $\dim K_B=\dim W_B+t$, including the case $t=0$.
			
			For any $U\subseteq B\setminus S$, consider the map from $E$ to
			$A_{S\cup U}/W_{S\cup U}$ induced by $j$. A vector in its kernel
			must have $x_1=\cdots=x_h=0$, because the observed lift meets $W_S$
			only in zero. Its $i$th coordinate for $i\in U$ must belong both to
			the chosen complement and to $W_{\{i\}}$, so its first $d_i$
			coordinates vanish as well. These conditions suffice, since they
			make the entire projection of $j(x)$ zero. The map therefore has
			rank $\max(h,e)$, where $e=\max_{i\in U}d_i$. Its image is
			$Y_{S\cup U}$ and $\dim Y_S=h$. The residual increment is
			$\max(0,e-h)=\max(0,e-a)$ because $e\le t$.
			The construction depends only on $S=P\cap B$. Take the direct sum
			over $B$ and invoke \cref{thm:min} once to realize this kernel
			by a minimum-order source. The same source attains the bound for all
			$J$ and all attained conditioning values.
			
			\item \emph{All observed sets simultaneously.} Fix any source attaining
			every minimum at $P=\varnothing$; the construction above supplies one.
			Applying the range formula with request $X\subseteq B$ shows that its
			kernel has
			\[
			\dim Y_X=u(X),\qquad u(X):=\max_{i\in X} d_i
			\quad(u(\varnothing)=0).
			\]
			In particular $\dim K_B=r(B)+t_B$ in every atom, so the source has
			minimum order. Write $c(S)=\dim A_S-r(S)$ for $S\subseteq B$.
			For observations $S=P\cap B$ and new requests $U=(J\setminus P)\cap B$,
			the residual increment is now fixed: it is
			$u(S\cup U)-u(S)=\max\{0,u(U)-u(S)\}$. The sharp residual in
			\eqref{eq:conditional} is $\max\{0,u(U)-c(S)\}$.
			Since $u(S)\le c(S)$, the two increments agree for every $U$ if
			$u(S)=t_B$ or $c(S)=u(S)$. Conversely, if $u(S)<t_B$, request a
			single index of deficit $t_B$ in $B\setminus S$. Equality of the two
			increments for this request forces $c(S)=u(S)$.
			
			Thus it suffices to consider $S\subseteq L_B$. Necessity of
			\eqref{eq:all-observations} follows by testing $S=L_B$ when
			$t_B>0$; when $t_B=0$ the equality holds on the empty set. For
			sufficiency, use
			\[
			c(L_B)=\sum_{i\in L_B}d_i+
			\operatorname{codim}_{\bigoplus_{i\in L_B}W_{\{i\}}} W_{L_B}.
			\]
			Both terms are nonnegative, and the first is at least $u(L_B)$.
			Equality with $u(L_B)$ means that at most one index of $L_B$ has
			positive deficit and $W_{L_B}$ is the full product of its coordinate
			images. These two properties pass to every subset $S\subseteq L_B$.
			Hence $c(S)=u(S)$ throughout $L_B$. Every other $S$ contains a chart
			with deficit $t_B$. This proves the test in every atom and hence the
			simultaneous assertion.\qedhere
		\end{enumerate}
	\end{proof}

	\begin{example}[Two queries with incompatible minimizers]
		For any finite prime field $k=\mathbb F_p$, take $Q=1$, no protected cuts, and $A_1=A_2=k$, $A_3=k^2$. Then
		$W=0$ and the minimum order is $p^2$. The source
		$\{(x,x,(x,y)):x,y\in k\}$ realizes $\{1\}\to2$. The source
		$\{(x,y,(x,y)):x,y\in k\}$ realizes $\{1,2\}\to3$.
		Any source satisfying both implications would make chart $1$
		determine chart $3$, contradicting chart coverage since $p<p^2$. The single atom has
		deficits $(1,1,2)$ and $W=0$. For $L_B=\{1,2\}$ the two sides
		of \eqref{eq:all-observations} are $2$ and $1$, respectively.
	\end{example}

	\paragraph{Families with a least observed set on each atom.}
	Let $\mathcal P$ be a nonempty family of observed sets. Suppose that
	in each atom $B$ the family $\{P\cap B:P\in\mathcal P\}$ has a least
	member $S_0$ under inclusion. Different atoms may obtain their least
	members from different observations. This hypothesis holds for every
	nested family and every family containing the empty observation.
	Within $B$ put $h(S)=\min\{t_B,\dim A_S-r(S)\}$ and
	$u(X)=\max_{i\in X}d_i$ with $u(\varnothing)=0$.
	One minimum-order source attains all requested sets for all
	$P\in\mathcal P$ if and only if in every atom
	\begin{equation}\label{eq:least-observation}
		h(S)=\max\{h(S_0),u(S\setminus S_0)\}
		\qquad (S=P\cap B,\ P\in\mathcal P).
	\end{equation}

	This criterion follows from the general observation-family test in
	\cref{cor:incomparable-cases}\textup{(i)}.

	For a fixed nested pair of codes, the inverse relative
	dimension/length profile minimizes a difference of punctured-code
	dimensions over coordinate sets of a given size
	\cite[Definition~1]{ZhuangDaiLuoHanVinck2014}. The quantity $\dim Y_X$ is that
	difference on the specified coordinates of $X$ after choosing bases
	in the chart kernels. Epperly, Govindarajan and Chandrasekaran
	\cite{EpperlyGovindarajanChandrasekaran2021} obtain simultaneous minimum ranks for overlapping
	submatrices by choosing a shared matrix block.
	In \cref{thm:conditional}, the observations and chart maps are fixed.
	The kernel $K$ varies subject to $W\le K$, the cut constraints, and surjectivity on each chart kernel.
	
	For instance take $Q=k^2$ and charts
	$E_1=C_{p^2}\times C_p\times k$ and
	$E_2=C_p\times C_{p^2}\times k^2$. Each quotient reads the first
	two coordinates modulo $p$. The $p$th powers of lifts of the two
	quotient basis vectors give independent kernel lines, one in each
	chart. These generate $W$ because the matching group is abelian
	and all its other coordinates have exponent $p$. Thus the forced
	ranks are one per chart and two jointly, and the residual
	dimensions are $d_1=1$ and $d_2=2$. With chart one observed, the
	forced contribution for chart two is one. With no cut, one
	residual coordinate can pass through the observed chart, leaving
	$p^2$ possible outputs in a minimizing source. Protecting the
	cut between the charts removes that capacity and leaves $p^3$
	possible outputs, exactly the size of the second chart's private fiber.

	\subsection{Conditional minima for arbitrary observation families}
	\label{subsec:incomparable}

	We extend \cref{thm:conditional} and the least-member test
	\eqref{eq:least-observation} to arbitrary families of whole-chart observations. Fix a cut atom $B$ and keep the
	notation $k=\mathbb F_p$, $r$, $d_i$, $t=t_B=\max_{i\in B}d_i$, $D_i=A_i/W_{\{i\}}$,
	$c(S)=\dim_kA_S-r(S)$, $h(S)=\min\{t,c(S)\}$ and
	$u(X)=\max_{i\in X}d_i$, with $u(\varnothing)=0$. For $X\subseteq Y\subseteq B$, the map
	$A_Y/W_Y\to A_X/W_X$ is the one induced by projection. A
	\emph{minimum-order kernel} on $B$ is a subspace $K_B$ with
	$W_B\le K_B\le A_B$, $\operatorname{pr}_{\{i\}}K_B=A_i$ for every
	$i\in B$, and $\dim K_B=r(B)+t$. For a nonempty family $\mathcal S$ of
	subsets of $B$, call $K_B$ \emph{$\mathcal S$-optimal} if
	\begin{equation}\label{eq:S-optimal}
		\dim\operatorname{pr}_{S\cup U}K_B-\dim\operatorname{pr}_SK_B
		=r(S\cup U)-r(S)+\max\{0,u(U)-h(S)\}
	\end{equation}
	for all $S\in\mathcal S$ and $U\subseteq B\setminus S$. By the proof of
	\cref{thm:conditional}, the left side is the exponent of the
	conditional range on $B$ for the observed set $S$ and the request $U$,
	and the right side is its sharp minimum. Kernels decompose over the
	atoms (\cref{thm:min}). Hence one minimum-order source attains every
	conditional minimum for every observed set in a family $\mathcal P$ if
	and only if, on every atom $B$, some minimum-order kernel is
	$\{P\cap B:P\in\mathcal P\}$-optimal. Put
	\[
	h_0=\min_{S\in\mathcal S}h(S),\qquad
	\mathcal S_0=\{S\in\mathcal S:h(S)=h_0\},\qquad
	T_0=\bigcup_{S\in\mathcal S_0}S .
	\]

	\begin{theorem}[Simultaneous optimality for arbitrary observation families]
		\label{thm:incomparable-test}
		An $\mathcal S$-optimal minimum-order kernel exists on $B$ if and only
		if both of the following hold.
		\begin{enumerate}[label=\textup{(\alph*)},nosep]
			\item \emph{Numerical condition.}
			$h(S)=u(S)$ for every $S\in\mathcal S\setminus\mathcal S_0$.
			
			\item \emph{Subspace condition.}
			The space $A_{T_0}/W_{T_0}$ contains an $h_0$-dimensional
			subspace $\Lambda$ that maps onto $D_i$ for every $i\in T_0$ and
			injectively into $A_S/W_S$ for every $S\in\mathcal S_0$.
		\end{enumerate}
		When they hold, the kernel constructed in the proof from any such
		$\Lambda$ is $\mathcal S$-optimal.
	\end{theorem}

	\begin{proof}
		Write $y(X)=\dim\operatorname{pr}_XK_B-r(X)$.
		Note first that $u(S)\le h(S)$ for every $S\subseteq B$, because $A_S/W_S$ maps onto each $D_i$ with $i\in S$. In particular $d_i\le h_0$ for every
		$i\in T_0$.
		
		Since $\dim K_B=r(B)+t$, \eqref{eq:S-optimal} with $U=B\setminus S$
		gives $y(S)=t-\max\{0,u(B\setminus S)-h(S)\}$.
		If $u(B\setminus S)=t$, this is $h(S)$.
		Otherwise a chart of deficit $t$ lies in $S$, so $u(S)=t=h(S)$ and $t\le y(S)\le y(B)=t$.
		In both cases $y(S)=h(S)$, and \eqref{eq:S-optimal} becomes
		\begin{equation}\label{eq:S-optimal-y}
			y(S\cup U)=\max\{h(S),u(U)\}\qquad
			(S\in\mathcal S,\ U\subseteq B\setminus S).
		\end{equation}
		Conversely, \eqref{eq:S-optimal-y} for all such pairs implies
		\eqref{eq:S-optimal}, by subtracting the instance $U=\varnothing$.
		\begin{itemize}
			\item \emph{Necessity.}
			Let $S\in\mathcal S\setminus\mathcal S_0$ and $S'\in\mathcal S_0$, and put $X=S\cup S'$.
			Applying \eqref{eq:S-optimal-y} from both sides,
			\[
			\max\{h_0,u(S\setminus S')\}=y(X)=\max\{h(S),u(S'\setminus S)\}
			\ge h(S)>h_0 .
			\]
			Hence $h(S)\le u(S\setminus S')\le u(S)\le h(S)$, which proves (a).
			For (b), put $\Lambda=\operatorname{pr}_{T_0}K_B/W_{T_0}$. For any
			$S\in\mathcal S_0$, \eqref{eq:S-optimal-y} with $U=T_0\setminus S$
			gives $\dim\Lambda=y(T_0)=h_0$, since $u(T_0)\le h_0$. Coverage of the
			charts makes $\Lambda$ map onto each $D_i$. Its image in $A_S/W_S$ has
			dimension $y(S)=h_0$, so that map is injective.
			
			\item \emph{Sufficiency.}
			Let $E=k^t$ with coordinates $x_1,\dots,x_t$, and recall $h_0\le t$. Choose a linear map $\lambda:k^{h_0}\to A_{T_0}$ lifting an isomorphism onto $\Lambda$.
			For $i\in B\setminus T_0$ choose a complement $C_i$ of $W_{\{i\}}$ in $A_i$ and an isomorphism $c_i:k^{d_i}\to C_i$.
			Define $j:E\to A_B$ by
			\[
			j(x)=\lambda(x_1,\dots,x_{h_0})
			+\sum_{i\in B\setminus T_0}c_i(x_1,\dots,x_{d_i}),
			\]
			where the first term is supported on $T_0$ and the $i$th summand on
			chart $i$. Put $K_B=W_B+j(E)$, so that $y(X)$ is the rank of the
			induced map $E\to A_X/W_X$.
			
			Call $X\subseteq B$ \emph{good} if it contains some $S\in\mathcal S_0$
			or a chart outside $T_0$ of deficit greater than $h_0$. We claim that
			for good $X$ the kernel of $E\to A_X/W_X$ is the subspace
			$x_1=\dots=x_m=0$, where $m=\max\{h_0,u(X\setminus T_0)\}$. Suppose
			$\operatorname{pr}_Xj(x)\in W_X$. For $i\in X\setminus T_0$ the
			component $c_i(x_1,\dots,x_{d_i})$ lies in $W_{\{i\}}\cap C_i=0$, so
			$x_1,\dots,x_{u(X\setminus T_0)}$ vanish. If $X\supseteq S\in\mathcal S_0$,
			then $\operatorname{pr}_S\lambda(x_1,\dots,x_{h_0})\in W_S$, and
			injectivity of $\Lambda\to A_S/W_S$ gives $x_1=\dots=x_{h_0}=0$. In the
			other case of goodness, $u(X\setminus T_0)>h_0$ already gives this.
			Conversely, if $x_1=\dots=x_m=0$, then $\operatorname{pr}_Xj(x)=0$.
			This proves the claim, so $y(X)=m$ for good $X$.
			
			Let $S\in\mathcal S_0$ and $U\subseteq B\setminus S$. Then $S\cup U$
			is good, and since $d_i\le h_0$ on $T_0$,
			$y(S\cup U)=\max\{h_0,u(U\setminus T_0)\}=\max\{h_0,u(U)\}$. Let
			$S\in\mathcal S\setminus\mathcal S_0$. By (a) and $d_i\le h_0<h(S)$ on
			$T_0$, some chart of deficit $h(S)$ lies in $S\setminus T_0$. Hence
			$S\cup U$ is good and
			$y(S\cup U)=\max\{h(S),u(U\setminus T_0)\}=\max\{h(S),u(U)\}$. This is
			\eqref{eq:S-optimal-y}. Coverage holds on $T_0$ because $\Lambda$
			maps onto each $D_i$, and outside $T_0$ by construction. Finally $B$ is
			good, so $y(B)=\max\{h_0,u(B\setminus T_0)\}$. If a chart of deficit
			$t$ lies outside $T_0$ this is $t$; otherwise $t\le h_0\le t$. Thus
			$\dim K_B=r(B)+t$.\qedhere
		\end{itemize}
	\end{proof}

	\begin{corollary}[Numerical criteria and pairwise obstructions to simultaneous optimality]
		\label{cor:incomparable-cases}
		In \cref{thm:incomparable-test} the following hold.
		\begin{enumerate}[label=\textup{(\roman*)},nosep]
			\item \emph{Least observed sets.}
			If $\mathcal S_0$ has a least member under inclusion, then (b) holds.
			If $\mathcal S$ has a least member $S_*$, then (a) is equivalent to \eqref{eq:least-observation}; thus the theorem recovers that test.
			
			\item \emph{Large fields.}
			If $p\ge1+|\mathcal S_0|+|\{i\in T_0:d_i>0\}|$, then (b) holds.
			Over such fields the test consists of (a) alone.
			
			\item \emph{Pairwise compatibility.}
			The pairwise condition $\max\{h(S_1),u(S_2\setminus S_1)\}=\max\{h(S_2),u(S_1\setminus S_2)\}$ for all $S_1,S_2\in\mathcal S$ is necessary but not sufficient.
		\end{enumerate}
	\end{corollary}

	\begin{proof}
		\leavevmode
		\begin{enumerate}[label=\textup{(\roman*)},leftmargin=*]
			\item \emph{Least observed sets.} Let $S_*$ be least in $\mathcal S_0$ and put
			$d=u(S_*)\le h_0\le c(S_*)$. As in the proof of
			\cref{thm:conditional}, combine surjections $k^d\to D_i$ for
			$i\in S_*$, lift through $A_{S_*}/W_{S_*}\to\bigoplus_{i\in S_*}D_i$, and
			extend to an injective map $\lambda_*:k^{h_0}\to A_{S_*}/W_{S_*}$ whose
			image maps onto every $D_i$ with $i\in S_*$. For
			$i\in T_0\setminus S_*$ send $x$ to $(x_1,\dots,x_{d_i})$ in a complement
			of $W_{\{i\}}$, which is possible because $d_i\le h_0$. The resulting
			map $k^{h_0}\to A_{T_0}/W_{T_0}$, composed with the projection to
			$A_{S_*}/W_{S_*}$, is $\lambda_*$. It is therefore injective, and its
			image $\Lambda$ maps onto $D_i$ for every $i\in T_0$. Every
			$S\in\mathcal S_0$ contains $S_*$, and the projection to
			$A_{S_*}/W_{S_*}$ factors through $A_S/W_S$, so $\Lambda\to A_S/W_S$ is
			injective.
			
			Now let $S_*$ be least in $\mathcal S$. Since $c$ is monotone,
			$h_0=h(S_*)$. For $S\in\mathcal S_0$, both (a) and
			\eqref{eq:least-observation} hold automatically, the latter because
			$u(S\setminus S_*)\le u(S)\le h(S)=h_0$. For $h(S)>h_0$, the test
			\eqref{eq:least-observation} reads $h(S)=u(S\setminus S_*)$, which
			implies $h(S)=u(S)$. Conversely $u(S)=h(S)>h_0\ge u(S_*)$ places a chart
			of deficit $u(S)$ in $S\setminus S_*$.
			
			\item \emph{Large fields.} For $h_0=0$ take $\Lambda=0$, using $d_i\le h_0$ on $T_0$. Let
			$h_0\ge1$ and choose a linear map $\varphi:k^{h_0}\to A_{T_0}/W_{T_0}$
			uniformly at random. Its composite with each surjection onto $D_i$ or
			onto $A_S/W_S$ is again uniformly distributed. A uniform linear map
			$k^a\to k^b$ fails to be onto with probability at most
			$\sum_{j<b}p^{j-a}<p^{b-a}/(p-1)$ when $b\le a$, and fails to be
			injective with probability at most $\sum_{j<a}p^{j-b}<p^{a-b}/(p-1)$
			when $a\le b$. Each required condition therefore fails with probability
			less than $1/(p-1)$, and the conditions for charts with $d_i=0$ always
			hold. Under the hypothesis some $\varphi$ satisfies every condition.
			It is injective because $\mathcal S_0\ne\varnothing$, and
			$\Lambda=\varphi(k^{h_0})$ satisfies (b).
			
			\item \emph{Pairwise compatibility.} Necessity follows from \eqref{eq:S-optimal-y} with
			$X=S_1\cup S_2$. \Cref{ex:field-dependence} shows that the condition
			is not sufficient.\qedhere
		\end{enumerate}
	\end{proof}
	
	\begin{example}[The test depends on the field]
		\label{ex:field-dependence}
		Take $Q=1$ and charts $A_1=\dots=A_4=k$ and $A_5=k^2$, all in one atom
		$B=\{1,\dots,5\}$. Then $G=A$, $P_p(G)=0$ and $W=0$, so $r=0$,
		$d=(1,1,1,1,2)$ and $t=2$. Let $\mathcal S$ consist of the six
		two-element subsets of $\{1,2,3,4\}$. Every $S\in\mathcal S$ has
		$c(S)=h(S)=2$, so $\mathcal S_0=\mathcal S$ and $T_0=\{1,2,3,4\}$.
		Condition (a) is vacuous and the pairwise condition holds. Condition
		(b) asks for a plane $\Lambda\le k^4$ on which the four coordinate
		functionals are nonzero and pairwise independent, that is, for four
		distinct points of the projective line over $\mathbb F_p$. Such a plane
		exists if and only if $p\ge3$. The numbers $r(X)$ and $\dim_kA_X$ are
		the same for every $p$. Hence no test expressed through these numbers
		alone is valid for all primes. The count of projective points proves both outcomes directly.
	\end{example}
	
	\begin{theorem}[NP-completeness of simultaneous optimality for observation families]
		\label{thm:incomparable-npc}
		Fix a prime $p$.
		\begin{enumerate}[label=\textup{(\roman*)},leftmargin=*]
			\item \emph{Decision complexity.} Given matrices for the chart kernels $A_i$ $(i\in B)$
			and for $W_B$, together with an explicit family $\mathcal S$, deciding
			whether an $\mathcal S$-optimal minimum-order kernel exists is
			NP-complete. Hardness holds already when $t_B=1$ and $h\equiv1$ on
			$\mathcal S$.
			
			\item \emph{Realization by central extensions.} Every subspace $W_B\le A_B$ arises from central
			extensions of an elementary abelian $p$-group $Q$ with abelian chart
			groups of exponent dividing $p^2$. Thus the hardness applies in the
			setting of \cref{thm:conditional} itself, already for such charts.
		\end{enumerate}
	\end{theorem}
	\begin{proof}
		\leavevmode
		\begin{enumerate}[label=\textup{(\roman*)},leftmargin=*]
			\item \emph{Decision complexity.} Condition (a) is checked by rank computations, and a basis of
			$\Lambda$ is a certificate for (b). For hardness, consider proper
			subspaces $H_1,\dots,H_s$ of $V=k^n$ with
			$\sum_{c'\ne c}H_{c'}=V$ for every $c$, and ask whether
			$V\ne\bigcup_cH_c$. This question is NP-hard. For $p\ge3$, take a simple graph
			on $[n]$ with at least three edges (adding disjoint edges if necessary) and the hyperplanes
			$H_{uv}=\{x_u=x_v\}$. Any two distinct hyperplanes span $V$, and a
			vector outside every $H_{uv}$ is a proper coloring with $p$ colors,
			so this is graph $p$-colorability
			\cite[Problem~GT4]{GareyJohnson1979}. For $p=2$, take a 3-CNF
			formula with no clause containing complementary literals. Introduce a homogenizing variable $z$ and append two positive
			clauses on disjoint triples of fresh variables. Either new clause
			can be satisfied independently of the original assignment. To a clause with
			literals $\ell_1,\ell_2,\ell_3$ assign the subspace on which the forms
			$\hat\ell_1,\hat\ell_2,\hat\ell_3$ vanish, where $\hat\ell=x_j$ for
			$\ell=x_j$ and $\hat\ell=x_j+z$ for $\ell=\neg x_j$. Add the hyperplane
			$\{z=0\}$. A vector outside all these subspaces has $z=1$ and encodes a
			satisfying assignment, and conversely. The annihilators of the two fresh clause subspaces and of
			$\{z=0\}$ are supported on disjoint sets of dual basis vectors.
			Any two of these three subspaces therefore span $V$. Deleting one
			member of the full family leaves at least two of the three. This is 3-SAT
			\cite{GareyJohnson1979}.
			
			Given such $H_c$, let $B=\{1,\dots,s,*\}$ with $A_c=H_c$ and $A_*=k$.
			Let $W_B=\ker\sigma\oplus0$, where $\sigma:\bigoplus_cH_c\to V$ is the
			sum map. For $T=\{1,\dots,s\}$ this gives $A_T/W_T\cong V$. It also
			gives $d_c=\operatorname{codim}\sum_{c'\ne c}H_{c'}=0$, $d_*=1$ and
			$t=1$. For $S_c=T\setminus\{c\}$ the map $A_T/W_T\to A_{S_c}/W_{S_c}$
			is the quotient map $V\to V/H_c$, so $h(S_c)=1$. With
			$\mathcal S=\{S_c\}$ one has $h_0=1$ and $T_0=T$, and (a) is vacuous.
			Condition (b) asks for a line meeting each $H_c$ only in zero,
			equivalently a nonzero vector outside their union. Bases of the
			$H_c$, the kernel of the sum map, and the list of sets $S_c$ can all
			be produced in polynomial time by Gaussian elimination.
			
			\item \emph{Realization by central extensions.} Let $A=A_B$, let
			$w_1,\ldots,w_r$ be a basis of the prescribed $W_B$, and let
			$C=(\mathbb Z/p^2\mathbb Z)^r$ with standard generators $e_\alpha$.
			Put $Q=C/pC$ and define the finite abelian groups
			\[
			G=(A\oplus C)/\langle(-w_\alpha,pe_\alpha):1\le\alpha\le r\rangle,
			\quad
			E_i=(A_i\oplus C)/\langle(-\operatorname{pr}_i w_\alpha,pe_\alpha):
			1\le\alpha\le r\rangle.
			\]
			The displayed relation subgroups project isomorphically onto $pC$,
			so $A$ embeds in $G$ and each $A_i$ embeds in $E_i$. Their quotients
			are $Q$. The coordinate projections define a homomorphism from $G$
			to the fiber product of the $E_i$ over $Q$. It induces the identity on
			$Q$ and on the common kernel $A$, hence is injective; both groups
			have order $|Q|\,|A|$, so it is an isomorphism. The $p$-multiples of
			the classes of the $e_\alpha$ are exactly $w_\alpha$, and $pA=0$.
			Consequently $A\cap P_p(G)=pG=W_B$ in additive notation. This gives
			a polynomial-size marked abelian presentation, including the case
			$r=0$. The NP-completeness assertion concerns the matrix input stated
			above, not explicitly enumerated group tables.\qedhere
		\end{enumerate}
	\end{proof}
	
	\subsection{Conditional minima for chart quotient readouts}
	\label{subsec:partial-chart-reads}
	
	A readout may report a quotient of the chart while the source must
	still cover every state of the original chart. Fix a subspace $N_i\le A_i$
	for each central extension in \cref{thm:min}. The readout $E_i\to\bar E_i=E_i/N_i$ still determines the
	common value in $Q$. The private kernel of $\bar E_i\to Q$ is
	$\bar A_i=A_i/N_i$. Both observed and
	requested label sets in this subsection refer to these fixed readouts.
	Write $\pi:A\to\bar A=\bigoplus_i\bar A_i$ for the coordinate
	quotient and $\bar W=\pi(W)$. For $X\subseteq B$ put
	\[
	\bar W_X=\operatorname{pr}_X\bar W,\quad
	\bar r(X)=\dim\bar W_X,\quad
	\bar D_i=\bar A_i/\bar W_{\{i\}},\quad
	\bar d_i=\dim\bar D_i,\quad
	\bar u(X)=\max_{i\in X}\bar d_i
	\]
	with $\bar u(\varnothing)=0$. On an atom $B$ set
	$t=\max_{i\in B}d_i$, $\bar t=\bar u(B)$,
	$\bar c(S)=\dim\bar A_S-\bar r(S)$ and
	$\bar h(S)=\min\{\bar t,\bar c(S)\}$.
	Thus $t$ measures the residual dimension needed for original chart coverage and $\bar t$ the dimension needed for readout coverage.

	\begin{theorem}[Sharp conditional ranges for quotient readouts]
		\label{thm:partial-chart-reads}
		\leavevmode
		\begin{enumerate}[label=\textup{(\roman*)},leftmargin=*]
			\item \emph{Fixed observation.}
			Among all chart-surjective sources exact at the protected cuts, the least conditional range of the quotient readouts on $J$ given the readouts on $P$ and the common boundary is $p^b$, where, writing $S=P\cap B$ and $U=(J\setminus P)\cap B$
			in each summand,
			\[
			b=\sum_{B\in\Pi}
			\bigl(
			\bar r(S\cup U)-\bar r(S)
			+\max\{0,\bar u(U)-\bar h(S)\}
			\bigr).
			\]
			For fixed $P$, a source of the original minimum order \eqref{eq:minimum} attains this bound for every $J$.
			
			\item \emph{Simultaneous attainment.}
			For a nonempty family $\mathcal P$ of observed sets, the following conditions are equivalent.
			\begin{enumerate}[label=\textup{(\alph*)},leftmargin=*]
				\item \emph{Attainment.}
				Some chart-surjective source exact at the protected cuts attains the bound in \textup{(i)} for every $P\in\mathcal P$ and every $J$.
				
				\item \emph{Attainment at minimum order.}
				A single source of the original minimum order attains all these bounds.
				
				\item \emph{Visible-space criterion.}
				On every atom, the visible spaces $\bar A_i,\bar W_B$ satisfy the two conditions of \cref{thm:incomparable-test} for
				$
				\mathcal S=\{P\cap B:P\in\mathcal P\}.
				$
			\end{enumerate}
		\end{enumerate}
	\end{theorem}

	\begin{proof}
		\leavevmode
		\begin{enumerate}[label=\textup{(\roman*)},leftmargin=*]
			\item \emph{Fixed observation.}
			
			\emph{Lower bound.}
			The full matching group $G$ maps onto the full matching group
			of the $\bar E_i$. Indeed, lifts of a matching tuple retain
			its common boundary value. The map carries $P_p(G)$ onto the
			subgroup generated by commutators and $p$th powers of the
			quotient matching group, and an element that maps into its
			kernel over $Q$ already lies in $A$. Thus $\bar W$ is its
			forced subspace. Projection preserves chart coverage and
			atomwise recombination. Applying \cref{thm:conditional} to
			the visible charts gives the stated lower bound for every
			source.
			
			\smallskip
			\emph{Attainment at the original minimum order.}
			We first determine which visible kernels can be obtained
			at minimum original source order. On a single atom, let
			$L\leq\bar A_B/\bar W_B$. A minimum original kernel $K_B$
			has visible residual image
			$
			L=\pi(K_B)/\bar W_B
			$
			if and only if $L$ maps onto every $\bar D_i$ and $\dim L\leq t$.
			
			Necessity follows from coverage and $\dim(K_B/W_B)=t$.
			For sufficiency, lift a surjection
			$k^t\to L$ linearly to $j':k^t\to\bar A_B$, and lift its
			coordinates further to maps $j_i:k^t\to A_i$. Their induced
			maps to $\bar D_i$ are onto. Let $V_i$ be the image of
			$N_i$ in $A_i/W_{\{i\}}$, and let $F_i$ be the kernel of
			the composite $k^t\to\bar D_i$. Then
			\[
			\dim F_i=t-\bar d_i
			\geq d_i-\bar d_i
			=\dim V_i.
			\]
			Fix a surjection $F_i\to V_i$. On a basis of $F_i$, adjust
			each $j_i$ value by a vector of $N_i$ so that its class
			modulo $W_{\{i\}}$ is the prescribed value in $V_i$.
			Leave $j_i$ unchanged on a complementary basis. The new
			map onto $A_i/W_{\{i\}}$ is surjective and still induces
			the prescribed map to $\bar A_i$.
			
			Set
			$
			K_B=W_B+(j_i)_{i\in B}(k^t).
			$
			It covers each chart coordinate, and at an index with
			$d_i=t$ the induced map onto $A_i/W_{\{i\}}$ is an
			isomorphism. Thus $\dim(K_B/W_B)=t$, and its visible image
			modulo $\bar W_B$ is exactly $L$. This proves the lifting
			criterion.
			
			Apply this criterion on each atom to the visible
			minimum-order kernel that attains the fixed-$P$ bound
			simultaneously for all $J$, supplied by
			\cref{thm:conditional}. By \cref{thm:min}, the resulting
			source has the original minimum order and attains all
			the bounds required in \textup{(i)}.
			
			\item \emph{Simultaneous attainment.}
			
			\emph{\textup{(c)} $\implies$ \textup{(b)}
				$\implies$ \textup{(a)}.}
			Condition \textup{(c)} yields visible minimum-order
			kernels attaining all the required bounds, by
			\cref{thm:incomparable-test}. The lifting criterion
			established in \textup{(i)} gives \textup{(b)}, which
			immediately implies \textup{(a)}.
			
			\smallskip
			\emph{\textup{(a)} $\implies$ \textup{(c)}.}
			Take any simultaneous optimizer and fix an atom.
			Write
			\[
			L=\pi(K_B)/\bar W_B,\qquad s=\dim L,
			\]
			and let $y(X)$ be the dimension of the image of $L$
			in $\bar A_X/\bar W_X$. Coverage gives $s\geq\bar t$.
			The residual part of the range exponent on this atom is
			$y(S\cup U)-y(S)$. Requests supported on one atom test
			its bound individually.
			
			If $s=\bar t$, the visible kernel has minimum order
			on this atom, so \cref{thm:incomparable-test} gives
			\textup{(c)}.
			
			Suppose $s>\bar t$.
			If $\bar c(S)<\bar t$ for some $S\in\mathcal S$, then
			$
			\bar u(S)\leq\bar c(S)<\bar t,
			$
			so a chart of deficit $\bar t$ lies outside $S$.
			For $U=B\setminus S$, the optimal residual increment
			gives
			$
			y(S)=s-\bar t+\bar c(S)>\bar c(S),
			$
			exceeding the dimension of $\bar A_S/\bar W_S$.
			Thus $\bar h(S)=\bar t$ for every $S\in\mathcal S$,
			and every sharp residual increment is zero. Taking
			$U=B\setminus S$ shows that each projection
			$
			L\longrightarrow\bar A_S/\bar W_S
			$
			is injective.
			
			We now construct a visible optimizer of dimension $\bar t$ modulo $\bar W_B$.
			Choose any $\bar t$-dimensional subspace $F\leq L$.
			Choose a complement to $\bar W_{\{i\}}$ in each $\bar A_i$
			and identify it with $\bar D_i$.
			For every $X\subseteq B$, these choices give compatible decompositions
			\[
			\bar A_X/\bar W_X
			=
			\left(\bigoplus_{i\in X}\bar W_{\{i\}}\right)\!/\bar W_X
			\ \oplus\!
			\bigoplus_{i\in X}\bar D_i.
			\]
			Keep the first component of the inclusion $F\to\bar A_B/\bar W_B$ fixed.
			For each component $a_i:F\to\bar D_i$, choose a surjection
			$a_i':F\to\bar D_i$ such that
			$
			\ker a_i'\subseteq\ker a_i.
			$
			To construct it, keep $a_i$ unchanged on a complement of its kernel and map its kernel onto a complement of its image.
			This is possible because $\dim F=\bar t\geq\dim\bar D_i$.
			
			For each $S\in\mathcal S$, any vector sent to zero
			by the modified map to $\bar A_S/\bar W_S$ was also
			sent to zero before modification: the first component
			is unchanged, and every coordinate kernel in the
			second component has only decreased. The original
			projection was injective on $F$, so the modified
			projection is injective as well.
			
			Since $\mathcal S$ is nonempty, the modified map on
			$B$ is also injective. Its image maps onto every
			$\bar D_i$, has dimension $\bar t$, and injects into
			the observation space at every $S\in\mathcal S$.
			Its full preimage in $\bar A_B$ is therefore a visible
			minimum-order kernel attaining every required zero
			residual increment. Applying
			\cref{thm:incomparable-test} proves \textup{(c)}.
			This argument includes $\bar t=0$, when $F=0$.
			\qedhere
		\end{enumerate}
	\end{proof}

	\paragraph{Different quotient maps in the premise and conclusion.}
	A query can also observe a quotient of a chart and request the full
	chart, or use different quotients of the same chart in different
	queries. The following finite count tests such queries together and
	distinguishes minimization at a fixed source order from minimization
	over all orders.
	
	Fix a finite list of queries. In each query $j$ let $O_j$ be the tuple
	of input readouts and $F_j$ the tuple consisting of those inputs and all
	requested outputs. Each readout is a chart quotient $E_i\to E_i/N$ with
	$N\le A_i$. Taking $N=0$ gives the full chart. Condition on an
	attained common boundary value and attained input values.
	On an atom $B$, each tuple induces a linear map $P$ on $A_B$ by taking
	the indicated coordinates in $B$ modulo their respective $N$. Include
	the full individual chart projections in the list of maps as coverage
	tests. Put $V=A_B/W_B$, $n=\dim V$ and, for each listed $P$, put
	\[
	r_P=\dim P(W_B),\qquad
	L_P=\ker\bigl(V\longrightarrow P(A_B)/P(W_B)\bigr).
	\]
	The map for an empty tuple is zero. A list may contain equal kernels.

	\begin{corollary}[A sharp range for different input and output readouts]
		\label{cor:mixed-query-minimum}
		Fix one of the queries and one atom. Write $F$ for the joint map on
		this atom and $O$ for its input part. Let $F_i,O_i$ denote their
		coordinate maps on $A_i$, taking the zero map at an unused chart. Set
		\[
		t_F=\max_{i\in B}\bigl(\dim F_i(A_i)-\dim F_i(W_{\{i\}})\bigr),
		\qquad c_O=\dim O(A_B)-r_O .
		\]
		The least contribution of this atom to the conditional range exponent is
		\[
		r_F-r_O+\max(0,t_F-c_O).
		\]
		Adding over atoms gives the sharp exponent for the query. A source of
		the minimum original order attains it.
	\end{corollary}

	\begin{proof}
		Work above a fixed boundary value. For an admissible kernel $K_B$, put
		$s_F=\dim F(K_B)-r_F$ and $s_O=\dim O(K_B)-r_O$. Coverage gives
		$s_F\ge t_F$, and projection gives $s_O\le \min(s_F,c_O)$. The
		conditional exponent is $r_F-r_O+s_F-s_O$, proving the bound.
		
		We attain it in two lifts. Take the visible coordinate spaces to be
		$O_i(A_i)$, with forced subspace $O(W_B)$. Write $d_i^O$ for the
		codimension of $O_i(W_{\{i\}})$ in $O_i(A_i)$ and set $h=\min(t_F,c_O)$.
		Each $d_i^O\le h$, since both the combined residual at $i$ and
		$O(A_B)/O(W_B)$ map onto the observed residual at $i$.
		There is an $h$-plane in $O(A_B)/O(W_B)$ covering every observed
		residual coordinate. To construct it, take the image of simultaneous
		surjections from a space of dimension $\max_i d_i^O$ onto those
		coordinates. One has full rank, so their joint image has that
		dimension. Lift the image along the surjection from $O(A_B)/O(W_B)$
		onto the product of the residual coordinates and extend it to an
		$h$-plane.
		
		Apply the coordinate lifting argument of \cref{thm:partial-chart-reads}
		first to $F_i(A_i)\to O_i(A_i)$, with $t_F$ parameters. It gives a
		covering combined kernel of residual dimension $t_F$ whose observed
		image has residual dimension $h$. Apply the same argument to
		$A_i\to F_i(A_i)$, now with $\max_i d_i$ parameters. This realizes the
		chosen combined image at minimum original kernel dimension. These maps
		are chart quotients on each boundary fiber: their kernels are contained
		in $A_i$, and the joint quotient for several readouts has the
		intersection of their kernels. Taking the direct sum of the resulting
		kernels and applying \cref{thm:min} completes the construction.
	\end{proof}
	
	\begin{proposition}[Intersection counts for mixed quotient queries]
		\label{prop:mixed-quotient-count}
		\leavevmode
		\begin{enumerate}[label=\textup{(\roman*)},leftmargin=*]
			\item \emph{Intersection counts.} For $0\le s\le n$ and a tuple of integers
			$0\le a_P\le\min(s,\dim L_P)$, define
			\[
			M_s(a)=
			\sum_{\substack{U_P\le L_P\ (\text{all }P)\\ \dim U_P=a_P}}
			\gaussian{n-\dim(\sum_P U_P)}{s-\dim(\sum_P U_P)} .
			\]
			Here \(\gaussian{v}{u}\) denotes the number of $u$-planes in
			$\mathbb F_p^v$ and is zero for $u<0$ or $u>v$. Starting with the
			larger tuples in the componentwise order, compute the integers
			\[
			N_s(a)=M_s(a)-
			\sum_{x\ge a,\,x\ne a} N_s(x)\prod_P\gaussian{x_P}{a_P} .
			\]
			
			\item \emph{Attainable exponents on an atom.} Retain the pairs $(s,x)$ with $N_s(x)>0$ for which
			$s-x_{P_i}=d_i$ at each full chart projection $P_i$, $i\in B$.
			These pairs give exactly the attainable vectors of conditional range
			exponents on $B$. Their entry at query $j$ is
			\[
			r_{F_j}-r_{O_j}+x_{O_j}-x_{F_j}.
			\]
			
			\item \emph{Simultaneous minima across atoms.} Choose one retained pair per atom and add the entry vectors. The
			coordinatewise minima of the resulting vectors give the sharp range
			exponents, and a single source attains them all exactly when their
			vector occurs. Restricting to $s=\max_{i\in B}d_i$ on each atom gives
			the same tests at minimum source order. To test whether all unrestricted
			minima are attained at minimum order, compare their vector with the
			vectors from this restriction.
		\end{enumerate}
	\end{proposition}
	
	\begin{proof}
		\leavevmode
		\begin{enumerate}[label=\textup{(\roman*)},leftmargin=*]
			\item \emph{Intersection counts.} For each $s$-plane $Y\le V$ let $x_P=\dim(Y\cap L_P)$. Count the
			tuples $(Y,(U_P)_P)$ with $U_P\le Y\cap L_P$ and $\dim U_P=a_P$.
			Fixing the $U_P$ first means choosing $Y$ through their sum. The
			quotient by that sum gives the Gaussian coefficient in $M_s(a)$.
			Fixing $Y$ first gives $\prod_P\gaussian{x_P}{a_P}$ choices.
			The coefficient at $x=a$ is one. Descending induction thus shows that
			$N_s(x)$ counts precisely the planes with intersection dimensions $x$.
			The coefficients can be computed by counting ordered independent
			vectors and dividing by the number of ordered bases of a plane.
			
			\item \emph{Attainable exponents on an atom.} The inverse image $K_B$ of $Y$ in $A_B$ has
			$\dim P(K_B)=r_P+s-x_P$. Applied to the full chart projection, this
			identity gives precisely the stated coverage test. In a fiber above
			the common boundary, fixing the inputs leaves a translate of their
			kernel on $K_B$. The dimension of its output image is the rank of the
			combined map $F_j$ minus the rank of $O_j$, which gives the displayed
			entry.
			
			\item \emph{Simultaneous minima across atoms.} Outputs and inputs are labeled chart coordinates, so these
			ranks add over atoms even if a chart has several readouts. Finally,
			\cref{thm:min} realizes every such direct sum of covering kernels
			and gives its order. This proves both attainment claims.\qedhere
		\end{enumerate}
	\end{proof}
	
	\subsection{Realizing representable rank functions by nonabelian charts}
	
	The following construction realizes subspace rank functions as cut costs of directly indecomposable charts with a common nonabelian quotient.
	Let $p$ be an odd prime, $k=\mathbb F_p$, and $I$ a finite nonempty
	set of chart labels. A subspace representation consists of finite-dimensional $k$-spaces
	$T_i\le T$ with $T=\sum_{i\in I}T_i$. Its rank function is
	\[
	r(S)=\dim_k\sum_{i\in S}T_i\qquad(S\subseteq I).
	\]
	Write $\Pi$ for the cut atoms. Use the complete chart graph, so every nontrivial bipartition has
	connected sides. When $|I|=1$ the single chart carries its designated
	quotient and there are no nontrivial cuts.
	
	\begin{theorem}[Realization of representable polymatroid ranks by group-atlas source minima]
		\label{thm:realization}
		Every such representation admits a finite class-two exponent-$p$ atlas
		with a common nonabelian quotient $Q$ on the edges, and a
		chart-surjective subgroup $H_0$ of its full matching group $G$, with
		the following properties.
		\begin{enumerate}[label=\textup{(\roman*)},leftmargin=*]
			\item \emph{Structure of the charts and seed.} Each chart and the group $H_0$ are directly indecomposable. The
			kernel of each chart map to $Q$ lies in the derived subgroup of the
			chart. The map has no group section whenever $r(\{i\})>0$.
			\item \emph{Minimum order and equality cases.} Among all chart-surjective subgroups of $G$ that recombine
			exactly at every protected cut, the least possible order is
			\begin{equation}\label{eq:realized-rank}
				|H_0|\,p^{\delta(\Pi)},\qquad
				\delta(\Pi)=\sum_{B\in\Pi} r(B)-r(I).
			\end{equation}
			A unique minimizer contains $H_0$. All minimizers are carried into
			one another by automorphisms of the charts that induce the identity
			on every edge.
		\end{enumerate}
	\end{theorem}
	
	\begin{proof}
		Put $P=k e_0\oplus T^*\oplus k a$ and $V=P\oplus P^*$. The extra
		line $ka$ ensures $\dim V\ge4$, including for the zero rank function.
		Write $\omega((x,\xi),(y,\eta))=\eta(x)-\xi(y)$ for the standard
		symplectic form. Extend the coordinate $e_0^*$ and each
		$\phi\in T=(T^*)^*$ to covectors on $V$ by zero on the other summands.
		Define $b:\bigwedge^2 V\to T^*$ by
		\[
		b(v,w)(\phi)=(e_0^*\wedge\phi)(v,w).
		\]
		The joint map $(\omega,b)$ surjects onto $k\oplus T^*$. The pair
		$(e_0,e_0^*)$ gives $(1,0)$, while the pairs $(e_0,x)$ for
		$x\in T^*\subseteq P$ give $(0,x)$.
		
		For an alternating map $\gamma:\bigwedge^2 V\to C$ use the group
		on $V\oplus C$ with product
		$(v,c)(w,d)=(v+w,c+d+\tfrac12\gamma(v,w))$.
		Take $Q$ from $\gamma=\omega$, and take the chart $E_i$ from
		$\gamma_i=(\omega,b|_{T_i})$ with center $k\oplus T_i^*$, where
		$b|_{T_i}$ restricts each value of $b$ to $T_i$. Every incidence
		map $E_i\to Q$ forgets $T_i^*$. Thus the matching group has underlying
		space $V\oplus k\oplus J$ with $J=\bigoplus_i T_i^*$. Its commutator
		has image $k\oplus W$, where
		\[
		W=\{(x|_{T_i})_i:x\in T^*\}\le J,\qquad
		\dim \operatorname{pr}_S W=r(S).
		\]
		Indeed the kernel of the displayed restriction to coordinates in
		$S$ is the annihilator in $T^*$ of $\sum_{i\in S}T_i$.
		Set $H_0=V\oplus k\oplus W$. Surjectivity of $T^*\to T_i^*$ gives
		chart coverage.
		\begin{enumerate}[label=\textup{(\roman*)},leftmargin=*]
			\item \emph{Structure of the charts and seed.} Each $E_i$ has derived subgroup $k\oplus T_i^*$,
			so its edge kernel is contained there. For $T_i\ne0$, choose
			$x\in T^*$ with $x|_{T_i}\ne0$. The primary vectors $e_0,x$ commute
			in $Q$, while any lifts in $E_i$ have commutator $(0,x|_{T_i})$.
			This excludes a group section.
			
			For both $E_i$ and $H_0$ the center equals the derived subgroup,
			since $\omega$ has zero radical. The scalar forms of their
			commutator maps are respectively the families
			\[
			c\omega+e_0^*\wedge\phi,\qquad
			c\in k,\quad \phi\in T_i\ \text{or}\ T.
			\]
			If $c\ne0$, testing a radical vector first against $P^*$ forces its
			$P$-coordinate to vanish, and testing against $P$ then kills its
			$P^*$-coordinate. If $c=0$ the rank is at most two. Thus the
			singular functionals form the proper hyperplane $c=0$.
			Suppose one of these groups were $X\times Y$ with both factors
			nontrivial. The equality of center and derived subgroup forces
			both factors to be nonabelian. Every functional on $X'\oplus Y'$
			annihilating $X'$ gives a singular form, since its radical contains
			the nonzero primary subspace $X/Z(X)$. Functionals annihilating
			$Y'$ are singular for the same reason. These two annihilators span
			the entire dual, contradicting containment in a proper hyperplane.
			This proves direct indecomposability.
			
			\item \emph{Minimum order and equality cases.} Let $F\le G$ cover every chart. Its primary projection is all of $V$,
			so $[F,F]=k\oplus W$. Odd characteristic makes $F$ an additive
			subspace in these coordinates. Indeed products corrected by half
			their commutators give sums, and powers give scalar multiples.
			In particular $F$ contains the common central line. On putting
			$D=F\cap J$, every such $F$ has the unique description
			\begin{equation}\label{eq:realized-twists}
				F_{D,h}=\{(v,c,z): z+D=h(v)\},\qquad
				W\le D\le J,\quad h\in\operatorname{Hom}_k(V,J/D).
			\end{equation}
			Conversely each displayed subgroup covers every chart, since
			$\operatorname{pr}_i W=T_i^*$.
			
			By \cref{thm:min}, the complete cuts require
			$D=\bigoplus_{B\in\Pi}\operatorname{pr}_B D$. Since $W$ already
			projects onto every private chart space, the least admissible $D$ is
			\[
			D_*=\bigoplus_{B\in\Pi}\operatorname{pr}_B W.
			\]
			It has dimension $\sum_{B\in\Pi}r(B)$, and $|F_{D,h}|=|Q|p^{\dim D}$.
			This proves \eqref{eq:realized-rank}. The minimizers are
			$F_{D_*,h}$, of which precisely $F_{D_*,0}$ contains $H_0$.
			A lift $\widetilde h:V\to J$ of $h$ acts on chart $i$ by
			$(v,c,z_i)\mapsto(v,c,z_i+\operatorname{pr}_i\widetilde h(v))$.
			These are group automorphisms fixing each edge value, and jointly
			they carry $F_{D_*,0}$ to $F_{D_*,h}$.\qedhere
		\end{enumerate}
	\end{proof}
	
	\noindent\textit{An abelian realization of the cut costs.}
	Let $p$ be any prime, put $d=\dim T$, and write
	$\widetilde H=(\mathbb Z/p^2\mathbb Z)^d$ additively.
	Identify $p\widetilde H$ with $T^*$ and take
	$\widetilde Q=\widetilde H/p\widetilde H$. For each $i$, let
	$N_i\le p\widetilde H$ be the annihilator of $T_i$ and form the
	chart $\widetilde H/N_i\to\widetilde Q$, whose kernel is $T_i^*$.
	The intersection of the $N_i$ is zero because the $T_i$ span $T$.
	Thus $\widetilde H$ embeds diagonally as a chart-surjective seed in
	the matching group $\widetilde G$. With $A=\bigoplus_i T_i^*$ we have
	$\widetilde G=\widetilde H+A$ and $pA=0$. Since $\widetilde G$ is
	abelian, its forced space in \eqref{eq:support-atoms} is
	$p\widetilde G=p\widetilde H\le A$. This is the image of the
	restriction map $T^*\to\bigoplus_i T_i^*$. Its projection on $B$
	has dimension $r(B)$ and every coordinate deficit vanishes. Hence
	\cref{thm:min} gives
	\[
	\frac{\min_F |F|}{|\widetilde H|}
	= p^{\sum_{B\in\Pi}r(B)-r(I)} ,
	\]
	where $F$ ranges over the chart-surjective subgroups exact at the
	protected cuts. Thus abelian charts realize the same relative costs as those in \cref{thm:realization}.
	
	\subsection{Recovery and the equality cases for scalar charts}
	
	On a fixed boundary fiber, whole-chart determination of scalar private
	coordinates is linear dependence. In \cref{thm:realization},
	each complement of $W$ in $J$ is a central direct factor of the full matching group, with complementary factor $H_0$.
	
	Suppose now that the $T_i$ are nonzero lines. Choose their generators
	$\lambda_i$ and identify each $T_i^*$ with $k$ by evaluation at
	$\lambda_i$. Write $M$ for the represented loopless matroid on $I$.
	For a subgroup $F$ write $\operatorname{cl}_F(S)$ for the labels whose
	whole-chart values are determined by the chart values on $S$.
	Denote the direct sum of the matroid restrictions $M|B$ by
	$M_\Pi=\bigoplus_{B\in\Pi}M|B$. Matroid representations by partitions express rank through joint
	observations \cite{Matus1999}. Finite groups and their subgroups
	also give information functions through coset observations
	\cite{ChanYeung2002}. Chen, Cheng and Bai use arrays of joint
	observations to characterize entropy functions proportional to the
	rank of a regular matroid \cite{ChenChengBai2025}. Here the charts
	and their quotient maps are fixed. We ask which subgroups supplying those observations attain
	the largest possible recovery closure.
	
	\begin{corollary}[Sharp tradeoff between recombination and recovery]
		\label{cor:recovery}
		\leavevmode
		\begin{enumerate}[label=\textup{(\roman*)},leftmargin=*]
			\item \emph{Recovery closures.} The seed satisfies $\operatorname{cl}_{H_0}=\operatorname{cl}_M$.
			For every chart-surjective, protected-cut-exact subgroup $F$,
			\[
			\operatorname{cl}_F(S)\subseteq\operatorname{cl}_{M_\Pi}(S)
			\qquad(S\subseteq I).
			\]
			
			\item \emph{Equality and compatibility with the cuts.} Equality for all $S$ holds precisely for the order minimizers in
			\cref{thm:realization}. In particular a chart-surjective subgroup
			realizes both the seed's exact whole-chart recovery and every protected
			cut if and only if every circuit of $M$ lies in one block of $\Pi$.
		\end{enumerate}
	\end{corollary}
	
	\begin{proof}
		\leavevmode
		\begin{enumerate}[label=\textup{(\roman*)},leftmargin=*]
			\item \emph{Recovery closures.} A nonempty set of whole charts reads the common value $(v,c)$.
			Within its fiber, differences of private coordinates range over
			$D$ from \eqref{eq:realized-twists}. Thus $S$ determines $i$
			exactly when the $i$th coordinate functional on $D$ belongs to the
			span of the coordinate functionals indexed by $S$. The empty set
			determines no chart, and the coordinate matroid of $D$ has no loops
			because $D\ge W$ projects onto every coordinate. For $D=W$ these
			functionals represent $M$. For $D=D_*$ the coordinates on each block
			represent $M|B$ and the different block spaces are independent, giving
			$M_\Pi$. Every other admissible $D$ contains $D_*$. Each determination
			on $D$ restricts to one on $D_*$.
			
			\item \emph{Equality and compatibility with the cuts.} If the two closures coincide, the
			two coordinate matroids have equal total rank. Their ranks are
			$\dim D$ and $\dim D_*$, forcing $D=D_*$. \Cref{cor:circuit} gives the circuit criterion for $H_0$.
			For a crossing circuit $C$ and $i\in C$, the set $C\setminus\{i\}$
			spans $i$ in $M$ and fails to span it in $M_\Pi$. The closure
			inclusion above transfers this failure to every admissible $F$.\qedhere
		\end{enumerate}
	\end{proof}
	
	For scalar charts $r(\{i\})=1$ for every label. Applying
	\cref{cor:rank} to the realized rank function gives
	\[
	r(S)=|S|+\delta(\Pi_S)-\delta(\Pi_{\rm disc})
	\qquad(S\ne\varnothing).
	\]
	Thus the increases for families of singleton cuts already determine
	the labeled matroid.
	The increase relative to the unprotected minimum at a single cut can agree for distinct scalar examples. On three
	labels take three identical nonzero columns for $U_{1,3}$, and take
	$(1,0),(0,1),(1,1)$ for $U_{2,3}$. Every nontrivial bipartition has
	exponent one for both matroids. Protecting two distinct singleton cuts
	gives the singleton partition, with exponents two and one, respectively.
	
	For $U_{2,3}$ the private readings in each fiber of the constructed $H_0\to Q$ are the triples $(a,b,a+b)$ from \cref{sec:recovery}. The same failure of determination now occurs with directly indecomposable chart groups and a common nonabelian quotient. More generally, for an individual requested chart with a scalar private reading the term $r(S_B\cup U_B)-r(S_B)$ of \eqref{eq:conditional} tests whether the readings in its own atom span it. With general central kernels the observed charts also have the capacity $a_B$ to convey private information beyond the forced subspace $W_{S_B}$.
	
	\section{Limits and extensions of the common-boundary criteria}
	\label{sec:limits}
	Two changes to the common-boundary problem affect how source minima
	are calculated. For a semisimple action on the kernels,
	\cref{thm:semisimple} gives the minimum one simple constituent at
	a time. Prescribing a second quotient and the kernel of the source
	over that quotient permits every downward-closed family of feasible
	cuts, even on a cycle of identical charts
	(\cref{thm:complex}). \Cref{thm:fixed-kernel-span-cost} computes
	the chart cost of this construction through a linear span program.
	
	\subsection{Semisimple kernels with nontrivial action}
	
	Keep the common quotient $Q$, chart coverage, and complete cuts of
	\cref{thm:min}. Allow each chart kernel $A_i$ to be an elementary abelian normal subgroup on which $Q$ acts by conjugation.
	Put $k=\mathbb F_p$ and suppose that the $kQ$-modules $A_i$ are
	semisimple. For each occurring simple module $S$, its endomorphism
	ring $D_S=\operatorname{End}_{kQ}(S)$ is a division ring by
	Schur's lemma and hence a field by Wedderburn's finite division ring
	theorem \cite{Lam2001Rings}. Using
	$M_{i,S}=\operatorname{Hom}_{kQ}(S,A_i)$, evaluation identifies
	$A_i$ with $\bigoplus_S S\otimes_{D_S}M_{i,S}$.
	Write $M_S=\bigoplus_i M_{i,S}$. The $S$-component of the extension
	class of the matching group belongs to
	$H^2(Q,S)\otimes_{D_S}M_S$. Let $W_S\le M_S$ be the image of
	this tensor under contraction by the $D_S$-dual of $H^2(Q,S)$.
	
	\begin{theorem}[Minimum source order with semisimple kernels]\label{thm:semisimple}
		For each atom $B\in\Pi$, let $r_S(B)=\dim_{D_S}\operatorname{pr}_B W_S$
		and $d_{i,S}=\dim_{D_S} M_{i,S}-r_S(\{i\})$.
		The least order of a chart-surjective subgroup exact at all protected
		cuts is $|Q|p^m$, where
		\[
		m=\sum_S (\dim_k S)\sum_{B\in\Pi}
		\left(r_S(B)+\max_{i\in B}d_{i,S}\right).
		\]
	\end{theorem}
	\begin{proof}
		Let $G$ be the matching group, $A=\bigoplus_iA_i$, and $H\le G$
		map onto $Q$. Since $A$ is abelian, conjugation by $H$ induces the
		given $Q$-action on $A$. Thus $K=H\cap A$ is a $kQ$-submodule
		and is normal in $G$. Conversely, a prescribed submodule $K$ is
		realized by such an $H$ precisely when $G/K\to Q$ splits, by
		taking the inverse image of a section. Write
		$K=\bigoplus_S S\otimes_{D_S} L_S$. Cohomology commutes with finite
		direct sums, and passage to $A/K$ acts on the $S$-component by
		$M_S\to M_S/L_S$. The quotient extension splits precisely when
		$W_S\le L_S$ for every $S$.
		
		The fibers of $H\to Q$ are translates of $K$, so exactness requires
		and suffices that $K$ factor over the cut atoms, and coverage requires
		and suffices that every coordinate projection of $K$ be onto.
		Both conditions respect the isotypic decomposition. On atom $B$, the
		corresponding multiplicity space must therefore contain
		$\operatorname{pr}_B W_S$ and project onto every $M_{i,S}$ for
		$i\in B$. The dimension bound and simultaneous supplement construction
		of Theorem~\ref{thm:min}, now over $D_S$, give the exact minimum
		$r_S(B)+\max_{i\in B} d_{i,S}$. A multiplicity of one contributes
		$\dim_k S$ to $\dim_k K$. Summing gives the claimed exponent.
	\end{proof}
	
	\begin{corollary}[Conditional range minima with semisimple kernels]
		\label{cor:semisimple-query-minima}
		Under the hypotheses of \cref{thm:semisimple}, let $b_S(P,J)$
		denote the expression in \eqref{eq:conditional} for the spaces
		$M_{i,S}$ and the forced subspace $W_S$, with all dimensions over
		$D_S$.
		\begin{enumerate}[label=\textup{(\roman*)},leftmargin=*]
			\item \emph{A fixed observation.} The smallest conditional range on $J$ given $P$ and $Q$ is
			\[
			p^{\sum_S (\dim_k S)b_S(P,J)}.
			\]
			For fixed $P$, one minimum-order source attains this range for every
			$J$.
			
			\item \emph{Simultaneous observations.} For any nonempty family $\mathcal P$ of observed sets, one
			minimum-order source works for all $P\in\mathcal P$ and all $J$ if and only
			if the criterion of \cref{thm:incomparable-test} holds in each atom
			of every multiplicity space. When $\mathcal P=2^I$, this reduces to
			\eqref{eq:all-observations} in each such space.
		\end{enumerate}
	\end{corollary}
	\begin{proof}
		\leavevmode
		\begin{enumerate}[label=\textup{(\roman*)},leftmargin=*]
			\item \emph{A fixed observation.} Fixing the boundary and the observed values leaves a translate of
			$K\cap\ker\operatorname{pr}_P$. Projection onto $J$ commutes with the
			isotypic decomposition. A multiplicity-space image of dimension $d$
			over $D_S$ contributes $d$ copies of $S$, hence $d\dim_k S$ to
			the dimension over $k$. The rank calculations and subspace
			constructions of
			\cref{thm:conditional,thm:incomparable-test} work over each
			finite field $D_S$. In particular, the fixed-$P$ choices can be made
			independently in the different summands and atoms. The splitting
			criterion in \cref{thm:semisimple} realizes their direct sum.
			
			\item \emph{Simultaneous observations.} Both source dimension and conditional exponents are sums of the
			individual contributions with strictly positive weights. Equality
			with a sum of lower bounds holds exactly when every summand attains
			its own bound. This gives the simultaneous criteria.\qedhere
		\end{enumerate}
	\end{proof}
	
	\paragraph{Failure of least support for a nonsplit coefficient module.}
	Take $Q=C_p=\langle g\rangle$ and $M=kv\oplus kQ$, with $v$ fixed.
	Let $N=1+g+\cdots+g^{p-1}$ and $a=N\cdot1$ in the regular summand.
	The cyclic group cohomology formula $H^2(Q,M)=M^Q/NM$ \cite{BrownCohomology1982} gives a
	nonzero class $[v]$. Both $kv$ and $k(v+a)$ are invariant lines.
	Quotienting by the first kills $v$, and quotienting by the second
	makes $v=-a$ a norm. Thus both quotients kill the class while
	quotienting by the intersection of their kernels leaves it nonzero.
	In particular the coefficient submodules whose quotients split the
	extension can have incomparable minimal elements.
	
	\subsection{Which families of actual cuts can one source satisfy?}
	
	Keep the atlas and the common boundary used for recombination fixed.
	Also fix a quotient map $q$ on the full matching group $G$ and a
	subgroup $K\le\ker q$. This map is additional to the chart boundary
	map. A candidate source $H$ must have $q(H)=q(G)$ and
	$H\cap\ker q=K$. The source may vary with the family of cuts. At each
	selected cut it must equal the fiber product of its two complete-side
	images over their common complete boundary.
	
	\begin{theorem}[Universality of fixed-kernel cut feasibility]
		\label{thm:complex}
		Fix an odd prime $p$. Let $\mathcal K\subseteq 2^{[m]}$ be any
		downward-closed family containing $\varnothing$. There exist a cycle of
		identical group charts of order $p^5$, identical nonabelian interfaces of
		order $p^3$, a quotient $q$ of its full matching group $G$, a subgroup
		$K\leq\ker q$, and $m$ actual cuts with both sides connected. These data have the following properties.
		\begin{enumerate}[label=\textup{(\roman*)},leftmargin=*]
			\item \emph{Group structure.} Both the charts and $G$ are directly indecomposable.
			The chart group, the interface group, and the incidence homomorphism depend only on $p$.
			
			\item \emph{Feasible cut families.} For every $J\subseteq[m]$, membership $J\in\mathcal K$ is equivalent
			to the existence of a subgroup $H\leq G$ exact at all cuts in $J$
			and satisfying
			\[
			q(H)=q(G),\qquad H\cap\ker q=K.
			\]
			
			\item \emph{Coverage and recovery.} Every subgroup satisfying the displayed quotient and kernel conditions
			covers every chart and every attainable complete boundary. Its
			whole-chart recovery closure is the identity, so a family of charts
			determines exactly its own members.
		\end{enumerate}
	\end{theorem}
	
	\begin{proof}
		Write $k=\mathbb F_p$. At each vertex take the class-two exponent-$p$
		group with primary coordinates $(s,y,t)$, central coordinates $(c,d)$,
		and basic brackets $[s,y]=c$ and $[s,t]=d$. Here and below the brackets
		specify the centered group law
		$(v,z)(v',z')=(v+v',z+z'+[v,v']/2)$. At both ends of every edge use
		the quotient that retains $(s,y,c)$. Its image $Q$ is the Heisenberg
		group of order $p^3$. Thus $G$ has common coordinates $(s,y,c)$ and
		independent private coordinates $(t_b,d_b)$ at each vertex $b$.
		Up to skew symmetry its nonzero basic brackets are $[s,y]=c$ and $[s,t_b]=d_b$.
		For any of the cuts below, the two crossing edges record the diagonal
		copy of $Q$.
		\begin{enumerate}[label=\textup{(\roman*)},leftmargin=*]
			\item \emph{Group structure.} In each chart and in $G$, the bracket on the primary quotient has the
			form $[(s,w),(s',w')]=sw'-s'w$ after identifying the $w$-space with
			the displayed central space. Its radical is zero and its image spans
			the center, so the center equals the derived subgroup. A vector with
			$s\ne0$ has primary centralizer exactly its own line. In a direct
			product decomposition the primary quotients of the factors are
			complementary and commute. One factor contains a primary vector with $s\ne0$, since the
			two primary quotients span the entire primary space. The other
			quotient lies in that vector's line and has zero intersection with
			the first quotient, so it is zero. A central
			direct factor must be trivial because the center equals the derived
			subgroup. This proves direct indecomposability.
			
			\item \emph{Feasible cut families.} Write $\mathcal F$ for the set of minimal nonfaces of $\mathcal K$.
			Arrange the cycle as a root $0$, an unprotected anchor $*$, and
			successive vertex blocks $B_1,\dots,B_m$. In $B_i$ put one vertex
			$(i,F)$ for each $F\in\mathcal F$ containing $i$. Give each empty
			block one dummy vertex. Add an unprotected dummy if needed to give
			the cycle at least three vertices. Cut $i$ separates the consecutive
			block $B_i$ from its complement. Write $r=t_0$, let $T$ be the space of
			private primary coordinates away from the root, and denote its standard
			basis by $e_b$. Define
			\[
			u_F=\sum_{i\in F}e_{(i,F)},\qquad
			U=\langle e_*,u_F\ (F\in\mathcal F)\rangle\leq T,\qquad
			\kappa(e_*)=\kappa(u_F)=1.
			\]
			These generators have disjoint nonempty supports, hence specify a
			linear functional $\kappa$ on $U$. On a matching tuple with nonroot
			private primary part $x\in T$, put $q=(s,y,x+U)\in k^2\oplus T/U$.
			This is a surjective homomorphism. Inside $\ker q$ prescribe $K$ by
			$r=\kappa(x)$, with all central coordinates free. Both $q$ and $K$
			are fixed before any cuts are selected.
			
			Let $Z=Z(G)$. Every candidate $H$ contains $Z$, since $Z\leq K$.
			An element of $H/Z$ with $(s,y,x)=0$ lies in $K/Z$ and has
			$r=\kappa(0)=0$, so the projection to $(s,y,x)$ is injective.
			For surjectivity, first lift $(s,y,x+U)$ using $q(H)=q(G)$ and
			then correct its $x$-coordinate by an element of $K/Z$. Consequently its equation, and that of
			every possible candidate, is
			\[
			r=f(s,y)+\lambda(x),\qquad f\in(k^2)^*,\quad
			\lambda\in T^*,\quad \lambda|_U=\kappa .
			\]
			Conversely the full central preimage of each such graph has the
			required quotient and kernel.
			
			Cut $i$ is exact precisely when $\lambda$ vanishes on all coordinates
			of $B_i$. Indeed, in that case the graph equation depends only on the
			complement and the shared boundary, so every boundary-matching pair
			of side states recombines. Otherwise take $x$ supported in $B_i$
			with $\lambda(x)\ne0$ and all shared coordinates zero. Its block side
			can be matched with the identity complement, but the resulting tuple
			violates the root equation. Thus feasibility at $J$ asks for
			coefficients $\lambda(e_b)$ vanishing on the selected blocks and
			satisfying
			\[
			\lambda(e_*)=1,\qquad
			\sum_{i\in F}\lambda(e_{(i,F)})=1\quad(F\in\mathcal F).
			\]
			If some $F\subseteq J$, its equation is impossible. Otherwise choose
			one $i\in F\setminus J$ per $F$, give its coefficient value one,
			and give the other nonanchor coefficients value zero. Disjoint supports
			make all equations hold simultaneously. This includes
			$\mathcal F=\varnothing$, and proves the equivalence.
			
			\item \emph{Coverage and recovery.} Every candidate covers each nonroot
			chart since $x$ is unrestricted, and the root chart since
			$\lambda(e_*)=1$. It covers the complete boundary $Q$.
			Varying $d_b$ alone changes precisely chart $b$, which proves free
			recovery on the whole-chart readouts.\qedhere
		\end{enumerate}
	\end{proof}
	
	Taking $\mathcal K$ to consist of all proper subsets of $[m]$ gives
	an obstruction of arbitrarily large minimal size. The transmitted
	commutator $[s,y]=c$ belongs to every actual interface and the groups
	are directly indecomposable. The feasibility equations descend to $G/Z(G)$.
	
	This construction uses at most
	$3+m+\sum_{F\in\mathcal F}|F|$ charts. Describing the prescribed
	kernel by the minimal nonfaces can therefore require exponentially
	many coordinates as $m$ grows. The local group order stays fixed.
	
	In \cref{thm:cyclic-recovery-compatibility-universality} the formation data vary.
	Here the incidence maps are fixed and the protected cuts vary. On $U$,
	a cut family $J$ is feasible exactly when $\kappa$ lies in the span of the
	coordinate functionals outside the selected blocks, restricted to $U$ and including
	the unprotected anchor. This is the acceptance criterion for a linear span program
	\cite{KarchmerWigderson1993}. The cuts are connected and preserve all whole-chart recovery witnesses.
	
	\paragraph{Chart cost in the cycle model.}
	Fix $p$ and $m\geq1$. Retain exactly one root,
	one anchor and $m$ nonempty consecutive blocks, the same chart and
	interface laws, and the same form of $q$ and $K$, but permit arbitrary
	$U\leq T$ containing $e_*$ and arbitrary $\kappa\in U^*$ with
	$\kappa(e_*)=1$. Count all vertices, including the root and anchor.
	For a downward-closed family $\mathcal K$ containing $\varnothing$,
	write $g(A)=1$ if and only if $[m]\setminus A\in\mathcal K$.
	A linear span representation over $k=\mathbb F_p$ consists of a
	finite-dimensional space $V$, a target $v\in V$, and subspaces
	$L_1,\ldots,L_m\leq V$ such that
	\[
	g(A)=1\quad\Longleftrightarrow\quad v\in\sum_{i\in A}L_i.
	\]
	Allow $v=0$ (in particular for the constant-true function), and
	assign the representation block cost $\sum_i\max(1,\dim L_i)$.
	
	\begin{theorem}[Exact chart cost in the fixed-kernel cycle model]
		\label{thm:fixed-kernel-span-cost}
		\leavevmode
		\begin{enumerate}[label=\textup{(\roman*)},leftmargin=*]
			\item \emph{Minimum number of charts.} In the restricted model just specified, the minimum number of charts
			realizing $\mathcal K$ equals
			\[
			2+\min_{(V,v,(L_i))}\sum_{i=1}^m\max(1,\dim L_i),
			\]
			where the minimum is over the linear span representations of $g$.
			
			\item \emph{Coverage and recovery.} Every candidate source in the realizing atlas covers the charts and
			complete boundaries and has free whole-chart recovery as in
			\cref{thm:complex}.
		\end{enumerate}
	\end{theorem}
	\begin{proof}
		\leavevmode
		\begin{enumerate}[label=\textup{(\roman*)},leftmargin=*]
			\item \emph{Minimum number of charts.} The graph description in the proof of \cref{thm:complex} used
			only the stated conditions on $U$ and $\kappa$; hence it applies to
			every atlas in this model. Denote the coordinate functional $e_b^*$
			restricted to $U$ by $\rho_b$. A selected family $J$ is feasible
			if and only if
			\[
			\kappa\in k\rho_*+\sum_{i\notin J}\langle\rho_b:b\in B_i\rangle,
			\]
			since a functional on $T$ vanishing on the selected blocks is a
			linear combination of the other coordinate functionals. In
			$V=U^*/k\rho_*$ take $v$ to be the image of $\kappa$ and $L_i$ the
			span of the images of the coordinate functionals of $B_i$.
			This gives a representation of $g$. Since $B_i$ is nonempty and
			$\dim L_i\leq |B_i|$, its block cost is at most the number of
			nonroot, nonanchor charts. This gives the lower bound.
			
			Conversely choose any span representation and replace $V$ by
			$\sum_i L_i$; it contains $v$ because $g([m])=1$. For each $i$ use
			$\max(1,\dim L_i)$ vertices, assigning to them basis vectors $v_b$
			of $L_i$ (one zero vector if $L_i=0$). Let $T'$ be the coordinate
			space of these vertices. The map
			\[
			j:V^*\longrightarrow T',\qquad
			j(\eta)=(\eta(v_b))_b
			\]
			is injective, since the $v_b$ span $V$. Set $T=k e_*\oplus T'$ and
			\[
			U=k e_*\oplus j(V^*),\qquad
			\kappa(ae_*+j(\eta))=a+\eta(v).
			\]
			These data satisfy $\kappa(e_*)=1$. Under the natural identification
			$U^*\cong k\oplus V$, the anchor functional is $(1,0)$, each block
			coordinate functional is $(0,v_b)$, and $\kappa=(1,v)$. The
			feasibility test above is therefore exactly the defining test for
			$g$.
			
			This realizes the chosen cost. Existence of a finite representation follows by applying the lower-bound construction to \cref{thm:complex}. The possible costs are nonnegative integers, so a least cost exists.
			
			\item \emph{Coverage and recovery.} The graph proof also gives the asserted coverage and recovery,
			independently of the choice of $U$.\qedhere
		\end{enumerate}
	\end{proof}
	
	\begin{corollary}[Exponential chart requirement in the fixed-kernel cycle model]
		\label{cor:fixed-kernel-chart-lower-bound}
		For each fixed odd prime $p$ and each $m\geq2$, some downward-closed
		family containing $\varnothing$ requires at least
		\[
		\sqrt{\frac{\binom{m}{\lfloor m/2\rfloor}}{1+2\log_2 p}}
		\ \geq\ \frac{2^{m/2}}{\sqrt{(m+1)(1+2\log_2 p)}}
		\]
		charts in the model of \cref{thm:fixed-kernel-span-cost}.
	\end{corollary}
	\begin{proof}
		Choose independently which subsets of size $\lfloor m/2\rfloor$ to
		include, include all smaller subsets and no larger subsets. These
		choices give $2^{\binom{m}{\lfloor m/2\rfloor}}$ different families
		containing $\varnothing$. Let $N$ be the largest minimum chart count
		among these families, which is finite by universality. Each family
		then admits data with at most $N$ charts. The $m$ positive block
		sizes have at most $(N+1)^m$ choices. For fixed sizes, $\dim T<N$;
		there are at most $p^{N^2}$ subspaces of $T$ (enumerate spanning
		lists padded with zeros), and at most $p^N$ choices of $\kappa$ for
		each. The group law and the root, anchor, ordered blocks and cuts
		are fixed by these data, and each choice determines one family.
		Consequently
		\[
		\binom{m}{\lfloor m/2\rfloor}\leq
		m\log_2(N+1)+(N^2+N)\log_2 p
		\leq (1+2\log_2 p)N^2,
		\]
		where $m\leq N$ and $\log_2(N+1)\leq N$. The central binomial
		coefficient is at least the average $2^m/(m+1)$ of the $m+1$
		binomial coefficients. This proves the bound.
	\end{proof}
	
	In this cycle model the extra vertices account for the root, anchor,
	and empty subspaces that still require a block. Accordingly, the block
	cost exceeds the row count of a smallest span representation by at
	most $m$.
	
	\section{Native extension lifting and preservation}
	
	Fix an intersection with the kernel of a group extension and an image
	in the quotient. Which subgroups realize both, and which of them
	preserve the chart observations? \Cref{thm:chief-normalizer-complement}
	identifies the lifts with complements in a normalizer quotient.
	For products of nonabelian simple groups, the diagonal-strip
	decomposition gives a uniqueness criterion. For abelian normalizer
	kernels, \cref{thm:cohomological-contract-enumerator,cor:abelian-complement-contracts}
	count the lifts satisfying the chart, recovery, and cut conditions
	simultaneously.
	
	\subsection{Nonsplit chief-layer lifts and contract-filtered outer monodromy}
	\label{subsubsec:nonsplit-chief-monodromy}
	
	A lift with prescribed kernel intersection exists when its conjugacy class is invariant and the resulting normalizer extension splits.
	
	\begin{theorem}[Normalizer--complement criterion for prescribed group-extension lifts]
		\label{thm:chief-normalizer-complement}
		Let
		\[
		1\longrightarrow N\longrightarrow G
		\overset{\pi}{\longrightarrow}Q\longrightarrow1
		\]
		be an extension of finite groups, not necessarily split.
		Fix $T\le Q$ and $K\le N$, and put
		\[
		G_T=\pi^{-1}(T),\qquad M_T=N_{G_T}(K),\qquad
		A_K=N_N(K)/K.
		\]
		The subgroup $T$ of the quotient acts on the set of $N$-conjugacy classes of
		subgroups of $N$ by conjugation using any representatives in $G$.
		Then:
		\begin{enumerate}[label=\textup{(\roman*)}]
			\item \emph{Invariance of the kernel class.}
			$\pi(M_T)=T$ if and only if the $N$-conjugacy class
			$[K]_N$ is fixed by $T$.
			\item \emph{Lifts and complements.}
			If this holds, there is a canonical extension
			\[
			1\longrightarrow A_K\longrightarrow M_T/K
			\longrightarrow T\longrightarrow1.
			\]
			Subgroups $H\le G$ satisfying
			$\pi(H)=T$ and $H\cap N=K$ correspond bijectively to
			complements of $A_K$ in $M_T/K$.
			\item \emph{Uniqueness and index.}
			If $N_N(K)=K$, such an $H$ exists precisely when
			$[K]_N$ is $T$-fixed; it is then unique and equals $M_T$.
			For every lift in \textup{(ii)}, the index is
			\[
			[G:H]=[Q:T]\,[N:K].
			\]
		\end{enumerate}
	\end{theorem}
	
	\begin{proof}
		\leavevmode
		\begin{enumerate}[label=\textup{(\roman*)},leftmargin=*]
			\item \emph{Invariance of the kernel class.}
			If $q\in T$ is represented by $g\in G_T$, then
			$q[K]_N=[K]_N$ precisely when $gKg^{-1}$ is conjugate
			to $K$ by an element of $N$.  Equivalently, some element of
			$Ng$ normalizes $K$.
			\item \emph{Lifts and complements.}
			Every subgroup $H$ in (ii) normalizes $K$, because
			$K=H\cap N\trianglelefteq H$; hence $H\le M_T$ and
			$H/K$ is a subgroup of $M_T/K$.  The kernel of
			$M_T/K\to T$ is $N_N(K)/K=A_K$.  The equalities
			$\pi(H)=T$ and $H\cap N=K$ say exactly that
			$H/K$ maps isomorphically onto $T$, which is the
			complement condition.  Conversely every complement has
			full image $T$ and preimage in $M_T$ intersecting $N$
			exactly in $K$.
			\item \emph{Uniqueness and index.}
			If $N_N(K)=K$ and $[K]_N$ is $T$-fixed, then $A_K=1$, so $M_T/K\cong T$ and
			the only complement is the whole group.  Hence $H=M_T$.
			Finally $|H|=|K||T|$, whereas $|G|=|N||Q|$. This calculation applies to every lift in (ii).\qedhere
		\end{enumerate}
	\end{proof}
	
	\begin{corollary}[Outer-monodromy criterion for subdirect chief layers]
		\label{cor:chief-subdirect-monodromy}
		Suppose $N=S^D$ for a finite nonabelian simple group $S$,
		with $D$ finite.  Let $K\le N$ project onto every simple
		coordinate.
		\begin{enumerate}[label=\textup{(\roman*)},leftmargin=*]
			\item \emph{Diagonal strips.}
			$K$ is a product of full twisted diagonal
			strips on a partition $\Pi$ of $D$, and
			\[
			N_N(K)=K,\qquad [N:K]=|S|^{|D|-|\Pi|}.
			\]
			\item \emph{Lifting and index.}
			For every $T\le Q$, a subgroup $H\le G$ satisfying
			$\pi(H)=T$ and $H\cap N=K$ exists if and only if
			the $N$-conjugacy class of the decorated strip partition
			is $T$-invariant; when it exists,
			\[
			H=N_{\pi^{-1}(T)}(K),\qquad
			[G:H]=[Q:T]|S|^{|D|-|\Pi|}.
			\]
			In particular, for $T=Q$ there is a bijection between
			$N$-conjugacy classes of such source subgroups and
			$Q$-fixed $N$-conjugacy classes of subdirect strip products.
			No splitting of $G\to Q$ is assumed.
		\end{enumerate}
	\end{corollary}
	
	\begin{proof}
		\leavevmode
		\begin{enumerate}[label=\textup{(\roman*)},leftmargin=*]
			\item \emph{Diagonal strips.}
			Scott's diagonal-strip lemma \cite{Scott1980}, in the form of
			\cite[Lemma~2.7]{JaikinPyber2011},
			gives the partition representation.  A full twisted
			diagonal in a direct power of a centerless simple group
			is self-normalizing within that power: if a coordinate
			conjugation preserves every equality
			$x_d=\alpha_d(x_{d_0})$ in the strip, then the conjugating
			coordinates differ by exactly the same twisted diagonal
			relation, because $Z(S)=1$.  A product of strips on disjoint
			supports is also self-normalizing in $N$, since inner
			automorphisms of $N$ do not permute the supports.
			Thus $N_N(K)=K$, and the strip representation gives $|K|=|S|^{|\Pi|}$.
			\item \emph{Lifting and index.}
			Apply \cref{thm:chief-normalizer-complement} using the normalizer and index from (i). For $T=Q$, the unique lift is equivariant under $N$-conjugation, and intersection with $N$ recovers $K$.\qedhere
		\end{enumerate}
	\end{proof}
	
	The $Q$-action on $[K]_N$ depends only on the outer action. Changing a
	lift by an element of $N$ changes the conjugate of $K$ only by
	$N$-conjugacy. The partition and relative outer-automorphism labels
	of the twisted identifications within each block encode $[K]_N$ for the fixed-point test.
	
	\begin{theorem}[Preservation tests and least index for chief-subdirect sources]
		\label{thm:chief-monodromy-full-contract}
		Assume the extension in
		\cref{cor:chief-subdirect-monodromy} is the full native
		matching group of a finite group atlas, and native
		realization permits unrestricted finite subgroup generation.
		Let the chart readouts $c_i:G\to X_i$ and descriptor
		readouts $\rho_d:G\to Y_d$ be homomorphisms, with
		finite descriptor set $\mathcal D$ and prescribed closure
		operator $c_0$ on $\mathcal D$. Take $X_i=c_i(G)$; if full
		coverage of a prescribed chart is required, first test that its
		state group equals this image.
		For $A\subseteq\mathcal D$ write
		$\rho_A=(\rho_d)_{d\in A}$.
		
		For each protected complete cut $t$ let
		$\ell_t:G\to L_t$, $r_t:G\to R_t$, and
		$b_t:G\to B_t$ be homomorphisms, with $b_t$ factoring
		through both side maps, such that the joint side map
		identifies $G$ with the full fiber product of
		$\ell_t(G)$ and $r_t(G)$ over $b_t(G)$.
		Write $D_t^L=\ker\ell_t$ and $D_t^R=\ker r_t$.
		For each proper meet-irreducible $c_0$-closed set
		$M$, let $M^+$ be its unique cover among the closed sets and choose
		$d_M\in M^+\setminus M$.
		\begin{enumerate}[label=\textup{(\roman*)},leftmargin=*]
			\item \emph{Preservation tests.}
			Let $[K]_N$ range over $Q$-fixed $N$-conjugacy
			classes of subdirect strip products in $N$, put
			$H_K=N_G(K)$, and require the following tests:
			\begin{enumerate}[label=\textup{(\alph*)}]
				\item \emph{Chart coverage:}
				$H_K\ker c_i=G$ for every $i$.
				\item \emph{Positive dependencies:}
				$H_K\cap\ker\rho_A\le\ker\rho_d$
				for all $A\subseteq\mathcal D$ and $d\in c_0(A)$.
				\item \emph{Exact recovery:}
				$H_K\cap\ker\rho_M\not\le\ker\rho_{d_M}$
				for every proper meet-irreducible $M$.
				\item \emph{Complete cuts:}
				\[
				H_K=(H_KD_t^L)\cap(H_KD_t^R),
				\qquad H_K\ker b_t=G
				\]
				for every named cut $t$.
			\end{enumerate}
			These tests are independent of the representative $K$ of
			its $N$-conjugacy class.  They are simultaneously necessary
			and sufficient for $H_K$ to preserve all the named chart,
			positive-dependency, exact-recovery, and full-range
			complete-cut contracts.
			\item \emph{Least index.}
			Among \emph{proper} sources with image $Q$ and
			subdirect intersection with $N$, the exact least index is
			\begin{equation}\label{eq:chief-subdirect-absolute-index}
				\iota_{\mathrm{chief\mbox{-}sd}}
				=\begin{cases}
					|S|^{\displaystyle\min_{[K]_N}(|D|-|\Pi_K|)},
					&\text{if a proper $Q$-fixed class passes (a)--(d)},\\
					+\infty,&\text{otherwise},
				\end{cases}
			\end{equation}
			where the minimum is taken only over passing classes
			with $K\ne N$.
		\end{enumerate}
	\end{theorem}
	
	\begin{proof}
		\leavevmode
		\begin{enumerate}[label=\textup{(\roman*)},leftmargin=*]
			\item \emph{Preservation tests.}
			By \cref{cor:chief-subdirect-monodromy}, every source in
			the stated class is $H_K=N_G(K)$ for a $Q$-fixed class
			$[K]_N$, and conversely every such class lifts uniquely.
			Conjugating $K$ by an element $n\in N$ conjugates $H_K$
			by $n$.  Since the kernels of all declared homomorphisms
			are normal in $G$, the four tests are unchanged under
			this conjugation.
			
			For any homomorphism $f$ on $G$,
			$f(H_K)=f(G)$ if and only if $H_K\ker f=G$.
			This proves (a) and the boundary-range part of (d).
			On $H_K$, the functional dependency $A\to d$ holds
			precisely when
			$H_K\cap\ker\rho_A\le\ker\rho_d$.
			Therefore (b) is exactly $c_0\le c_{H_K}$.
			For a proper meet-irreducible closed $M$, the set
			$c_{H_K}(M)$ is $c_0$-closed once (b) holds.
			If it properly contains $M$, it contains the unique
			cover $M^+$, and therefore $d_M$.
			Consequently (c) is equivalent to
			$c_{H_K}(M)=M$.  Every closed set of a finite closure
			system is an intersection of meet-irreducible closed
			sets containing it.  Hence (b) and (c) together are
			equivalent to $c_{H_K}=c_0$.
			
			For a complete cut, the preimage in $G$ of
			$\ell_t(H_K)$ is $H_KD_t^L$, and similarly the
			preimage of $r_t(H_K)$ is $H_KD_t^R$.
			Since the full ambient group is the fiber product over
			the complete boundary, the full fiber product of the
			\emph{actual} side images of $H_K$ is represented in $G$
			by $(H_KD_t^L)\cap(H_KD_t^R)$.
			Thus the first equality in (d) is precisely represented
			common-information exactness; the second is full
			boundary range.
			\item \emph{Least index.}
			The index follows from
			$[G:H_K]=[N:K]=|S|^{|D|-|\Pi_K|}$.
			All groups involved are finite, so the unrestricted
			native law admits a finite generating set for $H_K$.
			Minimizing over the passing proper classes gives \eqref{eq:chief-subdirect-absolute-index}.\qedhere
		\end{enumerate}
	\end{proof}
	
	\begin{proposition}[Sharp nonsplit lifting and nonlifting examples]
		\label{prop:chief-nonsplit-sharpness}
		There exist nonsplit extensions with a nonabelian minimal normal
		subgroup $N$ for which the following possibilities occur:
		\begin{enumerate}[label=\textup{(\roman*)}]
			\item \emph{Unique lift.}
			A proper subdirect diagonal intersection has a
			unique subgroup lift onto the quotient.
			\item \emph{Obstruction from outer labels.}
			An otherwise identical diagonal type has no such
			lift because of its outer-automorphism labels.
			\item \emph{Obstruction from nonsplitting.}
			A quotient-stable intersection fails to lift when
			the self-normalizing hypothesis is removed.
		\end{enumerate}
	\end{proposition}
	
	\begin{proof}
		We use the nonsplit extension
		$1\to A_6\to E\to C_2\to1$ with $E=M_{10}$;
		no element of the outer coset of $A_6$ in $E$ is an
		involution \cite[p.~354]{GuralnickMalle2002}; see also the entry for
		$A_6$ in the \emph{Atlas} \cite{ConwayEtAl1985Atlas}.
		\begin{enumerate}[label=\textup{(\roman*)},leftmargin=*]
			\item \emph{Unique lift.}
			Form
			\[
			G_+=E\times_{C_2}E
			=\{(x,y)\in E^2:xA_6=yA_6\},\quad
			N=A_6^2,\quad Q=C_2.
			\]
			For $K=\operatorname{diag}(A_6)$, the quotient acts
			by the same outer automorphism on both simple factors,
			so $[K]_N$ is fixed.  The unique lift is
			$H_+=\operatorname{diag}(E)$ and
			$[G_+:H_+]=|A_6|=360$.
			The extension $G_+\to C_2$ is nonsplit, since an
			involution outside $N$ would project to an outer
			involution in $E$.
			
			To make $N$ minimal normal, adjoin the involution $\tau$ interchanging
			the two coordinates. Every normal subgroup of $A_6^2$ is a product of
			factors, since a nontrivial projection onto a factor yields that
			factor by taking commutators and using simplicity and trivial center.
			The interchange acts transitively on the two factors. Hence $N$ is
			minimal normal in the enlarged group. A section of the enlarged
			extension would restrict to a section over the old quotient, so the
			extension remains nonsplit. The diagonal remains fixed, with unique
			lift $\operatorname{diag}(E)\times\langle\tau\rangle$.
			\item \emph{Obstruction from outer labels.}
			Let $G_-=E\times E$, with the same normal
			layer $N=A_6^2$ and quotient $Q=C_2^2$.
			Every proper subdirect subgroup of $N$ is a twisted
			full diagonal, the graph of some automorphism of $A_6$.
			An element acting outerly on the first factor and
			trivially on the second changes its relative outer
			automorphism label.  This change cannot be achieved
			by conjugation inside $N$.
			Thus no proper subdirect $K$ has $Q$-fixed conjugacy
			class, and no $H\le G_-$ with $\pi(H)=Q$ has
			proper subdirect intersection with $N$.
			This extension is also nonsplit: a complement would
			contain an involution mapping to $(1,0)\in C_2^2$,
			whose first projection would be an outer involution in $E$.
			
			Adjoin the coordinate interchange as in (i). The resulting extension
			has minimal normal subgroup $N$ and remains nonsplit. Any fixed diagonal
			class for the enlarged quotient would already be fixed by the old quotient.
			\item \emph{Obstruction from nonsplitting.}
			Take the same nonsplit extension $E\to C_2$
			but choose $K=1\le A_6$.  The class $[1]_{A_6}$ is
			fixed by $C_2$, yet a source meeting $A_6$ trivially
			and surjecting onto $C_2$ would split the extension.
			No such source exists.  Here $N_N(K)/K=A_6$ is
			nontrivial, and the normalizer--complement obstruction
			in \cref{thm:chief-normalizer-complement} is essential.\qedhere
		\end{enumerate}
	\end{proof}
	
	\begin{example}[A faithful nonsplit atlas with a protected complete cut]
		\label{ex:chief-nonsplit-faithful-cut}
		With $E=M_{10}$ and $A=C_3$, consider three charts
		\[
		X_1=E\times A,\qquad X_2=E\times A,\qquad
		X_3=E\times A
		\]
		on a path $1-2-3$.  On $12$ the interface is the
		whole $E$-coordinate; on $23$ it is $C_2\times A$,
		read by the outer coset in $E$ and the $A$-coordinate.
		With global coordinates $(x,y,u,v)$ the full matching
		group is
		\[
		G=\{(x,y,u,v)\in E^2\times A^2:
		xA_6=yA_6\},
		\]
		and the chart maps are
		$(x,y,u,v)\mapsto(x,u),(x,v),(y,v)$.
		The extension of $Q=C_2\times A^2$ by
		$N=A_6^2$ is nonsplit.  The proper native source
		\[
		H=\{(x,x,u,v):x\in E,\ u,v\in A\}
		\]
		is the unique subgroup over $Q$ with intersection
		$\operatorname{diag}(A_6)$ with $N$.
		It covers every individual chart and has index $|A_6|$.
		
		Protect the complete cut $1\mid\{2,3\}$, whose
		entire boundary is the $E$-state on edge $12$.
		The boundary image of $H$ equals $E$, and the fiber
		product of its two actual side images is exactly $H$.
		For the two declared whole-state descriptors
		$d_1=(x,u)$ on chart $1$ and $d_2=(y,v)$ on chart $3$,
		both the full source and $H$ have free recovery closure:
		when $d_1$ is fixed, $v$ can vary; when $d_2$ is fixed,
		$u$ can vary.  Thus $H$ satisfies chart coverage,
		exact recovery and the named full-range complete-cut
		contract in a nonsplit native atlas.
		
		The other complete cut $\{1,2\}\mid3$ is \emph{not}
		exact for $H$: it records only the common outer coset
		and the $A$-state, and does not transmit the equality
		$x=y$.  This illustrates why the normalizer-lift test
		must be followed by the separate cut contracts in
		\cref{thm:chief-monodromy-full-contract}.
	\end{example}
	
	\begin{remark}[Choice of the layer]
		\label{rem:chief-nonsplit-provenance}
		The strip criterion applies to any normal layer $N=S^D$ of the stated
		type. For $D$ nonempty, $N$ is minimal normal precisely when $G$ acts
		transitively on its simple factors. Indeed, the products on the factor
		orbits are its minimal nontrivial $G$-invariant normal subgroups.
	\end{remark}

	\subsection{Cohomological obstruction and simultaneous preservation for abelian normalizer kernels}
	\label{subsubsec:abelian-normalizer-contracts}

	When the normalizer kernel is abelian, several complements may lift the
	same subgroup of the quotient. The extension class decides whether a
	lift exists. To count the lifts that preserve several observations at
	once, we record their intersections with the observation kernels in one
	polynomial. Elementary abelian kernels give a formula in matrix ranks.

	\begin{lemma}[Kernel sizes for a lifted complement]
		\label{lem:normalizer-complement-kernel-count}
		Let $1\to N\to G\overset{\pi}{\to}Q\to1$ be an extension of
		finite groups, let $T\le Q$ stabilize the $N$-conjugacy class
		of $K\le N$, and put
		\[
		M=N_{\pi^{-1}(T)}(K),\qquad
		A=N_N(K)/K,\qquad E=M/K.
		\]
		Then $1\to A\to E\to T\to1$ is exact.  For every complement
		$C\le E$, let $H_C\le M$ denote its full preimage under $M\to E$.
		For a normal subgroup $D\trianglelefteq G$ define
		\[
		\kappa_D=|K\cap D|,\quad
		\overline D=(KD\cap M)/K\trianglelefteq E,\quad
		m_D(C)=|C\cap\overline D|.
		\]
		The image of $H_C\cap D$ under $M\to E$ is $C\cap\overline D$,
		with kernel $K\cap D$. In particular,
		\begin{equation}\label{eq:complement-kernel-cardinality}
			|H_C|=|K||T|,\qquad
			|H_C\cap D|=\kappa_D\,m_D(C).
		\end{equation}
	\end{lemma}

	\begin{proof}
		The existence of the exact normalizer sequence is
		\cref{thm:chief-normalizer-complement}.  The image of
		$D\cap M$ in $E$ is $\overline D$, since $K\le M$.
		The restriction $H_C\cap D\to C\cap\overline D$ is onto:
		if $h\in H_C$ and $d\in D\cap M$ have the same coset
		modulo $K$, then $d=hk$ for some $k\in K\le H_C$.
		Its kernel is $K\cap D$.  This proves the second formula;
		the first follows because $H_C/K\cong C\cong T$.
	\end{proof}

	\begin{theorem}[Cohomological contract enumerator for simultaneous native lifts]
		\label{thm:cohomological-contract-enumerator}
		Retain the setup of \cref{lem:normalizer-complement-kernel-count} and assume
		$A$ is an elementary abelian $p$-group.
		Choose a normalized set section $s:T\to E$ and let $\alpha\in Z^2(T,A)$ be its factor set for the induced $T$-action on $A$.		
		\begin{enumerate}[label=\textup{(\roman*)}]
			\item \emph{Existence and conjugacy.}
			A subgroup $H\le G$ with $\pi(H)=T$ and $H\cap N=K$	exists if and only if $[\alpha]=0$ in $H^2(T,A)$.
			When the class vanishes, the lifts are in bijection with the following affine space of complement cocycles
			\[
			\mathscr Z_\alpha
			=\{z:T\to A:\ z(1)=0,\quad
			z(tu)=z(t)+t\cdot z(u)+\alpha(t,u)\}.
			\]
			Conjugacy by $A$ identifies cocycles modulo $B^1(T,A)$.
			Consequently the set of $A$-conjugacy classes of lifts is a torsor for $H^1(T,A)$.
			Here $A$ acts on the complements in $E$.
			On the corresponding subgroups of $M$, this is conjugation by
			representatives in $N_N(K)$.
			
			\item \emph{Counting complement classes.}
			Let $\mathscr D$ be any finite set of distinct normal subgroups of $G$.  It may include kernels of chart, descriptor, and protected-cut homomorphisms.
			Write $h=\dim_{\mathbb F_p}H^1(T,A)$.
			When the extension splits, for each $D\in\mathscr D$ and $t\in T$ the condition
			\begin{equation}\label{eq:normalizer-affine-atom}
				z(t)s(t)\in\overline D
			\end{equation}
			becomes a (possibly inconsistent) affine linear system
			$L_{D,t}v=b_{D,t}$ in $h$ coordinates $v\in\mathbb F_p^h$
			for complement classes.
			If $[\alpha]=0$, define the intrinsic polynomial
			\begin{equation}\label{eq:cohomological-contract-polynomial}
				\mathscr F_{K,T}(\mathbf X)
				=\sum_{[C]}\prod_{D\in\mathscr D}X_D^{m_D(C)},
			\end{equation}
			where the sum runs over $A$-conjugacy classes of complements.
			It has the following rank-only expression, which does not enumerate those complements:
			\begin{equation}\label{eq:cohomological-rank-expansion}
				\boxed{\displaystyle
					\mathscr F_{K,T}(\mathbf X)
					=\sum_{S\subseteq\mathscr D\times T}
					\nu(S)\prod_{(D,t)\in S}(X_D-1),\qquad
					\nu(S)=
					\begin{cases}
						p^{h-\operatorname{rank}L_S},
						&L_Sv=b_S\text{ is consistent},\\
						0,&\text{otherwise}.
				\end{cases}}
			\end{equation}
			Here $L_Sv=b_S$ is the combined affine system of the
			atoms in $S$, with the empty system contributing $p^h$.
			If $[\alpha]\ne0$, set $\mathscr F_{K,T}=0$.
			The polynomial is independent of the set section, the
			cocycle coordinates, and the chosen affine gauge.
			
			\item \emph{Testing the preservation contracts.}
			Assume that $G$ is the full native matching group of
			a finite atlas, native realization allows unrestricted finite
			generation, and every declared observation is a group
			homomorphism.  Assume that the listed positive dependencies
			axiomatize the prescribed target closure $c_0$ (equivalently,
			one may test all $P\to d$ with $d\in c_0(P)$).
			Include in $\mathscr D$ the kernels needed for
			all named contracts.  For a monomial of
			\eqref{eq:cohomological-contract-polynomial} with exponent
			vector $(m_D)$, put
			\[
			n=|K||T|,\qquad a_D=\kappa_Dm_D.
			\]
			Then that monomial represents a source satisfying all the
			following contracts if and only if the displayed tests hold:
			\begin{align}
				\text{chart coverage:}\quad&
				a_{\ker c_i}=n/|c_i(G)| &&\text{for each }i;\label{eq:coho-chart}\\
				\text{positive dependencies:}\quad&
				a_{\ker\rho_P}=a_{\ker(\rho_P,\rho_d)}
				&&\text{for every named }P\to d;\label{eq:coho-positive}\\
				\text{exact recovery:}\quad&
				a_{\ker\rho_M}>a_{\ker(\rho_M,\rho_{d_M})}
				&&\text{for every proper meet-irreducible }M;\label{eq:coho-recovery}\\
				\text{exact complete cut:}\quad&
				a_{\ker b_t}\,a_{\ker(\ell_t,r_t)}
				=a_{\ker\ell_t}\,a_{\ker r_t}
				&&\text{for every protected }t;\label{eq:coho-cut}\\
				\text{full boundary range:}\quad&
				a_{\ker b_t}=n/|b_t(G)|
				&&\text{for every protected }t.\label{eq:coho-range}
			\end{align}
			The recovery test presupposes the positive tests and uses
			$d_M\in M^+\setminus M$ as in
			\cref{thm:chief-monodromy-full-contract}.
			Here coverage is measured against the declared ranges; if a chart
			must cover a fixed $X_i$, first check $c_i(G)=X_i$. At each cut assume
			that the boundary factors through each side. When the sides jointly
			separate $G$, the factor $a_{\ker(\ell_t,r_t)}$ equals one.
			The sum of coefficients of
			\eqref{eq:cohomological-contract-polynomial} over monomials
			passing \eqref{eq:coho-chart}--\eqref{eq:coho-range}
			counts the $A$-conjugacy classes of simultaneous faithful
			lifts.  In particular, a faithful lift exists exactly when
			this sum is positive.
		\end{enumerate}
	\end{theorem}

	\begin{proof}
		\leavevmode
		\begin{enumerate}[label=\textup{(\roman*)},leftmargin=*]
			\item \emph{Existence and conjugacy.}
			Choose the normalized section and write
			$s(t)s(u)=\alpha(t,u)s(tu)$ in additive notation for $A$.
			The conjugation action of $T$ on $A$ is independent of the
			chosen representatives because $A$ is abelian.
			A putative complement has a unique representative
			$z(t)s(t)$ above each $t\in T$; closure under multiplication
			is equivalent to the displayed affine cocycle equations.
			They admit a solution precisely when $\alpha$ is a
			coboundary, i.e., $[\alpha]=0$.  Their differences form
			$Z^1(T,A)$.  Conjugation by $a\in A$ changes $z$ by a
			$1$-coboundary, so $A$-conjugacy classes form a torsor
			for $H^1(T,A)$.
			
			\item \emph{Counting complement classes.}
			Fix a solution $z_0$ and linear representatives
			for $Z^1/B^1$.  Each class then has a representative
			$z_0+z_v$ depending affinely on $v\in\mathbb F_p^h$.
			For fixed $D,t$, the set of $a\in A$ with
			$a s(t)\in\overline D$ is either empty or a coset of
			the vector subgroup $A\cap\overline D$.
			Thus \eqref{eq:normalizer-affine-atom} imposes affine
			linear conditions on $v$.
			The condition is invariant under $A$-conjugation because
			$\overline D\trianglelefteq E$.
			For each $v$, let $I_{D,t}(v)$ be its zero--one
			indicator.  Then
			\[
			\prod_D X_D^{m_D(v)}
			=\prod_{(D,t)\in\mathscr D\times T}
			\bigl(1+(X_D-1)I_{D,t}(v)\bigr).
			\]
			Expanding and summing over $v$ gives
			\eqref{eq:cohomological-rank-expansion}, since a consistent
			affine system of rank $r$ in $h$ variables has exactly
			$p^{h-r}$ solutions.  The original definition
			\eqref{eq:cohomological-contract-polynomial} proves
			independence of all choices.
			
			\item \emph{Testing the preservation contracts.}
			Apply
			\cref{lem:normalizer-complement-kernel-count} to every
			normal observation kernel.  A homomorphism $f$ has
			$|f(H_C)|=n/a_{\ker f}$, proving the chart and boundary
			range tests.  The dependency $P\to d$ holds on $H_C$
			exactly when
			$H_C\cap\ker\rho_P=H_C\cap\ker(\rho_P,\rho_d)$,
			which for nested finite subgroups is equality of their
			orders.  Once all positive dependencies hold, exact
			recovery is equivalent to a nontrivial element in
			$(H_C\cap\ker\rho_M)\setminus\ker\rho_{d_M}$
			for every proper meet-irreducible $M$.  This proves
			\eqref{eq:coho-positive} and \eqref{eq:coho-recovery}.
			
			The joint side image of $H_C$ embeds in the fiber product of
			$\ell_t(H_C)$ and $r_t(H_C)$ over $b_t(H_C)$. Its size is $n/a_{\ker(\ell_t,r_t)}$.
			The two boundary maps are onto that actual range and are
			homomorphisms, so their fibers have constant size. The fiber
			product therefore has order
			\[
			\frac{|\ell_t(H_C)|\,|r_t(H_C)|}{|b_t(H_C)|}
			=\frac{n\,a_{\ker b_t}}{a_{\ker\ell_t}a_{\ker r_t}}.
			\]
			The embedding is onto precisely when the two sizes agree, which
			gives \eqref{eq:coho-cut}. Equivalently, the product of the two normal
			side kernels on $H_C$ is the boundary kernel: containment follows
			from factorization, and equality of orders is exactly  \eqref{eq:coho-cut}.
			Functions of both sides are precisely those constant on cosets of
			this product, giving represented common-information exactness.
			Since every source in the specified class comes from a complement, filtering the nonnegative monomial coefficients gives necessity and sufficiency.\qedhere
		\end{enumerate}
	\end{proof}

	\begin{corollary}[Contract counts for finite abelian normalizer kernels]
		\label{cor:abelian-complement-contracts}
		Let $A$ in \cref{lem:normalizer-complement-kernel-count} be any finite
		abelian group.
		\begin{enumerate}[label=\textup{(\roman*)},leftmargin=*]
			\item \emph{Obstruction and preservation tests.}
			The obstruction in $H^2(T,A)$, the parametrization by
			$H^1(T,A)$ when it vanishes, and the tests
			\eqref{eq:coho-chart}--\eqref{eq:coho-range} remain valid.
			
			\item \emph{Coefficients for a split extension.}
			For a split extension, the expansion in
			\eqref{eq:cohomological-rank-expansion} has the following coefficients.
			If any $(D,t)\in S$ has $t\notin\operatorname{im}(\overline D\to T)$,
			put $\nu(S)=0$. Otherwise define
			\[
			e_S:H^1(T,A)\longrightarrow\prod_{(D,t)\in S} A/(A\cap\overline D),
			\qquad [w]\longmapsto
			\big(w(t)\bmod (A\cap\overline D)\big)_{(D,t)\in S}.
			\]
			Fix a homomorphic section $s_*:T\to E$. For each atom choose
			$a_{D,t}\in A$ with $a_{D,t}s_*(t)\in\overline D$, and write $b_S$
			for the tuple of their residue classes. Then $\nu(S)=|\ker e_S|$ if
			$b_S\in\operatorname{im}e_S$, and $\nu(S)=0$ otherwise.
		\end{enumerate}
	\end{corollary}

	\begin{proof}
		\leavevmode
		\begin{enumerate}[label=\textup{(\roman*)},leftmargin=*]
			\item \emph{Obstruction and preservation tests.}
			The cocycle and conjugation arguments in
			\cref{thm:cohomological-contract-enumerator}(i) use only commutativity
			of $A$. The preservation tests use only the kernel orders in
			\cref{lem:normalizer-complement-kernel-count}.
			\item \emph{Coefficients for a split extension.}
			For a split extension, a complement is therefore represented by $t\mapsto w(t)s_*(t)$
			for a crossed homomorphism $w$. For each possible atom take
			$g=a_{D,t}s_*(t)\in\overline D$. For $a\in A$ the commutator
			$aga^{-1}g^{-1}$ lies in $A\cap\overline D$ by normality. In additive
			coordinates it equals $a-t\cdot a$. Thus coboundaries evaluate to zero
			in the indicated quotients, and $e_S$ is a well-defined homomorphism.
			An atom holds exactly when the corresponding component of $e_S([w])$
			equals its specified residue. Consequently the solution classes form
			the empty set or one coset of $\ker e_S$. Expanding the same product
			of indicators gives the coefficient formula.\qedhere
		\end{enumerate}
	\end{proof}
	
	\begin{example}[A split lift with incompatible complete-cut contracts]
		\label{ex:coho-incompatible-cuts}
		Let $S=A_5$ and use additive coordinates $x,y,b\in\mathbb F_2$.
		Put
		\[
		G=S\times\mathbb F_2^3,\quad
		N=S\times\mathbb F_2^{\{y,b\}},\quad
		K=S\times\mathbb F_2^{\{b\}},\quad
		T=\mathbb F_2^{\{x\}}.
		\]
		Then $A=N_N(K)/K\cong\mathbb F_2^{\{y\}}$ and the
		normalizer extension splits.  Its two complement lifts are
		\[
		H_0=\{(s,x,0,b)\},\qquad
		H_1=\{(s,x,x,b)\}.
		\]
		Both cover the chart $c(s,x,y,b)=(s,x,b)$ and preserve
		exact free recovery on the descriptor $d=x$.
		There are two complete-cut presentations with the same
		nontrivial full boundary $b$:
		\[
		\begin{array}{c|cc|c}
			&\ell(s,x,y,b)&r(s,x,y,b)&b(s,x,y,b)\\ \hline
			t_1&(s,x,b)&(y,b)&b\\
			t_2&(s,x,b)&(x+y,b)&b.
		\end{array}
		\]
		Each joint side map identifies $G$ with the full fiber
		product over $b$.  The first cut is exact on $H_0$ and
		not on $H_1$; the second is exact on $H_1$ and not on
		$H_0$.  Indeed, $|H_0|=|H_1|=4|S|$, and the two
		actual-side fiber-product cardinalities, in order
		$(t_1,t_2)$, are
		\[
		H_0:\ (4|S|,8|S|),\qquad
		H_1:\ (8|S|,4|S|).
		\]
		Both preserve the full boundary range $\mathbb F_2$,
		but neither is simultaneously faithful for both cuts.
		Thus $[\alpha]=0$ and every individual cut has a
		faithful complement, while the contract-filtered
		cohomological support is empty for the joint contract.
		The two cuts here are homomorphic presentations on the same group.
	\end{example}
	
	\begin{example}[A nonzero splitting obstruction with a nonabelian core]
		\label{ex:coho-nonsplit-core}
		Let $S=A_5$, take $G=S\times C_4$,
		$N=S\times C_2$ (where $C_2$ is the unique
		order-two subgroup of $C_4$), and $K=S\times1$.
		Then $T\cong C_2$, $A=N_N(K)/K\cong C_2$, and
		$N_G(K)/K\cong C_4$.  The resulting class in
		$H^2(C_2,C_2)$ is nonzero.  Thus no subgroup
		$H$ has $\pi(H)=T$ and $H\cap N=K$, even before
		any preservation contracts are imposed.
	\end{example}
	
	\begin{remark}[Counting complements subject to the contract]
		\label{rem:coho-contract-enumerator-scope}
		The $H^2$ obstruction and $H^1$ parametrization are classical
		group cohomology \cite{BrownCohomology1982}. For fixed $K,T$ and an elementary abelian normalizer quotient,
		the inclusion--exclusion count can have exponentially many terms in $|T||\mathscr D|$.
		Minimizing over $K,T$ also requires choosing the intersection and quotient image.
	\end{remark}
	
	\section{Contract-filtered composition and coupled extension spectra}
	
	Two families of group fragments are equivalent if every allowed
	continuation accepts both or rejects both. This section gives complete
	invariants for two formation laws. \Cref{thm:contract-filtered-spectrum}
	classifies independent choices by the sets of records they attain.
	\Cref{thm:coupled-extension-spectrum} shows that prescribing a subgroup
	of the joint kernel exposes the marked extension and its readout
	kernels. For central elementary abelian kernels,
	\cref{thm:affine-coupled-spectrum} computes these classes by linear
	algebra.
	
	\subsection{Complete invariants for coverage-preserving, contract-filtered composition}
	\label{subsec:contract-filtered-invariants}
	
	A component can have several admissible realizations, and a formation
	must cover the states of each chosen input. We use quotient observations
	on the actual boundary image to distinguish the resulting families
	while preserving that coverage.
	
	\begin{lemma}[Recovery of continuation records by coverage-preserving quotient tests]
		\label{lem:covered-normal-extension}
		Let $G$ be a marked fragment as in \cref{def:fragment}, with
		$Q=p(G)\leq P$. Write $\mathcal M(Q)$ for the proper meet-irreducible
		normal subgroups of $Q$, that is, the proper normal subgroups with a
		unique cover in the finite normal-subgroup lattice. For each
		$N\in\mathcal M(Q)$ name the quotient readout $q_N:Q\to Q/N$.
		Attach the chart $Q\hookrightarrow P$ by equality of the complete
		boundary, with these new readouts; write $J$ for their labels.
		Then the following statements hold.
		\begin{enumerate}[label=\textup{(\roman*)},leftmargin=*]
			\item \emph{Mutual coverage.}
			The matching group $G\times_P Q$ is canonically isomorphic to $G$;
			both chart projections are onto. More generally the same attachment to
			any fragment with image $Q'$ covers both charts exactly when $Q'=Q$.
			\item \emph{Recovery closure.}
			The recovery closure $c$ of the enlarged readout family $D\sqcup J$
			recovers the entire qualitative record by
			\begin{equation}\label{eq:covered-extension-decoding}
				B_G(T)=\bigcap_{\substack{N\in\mathcal M(Q)\\q_N\in c(T)}}N,
				\qquad V_G(T)=c(T\cup J)\cap D\qquad(T\subseteq D).
			\end{equation}
			The intersection of an empty family is $Q$.
			Conversely the qualitative record and $Q$ determine this closure.
			\item \emph{Numerical record.}
			The numerical record is recovered by additionally specifying
			$|G|$ and the joint range sizes $r(T\cup J)$ on the enlarged family:
			\begin{equation}\label{eq:covered-extension-count}
				m_G(T)=\frac{|G|}{r(T\cup J)}.
			\end{equation}
			\item \emph{Minimal library.}
			Among libraries of quotient observations of $Q$ used to distinguish
			all normal subgroups by their individual containment tests, the unique
			inclusion-minimal library is $\mathcal M(Q)$ (identifying quotient maps
			with the same kernel). No proper sublibrary of $\mathcal M(Q)$ suffices.
		\end{enumerate}
	\end{lemma}
	\begin{proof}
		\leavevmode
		\begin{enumerate}[label=\textup{(\roman*)},leftmargin=*]
			\item \emph{Mutual coverage.}
			The matching group consists of $(g,p(g))$ for $p(g)\in Q$.
			For a general $Q'$ the reachable boundary is $Q\cap Q'$, so mutual
			coverage is equivalent to $Q=Q'$.
			\item \emph{Recovery closure.}
			The readout $q_Np$ is determined by $T$ if and only if
			$K_G(T)\leq\ker(q_Np)$, equivalently $B_G(T)\leq N$.
			The image of the normal subgroup $K_G(T)$ under $G\twoheadrightarrow Q$
			is normal in $Q$. Every element of a finite lattice is the meet of the meet-irreducibles
			above it. Indeed, the top is the empty meet, and every proper
			meet-reducible element is a meet of two strictly greater elements.
			Downward induction gives the claim. In particular the joint
			kernel of the $q_N$ on $Q$ is trivial (also when $Q=1$ and $J$ is empty).
			Thus agreement on $J$ is precisely agreement on the complete boundary.
			These facts give (ii); its converse follows from the qualitative
			composition rule of \cref{thm:complete-compositional-invariants}.
			\item \emph{Numerical record.}
			The joint agreement kernel on $T\cup J$ is $K_G(T)\cap\ker p$;
			the homomorphism theorem gives (iii).
			\item \emph{Minimal library.}
			If a meet-irreducible $N$ with unique cover $N^+$ is omitted,
			then every normal subgroup strictly containing $N$ contains $N^+$.
			Consequently every individual containment test by a kernel other than
			$N$ has the same answer on $N$ and $N^+$. The library of all
			meet-irreducibles suffices by (ii). Thus every separating library of individual quotient tests contains
			all the meet-irreducibles.\qedhere
		\end{enumerate}
	\end{proof}

	\begin{definition}[Contract-filtered fragment families and their formation language]
		\label{def:contract-filtered-fragment}
		Fix the finite group-table realization context and the marked boundary
		and readout types of \cref{def:fragment}.
		A \textbf{contract-filtered fragment family} is a finite, possibly empty
		family of admissible realizations of the specified marked type.
		Its admissible inputs, formations, and tests are specified as follows.
		\begin{enumerate}[label=\textup{(\roman*)},leftmargin=*]
			\item \emph{Admissible input families.}
			An input component supplies such a family $\mathscr H$.
			In particular, a primitive component may be specified by
			an ambient marked group $G$, a named subgroup $K\leq\ker p$, and a
			boundary-image request $Q\leq P$, and supply the family
			\begin{equation}\label{eq:contract-fixed-k-family}
				\{H\leq G:p(H)=Q,\ H\cap\ker p=K\},
			\end{equation}
			with the maps restricted to $H$ and any stated local filters applied.
			A finite library of such requests is also allowed. A filter on a
			primitive realization may, for example, require coverage of prescribed
			local chart groups; such a filter is evaluated before the family is
			supplied. An ordinary fragment is the singleton family.
			
			\item \emph{Formation and recovery tests.}
			A formation chooses one realization from each occurrence, constructs
			its full matching group using the subgroup constraint $C$, and retains
			or introduces readouts and a boundary $\rho$ as in
			\cref{def:fragment}. It is required to cover every chosen input
			realization. Additional gates are exact recovery of a specified
			closure on a retained vocabulary, or specified positive and negative
			recovery judgments on that vocabulary. A gate may also test equality
			to the least closure generated by the local recovery laws on the
			specified input vocabularies (and any declared rules on new labels).
			When a closure is prescribed independently, its target is part of the
			formation data, not computed from the chosen realization. A gate rejects
			a choice rather than changing its matching group. All successful
			choices form the output family; failure of every choice gives the empty
			family. Each occurrence is an independent choice, not a shared unknown.
			The formation expressions and their gates are well typed independently
			of their success, and are closed under substitution.
			
			\item \emph{Numerical tests.}
			In the numerical language one can additionally test exact source order
			and specified joint or conditional range sizes. A chart map retained
			as a readout can be required to cover its declared target. At a named
			cut, two complete side readouts and a boundary factoring through both
			may be retained and subjected to exact represented common-information
			tests. These extra tests are on named maps of the chosen realization.
			
			\item \emph{Terminal tests.}
			The terminal test is existence of a successful choice (including the
			identity continuation). The language includes the two-chart attachments
			of \cref{lem:covered-normal-extension} with a fixed target closure on the
			enlarged vocabulary; in the numerical language it also allows a
			simultaneous list of prescribed range sizes and a source order. The targets of this test are fixed in advance, as allowed by
			\cref{subsec:framework-contract}.
		\end{enumerate}
	\end{definition}

	For $\epsilon\in\{\mathrm{rec},\mathrm{num}\}$ define the
	\textbf{attainable record spectrum}
	\begin{equation}\label{eq:contract-spectrum}
		\mathcal S_\epsilon(\mathscr H)
		=\{\mathcal I_\epsilon(H):H\in\mathscr H\}.
	\end{equation}
	Each member records values attained together on one realization.
	Repeated realizations of a record contribute a single member.
	Write $\equiv_{\Phi,\epsilon}$ for equality of Boolean success in
	every common continuation of the respective language.

	\begin{theorem}[Complete contract-filtered compositional spectra]
		\label{thm:contract-filtered-spectrum}
		Fix one of the languages of \cref{def:contract-filtered-fragment}. These languages permit independently prescribed full recovery targets (and, in the numerical variant, the specified numerical gates).
		\begin{enumerate}[label=\textup{(\roman*)},leftmargin=*]
			\item \emph{Invariance and composition.}
			The spectra in \eqref{eq:contract-spectrum} are invariant under marked isomorphisms transporting the admissible families. For each tuple of input records compute the output record by
			\cref{thm:complete-compositional-invariants}, rejecting a tuple
			unless, with $Q_i=B_{H_i}(\varnothing)$ and
			$M=C\cap\prod_i Q_i$,
			\begin{equation}\label{eq:contract-spectrum-coverage}
				\operatorname{pr}_i M=Q_i\qquad\text{for every }i.
			\end{equation}
			Reject also on failure of any declared gate. Collect precisely the
			remaining output records, with repetitions removed. This rule works
			for arbitrary finite matching graphs, including cyclic graphs, and
			arbitrary finite nested expressions in the language.
			
			\item \emph{Completeness and separation.}
			The spectra	in \eqref{eq:contract-spectrum} satisfy
			\begin{equation}\label{eq:contract-spectrum-completeness}
				\boxed{\mathscr H\equiv_{\Phi,\epsilon}\mathscr H'
					\quad\Longleftrightarrow\quad
					\mathcal S_\epsilon(\mathscr H)=\mathcal S_\epsilon(\mathscr H').}
			\end{equation}
			For each marked type the attained spectra form the coarsest exact
			summary for the stated tests. Every other summary sufficient for
			all the tests factors uniquely onto the spectrum image on its
			attained codes. A single two-chart test, with fixed targets and
			mutual coverage, separates any two distinct spectra. In the qualitative
			case the test need only add at most $|\mathcal M(Q)|$ readouts for the
			boundary image $Q$ of the separating record.
		\end{enumerate}
	\end{theorem}

	\begin{proof}
		\leavevmode
		\begin{enumerate}[label=\textup{(\roman*)},leftmargin=*]
			\item \emph{Invariance and composition.}
			Invariance follows by transporting every chosen realization and its maps.
			For a tuple of choices, every element of the matching group maps to $M$,
			and each member of $M$ has lifts obtained by independently lifting its
			coordinates. An arbitrary state of $H_i$ extends exactly when its
			boundary value belongs to $\operatorname{pr}_i M$. This gives
			\eqref{eq:contract-spectrum-coverage}; surjectivity along successive
			formation maps also preserves coverage of the chosen earlier charts.
			The kernel computations and recovery tests in
			\cref{thm:complete-compositional-invariants} depend only on the
			records and declared wiring. In particular each local closure and the
			least closure generated by their stipulated rules can be computed
			without knowing a representative group. Thus qualitative gates are
			record-determined. Both validity and failure of each dependency are
			retained, not merely positive dependencies.
			
			For the numerical gates, all joint range sizes are ratios of kernel
			orders, so a retained chart readout with prescribed target is onto if
			and only if its image cardinality equals that target's order. For the
			cut gate write $\ell,r,b$ for the two side and boundary homomorphisms
			on the candidate $W$, with $b$ factoring through each side. Put
			$A=\ker\ell$, $B=\ker r$, and $L=\ker b$. A function factors through
			both sides exactly when it is invariant under both normal kernels,
			equivalently under $AB$. Thus the intersection of the two represented
			function algebras is the represented boundary algebra exactly when
			$AB=L$. Since $AB\leq L$ and the groups are finite, this is equivalent to
			\begin{equation}\label{eq:contract-spectrum-cut}
				|L|\,|A\cap B|=|A|\,|B|.
			\end{equation}
			All four orders follow from the numerical record. If full boundary
			range is required, compare $|W|/|L|$ with the prescribed target order.
			The table law realizes every matching group. For a full matching graph
			and its complete reachable crossing boundary, the represented cut
			identity also follows directly from
			\cref{thm:arbitrary-graph-cut-information}.
			
			Consequently a tuple of records passes a gate exactly when any tuple
			of its realizing choices does. Every tuple of attained input records
			has simultaneous representatives because choices at different
			occurrences are independent. This proves the stated composition rule,
			including for empty families.
			
			\item \emph{Completeness and separation.}
			The implication from equality of spectra to contextual equivalence follows by induction on expressions using (i). For the converse, suppose a record $a$ belongs to exactly one spectrum;
			choose a realization of $a$ with boundary image $Q$. Attach the singleton
			normal-quotient chart of \cref{lem:covered-normal-extension}, retaining the	original readouts and the quotient observations. Prescribe the exact enlarged closure of that realization as a gate.
			In the numerical case also prescribe its order and the range sizes on $T\cup J$ for all $T\subseteq D$.
			These are fixed	constants of the continuation.
			Every choice with boundary image other than $Q$ fails mutual coverage. Among choices with image $Q$, the reconstruction formulas of \cref{lem:covered-normal-extension} imply that acceptance is equivalent to having precisely the record $a$.
			Conversely that record determines all prescribed values, so it passes.
			The	continuation therefore succeeds on exactly the family containing $a$.
			It is a single full matching attachment, and all its successful
			realizations retain every state of both chosen charts.
			This proves necessity and the separating claim. Finally any exact
			summary can identify only equivalent families. Sending its attained
			code to the spectrum of a representative is therefore well defined,
			onto the spectrum image, and unique. Every code in the domain of this factorization is attained.\qedhere
		\end{enumerate}
	\end{proof}

	\begin{proposition}[Rigidity of the independent fixed-kernel contract]
		\label{prop:independent-fixed-kernel-rigidity}
		Let $G_i$ be finite marked groups with boundary maps $p_i$ and fixed
		$K_i\leq\ker p_i$. Let $H\leq\prod_iG_i$ have projections $H_i$, each
		satisfying $H_i\cap\ker p_i=K_i$, and let $M=(\prod_i p_i)(H)$.
		\begin{enumerate}[label=\textup{(\roman*)},leftmargin=*]
			\item \emph{Independent kernels.}
			The following equivalence holds.
			\begin{equation}\label{eq:independent-fixed-kernel-rigidity}
				H\cap\ker(\textstyle\prod_i p_i)=\prod_i K_i
				\quad\Longleftrightarrow\quad
				H=\{(h_i)\in\textstyle\prod_i H_i:(p_i(h_i))_i\in M\}.
			\end{equation}
			
			\item \emph{Uniqueness over the prescribed images.}
			In particular, if a prescribed subgroup $C$ of the joint boundary
			is to be fully attained, this source is exactly the full matching
			source with constraint $C=M$. It is the unique such source over the
			chosen projections and prescribed boundary image. Conversely any such
			full matching with mutual coverage has the displayed independent
			kernel. If its new boundary $\rho$ is injective on $M$, its kernel over
			that new boundary is also $\prod_i K_i$.
		\end{enumerate}
	\end{proposition}

	\begin{proof}
		\leavevmode
		\begin{enumerate}[label=\textup{(\roman*)},leftmargin=*]
			\item \emph{Independent kernels.}
			The right-hand source is the preimage of $M$ under the joint map on
			$\prod_i H_i$, and its vertical kernel is $\prod_i K_i$. If the
			left-hand condition holds, any tuple in that preimage has the same
			boundary as some element of $H$. Their quotient lies in $\prod_i K_i$,
			which is contained in $H$, so the tuple belongs to $H$. The reverse
			implication follows by taking the kernel of the preimage.
			\item \emph{Uniqueness over the prescribed images.}
			Part (i) identifies the unique source on the prescribed projections and boundary image. Injectivity
			of $\rho$ on $M$ ensures that it introduces no new vertical elements.\qedhere
		\end{enumerate}
	\end{proof}
	
	\paragraph{Construction from fixed-kernel native data.}
	Given a finite group by its table and maps, one can compute its
	spectrum by enumerating the local subgroups, testing \eqref{eq:contract-fixed-k-family} and the local
	filters, and computing the agreement records of the accepted subgroups.
	The following formula gives a more structural construction, avoiding
	subgroup or complement enumeration in the elementary abelian
	normalizer case. Unlike an enumerator of kernel orders alone, it
	retains the simultaneous actual boundary images for every premise.
	
	\begin{theorem}[Joint fixed-kernel spectrum by a normalizer and rank transform]
		\label{thm:joint-fixed-kernel-spectrum}
		Fix $G,p,Q,K$ as in \eqref{eq:contract-fixed-k-family}. Put
		$N_0=\ker p$, $M=N_{p^{-1}(Q)}(K)$, and $D_T=\bigcap_{d\in T}\ker f_d$,
		where $D_\varnothing=G$.
		\begin{enumerate}[label=\textup{(\roman*)},leftmargin=*]
			\item \emph{Complements and their records.}
			If $p(M)\ne Q$ the family is empty. Otherwise
			put $E=M/K$, $A=(M\cap N_0)/K$, and write $\bar p:E\twoheadrightarrow Q$.
			Then the family is in bijection with complements to $A$ in $E$.
			For the homomorphic section $s:Q\to E$ of a complement, its inverse
			image $H$ has record
			\begin{equation}\label{eq:fixed-kernel-spectrum-record}
				\begin{aligned}
					B_H(T)&=\{q\in Q:s(q)\in (M\cap KD_T)/K\},\\
					m_H(T)&=|K\cap D_T|,\\
					V_H(T)&=\{d:K\cap D_T\leq\ker f_d\}.
				\end{aligned}
			\end{equation}
			In particular the vertical part of either record is fixed by $K$;
			all variability of the record is in the simultaneous boundary images.
			\item \emph{Enumeration by linear ranks.}
			Suppose additionally that $A$ is elementary abelian over $\mathbb F_l$
			for a prime $l$ (including the trivial group). The complete numerical
			spectrum, and hence the qualitative spectrum, admits a rank-only
			finite construction. Choose a normalized set section $s_0:Q\to E$;
			write $s(q)=z(q)s_0(q)$, using additive coordinates on $A$. Let
			$\Omega=\{(T,q):T\subseteq D,\ q\in Q\}$ index membership atoms.
			For $U\subseteq\Omega$, let $L_U z=b_U$ be the affine linear system
			consisting of the section equations and the requirements
			$s(q)\in (M\cap KD_T)/K$ for $(T,q)\in U$. An impossible atom makes
			the system inconsistent. With $n=|Q|\dim_{\mathbb F_l}A$, put
			\begin{equation}\label{eq:fixed-kernel-spectrum-rank}
				a(U)=\begin{cases}
					l^{n-\operatorname{rank}L_U},& L_Uz=b_U\text{ is consistent},\\
					0,&\text{otherwise}.
				\end{cases}
			\end{equation}
			For each proposed joint membership pattern $S\subseteq\Omega$, its
			multiplicity is
			\begin{equation}\label{eq:fixed-kernel-spectrum-mobius}
				b(S)=\sum_{U\subseteq\Omega\setminus S}(-1)^{|U|}a(S\cup U).
			\end{equation}
			A pattern belongs to the spectrum if and only if $b(S)>0$, with
			$B_H(T)=\{q:(T,q)\in S\}$ and the fixed vertical data of
			\eqref{eq:fixed-kernel-spectrum-record}. Multiplicity here counts
			actual subgroups, not conjugacy classes; the spectrum forgets it.
			All counts and the spectrum are independent of section and coordinates.
			Local filters determined by the record may be applied directly to the
			resulting patterns.
		\end{enumerate}
	\end{theorem}
	\begin{proof}
		\leavevmode
		\begin{enumerate}[label=\textup{(\roman*)},leftmargin=*]
			\item \emph{Complements and their records.}
			Apply the general normalizer--complement theorem
			(\cref{thm:chief-normalizer-complement}) to $p:G\to p(G)$ and the
			requested subgroup $Q$. If $Q\not\le p(G)$ or $p(M)\ne Q$ there is no
			lift; otherwise that theorem gives the stated bijection. For the
			normal observation kernel $D_T$, the image identity in
			\cref{lem:normalizer-complement-kernel-count} identifies
			$(H\cap D_T)K/K$ with $s(Q)\cap (KD_T\cap M)/K$. Projecting to $Q$
			gives $B_H(T)$. The vertical agreement kernel is
			$H\cap D_T\cap\ker p=K\cap D_T$, giving both remaining entries of
			\eqref{eq:fixed-kernel-spectrum-record}.
			\item \emph{Enumeration by linear ranks.}
			When $A$ is elementary abelian, conjugation on $A$ gives a well-defined
			$Q$-action, and the factor set $\alpha$ of $s_0$ gives the affine section
			equations $z(1)=0$ and
			\[
			z(qu)=z(q)+q\cdot z(u)+\alpha(q,u)\qquad(q,u\in Q).
			\]
			For a fixed $q$, the intersection of the subgroup $(M\cap KD_T)/K$
			with the fiber above $q$ is either empty or a coset of its intersection
			with $A$. Under the coordinates $z(q)s_0(q)$ this is an affine
			subspace of $A$. Thus every positive atom is an affine linear system,
			and \eqref{eq:fixed-kernel-spectrum-rank} counts precisely the sections
			satisfying at least the atoms in $U$. Inclusion-exclusion over all
			atoms outside $S$ gives exactly the sections satisfying precisely $S$,
			which proves \eqref{eq:fixed-kernel-spectrum-mobius}. Each such section
			corresponds to one complement and hence one $H$. This interpretation
			shows coordinate independence. The consistency and rank tests require
			only the native data and the finite normalizer extension. They refine
			the order-based complement counting in
			\cref{thm:cohomological-contract-enumerator} by recording the actual
			joint image subgroups, rather than just their sizes.\qedhere
		\end{enumerate}
	\end{proof}
	
	\begin{example}[The ambient record and marginal spectra are insufficient]
		\label{ex:contract-spectrum-obstructions}
		Let $P=Q=C_2$, with no private descriptors, and compare $C_4$ and
		$C_2\times C_2$ with their quotient maps onto $C_2$ of kernel order two.
		Their ambient numerical records coincide: the only boundary image is
		$C_2$ and the vertical multiplicity is two. For the fixed trivial
		kernel request $K=1$, the former has no lift onto $Q$, since its only
		order-two subgroup lies in the kernel; the latter has such lifts.
		The spectra are respectively empty and nonempty. Thus even the
		ambient numerical invariant cannot classify the fixed-kernel contract.
		
		For a different failure, use vector spaces over $\mathbb F_2$, with
		$Q=\mathbb F_2$, $K=0$, and $D=\{1,2\}$. In the first group
		$G=\mathbb F_2^2$, put $p(x,y)=x$ and $f_1(x,y)=f_2(x,y)=y$.
		In the second group $G'=\mathbb F_2^3$, put $p'(x,y,z)=x$,
		$f'_1(x,y,z)=y$, $f'_2(x,y,z)=z$. Lifts are graphs of linear maps
		from $Q$ to the vertical space. In both cases $m_H(T)=1$ and
		$V_H(T)=D$ for every $T$. For each nonempty premise $T$ the separately
		attainable $B_H(T)$ are precisely $0$ and $Q$, while the empty premise
		always has image $Q$. Nevertheless in the first family
		$B_H(\{1\})=B_H(\{2\})$ for every lift; in the second a lift can have
		$B_H(\{1\})=0$ and $B_H(\{2\})=Q$. The simultaneous spectra differ,
		and \cref{thm:contract-filtered-spectrum} supplies a coverage-preserving
		Boolean test separating them. Independent marginal feasibility does
		not imply joint feasibility of recovery targets.
	\end{example}
	
	\subsection{Complete invariants for literally coupled hidden kernels}
	\label{subsec:coupled-kernel-invariants}
	
	We now prescribe a subgroup of the product of the hidden kernels and
	allow this subgroup to be proper. To specify it unambiguously, we fix
	an identification of each kernel with a named group. With these
	markings a diagonal coupling detects the extension itself and the
	kernels of all retained observations. For central elementary abelian
	kernels the attainable records form an affine space before filtering.
	
	\begin{definition}[Marked kernel coupling language]
		\label{def:marked-kernel-coupling}
		Use the finite group-table context of \cref{def:fragment} and the
		family and gate conventions of \cref{def:contract-filtered-fragment}.
		The \textbf{marked kernel coupling language} is defined by the following
		realizations, coupling operations, output markings, and continuations.
		\begin{enumerate}[label=\textup{(\roman*)},leftmargin=*]
			\item \emph{Marked realizations.}
			In addition to $P$ and the labeled readout types, fix a finite group
			$K$ as part of the marked type. Each admissible realization consists of
			$G,p,(f_d)_{d\in D}$ together with an isomorphism
			$\iota:K\to\ker p$. Its image $Q=p(G)\leq P$ need not be prescribed
			in advance. The readouts are homomorphisms; their private values are
			used only through equality, cardinalities and coverage of their declared
			targets, not through tests of their actual labels.
			\item \emph{Permitted couplings.}
			A formation chooses independently one realization from each input
			family. Write $K_i$ for the named kernels and $Q_i=p_i(G_i)$.
			It supplies a boundary constraint $C\leq\prod_i P_i$ and a literal
			subgroup $L\leq\prod_i K_i$ projecting onto every $K_i$. Put
			$M=C\cap\prod_i Q_i$ and $W=(\prod_i p_i)^{-1}(M)$ in $\prod_i G_i$.
			The formation now ranges over \emph{all} subgroups $H\leq W$ satisfying
			\begin{equation}\label{eq:coupled-language-contract}
				(\textstyle\prod_i p_i)(H)=M,\qquad
				H\cap\textstyle\prod_i\ker p_i=(\textstyle\prod_i\iota_i)(L),
				\qquad \operatorname{pr}_i H=G_i\quad\text{for all }i.
			\end{equation}
			Thus the joint boundary must be fully attained, and even when $L$ is
			proper all local vertical kernels are fully represented. Any desired
			joint image can be named by taking $C$ equal to it. The gates of
			\cref{def:contract-filtered-fragment} are applied to each candidate.
			\item \emph{Output readouts and markings.}
			Readouts may be retained or forgotten and homomorphic readouts of the
			joint boundary may be introduced. A new boundary $\rho:C\to P'$ may
			be prescribed. Its marked kernel group $K'$ is fixed in the output
			type; for each successful source retain all isomorphisms
			$\iota':K'\to\ker(\rho\circ(\prod_i p_i)|_H)$, subject to any prescribed
			compatibility with the inclusion of $L$. More precisely, the formation
			may specify an injective homomorphism $\lambda:L\to K'$ and require $\iota'\lambda$
			to be the indicated inclusion of $L$ into $H$. When $\rho$ is injective
			on $M$, one can take $K'=L$ with its inherited marking. If no permitted
			mark exists the choice fails. The empty family represents failure.
			\item \emph{Continuations and equivalence.}
			The language includes finite nested formations of this kind, the
			identity continuation, and the diagonal tests below, with arbitrary
			finite reference fragments of the same marked type. It permits both
			independently prescribed recovery targets and targets generated by
			local laws, as in the preceding language. Equivalence means equality
			of Boolean success in every common continuation.
		\end{enumerate}
	\end{definition}
	
	The prescribed kernels belong to the chosen native realizations.
	Full matching with independent vertical kernels corresponds to $L=\prod_i K_i$.
	
	\begin{definition}[Decorated extension records and their attainable spectrum]
		\label{def:decorated-extension-record}
		For a marked realization put $N_d=\ker f_d$.
		\begin{enumerate}[label=\textup{(\roman*)},leftmargin=*]
			\item \emph{Section coordinates.}
			Choose a normalized set
			section $s:Q\to G$ of $p:G\to Q$, so $s(1)=1$. Identifying $K$ with its
			marked image, define
			\[
			\theta_q(a)=s(q)a s(q)^{-1},\qquad
			\alpha(q,r)=s(q)s(r)s(qr)^{-1},\qquad
			F_d(q)=\{a\in K:a s(q)\in N_d\}.
			\]
			Here $\theta_q\in\operatorname{Aut}(K)$, $\alpha(q,r)\in K$, and each
			$F_d(q)$ is empty or a coset of $K\cap N_d$.
			\item \emph{Record, change of section, and attainable spectrum.}
			The \textbf{decorated extension record} consists of
			$Q$ and the joint orbit of $(\theta,\alpha,(F_d)_{d\in D})$ under
			changes of normalized section. Explicitly, for a map $u:Q\to K$ with
			$u(1)=1$, the replacement $s'(q)=u(q)s(q)$ gives
			\begin{equation}\label{eq:decorated-extension-gauge}
				\begin{aligned}
					\theta'_q&=\operatorname{Inn}(u(q))\theta_q,\\
					\alpha'(q,r)&=u(q)\theta_q(u(r))\alpha(q,r)u(qr)^{-1},\\
					F'_d(q)&=F_d(q)u(q)^{-1}.
				\end{aligned}
			\end{equation}
			Write $\mathcal E(G)$ for this decorated record and
			$\mathcal S_{\rm ext}(\mathscr H)=\{\mathcal E(G):G\in\mathscr H\}$
			for its attainable spectrum. The orbit is a finite equivalence class
			of concrete data, not an equivalence defined by continuations.
		\end{enumerate}
	\end{definition}
	
	\begin{theorem}[Diagonal separation and complete coupled-kernel spectra]
		\label{thm:coupled-extension-spectrum}
		For the language of \cref{def:marked-kernel-coupling}, including arbitrary marked reference fragments and literal boundary and kernel diagonals, the following hold.
		\begin{enumerate}[label=\textup{(\roman*)},leftmargin=*]
			\item \emph{Equality of records.}
			Two realizations of the same marked type have equal decorated
			extension records if and only if there is a group isomorphism $\varphi$
			between them commuting with the boundary maps, intertwining the kernel
			markings, and transporting each named readout kernel to its counterpart.
			Such an isomorphism need not preserve actual private readout values.
			\item \emph{Composition.}
			The attainable decorated extension spectra have an explicit
			composition operation. For each tuple of records, reconstruct the
			marked extensions and readout kernels from their coordinates. Reject
			if $M$ does not project onto every $Q_i$. In the reconstructed $W$
			identify $L$ with its marked image, and put
			\begin{equation}\label{eq:coupled-normalizer}
				U=N_W(L),\qquad E=U/L,\qquad A=(U\cap\textstyle\prod_i K_i)/L.
			\end{equation}
			If the joint boundary map on $U$ is not onto $M$, reject. Otherwise
			form the complements to $A$ in $E$; their inverse images in $U$ are
			precisely the sources of \eqref{eq:coupled-language-contract}.
			Evaluate the gates, the output boundary and permitted output markings,
			and collect the decorated records of the successful outputs. The
			operation is independent of representatives and applies to arbitrary
			finite matching graphs and nested formations in the stated language.
			\item \emph{Separation by a diagonal test.}
			There is exact Boolean classification
			\begin{equation}\label{eq:coupled-spectrum-completeness}
				\boxed{\mathscr H\equiv_{\rm coupled}\mathscr H'
					\quad\Longleftrightarrow\quad
					\mathcal S_{\rm ext}(\mathscr H)=\mathcal S_{\rm ext}(\mathscr H').}
			\end{equation}
			A single diagonal test with mutual chart coverage separates distinct
			spectra. Only the prescribed boundary and hidden-kernel diagonals and
			mutual determination of corresponding readouts are needed; numerical
			gates and independently prescribed full recovery targets are unnecessary
			for this separation. In particular, for this language the qualitative
			and numerical variants of the Boolean equivalence coincide.
			Every sufficient behavioral summary factors uniquely onto the attained
			decorated spectra. The classification uses native extension data and
			finite subgroup computations, not enumeration of continuations.
		\end{enumerate}
	\end{theorem}
	\begin{proof}
		\leavevmode
		\begin{enumerate}[label=\textup{(\roman*)},leftmargin=*]
			\item \emph{Equality of records.}
			In section coordinates $(a,q)\leftrightarrow a s(q)$, multiplication is
			\[
			(a,q)(b,r)=(a\theta_q(b)\alpha(q,r),qr),
			\]
			and the readout kernel $N_d$ is exactly the union of the fibers
			$F_d(q)\times\{q\}$. Thus the data reconstruct a marked extension with
			its distinguished normal subgroups. A change of section changes the
			coordinates by the indicated bijection and yields
			\eqref{eq:decorated-extension-gauge}. Conversely an isomorphism over
			the markings and the boundary transports a section to a section,
			and the two sections in the target differ by a unique normalized
			map into the named kernel. Transport of the named normal subgroups
			is exactly the last equation in \eqref{eq:decorated-extension-gauge}.
			
			\item \emph{Composition.}
			Write $p_*=\prod_i p_i:W\twoheadrightarrow M$.
			Apply \cref{thm:chief-normalizer-complement} to this extension with
			prescribed hidden kernel $L$ and full boundary image $M$.
			It gives exactly \eqref{eq:coupled-normalizer} and the asserted
			bijection, provided $p_*(U)=M$. The application is to the joint
			matching group, not to its factors separately.
			When both $M$ and $L$ project onto each corresponding local group,
			any $g_i\in G_i$ has the same boundary as the $i$th coordinate of some
			$h\in H$. Multiplying $h$ by a suitable element of $L$ gives $g_i$
			in that coordinate. Thus all these lifts cover the charts. The
			converse necessity of the boundary coverage is immediate; vertical
			coverage is imposed in the language itself.
			
			All retained readout kernels on a candidate are intersections of $H$
			with inverse images of the original $N_d$. New readout kernels and the
			output boundary are determined by the joint boundary map. For any
			homomorphic observations on a finite group $H$, write $N_T$ for the
			intersection of their kernels in a premise $T$ (and $N_\varnothing=H$).
			Determination $T\to d$ is equivalent to $N_T\leq N_d$; the conditional
			range of $J$ at any attained $T$-value is $|N_T|/|N_{T\cup J}|$.
			Coverage of a declared target is decided by its cardinality and the
			index of the corresponding kernel. For two side readouts with kernels
			$N_1,N_2$, the equivalence generated by their agreement relations is
			agreement modulo $N_1N_2$, since these subgroups are normal. Hence their
			represented common information equals that of a specified cut readout
			precisely when its kernel is $N_1N_2$. These tests, and the generated
			closure from local laws, are determined by the reconstructed data.
			Marked isomorphisms transport the joint matching group, $L$, all its
			candidate subgroups, their gates and all permitted output markings
			bijectively. This proves the stated composition and its invariance.
			Repeated application gives the result for nested formations.
			\item \emph{Separation by a diagonal test.}
			Equality of spectra implies equality after every formation,
			by (ii), and thus equality of terminal success. For the converse,
			choose a record in the symmetric difference and a realization $R$ of
			it. Attach the fixed reference $R$ to the hole by $C=\Delta P$ and
			$L=\Delta K$ in the two literal markings. Require both projections
			of $H$ onto its chosen input $G$ and onto $R$ to be surjective. If the
			two boundary images differ, mutual coverage is impossible. Otherwise
			the boundary image must be their diagonal. An element of $H$ in the
			kernel of either chart projection has both boundary coordinates
			trivial and lies in $\Delta K$. Its other coordinate is therefore
			trivial as well. Both chart projections are isomorphisms, and $H$
			is the graph of an isomorphism $G\to R$ commuting with the boundary
			and fixing the kernel markings. Conversely the graph of any such
			isomorphism satisfies the prescriptions. For every $d$, require the
			readout of $G$ and the corresponding readout of $R$ to determine
			one another on $H$. On this graph the two implications are exactly
			transport of $N_d$ onto its counterpart. The inherited diagonal
			kernel supplies the required output marking. Hence the test succeeds
			if and only if the input spectrum contains the chosen record, by (i).
			This separates unequal spectra, including the empty spectrum. It uses
			only recovery implications; no range sizes are required. The same
			argument shows that if two exact behavioral summaries agree then their
			spectra agree. Sending each attained summary code to the spectrum of
			a representative is well-defined, surjective and uniquely determined.\qedhere
		\end{enumerate}
	\end{proof}
	
	\begin{corollary}[Exact count of central extensions distinguished by kernel coupling]
		\label{cor:coupled-extension-noncollapse}
		For a fixed finite group $Q$ and a central named kernel
		$K=\mathbb F_l^k$, with boundary image $Q$ and no private readouts, the following hold.
		\begin{enumerate}[label=\textup{(\roman*)},leftmargin=*]
			\item \emph{Numerical record.}
			All marked central extensions have the same numerical record of
			\cref{thm:complete-compositional-invariants}.
			\item \emph{Number of distinguishable extensions.}
			Their singleton
			spectra in the coupled language have exactly
			\begin{equation}\label{eq:coupled-extension-count}
				l^{k\dim_{\mathbb F_l}H^2(Q,\mathbb F_l)}
			\end{equation}
			distinct Boolean equivalence classes. Any deterministic sufficient
			summary on this family must distinguish that many values, and the
			marked central extension class attains the bound.
		\end{enumerate}
	\end{corollary}
	\begin{proof}
		\leavevmode
		\begin{enumerate}[label=\textup{(\roman*)},leftmargin=*]
			\item \emph{Numerical record.}
			The numerical record has only the boundary image $Q$ and the vertical
			multiplicity $|K|$.
			\item \emph{Number of distinguishable extensions.}
			For a central extension the action $\theta$ is
			trivial; the orbit of the normalized cocycle $\alpha$ in
			\eqref{eq:decorated-extension-gauge} is its class in $H^2(Q,K)$,
			with trivial action. Conversely each normalized cocycle constructs
			an extension by the displayed multiplication, with the prescribed
			markings. Coordinatewise cochains identify $H^2(Q,K)$ with $k$ copies
			of $H^2(Q,\mathbb F_l)$. Apply \cref{thm:coupled-extension-spectrum}.\qedhere
		\end{enumerate}
	\end{proof}
	
	\paragraph{An affine composition law for central elementary abelian kernels.}
	The extension coordinates give the attainable marked extensions directly, supplementing the boundary images in \cref{thm:joint-fixed-kernel-spectrum}.
	
	\begin{theorem}[Affine decorated spectrum of a prescribed central coupling]
		\label{thm:affine-coupled-spectrum}
		Fix a tuple of marked realizations, whose named kernels $K_i$ are central
		and elementary abelian over the same field $k=\mathbb F_l$. Write
		$A=\bigoplus_i K_i$ (in this theorem $A$ denotes the entire vertical
		space). Fix a subdirect joint boundary $M\leq\prod_i Q_i$ and a subdirect
		subspace $L\leq A$. Retain the full joint boundary of the output, with
		kernel marking $L$. The retained private readouts are a specified subset
		of the occurrence-labeled input readouts, together with any specified
		new homomorphic observations of $M$. Choose normalized central cocycles
		$\alpha_i$ for the inputs and write $\alpha$ for their coordinatewise
		restriction to $M$. Let $j:A/L\to A$ be any linear section of the vector
		space quotient. The following statements hold.
		\begin{enumerate}[label=\textup{(\roman*)},leftmargin=*]
			\item \emph{Existence and parametrization.}
			The sources satisfying \eqref{eq:coupled-language-contract} are
			in bijection with the solutions of the normalized affine cochain system
			\begin{equation}\label{eq:coupled-cochain-system}
				z(1)=0,\qquad \delta z=-\bar\alpha,
				\qquad z:M\to A/L,
			\end{equation}
			where $\delta z(q,r)=z(q)+z(r)-z(qr)$ and the bar denotes passage to
			$A/L$. If this system is inconsistent the output family is empty.
			Otherwise its solution space is a translate of
			$\operatorname{Hom}(M,A/L)$.
			\item \emph{Output cocycle and readouts.}
			Write the input kernel of an observation $d$ at chart $i$ as
			fibers with support $B_{i,d}=p_i(N_d)$, vertical subspace
			$L_{i,d}=K_i\cap N_d$, and coset coordinate
			$v_{i,d}(q_i)\in K_i/L_{i,d}$ for $q_i\in B_{i,d}$.
			For each solution $z$, put $t=jz$. The output cocycle and retained
			readout data are
			\begin{equation}\label{eq:coupled-affine-record}
				\begin{aligned}
					\beta_z&=\alpha+\delta t\in C^2_0(M,L),\\
					B_d^H&=M\cap\operatorname{pr}_i^{-1}(B_{i,d}),\\
					L_d^H&=L\cap\operatorname{pr}_i^{-1}(L_{i,d}),\\
					v_d^H(q)&=\bar\pi_d^{-1}
					\big(v_{i,d}(q_i)-\operatorname{pr}_i t(q)\bmod L_{i,d}\big),
					\qquad q\in B_d^H,
				\end{aligned}
			\end{equation}
			where $\bar\pi_d:L/L_d^H\to K_i/L_{i,d}$ is the induced isomorphism.
			For a new readout $h$ on $M$, take $B_h^H=\ker h$, $L_h^H=L$ and
			$v_h^H=0$. The notation $C^2_0$ denotes normalized cochains; the
			resulting cocycle is automatically a normalized $2$-cocycle.
			\item \emph{Attainable records and multiplicities.}
			Before gates are applied, the output decorated records constitute
			an affine subspace of an explicit vector space quotient. With the
			$B_d^H,L_d^H$ fixed as above, form
			\begin{equation}\label{eq:coupled-affine-quotient}
				\mathcal V=\left(C^2_0(M,L)\oplus
				\bigoplus_d C^1_0(B_d^H,L/L_d^H)\right)
				\big/\{(\delta u,(-u|_{B_d^H}\bmod L_d^H)_d):u\in C^1_0(M,L)\}.
			\end{equation}
			Here $C^1_0$ denotes functions vanishing at the identity. The image
			of $z\mapsto(\beta_z,(v_d^H)_d)$ in $\mathcal V$ is exactly the
			attainable decorated spectrum for this tuple. If the system has a
			solution $z_0$ and $F$ denotes this affine map on the solution space,
			put
			\begin{equation}\label{eq:coupled-affine-rank}
				r=\operatorname{rank}\big(v\mapsto F(z_0+v)-F(z_0),\quad
				v\in\operatorname{Hom}(M,A/L)\big).
			\end{equation}
			Then precisely $l^r$ distinct records occur, and each occurs for
			$l^{\dim\operatorname{Hom}(M,A/L)-r}$ distinct native subgroups.
			The rank and the spectrum are intrinsic under the marked comparison.
			\item \emph{Preservation tests and counts.}
			All recovery and numerical gates on the output can be evaluated
			on this affine quotient. More explicitly, for a premise $T$ write
			$L_T=\bigcap_{d\in T}L_d^H$ (empty intersection $L$) and let
			$J_T:L\to\bigoplus_{d\in T}L/L_d^H$ be the joint quotient map. Its
			boundary agreement image is
			\begin{equation}\label{eq:coupled-affine-memberships}
				B_H(T)=\{q\in\textstyle\bigcap_{d\in T}B_d^H:
				(v_d^H(q))_{d\in T}\in\operatorname{im} J_T\},
				\qquad |K_H(T)|=|B_H(T)|\,|L_T|.
			\end{equation}
			The empty premise has support $M$. Every membership atom is affine
			linear (or impossible) in the solution coordinates and invariant
			under the gauge quotient. Thus the gates select a finite Boolean
			combination of affine linear loci in the attainable affine space.
			Exact counts of the surviving decorated records can be obtained by
			linear ranks and finite inclusion--exclusion, without enumerating
			complements or continuations. For families of inputs, take the union
			of the successful spectra over tuples of their decorated records.
			No polynomial bound on the size of this union or on the Boolean
			filter is asserted.
		\end{enumerate}
	\end{theorem}
	\begin{proof}
		\leavevmode
		\begin{enumerate}[label=\textup{(\roman*)},leftmargin=*]
			\item \emph{Existence and parametrization.}
			In the product section coordinates, $W=A\times_\alpha M$.
			Since the kernel is central, $L$ is normal in $W$. A source of the
			specified image and kernel is exactly the inverse image of a section
			of $W/L\to M$. Such a section has the form $q\mapsto(z(q),q)$ in
			quotient coordinates and is a homomorphism precisely under
			\eqref{eq:coupled-cochain-system}. Distinct sections give distinct
			complements and therefore distinct subgroups. The homogeneous
			solutions are the homomorphisms from $M$ to the additive group $A/L$.
			Surjectivity of $M$ and $L$ supplies chart coverage by the proof of
			\cref{thm:coupled-extension-spectrum}.
			\item \emph{Output cocycle and readouts.}
			For a solution, use the section of $H\to M$ given by $q\mapsto(t(q),q)$.
			The elements of $H$ are precisely $(a+t(q),q)$ for $a\in L$.
			Their multiplication in these new coordinates has cocycle
			$\alpha+\delta t$, which takes values in $L$ by the section equation.
			For a retained readout from chart $i$, membership in its kernel is
			exactly the condition $q_i\in B_{i,d}$ and
			$\operatorname{pr}_i a+\operatorname{pr}_i t(q)\equiv
			v_{i,d}(q_i)\pmod{L_{i,d}}$. Since $L$ surjects onto $K_i$ this has a
			solution $a$ for every such $q$; its solutions form a coset modulo
			$L_d^H$. This gives the asserted isomorphism and the readout formulas. Boundary
			readouts depend only on $q$, giving their stated data.
			\item \emph{Attainable records and multiplicities.}
			With fixed support and vertical subspaces the last two parts of
			\eqref{eq:decorated-extension-gauge} in additive notation identify
			exactly $(\beta,(v_d))$ and $(\beta+\delta u,(v_d-u|_{B_d}))$.
			The action is trivial. Thus the quotient in
			\eqref{eq:coupled-affine-quotient}, restricted to valid data, gives
			precisely the decorated equivalence classes. Formula
			\eqref{eq:coupled-affine-record} is affine linear in $z$ on the solution
			space. Its quotient image is affine and has dimension $r$; fibers of
			that affine map have constant size by rank--nullity. Changing any
			input section, coordinate system, or linear lift $j$ merely changes
			representatives of the same marked sources and output records. In
			particular the number and the spectrum are independent of all choices.
			
			\item \emph{Preservation tests and counts.}
			Simultaneous kernel membership at a boundary value is
			exactly the simultaneous coset system $J_T(a)=(v_d^H(q))_d$.
			Every nonempty solution fiber is a coset of $L_T$, proving the
			formula. A gauge change translates the coset vector by $J_T(-u(q))$
			and so leaves membership in its image invariant. The latter is a
			linear-subspace membership condition on an affine vector, and hence
			an affine linear locus also on the quotient image (equivalently, use
			any linear complement of the gauge subspace). Determination is
			$|K_H(T)|=|K_H(T\cup\{d\})|$, since the latter is contained in the former;
			conditional ranges and source order follow by the kernel formula in
			\cref{thm:coupled-extension-spectrum}. Coverage of declared targets and
			represented common information of named homomorphic readouts can
			also be tested by these sizes: when the cut observation factors
			through both sides, equality with the product of side kernels is
			equivalent to equality of orders, and the product order is
			$|N_1||N_2|/|N_1\cap N_2|$. If factorization is not stipulated it can be
			tested by the corresponding determinations. Local-law targets for
			this fixed tuple are known. Each gate is consequently determined by a
			finite membership pattern. Intersections of positive membership
			atoms with the attainable affine space are affine systems; their
			solution counts are powers of $l$ or zero by rank--nullity. Applying
			inclusion--exclusion to the absent atoms counts each complete pattern,
			as in \cref{thm:joint-fixed-kernel-spectrum}, but now on the quotient
			space of distinct decorated records. This proves (iv).\qedhere
		\end{enumerate}
	\end{proof}
	
	\begin{example}[A strictly richer spectrum with fixed central coupling]
		\label{ex:coupled-affine-separation}
		Let $k=\mathbb F_l$, and take two split charts $G_i=k^2$ with boundary
		$p_i(x,y)=x$, named kernel $y\in k$, and one readout $f_i(x,y)=y$.
		Require diagonal boundary and diagonal vertical kernel. The compatible
		sources are precisely
		\[
		H_a=\{((x,y),(x,y+ax)):x,y\in k\},\qquad a\in k.
		\]
		Indeed the complement equation for $A/L\cong k$ is linear with
		one-dimensional solution space. In output coordinates $(x,y)$ the
		readouts are $y$ and $y+ax$. A marked isomorphism between two such
		outputs fixes $x$ and the vertical coordinate at $x=0$, and has the
		form $(x,y)\mapsto(x,y+bx)$. Preservation of the first readout kernel
		forces $b=0$, and preservation of the second then forces the parameters
		to agree. The affine quotient therefore has rank one and exactly $l$
		records. In contrast, the numerical record of
		\cref{thm:complete-compositional-invariants} is the same for all
		$a\ne0$: for either singleton premise the boundary agreement image
		is $k$ and the vertical kernel is zero; for both readouts the boundary
		image is zero and the vertical kernel is zero. The empty-premise data
		are fixed. For $a=0$ the double-premise boundary image is $k$ instead.
		Thus for $l>2$ the $l-1$ distinct nonzero parameters collapse to one
		numerical record of the uncoupled language.
		A diagonal test as in \cref{thm:coupled-extension-spectrum} separates
		every pair. Even without readouts, the extensions $C_4\to C_2$ and
		$C_2\times C_2\to C_2$ with their full named kernels of order two give
		equal singleton spectra in the uncoupled language and different
		coupled spectra. The test
		with reference $C_4$ succeeds only on the first.
	\end{example}
	
	\begin{remark}[What the kernel marking determines]
		\label{rem:coupled-extension-scope}
		The diagonal test uses the kernel markings and a reference fragment of
		each type. Its graph argument is the isomorphism case of Goursat's lemma
		\cite{Goursat1889,BauerSenZvengrowski2015}. The marked records supplement
		the classical extension coordinates \cite{BrownCohomology1982} with the named readout kernels.
		Retaining the full joint boundary makes $L$ the output kernel in the affine formula.
		For a noninjective boundary map, \cref{thm:coupled-extension-spectrum} also accounts for the larger output kernel and its allowed markings.
	\end{remark}
	
	\section{Source minima and rigidity}
	\label{sec:faithful-tree}
	
	For a specified atlas, the minimum is taken over the sources
	admitted by its formation law and preservation contract. On a finite
	Mal'cev join tree, chart coverage and exact common information across
	each full tree cut force the entire matching algebra as the observable
	image (\cref{thm:native-source-rigidity}). The paired-witness test
	of \cref{thm:faithful-matching-minimum} determines whether that
	minimum also gives the prescribed recovery. Varying the source law
	or the required observations changes the minimum, even with the same
	charts.
	
	\paragraph{Admissible classes and optimization objectives.}
	\begin{enumerate}[label=(\roman*),nosep]
		
		\item \emph{Native algebraic sources over prescribed charts.}
		A finite algebra $\mathbf H$ comes with compatible homomorphisms
		$h_i:H\to A_i$ into the actual chart algebras. The cost is $|H|$.
		Depending on the theorem, the contract uses state exactness and
		full-cut separator exactness, or state exactness and recovery
		faithfulness under simple descriptor images.
		
		\item \emph{Closure-only representations.}
		A finite relation $R$ over arbitrary column alphabets represents
		a specified determination closure. Here $|R|$ counts distinct rows,
		and chartwise exactness is weaker than faithfulness to the least
		global closure.
		
		\item \emph{Subamalgams of fixed native spectra.}
		The actual local spectral state sets and overlap maps are fixed,
		and a candidate $R\subseteq X$ must retain declared local states
		and non-recovery witnesses. The cost is its number of global sectors;
		additional separator contracts are imposed separately when required.
		For a single complete-side cut, the weaker relational contract has
		an exact spanning-tree minimum in
		\cref{thm:one-cut-spanning-minimum}.
		
		\item \emph{Sources under a different native law.}
		The binary meet-semilattice example minimizes $|H|$ among
		sources preserving the meet operation and the prescribed chart maps.
		
	\end{enumerate}
	
	The same cyclic atlas has sharp costs $8$, $5$, and $4$ under native
	linear assembly, fixed-chart relational assembly, and closure-only
	representation, respectively (\cref{prop:cyclic-triangle-source-minima}).
	
	\subsection{Native source minima under full-cut preservation}
	\label{subsec:native-algebraic-minima}
	
	For finite Mal'cev algebras, the canonical map to the full join classifies sources that cover every chart and satisfy the exact common-information condition at every full tree cut.
	
	\subsubsection{Sharp full-cut source minimum and equality cases}
	\label{subsec:sharp-native-source-resources}
	
	Suppose the native equational class $\mathcal V$ has a Mal'cev term,
	a ternary term $m$ satisfying
	\[
	m(x,y,y)=x,\qquad m(x,x,y)=y.
	\]
	Groups, modules, and rings have such a term. Kernel congruences in
	$\mathcal V$ permute \cite[Theorems~II.5.9 and~II.12.2]{BurrisSankappanavar1981}.
	Fleischer's fiber-product lemma describes a subdirect product over a
	common quotient \cite{Fleischer1955} (see also
	\cite[Lemma~IV.10.1]{BurrisSankappanavar1981}). At each full tree cut,
	exact common information identifies this quotient with the prescribed
	separator. Propagating the identification across the tree determines
	the minimum source and its equality cases.
	
	Fix an atlas of finite algebras $\mathbf A_i\in\mathcal V$ with
	homomorphic descriptor readouts, equal actual separator image algebras
	$\mathbf Y_a$, and a join tree $J$, as in
	\cref{thm:equational-native-amalgamation}.  Let $I$ be its nonempty
	chart set and let $\mathbf W$ be its constructed matching algebra,
	with underlying set
	\[
	W=\Bigl\{(a_i)\in\prod_{i\in I}A_i:
	\text{all specified overlap readouts agree}\Bigr\}.
	\]
	A \emph{native source over this atlas} consists of a finite algebra
	$\mathbf H\in\mathcal V$ and compatible homomorphisms
	$h_i:\mathbf H\to\mathbf A_i$.  They induce the canonical
	homomorphism
	\begin{equation}\label{eq:native-source-comparison-map}
		\psi:\mathbf H\longrightarrow\mathbf W,
		\qquad \psi(h)=(h_i(h))_{i\in I}.
	\end{equation}
	Each such finite algebra has a native table presentation in the fixed
	presentation context.  The comparison retains the actual
	local algebras and maps.  The presentation lists and generator markings carry their own
	specified comparisons.
	
	For $A\subseteq I$, let
	\[
	\mathcal M_A(H)
	=\mathbb C[(h_i)_{i\in A}]\subseteq\mathbb C^H
	\]
	be the algebra of functions depending only on the indicated partial
	state.  It is pulled back from the actual image $h_A(H)$, which can be a
	proper subset of the full partial join.  For a cut
	edge $a$, write $I=L_a\sqcup R_a$ for its two components and
	$h_a:H\to Y_a$ for the shared separator map.  The exact separator
	contract is
	\begin{equation}\label{eq:native-group-full-cut-contract}
		\mathcal M_{L_a}(H)\cap\mathcal M_{R_a}(H)
		=\mathbb C[h_a]
		\qquad(a\in E(J)).
	\end{equation}
	All algebras in this identity are represented on the same actual
	source $H$.
	
	\begin{lemma}[Cut information, kernel joins, and rectangularity]
		\label{lem:full-cut-kernel-rectangle}
		Let $\mathbf T$ be a finite algebra, with surjective homomorphisms
		$a:\mathbf T\to\mathbf A$, $b:\mathbf T\to\mathbf B$, and
		$s:\mathbf T\to\mathbf Y$, where
		$s=p\circ a=q\circ b$ for surjective homomorphisms
		$p:\mathbf A\to\mathbf Y$ and $q:\mathbf B\to\mathbf Y$.
		Assume that $(a,b)$ identifies $\mathbf T$ with a subalgebra
		of $\mathbf A\times_{\mathbf Y}\mathbf B$. Put
		\[
		\alpha=\ker a,\qquad \beta=\ker b,\qquad \sigma=\ker s.
		\]
		\begin{enumerate}[label=\textup{(\roman*)},leftmargin=*]
			\item \emph{Common information.}
			For the finite function algebras on $T$,
			\begin{equation}\label{eq:full-cut-kernel-meet}
				\mathbb C[a]\cap\mathbb C[b]=\mathbb C[s]
				\quad\Longleftrightarrow\quad
				\alpha\vee\beta=\sigma.
			\end{equation}
			\item \emph{Full fiber product.}
			If $\alpha$ and $\beta$ permute, these conditions are further
			equivalent to
			\begin{equation}\label{eq:full-cut-fiber-criterion}
				\mathbf T=\mathbf A\times_{\mathbf Y}\mathbf B
				\quad\text{under }(a,b).
			\end{equation}
			In particular, the last equivalence holds in every Mal'cev variety.
		\end{enumerate}
	\end{lemma}
	\begin{proof}
		\leavevmode
		\begin{enumerate}[label=\textup{(\roman*)},leftmargin=*]
			\item \emph{Common information.}
			A function on $T$ belongs to $\mathbb C[a]$ (respectively,
			$\mathbb C[b]$) precisely when it is constant on $\alpha$-classes
			(respectively, $\beta$-classes). The intersection consists of
			functions constant on classes of the equivalence relation generated
			by $\alpha\cup\beta$, namely the congruence $\alpha\vee\beta$.
			Since $s$ is onto, finite functions separate its fibers, giving
			\eqref{eq:full-cut-kernel-meet}.
			\item \emph{Full fiber product.}
			If the two congruences permute, their join equals
			$\alpha\circ\beta$. For $u\in A$, $v\in B$ with $p(u)=q(v)$,
			choose $t,t'\in T$ with $a(t)=u$ and $b(t')=v$.
			Then $(t,t')\in\sigma=\alpha\circ\beta$, so there is a
			$z\in T$ with $t\mathrel\alpha z\mathrel\beta t'$.
			Thus $a(z)=u$ and $b(z)=v$; every compatible pair occurs.
			Conversely, if every compatible pair occurs, any two elements
			of the same $s$-fiber can be connected by first retaining the
			$a$-coordinate and then the $b$-coordinate. Hence
			$\sigma\subseteq\alpha\circ\beta\subseteq\alpha\vee\beta$;
			the reverse inclusion is automatic. The Mal'cev conclusion follows
			from congruence permutability.\qedhere
		\end{enumerate}
	\end{proof}
	
	\paragraph{The common quotient and Goursat's lemma.}
	The preceding proof also recovers the two-factor finite-group form of
	Goursat's lemma \cite{Goursat1889,BauerSenZvengrowski2015}.
	For a subdirect subgroup $T\leq A\times B$, let $a,b$ be its
	projections and $N_a,N_b$ their normal kernels. Take
	$J=T/(N_aN_b)$, with quotient map $s$. Since each kernel lies in
	$\ker s$, $s$ induces surjections $A\to J$ and $B\to J$.
	Their pullback is exactly $T$ by \cref{lem:full-cut-kernel-rectangle},
	because the join of the two kernel congruences corresponds to $N_aN_b$.
	In particular the oriented quotient identifications are
	$A/a(N_b)\cong J\cong B/b(N_a)$.
	The same proof identifies $\mathbb C[a]\cap\mathbb C[b]=\mathbb C[s]$:
	$J$ is the finest deterministic common readout, also represented by the
	connected components of the bipartite support graph
	(\cref{thm:proper-cut-incidence}; \cite{GacsKorner1973,Nair2013}).
	With any full-support probability law on $T$ this is the common
	variable of G\'acs--K\"orner, whose entropy depends on the law.
	The information identity does not force rectangles for arbitrary
	relations; permuting congruences are the extra hypothesis that makes
	it sufficient in the native rigidity theorem below.
	
	\begin{theorem}[Full-cut rigidity, quotient classification, and sharp native minimum]
		\label{thm:native-source-rigidity}
		
		Fix the finite Mal'cev join-tree atlas and native source
		$\mathbf H$ with compatible homomorphic chart maps $h_i$ as above.
		Let $\psi:\mathbf H\to\mathbf W$ be the canonical comparison map
		\eqref{eq:native-source-comparison-map}, and put
		$\mathbf T=\psi(\mathbf H)\leq\mathbf W$. For a cut $a$, denote the
		images of $T$ on its complete sides by $\mathbf T_{L_a}$ and
		$\mathbf T_{R_a}$.
		\begin{enumerate}[label=\textup{(\alph*)},leftmargin=*]
			\item \emph{Equivalent full-cut conditions.}
			The following conditions are equivalent:
			\begin{enumerate}[label=(\roman*),nosep]
				\item Every $h_i$ is surjective, and the represented common-information
				algebras satisfy
				\eqref{eq:native-group-full-cut-contract} at \emph{every full tree cut}.
				\item The canonical comparison map
				$\psi:\mathbf H\twoheadrightarrow\mathbf W$ is surjective.
				\item Every $h_i$ is surjective, and each full-cut image is the
				entire fiber product of its \emph{actual} complete-side images:
				\begin{equation}\label{eq:native-source-all-cut-rectangles}
					\mathbf T\;=\;
					\mathbf T_{L_a}\times_{\mathbf Y_a}\mathbf T_{R_a}
					\qquad(a\in E(J)),
				\end{equation}
				under the two side projections.
			\end{enumerate}
			\item \emph{Observable quotient and minimum.}
			Consequently, the finite native sources satisfying the state and
			full-cut contracts are exactly, at the level of their prescribed
			chart maps, the surjective homomorphic extensions of $\mathbf W$.
			They have the canonical observable quotient
			\begin{equation}\label{eq:native-source-quotient-classification}
				\mathbf H/\theta\cong\mathbf W,
				\qquad\theta=\ker\psi=\bigcap_{i\in I}\ker h_i,
			\end{equation}
			and their sharp minimum cardinality is
			\begin{equation}\label{eq:sharp-native-source-cost}
				\boxed{\min|H|=|W|}.
			\end{equation}
			Equality holds precisely when $\theta=\Delta_H$, equivalently
			when the chart maps jointly separate $H$. In that case $\psi$ is
			the unique isomorphism to $\mathbf W$ over all the prescribed
			chart maps.
			\item \emph{Computing the minimum.}
			The minimum $|W|$ is computed by the separator messages in \eqref{eq:local-pair-count-message} and the root sum in \eqref{eq:local-pair-source-size}, applied to the underlying chart and separator sets. These formulas allow nonuniform fibers.
			
			For finite groups, writing the chart groups as $G_i$, the quotient
			classification is a short exact sequence
			\begin{equation}\label{eq:native-group-extension-classification}
				1\longrightarrow K\longrightarrow H
				\overset{\psi}{\longrightarrow}W\longrightarrow1,
				\qquad K=\ker\psi,
			\end{equation}
			and the sharp cardinalities are
			\begin{equation}\label{eq:sharp-native-group-source-cost}
				|H|=|K|\,|W|,
				\qquad\min|H|=|W|
				=\frac{\prod_{i\in I}|G_i|}{\prod_{a\in E(J)}|Y_a|}.
			\end{equation}
			For finite modules over a fixed finite ring, the analogous
			module quotient and product-over-separators formula hold.
			\item \emph{Recovery faithfulness.}
			If the atlas also has relative pair lifting,
			for example by the hypotheses of
			\cref{thm:congruence-interface-saturation,thm:complete-port-native-descent},
			then $c_W=c_0$ and every source in \textup{(i)}--\textup{(iii)}
			is recovery-faithful. The same minimum and equality
			classification then hold for the full faithful native contract.
		\end{enumerate}
	\end{theorem}
	\begin{proof}
		\leavevmode
		\begin{enumerate}[label=\textup{(\alph*)},leftmargin=*]
			\item \emph{Equivalent full-cut conditions.}
			\emph{\textup{(ii)} $\implies$ \textup{(i)} and \textup{(iii)}.} The join-tree construction
			\cref{thm:equational-native-amalgamation} makes $W$ project onto
			every chart and identifies it with the full fiber product of its
			complete sides at every cut. Surjective pullback along $\psi$
			preserves all the finite function-algebra intersections; therefore
			\textup{(ii)} implies both \textup{(i)} and \textup{(iii)}.
			
			\emph{\textup{(iii)} $\implies$ \textup{(i)}.} The full rectangle in each separator fiber
			implies that a function factors through each complete-side projection exactly when it factors through the named separator. Pullback from $T$ to $H$ gives the same
			intersection identity. Thus \textup{(iii)} implies \textup{(i)}.
			
			\emph{\textup{(i)} $\implies$ \textup{(ii)}.} The surjection $H\to T$ identifies the
			represented algebras of chart and separator observations, so the
			contracts also hold on $T$. At every cut the projections from
			$T$ onto its complete-side images are homomorphisms. Their kernels
			permute, since $T$ lies in the Mal'cev variety. By
			\cref{lem:full-cut-kernel-rectangle}, the exact common-information
			identity makes $T$ the full fiber product of those images.
			
			We show that $T=W$ by induction on the number of charts. If there
			is one chart, state surjectivity gives the result. Otherwise remove
			a leaf $\ell$, and let $T_0$ be the image of $T$ on all remaining
			charts. It projects onto every remaining chart. The cut contracts
			for the remaining tree descend to $T_0$: a function on $T_0$
			factoring through each side projection of a remaining cut pulls back to a function on $T$ factoring through each complete-side projection of the corresponding original cut. It must therefore be a function of that separator,
			and surjectivity of $T\to T_0$ makes this conclusion descend.
			By induction $T_0=W_0$, the full remaining matching algebra.
			At the removed leaf cut, the preceding two-factor conclusion gives
			\[
			T=A_\ell\times_{Y_{\ell,p(\ell)}}T_0
			=A_\ell\times_{Y_{\ell,p(\ell)}}W_0=W.
			\]
			This proves \textup{(i)} implies \textup{(ii)}.
			\item \emph{Observable quotient and minimum.}
			The homomorphism theorem gives
			$H/\ker\psi\cong W$ and the formula for $\theta$.
			Since $\psi$ is a surjection of finite sets, $|H|\geq|W|$;
			$W$ itself satisfies \textup{(ii)} and attains the bound. Equality
			holds if and only if $\psi$ is injective. The unique comparison map
			is then an isomorphism over the charts.
			\item \emph{Computing the minimum.}
			Apply \cref{prop:local-pair-source-cost}\textup{(i)} to the underlying finite sets.
			
			For finite groups the homomorphism theorem yields
			\eqref{eq:native-group-extension-classification}; each fiber of a
			surjective group homomorphism has size $|K|$. Equal-size fibers
			of the group separator maps give the product formula
			\eqref{eq:sharp-native-group-source-cost}. The same argument on
			underlying additive groups applies to finite modules.
			\item \emph{Recovery faithfulness.}
			Relative pair lifting gives $c_W=c_0$ by
			\cref{thm:equational-native-amalgamation}, and surjective pullback gives
			$c_H=c_W$.\qedhere
		\end{enumerate}
	\end{proof}
	
	\begin{theorem}[The faithful-native dichotomy and sharp source classification]
		\label{thm:faithful-matching-minimum}
		Fix the finite Mal'cev equational join-tree atlas, native
		presentation law, realized marked chart presentations,
		and homomorphic descriptor readouts of
		\cref{thm:equational-native-amalgamation}. Suppose the actual
		separator images agree. Let $\mathbf W$ be the full matching algebra,
		and $c_0$ the closure generated by the local chart theories.
		Consider the class $\mathfrak S_{\mathrm{cut}}$ of finite native
		algebras over the fixed chart maps that satisfy state coverage
		$\mathsf S$ and exact common information $\mathsf E$ at every
		\emph{complete} tree cut. Exactly one of the following regimes holds:
		\begin{enumerate}[label=(\roman*),nosep]
			\item \emph{Faithful regime.}
			Every ordered-pair root test for
			$S\subseteq V$, $d\notin c_0(S)$ succeeds, equivalently $c_W=c_0$.
			Then every source in $\mathfrak S_{\mathrm{cut}}$ is recovery-faithful;
			more precisely, $\mathfrak S_{\mathrm{cut}}$ consists exactly of the
			finite surjective homomorphic extensions
			$\psi:\mathbf H\twoheadrightarrow\mathbf W$ with
			$h_i=\pi_i\psi$. Each admits a native table presentation, so
			$\mathsf S,\mathsf F,\mathsf E,\mathsf N$ hold simultaneously.
			Its canonical quotient and sharp minimum are
			\[
			\mathbf H/\ker\psi\cong\mathbf W,
			\qquad \min_{\mathbf H\in\mathfrak S_{\mathrm{cut}}}|H|=|W|.
			\]
			The unique minimum, up to isomorphism over every prescribed
			chart map, is $\mathbf W$. The minimum has the tree recursion
			\cref{eq:local-pair-count-message,eq:local-pair-source-size} and, for finite groups and
			modules, the product formula \cref{eq:sharp-native-group-source-cost}.
			
			\item \emph{Obstructed regime.}
			Some ordered-pair root test fails, equivalently
			$c_W\ne c_0$. The rejected test specifies
			$S\subseteq V$ and $d\notin c_0(S)$ such that
			$d\in c_W(S)$. This additional determination holds in
			\emph{every} compatible source over the given charts, with or
			without native operations or state coverage. Hence there is
			no source satisfying the full faithful native contract.
			Nevertheless $\mathbf W$ remains a native source satisfying
			$\mathsf S,\mathsf E,\mathsf N$; failure concerns
			$\mathsf F$, not mere algebraic existence.
		\end{enumerate}
		
		Relative pair lifting guarantees the first regime. The pair test
		also decides the alternative directly for each given atlas.
	\end{theorem}
	\begin{proof}
		The finite pair test in
		\cref{thm:amalgamation-existence-criterion} gives the exhaustive
		alternative $c_W=c_0$ or a specified missing-pair obstruction.
		Independently of that alternative,
		\cref{thm:native-source-rigidity} identifies every element of
		$\mathfrak S_{\mathrm{cut}}$ with a surjective homomorphic
		extension of $W$ and conversely; its full-cut proof does not
		use recovery faithfulness.
		\begin{enumerate}[label=\textup{(\roman*)},leftmargin=*]
			\item \emph{Faithful regime.}
			A surjection $\psi:H\twoheadrightarrow W$
			preserves both positive and negative descriptor dependencies,
			so $c_H=c_W=c_0$. The finite native law supplies a presentation
			for $H$. The quotient and equality classification are exactly
			those of \cref{thm:native-source-rigidity}, which also provides
			the cardinality formulas. Relative pair lifting gives $c_W=c_0$ by \cref{thm:equational-native-amalgamation}.
			\item \emph{Obstructed regime.}
			Monotonicity under the chart tuple map
			$q:H\to W$ yields $c_W\subseteq c_H$ for any compatible $H$.
			Thus the extra dependency supplied by the rejected test persists
			on every such source. Finally, the unconditional matching
			construction \cref{thm:equational-native-amalgamation} still
			realizes $W$ natively with state coverage and exact cut algebras.\qedhere
		\end{enumerate}
	\end{proof}
	
	The construction uses acyclic joins and closure under products and
	subalgebras \cite{Beeri1983,Birkhoff1935}. Its paired witnesses
	certify nonimplications as in Armstrong relations \cite{Armstrong1974,BeeriDFS1984}.
	
	\begin{remark}[Meet sources and adjacent-chart tests]
		\label{rem:native-group-adjacent-insufficient}
		First, omit the Mal'cev hypothesis and use the fixed variety of
		meet-semilattices, with its binary meet operation and no required top
		constant.  Two one-descriptor charts have local algebra $\{0,1\}$
		and an empty separator.  Their full join has four states, but
		\[
		T=\{(0,0),(0,1),(1,0)\}
		\]
		is a proper subalgebra under coordinatewise meet.  Both projections
		are onto, its descriptor closure is free, and its two function
		algebras intersect only in the constants: the three states give a
		connected bipartite incidence path.  Thus all the state, recovery,
		and separator requirements hold on a smaller native source.
		Its two projection kernels do not permute.
		
		Second, retaining groups but checking only adjacent-chart algebras
		is insufficient.  Take three disjoint one-descriptor charts with
		local group $G_i=\mathbb F_2^2$ and identity descriptor
		$d_i=(u_i,z_i)$, arranged along a path with empty separators.  The
		proper subgroup
		\[
		T=\{((u_1,z_1),(u_2,z_2),(u_3,z_3)):
		u_1+u_2+u_3=0\}
		\]
		has $32$ states instead of $64$.  It projects onto every pair of
		local groups, so every adjacent-chart function-algebra intersection
		is constant.  Its descriptor closure remains free: whenever $d_i$
		is absent from a premise, varying $z_i$ refutes its determination.
		Nevertheless, across the cut separating chart $1$ from charts
		$2,3$, the bit $u_1=u_2+u_3$ is common information beyond the empty
		separator.
	\end{remark}
	
	\subsection{Native source rigidity from recovery faithfulness}
	\label{subsec:recovery-source-rigidity}
	
	The preceding minimum uses exact information across full tree cuts.
	For simple descriptor images in a Mal'cev variety, a proper
	state-exact subalgebra of the matching algebra adds a dependency.
	Thus recovery faithfulness and state coverage suffice for rigidity
	on any finite matching diagram, including cyclic diagrams.
	For scalar linear observations, recovery faithfulness also forces
	state coverage under the hypotheses below.
	
	\subsubsection{Recovery rigidity on arbitrary finite diagrams}
	\label{subsec:simple-descriptor-rigidity}
	
	Keep the fixed finite-signature Mal'cev variety $\mathcal V$ and its
	presentation context. A nontrivial algebra $\mathbf S$ is
	\emph{simple} if its only congruences are $\Delta_S,\nabla_S$.
	Consider any finite family of finite charts $\mathbf A_i\in\mathcal V$
	with descriptor sets $D_i$, whose union is $D$. At each chart require
	jointly separating onto homomorphisms
	$f_{i,d}:\mathbf A_i\twoheadrightarrow\mathbf S_d$, $d\in D_i$,
	where the finite algebras $\mathbf S_d$ are simple. Shared descriptors
	have the same actual target algebra. One-element descriptor images
	can be omitted and reinstated as constant observations.
	No join tree is assumed. Construct the full matching set
	\[
	W=\{(a_i)\in\prod_i A_i:
	f_{i,d}(a_i)=f_{j,d}(a_j)\text{ whenever }d\in D_i\cap D_j\}.
	\]
	If nonempty, it is a native subalgebra. Write $f_d:W\to S_d$ for
	its induced readouts and $c_0$ for the least closure generated by
	the local theories. A candidate source is a finite nonempty native
	algebra $\mathbf H$ with compatible homomorphisms $h_i:H\to A_i$;
	its comparison is $\psi=(h_i):H\to W$ and its observations are
	$f_d\psi$. State exactness means that every $h_i$ is onto.
	
	\begin{lemma}[Coordinate bases for subdirect products of simple algebras]
		\label{lem:simple-coordinate-basis}
		Let a finite algebra $\mathbf Z\in\mathcal V$ have jointly separating
		onto homomorphisms $g_d:Z\twoheadrightarrow S_d$ to finite simple
		algebras.
		\begin{enumerate}[label=\textup{(\roman*)},leftmargin=*]
			\item \emph{Constructing a coordinate basis.}
			An inclusion-minimal set $B$ with $c_Z(B)=D$ is obtained by scanning descriptors and retaining a
			descriptor exactly when the previously retained ones do not recover it.
			\item \emph{Isomorphism for every minimal basis.}
			If $B\subseteq D$ is an inclusion-minimal set with
			$c_Z(B)=D$, then the actual coordinate map is an isomorphism
			\begin{equation}\label{eq:simple-coordinate-basis}
				g_B:\mathbf Z\xrightarrow{\ \cong\ }
				\prod_{d\in B}\mathbf S_d.
			\end{equation}
		\end{enumerate}
	\end{lemma}
	\begin{proof}
		\leavevmode
		\begin{enumerate}[label=\textup{(\roman*)},leftmargin=*]
			\item \emph{Constructing a coordinate basis.}
			The two-factor step is Fleischer's lemma, and its iteration is the
			classical Foster--Pixley decomposition
			\cite[Lemma~IV.10.1 and Corollary~IV.10.2]{BurrisSankappanavar1981}.
			We record why the isomorphism uses the specified coordinates.
			Let $\mathbf R\leq\mathbf P\times\mathbf S$ be subdirect with
			$\mathbf S$ simple. On $R$, the coordinate kernels $\alpha,\beta$
			satisfy $\alpha\cap\beta=\Delta_R$, and $\beta$ is maximal proper.
			If $\alpha\leq\beta$, then $\alpha=\Delta_R$ and $R$ is the graph
			of a homomorphism $P\to S$. Otherwise
			$\alpha\vee\beta=\nabla_R$; permutability gives
			$\alpha\circ\beta=\nabla_R$, so $R=P\times S$.
			
			Apply this dichotomy to the image of the already retained coordinates
			and the next $g_d$. The graph case is exactly recovery of $d$ from
			those coordinates. In the other case the new coordinate supplies a
			full independent factor $S_d$. The retained coordinate map is therefore
			onto the indicated product at every stage. At termination all readouts
			factor through it, and joint separation makes it injective. Each
			retained coordinate is indispensable because its nontrivial factor
			can vary with all other retained coordinates fixed.
			\item \emph{Isomorphism for every minimal basis.}
			Scan
			any given minimal generating set $B$ first. If a label of $B$ were
			skipped, it would follow from earlier labels and $B$ would generate
			without it. After all of $B$ is retained, every later label is
			determined. This proves
			\eqref{eq:simple-coordinate-basis} for every such $B$.\qedhere
		\end{enumerate}
	\end{proof}
	
	\begin{theorem}[Faithful source rigidity for simple descriptor algebras]
		\label{thm:simple-diagram-source-rigidity}
		Fix the atlas above.
		\begin{enumerate}[label=\textup{(\roman*)},leftmargin=*]
			\item \emph{Existence.}
			A state-exact faithful native source exists
			if and only if $W$ is nonempty, all projections $W\to A_i$ are onto,
			and $c_W=c_0$.
			\item \emph{Certificate of a proper image.}
			Assume the existence conditions and let $B$ be any inclusion-minimal set with $c_0(B)=D$. If $Z=\psi(H)\ne W$ for a state-exact candidate,
			then there are $A\subsetneq B$ and $e\in B\setminus A$
			with
			\begin{equation}\label{eq:simple-source-defect-certificate}
				e\in c_H(A)\setminus c_0(A),\qquad
				|Z|=\prod_{d\in A}|S_d|,
				\qquad |A|\leq\lfloor\log_2|Z|\rfloor.
			\end{equation}
			They can be found by the finite coordinate scan on $Z$.
			\item \emph{Classification and minimum.}
			Under the existence conditions in (i), every state-exact candidate
			source satisfies
			\begin{equation}\label{eq:simple-diagram-faithfulness}
				c_H=c_0\quad\Longleftrightarrow\quad
				\psi:\mathbf H\to\mathbf W\text{ is onto}.
			\end{equation}
			Consequently all faithful candidates have the same observable quotient
			and the sharp source cost
			\begin{equation}\label{eq:simple-diagram-minimum}
				\mathbf H/\ker\psi\cong\mathbf W,
				\qquad
				\min|H|=|W|=\prod_{d\in B}|S_d|,
			\end{equation}
			where $B$ is any inclusion-minimal set with $c_0(B)=D$.
			Every minimum is uniquely isomorphic to $\mathbf W$ over the actual
			charts, by its comparison map. Thus the closure and the specified
			descriptor-image cardinalities determine the minimum, while the actual
			chart maps determine the equality isomorphism.
		\end{enumerate}
	\end{theorem}
	\begin{proof}
		\leavevmode
		\begin{enumerate}[label=\textup{(\roman*)},leftmargin=*]
			\item \emph{Existence.}
			If a faithful state-exact $H$ exists, its nonempty image lies in $W$
			and covers every chart. Hence $W$ covers every chart as well. Every
			local implication holds on $W$, and restriction to the image of $H$
			can only add implications. Therefore
			$c_0\subseteq c_W\subseteq c_H=c_0$. This proves necessity of the
			existence criterion. Conversely, a nonempty $W$ satisfying that criterion
			is itself an admissible native source; its finite table presentation is
			provided by \cref{lem:finite-table-presentation}.
			\item \emph{Certificate of a proper image.}
			Now assume the criterion. The global readouts on $W$ jointly separate
			states and are onto their simple targets. By
			\cref{lem:simple-coordinate-basis}, $f_B$ identifies $W$ with
			$\prod_{d\in B}S_d$. If a candidate has proper image $Z$, state
			exactness makes each restricted coordinate $f_d|_Z$, $d\in B$,
			onto. They still jointly separate $Z$. Applying the lemma to $Z$
			gives $A\subseteq B$ for which $Z\cong\prod_{d\in A}S_d$ via
			those actual coordinates. Since $Z$ is proper and finite, $A\ne B$.
			For any $e\in B\setminus A$, the $A$-coordinates determine all of
			$Z$, hence determine $e$ on $H$. In the full product $W$, they do
			not determine the nontrivial independent factor $S_e$. This proves
			the extra implication and rules out faithfulness. Each $|S_d|\geq2$
			gives the bound in \eqref{eq:simple-source-defect-certificate}.
			\item \emph{Classification and minimum.}
			By (ii), faithfulness forces surjectivity. If $\psi$ is onto, pullback preserves the whole determination theory,
			so $c_H=c_W=c_0$. The homomorphism theorem gives the observable
			quotient. Counting the fibers of $\psi$ proves $|H|\geq|W|$;
			equality holds exactly when $\psi$ is bijective, hence an isomorphism.
			Its chart components specify it uniquely. Finally the coordinate
			isomorphism for $B$ gives the product formula. This counts $W$ directly
			using the prescribed descriptor images.\qedhere
		\end{enumerate}
	\end{proof}
	
	For finite groups and
	modules the surjection gives $|H|=|\ker\psi|\,|W|$.
	
	\begin{proposition}[Sharp descriptor class for recovery-based rigidity]
		\label{prop:simple-descriptor-boundary}
		Fix a Mal'cev variety and a family $\mathcal D$ of allowed finite
		nontrivial descriptor algebras. The conclusion
		\eqref{eq:simple-diagram-faithfulness} holds for every atlas with
		onto, jointly separating $\mathcal D$-valued readouts and every
		state-exact native candidate, whenever $W$ is state-exact and faithful,
		if and only if every algebra in $\mathcal D$ is simple.
	\end{proposition}
	\begin{proof}
		Sufficiency is \cref{thm:simple-diagram-source-rigidity}. If an allowed
		$\mathbf Q$ has $\Delta_Q<\theta<\nabla_Q$, use two disjoint
		one-descriptor charts, both with algebra $\mathbf Q$ and identity
		readout. Their matching algebra is $Q^2$ with free two-descriptor
		closure. The proper native source
		\[
		Z=\{(x,y)\in Q^2:x/\theta=y/\theta\}
		=\mathbf Q\times_{\mathbf Q/\theta}\mathbf Q
		\]
		has both projections onto. A nontrivial $\theta$-class supplies pairs that vary each
		descriptor while fixing the other. Thus $c_Z=c_{Q^2}$, although
		\[
		|Z|=\sum_{C\in Q/\theta}|C|^2<|Q|^2.
		\]
		This violates the asserted universal conclusion in the same native class.
	\end{proof}
	
	\paragraph{Counterexamples for state coverage and permutability.}
	State exactness cannot generally be dropped: the single group chart
	$A_5$ with its identity descriptor has the same one-descriptor closure
	on the proper subgroup $\langle(12)(34)\rangle$, of order two.
	The three-element implication algebra in \cref{ex:goursat-faithful-nonrigid}
	gives a proper, state-exact, faithful source even with simple descriptor images.
	
	\subsubsection{Rigidity for nonsimple charts with named quotients}
	\label{subsec:minimal-congruence-layer-rigidity}
	
	A proper source can preserve whole-chart recovery while changing an
	unobserved congruence layer. We first treat charts with one independently
	variable minimal layer. Its quotient detects changes in the coarse
	coordinates, and the whole-chart observations detect further loss within
	a coarse fiber.
	
	\begin{theorem}[Native rigidity with exposed minimal congruence layers]
		\label{thm:minimal-congruence-layer-rigidity}
		Fix a finite Mal'cev atlas in the presentation context with
		nonempty matching algebra $W\leq\prod_i A_i$ projecting onto each chart.
		Assume the observation language includes the whole chart readouts
		$f_i:W\twoheadrightarrow A_i$ and their coarsenings
		$g_i=q_i f_i$, where $q_i:A_i\twoheadrightarrow A_i/\alpha_i$ and
		$\alpha_i$ is a minimal nonzero congruence of $A_i$ with nontrivial
		simple quotient. Suppose $W$ is saturated under coordinatewise
		$\prod_i\alpha_i$ in the ambient product, and $c_W=c_0$ for the declared
		observation language.
		\begin{enumerate}[label=\textup{(\roman*)},leftmargin=*]
			\item \emph{Surjectivity criterion.}
			Every state-exact native source $H$ over
			this atlas with comparison $\psi:H\to W$ satisfies
			\[
			c_H=c_0\quad\Longleftrightarrow\quad \psi\text{ is onto}.
			\]
			\item \emph{Quotient and minimum.}
			In particular, the observable quotient is $H/\ker\psi\cong W$ for every
			faithful state-exact candidate. The sharp minimum is $|W|$ and the
			comparison of every minimum is its unique isomorphism to $W$ over the
			specified charts.
		\end{enumerate}
	\end{theorem}
	\begin{proof}
		\leavevmode
		\begin{enumerate}[label=\textup{(\roman*)},leftmargin=*]
			\item \emph{Surjectivity criterion.}
			Surjectivity implies faithfulness by pullback of readouts. Conversely, assume faithfulness and
			put $Z=\psi(H)$. Then each $f_i|_Z$ is onto $A_i$ and $c_Z=c_W$.
			For a fixed chart $i$, the premise consisting of $g_i$ and all $f_j$
			with $j\ne i$ fails to recover $f_i$ on $W$. To see this, take two
			different $\alpha_i$-equivalent values, extend one of them to a row by
			chart surjectivity, and substitute the other using saturation.
			Faithfulness supplies a pair in $Z$ with the same property.
			
			On $Z$, let $\eta_i$ be the intersection of the kernels of that premise.
			Its direct image under $f_i|_Z$ is a congruence on $A_i$. Indeed, for
			$\kappa_i=\ker(f_i|_Z)$, congruence permutability identifies that image
			with $(\eta_i\vee\kappa_i)/\kappa_i$ under $Z/\kappa_i\cong A_i$.
			The image is contained in $\alpha_i$ and is nonzero by the witness, so
			minimality makes it $\alpha_i$.
			If a row $z\in Z$ has $i$-entry $u$ and $u\equiv v\pmod{\alpha_i}$,
			choose an $\eta_i$-pair $z',z''$ whose $i$-entries are $u,v$.
			Applying the Mal'cev term $m(z,z',z'')$ leaves every other coordinate
			of $z$ unchanged and replaces $u$ by $v$. Consequently $Z$ too is
			saturated under $\prod_i\alpha_i$.
			
			Let $\overline W$ and $\overline Z$ be the images in
			$\prod_i A_i/\alpha_i$. Their coarse coordinates are onto simple
			algebras and jointly separate their rows. Equality of the recovery
			closures on $W,Z$ restricts to equality for these coarse coordinates.
			Choose an inclusion-minimal generating subset $B$ of the coarse
			coordinates. By \cref{lem:simple-coordinate-basis}, projection to $B$
			is an isomorphism to the same specified product of simple algebras
			on both $\overline W$ and $\overline Z$. Thus $\overline Z=\overline W$.
			Saturation now gives $Z=W$.
			\item \emph{Quotient and minimum.}
			The homomorphism theorem gives the quotient,
			and $|H|\geq|W|$ is attained by $W$ itself in the fixed finite-table
			presentation context. At equality, the onto comparison is bijective
			and its prescribed chart components determine it uniquely.\qedhere
		\end{enumerate}
	\end{proof}
	
	For finite groups one can take $\alpha_i$ from a minimal nontrivial
	normal subgroup $N_i\trianglelefteq A_i$. Saturation then amounts to
	$\prod_i N_i\leq W$, and the minimum has the explicit form
	\[
	|W|=\left(\prod_i|N_i|\right)\prod_{i\in B}|A_i/N_i|,
	\]
	where $B$ is a minimal generating set of coarse coordinates.
	For instance, take $A_i=S_5$ and $N_i=A_5$. Partition the charts into
	groups, with one shared sign descriptor per group and a separate
	whole-chart descriptor per chart. Then $W$ is the product of the
	corresponding sign fiber products. The local implications generate the recovery on $W$. The sign in
	each group varies freely if the premise contains neither that sign
	nor any whole-chart descriptor from that group. The independent $A_5$ factors
	prevent any further whole-chart determinations.
	
	\subsubsection{Quotient-complete observations and coupled-layer rigidity}
	\label{subsec:quotient-complete-source-rigidity}
	
	When the named observations separate the chart congruences, the hidden
	layers may be coupled. For a finite algebra $A_i$, let
	$q_{i,s}:A_i\twoheadrightarrow D_{i,s}$, $s\in S_i$, be named
	homomorphic readouts. Call this finite family \emph{quotient-complete} if
	\begin{equation}\label{eq:chart-quotient-completeness}
		\theta=\bigcap_{\{s\in S_i:\theta\leq\ker q_{i,s}\}}
		\ker q_{i,s}
		\qquad(\theta\in\Con(A_i)),
	\end{equation}
	where an empty intersection is the universal congruence. In particular,
	these readouts jointly separate chart states. The condition concerns
	the charts and allows arbitrary interface signatures.
	
	\begin{theorem}[Recovery rigidity from quotient-complete chart observations]
		\label{thm:quotient-complete-source-rigidity}
		Let $W\leq\prod_{i\in I}A_i$ be a nonempty finite matching algebra in a
		Mal'cev variety, in the fixed presentation context and
		projecting onto every chart. Suppose the observation
		language contains a quotient-complete family at each chart, pulled back
		along its chart projection. Other declared observations may also be
		retained. For every finite native source $H$ with surjective compatible
		chart maps and comparison $\psi:H\to W$, equip $H$ with the pullbacks
		of all these observations.
		\begin{enumerate}[label=\textup{(\roman*)},leftmargin=*]
			\item \emph{Certificate of a proper image.}
			If $Z=\psi(H)<W$, there are a chart subset
			$J\subseteq I$ with $|J|\geq2$, a chart $i\in J$, and a named readout
			$q_{i,s}$ such that all the named readouts of the charts in
			$J\setminus\{i\}$ determine $q_{i,s}$ on $Z$ but do not determine it
			on $W$. Thus every proper state-exact subsource has a witnessed extra
			functional dependency in the declared language.
			\item \emph{Faithfulness, quotient, and minimum.}
			The comparison satisfies
			\[
			c_H=c_W\quad\Longleftrightarrow\quad\psi\text{ is onto}.
			\]
			No restriction on the matching graph or saturation under independent
			coordinate congruences is required. In particular, if $c_W=c_0$, this is
			an exact criterion for recovery-faithful native sources over the atlas.
			Their observable quotient is $H/\ker\psi\cong W$, their sharp minimum
			size is $|W|$, and every minimum has its comparison map as its unique
			isomorphism to $W$ over the prescribed charts.
		\end{enumerate}
	\end{theorem}
	
	\begin{proof}
		\leavevmode
		\begin{enumerate}[label=\textup{(\roman*)},leftmargin=*]
			\item \emph{Certificate of a proper image.}
			We first record the two-coordinate fact used in the proof. Let
			$T\leq B\times A$ project onto $B$ and $A$, and define
			\[
			\beta_T=\{(u,v)\in A^2:\text{some }b\in B\text{ has }
			(b,u),(b,v)\in T\}.
			\]
			This is a congruence of $A$. Reflexivity and symmetry are immediate,
			and compatibility follows by applying the operations to pairs with
			common first coordinate. For transitivity, if $(b,u),(b,v)$ and
			$(d,v),(d,w)$ belong to $T$, then the Mal'cev term gives
			$m((b,u),(b,v),(d,v))=(d,u)$, so $(u,w)\in\beta_T$.
			Moreover, every nonempty fiber $T_b$ is an entire $\beta_T$-class:
			if $(b,u)\in T$ and $(d,u),(d,v)\in T$, then
			$m((b,u),(d,u),(d,v))=(b,v)\in T$.
			
			Now suppose $Z<W$ and choose an inclusion-minimal $J\subseteq I$
			with $Z_J\ne W_J$, where subscripts denote coordinate images. Chart
			surjectivity gives $|J|\geq2$. Fix $i\in J$ and put
			$B=W_{J\setminus\{i\}}=Z_{J\setminus\{i\}}$. Both
			$T=Z_J$ and $V=W_J$ are subdirect products of $B$ and $A_i$.
			Their congruences above satisfy $\beta_T\leq\beta_V$. The inequality
			is strict: if the congruences were equal, each nonempty $T_b\subseteq V_b$
			would be the same congruence class, giving $T=V$.
			
			By \eqref{eq:chart-quotient-completeness}, there is a named $q_{i,s}$
			with $\beta_T\leq\ker q_{i,s}$ and
			$\beta_V\not\leq\ker q_{i,s}$. On $Z$, fixing the entire
			$J\setminus\{i\}$ state therefore fixes $q_{i,s}$. On $W$, a pair
			in the definition of $\beta_V$ refutes that determination. The named
			readouts on each other chart jointly separate its states, so replacing
			those entire chart states by their named readouts leaves both assertions
			unchanged. Surjective pullback from $Z$ to $H$ gives the same unintended
			dependency on $H$.
			\item \emph{Faithfulness, quotient, and minimum.}
			Consequently $c_H=c_W$ excludes a proper image.
			
			Conversely, a surjective comparison preserves every implication and
			every nonimplication. The homomorphism theorem gives the quotient, and
			a finite surjection gives $|H|\geq|W|$. The matching algebra $W$ itself
			attains the bound in the fixed finite-table realization context. At
			equality the comparison is an isomorphism; its prescribed chart
			components determine it uniquely.\qedhere
		\end{enumerate}
	\end{proof}
	
	\begin{corollary}[Recovery rigidity from meet-irreducible chart quotients]
		\label{cor:meet-irreducible-chart-rigidity}
		In \cref{thm:quotient-complete-source-rigidity}, it suffices to name the
		quotient $A_i/\mu$ for every proper meet-irreducible congruence
		$\mu\in\Con(A_i)$. For finite groups these are normal-subgroup
		quotients. Thus arbitrary finite group charts admit a finite observation
		enrichment giving recovery rigidity on every state-exact matching
		algebra, including those with coupled normal layers.
	\end{corollary}
	\begin{proof}
		In a finite lattice every element is the meet of the meet-irreducible
		elements above it. A proper meet-reducible element is the meet of two strictly larger
		elements. Repeating this decomposition terminates in a finite lattice. Apply this to $\Con(A_i)$ to obtain
		\eqref{eq:chart-quotient-completeness}, then apply the theorem.
	\end{proof}
	
	\begin{example}[A cyclic coupled layer detected by a coarse quotient]
		\label{ex:coupled-layer-quotient-rigidity}
		Take three group charts on the coordinate pairs $(x,y)$, $(y,z)$ and
		$(z,x)$, with $x,y,z\in\mathbb Z/4\mathbb Z$, and match the shared
		coordinates on the triangle. Then $W\cong(\mathbb Z/4\mathbb Z)^3$.
		No nontrivial change supported in just one chart remains in $W$, because
		its two neighboring charts fix both of its coordinates. Hence $W$ is
		not saturated under a product of nonzero chart congruences.
		The proper subgroup
		\[
		Z=\{(x,y,z):x+y+z\equiv0\pmod2\}
		\]
		projects onto all three pair charts. The whole-chart observations give $Z$ and $W$ the same closure.
		Any two charts determine the third. Given only one chart, the
		missing coordinate still has two possible values on $Z$.
		Name the quotient observation $x\bmod2$ at a chart containing $x$.
		On $Z$ it is a function of $(y,z)$. On $W$ it varies for fixed $(y,z)$. A
		quotient-complete family at every chart rules out every proper
		state-exact subgroup by \cref{thm:quotient-complete-source-rigidity}. The resulting sharp source minimum is
		$|W|=64$, for the recovery contract of the enriched full matching source.
	\end{example}
	
	Quotient-complete observations replace the independent-layer saturation
	and simple-coarse-quotient assumptions of
	\cref{thm:minimal-congruence-layer-rigidity}. For explicitly given $W$
	and $Z$, compute an extra-dependency certificate by deleting charts
	while the projections remain unequal, then scanning the named quotient
	readouts. This uses Mal'cev rectangularity and the finite congruence lattice.
	
	\subsubsection{Scalar-linear specialization}
	
	For scalar linear observations, a descriptor basis detects every loss of
	global rank. Recovery faithfulness therefore excludes every proper
	image in the matching source, even when the descriptor list omits
	projective directions.
	
	\begin{lemma}[Bounds for signed requirements on scalar linear observations]
		\label{lem:linear-access}
		Let $V$ be a finite descriptor-label set with scalar linear forms $v_a$ on a finite $\mathbb F_q$-space.
		Put	$R=\dim\operatorname{span}\{v_a:a\in V\}$,
		$c(T)=\{a\in V:v_a\in\operatorname{span}\{v_t:t\in T\}\}$,
		and $r(T)=\dim\operatorname{span}\{v_t:t\in T\}$.
		Fix an exact interface trace $D\in\operatorname{Fix}(c|_Q)$ on $Q\subseteq V$.
		A signed requirement $P^+\cup N^-$ is feasible when there exists a closed set $C$
		with $C\cap Q=D$, $P\subseteq C$ and $C\cap N=\varnothing$.
		\begin{enumerate}[label=\textup{(\roman*)},leftmargin=*]
			\item \emph{Reducing the positive requirements.}
			For every $P$ there is $G\subseteq P$ such that
			\[
			|G|=r(D\cup P)-r(D)\le R-r(D),\qquad
			c(D\cup G)=c(D\cup P).
			\]
			The least completion and all signed-feasibility decisions are unchanged on replacing $P$ with $G$.
			\item \emph{Minimal infeasible families.}
			Every inclusion-minimal infeasible family in this fixed trace has at most $R-r(D)+1$
			literals, and a positive-only one has at most $R-r(D)$ literals.
			\item \emph{Surjective pullback.}
			The conclusions remain true after pullback along a surjective linear source with invisible kernel.
		\end{enumerate}
	\end{lemma}
	
	\begin{proof}
		\leavevmode
		\begin{enumerate}[label=\textup{(\roman*)},leftmargin=*]
			\item \emph{Reducing the positive requirements.}
			Extend a basis of $\operatorname{span}(v_D)$ by forms drawn from $P$
			to a basis of $\operatorname{span}(v_{D\cup P})$; the added labels form
			$G$. Equality of the spans implies equality of the generated closure,
			so it preserves the canonical least completion and all tests for
			presence and absence of labels. Such a completion has the prescribed
			trace and avoids $N$ exactly when
			$c(D\cup P)\cap((Q\setminus D)\cup N)=\varnothing$.
			\item \emph{Minimal infeasible families.}
			If an infeasible family is inclusion-minimal, choose a forbidden label
			in this intersection. If it lies in $Q\setminus D$, the positive
			requirements alone are already infeasible, so no negative literal
			occurs. Otherwise a single negative literal witnesses failure.
			Deleting any positive premise must remove the forcing; hence the
			positives are linearly irredundant modulo the baseline span and number
			at most $R-r(D)$. There is at most one additional negative literal.
			\item \emph{Surjective pullback.}
			A surjective linear pullback induces an injection of dual spaces,
			so its pulled-back descriptor spans have identical ranks and
			closure; its kernel only repeats states and does not change the proof.\qedhere
		\end{enumerate}
	\end{proof}
	
	\begin{corollary}[Complete rigidity for faithful scalar-linear sources]
		\label{cor:native-linear-source-rigidity}
		Consider the fixed actual linear atlas of
		\cref{thm:linear-native-amalgamation}, satisfying separator visibility,
		with constructed source $X$ and
		\[
		R=\dim X=\sum_i r_i-\sum_a s_a.
		\]
		Let $H$ be any finite-dimensional $\mathbb F_q$-vector space with
		compatible linear maps $h_i:H\to X_i$, and let
		$\psi:H\to X$ be their joint map.  Give $H$ the pulled-back scalar
		descriptor observations.
		\begin{enumerate}[label=\textup{(\roman*)},leftmargin=*]
			\item \emph{Faithfulness and minimum size.}
			The comparison satisfies
			\begin{equation}\label{eq:linear-faithfulness-is-surjectivity}
				c_H=c_0
				\quad\Longleftrightarrow\quad
				\psi\text{ is onto}.
			\end{equation}
			Thus every faithful native linear source, including sources whose
			observations need not separate all their states, has the form
			\begin{equation}\label{eq:linear-source-complete-classification}
				H\cong X\oplus K\quad\text{over the atlas},
				\qquad K=\ker\psi,\quad k=\dim K,
				\qquad |H|=q^{R+k}.
			\end{equation}
			The splitting is not canonical.  The minimum $q^R$ is attained
			exactly when $K=0$, and then the comparison with $X$ is the unique
			isomorphism over all actual local charts.  Every faithful linear
			subamalgam of $X$ is therefore $X$ itself, and every faithful linear
			source with jointly separating descriptor observations is this minimum
			source up to the specified isomorphism.
			Without compatible maps to this fixed atlas, every full linear
			representation of the same closure still has at least $q^R$ states
			by its matroid rank.  The isomorphism classification above additionally
			uses the specified actual chart maps.
			\item \emph{Certificate of a proper image.}
			If $\psi$ is not onto, Gaussian elimination produces an unintended
			implication $A\to e$ with $|A|\leq\dim\psi(H)$, using only labels
			from one descriptor basis of $X^*$.  This gives an explicit recovery
			certificate for failure of a proposed smaller native source.
			\item \emph{Readout algebras and fibers.}
			Write $r(T)=\dim\operatorname{span}(v_T)$, as before.  For any faithful
			source in \eqref{eq:linear-source-complete-classification}, its
			diagonal algebra $M_H(T)$ of functions of the $T$-readout satisfies
			\begin{equation}\label{eq:linear-source-visible-hidden-cost}
				\dim_{\mathbb C}M_H(T)=q^{r(T)},\qquad
				|f_T^{-1}(y)|=q^{R+k-r(T)}
				\quad\bigl(y\in f_T(H)\bigr).
			\end{equation}
			Accordingly the exact summaries and signed access bounds of
			\cref{lem:linear-access} use the same intrinsic budget
			$R-r(D)$ for every faithful native realization of this atlas.  An additional invisible kernel increases the number of source states
			per observed value.
		\end{enumerate}
	\end{corollary}
	
	\begin{proof}
		\leavevmode
		\begin{enumerate}[label=\textup{(\roman*)},leftmargin=*]
			\item \emph{Faithfulness and minimum size.}
			Separator visibility gives $c_X=c_0$, so
			\cref{cor:scalar-contract-equivalence} proves
			\eqref{eq:linear-faithfulness-is-surjectivity}.
			A surjective linear map has a linear section, giving the displayed
			splitting and size formula.  All readouts factor through $\psi$,
			so the kernel factor is invisible to them.  The equality and
			uniqueness assertions follow as in
			\cref{thm:native-source-rigidity}.  Any scalar linear representation of $c_0$ has matroid rank $R$,
			so its source dimension is at least $R$. The classification over
			the original charts also requires the comparison map.
			\item \emph{Certificate of a proper image.}
			For the explicit certificate,
			take a descriptor basis $B$ of $X^*$, choose a subset $A$ whose
			restrictions form a basis of $\psi(H)^*$, and take $e\in B\setminus A$.
			Then $|A|=\dim\psi(H)$ and $e\in c_H(A)\setminus c_0(A)$.
			\item \emph{Readout algebras and fibers.}
			The $T$-readout
			has rank $r(T)$ and every realized fiber in $H$ has dimension
			$R+k-r(T)$.  Counting its values and fibers proves
			\eqref{eq:linear-source-visible-hidden-cost}.  The generated source
			algebras have the same inclusions and the same strict rank budget
			as those on $X$, so the stated access bounds apply unchanged.\qedhere
		\end{enumerate}
	\end{proof}
	
	The vector-valued readouts in \cref{rem:native-group-adjacent-insufficient}
	admit a faithful proper subspace. Dropping linear closure can also reduce
	the minimum, as the spectral subrelations in \cref{ex:linear-source-cost-scope} show.
	
	\subsection{Exact source minima beyond Mal'cev rigidity}
	\label{subsec:nonmalcev-source-minima}
	
	A proper subalgebra can preserve both chart coverage and common information.
	\Cref{thm:forest-native-minimum} computes the smallest such image by
	closing a spanning tree in each separator fiber under the native operations.
	We then characterize rigidity on the finite members of a variety and
	compute the cost for transitive group actions from their diagonal orbits.
	
	\subsubsection{Minimum size for one relational cut}
	\label{subsec:proper-cut-connectivity}
	
	For finite observations, the finest common readout is the component of
	the incidence graph \cite{GacsKorner1973,Nair2013}. Connectivity in each
	separator fiber therefore determines the least size of a relation with
	prescribed chart images and exact separator information.
	
	\begin{theorem}[Common information as bipartite incidence connectivity]
		\label{thm:proper-cut-incidence}
		Let $R\subseteq U\times V$ be a nonempty finite relation, and
		write $U_R=\pi_U(R)$, $V_R=\pi_V(R)$. Suppose a specified
		readout $h:R\to Y$ factors through each projection separately.
		Form the bipartite incidence graph $\Gamma_R$ with vertex parts
		$U_R,V_R$ and one edge $(u,v)$ for each $(u,v)\in R$.
		Let $M_U(R),M_V(R)$ be the multiplication algebras of all
		functions of the respective projections, represented in
		$\mathbb C^R$.
		\begin{enumerate}[label=\textup{(\roman*)},leftmargin=*]
			\item \emph{Common functions.}
			The common-function algebra satisfies
			\begin{equation}\label{eq:proper-cut-components}
				M_U(R)\cap M_V(R)
				\cong \mathbb C^{\pi_0(\Gamma_R)},
			\end{equation}
			where the isomorphism is the actual pullback of the component
			index of the edge containing a state.
			\item \emph{Exactness of the specified readout.}
			For the specified readout $h$,
			\begin{equation}\label{eq:proper-cut-connectivity}
				M_U(R)\cap M_V(R)=\mathbb C[h]
				\quad\Longleftrightarrow\quad
				\Gamma_R[h^{-1}(y)]\text{ is connected for every }
				y\in h(R).
			\end{equation}
			Here $\Gamma_R[h^{-1}(y)]$ consists of the edges of label $y$
			and their incident vertices. When equality fails, the indicator
			of one component within a split $h$-fiber is an explicit extra
			common-side observable not generated by $h$.
		\end{enumerate}
	\end{theorem}
	\begin{proof}
		\leavevmode
		\begin{enumerate}[label=\textup{(\roman*)},leftmargin=*]
			\item \emph{Common functions.}
			A function common to the two side algebras has the form
			$g(u)=k(v)$ at every edge $(u,v)\in R$. These equalities hold
			if and only if its side values are constant along paths of
			$\Gamma_R$, giving \eqref{eq:proper-cut-components}.
			\item \emph{Exactness of the specified readout.}
			Since $h$
			factors through both side coordinates, two incident edges have
			the same $h$-value, and hence every connected component has a
			single $h$-value. The component partition of the edges equals
			the $h$-fiber partition exactly when each nonempty fiber
			induces one connected graph. If a fiber splits, assigning value
			one to one of its components and zero to the other components
			defines a common-side function not constant on that $h$-fiber.\qedhere
		\end{enumerate}
	\end{proof}
	
	\begin{theorem}[Exact one-cut source minimum and all equality cases]
		\label{thm:one-cut-spanning-minimum}
		Let $p:U\twoheadrightarrow Y$ and $q:V\twoheadrightarrow Y$
		be surjections of nonempty finite sets, and let
		$W=U\times_Y V$ be their full fiber product. For $y\in Y$ put
		$a_y=|p^{-1}(y)|$, $b_y=|q^{-1}(y)|$.
		\begin{enumerate}[label=\textup{(\roman*)},leftmargin=*]
			\item \emph{Minimum size.}
			Among all nonempty
			subrelations $R\subseteq W$ satisfying
			\begin{equation}\label{eq:one-cut-preservation-contract}
				\pi_U(R)=U,\qquad \pi_V(R)=V,\qquad
				M_U(R)\cap M_V(R)=\mathbb C[p\circ\pi_U|_R],
			\end{equation}
			the exact minimum number of joint states is
			\begin{equation}\label{eq:one-cut-minimum}
				\boxed{\quad \min|R|=\sum_{y\in Y}(a_y+b_y-1).\quad}
			\end{equation}
			\item \emph{Equality cases and their count.}
			Equality holds precisely when, for each $y$, the pairs retained
			by $R$ form a spanning tree of the complete bipartite graph
			$K_{a_y,b_y}$. Thus these are \emph{all} equality cases over
			the specified side states, and their number is
			\begin{equation}\label{eq:one-cut-minimizer-count}
				\prod_{y\in Y}a_y^{\,b_y-1}b_y^{\,a_y-1}.
			\end{equation}
			\item \emph{Application to a matching graph.}
			For any finite matching graph, apply this result with
			$U=U_L^Y$, $V=U_R^Y$ and $Y=Y_{\partial}$ from
			\cref{thm:arbitrary-graph-cut-information}. This optimizes a
			single named cut over the actual side images of the full source.
		\end{enumerate}
	\end{theorem}
	\begin{proof}
		\leavevmode
		\begin{enumerate}[label=\textup{(\roman*)},leftmargin=*]
			\item \emph{Minimum size.}
			Under the two projection requirements every fiber graph has
			all its $a_y+b_y$ vertices. By
			\cref{thm:proper-cut-incidence}, separator exactness is
			equivalent to connectivity in each fiber. A connected graph
			on $a_y+b_y$ vertices has at least $a_y+b_y-1$ edges, with
			equality exactly for a tree. Conversely the full fiber product
			contains every possible edge between the two sides of the same
			boundary value. Choose a spanning tree of each complete
			bipartite fiber. Their union satisfies all three conditions in
			\eqref{eq:one-cut-preservation-contract} and attains the sum.
			\item \emph{Equality cases and their count.}
			Equality requires equality in each fiber, so every fiber graph must be a spanning tree. To count all choices, the matrix--tree theorem \cite[Chapter~5]{Moon1970Counting} for
			$K_{a,b}$ gives $a^{b-1}b^{a-1}$ spanning trees: the
			nonzero Laplacian eigenvalues are $a$ with multiplicity $b-1$,
			$b$ with multiplicity $a-1$, and $a+b$ once, and the product
			divided by $a+b$ is the tree count. Multiply over the
			independent boundary fibers.
			\item \emph{Application to a matching graph.}
			The final specialization follows
			from the full-rectangle identity
			\eqref{eq:arbitrary-graph-cut-fiber}.\qedhere
		\end{enumerate}
	\end{proof}
	
	The diagonal access--commutant construction of \cref{thm:local-pair-native}
	realizes each minimizing relation. Requiring closure under given algebraic
	operations changes the feasible sets.
	
	\begin{example}[A three-state relational minimum]\label{ex:one-cut-native-gap}
		With two binary charts and a singleton separator, the relation
		$\{(0,0),(0,1),(1,0)\}$ is a minimizing tree. In the group
		$\mathbb F_2^2$ its pairs generate the full four-state product. The
		only proper subdirect subgroup is the diagonal, whose common
		information is nonconstant.
	\end{example}
	
	\subsubsection{Spanning forests and minimum native sources}
	\label{subsubsec:nonmalcev-forest-forcing}
	
	For two charts the spanning-forest construction yields all minimum
	sources. On a general join tree it yields a simultaneous certificate for
	the full cuts. In both cases the source uses the fixed finite native
	table-presentation law.
	
	Let $\mathbf A,\mathbf B,\mathbf Y$ be nonempty finite algebras of one
	fixed finite finitary signature, with surjective homomorphisms
	$p:\mathbf A\twoheadrightarrow\mathbf Y$ and
	$q:\mathbf B\twoheadrightarrow\mathbf Y$. Set
	\[
	\mathbf W=\mathbf A\times_{\mathbf Y}\mathbf B,
	\qquad A_y=p^{-1}(y),\quad B_y=q^{-1}(y).
	\]
	A source is \emph{state-and-cut exact} if its maps onto $A$ and $B$ are
	surjective and the two pulled-back finite function algebras intersect
	precisely in the functions of the common separator:
	\begin{equation}\label{eq:forest-exact-cut}
		\mathbb C[A]\cap\mathbb C[B]=\mathbb C[Y]
		\quad\text{inside }\mathbb C^H.
	\end{equation}
	The notation refers to the actual readouts on $H$, as in
	\eqref{eq:native-group-full-cut-contract}.
	For each $y\in Y$ let $K_y$ be the complete bipartite graph on
	$A_y\sqcup B_y$; its edges are exactly the pairs of $W$ above $y$.
	Write $\mathscr F$ for the finite family of edge sets
	\[
	F=\bigcup_{y\in Y}F_y\subseteq W,
	\qquad F_y\text{ is a spanning tree of }K_y,
	\]
	and let $\langle F\rangle_{\mathbf W}$ be the subalgebra of $\mathbf W$
	generated by $F$, including any nullary operations of the signature.
	
	For a subalgebra $T\leq W$ covering $A$ and $B$, apply
	\cref{thm:proper-cut-incidence} with $h$ the common separator readout.
	Then \eqref{eq:forest-exact-cut} holds precisely when the selected edges
	in each $K_y$ form a connected spanning graph.
	
	\begin{theorem}[Exact two-chart native minimum and forest rigidity test]
		\label{thm:forest-native-minimum}
		Let $p:A\twoheadrightarrow Y$ and $q:B\twoheadrightarrow Y$ be the
		finite algebra maps above. The following statements hold for their fiber product.
		\begin{enumerate}[label=\textup{(\roman*)}]
			\item \emph{Subalgebras generated by forests.} For every $F\in\mathscr F$, the subalgebra
			$T_F=\langle F\rangle_{\mathbf W}$ is state-and-cut exact.
			\item \emph{Forests in each feasible image.} Every state-and-cut exact native source $H$ has a joint image
			$T\leq\mathbf W$ containing $T_F$ for at least one $F\in\mathscr F$.
			\item \emph{Minimum sources.} The sharp native-source minimum is
			\begin{equation}\label{eq:forest-native-minimum}
				\boxed{\displaystyle N_{\rm native}(A,B;Y)
					=\min_{F\in\mathscr F}|\langle F\rangle_{\mathbf W}|.}
			\end{equation}
			The minimum sources over the prescribed chart maps are exactly the
			distinct generated subalgebras $T_F$ attaining this minimum;
			different spanning forests may generate the same minimum source.
			\item \emph{Full-image rigidity.} The full matching joint image is forced by the state-and-cut contract
			for \emph{all} finite native sources over this fixed atlas if and only if
			\begin{equation}\label{eq:forest-rigidity-test}
				\langle F\rangle_{\mathbf W}=\mathbf W
				\quad\text{for every }F\in\mathscr F.
			\end{equation}
		\end{enumerate}
	\end{theorem}
	\begin{proof}
		\leavevmode
		\begin{enumerate}[label=\textup{(\roman*)},leftmargin=*]
			\item \emph{Subalgebras generated by forests.}
			Each $T_F$ is a subalgebra of $W$ because $W$ is a fiber product of
			homomorphisms. It contains every edge of the spanning tree $F_y$;
			therefore it projects onto $A_y$ and $B_y$ and has a connected graph
			in that fiber. No elements outside $W$ arise by algebraic closure.
			\Cref{thm:proper-cut-incidence} proves (i).
			\item \emph{Forests in each feasible image.}
			For any source $H$, its comparison image $T\leq W$ has the same
			function-cut equality as $H$, since pullback along $H\twoheadrightarrow T$
			is injective on functions and identifies the displayed observation
			algebras. State exactness likewise descends. Hence every $G_y(T)$ is
			connected, where $G_y(T)$ has edge set $T\cap(A_y\times B_y)$
			and vertex parts $A_y,B_y$. Choose a spanning tree $F_y$ in $G_y(T)$ for each
			$y$. Their union $F$ lies in $T$, whence the generated subalgebra
			$T_F$ lies in $T$.
			\item \emph{Minimum sources.}
			Every $T_F$ is an allowed native source by the fixed finite
			presentation law; (i)--(ii) and $|H|\geq |T|\geq |T_F|$ prove the
			minimum formula. Equality forces both $H\to T$ and $T_F\subseteq T$
			to be bijective; thus $H\cong T_F$ over the chart maps. Conversely
			all minimizing $T_F$ attain the bound. 
			\item \emph{Full-image rigidity.}
			If every $T_F=W$, then (ii) forces every feasible image $T$ to
			contain $W$, and hence to equal $W$. If some $T_F\ne W$, part (i)
			itself furnishes a counterexample to full rigidity.\qedhere
		\end{enumerate}
	\end{proof}
	
	\begin{corollary}[Bounds and equality cases for two-chart native source costs]
		\label{cor:forest-cost-envelope}
		For every fixed finite two-chart atlas as above,
		\begin{equation}\label{eq:forest-cost-envelope}
			\sum_{y\in Y} (|A_y|+|B_y|-1)
			\ \le\ N_{\rm native}(A,B;Y)\ \le\ |W|.
		\end{equation}
		\begin{enumerate}[label=\textup{(\roman*)},leftmargin=*]
			\item \emph{Lower bound.}
			The lower bound is attained if and only if some
			$F\in\mathscr F$ is already a subalgebra of $W$.
			\item \emph{Upper bound.}
			The upper bound is attained if and only if every
			$F\in\mathscr F$ generates $W$.
		\end{enumerate}
	\end{corollary}
	\begin{proof}
		\leavevmode
		\begin{enumerate}[label=\textup{(\roman*)},leftmargin=*]
			\item \emph{Lower bound.}
			A connected graph on $|A_y|+|B_y|$ vertices has at least
			$|A_y|+|B_y|-1$ edges, with equality precisely for a spanning tree.
			Summation gives the lower bound, which is attained by a native
			source precisely when its fiberwise tree edge set is closed under
			all native operations.
			\item \emph{Upper bound.}
			The upper bound follows from the feasibility
			of $W$, and its equality criterion is
			\eqref{eq:forest-rigidity-test}.\qedhere
		\end{enumerate}
	\end{proof}
	
	The same reduction preserves negative recovery witnesses. Let $D$ be
	the finite set of named descriptors of the two charts, $c_0$ the
	closure generated by the actual local determination theories, and
	$c_W$ the descriptor closure on $W$. Define the finite set of
	nonimplication requirements
	\[
	\mathscr N_0=\{(S,d):S\subseteq D,\ d\in D\setminus c_0(S)\}.
	\]
	A \emph{witness selector} $P$ chooses, for each $(S,d)\in\mathscr N_0$,
	a pair $(u_{S,d},v_{S,d})\in W^2$ agreeing on $S$ but differing on $d$;
	write $V(P)$ for the set of endpoints of the selected pairs, with
	$V(P)=\varnothing$ if there are no requirements.
	
	\begin{theorem}[Faithful two-chart minimum from forests and pairs]
		\label{thm:forest-faithful-minimum}
		\leavevmode
		\begin{enumerate}[label=\textup{(\roman*)},leftmargin=*]
			\item \emph{Existence.}
			A state-and-cut exact, recovery-faithful native source exists
			if and only if $c_W=c_0$.
			\item \emph{Minimum and equality cases.}
			When this equality holds, witness
			selectors exist, and the sharp faithful minimum is
			\begin{equation}\label{eq:forest-faithful-cost}
				\boxed{\displaystyle
					N_{\rm faithful}(A,B;Y,c_0)
					=\min_{F\in\mathscr F,\ P}
					\bigl|\langle F\cup V(P)\rangle_{\mathbf W}\bigr|.}
			\end{equation}
			All equality cases over the named charts are precisely the distinct
			minimizing generated subalgebras in \eqref{eq:forest-faithful-cost}.
			\item \emph{Meet-irreducible witnesses.}
			The same minimum and equality cases result if $P$ instead chooses one
			pair with exact agreement set $M$ for each proper meet-irreducible
			$c_0$-closed set $M$.
		\end{enumerate}
	\end{theorem}
	\begin{proof}
		\leavevmode
		\begin{enumerate}[label=\textup{(\roman*)},leftmargin=*]
			\item \emph{Existence.}
			Any compatible source maps into $W$, and local dependencies hold
			on $W$, giving $c_0\subseteq c_W\subseteq c_H$. Faithfulness thus
			forces $c_W=c_0$. Conversely, under this equality $W$ itself is a
			faithful native source with the full state-and-cut contract.
			\item \emph{Minimum and equality cases.}
			Choose $F$ and $P$. The generated algebra
			$U=\langle F\cup V(P)\rangle_W$ contains $T_F$, so it is state-and-cut
			exact by \cref{thm:proper-cut-incidence}. Since $U\subseteq W$ it satisfies
			all dependencies in $c_0$; the chosen pairs refute every
			nonimplication outside $c_0$. Therefore $c_U=c_0$.
			
			For any faithful feasible source, let $T$ be its image in $W$.
			Choose a fiberwise spanning forest inside $T$ by
			\cref{thm:forest-native-minimum}, and choose the required witness
			pairs inside $T$. Their generated subalgebra lies in $T$, remains
			feasible and faithful, and has no more states. The minimum and
			classification follow because equality forces the map onto the joint
			image and its inclusion of the generated subalgebra both to be bijections.
			\item \emph{Meet-irreducible witnesses.}
			For the shorter selector, every pair agreement set in a subrelation of
			$W$ is $c_0$-closed. In a faithful source, a proper meet-irreducible $M$
			is an intersection of the agreement sets containing $M$, so one of them
			equals $M$. Conversely, every $c_0$-closed set is an intersection of
			proper meet-irreducibles above it, allowing the empty intersection for
			the top. If $d\notin c_0(S)$, a selected $M$ contains $S$ and omits
			$d$, and its pair witnesses the failed implication
			\cite[Theorem~6.1]{BeeriDFS1984}. Thus the shorter selection has precisely
			the same feasibility and equality cases.\qedhere
		\end{enumerate}
	\end{proof}
	
	\begin{corollary}[Minimum source sizes for meet-semilattices and cyclic groups]
		\label{cor:forest-meet-group-separation}
		\leavevmode
		\begin{enumerate}[label=\textup{(\roman*)},leftmargin=*]
			\item \emph{Meet-semilattices.}
			Let $A,B$ be nonempty finite meet-semilattices with respective
			least elements $0_A,0_B$, with singleton separator. Then
			\begin{equation}\label{eq:forest-meet-native-cost}
				N_{\wedge}(A,B;1)=|A|+|B|-1.
			\end{equation}
			The equality cases are exactly the meet-closed spanning trees of
			$K_{|A|,|B|}$, viewed as edge subsets of $A\times B$.
			In particular
			\[
			T_{\star}=(A\times\{0_B\})\cup(\{0_A\}\times B)
			\]
			is an equality case. If both chart descriptors are their full states
			and $|A|,|B|\geq2$, these minima also have the free two-descriptor
			recovery closure.
			\item \emph{Comparison with cyclic groups.}
			On the same labeled $m$-state chart carriers, for
			$m\geq2$, take either cyclic-group operations or chain-meet
			operations, with singleton separator and full-state descriptors.
			Then
			\begin{equation}\label{eq:forest-law-gap}
				N_{\rm group}=m^2,\qquad
				N_{\wedge}=2m-1,\qquad
				\frac{N_{\rm group}}{N_{\wedge}}\longrightarrow\infty.
			\end{equation}
		\end{enumerate}
	\end{corollary}
	\begin{proof}
		\leavevmode
		\begin{enumerate}[label=\textup{(\roman*)},leftmargin=*]
			\item \emph{Meet-semilattices.}
			By \cref{thm:proper-cut-incidence}, every feasible relation is a
			connected bipartite graph on $|A|+|B|$ vertices, hence has at least
			$|A|+|B|-1$ edges. The displayed $T_{\star}$ attains that bound and
			is closed under componentwise meet: meets within either arm remain
			in that arm, and a cross-arm meet is $(0_A,0_B)$. At equality the
			graph is a spanning tree, and native closure is exactly meet-closure;
			this proves the claimed exhaustive equality classification.
			For two nonconstant full-state descriptors, connectedness and
			surjective projections ensure a vertex of degree at least two on
			each side, witnessing both failed cross-descriptor determinations;
			state surjectivity witnesses both failed constant determinations.
			\item \emph{Comparison with cyclic groups.}
			Finally, a subgroup of $C_m\times C_m$ with connected incidence
			graph must equal the full product by Fleischer's lemma (or by the
			group rectangle identity), giving $m^2$ states. For a chain meet
			law the first part gives $2m-1$ states. Both have the same underlying
			full-state observation tables and the same separator contract.\qedhere
		\end{enumerate}
	\end{proof}
	
	\begin{theorem}[Intermediate native cost from the generating rank of a regular group action]
		\label{thm:regular-action-source-rank}
		Let $G$ be a finite group of order $n$, let $d(G)$ be its least
		number of group generators (with $d(1)=0$), and give both chart
		carriers the \emph{left regular $G$-set} algebra: its unary
		operations are $L_g(x)=gx$, for each $g\in G$. The separator is
		one point, and the two full-state chart descriptors are retained.
		\begin{enumerate}[label=\textup{(\roman*)},leftmargin=*]
			\item \emph{Minimum sources for a regular action.}
			The sharp state-and-cut exact native-source cost is
			\begin{equation}\label{eq:regular-action-source-rank}
				\boxed{N_{G\text{-}\mathrm{set}}(G,G;1)
					=n\bigl(d(G)+1\bigr).}
			\end{equation}
			All minimum sources over the named charts are exactly
			\[
			T_D=\{(x,xd):x\in G,\ d\in D\}\subseteq G\times G,
			\qquad |D|=d(G)+1,
			\quad\langle DD^{-1}\rangle=G.
			\]
			Here $\langle DD^{-1}\rangle$ denotes the subgroup generated by
			$\{de^{-1}:d,e\in D\}$. For $n\ge2$ these sources also preserve
			the free two-descriptor recovery closure.
			\item \emph{Cyclic actions and comparison of costs.}
			If $G=C_n$ is cyclic, $n\ge2$, the same result holds in the
			\emph{fixed one-unary-operation signature} with $s(x)=x+1\pmod n$,
			since closure under $s$ entails closure under all its powers.
			Then $N_{G\text{-}\mathrm{set}}=2n$,
			and the distinct minimum sources number
			$n\varphi(n)/2$, where $\varphi$ is Euler's totient function.
			Consequently, on the same two $n$-state chart tables, with full-state
			readouts and singleton separator,
			\begin{equation}\label{eq:regular-action-three-costs}
				N_{\rm relational}=2n-1,\qquad
				N_{C_n\text{-}\mathrm{set}}=2n,\qquad
				N_{\rm group}=n^2.
			\end{equation}
			For $n\ge3$ these are three distinct sharp costs.
		\end{enumerate}
	\end{theorem}
	\begin{proof}
		\leavevmode
		\begin{enumerate}[label=\textup{(\roman*)},leftmargin=*]
			\item \emph{Minimum sources for a regular action.}
			A subset of $G\times G$ is a subalgebra of the two regular
			$G$-sets precisely when it is invariant under every diagonal left
			translation $(x,y)\mapsto(gx,gy)$. These translations preserve the
			difference $x^{-1}y$, and the orbits are exactly
			$O_d=\{(x,xd):x\in G\}$, each of size $n$. Hence every nonempty
			native joint image is $T_D=\bigcup_{d\in D}O_d$ for a unique
			nonempty set $D\subseteq G$, and $|T_D|=n|D|$.
			
			Each $O_d$ is a perfect matching, so $T_D$ projects onto both
			charts. Two left vertices $x,z$ share a right neighbor precisely
			when $z=xde^{-1}$ for some $d,e\in D$.
			Thus the left vertices connected to $x$ in the bipartite graph
			are exactly $x\langle DD^{-1}\rangle$; all right vertices have
			neighbors. By \cref{thm:proper-cut-incidence}, the cut identity
			holds precisely when $\langle DD^{-1}\rangle=G$.
			
			Fix $d_0\in D$. The subgroup $\langle DD^{-1}\rangle$ is generated
			by the $|D|-1$ differences $dd_0^{-1}$ for
			$d\in D\setminus\{d_0\}$. Hence connectivity forces
			$|D|\ge d(G)+1$. Equality is attained by
			$D=\{1,g_1,\dots,g_{d(G)}\}$ for a minimal generating tuple of
			$G$. This proves the cost and the complete equality classification.
			Full-state recovery faithfulness for $n\ge2$ follows, as in
			\cref{cor:forest-meet-group-separation}, from connectivity and
			surjectivity.
			\item \emph{Cyclic actions and comparison of costs.}
			For a nontrivial cyclic group, a minimizing orbit-label set $D$
			consists of two elements $\{a,b\}$ with $ab^{-1}$ a generator. Counting ordered
			pairs gives $n\varphi(n)$ and dividing by two gives the number of
			minimum sources. Finally the relational minimum is the bipartite
			tree bound, and the group-native minimum is $n^2$ by the full-cut
			Mal'cev rigidity theorem. This proves \eqref{eq:regular-action-three-costs}.\qedhere
		\end{enumerate}
	\end{proof}
	
	\begin{proposition}[Congruence permutability and universal two-chart source rigidity]
		\label{prop:forest-malcev-boundary}
		For a variety $\mathcal V$, the following are equivalent:
		\begin{enumerate}[label=\textup{(\roman*)}]
			\item Every finite algebra of $\mathcal V$ has permuting congruences.
			\item For every two-chart atlas of finite algebras in $\mathcal V$
			with surjective separator maps, every state-and-cut exact native
			source has joint image equal to the full fiber product $W$.
		\end{enumerate}
		If $\mathcal V$ is locally finite, these are further equivalent
		to the existence of a Mal'cev term for $\mathcal V$.
		
	\end{proposition}
	\begin{proof}
		\leavevmode
		\begin{itemize}
			\item \emph{\textup{(i)} $\implies$ \textup{(ii)}.}
			Under (i), the two projection kernels of any finite joint image
			$T$ permute. Exact cut information identifies their join with the
			kernel of $T\to Y$, so Fleischer's rectangle argument makes $T=W$.
			\item \emph{\textup{(ii)} $\implies$ \textup{(i)}.}
			Conversely fix a finite algebra $\mathbf C\in\mathcal V$ and any
			$\alpha,\beta\in\operatorname{Con}(\mathbf C)$. Take
			$A=C/\alpha$, $B=C/\beta$, $Y=C/(\alpha\vee\beta)$, with the
			canonical quotient maps, and let $T$ be the joint image of
			$c\mapsto(c/\alpha,c/\beta)$. Its projections are onto and its
			common-function partition is $\alpha\vee\beta$, the partition
			of the prescribed separator; thus the cut contract holds.
			By (ii), $T=A\times_Y B$. Accordingly whenever
			$a\equiv b\pmod{\alpha\vee\beta}$ there is $c\in C$ with
			$a\equiv c\pmod\alpha$ and $c\equiv b\pmod\beta$.
			Hence $\alpha\circ\beta=\alpha\vee\beta$; taking inverses gives
			$\beta\circ\alpha=\alpha\vee\beta$. Congruences permute.
			\item \emph{Locally finite varieties.}
			If $\mathcal V$ is locally finite, its free algebra on
			$x,y,z$ is finite. Apply (i) to the congruences generated by
			$(x,y)$ and $(y,z)$. Their permutability gives a term
			$t(x,y,z)$ with $x\equiv t\pmod{\operatorname{Cg}(y,z)}$
			and $t\equiv z\pmod{\operatorname{Cg}(x,y)}$. Freeness yields
			$t(x,y,y)=x$ and $t(x,x,y)=y$, a Mal'cev term. Conversely, suppose $a\mathrel{\alpha}b\mathrel{\beta}c$.
			Compatibility gives
			$a=t(a,b,b)\mathrel{\beta}t(a,b,c)\mathrel{\alpha}t(b,b,c)=c$.
			Interchanging $\alpha$ and $\beta$ proves that they permute.\qedhere
		\end{itemize}
	\end{proof}
	
	\begin{corollary}[Congruence permutability and universal join-tree source rigidity]
		\label{cor:finite-join-tree-rigidity-boundary}
		For a fixed equational class $\mathcal V$, the following are equivalent:
		\begin{enumerate}[label=\textup{(\roman*)},nosep]
			\item Every finite algebra in $\mathcal V$ has permuting congruences.
			\item For every finite equational join-tree atlas in $\mathcal V$,
			every finite native source with surjective chart maps and exact
			represented common information at \emph{all complete tree cuts}
			maps onto the full native matching algebra.
		\end{enumerate}
		The converse in \textup{(ii)} requires only two-chart atlases.
		For locally finite $\mathcal V$, either condition is also equivalent to
		the existence of a Mal'cev term.
	\end{corollary}
	\begin{proof}
		\leavevmode
		\begin{itemize}
			\item \emph{\textup{(i)} $\implies$ \textup{(ii)}.}
			Under \textup{(i)}, the two projection kernels of each finite
			complete-side image permute. Apply
			\cref{lem:full-cut-kernel-rectangle} at every cut, followed by the
			leaf-removal argument of \cref{thm:native-source-rigidity}.
			\item \emph{\textup{(ii)} $\implies$ \textup{(i)}.}
			Conversely, specialize \textup{(ii)} to a single edge and apply
			\cref{prop:forest-malcev-boundary}. The same proposition proves the locally finite assertion.\qedhere
		\end{itemize}
	\end{proof}
	
	\begin{example}[Faithful attachment without full-cut rigidity in a Goursat variety]
		\label{ex:goursat-faithful-nonrigid}
		Use the three-element implication algebra
		$H=\{(0,1),(1,0),(1,1)\}$ from the example following
		\cref{prop:finite-goursat-boundary}, with its two coordinate projections,
		private full-state readouts, and a singleton separator. Its joint image
		consists of three edges forming a spanning path in $K_{2,2}$. Both
		coordinates occur surjectively, and by \cref{thm:proper-cut-incidence}
		their only common functions are constants. The two degree-two vertices
		witness failure of recovery in both directions. Thus the source preserves
		the local closure, states, and complete cut, while omitting one of the
		four matching states. It has the native table presentation of
		\cref{prop:equational-levelii-bridge}. The variety is three-permutable
		\cite{Mitschke1971}.
	\end{example}
	
	\paragraph{Transfer thresholds and native costs.}
	By \cref{prop:finite-goursat-boundary}, automatic
	relative pair lifting for every finite chart over every congruence-complete
	separator signature is equivalent to three-permutability of all finite
	members. In turn,
	\cref{cor:finite-join-tree-rigidity-boundary} makes universal rigidity under
	state coverage and exact full-cut information equivalent to permutability
	of all finite members. Permuting congruences imply three-permutability, and
	\cref{ex:goursat-faithful-nonrigid} separates the two finite-member
	conditions \cite{CarboniEtAl1993,Gran2020Regular}.
	
	For a fixed observed atlas, let $\mathscr T_{\Xi}$ be the nonempty joint
	images admitted as sources on their own state sets by a given native law.
	In the finite equational table contexts used here, all homomorphic joint
	images have such realizations. In any context admitting every compatible
	native source's joint image as a source over the same charts, the minimum
	for a contract consisting only of the
	observation clauses of \cref{def:preservation-contract} together with that
	native realization is the least $|T|$ among the members of
	$\mathscr T_{\Xi}$ passing the observation tests. This follows from the
	joint-image transfer and $|H|\geq|\psi(H)|$; each eligible $T$ is itself
	realized. If the families of attainable images for two such laws are nested,
	the more restrictive family has a weakly larger minimum (with $+\infty$
	for infeasibility); if every passing image in the larger family is $W$,
	the same holds in the smaller family. This comparison concerns fixed
	observations and actual image inclusion, not arbitrary changes of signature
	or additional native constraints. The forest-and-pair formula
	(\cref{thm:forest-faithful-minimum}) and the group, meet and regular-action
	costs (\cref{cor:forest-meet-group-separation,thm:regular-action-source-rank})
	exhibit the different native costs of the same information requirements.
	
	\begin{proposition}[Simultaneous forest compression at all full cuts]
		\label{prop:forest-allcuts-compression}
		Let $J$ be a join tree with at least two charts, and let
		$\mathbf T\leq\mathbf W$ be any finite joint image that is
		state-exact and satisfies the common-function contract at every
		\emph{full} tree cut; no Mal'cev assumption is made.
		For a cut $e$ with sides $L_e,R_e$, let
		$T_{L_e},T_{R_e}$ be the \emph{actual partial images} of $T$.
		For every separator state $y\in Y_e$, choose a spanning tree
		$F_{e,y}\subseteq T$ of the incidence graph on
		$(T_{L_e})_y\sqcup(T_{R_e})_y$.
		Put $F=\bigcup_{e,y}F_{e,y}$ and
		$\mathbf U=\langle F\rangle_{\mathbf W}$.
		\begin{enumerate}[label=\textup{(\roman*)}]
			\item \emph{Partial images.} $U\subseteq T$, and for each cut $e$ the projections of $U$
			onto $L_e$ and $R_e$ coincide with those of $T$.
			\item \emph{State and cut exactness.} $U$ is state-exact and satisfies every full-cut
			common-information identity.
			\item \emph{Number of generators.} $U$ is generated by no more than
			\begin{equation}\label{eq:forest-allcuts-bound}
				\sum_{e\in E(J)}
				\bigl(|T_{L_e}|+|T_{R_e}|-|Y_e|\bigr)
			\end{equation}
			source states.
			\item \emph{Preserving recovery.} If in addition $c_T=c_0$, adjoining one witness pair for each
			$(S,d)\in\mathscr N_0$ before generating gives an algebra with the same state and cut properties and recovery closure $c_0$, using at most
			$2|\mathscr N_0|$ additional generators. Together these generators certify state coverage, all full cuts, and
			recovery on the prescribed finite tables.
		\end{enumerate}
	\end{proposition}
	\begin{proof}
		\leavevmode
		\begin{enumerate}[label=\textup{(\roman*)},leftmargin=*]
			\item \emph{Partial images.}
			The finite-function argument of
			\cref{thm:proper-cut-incidence} applies to each full cut, using
			its actual partial images, so every separator-fiber incidence graph
			of $T$ is connected. Each chosen $F_{e,y}$ lies in $T$; hence its
			algebraic closure $U$ remains inside $T$. For fixed $e$, the
			selected tree edges visit \emph{every} vertex of $T_{L_e}$ and
			$T_{R_e}$. Therefore the two side projections of $U$ equal those
			of $T$.
			\item \emph{State and cut exactness.}
			The incidence graph of $U$ on these side images contains
			$F_{e,y}$ in every fiber and is connected there. This gives the
			same cut identity. As each chart occurs on one side of at least
			one cut, individual chart surjectivity follows.
			\item \emph{Number of generators.}
			A tree on the vertices of $(T_{L_e})_y\sqcup(T_{R_e})_y$ has
			$|(T_{L_e})_y|+|(T_{R_e})_y|-1$ edges. Summing over $y$ and over
			all cuts proves \eqref{eq:forest-allcuts-bound}; duplicate selected
			states only improve the bound.
			\item \emph{Preserving recovery.}
			For recovery, choose witness pairs
			in $T$ for every failed $c_0$-implication. Adding them does not
			leave $T$ and preserves the already proved cut and state contracts;
			as in \cref{thm:forest-faithful-minimum}, their persistence gives the enlarged algebra recovery closure $c_0$.\qedhere
		\end{enumerate}
	\end{proof}
	
	\begin{example}[Full-cut exactness need not fill partial joins]
		\label{ex:forest-proper-partial-joins}
		Take a three-chart path with empty separators, with each chart the
		binary meet algebra $\{0,1\}$ and its full-state descriptor. The
		meet-subalgebra
		\[
		T=\{000,100,010,001\}\subsetneq\{0,1\}^3
		\]
		is state-exact and has free descriptor recovery. At each full cut
		its bipartite incidence graph is a spanning tree, so both cut
		identities hold. Nevertheless its image on the first two charts
		omits $11$, whereas their full matching join is $\{0,1\}^2$.
		Thus the partial vertex sets in
		\cref{prop:forest-allcuts-compression} must be the actual images of
		$T$, not the full partial joins of $W$.
	\end{example}
	
	\subsubsection{Double cosets and orbit costs over nontrivial separators}
	\label{subsubsec:nonmalcev-coset-cost}
	
	In the unary algebras of a finite group action, native source relations
	are unions of diagonal orbits. These orbits replace the spanning
	forests of \cref{thm:forest-native-minimum} in the cost calculation.
	They give a test by deletion of one orbit and a closed formula when
	both local stabilizers are normal.
	
	Fix a finite group $G$ and subgroups $H,K\le L<G$. Equip
	\[
	A=G/H,\qquad B=G/K,\qquad Y=G/L
	\]
	with their left $G$-set operations, named as unary operations
	$x\mapsto gx$ for $g\in G$. The chart maps
	$p(gH)=gL$ and $q(gK)=gL$ are surjective homomorphisms.
	As before, a source is \emph{state-and-cut exact} when its joint
	image onto $A,B$ is surjective and its pulled-back finite
	function algebras have intersection exactly $\mathbb C[Y]$.
	Put $W=A\times_Y B$ and write $r=[G:L]>1$.
	For $C=H\ell K\in H\backslash L/K$ define
	\begin{equation}\label{eq:coset-orbit-source}
		O_C=\{(gH,g\ell K):g\in G\}\subseteq W,
		\qquad w_C=|O_C|=[G:H\cap\ell K\ell^{-1}].
	\end{equation}
	Both $O_C$ and $w_C$ are independent of the representative
	$\ell\in C$. For a family $\mathcal D\subseteq H\backslash L/K$
	set $T_{\mathcal D}=\bigcup_{C\in\mathcal D}O_C$ and
	$E_{\mathcal D}=\bigcup_{C\in\mathcal D}C\subseteq L$.
	
	The following is the bi-coset connectivity criterion of
	\cite{DuXu2000}; see also \cite[Proposition~2.1(iii)]{FengWang2023BiCoset}.
	We give its proof in the present left-coset notation.
	
	\begin{lemma}[The bi-coset connectivity criterion]
		\label{lem:coset-incidence-connectivity}
		For a nonempty family $\mathcal D$, the separator-fiber incidence
		graph of $T_{\mathcal D}$ is connected if and only if
		\begin{equation}\label{eq:coset-generation-criterion}
			\langle E_{\mathcal D}E_{\mathcal D}^{-1}\rangle=L.
		\end{equation}
		The statement holds in every separator fiber.
	\end{lemma}
	\begin{proof}
		Consider the fiber over $L\in G/L$. Its left vertices are the
		cosets $L/H$, its right vertices are $L/K$, and its edges are
		\[
		(uH,u eK),\qquad u\in L,\quad e\in E_{\mathcal D}.
		\]
		Since $E_{\mathcal D}$ is a union of $(H,K)$-double cosets,
		this edge relation is well defined on the cosets.
		A two-edge walk between left vertices has the form
		\[
		uH\ --\ ueK\ --\ uee'^{-1}H,
		\qquad e,e'\in E_{\mathcal D}.
		\]
		Consequently the left vertices reached from $H$ are exactly the
		cosets $mH$ with $m\in\langle E_{\mathcal D}
		E_{\mathcal D}^{-1}\rangle$. To see the converse inclusion in this
		reachability assertion, concatenate the displayed two-edge walks;
		the set $EE^{-1}$ is symmetric, and $H$ is contained in its generated
		subgroup because $E$ is left $H$-invariant and nonempty.
		Every right vertex has a neighbor since the edge set is nonempty
		and left $L$-invariant. Hence the graph is connected exactly under
		\eqref{eq:coset-generation-criterion}. Translation by $G$ identifies
		all separator fibers, proving the last assertion.
	\end{proof}
	
	\begin{theorem}[Minimum source size and rigidity for finite coset-action atlases]
		\label{thm:coset-orbit-native-minimum}
		\leavevmode
		\begin{enumerate}[label=\textup{(\roman*)},leftmargin=*]
			\item \emph{Minimum and equality cases.}
			For the above finite $G$-set atlas, the minimum cardinality of a
			state-and-cut exact native source is
			\begin{equation}\label{eq:coset-orbit-native-minimum}
				\boxed{\displaystyle
					N_{G\text{-set}}(G/H,G/K;G/L)
					=\min_{\substack{\varnothing\ne\mathcal D\subseteq H\backslash L/K\\
							\langle E_{\mathcal D}E_{\mathcal D}^{-1}\rangle=L}}
					\sum_{C\in\mathcal D}[G:H\cap\ell_CK\ell_C^{-1}].}
			\end{equation}
			All minimum sources, up to isomorphism over the fixed chart maps,
			are exactly the distinct $T_{\mathcal D}$ attaining this minimum.
			\item \emph{Single-orbit rigidity test.}
			Full matching-image rigidity holds if and only if for every
			$C\in H\backslash L/K$ the relation
			\begin{equation}\label{eq:coset-single-orbit-deletion}
				W\setminus O_C
			\end{equation}
			is \emph{not} state-and-cut exact. Equivalently, for every $C$ either
			$E_{(H\backslash L/K)\setminus\{C\}}=\varnothing$ or
			\[
			\left\langle E_{(H\backslash L/K)\setminus\{C\}}
			E_{(H\backslash L/K)\setminus\{C\}}^{-1}\right\rangle\ne L.
			\]
			Thus testing full rigidity requires at most $|H\backslash L/K|$
			subgroup-generation tests, without listing spanning trees.
		\end{enumerate}
	\end{theorem}

	\begin{proof}
		\leavevmode
		\begin{enumerate}[label=\textup{(\roman*)},leftmargin=*]
			\item \emph{Minimum and equality cases.}
			A subalgebra of the $G$-set $W$ is exactly a $G$-invariant subset.
			The $G$-orbits of $W$ are the sets $O_C$ in
			\eqref{eq:coset-orbit-source}: the invariant $H g^{-1}g'K$
			records the double coset of a pair $(gH,g'K)\in W$.
			The stabilizer of $(H,\ell K)$ is $H\cap\ell K\ell^{-1}$,
			so orbit-stabilizer gives the stated weight.
			
			Every nonempty orbit already projects onto both transitive charts.
			By \cref{thm:proper-cut-incidence,lem:coset-incidence-connectivity}, a union $T_{\mathcal D}$ obeys
			the exact cut contract if and only if the subgroup-generation
			condition holds. Any compatible source has a joint image
			$T_{\mathcal D}$, and a minimum source is isomorphic to its image,
			just as in \cref{thm:forest-native-minimum}.
			This proves \eqref{eq:coset-orbit-native-minimum} and its equality cases.
			\item \emph{Single-orbit rigidity test.}
			If a proper union of orbits is feasible, it omits some orbit $O_C$.
			The larger relation $W\setminus O_C$ remains feasible because
			connectivity is monotone under adding edges; it is a proper native
			source. Conversely, a feasible single-orbit deletion itself
			disproves full rigidity. The connectivity lemma gives the equivalent
			subgroup tests.\qedhere
		\end{enumerate}
	\end{proof}

	\begin{proposition}[Cost bounds and rigidity regimes for finite coset-action atlases]
		\label{prop:coset-cost-trichotomy}
		Write $a=[L:H]$, $b=[L:K]$ and
		$
		R=r(a+b-1),\quad U=rab=|W|.
		$
		The following hold.
		\begin{enumerate}[label=\textup{(\roman*)}]
			\item \emph{Relational lower bound.} The relational lower bound $N_{G\text{-set}}=R$ is attained
			if and only if $a=1$ or $b=1$. In that case $R=N_{G\text{-set}}=U$.
			
			\item \emph{Intermediate cost and rigidity.} If $a,b>1$, then $R<N_{G\text{-set}}\le U$. The minimum is
			strictly intermediate if and only if at least one single-orbit
			deletion \eqref{eq:coset-single-orbit-deletion} passes the subgroup
			generation test. Otherwise $N_{G\text{-set}}=U$ and the full
			matching algebra is forced.
			
			\item \emph{Divisibility bound.} Every feasible source has at least
			\begin{equation}\label{eq:coset-biregular-lower}
				r\,\operatorname{lcm}(a,b)
				\left\lceil\frac{a+b-1}{\operatorname{lcm}(a,b)}\right\rceil
			\end{equation}
			states. In particular, if $\gcd(a,b)=1$, full rigidity holds.
		\end{enumerate}
	\end{proposition}

	\begin{proof}
		\leavevmode
		\begin{enumerate}[label=\textup{(\roman*)},leftmargin=*]
			\item \emph{Relational lower bound.}
			If $a=1$, each right chart state has only one compatible left
			state, so state exactness forces all of $W$; similarly for $b=1$.
			This also gives $R=U$.
			
			Suppose $a,b>1$. In each separator fiber, the group $L$ acts
			transitively on both parts of the incidence graph of any native
			source. Thus all left vertices have one common degree and all
			right vertices have one common degree. The edge count $e$ is
			divisible by both $a$ and $b$. A connected graph with exactly
			$a+b-1$ edges would be a tree, so we would have
			$a\mid(b-1)$ and $b\mid(a-1)$, which is impossible for
			$a,b>1$. The strict lower-bound statement follows.
			\item \emph{Intermediate cost and rigidity.}
			The single-orbit test of \cref{thm:coset-orbit-native-minimum}
			then distinguishes intermediate cost from full rigidity.
			\item \emph{Divisibility bound.}
			In general, connectivity requires $e\ge a+b-1$, and transitivity
			makes $e$ a positive multiple of $\operatorname{lcm}(a,b)$.
			Multiplying the resulting bound by $r$ proves
			\eqref{eq:coset-biregular-lower}. If $\gcd(a,b)=1$, the number
			of edges per fiber is divisible by $ab$ but is at most $ab$,
			hence is exactly $ab$.\qedhere
		\end{enumerate}
	\end{proof}
	
	\begin{theorem}[Closed-form minimum for normal separator stabilizers]
		\label{thm:coset-normal-stabilizer-rank}
		Assume in addition that $H$ and $K$ are normal subgroups of $L$.
		Let
		\[
		Q=L/(HK),\qquad d(Q)=\text{minimum number of generators of }Q,
		\qquad d(1)=0.
		\]
		\begin{enumerate}[label=\textup{(\roman*)},leftmargin=*]
			\item \emph{Minimum sources.}
			The minimum is
			\begin{equation}\label{eq:coset-normal-stabilizer-cost}
				\boxed{\displaystyle
					N_{G\text{-set}}(G/H,G/K;G/L)
					=[G:H\cap K]\bigl(d(Q)+1\bigr).}
			\end{equation}
			The minimum joint images are exactly the unions of
			$d(Q)+1$ distinct double-coset orbits whose corresponding
			subset $D\subseteq Q$ satisfies $\langle DD^{-1}\rangle=Q$.
			\item \emph{Full rigidity.}
			The full matching algebra is forced if and only if
			\begin{equation}\label{eq:coset-normal-full-rigidity}
				|Q|\le2.
			\end{equation}
			For $a,b>1$, it is strictly intermediate precisely when
			$|Q|\ge3$; when $a=1$ or $b=1$, the relational lower bound
			and full-matching bound coincide.
		\end{enumerate}
	\end{theorem}
	\begin{proof}
		\leavevmode
		\begin{enumerate}[label=\textup{(\roman*)},leftmargin=*]
			\item \emph{Minimum sources.}
			Normality gives $H\ell K=\ell HK$, so the double cosets
			are identified with the elements of $Q=L/(HK)$.
			The chosen representative $(H,\ell K)$ of each orbit has stabilizer
			$H\cap K$, so every orbit has weight $[G:H\cap K]$.
			Under the quotient map $L\twoheadrightarrow Q$, the criterion
			$\langle EE^{-1}\rangle=L$ is equivalent to
			$\langle DD^{-1}\rangle=Q$. Indeed $E$ is nonempty and is a union of
			$HK$-cosets. For any $e\in E$, the set $Ee^{-1}\subseteq EE^{-1}$
			contains $HK$. Hence the generated subgroup contains the quotient kernel.
			
			Fix $d_0\in D$. Every element $dd'^{-1}$ is a product of
			$dd_0^{-1}$ and $(d'd_0^{-1})^{-1}$, so
			$\langle DD^{-1}\rangle$ is generated by the at most
			$|D|-1$ differences $dd_0^{-1}$. Hence
			$|D|\ge d(Q)+1$. Equality is attained by
			$D=\{1,q_1,\dots,q_{d(Q)}\}$ for a minimum generating
			family of $Q$. This proves the formula and its equality cases.
			\item \emph{Full rigidity.}
			The whole matching source has exactly $|Q|$ double-coset
			orbits. Therefore the minimum is the whole source if and only if
			$d(Q)+1=|Q|$. This equality holds for groups of orders one and
			two. For every finite group of order at least three, a chain of
			strictly increasing generated subgroups at least doubles its size
			at each step, yielding $d(Q)\le\lfloor\log_2|Q|\rfloor$
			and $d(Q)+1<|Q|$. The remaining claims follow from
			\cref{prop:coset-cost-trichotomy}.\qedhere
		\end{enumerate}
	\end{proof}
	
	\begin{example}[Nontrivial separator with a full native source]
		\label{ex:coset-nonmalcev-full}
		Let $L=C_2\times C_2=\langle u\rangle\times\langle v\rangle$,
		$G=L\times C_2$, $H=\langle u\rangle$,
		$K=\langle v\rangle$. Then $Y=G/L$ has two states and
		every separator fiber has two left and two right chart states.
		Since $HK=L$, the quotient $Q$ is trivial and
		\cref{thm:coset-normal-stabilizer-rank} gives
		$N_{G\text{-set}}=|W|=8$, while the relational lower bound
		is $2(2+2-1)=6$. This is full rigidity in a unary,
		non-Mal'cev signature; it follows here from transitivity on the fiber pairs.
	\end{example}
	
	\begin{corollary}[Source minima for regular actions over subgroup coset separators]
		\label{cor:coset-principal-bundle-rank}
		\leavevmode
		\begin{enumerate}[label=\textup{(\roman*)},leftmargin=*]
			\item \emph{Principal bundles.}
			In \cref{thm:coset-normal-stabilizer-rank} take $H=K=1$.
			Then $Q=L$ and
			\begin{equation}\label{eq:coset-principal-bundle-cost}
				N_{G\text{-set}}(G,G;G/L)=|G|\bigl(d(L)+1\bigr).
			\end{equation}
			\item \emph{Cyclic separator.}
			In particular, for integers $q\ge2$ and $m\ge3$, identify
			$G=C_{qm}$ and its order-$m$ subgroup $L$, with separator
			$G/L\cong C_q$. The chart and separator observations are
			$x\mapsto x\bmod q$ on two copies of $\mathbb Z/(qm)$.
			With the \emph{single unary operation}
			$s(x)=x+1\pmod{qm}$, the minimum is
			\[
			N_{\mathrm{unary}}=2qm,
			\qquad |W|=qm^2,
			\qquad R=q(2m-1).
			\]
			This uses a fixed one-unary-operation signature as $q,m$ vary.
		\end{enumerate}
	\end{corollary}
	\begin{proof}
		\leavevmode
		\begin{enumerate}[label=\textup{(\roman*)},leftmargin=*]
			\item \emph{Principal bundles.}
			The first statement is \eqref{eq:coset-normal-stabilizer-cost}
			with $H=K=1$.
			\item \emph{Cyclic separator.}
			For cyclic $L=C_m$, $d(L)=1$.
			A source closed under the permutation $s$ is closed under all
			its powers, and hence is a union of the diagonal $C_{qm}$-orbits.
			Thus the one-operation and full-action classes have the same
			subalgebras and the same sharp minimum.\qedhere
		\end{enumerate}
	\end{proof}
	
	\begin{proposition}[Three sharp native source costs over the same nontrivial separator]
		\label{prop:coset-three-law-spectrum}
		For $q\ge2$ and $m\ge3$, fix the same two finite chart-state
		tables $X=\mathbb Z/(qm)$ and the same separator observation
		$p:X\twoheadrightarrow\mathbb Z/q$, $p(x)=x\bmod q$.
		\begin{enumerate}[label=\textup{(\roman*)},leftmargin=*]
			\item \emph{Three minimum costs.}
			Under the state-and-cut preservation contract, the minima for
			three different permitted native laws are
			\begin{equation}\label{eq:coset-three-law-spectrum}
				\boxed{\displaystyle
					N_{\mathrm{meet}}=q(2m-1)
					\;<\; N_{\mathrm{unary}}=2qm
					\;<\; N_{\mathrm{group}}=qm^2.}
			\end{equation}
			Here the meet-native law is obtained by identifying
			$X\cong\{0,\ldots,q-1\}\times\{0,\ldots,m-1\}$
			via $x=y+qj$, and using the coordinatewise minimum operation;
			its separator has the minimum operation on
			$\{0,\ldots,q-1\}$.
			The unary law and the group law are respectively
			$s(x)=x+1\pmod{qm}$ and addition modulo $qm$;
			with these laws the separator has the corresponding cyclic
			operation or quotient-group operation.
			These comparisons retain \emph{identical underlying chart-state
				and separator observation maps} but allow different native laws.
			\item \emph{Recovery closure.}
			With descriptors naming the two full chart states and their
			common separator, all three examples also preserve the same
			local-consequence functional-dependency closure.
		\end{enumerate}
	\end{proposition}
	\begin{proof}
		\leavevmode
		\begin{enumerate}[label=\textup{(\roman*)},leftmargin=*]
			\item \emph{Three minimum costs.}
			For the meet law, the matching algebra is the product
			of the separator chain and two $m$-element meet chains.
			Inside it, the subset
			\[
			T_{\wedge}=\left\{(y,i,j):
			i=0\text{ or }j=0\right\}
			\]
			is closed under coordinatewise minimum: a pair of points in the
			same arm remains in that arm, while a pair in opposite arms
			has both fiber coordinates zero. Over each $y$ its incidence
			graph is the double star with $2m-1$ edges, giving exact cut
			information and attaining the relational lower bound.
			The unary minimum is the preceding corollary.
			The cyclic group has a Mal'cev term, so the full-cut rigidity
			theorem \cref{thm:native-source-rigidity} forces the whole
			matching subgroup $W$ with $qm^2$ elements.
			\item \emph{Recovery closure.}
			For the descriptor assertion, in every one of the three sources
			both chart readouts determine the separator, but the separator
			determines neither chart state; neither chart readout determines
			the other. In each fiber the meet source has a vertex of degree $m$
			on each side, the unary minimum has degree two on each side,
			and the group source is a full rectangle. As $q\ge2$, the separator is not constant.
			These exhaust the functional dependencies among the three
			named descriptors and agree with the closure on $W$.\qedhere
		\end{enumerate}
	\end{proof}
	
	For the double-coset description of orbits and connectivity of bi-coset graphs, see
	\cite{FengWang2023BiCoset,GyurkiEtAl2024BiCoset}.
	
	\subsection{Logical implications among preservation contracts}
	\label{subsec:exact-contract-logic}
	
	Full-cut rigidity determines which preservation clauses imply the
	remaining clauses in a faithful native tree atlas.
	
	\begin{theorem}[Complete contract logic for faithful native tree atlases]
		\label{thm:complete-contract-logic}
		Consider finite Mal'cev join-tree atlases with onto interface maps and
		recovery-faithful full matching source $W$, so $c_W=c_0$.
		Candidates are finite native algebras with compatible homomorphic
		chart maps, realized in the fixed presentation law.
		\begin{enumerate}[label=\textup{(\roman*)},leftmargin=*]
			\item \emph{Universal implications.}
			Across this class the exact possible truth values of
			$(\mathsf C,\mathsf S,\mathsf F,\mathsf E)$ are the nine rows below.
			Equivalently, all universally valid logical implications among
			these four clauses follow from
			\begin{equation}\label{eq:complete-contract-logic}
				\mathsf S\Rightarrow\mathsf C,\qquad
				\mathsf F\Rightarrow\mathsf C,\qquad
				\mathsf S\wedge\mathsf E\Rightarrow\mathsf F.
			\end{equation}
			\item \emph{Realization on a fixed atlas.}
			Every row occurs on one fixed two-chart atlas in the native
			$\mathbb F_q$-module law, for any fixed finite field.
			\begin{center}
				\begin{tabular}{@{}cccc@{\qquad}cccc@{}}
					\toprule
					$\mathsf C$&$\mathsf S$&$\mathsf F$&$\mathsf E$&
					$U$&$U'$&$T$&Native clause\\
					\midrule
					0&0&0&0&$Z$&$P$&$Y$&$\mathsf N$\\
					0&0&0&1&$Z$&$A$&$Y$&$\mathsf N$\\
					1&0&0&0&$A$&$P$&$T_0$&$\mathsf N$\\
					1&0&0&1&$A$&$A$&$T_0$&$\mathsf N$\\
					1&0&1&0&$L$&$P$&$Y$&$\mathsf N$\\
					1&0&1&1&$L$&$A$&$Y$&$\mathsf N$\\
					1&1&0&0&$D$&$A$&$Y$&$\mathsf N$\\
					1&1&1&0&$P$&$A$&$Y$&$\mathsf N$\\
					1&1&1&1&$A$&$A$&$Y$&$\mathsf N$\\
					\bottomrule
				\end{tabular}
			\end{center}
			The sources denoted by $(U,U',T)$ are defined in the proof.
			In particular, the clauses $\mathsf S,\mathsf F,\mathsf E$ are
			pairwise incomparable in this class. The conjunction
			$\mathsf S\wedge\mathsf E$ forces surjectivity onto $W$ and all four clauses.
		\end{enumerate}
	\end{theorem}
	\begin{proof}
		\leavevmode
		\begin{enumerate}[label=\textup{(\roman*)},leftmargin=*]
			\item \emph{Universal implications.}
			State surjectivity preserves every local counterexample pair, giving
			$\mathsf S\Rightarrow\mathsf C$. Since the full matching source
			is state-exact, its closure $c_W=c_0$ has the prescribed local
			traces. Thus $\mathsf F\Rightarrow\mathsf C$.
			By \cref{thm:native-source-rigidity},
			$\mathsf S\wedge\mathsf E$ is equivalent to surjectivity of the
			joint chart map onto $W$. Surjective pullback preserves its full
			determination closure, proving the last implication.
			The displayed implications allow exactly the nine truth assignments
			in the table.
			\item \emph{Realization on a fixed atlas.}
			We realize them all without changing the atlas or native law.
			Put $Y=\mathbb F_q^3$, with coordinates $s_1,s_2,s_3$, and use
			the two chart state spaces
			\[
			X_1=\mathbb F_q^2\times\mathbb F_q^2\times Y,\qquad
			X_2=\mathbb F_q^2\times\mathbb F_q^2\times Y.
			\]
			The first chart has whole-block readouts $u,u',s$ and the scalar
			readout $a=s_1$; the second has $v,v',s$ and $b=s_2$.
			Their only shared descriptor is the whole-state readout $s$.
			The matching source is
			$W=(\mathbb F_q^2)^4\times Y$.
			Its closure is generated by $s\to a,b$: the four private blocks
			are independent, and $a,b$ together do not determine $s$.
			The same rules are generated by the two local theories, so
			$c_W=c_0$.
			
			For one pair $(u,v)$ of private blocks define
			\[
			\begin{array}{ll}
				A=\mathbb F_q^2\times\mathbb F_q^2,&
				Z=\{0\}\times\mathbb F_q^2,\\
				L=\langle(1,0)\rangle\times\mathbb F_q^2,&
				P=\{(u,v):u_1=v_1\},\\
				D=\{(u,v):u=v\}.&
			\end{array}
			\]
			Use the same notation for $(u',v')$, and put
			$T_0=\{s\in Y:s_1=s_2\}$.
			Each row specifies the submodule $H=U\times U'\times T$ of $W$.
			Every such source has its finite-table native presentation and
			the unchanged coordinate chart maps.
			
			For a private pair, $A$ has full coverage, free two-readout
			closure, and trivial common information. The relation $Z$ makes its left observation constant. In $L$ both
			observations vary independently, and some left chart states are missing.
			The relation $P$ covers both whole blocks and keeps the
			two-readout closure free, since fixing either block leaves
			$q$ values of the other. Its common information is the
			nonconstant coordinate $u_1=v_1$, so it fails exactness over
			the empty private separator. The diagonal $D$ covers both blocks. Each block determines
			the other, so they share nonconstant information.
			
			On $T_0$, each local theory on $\{s,a\}$ or $\{s,b\}$ is
			unchanged: the scalar is nonconstant, is determined by $s$,
			and does not determine the two-dimensional actual image of
			$s$. Globally $a=b$, introducing an extra dependency.
			The represented separator is still exact, because both
			sides see the entire actual $s$-value.
			
			The three factors are independent. For disjoint readout
			families the determination closure of a product is generated
			by the factor closures: a failed factor implication retains
			its counterexample with all other factors fixed. Chart
			surjectivity and chart-closure exactness likewise hold
			exactly when they hold in every factor. For linear left,
			right, and separator readouts $\ell,r,s$, cut exactness is
			$\ker\ell+\ker r=\ker s$. This equality splits over the
			independent factors. Thus separator exactness is also the
			conjunction of the factor tests.
			These observations give precisely the table entries.
			Since all allowed truth assignments occur, any Boolean
			implication valid throughout the class is a logical
			consequence of \eqref{eq:complete-contract-logic}.\qedhere
		\end{enumerate}
	\end{proof}
	
	\begin{corollary}[Scalar observations collapse the full contract to faithfulness]
		\label{cor:scalar-contract-equivalence}
		For a recovery-faithful finite-field join-tree atlas whose named
		readouts are scalar linear forms and jointly separate the chart
		states, every compatible linear native source satisfies
		\begin{equation}\label{eq:scalar-contract-equivalence}
			\mathsf F
			\quad\Longleftrightarrow\quad
			\psi(H)=W
			\quad\Longleftrightarrow\quad
			\mathsf S\wedge\mathsf E.
		\end{equation}
		Hence $\mathsf F$ alone is a sufficient and necessary observation
		test for the full faithful native contract in this class.
	\end{corollary}
	\begin{proof}
		The scalar readouts span $W^*$, so their matroid closure on $W$
		has rank $\dim W$. Equality of determination closures preserves
		matroid rank. On $\psi(H)\le W$ the restricted readouts span
		$\psi(H)^*$, so $\mathsf F$ forces
		$\dim\psi(H)=\dim W$ and therefore $\psi(H)=W$.
		Surjective pullback proves the converse, and
		\cref{thm:native-source-rigidity} proves the other equivalence.
		The given source already has its native presentation.
	\end{proof}
	
	\subsection{Faithfulness, state coverage and separators beyond trees}
	\label{subsec:contract}
	
	We use the setting of \cref{thm:simple-diagram-source-rigidity}: a finite
	atlas of algebras $A_i$ in a Mal'cev variety $\mathcal V$ on an arbitrary
	finite diagram, with onto descriptor readouts
	$f_{i,d}:A_i\twoheadrightarrow S_d$ that jointly separate each chart, a
	nonempty full matching set $W$, and candidate sources $H$ with
	compatible homomorphic chart maps and comparison $\psi:H\to W$.
	The comparison image may omit local states. Let $\mathcal D$ be a
	class of finite algebras in $\mathcal V$, each having at least two elements.
	Use the matching and native presentation conventions of the cited
	rigidity theorem; in particular candidates are nonempty and the global
	descriptor readouts jointly separate $W$.
	
	\begin{theorem}[Signatures for which faithfulness forces the full source]
		\label{thm:faithful-forcing}
		The following are equivalent.
		\begin{enumerate}[label=\textup{(\roman*)},nosep]
			\item Every $D\in\mathcal D$ is simple, and every subalgebra of $D$
			with at least two elements equals $D$.
			\item For every atlas as above with $\mathcal D$-valued readouts and
			every candidate $H$, the equality $c_H=c_W$ implies that $\psi$ is
			onto.
		\end{enumerate}
		Under \textup{(i)}, every candidate on every finite diagram satisfies
		\begin{equation}\label{eq:F-forces}
			\mathsf F\quad\Longleftrightarrow\quad
			\bigl(c_W=c_0\ \text{and}\ \psi(H)=W\bigr).
		\end{equation}
		For a finite Mal'cev join-tree atlas with onto interfaces, $c_W=c_0$,
		and all full edge cuts protected, it follows that
		$\mathsf F\Leftrightarrow\psi(H)=W\Leftrightarrow\mathsf S\wedge\mathsf E$.
	\end{theorem}
	\begin{proof}
		\leavevmode
		\begin{itemize}
			\item \emph{\textup{(i)} $\implies$ \textup{(ii)}.}
			Assume (i) and $c_H=c_W$, and put $Z=\psi(H)$, so that $c_Z=c_H$. Let
			$D'$ be the set of descriptors that are nonconstant on $W$. For
			$d\in D'$ the empty premise does not determine $d$ on $W$, hence not on
			$Z$. Thus $f_d(Z)\subseteq f_d(W)$ are subalgebras of $S_d$ with at
			least two elements, and both equal $S_d$ by (i). The readouts in $D'$
			still jointly separate $W$ and $Z$. Choose $B\subseteq D'$
			inclusion-minimal with $D'\subseteq c_W(B)$; since $c_Z=c_W$, it is also
			minimal for $Z$. By \cref{lem:simple-coordinate-basis}, $f_B$ maps both
			$W$ and $Z$ isomorphically onto $\prod_{d\in B}S_d$. As $Z\subseteq W$
			and $f_B$ is injective on $W$, we get $Z=W$.
			\item \emph{\textup{(ii)} $\implies$ \textup{(i)}.}
			If some $D\in\mathcal D$ is not simple, the atlas in the proof of
			\cref{prop:simple-descriptor-boundary} gives a proper $Z$ with
			$c_Z=c_W$. If $D$ is simple and has a proper subalgebra $D_0$ with
			$|D_0|\ge2$, use one chart $D$ with its identity readout and the
			candidate $D_0\hookrightarrow D$. Both closures make the single
			descriptor undetermined by the empty premise, so they agree, but $\psi$
			is not onto.
			\item \emph{Consequences for faithfulness.}
			For \eqref{eq:F-forces}, a surjective comparison gives $c_H=c_W$ by
			pullback. Conversely, $c_0\le c_W\le c_Z=c_H$ always holds, so
			$\mathsf F$ gives $c_W=c_0=c_H$ and then $\psi(H)=W$. On join trees,
			\cref{thm:native-source-rigidity} identifies $\psi(H)=W$ with
			$\mathsf S\wedge\mathsf E$.\qedhere
		\end{itemize}
	\end{proof}
	
	For modules over a finite ring, subalgebras are submodules, so condition
	(i) says that every descriptor module is simple. Scalar linear forms over
	$\mathbb F_q$ give \cref{cor:scalar-contract-equivalence}; a readout onto
	$\mathbb F_q^m$ with $m\ge2$ already violates (i). For groups, (i) holds
	exactly for cyclic groups of prime order. The example $A_5$ with
	$\langle(12)(34)\rangle$ after \cref{prop:simple-descriptor-boundary}
	realizes the second obstruction in the proof.
	
	The second equivalence of \cref{cor:scalar-contract-equivalence} does not
	survive on cycles, even for scalar readouts.
	
	\begin{proposition}[A separator-exact parity source on every cycle]
		\label{prop:parity-cycle}
		For every $n\ge3$ and every finite field $\mathbb F_q$ there is an $\mathbb F_q$-module atlas on the cycle
		with $n$ charts, with onto interfaces and scalar descriptor readouts,
		such that $W$ is state-exact and $c_W=c_0$, and a proper submodule
		$Z<W$ satisfies $\mathsf C$, $\mathsf S$ and $\mathsf E$ at every
		bipartition of the charts with complete crossing boundaries, but not
		$\mathsf F$.
	\end{proposition}
	\begin{proof}
		Number the charts $1,\dots,n$ cyclically and let $e_m$ join $m$ and
		$m+1$, indices modulo $n$, with interface $\mathbb F_q$. Its value is
		the descriptor $y_m$. Charts $1,2,3$ have state space $\mathbb F_q^2$:
		chart $1$ reads $(y_n,y_1)$, chart $2$ reads $(y_1,y_2)$ and chart $3$
		reads $(y_2,y_3)$. For $4\le m\le n$, chart $m$ is a relay with state
		space $\mathbb F_q$ and identity interface maps, so $y_{m-1}=y_m$. Then
		$W\cong\mathbb F_q^3$ with coordinates $a=y_1$, $b=y_2$ and
		$c=y_3=\dots=y_n$. It is state-exact, and its closure is generated by
		the relay equalities, which are local rules, so $c_W=c_0$. Let
		$Z=\{a+b+c=0\}$. Each corner chart reads two of $a,b,c$, which are
		independent on $Z$, so $\mathsf S$ holds and hence $\mathsf C$. On $Z$
		the premise $\{y_1,y_2\}$ determines $y_3$, which fails on $W$, so
		$\mathsf F$ fails.
		
		Let $L\sqcup R$ be a bipartition. Its crossing edges form an even,
		positive number of edges of the cycle. Suppose two of them lie in
		different classes among $\{e_1\}$, $\{e_2\}$ and $\{e_3,\dots,e_n\}$.
		Then the boundary determines two of $a,b,c$, hence a point of $Z$, and
		every common function of the two sides is a function of the boundary.
		Otherwise all crossing edges lie in $\{e_3,\dots,e_n\}$. Then charts
		$1,2,3$ are on one side, and the other side consists of relays, whose
		states are functions of $c$, the boundary value. In both cases
		$\mathsf E$ holds.
	\end{proof}
	
	\begin{theorem}[Contract logic on arbitrary diagrams and cut families]
		\label{thm:contract-logic-general}
		Fix a finite connected simple graph $G$ with at least two vertices and
		any family $\mathcal K$ of bipartitions of its vertex set into two
		nonempty parts. Either part may be disconnected. Consider finite Mal'cev atlases with chart graph $G$, onto
		interface maps, state-exact $W$ and $c_W=c_0$; test $\mathsf E$ at the
		cuts in $\mathcal K$ with complete crossing boundaries; and use native
		candidates as in \cref{thm:complete-contract-logic}. The universally
		valid implications among $\mathsf C,\mathsf S,\mathsf F,\mathsf E$ are
		generated as follows.
		\begin{enumerate}[label=\textup{(\arabic*)},nosep]
			\item \emph{Empty cut family.} The generators are $\mathsf S\Rightarrow\mathsf C$, $\mathsf F\Rightarrow\mathsf C$
			and $\top\Rightarrow\mathsf E$, if $\mathcal K=\varnothing$.
			\item \emph{Tree with all fundamental cuts.} The generators are $\mathsf S\Rightarrow\mathsf C$, $\mathsf F\Rightarrow\mathsf C$
			and $\mathsf S\wedge\mathsf E\Rightarrow\mathsf F$, if $G$ is a tree
			and $\mathcal K$ contains the fundamental cut obtained by deleting
			each edge of $G$.
			\item \emph{Remaining cut families.} The generators are $\mathsf S\Rightarrow\mathsf C$ and
			$\mathsf F\Rightarrow\mathsf C$ in all other cases.
		\end{enumerate}
		In each case a single $\mathbb F_q$-module atlas on $G$ realizes all
		$5$, $9$ or $10$ permitted truth assignments, respectively.
	\end{theorem}
	\begin{proof}
		$\mathsf S\Rightarrow\mathsf C$ holds as in
		\cref{thm:complete-contract-logic}. Since $W$ is state-exact and every
		local rule holds on it, $c_0=c_W$ has the prescribed local traces,
		which gives $\mathsf F\Rightarrow\mathsf C$. In case (2),
		\cref{thm:native-source-rigidity} gives
		$\mathsf S\wedge\mathsf E\Rightarrow\psi(H)=W\Rightarrow\mathsf F$. In
		case (1), $\mathsf E$ is vacuous. Any additional cuts in case (2)
		also hold on the full matching source: after fixing all crossing interfaces,
		compatible tuples on opposite sides glue independently
		(\cref{thm:arbitrary-graph-cut-information}).
		
		For realization we use products of atlases on $G$ with disjoint
		descriptor families. As in the proof of
		\cref{thm:complete-contract-logic}, each clause holds for a product
		source exactly when it holds in every factor. Three factors are used.
		
		\emph{Tree factor.} Choose an edge $e=ab$ of $G$, crossing some cut of
		$\mathcal K$ when $\mathcal K\ne\varnothing$. Place the two-chart atlas
		of \cref{thm:complete-contract-logic} on $a$ and $b$, with trivial
		states elsewhere and trivial interfaces off $e$. At a cut of
		$\mathcal K$ crossed by $e$, the clause $\mathsf E$ is the two-chart
		separator clause; at any other cut the factor lies on one side and
		$\mathsf E$ holds. The nine sources of that theorem give its nine rows. Here the shared
		descriptor of the factor is read on the edge $e$; all other edge
		interfaces are one-element modules. Its local closure is unchanged by
		adjoining the trivial charts. The choice of $e$ is possible because a
		nontrivial cut in a connected graph has a crossing edge.
		
		\emph{Coupling factor.} Suppose $G$ is a tree and the fundamental cut of
		an edge $e'$ is absent from $\mathcal K$. Give every vertex in one
		component of $G-e'$ the scalar descriptor $u$ on its chart
		$\mathbb F_q$, and every vertex in the other component a scalar $v$.
		Every edge inside a component uses identity interface maps. Use the
		trivial interface at $e'$. The matching space is $\mathbb F_q^2$
		with independent descriptors $u,v$. Its diagonal $D=\{u=v\}$ covers
		every chart and introduces the global implication $u\to v$.
		Every protected cut crosses an edge distinct from $e'$: the only
		nontrivial bipartition crossing $e'$ alone is its fundamental cut.
		The boundary of each protected cut therefore reads $u$ or $v$, which
		determines the entire diagonal source. Its cut identity follows.
		Thus $D$ satisfies $\mathsf C,\mathsf S,\mathsf E$ and fails $\mathsf F$
		even for protected cuts with disconnected sides.
		
		\emph{Parity factor.} If $G$ contains a cycle, place the atlas of
		\cref{prop:parity-cycle} on it, with trivial states off the cycle
		and trivial interfaces on all other edges. The argument of that
		proposition applies whenever a bipartition of $G$ splits the cycle.
		If the cycle lies on one side, the other side has only trivial states
		in this factor and its cut identity also holds. Thus $Z$ satisfies $\mathsf C$, $\mathsf S$ and $\mathsf E$ but not
		$\mathsf F$.
		\begin{enumerate}[label=\textup{(\arabic*)},leftmargin=*]
			\item \emph{Empty cut family.}
			The tree factor alone gives the five assignments, all with
			$\mathsf E=1$.
			\item \emph{Tree with all fundamental cuts.}
			The tree factor gives the nine rows.
			\item \emph{Remaining cut families.}
			Here $\mathcal K\ne\varnothing$, and either $G$ is a tree with an edge cut
			outside $\mathcal K$ or $G$ contains a cycle. Multiply the tree factor
			by the coupling or the parity factor. The products of the nine tree
			sources with the full source of the second factor give the nine rows.
			The product of the full tree source with $D$ or $Z$ gives
			$(\mathsf C,\mathsf S,\mathsf F,\mathsf E)=(1,1,0,1)$. These are the
			ten assignments allowed by $\mathsf S\Rightarrow\mathsf C$ and
			$\mathsf F\Rightarrow\mathsf C$.\qedhere
		\end{enumerate}
	\end{proof}
	
	The next result computes the universal logic when the native language
	has no operations. This case permits arbitrary finite sets, functions,
	and joint subrelations.
	
	\begin{corollary}[Universal contract logic for finite sets]
		\label{cor:unrestricted-contract-logic}
		Fix a graph and cut family as in \cref{thm:contract-logic-general}.
		\begin{enumerate}[label=\textup{(\roman*)},leftmargin=*]
			\item \emph{Universal implications.}
			For finite atlases in the empty signature with state-exact $W$ and
			$c_W=c_0$, the universal implications are generated by
			$\mathsf S\Rightarrow\mathsf C$ and $\mathsf F\Rightarrow\mathsf C$,
			together with $\top\Rightarrow\mathsf E$ when $\mathcal K$ is empty.
			\item \emph{Realization and sharpness.}
			One fixed set atlas on the prescribed graph realizes all five permitted
			assignments for an empty family, or all ten for a nonempty family.
			The same universal implications therefore hold over the class of all
			finite equational atlases with the stated matching and observation
			hypotheses when the native signature is allowed to vary.
		\end{enumerate}
	\end{corollary}
	\begin{proof}
		\leavevmode
		\begin{enumerate}[label=\textup{(\roman*)},leftmargin=*]
			\item \emph{Universal implications.}
			State coverage lifts local counterexamples and gives $\mathsf S\Rightarrow\mathsf C$. Since $W$ is state-exact and $c_W=c_0$, dependency faithfulness also gives the prescribed local closures. With no protected cuts, $\mathsf E$ holds vacuously.
			\item \emph{Realization and sharpness.}
			The module examples of \cref{thm:contract-logic-general}
			remain valid after forgetting the operations, with the same charts,
			readouts, source subsets, and closures. This proves the claim except on
			a tree when $\mathcal K$ contains every fundamental cut.
			
			For that case choose an edge $ab$. At $a$ take three states with binary
			readouts $(x,y)\in\{(0,0),(1,0),(1,1)\}$; at $b$ take one binary readout
			$z$ with both values. All other charts and all interfaces are singletons.
			The full product is state-exact. Neither $x$ nor $y$ determines the other,
			and $z$ varies independently, so $c_W=c_0$ is the free descriptor closure.
			Use the four triples
			\[
			Z=\{(0,0,0),(0,0,1),(1,0,0),(1,1,1)\}.
			\]
			Both charts are covered. At a cut separating $a,b$, the incidence graph
			of the two actual side images is a connected path on five vertices, so
			its common functions are constant by \cref{thm:proper-cut-incidence}.
			At every other cut one side is a singleton. Thus every cut is exact.
			On $Z$, however, $(x,z)$ determines $y$, which fails on $W$.
			
			Take the product of this atlas and the two-chart factor used in the
			preceding theorem, with disjoint descriptor names. Combining the full
			source of this new factor with each of the nine earlier sources gives
			the nine original assignments. Combining $Z$ with the full source of the
			earlier factor adds $(\mathsf C,\mathsf S,\mathsf F,\mathsf E)=(1,1,0,1)$.
			To check the cut clause for these products of set atlases, observe that
			two edges of a product incidence graph are connected precisely when
			their components in both factor graphs are connected. Indeed, a walk
			can change the edge in one factor at a time while keeping the other
			factor edge fixed, and projection preserves connectivity. Boundary
			fibers are likewise products. Thus the cut test is the conjunction
			of the factor tests. All ten assignments occur on this single product atlas. Since the empty-signature
			examples are themselves equational examples, no further implication holds
			uniformly over the larger class.\qedhere
		\end{enumerate}
	\end{proof}
	
	\subsection{Faithful descent and source counts on a tree}
	\label{subsec:uniform-descent}
	
	\Cref{thm:main-native-descent} combines the recovery and full-cut
	criteria for complete separators. To compute its minimum, we first
	identify when the tree count factors into the quotient of the chart
	and separator sizes.
	
	\subsubsection{Uniform and nonuniform separator fibers}
	
	The fiber sizes in \eqref{eq:local-pair-count-message}--\eqref{eq:local-pair-source-size} can depend on the separator value even for Mal'cev algebras.
	
	By \cref{prop:local-pair-source-cost}\textup{(ii)}, constant fibers at the child end of each edge suffice for the product formula. This holds for group and module homomorphisms, but can fail for other Mal'cev algebras.
	
	\begin{example}[Three-element nonuniform Mal'cev interface]\label{ex:nonuniform-join-count}
		Let $A=\{0,1,2\}$ and $q:A\twoheadrightarrow\ff_2$ send $0\mapsto0$ and $1,2\mapsto1$.
		Define a ternary operation
		\[
		m(a,b,c)=
		\begin{cases}
			a,& b=c,\\
			c,& a=b,\\
			s\bigl(q(a)+q(b)+q(c)\bigr),&\text{otherwise},
		\end{cases}
		\qquad s(0)=0,\quad s(1)=1.
		\]
		The first two cases give $m(a,b,b)=a$ and $m(a,a,c)=c$;
		they agree when all three entries coincide. In all three cases
		$q(m(a,b,c))=q(a)-q(b)+q(c)$ in $\ff_2$, so $q$ is a homomorphism.
		Its congruence classes have sizes $1$ and $2$.
		Consequently
		\[
		\bigl|A\times_{\ff_2}A\bigr|=1^2+2^2=5,
		\qquad \frac{|A|^2}{|\ff_2|}=\frac92.
		\]
		The recursion \eqref{eq:local-pair-count-message}--\eqref{eq:local-pair-source-size} gives the
		integer count for both uniform and nonuniform fibers.
		
	\end{example}
	
	\subsubsection{Descent at complete separators}
	
	When every separator is congruence-complete, faithful descent and full-cut rigidity determine the minimum over the fixed charts.
	
	\begin{corollary}[Faithful native minima at congruence-complete separators]\label{thm:main-native-descent}
		Let a finite join tree of finite algebras in a Mal'cev variety have
		the marked presentation law, jointly separating homomorphic
		readouts, and equal actual separator images. Assume every separator is
		congruence-complete, and let $W$ be the matching algebra.
		\begin{enumerate}[label=\textup{(\roman*)},leftmargin=*]
			\item \emph{Faithful minimum.}
			The source $W$ satisfies the full faithful native contract. Every state- and full-cut-exact compatible native source maps onto $W$, so $W$ is the unique minimum over the prescribed chart maps.
			\item \emph{Scalar linear charts.}
			For scalar linear readouts over a finite field, every compatible linear source that is recovery-faithful maps onto $W$, and equality of source cardinalities holds exactly when that map is an isomorphism.
		\end{enumerate}
	\end{corollary}
	\begin{proof}
		\leavevmode
		\begin{enumerate}[label=\textup{(\roman*)},leftmargin=*]
			\item \emph{Faithful minimum.}
			\Cref{thm:complete-port-native-descent} gives the faithful native realization. Apply \cref{thm:native-source-rigidity} for the minimum and uniqueness over the chart maps.
			\item \emph{Scalar linear charts.}
			Part~\textup{(i)} gives $c_W=c_0$. Apply \cref{cor:scalar-contract-equivalence} for surjectivity. A surjection between finite sets of equal cardinality is bijective.\qedhere
		\end{enumerate}
	\end{proof}
	
	\subsection{One cyclic atlas and three distinct sharp source costs}
	\label{subsec:cyclic-three-source-costs}
	
	The candidate class changes the optimum even on the smallest cyclic coordinate cover.
	Fix binary linear charts with descriptor sets $\{x_1,x_2\}$, $\{x_2,x_3\}$, and
	$\{x_3,x_1\}$, each having the full actual state space $\mathbb F_2^2$.
	Their matching source is $W=\mathbb F_2^3$, their generated closure is free, and the cover has no join tree.
	
	\begin{proposition}[All native and fixed-table minima for the binary triangle]
		\label{prop:cyclic-triangle-source-minima}
		For this atlas the following costs and equality descriptions hold.
		\begin{enumerate}[label=(\roman*)]
			\item \emph{Linear native sources.} Among native $\mathbb F_2$-linear sources with compatible chart
			maps, state exactness, and faithful recovery, the minimum is $8$.
			Every minimum is $W$ up to its unique isomorphism over the charts.
			
			\item \emph{Fixed binary charts.} Among subrelations $R\subseteq\{0,1\}^3$ with the same three
			actual two-coordinate chart images and faithful recovery, the minimum
			is $5$. Its eight equality cases are exactly one complete parity
			class together with one word of the opposite parity.
			
			\item \emph{Arbitrary alphabets.} Among arbitrary-alphabet representations of the free closure on
			three descriptors, with no prescribed chart images, the minimum is $4$.
		\end{enumerate}
	\end{proposition}
	
	\begin{proof}
		\leavevmode
		\begin{enumerate}[label=\textup{(\roman*)},leftmargin=*]
			\item \emph{Linear native sources.}
			The coordinate readouts are onto simple $\mathbb F_2$-modules, and
			$W$ is state-exact and faithful.
			Thus \textup{(i)} follows from \cref{thm:simple-diagram-source-rigidity}.
			\item \emph{Fixed binary charts.}
			Four or fewer rows cannot suffice. Every prescribed
			pair image has four values, so at least four rows are needed; with
			exactly four, any two-coordinate projection is bijective and recovers
			the third descriptor, contrary to free closure.
			
			For a five-row subset $R$, put $U=\{0,1\}^3\setminus R$.
			All three pair projections are onto exactly when $U$ contains no edge
			of the cube, since each omitted pair value corresponds to an entire
			two-word coordinate fiber. Every independent three-word subset of the
			cube lies in one parity class. Indeed, an even and an odd word are
			nonadjacent only when antipodal; no third word is nonadjacent to both.
			Thus $R$ must consist of one parity class and one extra word, giving
			$2\cdot4=8$ possibilities. Conversely these sets have the required pair
			images. Any five-word subset contains an edge in each coordinate direction.
			The four edges of that direction partition the eight words.
			Choosing at most one word from each edge gives at most four words.
			Those three kinds of pairs refute every nontrivial determination, so
			all eight sets are faithful.
			\item \emph{Arbitrary alphabets.}
			Free recovery requires, for each descriptor, a row
			pair differing only in that descriptor. Choose one such pair for each
			of the three descriptors. In any cycle of the chosen edges, each coordinate
			changed along the cycle would change exactly once. It could not then
			return to its initial value. The three distinct edges therefore form a
			forest on at least four vertices. The rows $000,100,010,001$ attain four.
			The binary case uses the same witness-forest argument
			\cite[ECCC revision~1, Claim~5.13, p.~25]{ByramjiImpagliazzo2025}.\qedhere
		\end{enumerate}
	\end{proof}
	
	Each five-row minimum is a relation with a commuting spectral
	realization. It fails the linear subspace requirement.
	The four-row construction preserves the same closure while omitting
	a value from each pair chart. These examples show that the same
	closure allows different chart images and overlap diagrams.
	
	\subsection{Native group and relational source minima}
	The one-cut spanning-tree minimum gives a comparison with the
	minimum under a fixed group formation law on the same local states.
	
	\begin{proposition}[Sharp group and relational source costs for the two $S_3$ charts]
		\label{prop:native-group-versus-spectral-cost}
		Let two charts have local groups $S_3$ and descriptors
		$\{a,s\}$ and $\{s,b\}$, where $a,b$ are the identity readouts
		of their respective $S_3$ states and $s$ is the sign homomorphism
		to $C_2$. Retain these actual local groups and actual readouts.
		Their full matching source is
		\[
		W=S_3\times_{C_2}S_3
		=\{(g,h):\operatorname{sgn}(g)=\operatorname{sgn}(h)\}.
		\]
		The following two minima are exact.
		\begin{enumerate}[label=(\roman*)]
			\item \emph{Group sources.} Among subgroup subamalgams $H\leq W$ whose projections onto
			both local groups are surjective and whose descriptor closure is
			exactly generated by $a\to s$ and $b\to s$, the minimum is
			$|H|=18$, attained only by $H=W$.
			The same lower bound of eighteen holds for any finite group that maps onto both
			local $S_3$ groups, has equal sign observations, and has the required
			descriptor closure. At eighteen elements its joint observation map
			is an isomorphism onto $W$ over both charts.
			\item \emph{Spectral state relations.} Among arbitrary subrelations $R\subseteq W$ preserving all
			local states, the same descriptor closure, and the exact separator
			identity
			\[
			\mathbb C[a]\cap\mathbb C[b]=\mathbb C[s]
			\quad\text{inside }\mathbb C^R,
			\]
			the minimum is $|R|=10$.
		\end{enumerate}
		The first minimum is a native group-formation requirement; the
		second permits finite spectral state relations without group
		closure. \end{proposition}
	\begin{proof}
		\leavevmode
		\begin{enumerate}[label=\textup{(\roman*)},leftmargin=*]
			\item \emph{Group sources.}
			Write $S_3=C_3\rtimes C_2$ with elements $(u,\epsilon)$,
			$u\in\mathbb F_3$, $\epsilon\in\mathbb F_2$, and multiplication
			\[
			(u,\epsilon)(v,\delta)
			=(u+(-1)^\epsilon v,\epsilon+\delta).
			\]
			Then
			\[
			W=\{(u,v,\epsilon):u,v\in\mathbb F_3,
			\epsilon\in\mathbb F_2\}
			\cong(C_3\times C_3)\rtimes C_2,
			\]
			where $C_2$ acts by simultaneous inversion. This is a nonabelian
			group of order $18$. Its binary singleton separator satisfies
			relative pair lifting, so the native construction is recovery
			faithful by \cref{thm:binary-separator-threshold} and
			\cref{thm:equational-native-amalgamation}.
			
			If $H\leq W$ projects onto a local group of order six, its order
			is divisible by six and divides eighteen; hence it is six or
			eighteen. At order six both projections are bijective, creating
			$a\to b$ and $b\to a$, which are absent from the generated local
			closure. Thus a faithful such $H$ has order eighteen and equals
			$W$.
			
			For the more general assertion, the joint readout of the finite
			group is a homomorphism into $W$, with a subgroup image $H$ which
			projects onto both local groups. A descriptor implication holds
			on the original group exactly when it holds on this image, since
			every image state has a preimage and all descriptors factor through
			the joint readout. Its faithful image is therefore $W$ by the
			preceding argument. Surjectivity onto $W$ gives the cardinality
			bound, and equality makes the joint homomorphism an isomorphism.
			\item \emph{Spectral state relations.}
			For a state-exact subrelation $R$, each sign fiber is a bipartite
			graph with three vertices on each side, all six vertices used.
			By \cref{thm:proper-cut-incidence}, exact separator intersection
			is equivalent to connectedness of each sign fiber. Each then needs at least five
			edges, giving $|R|\geq10$.
			
			For attainment, take
			\begin{equation}\label{eq:ten-state-spectral-group-comparison}
				R=\{(u,0,\epsilon):u\in\mathbb F_3,
				\epsilon\in\mathbb F_2\}
				\ \cup\
				\{(0,v,\epsilon):v\in\mathbb F_3,
				\epsilon\in\mathbb F_2\}.
			\end{equation}
			There are five states per sign, since the two displayed axes meet
			in $(0,0,\epsilon)$. Both local projections are onto, and each
			sign graph is a connected double star. The two local implications
			hold. At fixed $a=(0,\epsilon)$ all three $b$-values occur.
			The symmetric case gives the same variation of $a$ with $b$ fixed.
			Both private observations vary within each sign, and both sign
			values occur. These tests exhaust the nontrivial nonimplications
			on $\{a,s,b\}$, proving exact recovery faithfulness. The relation fails closure under multiplication. The product of
			$(1,0,0)$ and $(0,1,0)$ is $(1,1,0)$, outside the relation.\qedhere
		\end{enumerate}
	\end{proof}
	
	The group source in this example also has the compact native datum
	\[
	\left\langle r_1,r_2,t\ \middle|\
	r_1^3=r_2^3=t^2=1,\quad [r_1,r_2]=1,\quad
	tr_1t=r_1^{-1},\quad tr_2t=r_2^{-1}\right\rangle.
	\]
	The relations reduce every element to $r_1^u r_2^v t^\epsilon$,
	so the quotient has at most eighteen elements; the displayed
	semidirect product realizes all eighteen normal forms. Its two
	source maps to the standard $S_3$ presentation send respectively
	$(r_1,r_2,t)$ to $(r,1,t)$ and $(1,r,t)$. Thus the example can be
	formed directly in the group-presentation law.
	
	\begin{theorem}[Sharp group and relational source minima with an unbounded cost ratio]
		\label{thm:unbounded-native-resource-separation}
		Let $\pi:G\twoheadrightarrow Q$ be a surjective homomorphism of finite
		groups with $m=|\ker\pi|\ge2$, and put $W=G\times_QG$. Retain the
		whole-chart readouts $a,b$ and the common readout
		$s=\pi a=\pi b$. Require both chart images to be $G$, the determination
		closure of these three readouts to equal that on $W$, and the exact cut
		identity
		\begin{equation}\label{eq:resource-separation-contract}
			\mathbb C[a]\cap\mathbb C[b]=\mathbb C[s]
		\end{equation}
		inside the function algebra of the candidate source. Compare the
		following two source classes: finite groups with the specified chart
		homomorphisms in $\Xi_{\mathrm{Grp}}^{\mathrm{pres}}$, and arbitrary
		subrelations of $W$. Their sharp minimum cardinalities are
		\begin{equation}\label{eq:unbounded-native-resource-minima}
			N_{\mathrm{grp}}=|Q|m^2,\qquad
			N_{\mathrm{rel}}=|Q|(2m-1).
		\end{equation}
		\begin{enumerate}[label=\textup{(\roman*)},leftmargin=*]
			\item \emph{Minimum group sources.}
			Every minimum group source has its joint chart map as its unique
			isomorphism onto $W$ over the charts.
			\item \emph{Minimum relations.}
			Every minimum relation consists,
			in each $Q$-fiber, of a spanning tree on the two prescribed sets of
			$m$ chart states. There are exactly
			\begin{equation}\label{eq:unbounded-native-resource-equality-count}
				m^{\,2|Q|(m-1)}
			\end{equation}
			such minimum relations over the fixed state labels.
			\item \emph{Ratio of costs.}
			In particular, the exact native-to-relational cost ratio is
			\begin{equation}\label{eq:unbounded-native-resource-ratio}
				\frac{N_{\mathrm{grp}}}{N_{\mathrm{rel}}}
				=\frac{m^2}{2m-1}.
			\end{equation}
			This ratio is unbounded even with $Q$ fixed and the group-presentation law
			unchanged: take $G=Q\times C_m$ with its projection onto $Q$.
		\end{enumerate}
	\end{theorem}
	\begin{proof}
		\leavevmode
		\begin{enumerate}[label=\textup{(\roman*)},leftmargin=*]
			\item \emph{Minimum group sources.}
			Let $Z\le W$ be the joint image of a finite group source.
			Surjective pullback from $Z$ preserves the descriptor dependencies and
			the intersection in \eqref{eq:resource-separation-contract}.
			\Cref{thm:native-source-rigidity} gives $Z=W$. Hence $|H|\ge |Q|m^2$,
			with equality precisely when the joint map is an isomorphism over the
			two charts. The matching group attains this bound.
			\item \emph{Minimum relations.}
			For an arbitrary relation, each connected fiber has at least $2m-1$
			edges, with equality precisely for a spanning tree. Choosing one such
			tree independently in each fiber gives both chart surjectivities and
			the exact cut identity. These relations also have the required
			determination closure. The implications $a\to s$ and $b\to s$ hold.
			At fixed $s$, all $m$ values of each private chart state occur.
			Since $2m-1>m$, a fiber tree has a vertex of degree at least two
			on each side. Thus $a$ varies for some fixed value of $b$, and $b$
			varies for some fixed value of $a$. The only constant readout among
			these three is $s$ when $|Q|=1$, exactly as on $W$. These observations
			determine the full closure on $\{a,b,s\}$. Therefore every fiberwise
			spanning-tree relation is feasible, and every minimum has this form.
			
			The matrix--tree count for the labeled complete bipartite graph
			$K_{m,m}$ is $m^{2m-2}$, as in
			\cref{thm:one-cut-spanning-minimum}. Multiplication over the $|Q|$
			independent fibers gives \eqref{eq:unbounded-native-resource-equality-count}.
			Each finite relation also has the diagonal spectral realization of
			\cref{prop:spectral-observation-interpretation}.
			\item \emph{Ratio of costs.}
			Finally
			\eqref{eq:unbounded-native-resource-ratio} follows from the two minima,
			and $G=Q\times C_m$ realizes every integer $m\ge2$ with the same
			separator group. Its ratio is asymptotic to $m/2$.\qedhere
		\end{enumerate}
	\end{proof}
	
	\section{Terminal-tree representations}
	\label{sec:terminal-representations}
	
	At an internal vertex of a terminal tree, let all incident edge
	descriptors except one determine the remaining descriptor. We ask
	which global dependencies follow and how many rows realize exactly
	those dependencies. \Cref{thm:tree-rigidity} proves that a tree with
	$n$ terminals needs $n$ rows and determines the equality partition
	of every column at the minimum. Local realizations can be smaller:
	three rows suffice for a trivalent tree. We also count the dependencies
	lost at smaller sizes and construct a faithful source with the local
	spectra fixed.
	
	\subsection{The closure generated by terminal-tree rules}
	\label{subsec:terminal-tree-closure-minima}
	
	Terminal paths characterize nonimplications, and boundaries of connected
	sets of internal vertices give the minimal implications. These two
	descriptions will supply the witnesses used to compute source sizes.
	
	\subsubsection{Chart rules and representation minima}
	
	Let $T=(W,E)$ be a finite tree with terminal set $\Omega$ equal to its degree-one
	vertices, $n=|\Omega|\geq3$, and all other vertices of degree at least three.
	Edge labels are the descriptors.
	Write $J=W\setminus	\Omega$ and $\delta_T(U)$ for the edges with exactly one endpoint in $U$.
	For each $v\in J$, put $I_v=\delta_T(v)$ and impose all rules
	\begin{equation}\label{eq:tree-local-rules}
		I_v\setminus\{e\}\longrightarrow e\qquad(e\in I_v).
	\end{equation}
	Let $c_T$ be their least consequence closure.
	The local closure on $I_v$ is the circuit closure in which any $|I_v|-1$ descriptors
	determine the remaining one, with no other nontrivial determinations.
	
	For a trivalent tree, choosing a terminal as the external root orients
	the internal vertices as a rooted full binary tree. The choice leaves
	the local rules unchanged.
	
	A finite relation $R$, with nonempty row set and one column $f_e$ per edge,
	induces the determination closure
	\[
	c_R(S)=\{e\in E:\ \forall x,y\in R,\quad
	(f_s(x)=f_s(y)\text{ for all }s\in S)
	\Rightarrow f_e(x)=f_e(y)\}.
	\]
	Values in different columns may belong to unrelated sets. Repeated full rows are discarded.
	We call $R$	\emph{sound} if $c_T(S)\subseteq c_R(S)$ for every $S\subseteq E$;
	\emph{chart-exact} if $c_R|_{I_v}=c_T|_{I_v}$ for every $v\in J$; and \emph{faithful} if $c_R=c_T$.
	Chart-exactness implies soundness. A chart-exact relation can still
	introduce extra dependencies. The objectives prescribe closures
	and allow arbitrary numerical values for the observable operators.
	Here restriction of a closure to $I$ means $c|_I(S)=c(S)\cap I$ for $S\subseteq I$.
	Define
	\[
	d_{\rm loc}(T)=\min_{R\text{ chart-exact}}|R|,
	\qquad
	d_{\rm faith}(T)=\min_{R\text{ faithful}}|R|.
	\]
	For a closure $c$, an access state is a fixed point $C=c(C)$.
	The signs denote membership tests. The requirement $e^+$ means
	$e\in C$, and $e^-$ means $e\notin C$. The signed Helly number is the
	largest size of an inclusion-minimal infeasible family of these tests.
	Use the value one when no larger family occurs. This number measures
	feasibility of tests on selected closed sets in a given closure.
	
	Both costs count distinct rows, or equivalently distinct joint
	spectral sectors. Encoding a column by a diagonal self-adjoint
	operator gives the corresponding commuting realization, as in
	\cite[Theorem~5.7]{Ma2026}.
	
	\subsubsection{Terminal paths, closed sets, and complete obstruction geometry}
	
	For distinct terminals $a,b$, let $P_{ab}\subseteq E$ be their unique path and
	put $C_{ab}=E\setminus P_{ab}$. A \emph{coatom} is an inclusion-maximal proper
	closed set. The following elementary tree description is also a specialization
	of the stopping-set/coatom correspondence for hypergraph Horn functions
	\cite{BercziHH2024,Zhu2026}.
	
	\begin{lemma}[Terminal paths, recovery closure, and coatoms]\label{lem:tree-paths}
		For $F\subseteq E$, the following are equivalent:
		\begin{enumerate}[label=(\roman*),nosep]
			\item $c_T(F)=F$;
			\item in the forest $T\setminus F$, no internal vertex has degree one;
			\item every edge of $T\setminus F$ lies on a terminal-to-terminal path in that
			forest.
		\end{enumerate}
		Consequently every closed set is an intersection of the $C_{ab}$, and the
		$\binom n2$ distinct sets $C_{ab}$ are exactly the coatoms and exactly the
		proper meet-irreducible closed sets. Moreover
		\begin{equation}\label{eq:tree-path-test}
			e\notin c_T(S)
			\quad\Longleftrightarrow\quad
			\exists a\ne b:\ e\in P_{ab},\quad P_{ab}\cap S=\varnothing.
		\end{equation}
		The binary edge-cut descriptors on $\Omega$ realize $c_T$ faithfully.
		The identification of the coatom complements with terminal paths also
		follows directly from Zhu's frequency-two apex-graph
		result \cite[Proposition~5.6]{Zhu2026}: its graphic circuits are
		the terminal-to-terminal paths after identifying the terminals.
	\end{lemma}
	\begin{proof}
		\leavevmode
		\begin{itemize}
			\item \emph{\textup{(i)} $\Longleftrightarrow$ \textup{(ii)}.}
			All the local rules hold on $F$ exactly when each internal vertex
			has either zero or at least two remaining edges. This proves (i)$\Leftrightarrow$(ii).
			\item \emph{\textup{(ii)} $\Longleftrightarrow$ \textup{(iii)}.}
			In a finite tree satisfying (ii), each side of a remaining edge contains a
			leaf of its component, necessarily an original terminal. This proves (iii).
			Conversely a terminal path through an edge supplies a second remaining edge
			at each internal endpoint, proving (ii).
			\item \emph{Coatoms and the path test.}
			The complement of any terminal path satisfies (ii). The full path is the only nonempty subset of its edges satisfying
			(ii). Any nonempty proper subset leaves a nonterminal end. Thus $C_{ab}$ is a coatom.
			By (iii), the complement of any closed set is a union of terminal paths,
			so that closed set is an intersection of their complements. In a finite
			intersection system generated by its coatoms, every proper meet-irreducible
			is a coatom. Distinct terminal pairs have distinct paths. The same
			intersection description, applied to $c_T(S)$, proves \eqref{eq:tree-path-test}.
			For an edge-cut descriptor, two terminal rows disagree precisely on their
			connecting path. Hence its determination test is \eqref{eq:tree-path-test}.\qedhere
		\end{itemize}
	\end{proof}
	
	\begin{theorem}[Minimal recovery premises and signed obstructions on terminal trees]
		\label{thm:tree-boundary-cores}
		\leavevmode
		\begin{enumerate}[label=\textup{(\roman*)},leftmargin=*]
			\item \emph{Minimal premises.}
			For $e\notin P\subseteq E$, the set $P$ is an inclusion-minimal premise
			forcing $e$ if and only if
			\begin{equation}\label{eq:tree-core-boundary}
				P\cup\{e\}=\delta_T(U)
			\end{equation}
			for a nonempty connected set $U\subseteq J$.
			\item \emph{Signed obstructions.}
			The nontrivial minimal
			signed obstructions at the constant interface are exactly
			\[
			\bigl(\delta_T(U)\setminus\{e\}\bigr)^+\cup\{e^-\},
			\quad \varnothing\ne U\subseteq J\text{ connected},\quad e\in\delta_T(U).
			\]
			The additional minimal signed obstructions are the direct contradictions
			$\{e^+,e^-\}$.
			\item \emph{Boundary sizes and Helly number.}
			The boundary sizes and signed Helly number are
			\begin{equation}\label{eq:tree-core-sizes}
				|\delta_T(U)|=2+\sum_{v\in U}(\deg_T(v)-2),\qquad
				\hr(q_\varnothing;\mathcal A_\pm)=n.
			\end{equation}
			Here $q_\varnothing$ is the constant map on $\operatorname{Fix}(c_T)$ and
			$\mathcal A_\pm$ consists of the positive and negative membership requirements.
			The signed Helly number is the largest size of an inclusion-minimal
			infeasible signed family, with a minimum value of one by convention.
		\end{enumerate}
	\end{theorem}
	
	\begin{proof}
		\leavevmode
		\begin{enumerate}[label=\textup{(\roman*)},leftmargin=*]
			\item \emph{Minimal premises.}
			If $e\in c_T(P)$, no terminal path through $e$ avoids $P$. In the component
			of $T\setminus P$ containing $e$, at least one side of $e$ has no terminal.
			Let $U$ be the vertices on that side. Then $U$ is nonempty, connected, and
			internal, and $\delta_T(U)\subseteq P\cup\{e\}$. Deleting the boundary
			except $e$ leaves that same terminal-free side, so its boundary except $e$
			already forces $e$. Minimality of $P$ gives equality in
			\eqref{eq:tree-core-boundary}.
			
			Conversely, each component attached to a connected internal $U$ through a
			boundary edge contains a terminal. To see this, it is a finite tree and its
			only possible missing incidence is its attachment to $U$; an original
			internal vertex cannot be its terminal-free end. For $f\in\delta_T(U)-e$,
			choose terminals beyond $e$ and $f$. Their path crosses exactly $e,f$ among
			the boundary edges. It witnesses failure of determination after deleting $f$
			from the premise. Deleting all other boundary edges does force $e$, as the
			$U$-side contains no terminal. This proves minimality.
			\item \emph{Signed obstructions.}
			All-positive signed
			families are feasible, and a signed family with positive support $P$ is feasible exactly when
			$c_T(P)$ avoids every negatively requested label: $c_T(P)$ is the least
			possible state. An infeasible family therefore has an infeasible subfamily
			with just one negative. Minimality gives the stated classification.
			If its target is already one of its positives, it is the direct contradiction.
			\item \emph{Boundary sizes and Helly number.}
			Counting degrees in the connected tree on $U$ proves the size formula.
			All summands are positive, so the largest boundary is $\delta_T(J)$, the
			$n$ terminal edges. Its rooted obstruction has size $n$, and no classified
			obstruction is larger.\qedhere
		\end{enumerate}
	\end{proof}
	
	\begin{remark}[Matroid interpretation and the classical presentation optimum]
		\label{rem:tree-horn-minimum}
		Identify all terminals of $T$ to one vertex to obtain a multigraph $G_T$.
		The closure $c_T$ is the closure of the cographic matroid of $G_T$.
		Indeed, an edge outside $S$ belongs to its cographic closure precisely when
		it is a bridge in $G_T\setminus S$, equivalently when no terminal path in
		$T\setminus S$ crosses it. Its circuits are the bonds
		$\delta_T(U)$ above. This cographic matroid is simple and binary, and its
		chordless circuits are exactly the internal stars $I_v$, by the
		boundary argument in the proof of \cref{cor:tree-reconstruct}. Thus the local star presentation is already the unique minimum
		circuit generator and circuit-clause presentation supplied by
		\cite[Theorems~16--17, arXiv v2]{BercziMatroid2024}.
	\end{remark}
	
	\subsection{Rigidity and fidelity of terminal-tree representations}
	\label{subsec:terminal-tree-source-rigidity}
	
	For arbitrary finite column alphabets, \cref{thm:tree-rigidity}
	classifies all minimum faithful row partitions and identifies the tree
	data they determine. At smaller sizes, the number of failed coatom
	tests measures the loss of recovery. Chartwise exactness alone may
	permit fewer rows.
	
	\subsubsection{Rigidity of minimum representations over arbitrary alphabets}
	
	For two distinct rows $x,y$ put
	\[
	\operatorname{agr}(x,y)=\{e:f_e(x)=f_e(y)\},\qquad
	\Delta(x,y)=E\setminus\operatorname{agr}(x,y).
	\]
	The use of agreement sets to characterize Armstrong relations and to bound their minimum
	size is classical \cite[Theorems~6.1 and~6.5]{BeeriDFS1984}.
	The terminal paths determine how the pendant and internal columns can
	partition the rows of a minimum representation.
	The classical uniform family in \cref{thm:uniform-nonrigid} shows why the equality
	classification cannot be inferred from the numerical bound alone.
	
	\paragraph{Reduction to the classical pendant-column problem.}
	Let $E_\partial=\{e_a:a\in\Omega\}$ be the $n$ pendant descriptors.
	The path test gives, for $S\subseteq E_\partial$,
	\[
	c_T(S)\cap E_\partial=
	\begin{cases}
		S,&|S|\leq n-2,\\
		E_\partial,&|S|\geq n-1.
	\end{cases}
	\]
	Indeed, two missing terminal edges are joined by a path avoiding $S$;
	when at most one is missing, every terminal path meets $S$.
	This is the uniform closure $U_{n-1,n}$.
	Its row minimum $n$, attained by the identity matrix, is an instance of the classical theory of minimum matrix representations
	\cite{DemetrovicsKatona1981,DemetrovicsKatona1982,DemetrovicsFurediKatona1985}.
	The pendant descriptors jointly recover every edge, so projection onto them preserves distinct rows of a faithful table.
	Consequently every $n$-row faithful tree table restricts to a minimum uniform table.
	The rigidity argument determines the minimum pendant partitions
	and their extensions to the internal edges.
	
	\begin{lemma}[Coatom witnesses and minimum faithful size for terminal trees]
		\label{lem:tree-saturation}
		Let $R$ be a sound relation.
		\begin{enumerate}[label=\textup{(\roman*)},leftmargin=*]
			\item \emph{Coatom witnesses.}
			Every agreement set is $c_T$-closed.
			The relation $R$ is faithful if and only if every $C_{ab}$ occurs as the agreement set of two distinct rows.
			\item \emph{Size and equality cases.}
			If $R$ is faithful, then $|R|\geq n$.
			At equality, every row pair has disagreement set equal to one terminal path, and the
			correspondence between row pairs and terminal pairs is a bijection.
		\end{enumerate}
	\end{lemma}
	
	\begin{proof}
		\leavevmode
		\begin{enumerate}[label=\textup{(\roman*)},leftmargin=*]
			\item \emph{Coatom witnesses.}
			If two rows agree on $S$, soundness makes them agree on $c_T(S)$, so their agreement set is closed.
			In a faithful relation, any $e\notin C_{ab}$ has a
			pair of rows agreeing on $C_{ab}$ but not on $e$.
			Its proper closed agreement	set contains the coatom and therefore equals it.
			Conversely, for every $e\notin c_T(S)$, \cref{lem:tree-paths} supplies a coatom containing $S$
			but not $e$;
			a pair witnessing that coatom refutes $S\to e$ in $R$.
			\item \emph{Size and equality cases.}
			For a faithful relation there are $\binom n2$ different coatoms and only $\binom{|R|}2$ row pairs.
			When the two numbers agree, every pair is used once and has the stated form.\qedhere
		\end{enumerate}
	\end{proof}

	\begin{theorem}[Minimum size and rigidity of faithful terminal-tree representations]
		\label{thm:tree-rigidity}
		\leavevmode
		\begin{enumerate}[label=\textup{(\roman*)},leftmargin=*]
			\item \emph{Minimum size.}
			For every terminal tree above,
			$
			d_{\rm faith}(T)=n.
			$
			
			\item \emph{Rigidity at equality.}
			If a faithful relation $R$ has exactly $n$ rows, there is a bijection $\theta:\Omega\to R$ such that for every edge $e$,
			\begin{equation}\label{eq:tree-rigidity}
				f_e(\theta(a))=f_e(\theta(b))
				\quad\Longleftrightarrow\quad
				a,b\text{ lie on the same side of }T\setminus e.
			\end{equation}
			Thus every minimum representation has binary equality partitions and is the terminal-cut representation up to row permutation and independent renaming of the two values in each column.
			Two-element column images are a	conclusion, with no restriction on unused symbols in the ambient alphabets.
			In the native realization, the same assertion concerns the represented commutative descriptor algebras, not unique numerical eigenvalues.
		\end{enumerate}
	\end{theorem}

	\begin{proof}
		\leavevmode
		\begin{enumerate}[label=\textup{(\roman*)},leftmargin=*]
			\item \emph{Minimum size.}
			The lower bound is \cref{lem:tree-saturation}; the terminal-cut construction attains it.
			
			\item \emph{Rigidity at equality.}
			Suppose equality holds.
			For a terminal $a$, let $e_a$ be its pendant edge.
			Exactly $n-1$ terminal paths use $e_a$.
			The preceding lemma assigns each terminal path to exactly one row pair.
			Thus the partition of the $n$ rows given by $f_{e_a}$ separates
			exactly $n-1$ unordered pairs.
			A nonconstant partition of $n\geq3$ points separates at least $n-1$ pairs, with
			equality only for a singleton and its complement.
			Two blocks of sizes $k,n-k$ separate $k(n-k)\geq n-1$,
			with equality only at $k=1,n-1$; three or more blocks separate at least
			$2n-3>n-1$ pairs.
			
			Let $r_a$ be the isolated row. Distinct terminals isolate distinct rows,
			since their pendant edges occur in different sets of terminal paths and
			their columns therefore separate different sets of row pairs.
			Hence $a\mapsto r_a$ is a bijection.
			Among the pendant-edge columns, the pair $r_a,r_b$ disagrees
			on precisely $e_a$ and $e_b$; it may also disagree on internal-edge
			columns.
			The disagreement set is a terminal path. Its two pendant edges
			identify its endpoints.
			Therefore $\Delta(r_a,r_b)=P_{ab}$.
			This is exactly \eqref{eq:tree-rigidity} for every edge, including internal edges.\qedhere
		\end{enumerate}
	\end{proof}
	
	\begin{corollary}[Reconstruction of the unrooted terminal tree from its recovery closure]
		\label{cor:tree-reconstruct}
		For every terminal tree with all internal degrees at least three,
		the closure $c_T$ alone determines its unrooted edge-labeled tree,
		up to terminal relabeling.
		\begin{enumerate}[label=\textup{(\roman*)},leftmargin=*]
			\item \emph{Intrinsic star test.}
			Its internal stars are exactly those circuits
			$C$ for which there is no circuit $D$ with
			\begin{equation}\label{eq:intrinsic-tree-star-test}
				|D\setminus C|=1,\qquad |D|<|C|.
			\end{equation}
			\item \emph{Reconstruction.}
			A bijection of descriptor labels preserves the closure
			if and only if it is induced by an isomorphism of these trees.
			Every optimal source yields the same tree from its equality partitions.
			For trivalent trees, testing three-element circuits suffices.
			Choosing a terminal as the root and orienting the tree from that
			root gives the same closure and minimum source.
		\end{enumerate}
	\end{corollary}
	\begin{proof}
		\leavevmode
		\begin{enumerate}[label=\textup{(\roman*)},leftmargin=*]
			\item \emph{Intrinsic star test.}
			By \cref{thm:tree-boundary-cores}, the circuits are the boundaries
			$\delta_T(U)$ of nonempty connected internal vertex sets. If $|U|>1$,
			choose a leaf $v$ of the subtree on $U$. Its star $D=I_v$ has exactly
			one edge outside $C=\delta_T(U)$, and
			$|C|=2+\sum_{u\in U}(\deg_T(u)-2)>\deg_T(v)=|D|$.
			Conversely, for $C=I_v$, a different boundary $D=\delta_T(U)$
			has at least two edges outside $C$. If $v\notin U$, at most one
			star edge meets $U$ and $|D|\geq3$. If $v\in U$ and $U\ne\{v\}$,
			a nonempty branch of $U\setminus\{v\}$ has at least two boundary
			edges away from $v$, by the same degree count. This proves the test.
			\item \emph{Reconstruction.}
			Reconstruct one internal vertex for each recovered star. A descriptor
			belonging to two stars joins their vertices; one belonging to just one
			star receives a terminal endpoint. These are all descriptors and give
			exactly $T$. This construction proves the closure-isomorphism assertion.
			For a minimum source, \cref{thm:tree-rigidity} gives precisely the edge
			splits; their reconstruction also follows from Buneman's splits-equivalence theorem
			\cite{Buneman1971}; see \cite[Theorem~2.1]{Huber2014}.
			When internal degrees are three, the boundary-size formula makes
			the three-element circuits exactly the stars. Rechoosing a root or
			formation direction changes none of these undirected data.\qedhere
		\end{enumerate}
	\end{proof}
	
	\begin{proposition}[Subdivisions: all minimum sources and the exact reconstruction ambiguity]
		\label{prop:subdivided-tree-reconstruction}
		Allow a finite terminal tree $T$ with $n\geq3$ terminals to have
		internal vertices of degree two, and impose the same local rules
		\eqref{eq:tree-local-rules}. Suppress all degree-two vertices to obtain
		$\overline T$, and let $q:E(T)\to E(\overline T)$ map each maximal
		chain to its reduced edge.
		\begin{enumerate}[label=\textup{(\roman*)},leftmargin=*]
			\item \emph{Closure and chain classes.}
			The closure satisfies
			\begin{equation}\label{eq:subdivision-closure}
				c_T(S)=q^{-1}\!\left(c_{\overline T}(q(S))\right).
			\end{equation}
			The fibers of $q$ are exactly the mutual-recovery classes
			$e\sim f$ defined by $e\in c_T(\{f\})$ and
			$f\in c_T(\{e\})$.
			\item \emph{Minimum sources.}
			The minimum arbitrary-alphabet faithful source
			still has $n$ rows; every minimum has the terminal-cut partitions
			of $T$, with identical equality partitions along each such chain.
			\item \emph{Reconstruction ambiguity.}
			The recovery closure determines the reduced tree with each reduced
			edge labeled by its whole descriptor class, hence also the length
			of every subdivided chain. Among edge-labeled terminal trees, the
			only remaining undirected ambiguity is the order of labels along
			these chains, modulo tree isomorphisms. Roots and formation directions
			remain additional data.
		\end{enumerate}
	\end{proposition}
	\begin{proof}
		\leavevmode
		\begin{enumerate}[label=\textup{(\roman*)},leftmargin=*]
			\item \emph{Closure and chain classes.}
			At a degree-two vertex the two incident descriptors recover each other.
			Saturating these rules first saturates exactly the maximal chains.
			The rules at the remaining internal vertices then become exactly the
			star rules of $\overline T$, proving \eqref{eq:subdivision-closure}.
			Its reduced closure has no distinct mutually recovering descriptors:
			all its circuits have at least three elements by the boundary-size
			formula. Thus its chain classes are recovered intrinsically.
			\item \emph{Minimum sources.}
			In any faithful relation, mutual recovery is equality of column
			partitions. Keep one descriptor from each class. This preserves
			distinct rows and gives a faithful representation of
			$c_{\overline T}$. Conversely any faithful reduced representation
			extends by repeating its equality partitions along the chains.
			\Cref{thm:tree-rigidity} therefore gives the minimum $n$ and all
			equality cases.
			\item \emph{Reconstruction ambiguity.}
			Apply \cref{cor:tree-reconstruct} to the induced
			closure on classes to recover the reduced tree. Formula
			\eqref{eq:subdivision-closure} is insensitive to the order within
			each class, and every tree with this class-labeled reduction is
			obtained by ordering those labels along its chains. This proves
			the exact ambiguity statement.\qedhere
		\end{enumerate}
	\end{proof}
	
	\subsubsection{All minimum chartwise completions and an exponential gap}
	
	On a trivalent tree, the chartwise and faithful optima give the exact
	additional cost of preserving dependencies across charts.
	
	\begin{theorem}[Chartwise and faithful source minima for trivalent terminal trees]
		\label{thm:tree-two-costs}
		Let every internal vertex of $T$ have degree three.
		\begin{enumerate}[label=\textup{(\roman*)},leftmargin=*]
			\item \emph{Minimum sizes.}
			The two minima are
			\begin{equation}\label{eq:tree-two-costs}
				d_{\rm loc}(T)=3,\qquad d_{\rm faith}(T)=n.
			\end{equation}
			\item \emph{Colorings and global closures.}
			The column partitions of a three-sector chart-exact representation
			give a proper three-edge-coloring of $T$. Two such representations
			have the same closure exactly when their colorings differ by a
			single permutation of the colors. Consequently there are exactly
			\begin{equation}\label{eq:tree-local-count}
				2^{n-3}
			\end{equation}
			distinct global closures arising from minimum chart-exact representations.
			\item \emph{Perfect binary trees.}
			For a perfect rooted binary tree of height $h\ge1$, adjoining an external
			root terminal gives $n=2^h+1$. The faithful cost grows
			exponentially with height and linearly with the number of terminals,
			while the chartwise cost stays three.
		\end{enumerate}
	\end{theorem}
	\begin{proof}
		\leavevmode
		\begin{enumerate}[label=\textup{(\roman*)},leftmargin=*]
			\item \emph{Minimum sizes.}
			On three rows, there are exactly three nonconstant, nondiscrete partitions:
			the three singleton-versus-complement partitions. A chart with the
			pair-to-third closure must use these three different partitions. A constant
			column would be determined by the empty set, a discrete column would
			determine the others, and identical partitions would determine each other;
			all are forbidden by the chart closure. A two-row relation has only one
			nonconstant partition and cannot realize a chart, proving the lower bound
			three. Assigning its three singleton partitions as colors shows that a
			three-row chart-exact representation is exactly a proper three-edge-coloring.
			Every trivalent tree has such a coloring: color a first internal star and
			then, on reaching a new internal vertex, assign its other two edges the two
			remaining colors. Thus the lower bound is attained.
			The faithful minimum follows from \cref{thm:tree-rigidity}.
			\item \emph{Colorings and global closures.}
			There are $n-2$ internal vertices. The first star has $6$ colorings and each
			of the remaining $n-3$ internal vertices has two choices. Global color
			permutation acts freely, so the number of colorings modulo permutation is
			$6\,2^{n-3}/6$. A singleton descriptor in a resulting global closure
			determines exactly its color class; two distinct colors determine every
			label. Thus the color classes, and only those classes, determine the closure,
			and inequivalent colorings give distinct closures. \item \emph{Perfect binary trees.}
			For the rooted tree, the additional external terminal
			accounts for the $+1$.\qedhere
		\end{enumerate}
	\end{proof}
	
	\subsubsection{Sharp bounds on failed coatom tests at fixed source size}
	
	The test for each terminal pair $a,b$ asks whether $C_{ab}$ remains
	closed under the source closure. For a sound relation define
	\begin{equation}\label{eq:tree-loss}
		\mathcal E_T(R)=
		\bigl|\{\{a,b\}:c_R(C_{ab})=E\}\bigr|.
	\end{equation}
	Since $c_R(C_{ab})$ is $c_T$-closed, it equals either $C_{ab}$ or $E$.
	The test fails when $C_{ab}$ is a key for $R$.
	
	\begin{theorem}[Minimum number of failed coatom tests at each source size]
		\label{thm:tree-loss-frontier}
		For every integer $m\ge1$, allowing arbitrary column alphabets,
		\begin{equation}\label{eq:tree-loss-frontier}
			\min_{\substack{R\text{ sound}\\ |R|=m}}\mathcal E_T(R)
			=\min_{\substack{R\text{ sound}\\1\le |R|\le m}}\mathcal E_T(R)
			=\max\left\{0,\binom n2-\binom m2\right\}.
		\end{equation}
		The sharp lower bound on the failure probability of a uniformly
		sampled coatom test is this quantity divided by $\binom n2$.
	\end{theorem}
	\begin{proof}
		A sound relation passes the test at $C_{ab}$ exactly when some row pair has
		agreement set $C_{ab}$, by maximality and distinctness of full rows. One row
		pair cannot witness two different coatoms. Thus at most $\binom m2$ tests
		are passed. For $m\le n$, retain any $m$ terminal rows of the canonical
		representation. Deleting rows preserves all valid dependencies, so the
		result is sound. Its pair disagreements are exactly the paths between
		retained terminals, giving exactly $\binom m2$ different coatom witnesses.
		For $m>n$, take a Cartesian power of the canonical $n$-row
		relation with at least $m$ rows. Each descriptor reads the tuple
		of its values in the factors. All original implications hold
		coordinatewise, and the $n$ diagonal rows witness every coatom.
		Retaining those rows and any $m-n$ others gives a sound relation
		of exactly $m$ rows passing every test. The asserted probability
		is the number of failures divided by the number of tests.
	\end{proof}
	
	Requiring chartwise exactness at the same row count can force more failed coatom tests.
	
	\begin{theorem}[Coatom tests for minimum local models]
		\label{thm:tree-three-loss}
		Let $T$ be trivalent and $R$ a three-sector chart-exact representation. Let
		$n_i$ be the number of terminal edges of color $i$, for its associated
		three-coloring.
		\begin{enumerate}[label=\textup{(\roman*)},leftmargin=*]
			\item \emph{Number of passed tests.}
			The number of passed coatom tests is exactly
			\begin{equation}\label{eq:tree-three-witness}
				\binom n2-\mathcal E_T(R)
				=\bigl|\{\{i,j\}\subseteq\{1,2,3\}:n_i+n_j=2\}\bigr|.
			\end{equation}
			\item \emph{Extremal tree shapes.}
			For $n\geq5$, this number is at most one. If the subtree induced by the
			internal vertices branches, it is zero for every minimum local model.
			If that subtree is a path (the tree is a caterpillar), some minimum local
			model passes one test. Thus, for $n\geq5$, the maximum number of passed tests
			among minimum local models is one for caterpillars and zero otherwise.
		\end{enumerate}
	\end{theorem}
	\begin{proof}
		\leavevmode
		\begin{enumerate}[label=\textup{(\roman*)},leftmargin=*]
			\item \emph{Number of passed tests.}
			The row pair indexed by colors $i,j$ disagrees exactly on the edges colored
			$i$ or $j$. At every internal vertex this subgraph has degree two, so its
			nontrivial components are terminal-to-terminal paths. Its number of paths
			is $(n_i+n_j)/2$. It witnesses a coatom precisely when there is exactly one
			such path. This proves \eqref{eq:tree-three-witness}.
			\item \emph{Extremal tree shapes.}
			Two different color pairs with sum two would, after renaming, give
			$n_1+n_2=n_1+n_3=2$, hence $n=4-n_1\leq4$. Thus at most one is possible
			when $n\geq5$.
			
			Every two-color subgraph contains every internal vertex. If it is a single
			path, all internal vertices lie on one path of $T$, so their induced subtree
			cannot branch. Conversely, when they form a path, alternate two colors on
			its edges, continue these colors to one terminal edge at each end, and color
			every other terminal edge with the third color. This is proper, and its
			chosen two-color subgraph is a single terminal path. For $n\geq5$, the first
			part shows that it passes exactly one coatom test.\qedhere
		\end{enumerate}
	\end{proof}
	
	\begin{example}[A depth-three tree with complete local tests and failed global tests]
		A perfect rooted binary tree of height three gives $n=9$ terminals after
		adjoining the external root terminal. Its internal subtree branches. There
		are $2^6=64$ different three-sector chart-exact global closures, and every
		one turns all $\binom92=36$ coatoms of $c_T$ into keys. In contrast, the
		faithful optimum has nine sectors and is unique at the partition level.
		
	\end{example}
	
	\subsubsection{Signed requirements after conditioning}
	
	The tree structure of the minimum faithful source gives an exact bound
	on the size of a minimal infeasible signed family after conditioning
	on positive requirements.
	Fix $D\subseteq E$ and consider
	\[
	X_D=\{C\in\operatorname{Fix}(c_T):D\subseteq C\}.
	\]
	This is the fiber of the trace to $D$ at the full code $D$. It imposes a
	positive baseline only; there are no additional excluded interface labels.
	Put $F=c_T(D)$. Let $H_1,\ldots,H_k$ be the nontrivial components of
	$T\setminus F$, and let $n_i\geq2$ be their numbers of original terminals.
	Each internal vertex in such a component has degree at least two.
	
	\begin{lemma}[Bound for signed tests from spectral partitions]
		\label{lem:partition-budget}
		Let $V$ be a finite descriptor set with finite-valued functions
		$f_v$ on a set $\Omega$ of $t$ spectral sectors. Let $c$ be
		the functional-determination closure of those functions. Fix
		$Q\subseteq V$ and a consistent exact trace
		$D\in\operatorname{Fix}(c|_Q)$. The joint-value partition for $D$
		has $b(D)$ blocks. Among closed access states $C$ satisfying
		$C\cap Q=D$, every inclusion-minimal infeasible signed family
		has at most $t-b(D)+1$ literals.
	\end{lemma}
	\begin{proof}
		The joint partition of the baseline has $b(D)$ blocks. A signed
		family is infeasible exactly when the closure of its positives together
		with $D$ forces an excluded trace label in $Q\setminus D$ or a
		negative literal. Every positive
		literal indispensable to a forcing contradiction must strictly refine
		the partition already generated by the baseline and preceding positive
		literals; otherwise it could be removed without affecting feasibility.
		Each strict refinement increases the block count by at least one,
		while the maximum is $t$. Thus there are at most $t-b(D)$ indispensable
		positives. In an inclusion-minimal infeasible family, a forced excluded
		boundary label needs no negative literal. In the other case
		one negative literal suffices to witness failure. This proves the bound.
	\end{proof}
	
	\begin{theorem}[Conditioned terminal-tree closure and its signed Helly number]
		\label{thm:tree-conditioned}
		\leavevmode
		\begin{enumerate}[label=\textup{(\roman*)},leftmargin=*]
			\item \emph{Residual product.}
			The interval of closed sets above $F$ is the direct product of the residual
			closure systems on $E(H_i)$, generated by the all-but-one rules at their
			remaining internal vertices. Labels of $F$ are fixed present.
			\item \emph{Signed Helly number.}
			With signed
			edge requirements on $X_D$, the exact signed Helly number is
			\begin{equation}\label{eq:tree-conditioned-rank}
				h_T(D)=
				\begin{cases}
					1,&F=E,\\
					\max_{1\leq i\leq k} n_i,&F\ne E.
				\end{cases}
			\end{equation}
			\item \emph{Comparison with the partition bound.}
			In the canonical $n$-sector source, if $b(D)$ is the number of baseline
			spectral blocks, then
			\[
			n-b(D)=\sum_{i=1}^k(n_i-1).
			\]
			When components remain, compare the generic bound $1+\sum_i(n_i-1)$
			from \cref{lem:partition-budget} with the exact value $\max_i n_i$.
			They agree for one unresolved component. The exact value is strictly
			smaller for two or more unresolved components.
		\end{enumerate}
	\end{theorem}
	\begin{proof}
		\leavevmode
		\begin{enumerate}[label=\textup{(\roman*)},leftmargin=*]
			\item \emph{Residual product.}
			A local rule above $F$ has all its already-present edges removed from the
			premise. Its remaining edges belong to a single component $H_i$. There is
			no rule connecting two such components, so closed supersets of $F$ are
			exactly independent choices of closed sets in the residual systems.
			\item \emph{Signed Helly number.}
			A minimally infeasible signed family either contains a negative on a fixed
			label (a singleton), or is contained in one component. Indeed, otherwise
			one failing component would already give a proper infeasible subfamily.
			
			For a residual component with $n_i$ terminals, the edge-cut representation
			on those terminals remains exact by the proof of \cref{lem:tree-paths};
			internal degree two causes repeated cuts but no failure of that proof.
			Its signed rank is at most $n_i$: in a minimal forcing premise each
			positive descriptor strictly refines the current partition of these $n_i$
			rows, so at most $n_i-1$ positives are indispensable, followed by one negative. If it has
			an internal vertex, all terminal-incident edges except one minimally force
			the remaining one. After removing any one positive edge, the path between
			its terminal and the target terminal witnesses failure. This gives a core
			of size $n_i$. If the component is a single edge, then $n_i=2$, and its
			positive-negative contradiction has size two. These facts prove
			\eqref{eq:tree-conditioned-rank}. If $F=E$, the only state is $E$ and every
			infeasible family contains a negative singleton.
			\item \emph{Comparison with the partition bound.}
			Two terminal rows agree on $D$ precisely when their path avoids $D$.
			No further deletion of a forced edge changes their connectivity, so the
			baseline blocks are exactly the terminal sets of components of
			$T\setminus F$, including isolated terminals. If their total number is
			$b(D)$, subtracting one from the size of each block gives the displayed
			identity.\qedhere
		\end{enumerate}
	\end{proof}
	
	For example, a baseline splitting the terminals into two unresolved blocks
	of sizes five and five has exact signed Helly number five. In a trace fiber
	$C\cap S=D$, the extra exclusions $S\setminus D$ can remove
	otherwise available extremal cores.
	
	\subsubsection{Presentation minima and source minima}
	
	A minimum presentation counts implications, whereas a minimum faithful
	relation counts rows. The following examples calculate the latter for
	two closures with distinct minimum tables. The second closure has a
	unique minimum circuit-clause presentation.
	
	\begin{example}[A classical closure with several minimum sources]
		\label{thm:uniform-nonrigid}
		Let $n\ge4$, $r=\binom n2$, and let $c$ be the rank-two uniform
		closure on $r$ labeled descriptors $E$:
		\[
		c(\varnothing)=\varnothing,\qquad c(\{e\})=\{e\},\qquad
		c(S)=E\quad(|S|\ge2).
		\]
		Its minimum faithful row count is $n$, but the number of inequivalent
		minimum tables, with descriptor labels fixed and allowing row permutations
		and independent renamings of values in each column, is
		\begin{equation}\label{eq:uniform-count}
			\frac{r!}{n!}.
		\end{equation}
		For $n=4$, this gives $30$ inequivalent four-row minima on six descriptors.
		The construction and row minimum are classical
		\cite[Example~6.6]{BeeriDFS1984}; see also
		\cite{DemetrovicsFurediKatona1985}. The displayed count
		comes from the row-permutation orbits.
	\end{example}
	
	\begin{proof}
		In a faithful table, two distinct rows can agree on at most one descriptor:
		agreement on two would force agreement on every descriptor. For each
		$e$, faithfulness requires a pair agreeing on $e$ and disagreeing on any
		chosen $f\ne e$; its exact agreement set is therefore $\{e\}$.
		Thus the $r$ distinct singleton agreement sets require $r$ different row
		pairs, giving at least $n$ rows.
		
		Conversely, assign the $r$ labeled descriptors bijectively to the row
		pairs of an $n$-element row set. In column $e$, give its assigned pair one
		common value and give every other row a separate value. Every row pair
		then agrees precisely on its assigned descriptor. The resulting closure
		is $c$, so the lower bound is attained.
		
		In any $n$-row faithful table all $r$ pairs must be used once. In each
		column there is exactly one equal pair, so the preceding construction
		determines every minimum table up to column value renaming. There are
		$r!$ assignments. The row-permutation action is free, because for
		$n\ge3$ a permutation fixing every unordered pair fixes every row.
		Dividing by $n!$ proves \eqref{eq:uniform-count}.
	\end{proof}
	
	\begin{example}[Unique minimum circuit clauses without source rigidity]
		\label{ex:horn-source-nonrigid}
		On the five labeled descriptors $E=\{a,b,c,d,e\}$, let $c$ be generated
		by $ab\to c$, $ac\to b$, and $bc\to a$, with $d,e$ free. This is the
		simple binary matroid $U_{2,3}\oplus U_{1,1}\oplus U_{1,1}$, represented
		graphically by a triangle and two bridges. Its only circuit is
		$\{a,b,c\}$, so its three circuit clauses form the unique minimum
		circuit-clause presentation. Nevertheless the following two faithful
		tables are inequivalent minimum representations:
		\[
		\begin{array}{c|ccccc}
			R_*&a&b&c&d&e\\\hline
			&1&0&0&0&0\\
			&0&1&0&0&0\\
			&0&0&1&0&0\\
			&1&0&0&1&0\\
			&1&0&0&0&1
		\end{array}
		\qquad
		\begin{array}{c|ccccc}
			R'&a&b&c&d&e\\\hline
			&1&0&0&0&0\\
			&0&1&0&0&0\\
			&0&0&1&0&0\\
			&1&0&0&1&0\\
			&1&0&0&1&1
		\end{array}
		\]
		Both have closure $c$. The $d$-column partitions have block sizes
		$1,4$ and $2,3$, respectively, so row permutations and independent
		value renamings cannot identify the tables with descriptor labels fixed.
		This is a direct-product closure in the sense of
		\cite{DemetrovicsFurediKatona1985}. The agreement-set argument below
		gives its five-row minimum.
	\end{example}
	
	\begin{proof}
		The proper meet-irreducibles are
		\[
		\{a,d,e\},\quad\{b,d,e\},\quad\{c,d,e\},\quad
		\{a,b,c,e\},\quad\{a,b,c,d\}.
		\]
		In each displayed table the three pairs among the first three rows
		realize the first three sets. The remaining two are witnessed by row pairs $(1,4)$,
		$(1,5)$ in $R_*$ and $(1,4)$, $(4,5)$ in $R'$. Every pair agreement
		set is closed under the three rules, so the classical agreement-set
		criterion proves faithfulness.
		
		For a direct size check, suppose that a faithful table had at most four
		rows. The first three meet-irreducibles require three distinct pairs
		inside fibers of the joint readout $(d,e)$. This readout is nonconstant.
		On at most four rows, three such pairs are possible only with fiber
		sizes $3,1$. A witness for $\{a,b,c,e\}$ must cross these fibers and
		disagree only in $d$; a witness for $\{a,b,c,d\}$ must cross them and
		disagree only in $e$. The same two joint $(d,e)$ values cannot meet both
		requirements. Hence at least five rows are necessary, and the displayed
		tables attain the minimum.
	\end{proof}
	
	\subsection{Native spectral minima over fixed local spectra}
	\label{subsec:fixed-spectral-source-minima}
	
	Fix the local spectra, overlap identifications, and access--commutant
	realization law. On a terminal tree the opposite-polarity interfaces
	give a state-exact, faithful matching source with $n$ sectors, the
	same count as the closure-only row minimum. For a general fixed join
	tree, the state-exact faithful minimum is a finite witness-cover
	problem. Exactness at full cuts is an additional condition when required.
	
	\subsubsection{Terminal-tree spectral amalgamation and additive source cost}
	\label{subsec:tree-spectral-gluing}
	
	The three-row models realize the chart closures. Here the source must
	also cover the local spectra prescribed by the tree. The terminal
	geometry identifies the global states and gives the sharp spectral
	minimum, with equality cases in \cref{thm:tree-rigidity}.
	Let $T$ be a finite tree with $n\ge3$ terminals $\Omega$, internal vertices
	$I$, and all internal degrees at least three. Descriptor labels are the
	edges $E(T)$, and at every internal vertex $v$ the intended local
	formation rules are
	$\delta_T(v)\setminus\{e\}\to e$ for each incident edge $e$.
	The least generated global closure is the $c_T$ defined in
	\cref{subsec:terminal-tree-closure-minima}.
	
	For each $v\in I$ take the \emph{native local spectrum}
	$X_v=\delta_T(v)$: a state chooses one incident edge. The local
	observable for incident edge $e$ is its characteristic function
	$1_{\{e\}}$ or its complement. On a terminal edge use the
	characteristic function. On an internal edge $e=vw$, choose one of
	its endpoints, say $v$, and give the two local observables opposite
	polarities:
	\begin{equation}\label{eq:opposite-polarity}
		b_{v,e}(d)=\Ind[d=e],\qquad
		b_{w,e}(d')=1-\Ind[d'=e].
	\end{equation}
	The common interface is the two-state spectrum $Y_e=\{0,1\}$;
	\eqref{eq:opposite-polarity} makes both maps onto. The equality
	$b_{v,e}(d_v)=b_{w,e}(d_w)$ is equivalent to saying that
	\emph{exactly one endpoint chooses the internal edge $e$}.
	A local chart with $d=\deg_T(v)$ states has precisely the declared
	all-but-one recovery closure: $d-1$ of its indicator coordinates determine
	the last, while omitting two or more coordinates leaves each omitted
	coordinate undetermined. Complementing any binary
	coordinate leaves this local closure unchanged.
	
	\begin{theorem}[Faithful spectral amalgamation and rigidity for terminal-tree atlases]
		\label{thm:native-spectral-tree}
		For the local spectra and overlap maps \eqref{eq:opposite-polarity},
		write
		\begin{equation}\label{eq:tree-spectral-limit}
			X_T=\{(d_v)_{v\in I}\in\prod_{v\in I}X_v:
			b_{v,e}(d_v)=b_{w,e}(d_w)\text{ for every internal }e=vw\}.
		\end{equation}
		The following hold.
		\begin{enumerate}[label=(\roman*),nosep]
			\item \emph{Terminal sectors.} There is a canonical bijection $\Omega\xrightarrow{\sim}X_T$:
			a terminal $t$ assigns to each internal vertex the first edge on the
			path toward $t$. Hence $|X_T|=n$.
			\item \emph{State coverage.} Every local chart sector extends to a global sector. In fact, for
			each terminal $t$ the adjacent chart state choosing its terminal edge
			has $t$ as its \emph{unique} global preimage. Therefore no proper
			subset of $X_T$ is state-exact for this fixed native atlas.
			\item \emph{Faithfulness.} The globally glued edge observables are precisely the binary cut
			functions of terminal paths (up to independent value complements),
			and their recovery closure is exactly $c_T$. Thus $X_T$ is a
			faithful $n$-sector native source.
			\item \emph{Rigidity at the minimum.} Every faithful representation of $c_T$ has at least $n$
			distinct sectors, and every $n$-sector faithful representation is
			equivalent to this terminal-cut source up to sector permutation
			and descriptor-value renaming.
		\end{enumerate}
	\end{theorem}
	\begin{proof}
		\leavevmode
		\begin{enumerate}[label=\textup{(\roman*)},leftmargin=*]
			\item \emph{Terminal sectors.}
			The internal graph on $I$ is a tree, with $k=|I|$ vertices and $k-1$
			internal edges. Every internal edge is chosen at exactly one endpoint
			in a compatible global tuple, so the $k$ vertex choices consist of
			$k-1$ internal-edge choices and precisely one terminal-edge choice.
			If $t$ is that terminal, start at its incident internal vertex and
			propagate away from $t$: the first vertex chooses $t$, every other
			vertex must choose the first edge on the path to $t$. Conversely this
			assignment satisfies every interface equality. The terminal direction therefore has a unique global preimage.
			\item \emph{State coverage.}
			Every local direction at
			$v$ points into a component containing at least one terminal, so all
			local states occur. Omitting any terminal would omit its uniquely extended terminal direction.
			\item \emph{Faithfulness.}
			For any edge $e$, the glued binary observable separates terminals
			according to the two components of $T\setminus\{e\}$: at the terminal
			$t$, the direction toward $t$ lies on one side or the other of $e$.
			Thus two global states $a,b\in\Omega$ disagree on $e$ exactly when
			$e$ lies on the connecting path $P_{ab}$. A counterexample to recovery of $e$ from $S$ is a terminal path
			with $P_{ab}\cap S=\varnothing$ and $e\in P_{ab}$. Recovery holds
			exactly when the set of such paths is empty. By
			\cref{lem:tree-paths}, this is exactly $e\in c_T(S)$, proving (iii).
			\item \emph{Rigidity at the minimum.}
			The lower bound and classification at equality follow from
			\cref{thm:tree-rigidity}.\qedhere
		\end{enumerate}
	\end{proof}
	
	The sector count can also be obtained by exact local recursion. Cut any
	internal edge $e$, separating the terminal set into $\Omega_L$ and
	$\Omega_R$, with $n_L,n_R>0$. Let $U_L$ and $U_R$ denote the locally
	compatible partial-state spaces on the two sides, retaining the
	boundary bit on $e$. Each side has one extra sector recording the
	option that its unique chosen direction exits through $e$:
	\begin{equation}\label{eq:additive-fibers}
		\begin{gathered}
			|U_L|=n_L+1,\quad |U_R|=n_R+1,\\
			(|U_{L,0}|,|U_{L,1}|)=(1,n_L),\quad
			(|U_{R,0}|,|U_{R,1}|)=(n_R,1),
		\end{gathered}
	\end{equation}
	where the bits have been named so that the common interface values
	match the displayed fiber counts.
	
	\begin{theorem}[Exact separators and additive source size for terminal-tree atlases]
		\label{thm:exact-native-separator}
		\leavevmode
		\begin{enumerate}[label=\textup{(\roman*)},leftmargin=*]
			\item \emph{Partial states and source size.}
			For every internal edge cut as above,
			\begin{equation}\label{eq:exact-native-rectangle}
				X_T\cong U_L\times_{\{0,1\}}U_R,\qquad
				|X_T|=n_L+n_R.
			\end{equation}
			\item \emph{Common functions.}
			On $\ell^2(X_T)$ let $M_L$ and $M_R$ be the multiplication algebras
			of functions of the left and right partial states, respectively.
			For each edge $f$ of $T$, let $A_f$ denote the unital multiplication
			algebra generated by its binary observable. Then
			\begin{equation}\label{eq:exact-native-intersection}
				\boxed{M_L\cap M_R=A_e.}
			\end{equation}
			\item \emph{Native formation.}
			Consequently the common source $M=\mathbb C^{X_T}$ admits exact
			native recovery and formation:
			\begin{equation}\label{eq:native-commutant-promotion}
				\begin{gathered}
					M(S)=\operatorname{alg}^{*}\{A_e:e\in S\},\qquad
					K(S)=M(S)',\\
					K(S\cup U)=K(S)\cap K(U),\quad K(S)'=M(S).
				\end{gathered}
			\end{equation}
		\end{enumerate}
	\end{theorem}
	\begin{proof}
		\leavevmode
		\begin{enumerate}[label=\textup{(\roman*)},leftmargin=*]
			\item \emph{Partial states and source size.}
			The counting argument in the preceding proof applies separately to
			each partial tree, with the exposed boundary edge treated as another
			possible terminal direction. It gives $n_L+1$ and $n_R+1$ sectors;
			exactly one sector on each side chooses to exit, while the other
			$n_L$ or $n_R$ sectors choose terminals in that side. Opposite
			polarities make the exit choice on one side compatible with a terminal
			choice on the other, proving \eqref{eq:additive-fibers} and
			\eqref{eq:exact-native-rectangle}.
			\item \emph{Common functions.}
			For each interface bit $y$, the compatible joint fiber is the
			\emph{entire Cartesian rectangle} $U_{L,y}\times U_{R,y}$; both
			factors are nonempty. A function that factors through each of the two partial-state projections is constant on every separator rectangle. Conversely, every function of $y$ factors through either side.
			This gives \eqref{eq:exact-native-intersection}.
			\item \emph{Native formation.}
			Finally, in the faithful diagonal representation on
			$\ell^2(X_T)$, finite-dimensional bicommutant recovery gives
			$M(S)''=M(S)$. The identity
			$(M(S)\vee M(U))'=M(S)'\cap M(U)'$ proves
			\eqref{eq:native-commutant-promotion}.\qedhere
		\end{enumerate}
	\end{proof}
	
	The opposite-polarity interfaces also satisfy the local lifting criterion
	of \cref{thm:two-state-joint-separators}.
	
	\subsubsection{An exact minimum-size criterion for fixed local sources}
	
	The join $X$ in \eqref{eq:native-joint-spectrum} is the largest compatible
	fixed-atlas source. A smaller $R\subseteq X$ can preserve all required
	local sectors and all required nonimplications only if it keeps enough
	single-state and pair-state witnesses. For a chart state $x_i$ put
	$U_{i,x_i}=\{x\in X:\pi_i(x)=x_i\}$, and for a required nonimplication
	$e\notin c_0(S)$ let
	\[
	W_{S,e}=\{\{x,x'\}\subseteq X:\ x\ne x',\
	f_s(x)=f_s(x')\ (s\in S),\quad f_e(x)\ne f_e(x')\}.
	\]
	\begin{proposition}[Exact fixed-atlas faithful-source optimization]
		\label{prop:native-minimum-cover}
		\leavevmode
		\begin{enumerate}[label=\textup{(\roman*)},leftmargin=*]
			\item \emph{Existence.}
			For finite join-tree atlas data with equal realized separator images
			on every tree edge, a state-exact faithful subamalgam exists if and only
			if $c_X=c_0$. On any finite atlas, the same equivalence holds
			provided the full join $X$ projects onto every chart. If a chart
			state is absent from $X$, state-exact assembly is impossible.
			\item \emph{Minimum size.}
			Whenever a state-exact faithful subamalgam exists, its minimum number of
			distinct sectors is exactly
			\begin{equation}\label{eq:minimum-native-cover}
				\min\left\{|R|:
				\begin{gathered}
					R\subseteq X,\quad R\cap U_{i,x_i}\ne\varnothing\ \text{for every }i,x_i,\\
					\forall e\notin c_0(S)\ \exists\{x,x'\}\in W_{S,e}
					\text{ with }\{x,x'\}\subseteq R
				\end{gathered}\right\}.
			\end{equation}
			Every finite common source whose observables are generated by the
			fixed local sources factors through a subset of $X$, so this
			optimization is exact under the stated fixed-atlas contract.
			\item \emph{Meet-irreducible tests.}
			Equivalently, one may require a pair with exact agreement set $M$
			for each proper meet-irreducible closed set $M$ of $c_0$
			\cite{BeeriDFS1984}.
			\item \emph{Terminal-tree atlas.}
			For the oriented terminal-tree atlas,
			\eqref{eq:minimum-native-cover} equals $n$ and its only
			inclusion-minimal solution is $R=X_T$.
		\end{enumerate}
	\end{proposition}
	\begin{proof}
		\leavevmode
		\begin{enumerate}[label=\textup{(\roman*)},leftmargin=*]
			\item \emph{Existence.}
			On a join tree with equal separator images, \cref{prop:closure-descent}
			gives full chart projections. On any finite atlas, a subrelation
			of $X$ can cover only chart states already present in $X$.
			If $c_X\ne c_0$, some
			required $W_{S,e}$ is empty and no subset can be faithful. If $X$ covers all chart states and
			$c_X=c_0$, then $X$ itself is feasible.
			\item \emph{Minimum size.}
			Every compatible global state has well-defined local projections,
			so a source with no extraneous state distinctions is a subset of $X$.
			The surjectivity requirements are precisely the unary coverage
			conditions $R\cap U_{i,x_i}\ne\varnothing$. By definition, an FD
			$S\to e$ fails in $R$ precisely when $R$ contains both endpoints of a
			pair in $W_{S,e}$. Thus all and only the required nonimplications are
			preserved exactly under the pair conditions.
			\item \emph{Meet-irreducible tests.}
			For the meet-irreducible version, all row-pair
			agreement sets are $c_0$-closed. A proper meet-irreducible $M$
			has a unique cover $M^+$ in this finite closure lattice. For
			$e\in M^+\setminus M$, a witness against $M\to e$ must have
			agreement set exactly $M$, since every larger closed set contains
			$M^+$. Conversely, every closed set is an intersection of
			meet-irreducibles above it, so if $e\notin c_0(S)$, one such
			$M$ contains $S$ and excludes $e$. A witness pair for $M$
			therefore refutes $S\to e$.
			\item \emph{Terminal-tree atlas.}
			Finally, for every terminal $t$ in the
			native tree atlas, the local sector choosing its terminal edge has
			unique global preimage $t$. Every feasible $R$ therefore contains
			every terminal, proving the final assertion.\qedhere
		\end{enumerate}
	\end{proof}
	
	The optimization combines a coverage requirement for each chart state
	with a pair requirement for each proper meet-irreducible. Finite joins
	and witness propagation implement these tests
	\cite{Beeri1983,BeeriDFS1984,SagivShmueli1989,AtseriasKolaitis2025,AbramskyCaru}.
	
	\section{Witness descent and globalization obstructions}
	\label{sec:native-witness-descent}
	
	A pair of chart states witnessing a failed dependency may extend
	along each path separately and still fail to extend around a cycle.
	For bijective transport, the possible root pairs are the common
	fixed points of the cycle transports (\cref{thm:native-witness-holonomy}).
	For finite multivalued formation, \cref{thm:witness-globalization}
	tests one exact-agreement witness per proper meet-irreducible.
	Linear and group diagrams give tests for coupled cycles
	(\cref{thm:graph-relative,thm:joint-chord-witness}). We also treat
	lifting through a quotient and the cost of repairing one missing witness.
	
	For group matching formations, \cref{thm:complete-compositional-invariants}
	computes recovery from joint boundary agreement images. Here the lifting
	tests can also prescribe a preimage at one vertex.
	
	\subsection{Witness pairs and transport around cycles}
	\label{subsec:structural-graph-witness}
	
	Fix the required agreement pattern on the full native records.
	A local pair extends if its continuations match at every interface.
	
	\begin{definition}[Realized-source witness atlas]\label{def:witness-atlas}
		Let $G=(V,E)$ be a finite connected graph.
		A \textbf{realized-source witness atlas} on $G$ is specified by the
		following local records, interface maps, and readouts, together with
		the witness-pair relations determined by a target closure.
		\begin{enumerate}[label=(\roman*),leftmargin=*]
			\item \emph{Local native records.} At each vertex $v$,
			choose a family $\mathcal B_v$ of actual realized triples
			$\mathbb P=(\Xi,P,\mathfrak P)$ and a specified observation construction
			$O_v(\mathbb P)$ as in
			\cref{def:associated-observation-presentation}. A finite set of
			\emph{realized native records} is a specified subset
			\[
			Z_v\subseteq\bigsqcup_{\mathbb P\in\mathcal B_v}
			\{\mathbb P\}\times O_v(\mathbb P).
			\]
			Every record thus retains its native judgment and its chosen observation.
			The direct selected-state interpretation uses the tagged set
			$O_v(\mathbb P)=\bigsqcup_{K\in\Pkg_{\mathfrak P}}\RR_{\mathfrak P}(K)$;
			a spectral interpretation uses the observations of
			\cref{prop:spectral-observation-interpretation}.
			\item \emph{Interface matching.} The maps below act on the entire record. For an
			edge $e=\{v,w\}$ let $I_e$ be an interface and let
			$\rho_{v,e}:Z_v\to I_e$, $\rho_{w,e}:Z_w\to I_e$ be specified
			restriction maps.  The global matching set retains the full records.
			\[
			Z_G=\{(z_v)\in\prod_{v\in V}Z_v:
			\rho_{v,e}(z_v)=\rho_{w,e}(z_w)
			\text{ for every }e=\{v,w\}\}.
			\]
			\item \emph{Global readouts.} Let $D$ be a finite descriptor set, with every descriptor observed on
			at least one vertex.  For each $d\in D$, require that every $z=(z_v)\in Z_G$ has the same
			$d$-readout at all vertices observing $d$. This global consistency
			hypothesis makes $d(z)$ well defined even when its observing vertices
			are not adjacent. It follows, for example, if the observing vertices
			form a connected subgraph and each of its edge interface matches
			preserves the $d$-readout; edgewise agreement alone without that
			connectedness need not suffice.
			\item \emph{Witness pairs.} Fix a target closure $c_0$ on $D$. For a proper $c_0$-closed
			set $M\subsetneq D$, let
			$W_v(M)\subseteq Z_v\times Z_v$ consist of the ordered realized-source
			pairs that agree on every descriptor of $M$ observed at $v$ and
			disagree on every descriptor of $D\setminus M$ observed at $v$.
			For an edge $e=\{v,w\}$, its \emph{paired witness matching relation}
			$C_e(M)\subseteq W_v(M)\times W_w(M)$ is defined by
			\[
			((z,z'),(y,y'))\in C_e(M)
			\quad\Longleftrightarrow\quad
			\bigl(\rho_{v,e}(z),\rho_{v,e}(z')\bigr)
			=\bigl(\rho_{w,e}(y),\rho_{w,e}(y')\bigr).
			\]
			Matching uses the full records in $Z_v$. A projection to endpoints
			can be used first only if it preserves the required formation
			and interface data.
		\end{enumerate}
	\end{definition}
	
	\begin{lemma}[Exact witness descent]\label{lem:witness-exact-descent}
		With the preceding hypotheses, ordered pairs of global realized states
		in $Z_G$ whose \emph{exact} descriptor agreement set is $M$ correspond
		bijectively to families $a_v\in W_v(M)$ satisfying
		$(a_v,a_w)\in C_e(M)$ on every edge.  In particular, local
		nonemptiness and compatibility on each edge separately do not imply
		a global witness.
	\end{lemma}
	\begin{proof}
		Take the two componentwise projections of a matching witness family.
		Each separately satisfies all interface equations and hence belongs to
		$Z_G$.  The local readout conditions yield equality on $M$ and
		inequality on $D\setminus M$, since each descriptor is observed
		somewhere.  Conversely, a global pair with exact agreement $M$
		restricts to a pair in each $W_v(M)$ and satisfies all paired
		matching relations.
		The constructions are inverse.
	\end{proof}
	
	\paragraph{Separate and joint feasibility of paired constraints.}
	Lemma~\ref{lem:witness-exact-descent} requires a pair satisfying
	all constraints at once. On the affine heap $W=\mathbb F_2^2$ with
	coordinates $x,y$, consider the following pair relations.
	\[
	\begin{aligned}
		Q_x&=\{(u,v)\in W^2:u_y=v_y,\ u_x\ne v_x\},\\
		Q_y&=\{(u,v)\in W^2:u_x=v_x,\ u_y\ne v_y\}.
	\end{aligned}
	\]
	The relations are both nonempty and have empty intersection
	$Q_x\cap Q_y=\varnothing$. If
	$\pi_{\varnothing}$ forgets $x,y$ in both copies, then
	\[
	\pi_{\varnothing}(Q_x\cap Q_y)=\varnothing,
	\qquad
	\pi_{\varnothing}(Q_x)\cap\pi_{\varnothing}(Q_y)
	=\{(\langle\rangle,\langle\rangle)\}.
	\]
	Separate projections discard the identity of the pair.
	\begin{definition}[Bijective witness transport]\label{def:unique-witness-transport}
		For a fixed $M$, the atlas has \emph{bijective witness transport} if
		every $W_v(M)$ is nonempty and, for each oriented edge $v\to w$,
		$C_e(M)$ is the graph of a bijection
		$\tau_{vw}:W_v(M)\to W_w(M)$ with
		$\tau_{wv}=\tau_{vw}^{-1}$.
	\end{definition}
	
	\begin{theorem}[Common fixed points of cycle transport]\label{thm:native-witness-holonomy}
		Assume bijective witness transport for $M$.  Choose a spanning tree
		$T$ of $G$, a root $o$, and write $W=W_o(M)$.  Let
		$t_v:W\to W_v(M)$ be the unique transport along the tree from $o$
		to $v$.  For each oriented edge $e=(v,w)$ outside the tree, compose the
		transport maps around the corresponding cycle. This permutation
		of the root labels is its \emph{holonomy}.
		\[
		h_e=t_w^{-1}\circ\tau_{vw}\circ t_v\in\operatorname{Sym}(W).
		\]
		Then the global pairs with exact agreement set $M$ are in natural
		bijection with
		\begin{equation}\label{eq:native-holonomy-fixed}
			\quad\operatorname{Fix}_W\langle h_e:e\notin T\rangle
			=\bigcap_{e\notin T}\operatorname{Fix}_W(h_e).\quad
		\end{equation}
		In particular, all $|W|$ root witnesses extend on a tree; on an
		arbitrary graph, a global witness exists if and only if the common
		fixed-point set is nonempty.  This criterion is independent of the
		chosen spanning tree.
	\end{theorem}
	\begin{proof}
		An assignment compatible along the tree is uniquely determined by
		its root label $a\in W$, because all tree transports are bijections;
		it has value $t_v(a)$ at $v$.  Compatibility on a chord $e=(v,w)$
		is the equation $\tau_{vw}t_v(a)=t_w(a)$, equivalent to
		$h_e(a)=a$.  Thus exactly the common fixed points extend across
		all chords.  Every choice of $T$ describes the global matching space in
		\cref{lem:witness-exact-descent}. Its nonemptiness and cardinality
		are therefore independent of the spanning tree.
	\end{proof}
	
	\begin{theorem}[Sharp cycle bounds for failure of witness globalization]
		\label{thm:witness-cycle-core}
		In the setting of \cref{thm:native-witness-holonomy}, let $M\subsetneq D$.
		\begin{enumerate}[label=(\roman*),leftmargin=*]
			\item \emph{Reversal of witnesses.} Witness reversal $\iota(z,z')=(z',z)$ induces a fixed-point-free
			involution on $W=W_o(M)$, every holonomy commutes with $\iota$, and
			$m=|W|$ is even.
			\item \emph{Sharp cycle bound.} If no global exact-agreement witness exists,
			\begin{equation}\label{eq:native-half-cycle-core}
				\begin{gathered}
					\text{at most }m/2\text{ fundamental-cycle holonomies}\\
					\text{already have no common fixed point}.
				\end{gathered}
			\end{equation}
			The bound is sharp for every even $m\ge2$ within finite
			atlases of pairs of realized records, including examples with
			nonempty global state sets.
			\item \emph{Unrestricted permutations.} For arbitrary permutations of a set of $m$ labels with $m\ge2$, the sharp
			bound is $m-1$.
		\end{enumerate}
	\end{theorem}
	\begin{proof}\leavevmode
		\begin{enumerate}[label=(\roman*),leftmargin=*]
			\item \emph{Reversal of witnesses.} Choose a vertex observing a descriptor in $D\setminus M$, which exists
			by hypothesis. Every exact-agreement pair there has distinct entries,
			so reversal is free at that vertex. The edge matching relations are
			invariant under simultaneous reversal. By uniqueness of transport,
			$\tau_{vw}\iota=\iota\tau_{vw}$. Bijective transport along a path
			therefore carries freeness to the chosen root, even if the root itself
			does not observe a descriptor outside $M$. Each fundamental holonomy
			commutes with $\iota$.
			\item \emph{Sharp cycle bound.} The fixed-point set of each holonomy is a union of
			$\iota$-orbits.  There are $n=m/2$ such two-element orbits.
			If there is no common fixed point, select a holonomy that excludes
			at least one orbit.  If any orbits survive, choose a surviving
			witness and a further holonomy moving it; commutation with
			$\iota$ excludes its entire two-element orbit.  After at most
			$n$ such choices no orbit survives.
			
			For sharpness put $n=m/2$, and at every vertex of a connected graph
			with $n$ independent fundamental cycles use the realized state set
			$Z_v=\{\star\}\sqcup(\{1,\ldots,n\}\times\{0,1\})$.
			The shared descriptor $A$ records the block label, taking the value
			$\star$ at the distinguished state; vertex $v$ also has a private
			binary descriptor $d_v$ that reads the second bit (zero at $\star$).
			For $M=\{A\}$, the local exact-agreement witnesses are precisely
			the $2n$ ordered pairs within a common two-state block.  Take all
			interface state maps to be bijections of $Z_v$ preserving $A$.
			On the chosen spanning tree use identity maps at both incidences.
			At chord $j$, use the identity on one incidence and the map that
			flips the two states in block $j$ and fixes all others on the other.
			The induced holonomy reverses exactly the two ordered witnesses from
			block $j$ and fixes the others.  All $n$ chords together have no
			common fixed witness, but omitting any one chord leaves its block's
			ordered witnesses fixed.  The global realized-source space is
			nonempty because $\star$ survives every interface.  Thus $n$
			cycles are sometimes necessary.
			\item \emph{Unrestricted permutations.} For unrestricted permutation witnesses on $m$ labels, a nonidentity
			permutation moves at least two labels; the same greedy argument
			therefore gives $m-1$.  The transpositions $(0\ j)$,
			$1\le j<m$, show sharpness for unrestricted permutations.\qedhere
		\end{enumerate}
	\end{proof}
	
	\begin{corollary}[Faithful globalization under bijective witness transport]\label{cor:holonomy-faithful-Armstrong}
		Let $c_0$ be a closure operator on the finite descriptor set $D$, and
		suppose the global source $Z_G$ is nonempty and satisfies every
		dependency of $c_0$, so $c_0\le c_{Z_G}$ pointwise.  For every proper
		meet-irreducible closed set $M$ of $c_0$, form the exact-agreement
		witness atlas of \cref{def:witness-atlas}.  If each of these atlases
		has bijective witness transport, the following hold.
		\begin{enumerate}[label=(\roman*),leftmargin=*]
			\item \emph{Faithfulness.} \begin{equation}\label{eq:witness-holonomy-faithful}
				\boxed{\quad c_{Z_G}=c_0
					\quad\Longleftrightarrow\quad
					\bigcap_{e\notin T}\operatorname{Fix}(h_{e,M})\ne\varnothing
					\text{ for every proper meet-irreducible }M.\quad}
			\end{equation}
			\item \emph{Certificate size.} Every failure of faithful recovery has a certificate involving
			at most $|W_o(M)|/2$ fundamental cycles for at least one $M$.
			State coverage and exact separator intersections require separate tests.
		\end{enumerate}
	\end{corollary}
	\begin{proof}\leavevmode
		\begin{enumerate}[label=(\roman*),leftmargin=*]
			\item \emph{Faithfulness.} For a finite set of global states, its dependency closure is the
			intersection closure of their pairwise exact agreement sets.  Since
			$c_0\le c_{Z_G}$, all such agreement sets are $c_0$-closed.
			The agreement-set characterization \cite[Theorem~6.1]{BeeriDFS1984} uses the
			following fact. Proper meet-irreducibles generate all closed sets
			by intersection, and every family that does so contains them.
			Therefore $c_{Z_G}=c_0$ if and only if every such $M$ occurs as the
			exact agreement set of some global pair.  Apply
			\cref{lem:witness-exact-descent,thm:native-witness-holonomy}.
			\item \emph{Certificate size.} If faithfulness fails, choose a missing meet-irreducible witness $M$.
			Then \cref{thm:witness-cycle-core} gives the stated cycle bound.\qedhere
		\end{enumerate}
	\end{proof}
	
	\begin{remark}[Relational composition and cycle constraints]
		\label{rem:witness-relational-holonomy}
		For arbitrary finite paired matching relations $C_i$ on a single
		cycle, exact witness existence is the relational fixed-point test
		\[
		\Delta_{W_0}\cap (C_0;C_1;\cdots;C_{s-1})\ne\varnothing,
		\]
		where $;$ denotes relational composition and the composite is a
		relation on $W_0$. On a triangle with two identity transports and one flip
		on $\{0,1\}$, this intersection is empty. On a graph with several cycles all chord constraints
		must hold on the same tree-compatible witness family. For instance,
		in the four-witness construction of \cref{thm:witness-cycle-core},
		each block flip fixes the other block, but the two flips have no common
		fixed witness. For translations of an abelian group, the holonomy of
		each cycle is translation by its signed sum, so a global witness exists
		exactly when every such sum vanishes. For graph descent and its
		connection with sheaf constructions, see
		\cite{AbramskyBrandenburger,Stacks}.
	\end{remark}
	
	At a protected cut, \cref{thm:proper-cut-incidence} tests the additional
	identity $M_L\cap M_R=A_e$ on the incidence graph of the actual side images.
	\subsection{Witness extension with multivalued relations}
	\label{subsec:structural-multivalued-witness}
	
	Different internal sources may realize the same endpoint, giving a
	witness several continuations. Finite relations with multiple inputs
	and outputs supply additional Level-II data to the realized system of
	\cref{def:standalone-realized-system} and conditions (H1)--(H3).
	Constraint joins and valuation algebras provide the relational
	operations used to test recovery and repair on full records
	\cite{AbramskyCaru,BeeriDFS1984,AbramskyBrandenburger}.
	
	\begin{definition}[Finite multivalued native formation atlas]
		\label{def:multi-native-atlas}
		Let $V$ be a finite set of package positions, with finite realized-record
		sets $Z_v$ and their declared observation interpretations as in
		\cref{def:witness-atlas}.
		A \textbf{finite multivalued native formation atlas} on these records
		is specified by the following formation relations and descriptor
		readouts, together with their associated paired relations.
		\begin{enumerate}[label=(\roman*),leftmargin=*]
			\item \emph{Records and formation relations.} Each record retains an actual source--endpoint
			judgment; its selected state or associated observation has the stated
			separate type. Let $\mathcal F$ be a finite index set of formation
			constraints.  Every $\alpha\in\mathcal F$ has a nonempty scope
			$S_\alpha\subseteq V$ and an admissible relation
			\[
			F_\alpha\subseteq\prod_{v\in S_\alpha}Z_v.
			\]
			The global realized-source space consists of the tuples whose restrictions
			belong to every $F_\alpha$.
			\item \emph{Descriptor readouts.} Fix a finite descriptor set $D$ observed collectively
			by the positions and require the readouts of shared descriptors to agree
			on admissible global tuples.  Any shared-descriptor consistency
			requirements are treated as protected, non-editable constraints in a
			repair problem. Keep the constraints enforcing this agreement
			when forming a protected paired matching set $\Lambda_0(M)$.
			\item \emph{Paired relations.} For $M\subsetneq D$ let $W_v(M)$ be the
			ordered local realized-state pairs whose descriptor agreement pattern at
			$v$ equals the restriction of $M$.  (If a descriptor is not observed at
			$v$, it imposes no condition there.)  Put
			\[
			\mathcal P_M=\prod_{v\in V}W_v(M),\qquad
			\widetilde F_\alpha(M)=
			\left\{w\in\prod_{v\in S_\alpha}W_v(M):
			(\pi_1w_v)_{v\in S_\alpha},
			(\pi_2w_v)_{v\in S_\alpha}\in F_\alpha\right\}.
			\]
			Both $F_\alpha$ and $\widetilde F_\alpha(M)$ may be general relations.
			They can describe a binary interface, a formation with several inputs,
			or a restriction at one site.
		\end{enumerate}
	\end{definition}
	A tuple of matching local records gives an observational witness.
	Realizing it as a native endpoint requires the promotion construction
	of \cref{def:associated-observation-presentation}.
	
	\begin{theorem}[Exact multivalued native witness descent]
		\label{thm:multi-native-exact-descent}
		Fix a finite multivalued native atlas and $M\subsetneq D$.
		\begin{enumerate}[label=(\roman*),leftmargin=*]
			\item \emph{General descent.} The pairs of
			global realized states with \emph{exact} descriptor agreement $M$ are in
			natural bijection with
			\begin{equation}\label{eq:multi-native-join}
				\quad\Lambda(M)=
				\left\{w\in\mathcal P_M:
				w|_{S_\alpha}\in\widetilde F_\alpha(M)
				\text{ for all }\alpha\in\mathcal F\right\}.\quad
			\end{equation}
			\item \emph{Cycle test.} In a simple cyclic graph with paired edge relations
			$C_i(M)\subseteq W_i(M)\times W_{i+1}(M)$, this becomes the
			\emph{relational holonomy} criterion
			\begin{equation}\label{eq:multi-native-cycle}
				\quad\Lambda(M)\ne\varnothing
				\iff \Delta_{W_0(M)}\cap
				(C_0(M);C_1(M);\cdots;C_{k-1}(M))\ne\varnothing.\quad
			\end{equation}
			Here $;$ is ordinary relational composition, not presumed functional
			transport.  For a graph with several cycles, separately feasible
			cycle relations need not have a common global witness.
		\end{enumerate}
	\end{theorem}
	\begin{proof}\leavevmode
		\begin{enumerate}[label=(\roman*),leftmargin=*]
			\item \emph{General descent.} The two componentwise projections of a tuple in \eqref{eq:multi-native-join}
			are admissible for every native formation relation and therefore define
			global realized states.  Local equality and inequality tests yield exact
			agreement $M$, because every descriptor is observed somewhere and shared
			readouts are consistent.  Conversely, any global pair with exact agreement
			$M$ restricts to an element of every paired relation.  The constructions
			are mutually inverse.
			\item \emph{Cycle test.} In a simple cycle, relational composition
			expresses existence of the successive intermediate witnesses, and
			intersecting the composite with the diagonal closes the cycle.  On a
			graph with several cycles, different cycle composites may use different
			intermediate witnesses, so their separate diagonal tests are insufficient.\qedhere
		\end{enumerate}
	\end{proof}
	
	\begin{definition}[Exact witness descent and its local and global obstructions]
		\label{def:witness-globalization-obstruction}
		For a finite multivalued atlas and $M\subsetneq D$, the following cases describe exact agreement.
		\begin{enumerate}[label=(\roman*),leftmargin=*]
			\item \emph{Exact witness descent.}
			We say that \emph{exact-$M$ witness descent holds} if $\Lambda(M)\ne\varnothing$
			in \eqref{eq:multi-native-join}.
			\item \emph{Global obstruction.} If every local paired-record set
			$W_v(M)$ is nonempty but $\Lambda(M)=\varnothing$, we call $M$ a
			\emph{globalization obstruction} for the declared formation relations.
			\item \emph{Local obstruction.} If some $W_v(M)$ is empty, the failure is already a \emph{local witness
				obstruction}. These notions concern synchronized \emph{ordered pairs};
			a single-state extension obstruction, as studied below for linear and
			group diagrams, is a separate lifting question.
		\end{enumerate}
	\end{definition}
	
	\begin{theorem}[Sharp bounds on the number of constraints obstructing a witness]
		\label{thm:multi-native-core}
		Fix $M\subsetneq D$ and protect a subfamily
		$\mathcal F_0\subseteq\mathcal F$.  Suppose
		\[
		\Lambda_0(M)=\{w\in\mathcal P_M:
		w|_{S_\alpha}\in\widetilde F_\alpha(M)
		\text{ for }\alpha\in\mathcal F_0\}\ne\varnothing.
		\]
		\begin{enumerate}[label=(\roman*),leftmargin=*]
			\item \emph{Obstruction bound and sharpness.} If the full atlas has no global witness with exact agreement $M$,
			then already a subfamily of at most
			\begin{equation}\label{eq:multi-core-bound}
				\quad\min\left\{|\mathcal F\setminus\mathcal F_0|,
				\tfrac12|\Lambda_0(M)|\right\}\quad
			\end{equation}
			remaining formation constraints has no compatible witness together with
			$\mathcal F_0$.  This bound is sharp in the finite class, even when all
			formation constraints are binary and the global realized-state space is
			nonempty.
			\item \emph{Tree estimate.} If the protected factors form a spanning tree and each
			oriented tree edge permits at most $b_e$ paired continuations per
			incoming witness, then
			\begin{equation}\label{eq:multi-ambiguity-bound}
				|\Lambda_0(M)|\le |W_o(M)|\prod_{e\in T}b_e
			\end{equation}
			for a root $o$ and orientation away from $o$; in particular this gives
			an explicit ambiguity-dependent obstruction certificate.
		\end{enumerate}
	\end{theorem}
	\begin{proof}\leavevmode
		\begin{enumerate}[label=(\roman*),leftmargin=*]
			\item \emph{Obstruction bound and sharpness.} Simultaneously reversing the two realized-state components at every
			position is an involution on $\Lambda_0(M)$.  Since a position observes
			some descriptor outside $M$, no global paired assignment is fixed by
			this reversal.  Each native paired formation relation is invariant under
			reversal.  Consequently every remaining factor's admissible subset of
			$\Lambda_0(M)$ is a union of two-element reversal orbits.
			If their full intersection is empty, choose a constraint excluding an
			orbit that survives the previously chosen constraints, and continue.
			Each step removes at least one surviving orbit.  There are
			$|\Lambda_0(M)|/2$ orbits, giving the bound; the number of available
			constraints gives the second bound.  The block-flip example in \cref{thm:witness-cycle-core} attains the
			factor $1/2$. Its global state set contains the distinguished state.
			\item \emph{Tree estimate.} Start with a root pair and extend it
			along the oriented edges. At each edge there are at most $b_e$
			choices. Each nonroot vertex has one parent.\qedhere
		\end{enumerate}
	\end{proof}
	
	\paragraph{Native records in the sharpness family.}
	The finite state sets used in \cref{thm:witness-cycle-core}, and
	hence for sharpness in \cref{thm:multi-native-core}, have the
	declared native-record interpretation. At vertex $v$, take
	$M_v=\mathbb C^{Z_v}$ acting by multiplication on $\ell^2(Z_v)$
	and the realized datum $(\mathrm{com},M_v)$ with endpoint
	$\mathfrak P_{M_v'}$ in its access--commutant context.
	Tag each character of $M_v$ by this realized triple, as in
	\cref{prop:spectral-observation-interpretation}.
	On these records the descriptor functions and interface permutations
	induce the stated observations and paired relations.
	
	\subsection{Faithful witness globalization}
	\label{subsec:witness-globalization}
	
	Exact recovery requires a pair with exact agreement set $M$ for each
	proper meet-irreducible $M$ of the target closure.
	
	\begin{lemma}[One-cover certificate for an exact recovery witness]
		\label{lem:meetirreducible-cover-witness}
		Let $X$ be a nonempty finite realized source with coherent
		readouts indexed by a finite set $D$, and let $c_0\le c_X$ be a target
		determination closure.  For a proper meet-irreducible
		$c_0$-closed set $M$, let $M^+$ be its unique cover in the
		finite lattice of $c_0$-closed sets and fix any
		$d_M\in M^+\setminus M$.  An ordered pair of source states
		that agrees on $M$ has \emph{exact} agreement set $M$ if
		and only if it disagrees on $d_M$.
	\end{lemma}
	\begin{proof}
		The agreement set $E(x,x')$ of any two source states
		is $c_X$-closed and hence $c_0$-closed, since every
		implication valid in $X$ is respected by the pair.
		If $E(x,x')$ properly contains $M$, the unique-cover
		property in a finite closure lattice implies
		$M^+\subseteq E(x,x')$.  Thus disagreement on $d_M$
		rules out any proper enlargement of $M$.
		The converse is immediate.
	\end{proof}
	
	\begin{theorem}[Exact criterion for faithful witness globalization]
		\label{thm:witness-globalization}
		Let a finite multivalued native formation atlas be given as in
		\cref{def:multi-native-atlas}, and let $X$ be its nonempty set of global
		admissible realized-record tuples. Assume that all descriptors in the
		finite set $D$ have consistent global readouts. Write $c_X$ for the
		determination closure of those readouts and let $c_0$ be a target closure
		operator on $D$ satisfying $c_0\le c_X$ pointwise. Denote by
		$\mathcal M(c_0)$ the proper meet-irreducible closed sets of $c_0$.
		Then the following conditions are equivalent:
		\begin{enumerate}[label=(\roman*)]
			\item The global source is recovery-faithful: $c_X=c_0$.
			\item Every $M\in\mathcal M(c_0)$ occurs as the \emph{exact} agreement
			set of an ordered pair of global realized states.
			\item For every $M\in\mathcal M(c_0)$, the finite paired-formation
			matching set $\Lambda(M)$ of \eqref{eq:multi-native-join} is nonempty.
		\end{enumerate}
		If these conditions fail, an $M\in\mathcal M(c_0)$ with
		$\Lambda(M)=\varnothing$ certifies the failure. More precisely, either
		$\mathcal P_M=\varnothing$ (the obstruction is already local), or some
		subfamily of at most
		\begin{equation}\label{eq:faithful-globalization-core}
			\min\!\left\{|\mathcal F|,\frac{|\mathcal P_M|}{2}\right\}
		\end{equation}
		formation constraints excludes every candidate in $\mathcal P_M$.
		The bound counts constraints for explicitly listed pairs of records.
	\end{theorem}
	\begin{proof}
		For two global states $x,x'\in X$, let
		$E(x,x')=\{d\in D:d(x)=d(x')\}$. The dependency closure satisfies
		\[
		c_X(U)=\bigcap_{\substack{x,x'\in X\\U\subseteq E(x,x')}}E(x,x'),
		\]
		where the empty intersection is $D$. Each $E(x,x')$ is $c_X$-closed:
		the displayed pair witnesses $d\notin c_X(E(x,x'))$ for every
		$d\notin E(x,x')$. Since $c_0\le c_X$, it is also $c_0$-closed.
		\begin{itemize}
			\item \emph{\textup{(i)} $\implies$ \textup{(ii)}.} If $c_X=c_0$, the displayed intersection represents every proper
			$c_0$-closed set $M$ as an intersection of exact agreement sets
			containing $M$. When $M$ is meet-irreducible, at least one set in this
			finite intersection must equal $M$; hence (i) implies (ii).
			\item \emph{\textup{(ii)} $\implies$ \textup{(i)}.} Conversely, every closed set in a finite closure system is the
			intersection of the meet-irreducible closed sets containing it. If
			(ii) holds and $d\notin c_0(U)$, there is a proper meet-irreducible
			$M\supseteq c_0(U)$ with $d\notin M$, and an actual global pair with
			agreement set $M$. Thus $d\notin c_X(U)$. Together with
			$c_0\le c_X$, this proves (i).
			\item \emph{\textup{(ii)} $\Longleftrightarrow$ \textup{(iii)}.} By \cref{thm:multi-native-exact-descent}, $\Lambda(M)$ is in bijection with the global pairs having agreement set $M$.
			\item \emph{Obstruction certificate.} If (iii) fails, the claimed $M$ is an obstruction. If
			$\mathcal P_M\ne\varnothing$, apply \cref{thm:multi-native-core} with
			no protected formation constraints. Since $M\subsetneq D$ and every
			descriptor is observed at a position, reversal acts freely on
			$\mathcal P_M$, and the stated core bound follows. If
			$\mathcal P_M=\varnothing$, the local obstruction is immediate.\qedhere
		\end{itemize}
	\end{proof}
	
	\begin{corollary}[Recovery faithfulness and relative constraint obstruction ranks]
		\label{cor:synchronized-recovery-transfer}
		Retain the finite multivalued atlas, nonempty global source $X$,
		and target closure $c_0\le c_X$ of
		\cref{thm:witness-globalization}. Fix a protected
		subfamily $\mathcal F_0\subseteq\mathcal F$ containing the declared
		shared-readout consistency constraints, and form the protected
		paired matching set $\Lambda_0(M)$ for each proper meet-irreducible
		$M\in\mathcal M(c_0)$. For
		$J\subseteq\mathcal F\setminus\mathcal F_0$ write
		\[
		\Lambda_J(M)=\{w\in\Lambda_0(M):
		w|_{S_\alpha}\in\widetilde F_\alpha(M)
		\text{ for every }\alpha\in J\},
		\]
		and define the relative \emph{constraint obstruction rank}
		\begin{equation}\label{eq:synchronized-transfer-rank}
			\rho_{\mathcal F_0}(M)=
			\min\{|J|:J\subseteq\mathcal F\setminus\mathcal F_0,
			\ \Lambda_J(M)=\varnothing\},
		\end{equation}
		with minimum $+\infty$ if no such $J$ exists. The following hold.
		\begin{enumerate}[label=(\roman*),leftmargin=*]
			\item \emph{Faithfulness.} \[
			c_X=c_0\quad\Longleftrightarrow\quad
			\rho_{\mathcal F_0}(M)=+\infty
			\quad\text{for every }M\in\mathcal M(c_0).
			\]
			\item \emph{Finite obstruction ranks.} If recovery faithfulness fails, at least one $M$ has finite
			obstruction rank, bounded by
			\begin{equation}\label{eq:synchronized-transfer-bound}
				\rho_{\mathcal F_0}(M)\le
				\min\!\left\{|\mathcal F\setminus\mathcal F_0|,
				\frac{|\Lambda_0(M)|}{2}\right\}.
			\end{equation}
			In particular, $\rho_{\mathcal F_0}(M)=0$ exactly when the
			protected constraints already exclude every candidate for $M$.
		\end{enumerate}
	\end{corollary}
	\begin{proof}\leavevmode
		\begin{enumerate}[label=(\roman*),leftmargin=*]
			\item \emph{Faithfulness.} For a fixed $M$, imposing the entire remaining formation family
			produces exactly $\Lambda(M)$ of
			\eqref{eq:multi-native-join}. Since constraints only delete
			candidates, $\Lambda(M)$ is nonempty if and only if no subfamily
			$J$ makes $\Lambda_J(M)$ empty, equivalently
			$\rho_{\mathcal F_0}(M)=+\infty$. Applying
			\cref{thm:witness-globalization} for every proper
			meet-irreducible gives the first equivalence.
			\item \emph{Finite obstruction ranks.} If $\Lambda_0(M)$ is empty, $J=\varnothing$
			already certifies failure, giving $\rho_{\mathcal F_0}(M)=0$.
			Conversely rank zero means that this empty choice works. For a
			failed $M$ with $\Lambda_0(M)\ne\varnothing$, apply the protected
			form of \cref{thm:multi-native-core} to obtain
			\eqref{eq:synchronized-transfer-bound}.\qedhere
		\end{enumerate}
	\end{proof}
	
	\subsection{Coupled cycles and linear globalization tests}
	\label{subsec:graph-globalization}
	A choice of edge values may pass every individual interface test and
	still fail to extend around two cycles that share an edge. The first
	example exhibits a recovery implication missing from the separate
	chart laws. The graph criterion below decides whether a prescribed
	local value extends to a single compatible assignment. Its obstruction
	lies in a cokernel. A rooted spanning tree bounds its dimension even
	when the two incidence maps have different images. For a finite algebra
	given by tables, \cref{thm:finite-native-cyclic-certificate} also tests
	whether a smaller native source can meet the prescribed contracts.
	
	\begin{example}[Complete individual interfaces do not control interacting cycles]
		\label{ex:interacting-cycle-ports}
		Let a graph have vertices $0,1,2,3$ and oriented edges
		$01,12,20,13,30$.  Its triangles $0120$ and $0130$ share edge $01$,
		so its cycle rank is two.  For every edge $uv$, take a native
		$\mathbb F_q$-module chart with state $(s_u,s_v)\in\mathbb F_q^2$
		and scalar descriptors
		\[
		\{s_u,a_{uv},s_v\},\qquad a_{uv}=s_v-s_u.
		\]
		Identify actual shared vertex readouts.  Distinct incident edge
		charts meet in a single scalar descriptor, and every such interface is
		congruence-complete in the fixed $\mathbb F_q$-module class.
		Every chart has precisely the local triangle closure in which any
		two descriptors determine the third.
		
		The native full join is $\mathbb F_q^4$, with free vertex values
		$(s_0,s_1,s_2,s_3)$.  It projects onto every chart.  The intersection
		of any two chart function algebras is exactly the algebra of their
		shared vertex values, since the four vertex coordinates are
		independent.  Nevertheless the full join satisfies
		\[
		a_{01}+a_{12}+a_{20}=0,\qquad
		a_{01}+a_{13}+a_{30}=0.
		\]
		For example, $\{a_{01},a_{12}\}\to a_{20}$ holds globally.
		It is not generated by the local rules: starting with any set of
		edge labels $a_{uv}$ and no vertex label, each chart has at most one
		known descriptor, so no local triangle implication can fire.
		Thus all individual interface tests and one-state projections pass,
		yet recovery faithfulness fails.  To see that the atlas has no join tree, consider the chart
		nodes containing each of $s_0,s_1,s_2$. They would form three
		pairwise intersecting subtrees, since each pair occurs in an edge
		chart. Pairwise intersecting subtrees of a tree have a common node,
		but no chart contains all three readouts.
	\end{example}
	
	Let each vertex and edge of a finite graph carry a finite-dimensional
	\(\mathbb F_q\)-space \(X_i,Y_e\), with linear incidence maps. Choose
	an orientation and put
	\[
	\delta:\bigoplus_iX_i\to\bigoplus_eY_e,\qquad
	(\delta x)_e=p_{i,e}(x_i)-p_{j,e}(x_j),\qquad
	W=\ker\delta.
	\]
	Fix protected difference subspaces \(K_i^\lambda\le X_i\), including
	\(K_i^{\mathrm{state}}=X_i\), and restrict the coboundary to
	\(\delta_\lambda:\bigoplus_iK_i^\lambda\to
	\bigoplus_e(p_{i,e}K_i^\lambda+p_{j,e}K_j^\lambda)\).
	At a vertex \(r\), let \(\delta_{\lambda,-r}\) be the restriction to
	tuples whose \(r\)-coordinate vanishes, with the same codomain as
	\(\delta_\lambda\). Write \(s_{\lambda,r}\) for the contribution
	from the \(r\)-coordinate alone.

	\begin{theorem}[Linear globalization: exact sequence and obstruction dimension]
		\label{thm:graph-relative}
		\leavevmode
		\begin{enumerate}[label=(\roman*),leftmargin=*]
			\item \emph{Exact sequence and extension criterion.}
			There is an exact sequence
			\[
			\begin{gathered}
				0\to\ker\delta_{\lambda,-r}\to\ker\delta_\lambda
				\xrightarrow{\rho_r}K_r^\lambda
				\xrightarrow{\omega_{\lambda,r}}
				\operatorname{coker}\delta_{\lambda,-r}
				\to\operatorname{coker}\delta_\lambda\to0,
				\\
				\omega_{\lambda,r}(z)=[s_{\lambda,r}z].
			\end{gathered}
			\]
			Consequently a prescribed local state or protected difference
			\(z\in K_r^\lambda\) globalizes if and only if
			\(\omega_{\lambda,r}(z)=0\).
			A nonzero obstruction has a dual functional certificate.
			
			\item \emph{Dimension bound.}
			In a connected graph choose any spanning tree $T$ rooted at $r$.
			For a tree edge $e$ let $v(e)$ be its child, and for each edge put $L_e=p_{i,e}K_i^\lambda+ p_{j,e}K_j^\lambda$.
			Then
			\[
			\dim\operatorname{im}\omega_{\lambda,r}
			\le \dim\operatorname{coker}\delta_{\lambda,-r}
			\le \sum_{e\notin T}\dim L_e+
			\sum_{e\in T}\operatorname{codim}_{L_e}
			(p_{v(e),e}K_{v(e)}^\lambda).
			\]
			If the child map covers $L_e$ on every tree edge, the second sum
			vanishes.
			In particular, matching restricted images on a tree	guarantee that every prescribed root value extends.
		\end{enumerate}
	\end{theorem}

	\begin{proof}
		\leavevmode
		\begin{enumerate}[label=(\roman*),leftmargin=*]
			\item \emph{Exact sequence and extension criterion.} A vector \(z\) extends exactly when
			\(s_{\lambda,r}z\in\operatorname{im}\delta_{\lambda,-r}\),
			which proves exactness at \(K_r^\lambda\). The kernel of
			$\rho_r$ is the zero-root part of $\ker\delta_\lambda$.
			The map between cokernels is onto, and its kernel is
			$(\operatorname{im}\delta_{\lambda,-r}+s_{\lambda,r}K_r^\lambda)/
			\operatorname{im}\delta_{\lambda,-r}$, the image of
			$\omega_{\lambda,r}$. A functional on the cokernel separating
			a nonzero class from zero gives the dual certificate.
			\item \emph{Dimension bound.} For each tree edge choose a complement $J_e$ to the image of the
			child map in $L_e$. Given a tuple of edge values, start at the
			root with value zero and visit children in order of increasing
			depth. After the parent value is fixed, the child value can
			realize any desired tree-edge difference modulo $J_e$, since
			its signed image maps onto $L_e/J_e$. Later assignments leave
			that difference unchanged. Thus every cokernel class has a
			representative in
			$\bigoplus_{e\in T}J_e\oplus\bigoplus_{e\notin T}L_e$.
			Taking dimensions gives the bound.\qedhere
		\end{enumerate}
	\end{proof}

	\paragraph{A sharp tree-edge correction.}
	On a single edge over \(\mathbb F_2\), let both vertex spaces and the edge
	space be \(\mathbb F_2\), with incidence maps \(p_1=\mathrm{id}\) and
	\(p_2=0\), and take \(K_i=X_i\). For \(r=1\), the obstruction map is
	the identity into \(\operatorname{coker}\delta_{-1}\cong\mathbb F_2\).
	Its image has dimension one although there are no chords. The single tree-edge codimension in the bound is one, so the added term is attained.

	The finite Mal'cev analogue replaces linear cokernels by actual
	\emph{support subalgebras}.  On a single cycle with compatible
	subdirect edge relations, compose those relations around the cycle
	to obtain \(H\le A_0\times A_0\).  Fleischer's lemma gives
	onto homomorphisms \(\alpha,\beta:A_0\twoheadrightarrow D\)
	with \(H=\{(x,y):\alpha(x)=\beta(y)\}\).
	Hence
	$
	\pi_0(W)=\operatorname{Eq}(\alpha,\beta).
	$
	The same formula applies to compatible protected-pair algebras. For
	conjugation by a transposition in \(S_3\), the fixed subgroup is nonnormal,
	so this support need not be a congruence class.  In general Mal'cev graphs, distinct
	cycle constraints can share intermediate choices. When the spanning-tree
	transport is deterministic, intersecting the chord equalizers gives
	the root support. General finite instances can instead use the
	exact support relation and the classical Mal'cev-CSP decision method \cite{BulatovDalmau2006}.
	The join-tree criterion appears in \cref{thm:complete-port-native-descent}.
	
	\Cref{app:two-chord-theta-example} gives a two-chord linear
	example in which each chord preserves the free recovery theory
	and their conjunction introduces a four-descriptor circuit.

	\subsubsection{Finite-table certificates for native cyclic sources}
	\label{subsec:finite-native-cyclic-certificate}

	For a matching graph given by explicit finite tables, a complete test
	for a faithful native subsource can enumerate the matching states.
	One Boolean variable per state specifies the source. Closure under
	operations, agreement of descriptor values, and incidence at each
	named cut then give exact constraints on these variables.

	\begin{theorem}[Complete finite-table native globalization and minimum-source test]
		\label{thm:finite-native-cyclic-certificate}
		Take the finite equational matching data of
		\cref{cor:algebraic-arbitrary-graph-cut}, on any finite graph,
		with full matching algebra $\mathbf W$ and finite coherent
		descriptor readouts $(f_d)_{d\in D}$. Fix a target
		determination closure $c_0$. For each
		member $t$ of a specified finite family of cuts, fix a
		separator observation $h_t:W\to Y_t$ that factors through
		both complete side-state projections on $W$.
		For every $w\in W$ introduce a Boolean variable $x_w$ and
		put $R_x=\{w:x_w=1\}$. Conditions \textup{(i)--(iv)} below are
		necessary and sufficient for $R_x$ to be a faithful
		\emph{native subalgebra} preserving all original chart states,
		$c_0$, and the named separator algebras:
		\begin{enumerate}[label=(\roman*)]
			\item \emph{Native subalgebra.}
			$R_x\ne\varnothing$ and, for each basic $n$-ary operation $\omega$ of the fixed finite algebraic signature,
			\begin{equation}\label{eq:finite-csp-native}
				\bigwedge_{j=1}^n x_{w_j}=1
				\quad\Longrightarrow\quad
				x_{\omega^W(w_1,\ldots,w_n)}=1
				\qquad(w_1,\ldots,w_n\in W).
			\end{equation}
			For nullary operations the conclusion is required unconditionally.
			
			\item \emph{Chart coverage.}
			For every chart state $a\in X_i$,
			\begin{equation}\label{eq:finite-csp-coverage}
				\bigvee_{\{w\in W:\pi_i(w)=a\}}x_w=1.
			\end{equation}
			
			\item \emph{Recovery closure.}
			For each pair of distinct rows $w,w'\in W$ whose exact
			agreement set $A(w,w')=\{d:f_d(w)=f_d(w')\}$ is
			not $c_0$-closed, require $\neg(x_w\wedge x_{w'})$.
			For each $S\subseteq D$ and $d\notin c_0(S)$, require
			\begin{equation}\label{eq:finite-csp-recovery}
				\bigvee_{\substack{\{w,w'\}\subseteq W:\\
						f_s(w)=f_s(w')\ (s\in S),\\
						f_d(w)\ne f_d(w')}}(x_w\wedge x_{w'})=1.
			\end{equation}
			
			\item \emph{Named cut algebras.}
			For each named cut $t$ and every value $y\in h_t(R_x)$, the bipartite graph of retained complete-side state pairs with $h_t$-value $y$ is connected.
			Its left and right vertices are the two projections	of the retained rows in that fiber.
			
			\item \emph{Minimum embedded source.}
			If a preservation contract additionally requires full coverage of
			prescribed side values or values of $h_t$, add unary coverage
			conditions of the form \eqref{eq:finite-csp-coverage}. To require
			equality with an exact prescribed image, also require $x_w=0$
			for every row whose value lies outside that image. There is a proper embedded-native-subalgebra solution in
			this fixed class if and only if these finite constraints
			have a Boolean model with $\sum_wx_w<|W|$. Moreover, the
			exact minimum among all such embedded native sources is
			\begin{equation}\label{eq:finite-csp-minimum}
				\min\left\{\sum_{w\in W}x_w:
				\textup{\emph{(i)--(iv)}}\text{ and all prescribed
					additional coverage conditions hold}\right\},
			\end{equation}
			with minimum $+\infty$ if there is no model. Each minimizing Boolean model gives one minimum
			embedded source, and its indicator vector determines it uniquely.
			
			\item \emph{Minimum over compatible chart maps.}
			The same minimum holds over all nonempty finite algebras in the fixed variety equipped with compatible homomorphic chart maps, when the protected readouts factor through their joint map to $W$ and the contract consists of the stated tests.
			Each minimizer in this class is isomorphic over the charts to a minimizing embedded source.
		\end{enumerate}
	\end{theorem}

	\begin{proof}
		The natural join can be enumerated from the supplied finite
		tables; if it is empty, no nonempty source exists.
		\begin{enumerate}[label=(\roman*),leftmargin=*]
			\item \emph{Native subalgebra.} Condition
			\textup{(i)} is exactly nonempty closure under the prescribed
			operations, including constants, so it is equivalent to being
			a subalgebra of $\mathbf W$. The finite table presentation
			then supplies the native realization in the fixed context
			by \cref{lem:finite-table-presentation}.
			\item \emph{Chart coverage.} Condition
			\textup{(ii)} is exactly surjectivity onto each chart.
			\item \emph{Recovery closure.} A positive target implication fails on two retained rows exactly
			when their agreement set omits some member of its $c_0$-closure.
			The forbidden pairs in \textup{(iii)} therefore enforce all
			positive implications. This test is automatic if $c_0\le c_W$.
			A negative implication
			$S\not\to d$ survives precisely when two retained states
			agree on $S$ but disagree on $d$, which is
			\textup{(iii)}. Therefore \textup{(iii)} is equivalent to
			$c_{R_x}=c_0$.
			\item \emph{Named cut algebras.} As each $h_t$ factors through both sides
			on $W$, it also does on $R_x$, and
			\cref{thm:proper-cut-incidence} makes \textup{(iv)}
			equivalent to the exact named common-information algebra.
			\item \emph{Minimum embedded source.} The optional unary conditions cover the prescribed values and,
			when equality is required, exclude precisely the out-of-range rows. All constraints concern
			the chosen subset of the fixed $W$; conversely every
			permitted embedded native subsource gives precisely its
			indicator vector satisfying them. Minimizing that indicator
			sum gives the optimum over embedded sources.
			\item \emph{Minimum over compatible chart maps.} For a finite source $H$ in the stated variety, the compatible
			chart maps give a homomorphism $H\to W$. Its image $R$ is
			a subalgebra and $H\to R$ is onto. Thus $R$ has exactly the
			same chart images, relation of descriptor values and cut-side
			incidences as $H$. It meets the same tests and $|R|\le |H|$.
			Conversely every feasible $R$ supplies such a source by its
			chart projections. At the minimum the surjection $H\to R$
			is bijective, hence is an algebra isomorphism commuting with
			all chart maps.\qedhere
		\end{enumerate}
	\end{proof}

	\subsection{Native globalization with compatible and crossing cuts}

	\Cref{thm:compatible-cut-native-globalization} constructs the
	greatest faithful source with a given boundary image at each of a
	family of compatible cuts. For crossing cuts, one complete side
	image per cut gives the corresponding construction in
	\cref{thm:oriented-side-native-globalization}.

	\paragraph{Compatible named cuts.}
	Two bipartitions are \emph{compatible} if at least one of the four
	intersections of a side of one with a side of the other is empty.
	In particular, the edge splits of a tree are pairwise compatible.

	\begin{theorem}[Greatest sources with prescribed boundaries at compatible cuts]
		\label{thm:compatible-cut-native-globalization}
		Take the finite equational matching data of
		\cref{cor:algebraic-arbitrary-graph-cut} on any finite graph with
		nonempty full matching algebra $W$. Fix coherent descriptor
		readouts and a target closure $c_0\le c_W$. For each cut $t$ in a
		pairwise compatible family, write $b_t:W\to Y_t$ for the complete
		boundary map onto its reachable image, and let
		$h_t=\rho_t\circ b_t$ be the named observation.
		For each family of nonempty boundary subalgebras $D_t\le Y_t$
		on which $\rho_t$ is injective, take all rows with those allowed
		boundary values.
		\[
		W_D=\{w\in W:b_t(w)\in D_t\text{ for every }t\}.
		\]
		\begin{enumerate}[label=(\roman*),leftmargin=*]
			\item \emph{Existence criterion.} There exists a native subalgebra of $W$ with boundary images exactly
			$D_t$ that preserves all original chart states, has closure $c_0$, and has
			exact named common information at every cut if and only if $W_D$
			is nonempty, reaches each $D_t$, projects onto each chart, and
			has $c_{W_D}=c_0$.
			\item \emph{Greatest source and coverage.} In that case $W_D$ is the greatest faithful
			native subalgebra with those boundary images. Every faithful
			native subalgebra has such a boundary profile. Specified coverage
			of side-state values can be tested on $W_D$ as well, and the
			named observation ranges are exactly $\rho_t(D_t)$.
		\end{enumerate}
	\end{theorem}

	\begin{proof}
		\leavevmode
		\begin{enumerate}[label=(\roman*),leftmargin=*]
			\item \emph{Existence criterion.} Let $U$ be any faithful native subalgebra with the stated chart,
			recovery and cut contracts. Then $D_t=b_t(U)$ is a subalgebra.
			Every function of $b_t|_U$ factors through both complete sides.
			If two different values in $D_t$ had the same $\rho_t$-value,
			a function separating them would violate exact named common
			information on $U$. Thus $\rho_t$ is injective on $D_t$. Since
			$U\subseteq W_D\subseteq W$, the source $W_D$ reaches $D_t$ and
			retains chart coverage and any required side-state values already
			attained on $U$. Determination is antitone under inclusion of
			sources, so $c_0\le c_W\le c_{W_D}\le c_U=c_0$.
			
			Conversely, suppose $W_D$ meets the stated conditions. The preimage
			of each subalgebra $D_t$ under the homomorphism $b_t$ is closed
			under the native operations, so $W_D$ has a finite native
			presentation by \cref{lem:finite-table-presentation}. Fix a cut
			$t$ and a reached boundary value $y\in D_t$. The fiber
			$b_t^{-1}(y)$ in $W$ is a full rectangle of complete side states
			by \cref{thm:arbitrary-graph-cut-information}. For any other
			cut $s$, compatibility places one entire side of $s$ in one
			side of $t$. The boundary $b_s$ is readable from that side of
			$s$, so the restriction $b_s\in D_s$ in the rectangle depends
			on only one of its two factors. After all such restrictions the
			fiber in $W_D$ is still a full nonempty rectangle.
			A function of both complete sides is consequently constant on
			each such fiber, and any function of the boundary is common to
			both sides. The common algebra on $W_D$ is therefore precisely
			$\mathbb C[b_t|_{W_D}]=\mathbb C[h_t|_{W_D}]$, by injectivity
			of $\rho_t$ on $D_t$. Chart coverage and recovery were assumed.
			\item \emph{Greatest source and coverage.} Every source with that boundary profile is contained in $W_D$,
			so $W_D$ is greatest. Each faithful source has the profile obtained from its boundary images as above. Required side values attained by such a source remain attained in $W_D$, and $h_t(W_D)=\rho_t(D_t)$.\qedhere
		\end{enumerate}
	\end{proof}

	For matrix presentations of finite modules, the boundary inverse images
	and chart projections are row-space or
	intersection calculations. A linear $\rho_t$ whose full image must
	be preserved restricts the search at that boundary to complements
	in $Y_t$ of $\ker\rho_t$.

	\begin{example}[A proper native source restores a named separator]
		\label{ex:compatible-cut-boundary-trim}
		Take a four-cycle $1,2,3,4,1$ with one copy of $\mathbb F_2$ per
		edge interface. At each vertex take $\mathbb F_2^2$ with the two
		incident edge coordinates as its interfaces. The full matching space
		is $\mathbb F_2^4$, with coordinates $x_{12},x_{23},x_{34},x_{41}$.
		Protect the two descriptors $x_{12}$ and $x_{34}$ and their free
		closure. For the cut $\{1,2\}\mid\{3,4\}$ the complete boundary
		is $(x_{23},x_{41})$ and the named observation is $x_{23}$.
		The full matching space has two independent common coordinates.
		The diagonal boundary $x_{23}=x_{41}$ gives a three-dimensional
		saturation with the exact named separator and its full range.
		No vertex is incident to both crossing edges, so every local
		state survives. The two protected descriptors remain independent.
	\end{example}

	\begin{example}[Crossing cuts can invalidate boundary saturation]
		\label{ex:crossing-cut-boundary-saturation}
		On the same four-cycle, use nontrivial independent $\mathbb F_2$
		interfaces only on edges $12$ and $34$. Each chart has one nontrivial
		coordinate, and the full matching space has two independent
		coordinates $x_{12},x_{34}$. At $t=\{1,2\}\mid\{3,4\}$ the complete
		boundary is a singleton and the named separator is constant. At
		$s=\{1,3\}\mid\{2,4\}$ the complete boundary is
		$(x_{12},x_{34})$, and the named separator is $x_{12}$. Choose the
		diagonal profile at $s$ and the singleton at $t$. The saturation
		covers all chart states and has exact separator information at $s$.
		At $t$ both sides read the same nonconstant coordinate, so its
		named separator is inexact. In fact every chart-exact subsource
		with this boundary profile has this defect at $t$.
	\end{example}

	\paragraph{Crossing cuts with one prescribed side per cut.}
	For crossing cuts, record one complete side image at each cut. The
	required monotonicity holds for every finite relation.

	\begin{lemma}[One-sided stability of exact common information]
		\label{lem:one-sided-cut-stability}
		Let $U\subseteq R$ be nonempty finite sets of complete state tuples,
		with side readouts $\ell,r$ and a common readout $h$ factoring through
		both sides on $R$. Write $M_\ell$ and $M_r$ for the algebras of complex
		functions factoring through the indicated readouts. If $\ell(U)=\ell(R)$ and
		\[
		M_\ell(U)\cap M_r(U)=\mathbb C[h|_U],
		\]
		then
		\[
		M_\ell(R)\cap M_r(R)=\mathbb C[h|_R].
		\]
	\end{lemma}
	
	\begin{proof}
		Fix $y\in h(R)$. Because $h$ factors through $\ell$ and
		$\ell(U)=\ell(R)$, the $y$-fiber is already reached in $U$.
		Its incidence graph for $U$ is connected by
		\cref{thm:proper-cut-incidence}. Every edge of the corresponding
		graph for $R$ has a left vertex attained by some state of $U$;
		that state has the same $h$-value. Thus every added edge, and
		any new right vertex incident to it, attaches to the connected
		$U$-graph in this fiber. The entire $R$-graph is connected.
		Apply \cref{thm:proper-cut-incidence} again.
	\end{proof}
	
	\begin{theorem}[Faithful native globalization with prescribed side images at crossing cuts]
		\label{thm:oriented-side-native-globalization}
		Take the finite equational matching data of
		\cref{cor:algebraic-arbitrary-graph-cut}, with nonempty full matching
		algebra $W$, coherent descriptor readouts, and target closure
		$c_0\le c_W$. Let $\mathcal T$ be any finite family of named cuts,
		with no compatibility assumption. Orient each cut as
		$I=L_t\sqcup R_t$, put $\ell_t=\pi_{L_t}$ and $r_t=\pi_{R_t}$,
		and let $h_t$ be its named observation, factoring through both
		complete sides on $W$.		
		For each $t$, prescribe a nonempty subalgebra
		$A_t\le\ell_t(W)$ and form the oriented-side saturation
		\begin{equation}\label{eq:oriented-side-saturation}
			W_A=\{w\in W:\ell_t(w)\in A_t\text{ for every }t\in\mathcal T\}.
		\end{equation}
		There exists a native subalgebra $U\le W$ with
		$\ell_t(U)=A_t$ for every $t$, preserving all original chart states,
		$c_0$, and all named common-information algebras, if and only if:
		\begin{enumerate}[label=(\roman*)]
			\item \emph{Chart and side coverage.}
			$W_A\ne\varnothing$, $\ell_t(W_A)=A_t$ for every $t$, and
			$\pi_i(W_A)=X_i$ for every original chart $i$;
			
			\item \emph{Recovery closure.}
			$c_{W_A}=c_0$;
			
			\item \emph{Named common information.}
			For every $t$ and $y\in h_t(W_A)$, the incidence graph of complete side pairs in $W_A\cap h_t^{-1}(y)$ is connected.
		\end{enumerate}
		When these conditions hold, $W_A$ is the greatest faithful embedded
		native source with the prescribed oriented side images.
		Every faithful embedded native source is contained in such a greatest
		source, obtained by taking $A_t=\ell_t(U)$. Thus existence for arbitrary
		crossing cuts is characterized exactly by the successful oriented-side
		profiles. Additional required coverage of specified side values or named
		observation values can be tested directly on $W_A$.
	\end{theorem}

	\begin{proof}
		Every $\ell_t$ is a homomorphism, so $W_A$, when nonempty, is a
		subalgebra of $W$ and has the required finite native presentation.
		\begin{itemize}
			\item \emph{Necessity of (i)--(iii).} Suppose first that $U$ exists. Its side images are subalgebras,
			and $U\subseteq W_A\subseteq W$ by construction. Hence all
			required chart states and side-image values are still attained in
			$W_A$, while $\ell_t(W_A)\subseteq A_t$ gives equality. Inclusion
			of sources reverses determination, so
			\[
			c_0\le c_W\le c_{W_A}\le c_U=c_0.
			\]
			For each cut, $\ell_t(U)=\ell_t(W_A)$, and
			\cref{lem:one-sided-cut-stability} transfers exact common information
			from $U$ to $W_A$. This proves all three conditions.
			
			\item \emph{Sufficiency and greatestness.} Conversely,
			\textup{(i)--(iii)} state precisely the chart, recovery, and cut
			contracts for the native source $W_A$, using
			\cref{thm:proper-cut-incidence} for the last clause. Every source
			with the prescribed side images lies in $W_A$, which proves
			greatestness. Coverage requirements are retained under enlargement
			from $U$ to $W_A$ and can be imposed on the latter directly.\qedhere
		\end{itemize}
	\end{proof}

	\begin{corollary}[Recovery and cut obstructions inherited by subformations]
		\label{cor:oriented-side-obstructions}
		\label{app:full-side-hereditary}
		\begin{enumerate}[label=\textup{(\roman*)},leftmargin=*]
			\item \emph{Recovery under deletion.}
			If $U\subseteq R$ are nonempty finite sources with the same
			descriptor readouts, then $c_R\subseteq c_U$.
			In particular, if $c_0\subseteq c_R$ and $c_R\ne c_0$,
			then $c_U\ne c_0$.
	
			\item \emph{Inherited cut obstruction.}
			Let $\ell,r$ be side readouts on a nonempty finite source $R$,
			and let $h$ factor through both. If
			\[
			\mathbb C[\ell]\cap\mathbb C[r]\supsetneq\mathbb C[h]
			\quad\text{on }R,
			\]
			then the same strict inclusion holds on every nonempty
			$U\subseteq R$ retaining the full image of at least one side.
			More generally, it holds whenever $U$ retains a left vertex
			from each of two distinct components of the incidence graph
			in a single $h$-fiber.
			In particular, if the cut condition of
			\cref{thm:oriented-side-native-globalization}\textup{(iii)}
			fails on $W_A$ and $\ell_t(W_A)=A_t$, every subrelation
			$U\subseteq W_A$ with $\ell_t(U)=A_t$ inherits that failure.
	
			\item \emph{Finite cut certificate.}
			For explicit finite tables, a partition of the left vertices
			in one $h$-fiber into two nonempty unions of components,
			with no right vertex incident to both parts, certifies the
			cut obstruction. Verification uses the incidence table
			of that fiber.
		\end{enumerate}
	\end{corollary}
	
	\begin{proof}\leavevmode
		\begin{enumerate}[label=\textup{(\roman*)},leftmargin=*]
			\item \emph{Recovery under deletion.}
			Every pair of states in $U$ is a pair in $R$, so every
			functional dependency holding on $R$ also holds on $U$.
	
			\item \emph{Inherited cut obstruction.}
			If one full side image is retained, exact common information
			on $U$ would imply exact common information on $R$ by
			\cref{lem:one-sided-cut-stability}, with the sides interchanged
			if necessary. This contradicts the hypothesis.
			For the more general assertion, retaining a left vertex
			retains an incident edge and hence a state in its component.
			Deleting edges cannot join distinct components.
			A subrelation meeting two components in the same $h$-fiber
			therefore fails \cref{thm:proper-cut-incidence}.
			Retaining the complete left image forces both components
			to be met, proving the assertion for $W_A$.
	
			\item \emph{Finite cut certificate.}
			The stated partition is a separation certificate.
			Checking each edge verifies that no path joins its two parts.
			\qedhere
		\end{enumerate}
	\end{proof}

	\begin{example}[Boundary profiles can admit a faithful source while their saturation fails]
		\label{ex:crossing-cut-faithful-proper-source}
		Take six charts $X_i=\mathbb F_2$ with identity interfaces on three
		disjoint matching edges $12$, $34$, and $56$. Thus
		$W=\mathbb F_2^3$ has coordinates $(a,c,h)$ and chart tuple
		$(a,a,c,c,h,h)$. Use the crossing cuts
		\[
		t:\{1,2,5\}\mid\{3,4,6\},\qquad
		s:\{1,3,5,6\}\mid\{2,4\}.
		\]
		Their complete boundaries are respectively $b_t=h$ and
		$b_s=(a,c)$. Name the observations $h_t=h$ and $h_s=a$.
		Choose $D_t=\mathbb F_2$ and
		$D_s=\{(a,a):a\in\mathbb F_2\}$; the named observations are
		injective on these respective boundary images. The boundary
		saturation is
		\[
		W_D=\{(a,a,h):a,h\in\mathbb F_2\}.
		\]
		It preserves every chart state, but both complete sides at $t$
		read $(a,h)$, so their common algebra strictly contains
		$\mathbb C[h]$. Nevertheless the proper native subspace
		\[
		U=\{(a,a,a):a\in\mathbb F_2\}
		\]
		has the same two boundary images, preserves every chart state,
		and has exact named common information at both cuts. If the
		protected descriptor family consists of $h$, its determination
		closure on $U$ equals that on $W$. Thus even existence of a
		faithful source with a boundary profile does not make boundary
		saturation faithful for crossing cuts. Taking $A_t=\ell_t(U)$ and $A_s=\ell_s(U)$ gives $W_A=U$.
	\end{example}

	Given an oriented-side profile, test its cuts by incidence
	connectivity, using at most $|W_A|$ state edges per cut. For scalar
	linear descriptors, linear named observations, and target $c_0=c_W$,
	compute \eqref{eq:oriented-side-saturation}, check the required
	images, and apply the row-space tests of
	\cref{thm:arbitrary-linear-constraints}. A loss of rank makes a
	descriptor basis dependent. Retention of rank preserves every span
	relation, hence the scalar closures agree. For a supplied profile,
	these checks take polynomial time in the matrix input size.

	\begin{proposition}[Complexity of faithful-source existence from global and local tables]
		\label{prop:local-table-profile-hardness}
		\leavevmode
		\begin{enumerate}[label=(\roman*),leftmargin=*]
			\item \emph{Explicit matching table.} For an explicitly listed matching algebra $W$, operation tables, chart
			maps, readouts, named cuts, and a list of all $c_0$-closed sets with
			$c_0\le c_W$, the existence of a successful oriented-side profile belongs
			to $\mathrm{NP}$.
			
			\item \emph{Local chart tables.} When only the local charts and incidence maps are given by tables, the
			following restricted decision problem is $\mathrm{NP}$-complete. The
			algebraic signature is empty, there is one constant descriptor with its
			constant closure, and the only named cut separates a leaf from the other
			charts and records the complete leaf state via an identity interface.
			The question is whether a native source covers every chart state and
			satisfies this cut and recovery contract. Hardness persists with a
			specified matching tuple and at most seven states per chart.
		\end{enumerate}
	\end{proposition}

	\begin{proof}
		\leavevmode
		\begin{enumerate}[label=(\roman*),leftmargin=*]
			\item \emph{Explicit matching table.} Guess the indicator of a nonempty subset
			$U\subseteq W$. Use the operation and projection tables to check that
			$U$ is a subalgebra and covers every chart. At each named cut test
			connectivity in each reached separator fiber. Since $c_0\le c_W\le c_U$,
			equality $c_U=c_0$ is equivalent to every listed $c_0$-closed set $K$
			being $c_U$-closed. For each $d\notin K$, scan pairs in $U$ for two
			rows agreeing on $K$ and differing on $d$. All tests take polynomial
			time in the explicit input length. The profile of $U$ consists of its
			left-side images by \cref{thm:oriented-side-native-globalization}.
			\item \emph{Local chart tables.} For the local-table problem a certificate consists of one global matching
			tuple extending each state of each chart. Its size and verification time
			are polynomial in the local tables. Their union is a subalgebra in the
			empty signature. The leaf reads the named separator and the other side
			reads the same value, so the common algebra is exactly the leaf algebra.
			Every nonempty source has the prescribed constant-descriptor closure.
			
			Reduce from three-colorability of a finite graph
			\cite[Problem~GT4]{GareyJohnson1979}. We may take the graph connected:
			join successive components by bridges, permuting the colors within
			each new component to accommodate its bridge. For each graph vertex
			use a chart with states $\bot,1,2,3$. For each oriented graph edge
			use a chart with states $(\bot,\bot)$ and the six ordered pairs of
			distinct colors. Connect this chart to the endpoint charts by its two
			coordinate projections, matched against their identity maps.
			Attach a copy of one vertex chart as a leaf with an identity interface
			and name its state as the separator. The all-$\bot$ assignment is
			always a global tuple. Connectedness forces every other global tuple
			to use colors everywhere, and those tuples are exactly proper
			three-colorings, with their edge and leaf states determined.
			If a coloring exists, permutations of its three colors reach each
			color at every vertex and each ordered pair of distinct colors at
			each edge. Together with the all-$\bot$ tuple these give full chart
			coverage. Conversely, extending any colored vertex state gives a
			three-coloring. The construction has polynomial size.\qedhere
		\end{enumerate}
	\end{proof}

	\subsection{Linear cut tests and the forest boundary}

	Restricting linear readouts to a matching subspace amounts to
	taking their row spaces modulo its annihilator.
	\Cref{thm:arbitrary-linear-constraints} uses this observation to
	test coverage, recovery, and common information at a cut.
	\Cref{prop:forest-automatic-faithfulness} then identifies the
	finite simple graphs on which all binary linear edge-chart
	atlases preserve the local recovery rules.
	
	\begin{theorem}[Dual preservation criteria for homogeneous linear constraints]
		\label{thm:arbitrary-linear-constraints}
		Let $T$ be a finite-dimensional $\mathbb F_q$-space with original-chart
		surjections $\pi_i:T\twoheadrightarrow A_i$ and protected scalar
		descriptor forms $f_d\in T^*$, with $c_T=c_0$. Let
		$
		a:T\longrightarrow E
		$
		be any linear map to a finite-dimensional vector space, and put
		\[
		W=\ker a,\qquad C=\operatorname{im}a^*=W^\perp,
		\qquad V_i=\operatorname{im}\pi_i^*,\qquad
		V=\operatorname{span}\{f_d:d\in D\}.
		\]
		For each named cut, let $\ell:T\to L_0$ and $r:T\to R_0$ be its complete assembled-side linear readouts and let the preassigned linear separator $h_W:W\to H$ factor through both $\ell|_W$ and $r|_W$.
		For any linear extension $\widetilde h:T\to H$ set
		$U_L=\operatorname{im}\ell^*$, $U_R=\operatorname{im}r^*$, and
		$U_h=\operatorname{im}\widetilde h^*$.
		The following criteria hold.
		\begin{enumerate}[label=(\roman*),leftmargin=*]
			\item \emph{Chart coverage.}
			\begin{equation}\label{eq:arbitrary-linear-S}
				\pi_i(W)=A_i\ \text{for every } i
				\quad\Longleftrightarrow\quad V_i\cap C=0\ \text{for every } i.
			\end{equation}
			\item \emph{Recovery.}
			\begin{equation}\label{eq:arbitrary-linear-F}
				c_W=c_0\quad\Longleftrightarrow\quad V\cap C=0.
			\end{equation}
			
			\item \emph{Separator algebras and side images.}
			At each named cut,
			\begin{equation}\label{eq:arbitrary-linear-E}
				M_L(W)\cap M_R(W)=\mathbb C[h_W]
				\quad\Longleftrightarrow\quad (U_L+C)\cap(U_R+C)=U_h+C.
			\end{equation}
			Here $M_L(W),M_R(W)$ are the multiplication algebras of the
			complete linear side readouts restricted to $W$.
			The full assembled-side images are preserved,
			$\ell(W)=\ell(T)$ and $r(W)=r(T)$, precisely when
			\begin{equation}\label{eq:arbitrary-linear-side-images}
				U_L\cap C=U_R\cap C=0.
			\end{equation}
			
			\item \emph{Native realization.} The subspace $W$ has a native finite-module table
			presentation in the fixed formation law. Hence the conjunction of
			the displayed conditions at all protected charts and cuts is an
			exact $\mathsf S+\mathsf F+\mathsf E+\mathsf N$ test, with
			\eqref{eq:arbitrary-linear-side-images} added when complete
			assembled-side coverage is required.
			
			\item \emph{Successive constraints.} If $a$ is assembled from any finite collection of maps
			$a_j:T\to E_j$, then
			\begin{equation}\label{eq:arbitrary-linear-constraint-sum}
				W=\bigcap_{j=1}^m\ker a_j,\qquad
				C=\sum_{j=1}^m\operatorname{im}a_j^*,
			\end{equation}
			with the exact update when adding a constraint
			$C_{k+1}=C_k+\operatorname{im}a_{k+1}^*$.
			All the tests are polynomial-time row space calculations for
			explicitly supplied matrices and a listed family of cuts.
		\end{enumerate}
	\end{theorem}
	
	\begin{proof}
		\leavevmode
		\begin{enumerate}[label=\textup{(\roman*)},leftmargin=*]
			\item \emph{Chart coverage.}
			Restriction of forms gives an onto map $T^*\to W^*$ whose
			kernel is $W^\perp=C$. The dual of $\pi_i|_W$ is injective
			if and only if $\pi_i|_W$ is surjective; equivalently its
			row space loses no vector modulo $C$. This proves
			\eqref{eq:arbitrary-linear-S}.
			
			\item \emph{Recovery.}
			For scalar forms, functional dependence is exactly linear-span
			membership. If $V\cap C=0$, restriction is injective on $V$,
			so every span implication and nonimplication survives.
			Conversely, if $V\cap C\ne0$, take a subset of the protected
			descriptor forms that is a basis of $V$. Their restrictions
			become dependent, and at least one formerly independent form
			is recovered from the other basis descriptors on $W$.
			This implication was absent on $T$, proving necessity in
			\eqref{eq:arbitrary-linear-F}.
			
			\item \emph{Separator algebras and side images.}
			Apply the same argument with the surjections $T\to\ell(T)$ and $T\to r(T)$. Their row spaces are
			$U_L$ and $U_R$ even if the original codomains are larger than
			their images, giving \eqref{eq:arbitrary-linear-side-images}.
			
			The restricted side rowspaces in $W^*=T^*/C$ are
			$(U_L+C)/C$ and $(U_R+C)/C$. Functions of the respective
			linear side readouts coincide precisely on the quotient by
			$\ker(\ell|_W)+\ker(r|_W)$, whose dual is the intersection
			of these rowspaces: invariance under translations by each kernel
			is equivalent to invariance under translations by their sum.
			As $h_W$ factors through both sides, its
			function algebra is the common algebra exactly when its
			row space is this intersection. This is
			\eqref{eq:arbitrary-linear-E}. Different extensions of $h_W$
			have the same row space modulo $C$, so the test is intrinsic.
			
			\item \emph{Native realization.}
			Finally, $W$ is closed under the original module operations,
			so \cref{lem:finite-table-presentation} gives a native source.
			
			\item \emph{Successive constraints.}
			The annihilator of the intersection of kernels is the sum of
			their annihilators, giving \eqref{eq:arbitrary-linear-constraint-sum} and the update rule.\qedhere
		\end{enumerate}
	\end{proof}

	\subsubsection{Mixed obstruction directions for two linear constraints}
	\label{app:joint-dual}

	Retain the source $T$, protected descriptor space $V$, and chart
	readout spaces $V_i$ of \cref{thm:arbitrary-linear-constraints}.
	Two independent scalar chord equations can create a joint
	obstruction even when each equation separately preserves a readout space.

	\begin{corollary}[Unique mixed rank-two obstruction direction]
		\label{app:joint-mixed}
		Suppose $a,b\in T^*$ are independent scalar chord equations and a	protected readout space $U\le T^*$ satisfies
		$U\cap\langle a\rangle=U\cap\langle b\rangle=0$.
		Then $U\cap\langle a,b\rangle$ is zero or one-dimensional. If it
		is nonzero, it is the line spanned by a vector
		$\lambda a+\mu b$ with $\lambda\mu\ne0$. Exactly $q-1$
		projective mixed directions are possible as $U$ varies. Taking $U=V$ detects
		joint recovery loss, while $U=V_i$ detects joint chart-state loss.
	\end{corollary}
	
	\begin{proof}
		The two unmixed lines are disjoint from $U$, so the intersection
		cannot be the whole two-space. Every nonzero intersection vector
		uses both coefficients. The projective line has $q+1$ directions,
		two unmixed and $q-1$ mixed. Each mixed direction is attained
		by taking $U$ equal to that line.
	\end{proof}
	
	In \cref{app:two-chord-theta-example}, the binary theta graph has one mixed recovery direction.
	The protected descriptors fail the joint-separation hypothesis of \cref{thm:simple-diagram-source-rigidity}.

	\subsubsection{Forests as the sharp universal boundary}
	\label{subsec:automatic-faithfulness-graph-boundary}
	
	Even with complete interfaces and surjective chart projections, a
	cycle may introduce a global dependency. Binary linear edge charts
	suffice to characterize the graphs on which faithfulness is automatic.
	
	Let $G=(V,E)$ be a nonempty finite simple graph.
	An edge-chart atlas on $G$ has one finite-dimensional $\mathbb F_2$-space $A_v$ at each vertex.
	Its descriptor set is $B_v=\{e\in E:v\in e\}\sqcup P_v$, where the private sets $P_v$
	are disjoint from each other and from $E$.
	Each edge descriptor is an onto linear map to $\mathbb F_2$; private readouts are linear maps to finite-dimensional spaces.
	All readouts at each chart jointly separate its states.
	Thus every shared interface is a complete singleton signature on $\mathbb F_2$.
	Let $W$ be the full matching space and $c_0$ the closure generated by the local determination theories.
	
	\begin{proposition}[Forests are exactly the universally faithful edge-chart graphs]
		\label{prop:forest-automatic-faithfulness}
		For such a graph $G$, the following are equivalent:
		\begin{enumerate}[label=(\roman*),nosep]
			\item $G$ is a forest.
			
			\item Every edge-chart atlas on $G$ has onto local projections
			$W\to A_v$ and satisfies $c_W=c_0$.
			
			\item Every edge-chart atlas on $G$ whose local projections are onto
			satisfies $c_W=c_0$.
		\end{enumerate}
		If $G$ contains a simple cycle $C$, failure is witnessed with $A_v=\mathbb
		F_2^{E(v)}$, $W\cong\mathbb F_2^E$, coordinate edge readouts, and one scalar private
		readout at each vertex of $C$.
		All interfaces are nontrivial, complete, and onto, and every local state extends to the native matching source.
	\end{proposition}
	
	\begin{proof}\leavevmode
		\begin{itemize}
			\item \emph{\textup{(i)} $\implies$ \textup{(ii)}.} If $G$ is a forest, each component is a join tree and its singleton
			linear interfaces are congruence-complete. Apply
			\cref{thm:complete-port-native-descent} componentwise. The matching
			source of the whole forest is the product of its nonempty component
			sources. A counterexample pair in one component extends by using
			identical states in every other component, so the product adds no
			cross-component dependencies. This proves \textup{(ii)}.
			\item \emph{\textup{(ii)} $\implies$ \textup{(iii)}.} Restrict the assertion in \textup{(ii)} to atlases whose chart projections are onto.
			\item \emph{\textup{(iii)} $\implies$ \textup{(i)}.} We prove the contrapositive. Choose a simple cycle $C$ and use the chart spaces
			in the statement, where $E(v)$ is the set of edges incident to $v$.
			For $v\in V(C)$ with its two cycle edges $e_v^-,e_v^+$, add the
			private scalar readout
			\[
			a_v=x_{e_v^-}+x_{e_v^+}.
			\]
			There are no other private readouts. The shared edge coordinates
			separate each chart, and matching identifies exactly the two copies
			of each edge variable. Consequently $W=\mathbb F_2^E$ and every
			projection onto $A_v$ is onto. On $W$ one has
			\begin{equation}\label{eq:cycle-private-parity-obstruction}
				\sum_{v\in V(C)}a_v=0,
			\end{equation}
			since each cycle edge occurs twice. Fix $w\in V(C)$ and let
			$T=\{a_v:v\in V(C)\setminus\{w\}\}$.
			At a chart $v\ne w$ on the cycle, knowledge of the one sum $a_v$
			does not determine either summand, and does not determine any other
			incident edge variable. At $w$ and at all other vertices there are
			no premises and no constant edge readouts. Thus $T$ is closed under
			every local theory, so $a_w\notin c_0(T)$. Equation
			\eqref{eq:cycle-private-parity-obstruction} gives $a_w\in c_W(T)$,
			contradicting \textup{(iii)}. All spaces and maps have the fixed native
			presentation realizations of \cref{lem:finite-table-presentation}.\qedhere
		\end{itemize}
	\end{proof}

	The forest direction uses classical join-tree consistency \cite{Beeri1983}.

	\subsection{Group-native chord criteria for globalization}

	For groups, test all the chord constraints in one action.
	The orbit of the identity gives an exact extension test and bounds the index of the subgroup of root states that extend.

	Fix a connected finite graph, finite groups $A_i,B_e$, and homomorphisms
	$p_{i,e}:A_i\to B_e$. Choose a spanning tree $T$ with equal incidence
	images $L_e$ on every tree edge. The two images on a chord may differ. Let $S$ be the tree matching group and $W\le S$ the full
	graph matching group. For a root $r$, the projection $\rho:S\to A_r$
	is onto; write $N_r=\ker\rho$.

	On $Q=\prod_{e\notin T}B_e$ define a right $S$-action by
	$
	(d^s)_e=p_{i,e}(s_i)^{-1}d_ep_{j,e}(s_j).
	$
	Let $O=1^S$, $O_r=1^{N_r}$ and $R_T=\prod_{e\notin T}|B_e|$.

	\begin{theorem}[Joint-cycle extension criterion and index bound for finite group atlases]
		\label{thm:joint-chord-witness}
		\leavevmode
		\begin{enumerate}[label=(\roman*),leftmargin=*]
			\item \emph{Source cardinality.}
			\[
			W=\operatorname{Stab}_S(1),\qquad
			|W|=\frac{|S|}{|O|},\qquad |S|=\frac{\prod_i|A_i|}{\prod_{e\in T}|L_e|}.
			\]
			
			\item \emph{Local extension criterion.}
			A local state $a\in A_r$ extends globally if and only if
			$1^{\widetilde a}\in O_r$ for any tree lift $\widetilde a$.
			If $1^n=1^{\widetilde a}$ with $n\in N_r$, then
			$\widetilde a n^{-1}\in W$ is a global extension.
			The pointed quotient of $O$ by $N_r$-orbits is canonically the
			homogeneous $A_r$-set $\pi_r(W)\backslash A_r$; it need not be a group.
			
			\item \emph{Projection index.}
			The orbit sizes give
			\[
			[A_r:\pi_r(W)]=\frac{|O|}{|O_r|}\le R_T,
			\qquad |W\cap N_r|=\frac{|N_r|}{|O_r|}.
			\]
		\end{enumerate}
	\end{theorem}

	\begin{proof}
		\leavevmode
		\begin{enumerate}[label=(\roman*),leftmargin=*]
			\item \emph{Source cardinality.} The action law follows by direct multiplication in the correct order:
			$(d^s)^t=d^{st}$. Fixing the identity tuple is exactly the system of all
			chord equalities. Orbit--stabilizer proves the source formula, and
			successive tree leaf attachments prove the formula for $|S|$.
			
			\item \emph{Local extension criterion.} Normality of $N_r$ implies $(d^{N_r})^s=(d^s)^{N_r}$, so its orbits
			form blocks. Their induced transitive action factors through
			$S/N_r\cong A_r$. The stabilizer of $O_r$ in $S$ is $WN_r$:
			$1^s=1^n$ is equivalent to $sn^{-1}\in W$. Its image in $A_r$ is
			$\pi_r(W)$, proving the coset identification and lift criterion.
			
			\item \emph{Projection index.} All $N_r$-orbits in $O$ have the same size. Counting them gives the
			index ratio; $O\subseteq Q$ gives the bound. Applying orbit--stabilizer
			to $N_r$ gives the remaining cardinality.\qedhere
		\end{enumerate}
	\end{proof}

	\subsubsection{Orbit algorithm, sharp index, and nonnormal root support}
	\label{app:group-orbit-algorithm}
	\label{app:joint-cycle-diagnostics}
	
	The orbit algorithm uses the tree matching group $S$, its
	chord-value action on $Q$, the root kernel $N_r$, and the
	bound $R_T$ of \cref{thm:joint-chord-witness}.
	
	\paragraph{Generators for the matching group.}
	For full local multiplication tables, compute a generating family for
	$S$ without enumerating it: lift each root element outward on $T$, and
	for each nonroot vertex lift each element of its parent-interface kernel on
	its subtree, taking identity outside. Multiplication by inverse lifts
	in outward order reduces any matching tuple to identity, proving
	generation. Omitting the root lifts generates $N_r$.
	
	Traverse the orbits on $Q$ under these generators, storing transporters
	$t_d$ with $1^{t_d}=d$ and $t_1=1$. Include the inverses of the
	generators in the traversal. At most $R_T$ states are visited. The classical
	Schreier generators
	\[
	t_d g t_{d^g}^{-1}
	\]
	generate $W$; telescoping a generator word that fixes $1$ proves this.
	The $N_r$ traversal supplies the correction in
	Theorem~\ref{thm:joint-chord-witness}. This is polynomial in explicit input size and
	$R_T$. Store each transporter as a tuple of vertex-group elements.
	Each traversal step uses one table multiplication per vertex and the
	listed incidence maps. The output consists of subgroup generators.
	The orbit and stabilizer calculations use standard group action
	theory \cite{DixonMortimer1996}. For graphical realizations of groups
	see \cite{Forney2014}.
	
	\paragraph{Sharp index and nonnormal root support.}
	With $\beta$ chords, take $A_i=\mathbb F_3^\beta$, identity tree interfaces, and
	chord $j$ maps $x\mapsto x_j$ and $x\mapsto -x_j$ into $\mathbb F_3$.
	The full source is zero, so the index equals $3^\beta=R_T$.
	For a triangle of $S_3$ groups with identity tree maps and chord maps
	identity and conjugation by $(12)$, the root support is the nonnormal
	two-element centralizer.

	\subsubsection{A complete full-source test on one group-native cycle}

	Traversing a cycle returns a relation on its starting chart. Its
	diagonal detects the states that extend all the way around the
	cycle. Restricting the same relations to prescribed agreement sets
	tests recovery as well.
	
	Let $\Gamma$ be a simple cycle on $n\ge3$ vertices, with finite
	vertex groups $A_i$, finite edge groups $B_{ij}$ and incidence
	homomorphisms $p_{i,ij}:A_i\to B_{ij}$ having equal actual images
	at each edge.  All groups have marked native presentations under the
	fixed finite-group presentation law of
	\cref{def:equational-native-context}.  On each chart fix a finite
	family of homomorphisms as descriptor readouts.  Assume that every
	shared descriptor has the same value on all matching tuples, so that
	the global homomorphisms $f_d:W\to D_d$ are well-defined.  The local
	readouts need not jointly separate the vertex groups in this cycle test.
	Write $D$ for the descriptor set, $c_0$ for the closure generated by
	all local determination rules, and
	\[
	W=\{(a_i)\in\prod_i A_i:
	p_{i,ij}(a_i)=p_{j,ij}(a_j)\text{ for each edge }ij\}.
	\]
	Because all interface conditions are group equations, $W$ is a nonempty
	subgroup of the product, and $c_0\subseteq c_W$.

	\begin{proposition}[Coverage, recovery, cut information, and native realization on a group cycle]
		\label{app:cyclic-full-source-test}
		In the preceding fixed group-native class the following tests are
		necessary and sufficient for the corresponding clauses of the
		full-matching-source preservation contract.
		\begin{enumerate}[label=\textup{(\roman*)}]
			\item \emph{Chart coverage.}
			For an oriented edge $ij$ let
			$R_{ij}=\{(a,b)\in A_i\times A_j:p_{i,ij}(a)=p_{j,ij}(b)\}$.
			Compose these subgroup correspondences once around the cycle from
			vertex $i$ back to $i$, obtaining $H_i\le A_i\times A_i$.
			By Goursat's lemma \cite[Theorem~2.1]{BauerSenZvengrowski2015},
			choose onto homomorphisms
			$\alpha_i,\beta_i:A_i\twoheadrightarrow Q_i$ with
			$H_i=\{(a,b):\alpha_i(a)=\beta_i(b)\}$.
			Then
			\[
			\pi_i(W)=\operatorname{Eq}(\alpha_i,\beta_i),\qquad
			\pi_i(W)=A_i\ \Longleftrightarrow\ \alpha_i=\beta_i.
			\]
			Thus the state-exactness clause $\mathsf S$ holds precisely when the
			last equality holds at every vertex.
			
			\item \emph{Recovery faithfulness.}
			Let $\mathcal M(c_0)$ be the proper meet-irreducible closed sets
			of $c_0$.  For $M\in\mathcal M(c_0)$ put
			\[
			K_i(M)=\bigl\{a\in A_i:
			f_{i,d}(a)=1_{D_d}\ (d\in M\cap D_i),\quad
			f_{i,d}(a)\ne1_{D_d}\ (d\in D_i\setminus M)\bigr\},
			\]
			where $D_i$ denotes the descriptors observed on chart $i$.
			Restrict $R_{ij}$ to
			$R^M_{ij}=R_{ij}\cap(K_i(M)\times K_j(M))$.
			Then $\mathsf F$, namely $c_W=c_0$, holds if and only if,
			for every $M\in\mathcal M(c_0)$, the relational product of the
			$R^M_{ij}$ around the full cycle meets the diagonal of $K_i(M)^2$.
			This test is independent of the starting vertex $i$.
			
			\item \emph{Complete cut information.}
			For any prescribed bipartition of the vertex set $L\mid R$, form
			the partial groups $U_L$ and $U_R$ satisfying all edges internal
			to their own sides.  Collect the interfaces on all edges crossing the
			cut into the same product group $B_{\partial}$, with homomorphisms
			$b_L:U_L\to B_{\partial}$ and $b_R:U_R\to B_{\partial}$.
			Put $Y_{\partial}=\im b_L\cap\im b_R$, and restrict $U_L,U_R$
			to the inverse images of $Y_{\partial}$, calling the restricted
			groups $W_L,W_R$. Then
			\[
			W\cong W_L\times_{Y_{\partial}}W_R,
			\qquad \mathbb C[W_L]\cap\mathbb C[W_R]
			=\mathbb C[Y_{\partial}]\quad\text{inside }\mathbb C^W.
			\]
			If the \emph{preassigned} separator observation is
			$h=\rho\circ b_{\partial}$ for a fixed
			$\rho:Y_{\partial}\to Z$, its exactness
			$\mathbb C[W_L]\cap\mathbb C[W_R]=\mathbb C[h]$ holds if and
			only if $\rho$ is injective on the actual reachable joint boundary
			$Y_{\partial}$.  The same condition must be tested separately for
			every named cut.
			
			\item \emph{Native realization.}
			The finite table presentation $P[W]$ gives a native realization
			of the full matching group under the already fixed law, with
			native homomorphisms to the prescribed local groups.  Thus the
			native clause $\mathsf N$ holds for this constructed source.
		\end{enumerate}
		The conjunction of \textup{(i)}--\textup{(iv)} is an exact test for
		$\mathsf S+\mathsf F+\mathsf E+\mathsf N$ for this \emph{full}
		matching source.
	\end{proposition}
	\begin{proof}\leavevmode
		\begin{enumerate}[label=\textup{(\roman*)},leftmargin=*]
			\item \emph{Chart coverage.}
			Each $R_{ij}$ is a subdirect subgroup of $A_i\times A_j$ because
			the two incidence maps have the same actual image. To compose two
			such subgroups, take their fiber product over the common vertex group
			and project to the outer coordinates. The image is a subgroup, and
			surjectivity at both ends follows by extending across the two edges in
			either order.  A root value $a$
			extends around the cycle exactly when the return value can be chosen
			equal to $a$; thus $\pi_i(W)=\{a:(a,a)\in H_i\}$.  The Goursat
			presentation gives the equalizer formula and (i).
			
			\item \emph{Recovery faithfulness.}
			Two states of $W$ agree at descriptor $d$ exactly when
			their group difference $w^{-1}w'$ is killed by $f_d$. Consequently
			the exact agreement sets of ordered state pairs are precisely
			$\{d:f_d(g)=1\}$ for elements $g\in W$.  Every such agreement set
			is $c_0$-closed.  The closure induced by the entire group is the
			intersection closure of these agreement sets. In a finite closure
			system, its proper meet-irreducibles must each occur in any finite
			intersection-generating family of closed sets. Conversely, their
			occurrence generates the entire closed-set family and establishes
			$c_W=c_0$.  An element $g\in W$ has exact agreement set $M$ if and
			only if its vertex coordinates lie in all the $K_i(M)$ and satisfy
			the edge matches.  Relational composition around the cycle imposes
			exactly this common system of constraints, with equality of the
			initial and final root coordinate. This proves (ii).
			
			\item \emph{Complete cut information.}
			The only conditions connecting a partial left assignment
			to a partial right assignment are equality of all their crossing-interface
			values. Restricting to the intersection $Y_{\partial}$ of reachable
			joint boundary images removes precisely the nonextendible partial
			assignments. Both restricted boundary maps are onto $Y_{\partial}$,
			so each fiber of $W$ over $y$ is the full nonempty rectangle
			$W_{L,y}\times W_{R,y}$. A function on $W$ which depends only on
			the left side and only on the right is therefore constant on each
			such rectangle, and is an arbitrary function of $y$. Thus the common
			algebra is $\mathbb C[Y_{\partial}]$.  It equals
			$\mathbb C[h]$ exactly when the named $h$ separates all reachable
			boundary values, which is the injectivity condition on $\rho$.
			\item \emph{Native realization.}
			Closure of finite groups under finite products and
			subgroups makes $W$ a native finite group. Its table presentation
			and maps to the original marked presentations are supplied by
			\cref{lem:finite-table-presentation}.\qedhere
		\end{enumerate}
	\end{proof}

	\subsubsection{Joint witnesses for chord constraints}
	One group difference must witness each required failure of determination
	while satisfying every chord equation. The count below tests this for
	any finite family. Two chords already exhibit a failure of separate tests.
	
	Fix a connected graph with a chosen spanning tree $T$ and additional
	edges $e_1,\ldots,e_m$.  Assume equal actual interface images
	on the tree, and let $S$ be the corresponding tree matching group.
	Suppose the retained descriptor readouts are homomorphisms
	$f_d:S\to D_d$. These include private chart homomorphisms and shared
	readouts already identified by $T$.
	For $1\leq j\leq m$ define the chord-valid subgroups
	\[
	C_j=\bigl\{s\in S:
	p_{i_j,e_j}(s_{i_j})=p_{k_j,e_j}(s_{k_j})\bigr\},
	\qquad W=\bigcap_{j=1}^m C_j.
	\]
	For $M\subseteq D$ let
	\[
	S_M=\bigl\{s\in S:
	f_d(s)=1_{D_d}\ (d\in M),\quad
	f_d(s)\ne1_{D_d}\ (d\notin M)\bigr\}.
	\]
	\begin{proposition}[Exact joint-chord recovery certificate]
		\label{app:two-chord-certificate}
		Suppose $c_0\subseteq c_W$ is the least closure generated by the
		prescribed local rules.
		\begin{enumerate}[label=\textup{(\roman*)},leftmargin=*]
			\item \emph{Joint witnesses.}
			\begin{equation}\label{app:eq-joint-chord}
				c_W=c_0\quad\Longleftrightarrow\quad
				S_M\cap W\ne\varnothing
				\quad\text{for every }M\in\mathcal M(c_0).
			\end{equation}
			Already for $m=2$ the separate conditions $S_M\cap C_j\ne\varnothing$
			can hold with different witnesses and fail jointly.
			
			\item \emph{Counting witnesses.}
			An exact finite cardinality test is obtained by putting
			$K_M=W\cap\bigcap_{d\in M}\ker f_d$:
			\begin{equation}\label{app:eq-joint-chord-count}
				\left|S_M\cap W\right|
				=\sum_{J\subseteq D\setminus M}(-1)^{|J|}
				\left|K_M\cap\bigcap_{d\in J}\ker f_d\right|.
			\end{equation}
			Positivity tests the joint recovery condition. For chart-state
			coverage, use the simultaneous orbit test in \cref{thm:joint-chord-witness}.
		\end{enumerate}
	\end{proposition}
	\begin{proof}\leavevmode
		\begin{enumerate}[label=\textup{(\roman*)},leftmargin=*]
			\item \emph{Joint witnesses.}
			A source difference with exact descriptor agreement set $M$
			lies in $S_M$. It satisfies every chord equation precisely when
			it belongs to $W$.  Apply the finite
			meet-irreducible agreement-set criterion proved in
			\cref{app:cyclic-full-source-test} to $W$, yielding \eqref{app:eq-joint-chord}. The example below exhibits failure of the separate chord tests.
			\item \emph{Counting witnesses.}
			To count such differences, start from the subgroup $K_M$ and
			exclude the union of the kernels $\ker f_d$, $d\in D\setminus M$.
			Finite inclusion--exclusion gives \eqref{app:eq-joint-chord-count}.\qedhere
		\end{enumerate}
	\end{proof}
	
	\begin{example}[A two-cycle synchronization defect]
		\label{app:two-chord-theta-example}
		Work with finite $\mathbb F_2$-modules, regarded as elementary abelian
		groups under addition.  Let the simple theta graph
		have vertices $0,1,2,3$, spanning-tree edges
		$01,12,13$, and chords $02,03$.  Its two triangles
		$0120$ and $0130$ share the edge $01$.
		At each vertex $v$ use the full native group
		\[
		A_v=\mathbb F_2^{\delta(v)}
		\]
		with one independent coordinate $x_{v,e}$ for each incident graph
		edge $e$.  Both incidence maps at $e$ are the coordinate projection
		to $\mathbb F_2$.  The edge maps are onto. Their scalar values are interface observations.
		The four protected recovery descriptors are the private readouts
		$D=\{d_0,d_1,d_2,d_3\}$, with
		\[
		\begin{aligned}
			d_0&=x_{0,02}+x_{0,03}, &
			d_2&=x_{2,02}+x_{2,12},\\
			d_1&=x_{1,12}+x_{1,13}, &
			d_3&=x_{3,13}+x_{3,03}.
		\end{aligned}
		\]
		Each local descriptor is nonconstant and is the only descriptor in
		its local recovery vocabulary, so $c_0$ is the free closure on $D$.
		The private readouts do not jointly separate local states.
		\Cref{thm:complete-port-native-descent} requires joint separation.
		
		For any subset $F\subseteq\{02,03\}$, impose the tree equations
		and precisely the chord equations in $F$, while retaining the same
		vertex groups and readouts. Call the resulting matching source $W_F$.
		Every original chart projection $W_F\to A_v$ is onto, because the
		coordinates on distinct edges can be assigned independently and a
		matching equality only identifies the two incidences of its own edge.
		There are ten local incidence coordinates, and each imposed edge
		removes one independent dimension. Consequently
		\[
		\dim W_{\varnothing}=7,\qquad
		\dim W_{\{02\}}=\dim W_{\{03\}}=6,\qquad
		\dim W_{\{02,03\}}=5.
		\]
		In either one-chord source the four private descriptor forms are
		linearly independent. Indeed, when chord $02$ is absent the two
		incidences of $02$ remain independent; their coefficients force the
		coefficients of $d_0,d_2$ in any linear relation to vanish. The
		matched edges $12$ and $13$ then force the coefficients of
		$d_1,d_3$ to vanish. The argument with $03$ absent is symmetric.
		Hence
		\[
		c_{W_{\{02\}}}=c_{W_{\{03\}}}=c_0.
		\]
		With both chords imposed, write $x_e$ for their common incidence
		values. The retained descriptors become the four successive edge
		sums around the cycle $0,2,1,3,0$, and satisfy
		\[
		d_0+d_2+d_1+d_3=0.
		\]
		Deleting any one of the four sums leaves the three incidence
		vectors of a path. In a relation among them, each end coordinate
		forces its adjacent coefficient to vanish; the middle coefficient
		then vanishes too. The four forms therefore have rank three.
		Thus each three-element subset of $D$ determines the remaining
		descriptor in $W_{\{02,03\}}$, although $c_0$ is free:
		\[
		c_{W_{\{02,03\}}}\ne c_0.
		\]
		The matching groups have finite-module presentations. In each
		graph of imposed edges the complete tuple of crossing-edge values
		gives exact common information at every cut
		(\cref{thm:arbitrary-graph-cut-information}).
	\end{example}

	\subsubsection{Normal-chord recovery and coverage}
	A failed dependency survives simultaneous normal chord constraints when
	the normal images of its witnesses have the following commutator property.

	\begin{theorem}[Recovery preservation under normal subgroup intersections]
		\label{thm:perfect-monolith-recovery}
		Let $S$ be a group with globally consistent homomorphic descriptor
		readouts $f_d:S\to D_d$ and target determination closure $c_0$.
		Suppose that, in each nontrivial actual descriptor image
		$G_d=f_d(S)$, any two nontrivial normal subgroups have nontrivial
		commutator subgroup. For finite $G_d$ this means that it has
		a unique minimal nontrivial normal subgroup $M_d$ and that
		\begin{equation}\label{eq:perfect-monolith-condition}
			[M_d,M_d]=M_d.
		\end{equation}
		For any nonempty finite family of normal subgroups
		$C_1,\ldots,C_m\trianglelefteq S$ satisfying
		$c_{C_j}=c_0$ for every $j$, one has
		\begin{equation}\label{eq:perfect-monolith-recovery}
			c_{\bigcap_{j=1}^m C_j}=c_0.
		\end{equation}
	\end{theorem}

	\begin{proof}
		In a finite nontrivial group, distinct minimal normal subgroups
		intersect trivially and commute. The derived subgroup of a minimal
		normal subgroup is characteristic in it and is therefore either
		trivial or the whole subgroup. The commutator condition forces a
		unique perfect minimal normal subgroup $M$. Conversely, every
		nontrivial normal subgroup contains $M$ when it is unique, and
		$[M,M]=M$ then implies that any two such subgroups have a
		nontrivial commutator. This proves the finite characterization.
		
		Restrictions to subgroups preserve all positive implications,
		so it is enough to retain every nonimplication
		$e\notin c_0(U)$. Write
		\[
		H_U=\bigcap_{d\in U}\ker f_d,
		\qquad N_j=H_U\cap C_j.
		\]
		Each $N_j$ is normal in $S$. Faithfulness of $C_j$ implies
		$f_e(N_j)\ne1$: the difference of two states in $C_j$ that agree on
		$U$ and differ on $e$ lies in $N_j$ and has nonidentity $e$-image.
		Put $L_j=\bigcap_{i=1}^jN_i$. We prove by induction that
		$f_e(L_j)\ne1$. This holds for $j=1$. At the next step the
		images $f_e(L_j)$ and $f_e(N_{j+1})$ are nontrivial normal
		subgroups of $G_e=f_e(S)$. Their commutator subgroup is nontrivial
		by the hypothesis on $G_e$. Normality in $S$ gives
		\[
		1\ne[f_e(L_j),f_e(N_{j+1})]
		=f_e([L_j,N_{j+1}])\le f_e(L_{j+1}).
		\]
		An element of $L_m$ with nonidentity $e$-image is the required
		group difference in the joint source.
		Every required nonimplication survives, proving the equality
		of closures.
	\end{proof}

	The criterion applies, for example, to $S_5$ with perfect monolith $A_5$.
	The next two results give the corresponding coverage bound and
	the exact descriptor-image criterion for normal constraints.
	\begin{proposition}[Commutator protection for normal native chords]
		\label{app:joint-normal}
		Let $S$ be a group, $\pi:S\twoheadrightarrow A$ and
		$C_1,\ldots,C_m\trianglelefteq S$ with $\pi(C_j)=A$.
		Assume $m\ge1$. Write $\gamma_1(A)=A$ and
		$\gamma_{r+1}(A)=[\gamma_r(A),A]$.
		\begin{enumerate}[label=\textup{(\roman*)},leftmargin=*]
			\item \emph{Lower central series bound.}
			\[
			\gamma_m(A)\ \le\ \pi\Bigl(\bigcap_{j=1}^m C_j\Bigr).
			\]
			For two normal chords, the quotient of $A$ by the joint image is
			abelian. A perfect $A$ retains every state under any such finite family.
			
			\item \emph{Two quotient maps.}
			When $C_1C_2=S$, the two oriented Goursat quotient maps induce
			a surjection from $A$ to an abelian group with kernel
			$\pi(C_1\cap C_2)$.
		\end{enumerate}
	\end{proposition}

	\begin{proof}
		\leavevmode
		\begin{enumerate}[label=\textup{(\roman*)},leftmargin=*]
			\item \emph{Lower central series bound.}
			Choose preimages in the individual normal subgroups of arbitrary
			elements of $A$. A left-nested commutator of one element from each
			$C_j$ belongs to every $C_j$: normality places each successive
			commutator in each subgroup already used. Its image is the
			corresponding left-nested commutator in $A$. Such commutators
			generate $\gamma_m(A)$. For two chords the joint image is normal
			and contains $[A,A]$. If $A$ is perfect, its lower central series is constant, so the joint image equals $A$.
			\item \emph{Two quotient maps.}
			If $C_1C_2=S$, the map $S\to(S/C_1)\times(S/C_2)$ is onto and the image of
			$N=\ker\pi$ is a normal subdirect product. Goursat's lemma
			\cite[Theorem~2.1]{BauerSenZvengrowski2015}
			expresses it as the equality relation of two oriented maps to one
			common quotient $J$. Its image in $J\times J$ is the diagonal,
			which is normal because the two quotient maps are onto.
			Conjugating $(v,v)$ by $(u,1)$ gives $(uvu^{-1},v)$ on the
			diagonal; hence $J$ is abelian. The quotient maps $q_j:S\to S/C_j$ define, for $a=\pi(s)$,
			\[
			\delta(a)=\lambda_1(q_1(s))\lambda_2(q_2(s))^{-1}.
			\]
			This is a well-defined homomorphism because $J$ is abelian and
			$q(N)$ lies in the Goursat relation. It is onto since
			$(q_1,q_2):S\twoheadrightarrow(S/C_1)\times(S/C_2)$ is onto,
			and has kernel $\pi(C_1\cap C_2)$: the relation
			$\lambda_1(q_1(s))=\lambda_2(q_2(s))$ means
			$q(s)\in q(N)$, equivalently $s\in N(C_1\cap C_2)$.
			This yields the exact oriented obstruction.\qedhere
		\end{enumerate}
	\end{proof}

	\begin{corollary}[Normal-subgroup protection and its converse]
		\label{app:joint-normal-recovery}
		\leavevmode
		\begin{enumerate}[label=\textup{(\alph*)},leftmargin=*]
			\item \emph{Preservation under intersection.}
			Let $S$ be a group with globally consistent homomorphic descriptor readouts, with target closure $c_0$.
			Suppose every nontrivial actual descriptor image is nonabelian simple.
			If normal chord subgroups $C_j$	form a nonempty finite family and are each recovery-faithful ($c_{C_j}=c_0$), then their intersection
			is recovery-faithful.
			If each original chart is perfect and all $C_j$ individually cover it, the intersection preserves every original chart image as well.
			Exact cut information still requires the preassigned separator to detect the full reachable joint crossing-interface value.
			
			\item \emph{Criterion for protection of nontrivial images.}
			Fix any nontrivial group $G$.
			The following are	equivalent:
			\begin{enumerate}[label=\textup{(\roman*)},nosep]
				\item For every surjection $f:S\twoheadrightarrow G$, the
				intersection of any two normal subgroups of $S$ on each of which
				$f$ is nonconstant still has nonconstant $f$-image.
				
				\item $[N,L]\ne1$ for every two nontrivial normal subgroups
				$N,L$ of $G$, allowing $N=L$.
			\end{enumerate}
			Condition \textup{(i)} then holds for every nonempty finite number
			of subgroups. For finite $G$, it holds precisely when $G$ has exactly
			one minimal nontrivial normal subgroup $M$ and $[M,M]=M$. In
			this finite case one may restrict \textup{(i)} to finite $S$.
		\end{enumerate}
	\end{corollary}

	\begin{proof}
		\leavevmode
		\begin{enumerate}[label=\textup{(\alph*)},leftmargin=*]
			\item \emph{Preservation under intersection.}
			Every nonabelian simple group is perfect. Its only nontrivial
			normal subgroup is the group itself, so \cref{thm:perfect-monolith-recovery} gives recovery
			faithfulness. The chart statement follows from
			\cref{app:joint-normal} applied to each perfect chart image.
			The complete crossing-edge value is observable on both sides,
			so exact common information at a named coarsening requires
			that it separate all attained complete boundary values. For the
			full graph matching source, \cref{thm:arbitrary-graph-cut-information}
			also gives rectangular fibers over those values.
			
			\item \emph{Criterion for protection of nontrivial images.}
			For \textup{(ii)} $\implies$ \textup{(i)}, let $C,D\trianglelefteq S$ have nontrivial $f$-images. Condition \textup{(ii)} gives
			$1\ne[f(C),f(D)]\le f(C\cap D)$, and induction gives the
			finite-family assertion.
			For \textup{(i)} $\implies$ \textup{(ii)}, suppose nontrivial normal
			subgroups $N,L\trianglelefteq G$ commute.
			If $N\cap L=1$, take $f=\operatorname{id}_G$ and use these two subgroups in \textup{(i)}.
			Otherwise $M=N\cap L$ is a nontrivial abelian normal subgroup. Put
			\[
			S=G\times_{G/M}G,
			\qquad f(x,y)=x.
			\]
			The subgroups $M\times\{1\}$ and $\{(m,m):m\in M\}$ intersect
			trivially and both map onto $M$. The first is normal in $S$.
			For $(x,y)\in S$, the element $y^{-1}x$ lies in $M$ and hence
			centralizes $M$. Thus $xmx^{-1}=ymy^{-1}$ for every $m\in M$,
			proving normality of the second. This contradicts \textup{(i)}.
			
			The finite-group characterization is proved in
			\cref{thm:perfect-monolith-recovery}. Both counterexamples above
			are finite when $G$ is finite.\qedhere
		\end{enumerate}
	\end{proof}

	The same commutator condition characterizes groups that are equational
	domains for polynomial equations, by the congruence-permutable
	criterion of Aichinger, Behrisch and Rossi
	\cite[Theorem~4.8]{AichingerBehrischRossi2025}. Here it characterizes
	preservation of a nonconstant readout upon intersection in every group
	mapping onto the fixed image.

	\begin{example}[Perfect descriptor images alone do not suffice]
		\label{ex:perfect-image-not-monolith}
		Let $G=A_5$, $S=G\times G$, and take a single descriptor
		$f_e=\operatorname{id}_S$ with nonconstant target closure
		$c_0=c_S$. Each normal subgroup $C_1=G\times\{1\}$ and $C_2=\{1\}\times G$
		retains $c_{C_j}=c_0$. On their intersection $C_1\cap C_2=1$ the
		observation is constant, so the nonimplication
		$e\notin c_0(\varnothing)$ fails. The image $S$ is perfect and
		has two distinct minimal nontrivial normal subgroups.
		Thus merely requiring each descriptor image to be perfect
		would make \eqref{eq:perfect-monolith-recovery} false.
	\end{example}
	
	Normality is essential even for $G=A_5$. Choose $1\ne g\in G$.
	The graphs of the identity and of conjugation by $g$ both project
	onto the first factor of $G\times G$. Their intersection projects
	onto $C_G(g)$, which is proper since $G$ has trivial center. Neither
	graph is normal: conjugation in the second factor can change the
	second entry without changing the first.

	A protected-difference contract can use the group-orbit test on its
	subgroup diagrams when their restricted images match on the tree.

	\begin{proposition}[Vanishing extension fractions in finite pointed Mal'cev atlases]
		\label{prop:malcev-group-fraction-counterexample}
		Even with one cycle, a binary chord, and surjective interface maps, there is no
		positive lower bound on the fraction of local states extending globally
		for arbitrary finite pointed Mal'cev algebras.
	\end{proposition}

	\begin{proof}
		For $M\ge1$, choose a set $A$ mapping onto $\mathbb F_2^2$ by $f$, with each
		diagonal fiber of size one and each off-diagonal fiber of size $M$.
		Choose a representative $s(v)$ in every fiber, and define
		\[
		m(x,y,z)=
		\begin{cases}
			z,&x=y,\\
			x,&y=z,\\
			s(f(x)+f(y)+f(z)),&\text{otherwise}.
		\end{cases}
		\]
		The overlap of the first two cases is consistent. The Mal'cev identities
		hold and $f$ is a homomorphism to the affine Mal'cev operation. Its two
		coordinate maps $\alpha,\gamma:A\to\mathbb F_2$ are surjective homomorphisms.
		Use identity tree interfaces on a triangle, and chord maps $\alpha,\gamma$.
		The tree source is diagonal $A$, while the full source has
		$\alpha(a)=\gamma(a)$, hence only two states out of $2M+2$.
		The survival fraction $1/(M+1)$ tends to zero. A distinguished element
		in the zero fiber makes all maps pointed if required.
	\end{proof}
	
	\subsection{Abelian layers and group corrections}
	\label{subsec:derived-towers}
	
	A matching tuple in a quotient diagram lifts through an abelian
	layer exactly when a system of affine equations is consistent.
	We identify the obstruction and describe all lifts. For general
	group extensions the corresponding obstruction is an orbit.
	
	Fix a finite oriented matching graph in a Mal'cev variety. Suppose
	the vertex and edge algebras have quotient maps
	$A_i\to\overline A_i$ and $B_e\to\overline B_e$ by Abelian
	congruences, and the incidence homomorphisms $p_{i,e}$ descend to
	these quotients. For a prescribed matching tuple $\bar x$, its vertex
	lift fibers and its common edge lift fibers are Abelian heaps under
	the difference term \cite[Corollary~5.8]{FreeseMcKenzie1987}.
	Choose origins in these fibers and write the corresponding Abelian
	groups as $K_i$ and $K_e$. By the same result, the chosen Mal'cev term
	is affine in these coordinates. It has zero constant term by
	idempotence, and $m(x,y,y)=x$ and $m(x,x,z)=z$ force its three
	additive coefficients to be $\id,-\id,\id$. Thus it defines the same
	heap operation. Each incidence map preserves the heap
	operation, so on this fiber it has the form
	$k_i\mapsto c_{i,e}+\lambda_{i,e}(k_i)$ for a group homomorphism
	$\lambda_{i,e}:K_i\to K_e$.
	
	For each oriented edge $e=(i,j)$ put
	\[
	(\delta k)_e=\lambda_{i,e}(k_i)-\lambda_{j,e}(k_j),
	\qquad d_e=c_{j,e}-c_{i,e},\qquad
	\delta:\prod_i K_i\longrightarrow\prod_e K_e.
	\]
	Lifting $\bar x$ requires exactly $\delta k=d$. Thus $[d]$ in
	$\operatorname{coker}\delta$ is the $H^1$ obstruction and, when it
	vanishes, the lift set is a torsor for $H^0=\ker\delta$.
	Moving the vertex origins by $k^0$ replaces $d$ by $d-\delta k^0$
	under the canonical identifications of the direction groups. Moving
	an edge origin translates both of its offsets by the same amount.
	Hence the obstruction and the action on successful lifts are
	intrinsic to the chosen quotient tuple.
	
	For a compatible tower with Abelian layers this calculation applies
	at each step to a chosen matching tuple on the previous layer.
	In particular, for surjective incidences of solvable Mal'cev
	algebras use the derived congruences
	$\theta^{(0)}=\nabla$ and
	$\theta^{(j+1)}=[\theta^{(j)},\theta^{(j)}]$. The commutator quotient
	formula gives the Abelian kernel of
	$A/\theta^{(j+1)}\to A/\theta^{(j)}$ and makes the derived sequence
	compatible with each surjective incidence
	\cite[Proposition~4.4(1)]{FreeseMcKenzie1987}.
	At each layer the obstruction is evaluated at the particular matching lift
	chosen at the preceding layer.

	\paragraph{The group correction orbit.}
	For compatible group extensions
	$1\to K_i\to E_i\to\bar E_i\to1$ at vertices and
	$1\to K_e\to F_e\to\bar F_e\to1$ at edges, fix a quotient matching
	tuple and choose vertex lifts $u_i\in E_i$. If $e=(i,j)$, the
	mismatch $d_e=p_{i,e}(u_i)^{-1}p_{j,e}(u_j)$ lies in $K_e$.
	Writing $\lambda_{i,e}=p_{i,e}|_{K_i}$, right multiplication of the
	lifts by $k=(k_i)_i\in\prod_iK_i$ changes the mismatch to
	$
	(d^k)_e=\lambda_{i,e}(k_i)^{-1}d_e\lambda_{j,e}(k_j).
	$
	This is a right action, since $(d^k)^\ell=d^{k\ell}$. Any other
	initial lift differs by such a tuple $k$, so the orbit $[d]$ is
	well-defined. It contains the identity exactly when a simultaneous
	correction satisfies all edges. After one matching lift is found,
	all matching lifts are its right translates by the matching subgroup
	of $\prod_i K_i$. In the Abelian-kernel case, the correction equations
	reduce to the preceding coboundary equations.
	
	\subsubsection{Nonnormal correction orbits and oriented quotient data}
	\label{app:joint-orbit}
	\label{app:joint-structural}
	
	For nonnormal chord subgroups, the relative positions of lifts
	define an action on double cosets. Its distinguished stabilizer
	records the attainable root states.
	
	Fix a native short exact sequence
	\[
	1\longrightarrow N\longrightarrow S\overset{\pi}{\longrightarrow}
	A\longrightarrow1
	\]
	and two individually root-surjective native chord subgroups $C_1,C_2\le S$.
	Put $K_j=C_j\cap N$.
	Each $C_j$ determines $K_j$ and a section $\sigma_j:A\to N_S(K_j)/K_j$ of the homomorphism induced by $\pi$:
	the section sends $a$ to the coset of any lift of $a$ in $C_j$.
	Conversely such a kernel and section reconstruct $C_j$ by taking the full inverse image of the section subgroup.
	
	\begin{theorem}[Double-coset action and exact root-state coverage]
		\label{app:joint-double-coset}
		\begin{enumerate}[label=\textup{(\roman*)},leftmargin=*]
			\item \emph{Correction action.}
			The pointed set $\mathscr D=K_1\backslash N/K_2$ with basepoint	$d_0=K_1K_2$ carries a canonical left $A$-action.
			For lifts $c_j(a)\in C_j$ of $a\in A$ it is
			\[
			a\cdot(K_1tK_2)=K_1c_1(a)t c_2(a)^{-1}K_2.
			\]
			The action is independent of the lifts and determined by the two normalizer sections.
			\item \emph{Stabilizer and coverage.}
			For $W=C_1\cap C_2$,
			\[
			\pi(W)=\operatorname{Stab}_A(d_0),\qquad
			[A:\pi(W)]=|A\cdot d_0|,\qquad
			|W|=|K_1\cap K_2|\,|\operatorname{Stab}_A(d_0)|.
			\]
			Thus complete root-state coverage is equivalent to the basepoint being fixed under the entire $A$-action.
		\end{enumerate}
	\end{theorem}
	
	\begin{proof}\leavevmode
		\begin{enumerate}[label=\textup{(\roman*)},leftmargin=*]
			\item \emph{Correction action.}
			Individual surjectivity gives $S=NC_j$. The two left coset
			spaces $S/C_j$ are transitive $N$-sets identified with $N/K_j$.
			The diagonal $N$-orbit of $(C_1,tC_2)$ is indexed by
			$K_1tK_2$: translate the first coset to $C_1$ and the remaining
			freedom on the left is $K_1$. As $N$ is normal, the left $S$-action
			on this orbit space factors through $A=S/N$. For a lift
			$c_1(a)$, the translated pair is
			$(C_1,c_1(a)t c_2(a)^{-1}C_2)$, giving the stated formula.
			Each $c_j(a)$ normalizes $K_j$ and a different lift changes it
			by a factor in $K_j$, which gives the same double coset.
			\item \emph{Stabilizer and coverage.}
			The stabilizer in $S$ of the ordered coset pair $(C_1,C_2)$ is
			$C_1\cap C_2=W$. Its $N$-orbit has stabilizer $NW$, so the
			stabilizer of the distinguished point in $A$ is
			$NW/N=\pi(W)$. Orbit--stabilizer and the kernel of $\pi|_W$
			give the stated index and cardinality formulas.\qedhere
		\end{enumerate}
	\end{proof}
	
	\paragraph{Continuation states for two automorphism-graph chords.}
	For the root-reorientation language on $\mathscr D$, a state $d$
	is tested by membership of named elements in $\operatorname{Stab}_A(d)$.
	Reorientation by $a\in A$ updates this subgroup by conjugation by $a$.
	Equal stabilizers therefore give identical answers after every
	reorientation. A membership query in their symmetric difference
	distinguishes any two different stabilizers. This proves the exact
	summary for that language.
	Two pairs of chords with the same root support can respond differently
	to a further chord. Take $G=A_5$ and $S=G\times G$ with the first projection
	onto $A=G$. Then
	taking both initial chords to be the identity graph or both to
	be the graph of a nonidentity automorphism gives full root coverage
	in either case. Adding the identity graph retains full coverage
	for the first choice and restricts the second choice to the
	automorphism's fixed subgroup.
	
	The ordered pair of group coordinates defines the labeled interface
	for the following count. For
	$\alpha\in\operatorname{Aut}(G)$ let
	$\Gamma_\alpha=\{(x,\alpha(x)):x\in G\}$. Use the intersections
	$W_{\alpha,\beta}=\Gamma_\alpha\cap\Gamma_\beta$ as the initial
	states. Future attachments intersect the current relation with
	another $\Gamma_\gamma$.
	The initial intersection is the graph of $\alpha|_K$ on
	$K=\operatorname{Fix}(\beta^{-1}\alpha)$. Conversely, for a fixed
	$\alpha$ and any $\tau\in\operatorname{Aut}(G)$, the choice
	$\beta=\alpha\tau^{-1}$ gives $K=\operatorname{Fix}(\tau)$ and
	the restriction $\alpha|_K$.
	
	For $G=A_5$, every automorphism is conjugation by a permutation in
	$S_5$ \cite[$A_5$ entry]{ConwayEtAl1985Atlas}. Thus the possible domains are precisely
	$K=C_G(\tau)$ for $\tau\in S_5$. For a fixed $K$, two conjugations
	restrict to the same map on $K$ exactly when their relative
	permutation centralizes $K$ pointwise. The table records the
	separate counts, with fixed coordinate labels.
	\begin{center}
		\small
		\begin{tabular}{@{}lrrr@{}}
			\toprule
			Cycle type of $\tau$ & Distinct $K$ & $|C_{S_5}(K)|$ &
			Graphs on these domains \\
			\midrule
			Identity 					 & 1  & 1 & 120 \\
			Transposition 				 & 10 & 2 & 600 \\
			Two disjoint transpositions  & 5  & 4 & 150 \\
			Three-cycle, with or without a disjoint transposition
			& 10 & 6 & 200 \\
			Four-cycle 					 & 15 & 8 & 225 \\
			Five-cycle 					 & 6  & 5 & 144 \\
			\midrule
			Total & & & $1439$ \\
			\bottomrule
		\end{tabular}
	\end{center}
	For a transposition the even centralizer acts as $S_3$ on the
	remaining triple and by its sign on the transposed pair. Its orbits
	recover the pair and triple, and its centralizer consists of the two
	permutations on the pair acting trivially on the triple. For two
	disjoint transpositions the even centralizer is the regular Klein
	four group on the four moved letters, so its support specifies the
	group and its centralizer has order four. A three-cycle gives the
	cyclic subgroup on its triple, also obtained from a disjoint product
	with a transposition. Its centralizer is cyclic of order three on
	the triple with an independent permutation of the remaining two
	letters. A four-cycle gives the order-two subgroup generated by its
	square, with centralizer of order eight obtained by preserving the
	two transposed pairs up to interchange. There are $15$ such squares.
	A five-cycle gives its own cyclic group, with $24/4=6$ distinct
	domains and centralizer order five. Together with $K=G$ these
	families have distinct orders, so the table counts each domain once.
	
	If all labeled membership queries $(x,y)\in W$ are protected,
	two different partial graphs differ on one such query. The graph
	itself determines every future intersection and query, so the
	exact initial state count is $1439$. For the language that tests
	only whether the first projection is onto $G$, every proper
	partial graph stays non-surjective under every further graph
	intersection. These graphs give one class. A full graph
	$\Gamma_\alpha$ is distinguished from every other full graph
	by attaching $\Gamma_\alpha$ before testing surjectivity. The
	initial count in this language is therefore $120+1=121$.
	
	For split two-block sources one can replace the correction action by
	an explicit oriented quotient.
	
	\begin{proposition}[Joint chart images and common information for two group correspondences]
		\label{app:joint-two-block}
		Let $S=A\times B$ and let
		$C_j=\{(a,b):f_j(a)=g_j(b)\}$ be two subdirect group
		correspondences, with onto maps $f_j:A\twoheadrightarrow Q_j$ and
		$g_j:B\twoheadrightarrow Q_j$. Put
		$F=(f_1,f_2)$, $G=(g_1,g_2)$,
		$P=F(A)\cap G(B)$, and
		$J=B/(\ker g_1\ker g_2)$. The two induced onto maps
		$\lambda_j:Q_j\to J$ are the oriented Goursat maps of $G(B)$.
		\begin{enumerate}[label=\textup{(\roman*)},leftmargin=*]
			\item \emph{Chart images.}
			Then the joint native source has root support
			\[
			\pi_A(W)=\operatorname{Eq}(\lambda_1 f_1,\lambda_2 f_2)
			=F^{-1}(P),
			\qquad \pi_B(W)=G^{-1}(P).
			\]
			Both complete original chart images are preserved if and only if
			$F(A)=G(B)$.
			\item \emph{Fiber product and common information.}
			With the actual restricted side groups,
			\[
			\begin{aligned}
				W&\cong F^{-1}(P)\times_P G^{-1}(P),\\
				|W|&=|\ker F|\,|\ker G|\,|P|,\\
				\mathbb C[F^{-1}(P)]\cap\mathbb C[G^{-1}(P)]
				&=\mathbb C[P]\quad\text{inside }\mathbb C^W.
			\end{aligned}
			\]
			For a named observation $\rho:P\to Z$, define the common
			readout on $W$ by $h(a,b)=\rho(F(a))=\rho(G(b))$. It has exact
			common information precisely when $\rho$ is injective on $P$. The common quotient $J$
			may be nonabelian and the root equalizer nonnormal.
		\end{enumerate}
	\end{proposition}
	
	\begin{proof}\leavevmode
		\begin{enumerate}[label=\textup{(\roman*)},leftmargin=*]
			\item \emph{Chart images.}
			The image $G(B)$ is a subdirect subgroup of $Q_1\times Q_2$.
			Goursat's lemma \cite[Theorem~2.1]{BauerSenZvengrowski2015}
			identifies it with the equality relation of the
			induced maps $\lambda_j$ to $J$. Thus $a$ extends to a hidden
			state exactly when $F(a)\in G(B)$, giving the first root formula,
			and symmetrically for the other side. The two complete image
			conditions are $F(A)\subseteq G(B)$ and the reverse inclusion.
			\item \emph{Fiber product and common information.}
			Each fiber over $p\in P$ is the rectangle
			$F^{-1}(p)\times G^{-1}(p)$ of size $|\ker F|\,|\ker G|$. This gives the product,
			cardinality, and function-algebra intersection formulas. The
			named separator detects all common functions exactly when it is
			injective on $P$.\qedhere
		\end{enumerate}
	\end{proof}
	
	Take $G=A_5$, $A=G^2$, and $B=G$. With root state $(x,y)$ and hidden state $b$, the two equations $x=b$, $y=b$ give the nonnormal diagonal subgroup of $G^2$ as joint root image.
	Its index in $G^2$ is $60$, although each chord separately covers $G^2$.

	\subsection{Exact quotient lifting on cyclic diagrams}
	On trees, the local condition of \cref{thm:canonical-reduction-canonical}
	suffices even when one representative is prescribed. A cycle also
	requires the transported representatives to agree on return. We
	compare lifting with a prescribed representative and lifting with
	a free choice of representative, first for isomorphism transport
	and then for multivalued linear relations.
	
	\begin{proposition}[Exact cyclic lifting test for deterministic native quotients]
		\label{prop:deterministic-cyclic-quotient}
		Let $G$ be a finite connected graph. At each vertex $v$ fix a
		finite algebra $\mathbf X_v$ in one equational signature, and
		for each oriented edge $v\to w$ a specified algebra
		isomorphism $\tau_{vw}:X_v\to X_w$, with inverse
		$\tau_{wv}$. Let $q_v:\mathbf X_v\twoheadrightarrow
		\mathbf Y_v$ be congruence quotients such that every
		$\tau_{vw}$ induces an isomorphism
		$\bar\tau_{vw}:Y_v\to Y_w$ and
		$q_w\tau_{vw}=\bar\tau_{vw}q_v$.
		Choose a spanning tree $T$ and root $o$. Let
		$h_e\in\operatorname{Aut}(\mathbf X_o)$ and
		$\bar h_e\in\operatorname{Aut}(\mathbf Y_o)$ be the
		fundamental-cycle holonomies of the two transported
		diagrams, and set
		\[
		F=\bigcap_{e\notin T}\operatorname{Fix}(h_e),\qquad
		\bar F=\bigcap_{e\notin T}\operatorname{Fix}(\bar h_e).
		\]
		The full source and quotient matching sets are naturally identified with
		$F$ and $\bar F$, respectively, and their comparison is $q_o|_F$.
		\begin{enumerate}[label=(\roman*),leftmargin=*]
			\item \emph{Ordinary lifting.} Every quotient matching state has a lift exactly when
			\begin{equation}\label{eq:det-cyclic-unpinned}
				q_o(F)=\bar F.
			\end{equation}
			\item \emph{Lifting with a prescribed representative.} Every quotient matching state lifts
			with any prescribed representative at any one vertex exactly when
			\begin{equation}\label{eq:det-cyclic-pinned}
				q_o^{-1}(\bar F)\subseteq F.
			\end{equation}
			\item \emph{Native realization.} The lifting criteria are finite permutation tests.
			Whenever the original matching set is nonempty it is a subalgebra, with the
			native finite-table realization of the fixed equational context.
		\end{enumerate}
	\end{proposition}
	\begin{proof}\leavevmode
		\begin{enumerate}[label=(\roman*),leftmargin=*]
			\item \emph{Ordinary lifting.} The tree-edge constraints uniquely transport a root state
			$x\in X_o$ to every vertex. A non-tree edge is satisfied
			if and only if $h_e(x)=x$. The same argument on the
			quotients identifies quotient states with $\bar F$.
			Commutation of the quotient maps with the transports
			implies $q_oh_e=\bar h_eq_o$, so the matching comparison
			is exactly $q_o|_F$. This gives
			\eqref{eq:det-cyclic-unpinned}.
			\item \emph{Lifting with a prescribed representative.} For a prescribed
			representative at a vertex $v$, tree transport determines
			one and only one corresponding root representative in
			$q_o^{-1}(y)$, where $y\in\bar F$ is the quotient root
			state. Such a lift exists precisely when that
			representative belongs to $F$. Requiring this for all
			vertices, all quotient matching states, and all their
			representatives gives \eqref{eq:det-cyclic-pinned}.
			\item \emph{Native realization.} All matching constraints are preserved by the algebra
			operations, proving the final native assertion.\qedhere
		\end{enumerate}
	\end{proof}
	
	\begin{example}[A cyclic lift can fail at a prescribed representative]
		\label{ex:det-cyclic-pinned-gap}
		Take a triangle whose vertex algebras are
		$\mathbb F_3$-modules. Use identity transport along a
		spanning tree and $x\mapsto-x$ on the remaining edge.
		Quotient each vertex to the zero module. The original
		matching source is $F=\{0\}$ and the quotient matching
		source is $\bar F=\{0\}$, so
		\eqref{eq:det-cyclic-unpinned} holds. But a nonzero
		representative of the quotient state cannot extend around
		the cycle. Thus \eqref{eq:det-cyclic-pinned} fails even
		in a finite module-native class.
	\end{example}
	
	\paragraph{Canonical cyclic reduction from semisimple holonomy.}
	The lifting requirements separated by
	\cref{ex:det-cyclic-pinned-gap} each have a coarsest quotient under
	semisimple cycle transport. The proof uses invariant complements,
	which Maschke's theorem supplies in coprime characteristic
	\cite{SerreRepresentations1977}.
	
	\begin{theorem}[Semisimple holonomy and coarsest exact cyclic quotients]
		\label{thm:semisimple-cyclic-reduction}
		In \cref{prop:deterministic-cyclic-quotient}, suppose the vertex algebras
		are finite-dimensional $\mathbb F_q$-modules and all transports are linear.
		Let $t_v:X_o\to X_v$ be the tree transports, let
		$H=\langle h_e:e\notin T\rangle$, and put $F=X_o^H$.
		Protect linear observations with kernels $A_v\le X_v$, and set
		\[
		A=\bigcap_v t_v^{-1}(A_v),\qquad
		A_H=\bigcap_{h\in H}hA.
		\]
		Assume $X_o$ is semisimple as an $\mathbb F_q[H]$-module; in particular,
		this holds if $\operatorname{char}\mathbb F_q$ does not divide $|H|$.
		Among all observation-preserving coordinate quotients that commute with
		the transports, the following statements hold.
		\begin{enumerate}[label=\textup{(\roman*)}]
			\item \emph{Ordinary lifting.} Every quotient has ordinary exact lifting. The greatest root
			kernel is $N_o^*=A_H$, with $N_v^*=t_v(A_H)$, so this family gives the
			coarsest ordinary exact quotient.
			\item \emph{Lifting with a prescribed representative.} A quotient has exact lifting with an arbitrary representative
			prescribed at any one vertex if and only if its root kernel satisfies
			$N_o\le F$. Its greatest allowed root kernel is therefore
			\[
			N_o^*=A\cap F,\qquad N_v^*=t_v(A\cap F).
			\]
			The same criterion applies if a representative may be prescribed only at
			one specified vertex.
			\item \emph{Canonical kernels and computation.} These greatest kernels are intrinsic to the transported diagram and the
			protected observations: a change of root or spanning tree only changes
			their root coordinates. They are computable from generators of $H$ by
			linear algebra. The quotients and their factor maps are native in the
			fixed finite-module presentation law.
		\end{enumerate}
	\end{theorem}
	\begin{proof}\leavevmode
		\begin{enumerate}[label=(\roman*),leftmargin=*]
			\item \emph{Ordinary lifting.} Transport-compatible kernel families correspond exactly to
			$H$-invariant root subspaces $N_o$, with $N_v=t_vN_o$.
			The protected-observation conditions are precisely $N_o\le A$.
			Their greatest invariant solution is $A_H$: it is invariant, lies in
			$A$, and contains every invariant subspace of $A$.
			
			Write $N=N_o$ and $q:X_o\to X_o/N$. Semisimplicity gives an invariant
			complement $U$ to $N$. The restriction $q|_U$ is an $H$-module
			isomorphism onto $X_o/N$, hence each $H$-fixed quotient vector has an
			$H$-fixed lift in $U$. Thus
			\[
			(X_o/N)^H=q(F),\qquad
			q^{-1}\bigl((X_o/N)^H\bigr)=F+N.
			\]
			By \cref{prop:deterministic-cyclic-quotient}, the first equality gives
			ordinary exact lifting.
			\item \emph{Lifting with a prescribed representative.} The second equality gives pinned exact lifting
			precisely when $F+N\subseteq F$, equivalently $N\le F$.
			The subspace $A\cap F$ is already invariant, since $H$ acts trivially
			on $F$, so it is the greatest protected kernel satisfying this condition.
			Transport identifies a prescribed representative at any vertex with a
			prescribed root representative.
			\item \emph{Canonical kernels and computation.} The holonomy group is the group of transports around all root
			loops. Changing the tree changes each $t_v$ by a root holonomy, and
			changing the root conjugates these data. Greatestness therefore proves
			the canonical claims and the unique quotient factorizations.
			
			For computation, $F$ is the common kernel of $h-\id$ for a generating
			set of $H$. Starting at $A$, repeatedly replace a current subspace $B$
			by $B\cap\bigcap_h hB$, using the generators and their inverses.
			At a stable step $B$ is invariant; every invariant subspace of $A$
			survives all steps. Thus the process computes $A_H$, with at most
			$\dim X_o$ strict dimension decreases. The native assertion follows
			from the fixed module quotient construction.\qedhere
		\end{enumerate}
	\end{proof}
	
	\begin{corollary}[Coarsest exact quotients for coprime gain continuations]
		\label{cor:coprime-gain-continuations}
		Fix a finite module $X$ over $\mathbb F_q$ and a finite subgroup
		$H\le\operatorname{GL}(X)$ of order prime to
		$\operatorname{char}\mathbb F_q$. Allow all finite connected diagrams
		whose edge constraints have the form $x_w=h_ex_v$, with $h_e\in H$.
		Require a uniform quotient kernel $N\le X$ invariant under every named
		transport in $H$, and protect a fixed family of linear readouts with
		joint kernel $A$.
		\begin{enumerate}[label=(\roman*),leftmargin=*]
			\item \emph{Greatest kernels.} The coarsest quotient exact on every diagram has
			kernel $A_H=\bigcap_{h\in H}hA$. Under the additional one-representative
			lifting contract, its kernel is $A\cap X^H$.
			\item \emph{Substitution.} Both conclusions remain valid after replacing any part of a diagram by
			another finite diagram using the same transport group, including
			replacements that introduce cycles. Matching images and their retained
			boundary-readout relations are exact for every expanded diagram.
		\end{enumerate}
	\end{corollary}
	\begin{proof}\leavevmode
		\begin{enumerate}[label=(\roman*),leftmargin=*]
			\item \emph{Greatest kernels.} Every diagram's holonomy is a subgroup of $H$, so its order is prime to
			the characteristic. For an invariant $N$, a direct proof of ordinary
			lifting averages any representative of a fixed quotient vector over
			that holonomy subgroup. The averaged vector is fixed and has the same
			quotient value. If $N\le X^H$, all representatives of this fixed
			quotient vector are fixed as well, proving pinned lifting.
			Conversely, a triangle with identity transports on two edges and
			transport $h\in H$ on the third belongs to the language. Pinned
			lifting of its zero quotient state forces every $n\in N$ to be fixed
			by $h$. Requiring this for every $h$ forces $N\le X^H$.
			The displayed kernels are the greatest ones allowed by the observations
			and these conditions.
			\item \emph{Substitution.} An expanded diagram still uses gains in $H$,
			so the same arguments apply after substitution. Projecting the equality
			of full matching images proves the boundary-readout assertion.\qedhere
		\end{enumerate}
	\end{proof}
	
	\begin{example}[The modular obstruction persists for deterministic transport]
		\label{ex:modular-deterministic-no-coarsest}
		Over $\mathbb F_p$, let $X$ have basis $e_1,e_2,e_3$ and let
		\[
		h(e_1)=e_1,\qquad h(e_2)=e_2,\qquad h(e_3)=e_3+e_1.
		\]
		Use a triangle with identity tree transports and chord transport $h$,
		and protect the $e_3$ coordinate at every vertex. Here $|\langle h\rangle|=p$
		and $F=\operatorname{Fix}(h)=\langle e_1,e_2\rangle$.
		The uniform kernels $N^1=\langle e_2\rangle$ and
		$N^2=\langle e_1+e_2\rangle$ both preserve the observation and have
		pinned exact lifting: for either kernel,
		\[
		(h-\id)x\in N^j\quad\Longleftrightarrow\quad x\in F,
		\]
		because $N^j\cap\langle e_1\rangle=0$.
		Every protected common coarsening must have kernel $F$.
		But $h$ acts trivially on $X/F$, while the original matching source is
		$F$ and has zero image there. This common coarsening fails even ordinary
		lifting. This example has several optimal pinned-exact quotients.
		The following theorem classifies the entire optimal family and the
		coarsest-quotient criterion for every cyclic gain group, including
		modular characteristic.
	\end{example}
	
	\paragraph{All optimal quotients for cyclic transport.}
	Fix a finite-dimensional $\mathbb F_q$-space $X$, an invertible linear
	map $h$, and linear observations $o_j:X\to O_j$ with joint kernel
	$A=\bigcap_j\ker o_j$. The permitted continuation language consists
	of all finite connected gain diagrams with every vertex space equal
	to $X$ and every edge transport in the finite cyclic group
	$\langle h\rangle$. Reverse edges carry inverse transports.
	A uniform quotient uses the same $h$-invariant subspace $N$ at every
	vertex. It is \emph{universally pinned-exact} if it preserves the observations
	and satisfies this test on every permitted diagram. Take any matching
	quotient state and prescribe any representative at one variable.
	A lift must exist with that representative.
	
	\begin{theorem}[Minimum cyclic-gain quotients with prescribed-representative lifting]
		\label{thm:cyclic-optimal-quotients}
		In this setting put
		\begin{equation}\label{eq:cyclic-optimal-subspaces}
			\begin{gathered}
				F=\ker(h-\id),\qquad P=A\cap F,\qquad
				D=P\cap\im(h-\id),\\
				r=\dim P,\quad d=\dim D.
			\end{gathered}
		\end{equation}
		The universally pinned-exact kernels are exactly
		\begin{equation}\label{eq:cyclic-optimal-admissibility}
			N\le P,\qquad N\cap D=0.
		\end{equation}
		They are also exactly the protected transport-compatible kernels
		with pinned exact lifting on a triangle whose two tree gains are
		identities and whose chord gain is $h$.
		The following consequences hold in every characteristic.
		\begin{enumerate}[label=\textup{(\roman*)}]
			\item \emph{Minimum size.} The maximum kernel dimension is $r-d$, and the minimum number
			of states in each vertex quotient is
			\begin{equation}\label{eq:cyclic-optimal-size}
				\min_N |X/N|=q^{\dim X-r+d}.
			\end{equation}
			Every inclusion-maximal admissible kernel attains this minimum.
			The optimal kernels are precisely the complements of $D$ in $P$.
			\item \emph{Classification and enumeration.} There are exactly
			\begin{equation}\label{eq:cyclic-optimal-count}
				q^{d(r-d)}
			\end{equation}
			optimal kernels, hence this many minimum quotient maps up to
			isomorphism over $X$. Their family is canonically the affine space
			of linear sections of $P\twoheadrightarrow P/D$, a torsor for
			$\operatorname{Hom}_{\mathbb F_q}(P/D,D)$.
			\item \emph{Greatest kernel.} A greatest admissible kernel exists if and only if
			\begin{equation}\label{eq:cyclic-coarsest-criterion}
				d=0\quad\text{or}\quad d=r.
			\end{equation}
			It is $P$ when $d=0$ and $0$ when $d=r$. In either case its quotient
			is coarsest over $X$: every admissible quotient has a unique
			surjective factor map to it. If $0<d<r$, the optimal quotients form
			the nonsingleton family in \textup{(ii)}, and no admissible quotient
			is coarser than all of them.
			\item \emph{Substitution and computation.} Every admissible kernel remains pinned-exact after arbitrary
			finite substitutions using gains in $\langle h\rangle$. Matching
			images and their retained boundary-readout relations are exact on
			each expanded diagram. The subspaces in
			\eqref{eq:cyclic-optimal-subspaces}, one optimal kernel, and a
			parameterization of all optimal kernels are computable by linear
			algebra. All these quotients and their factor homomorphisms have
			native realizations in the fixed finite-module presentation law.
		\end{enumerate}
	\end{theorem}
	\begin{proof}
		For the specified triangle the root matching space is $F$.
		For an $h$-invariant kernel $N$, a quotient root value $x+N$ is
		fixed by the induced holonomy precisely when $(h-\id)x\in N$.
		By \cref{prop:deterministic-cyclic-quotient}, pinned exactness is
		therefore equivalent to
		\begin{equation}\label{eq:cyclic-optimal-root-test}
			(h-\id)^{-1}(N)=F.
		\end{equation}
		Since $N$ is $h$-invariant, its vectors satisfy the left-hand test;
		thus \eqref{eq:cyclic-optimal-root-test} forces $N\le F$.
		It also forces $N\cap\im(h-\id)=0$: if
		$z=(h-\id)x\in N$, then $x\in F$ and $z=0$.
		Conversely, that zero intersection makes
		$(h-\id)x\in N$ imply $(h-\id)x=0$.
		Observation preservation is $N\le A$. Once $N\le P$, invariance
		is automatic and $N\cap\im(h-\id)=N\cap D$.
		This proves \eqref{eq:cyclic-optimal-admissibility} for the triangle.
		
		Now consider any permitted diagram and any quotient matching state.
		Root the diagram at the variable where a representative is prescribed
		and choose a spanning tree. Tree transport gives concrete values at
		all vertices and the required quotient values. Every root holonomy
		is $h^t$ for an integer $t$. Since $h$ has finite order, choose
		$t\ge0$. The identity
		\[
		h^t-\id=(h-\id)(\id+h+\cdots+h^{t-1})
		\]
		(with zero right side for $t=0$) gives
		$\im(h^t-\id)\subseteq\im(h-\id)$.
		The quotient matching equations say $(h^t-\id)x\in N$ for the
		prescribed root representative $x$. Their intersection with
		$\im(h-\id)$ is zero by
		\eqref{eq:cyclic-optimal-admissibility}, so every concrete holonomy
		equation holds. This gives the desired lift. Necessity for the whole
		language follows because it contains the specified triangle.
		\begin{enumerate}[label=(\roman*),leftmargin=*]
			\item \emph{Minimum size.} Projection $P\to P/D$ is injective on an admissible $N$, so
			$\dim N\le r-d$. Extend a basis of $D$ together with a basis of $N$
			to a basis of $P$. The added vectors enlarge $N$ to a complement of
			$D$. Therefore the inclusion-maximal admissible kernels are exactly
			the complements, and they all have dimension $r-d$. Hence quotient size is minimized at dimension $r-d$.
			\item \emph{Classification and enumeration.} Each complement is the image of a unique linear section
			of $P\to P/D$, obtained by inverting the restricted projection.
			The difference of two sections is a map $P/D\to D$, and adding any
			such map to a section yields another section. Thus the section
			family is the stated torsor and has $q^{d(r-d)}$ members. Two
			quotient maps over $X$ are isomorphic exactly when their kernels
			coincide, proving \textup{(ii)}.
			\item \emph{Greatest kernel.} If $d=0$, every subspace of $P$ is admissible and $P$ is greatest.
			If $d=r$, the only admissible kernel is $0$. If $0<d<r$, the count
			in \textup{(ii)} gives at least two distinct complements of the same
			maximum dimension. An admissible kernel containing both would
			have larger dimension, a contradiction. Hence there is no greatest
			kernel. When a greatest kernel $N^*$ exists, the unique factor map
			is $x+N\mapsto x+N^*$; it retains the observations and commutes
			with the transports. It is itself pinned-exact on each permitted
			diagram: lift a coarse witness to $X$ with a representative of the
			prescribed intermediate class, and then project modulo $N$.
			This proves \textup{(iii)}.
			\item \emph{Substitution and computation.} An expanded diagram after substitution is again in the same language;
			the lifting proof and then projection onto its boundary prove
			\textup{(iv)}'s exactness claims.
			Kernels, images, intersections, and a basis complement are found by
			Gaussian elimination. After choosing one section, a basis of
			$\operatorname{Hom}(P/D,D)$ parameterizes the entire optimal family.
			Each $N$ is a subspace and hence a native module congruence; its
			finite quotient has the presentation and source maps supplied by
			\cref{lem:finite-table-presentation,subsec:canonical-reduction-native-verification}.\qedhere
		\end{enumerate}
	\end{proof}
	
	In \cref{ex:modular-deterministic-no-coarsest},
	$P=\langle e_1,e_2\rangle$ and $D=\langle e_1\rangle$.
	The complete list of optimal kernels is therefore
	\[
	N_a=\langle e_2+a e_1\rangle\qquad(a\in\mathbb F_p).
	\]
	This gives exactly $p$ minimum vertex quotients with $p^2$ states,
	each exact for the entire cyclic gain language. Here $r=2$ and $d=1$.
	In the semisimple cyclic case, $D=0$ gives the unique kernel $A\cap F$
	of \cref{thm:semisimple-cyclic-reduction}.
	
	\phantomsection\label{par:all-gain-arrangement}
	\paragraph{Optimal quotients for an arbitrary finite gain group.}
	For several generators the single forbidden subspace in
	\cref{thm:cyclic-optimal-quotients} becomes a family of subspaces.
	The forbidden lines determine both the minimum quotient size and
	whether a greatest kernel exists.
	
	\begin{proposition}[Exact quotient kernels and minimum size for finite gain groups]
		\label{prop:finite-gain-kernels}
		Allow all finite connected gain diagrams with gains in a finite group
		$H\le\operatorname{GL}(X)$ over $\mathbb F_q$, in any characteristic.
		Keep the same observations with joint kernel $A$ and the same
		one-representative lifting requirement. Put
		\[
		P=A\cap X^H,\qquad D_g=P\cap\operatorname{im}(g-\id)
		\quad(g\in H).
		\]
		\begin{enumerate}[label=(\roman*),leftmargin=*]
			\item \emph{Admissible kernels.} A uniform kernel is admissible on every diagram if and only if
			$N\le P$ and $N\cap D_g=0$ for every $g\in H$.
			In particular, it suffices to test one triangle for each group element.
			\item \emph{Counting and minimum size.} Let $r=\dim P$ and let $\mathcal L$ be the set of lines in $P$ contained
			in at least one $D_g$. Maximizing the admissible dimension is now the classical
			critical problem for this set of projective points \cite{CrapoRota1970Critical}.
			The number of admissible $k$-dimensional kernels is
			\[
			b_k=\sum_{J\subseteq\mathcal L}(-1)^{|J|}
			\genfrac{[}{]}{0pt}{}{r-s_J}{k-s_J}_{q},\qquad s_J=\dim\sum_{L\in J}L .
			\]
			Here $\genfrac{[}{]}{0pt}{}{u}{v}_{q}$ counts $v$-planes in a space of dimension $u$,
			with value zero for $v<0$ or $v>u$.
			With $k_*=\max\{k:b_k>0\}$, the minimum quotient size is
			$q^{\dim X-k_*}$ and exactly $b_{k_*}$ kernels attain it.
			\item \emph{Greatest kernel.} Let $P_0$
			be the span of all lines of $P$ outside $\mathcal L$, with empty span
			zero. A greatest kernel exists precisely when $P_0\cap D_g=0$ for
			every $g$, in which case it equals $P_0$.
		\end{enumerate}
	\end{proposition}
	\begin{proof}\leavevmode
		\begin{enumerate}[label=(\roman*),leftmargin=*]
			\item \emph{Admissible kernels.} On the triangle with chord gain $g$ and identity tree gains,
			the zero quotient state must lift with every prescribed $n\in N$.
			This forces $gn=n$. Also, if $(g-\id)x\in N$, the fixed quotient
			state $x+N$ must lift with representative $x$, so $(g-\id)x=0$.
			These tests give the two necessary conditions, along with $N\le A$
			for the observations. Conversely, $H$ fixes such an $N$ pointwise.
			Root any diagram at the vertex with prescribed representative $x$.
			Its quotient state satisfies $(g-\id)x\in N$ for every fundamental
			holonomy $g$. The intersection condition makes all these differences
			zero. Transporting $x$ along a spanning tree thus gives the required
			lift. The same proof applies to every substitution within this language.
			\item \emph{Counting and minimum size.} A subspace of $P$ is inadmissible precisely when it contains a
			forbidden line. The $k$-planes containing the lines in $J$ correspond
			to the $(k-s_J)$-planes in the quotient by their span. Inclusion--exclusion
			gives $b_k$, and maximizing $k$ minimizes the quotient size.
			\item \emph{Greatest kernel.} Every admissible line lies in any greatest kernel, and every admissible
			kernel is spanned by admissible lines. Their span $P_0$ is therefore
			greatest exactly when it avoids all the forbidden lines.\qedhere
		\end{enumerate}
	\end{proof}
	
	\paragraph{Quotient lifting for multivalued linear relations.}
	For bijective transport, \cref{prop:deterministic-cyclic-quotient}
	gives the fixed-point test. In a multivalued linear diagram the
	quotient relation at each edge is the image of its given relation
	under the coordinate quotients.
	
	\begin{theorem}[Exact lifting for multivalued linear quotients]
		\label{thm:multivalued-linear-quotient}
		Let a finite relation diagram have variable spaces $X_i$ over
		$\mathbb F_q$. Write $X=\bigoplus_i X_i$. For each relation occurrence
		$a$ let $L_a:X\to E_a$ give its homogeneous equations, with $L_a$
		depending only on the incident variables. Each local relation is
		taken on its set of incident variables, so any identifications of
		repeated arguments are imposed before its quotient image is formed. Put
		$L=(L_a)_a:X\to\bigoplus_a E_a$ and $W=\ker L$.
		For quotient kernels $N_i\le X_i$ set $N=\bigoplus_i N_i$ and let
		$q:X\to X/N$ be the coordinate quotient. Let $\overline W$ be the
		matching space of the individual quotient relation images. Define
		\[
		B_N=\bigoplus_a L_a(N),\qquad
		K=\operatorname{im} L\cap B_N,\qquad
		N^{(r)}=\{n\in N:n_r=0\}.
		\]
		Then the following hold.
		\begin{enumerate}[label=\textup{(\roman*)}]
			\item \emph{Ordinary lifting.} There is a canonical isomorphism
			$\overline W/q(W)\cong K/L(N)$. Thus every quotient matching
			state has a lift exactly when $K=L(N)$.
			\item \emph{One prescribed representative.} Every quotient matching state lifts with any prescribed
			concrete representative at variable $r$ exactly when
			$K=L(N^{(r)})$. Equivalently, ordinary lifting holds and
			$\pi_r(W\cap N)=N_r$. Requiring this equality for every $r$
			gives exact lifting with a representative prescribed at any one
			variable.
			\item \emph{Several prescribed representatives.} To prescribe representatives simultaneously
			at a set $J$ of variables, replace $N^{(r)}$ by
			$N^{(J)}=\{n\in N:n_j=0\ (j\in J)\}$. Exact lifting for all
			such requests is equivalent to $K=L(N^{(J)})$, or equivalently
			to ordinary lifting together with
			$\pi_J(W\cap N)=\prod_{j\in J}N_j$.
			\item \emph{Computation and realization.} These conditions and the dimension of $\overline W/q(W)$ are
			computable by Gaussian elimination from the local equation matrices
			and bases of the $N_i$. The matching and quotient matching spaces
			are native finite modules under the fixed finite-module
			presentation law of \cref{def:equational-native-context}.
		\end{enumerate}
	\end{theorem}
	\begin{proof}\leavevmode
		\begin{enumerate}[label=(\roman*),leftmargin=*]
			\item \emph{Ordinary lifting.} A tuple $q(x)$ satisfies the quotient image of relation $a$ exactly
			when some $n\in N$ has $L_a(x+n)=0$. The correction may be chosen
			separately at each occurrence, so $q^{-1}(\overline W)=L^{-1}(B_N)$.
			The map $q(x)\mapsto Lx+L(N)$ is therefore well defined and onto
			$K/L(N)$. Its kernel is $q(W)$, since $Lx\in L(N)$ exactly when
			one correction makes all equations hold at once. This proves
			\textup{(i)}.
			\item \emph{One prescribed representative.} Choose $x$ with its $r$-coordinate equal to the
			prescribed representative. A simultaneous correction that keeps
			that coordinate fixed exists exactly when $Lx\in L(N^{(r)})$.
			Every element of $K$ arises from such a lifting request by taking
			any preimage $x$ under $L$ and prescribing $x_r$. This proves the
			first equivalence. Under ordinary lifting, $L(N^{(r)})=L(N)$
			exactly when each $n\in N$ differs from some $n'\in N^{(r)}$
			by an element of $W\cap N$. This is the stated projection
			condition.
			\item \emph{Several prescribed representatives.} With prescribed representatives at all variables in
			$J$, first choose the corresponding concrete coordinates of $x$.
			The permitted corrections are exactly $N^{(J)}$. As before,
			every $Lx\in K$ occurs in a request. Under ordinary lifting,
			$L(N^{(J)})=L(N)$ holds precisely when
			$N=N^{(J)}+(W\cap N)$, which is equivalent to the stated
			surjectivity under $\pi_J$.
			\item \emph{Computation and realization.} All tests use linear images,
			sums and intersections.
			A linear functional certifies a failed lift if it vanishes on
			the allowed correction image and takes a nonzero value on $Lx$. Finally, the equation kernels and quotient images are
			linear subspaces, and \cref{lem:finite-table-presentation} supplies
			their native realizations in the declared law.\qedhere
		\end{enumerate}
	\end{proof}
	
	By \textup{(i)}, the dimension of $K/L(N)$ is independent of the
	equation matrices. In \cref{thm:graph-relative} the prescribed input is
	a local value. Here the quotient edge relations specify the states
	requiring simultaneous correction.
	
	\begin{theorem}[Joins of exact linear quotient kernels]
		\label{thm:linear-quotient-joins}
		In the setting of \cref{thm:multivalued-linear-quotient}, put
		$I=\operatorname{im}L$ in $E=\bigoplus_aE_a$. Let $P$ be a set of
		variable indices. A coordinate kernel $N=\bigoplus_iN_i$ is
		\emph{$P$-exact} if it has ordinary lifting and satisfies the
		prescribed-representative test at each index in $P$. Each test
		prescribes a representative at one variable at a time.
		\begin{enumerate}[label=(\roman*),leftmargin=*]
			\item \emph{Lifting for a sum of kernels.} Let $N^1,\ldots,N^s$ be $P$-exact and take their coordinatewise sum
			$N=\sum_jN^j$. Then $N$ is $P$-exact if and only if
			\begin{equation}\label{eq:linear-join-distributivity}
				I\cap\Bigl(\sum_j B_{N^j}\Bigr)
				=\sum_j\bigl(I\cap B_{N^j}\bigr).
			\end{equation}
			The quotient of the left-hand side of \eqref{eq:linear-join-distributivity} by its
			right-hand side is canonically the ordinary lifting defect of the joined kernel. If
			$P\ne\varnothing$, it is also canonically the pinned lifting defect at each $r\in P$.
			It and its dimension can be computed by Gaussian elimination from the given matrices
			and kernel bases.
			\item \emph{Greatest kernel under a splitting condition.} Fix a vector-space complement $C$ of $I$ in $E$ and subspaces $A_i\le X_i$
			representing protected observations. Consider the $P$-exact coordinate kernels $N$
			satisfying $N_i\le A_i$ for all $i$ and
			\begin{equation}\label{eq:linear-join-split}
				B_N=(B_N\cap I)\oplus(B_N\cap C).
			\end{equation}
			This class contains the zero kernel and is closed under coordinatewise sums.
			Its greatest member gives the coarsest quotient in the class.
		\end{enumerate}
	\end{theorem}
	\begin{proof}\leavevmode
		\begin{enumerate}[label=(\roman*),leftmargin=*]
			\item \emph{Lifting for a sum of kernels.} Linearity gives $B_{\sum_jN^j}=\sum_jB_{N^j}$ and $L(\sum_jN^j)=\sum_jL(N^j)$. For
			any coordinate $r$, the coordinatewise nature of the kernels also gives
			\[
			\left(\sum_jN^j\right)^{(r)}=\sum_j(N^j)^{(r)}.
			\]
			By \cref{thm:multivalued-linear-quotient}, $P$-exactness of each $N^j$ implies $I\cap
			B_{N^j}=L(N^j)$ and, for every $r\in P$, also $I\cap B_{N^j}=L((N^j)^{(r)})$. Thus
			the denominator in the asserted quotient is both $L(N)$ and, for each $r\in P$,
			$L(N^{(r)})$. The lifting criterion in that theorem now proves the equivalence and
			identifies the ordinary defect. For pinned lifting, the same proof of that theorem
			identifies the defect with $(I\cap B_N)/L(N^{(r)})$. The equalities also give the
			elimination test.
			\item \emph{Greatest kernel under a splitting condition.} If each $B_{N^j}$ satisfies \eqref{eq:linear-join-split}, then their sum is the
			direct sum of its intersections with $I$ and $C$, and its intersection with $I$ is
			$\sum_j(I\cap B_{N^j})$. Hence \eqref{eq:linear-join-distributivity} holds and the
			sum is $P$-exact. It also satisfies \eqref{eq:linear-join-split}, and the
			protected-observation containments are preserved by sums. The zero kernel satisfies
			all these conditions. There are only finitely many coordinate kernels over the finite
			field, so their sum belongs to the class and contains every member. This is its
			greatest kernel, giving the coarsest quotient in the class.\qedhere
		\end{enumerate}
	\end{proof}
	
	\begin{example}[Pinned-exact linear quotients may have no coarsest member]
		\label{ex:multivalued-linear-no-coarsest}
		Take two variable spaces $X_i=\mathbb F_2^3$ with coordinates
		$(x_i,z_i,t_i)$ and protect the observations $x_i$. Use four
		relation occurrences joining both variables, each given by one of
		\[
		z_1+z_2=0,\qquad z_1+z_2+x_1=0,\qquad
		t_1+t_2=0,\qquad t_1+t_2+x_1=0.
		\]
		Each individual relation has four compatible values at one variable
		for every value at the other. Together they force
		$x_1=0$, $z_1=z_2$, and $t_1=t_2$. The quotient forgetting $z_i$
		at both variables is exact with either representative prescribed:
		the last two equations still hold in the quotient, and the dropped
		coordinates can both equal the prescribed $z_i$. The quotient
		forgetting $t_i$ is exact by the same argument using the first two
		equations. Any observation-preserving common coarsening has to
		forget both $z_i$ and $t_i$ at each variable. Its four individual
		relation images impose no equations on $(x_1,x_2)$, while the image
		of the original matching space has $x_1=0$. Its obstruction space
		has dimension one and there is no coarsest quotient with the
		one-representative lifting contract for this diagram.
	\end{example}

	\subsection{Target-specific repair of globalization obstructions}
	\label{subsec:targeted-globalization-repair}

	A missing exact-agreement witness consists of two global states.
	Each editable factor must admit both of their local tuples.
	\Cref{thm:multi-native-repair} minimizes the number of insertions
	over all pairs compatible with the protected factors.

	\begin{definition}[Factor and tuple costs for exact witness repair]
		\label{def:multi-native-repair}
		Fix a finite atlas as in \cref{def:multi-native-atlas},	a proper agreement pattern $M\subsetneq D$, and protected factors $\mathcal F_0$.
		At each remaining factor any absent tuple on its given scope may be inserted.
		The realized-state sets, readouts, and all factors enforcing consistency of shared readouts are fixed.
		For $w=((a_v,b_v))_{v\in V}\in\Lambda_0(M)$ let
		$
		A_\alpha(w)=(a_v)_{v\in S_\alpha},
		\quad
		B_\alpha(w)=(b_v)_{v\in S_\alpha}.
		$
		Define the \textbf{violated-factor cost} and \textbf{missing-tuple cost}
		by \textup{(i)} and \textup{(ii)}, respectively.
		Both minima are	$+\infty$ if $\Lambda_0(M)$ is empty.
		\begin{enumerate}[label=(\roman*),leftmargin=*]
			\item \textbf{Violated-factor cost.}
			Set
			\[
			q(w)=\sum_{\alpha\notin\mathcal F_0}
			\mathbf 1\!\left[\{A_\alpha(w),B_\alpha(w)\}\nsubseteq F_\alpha\right],
			\qquad \sigma_0(M)=\min_{w\in\Lambda_0(M)}q(w).
			\]
			
			\item \textbf{Missing-tuple cost.}
			Set
			\[
			d(w)=\sum_{\alpha\notin\mathcal F_0}
			\left|\{A_\alpha(w),B_\alpha(w)\}\setminus F_\alpha\right|,
			\qquad \Delta_0(M)=\min_{w\in\Lambda_0(M)}d(w).
			\]
		\end{enumerate}
	\end{definition}

	\begin{theorem}[Minimum witness-repair cost and cycle-rank bounds]
		\label{thm:multi-native-repair}
		\leavevmode
		\begin{enumerate}[label=(\roman*),leftmargin=*]
			\item \emph{Minimum repair cost.}
			The minimum total number of native formation tuples that must be
			inserted into the unprotected factors to create at least one
			global exact-$M$ recovery witness is exactly
			\begin{equation}\label{eq:multi-native-repair-cost}
				\quad\Delta_0(M)=\min_{w\in\Lambda_0(M)}
				\sum_{\alpha\notin\mathcal F_0}
				|\{A_\alpha(w),B_\alpha(w)\}\setminus F_\alpha|.\quad
			\end{equation}
			If $\Lambda_0(M)$ is empty, no repair produces the witness.
			
			\item \emph{Bounds and sharpness.}
			If $\Lambda_0(M)$ is nonempty, then
			\begin{equation}\label{eq:multi-native-repair-sandwich}
				\quad\sigma_0(M)\le\Delta_0(M)\le2\sigma_0(M)
				\le2|\mathcal F\setminus\mathcal F_0|.\quad
			\end{equation}
			For a connected graph with protected spanning tree $T$, the number
			of chords is $\beta=|E|-|V|+1$. Write $\Delta_T(M)$ for the minimum
			when the tree edges are the protected factors.
			Whenever the tree has a paired witness, $\Delta_T(M)\le2\beta$.
			In particular, one chord requires at most two insertions under this hypothesis.
			The factors $1$ and $2$ are both attainable with a single unprotected interface. 
		\end{enumerate}
	\end{theorem}

	\begin{proof}
		\leavevmode
		\begin{enumerate}[label=(\roman*),leftmargin=*]
			\item \emph{Minimum repair cost.}
			Any repaired global pair restricts to some $w\in\Lambda_0(M)$.
			For every unprotected factor, both tuples $A_\alpha(w)$ and
			$B_\alpha(w)$ must be admitted, so at least the number counted by
			$d(w)$ of previously absent tuples must have been inserted.
			If $\Lambda_0(M)$ is empty, the first observation rules out every
			repair. Otherwise choose a minimizer $w$ and insert exactly those missing
			tuples into their indicated factors.  Both componentwise native-state
			assignments then satisfy every factor, and their exact descriptor
			agreement is $M$.  This proves \eqref{eq:multi-native-repair-cost}.
			
			\item \emph{Bounds and sharpness.}
			Every violated factor requires either one or two distinct missing
			tuples, giving $q(w)\le d(w)\le2q(w)$ for every $w$ and hence the
			claimed inequalities after minimization. The cyclomatic estimate uses
			$|E\setminus T|=\beta$. For sharpness, take two positions each with
			state set $\{*,0,1\}$; let a shared descriptor $A$ equal $1$ at $*$
			and $0$ at $0,1$, and let the private descriptor $d_v$ equal $1$
			only at state $1$. With $D=\{A,d_1,d_2\}$ and $M=\{A\}$, take
			the protected binary relation $F_0=\{(*,*),(0,0),(1,1)\}$. The only
			protected exact-$M$ witness has components $(0,0),(1,1)$, in either
			order. For the unprotected factor use either
			$F_1=\{(*,*),(0,0)\}$ or $F_1=\{(*,*)\}$. Its allowed extension
			contains all three protected tuples. The minimum number of inserted
			$F_1$-tuples is respectively one or two, while $(*,*)$ is initially
			a global state in both cases.\qedhere
		\end{enumerate}
	\end{proof}

	\section{Simultaneous repair and compatibility}
	\label{sec:simultaneous-faithful-repair}

	A repair for one missing witness can introduce a pair that violates
	a required dependency. Given the admissible tuple insertions,
	\cref{thm:joint-faithful-repair,thm:selected-safe-repair} characterize
	the selections that supply all required witnesses while preserving the
	positive dependencies.

	\subsection{Certificates for simultaneous faithful repair}
	\label{subsec:structural-repair}

	Each recovery witness, local state, or separator edge has an edit support.
	A repair must contain the supports of its chosen certificates.

	\begin{definition}[Allowed formation envelope, edit supports, and costs]
		\label{def:safe-formation-envelope}
		In a finite native formation atlas as in
		Definition~\ref{def:multi-native-atlas}, specify allowed relations
		$F_\alpha\subseteq\widehat F_\alpha$ on the same scopes, with
		protected factors required to satisfy $\widehat F_\alpha=F_\alpha$.
		\begin{enumerate}[label=\textup{(\roman*)},leftmargin=*]
			\item \textbf{Allowed edits and formation envelope.}
			Let $\mathcal U$ be the disjoint set of allowed missing native
			tuples $(\alpha,t)$, where	$t\in\widehat F_\alpha\setminus F_\alpha$.
			For $E\subseteq\mathcal U$ write $F_\alpha[E]$ for the original
			relation together with its tuples selected by $E$.
			Let
			\[
			\Omega=\left\{x\in\prod_{v\in V}Z_v:
			x|_{S_\alpha}\in\widehat F_\alpha\ \text{for every }\alpha\right\}
			\]
			be the allowed full formation envelope.
			Protected readout consistency is assumed on $\Omega$, so each global descriptor $d\in D$	has a single well-defined value $d(x)$.
			
			\item \textbf{Activation supports and costs.}
			For $x\in\Omega$ put
			\[
			\sigma(x)=\{(\alpha,x|_{S_\alpha})\in\mathcal U:\alpha\in\mathcal F\},
			\qquad \sigma(C)=\bigcup_{x\in C}\sigma(x).
			\]
			The repaired realized-source space has the exact activation law
			\begin{equation}\label{eq:activation}
				\quad X(E)=\{x\in\Omega:\sigma(x)\subseteq E\}.\quad
			\end{equation}
			Assign each allowed edit $t\in\mathcal U$ a nonnegative cost $w_t$
			and use $w(E)=\sum_{t\in E}w_t$.
		\end{enumerate}
	\end{definition}

	When the contract permits a constructed commutant output, every nonempty
	$X(E)$ has the full source $M_E=\mathbb C^{X(E)}$ on
	$\ell^2(X(E))$ and the native judgment
	\[
	\Xi_{\ell^2(X(E))}^{\mathrm{com}}\vdash
	(\mathrm{com},M_E)\rightsquigarrow\mathfrak P_{M_E'}.
	\]
	The recovered source is $M_E''=M_E$. Its characters are the states
	of $X(E)$, including states with identical descriptor values.
	The descriptor and separator maps act by multiplication. Other native
	laws specify the realization of the repaired source separately.
	
	For each commuting separator test $i$, fix a finite nonempty
	separator range $H_i$ and prescribed left and right endpoint
	alphabets $L_i,R_i$. Fix maps $p_{L,i}:L_i\to H_i$ and
	$p_{R,i}:R_i\to H_i$ with every fiber nonempty. The endpoint and
	separator readouts on the whole allowed envelope are required to satisfy
	\[
	\ell_i:\Omega\to L_i,\qquad r_i:\Omega\to R_i,\qquad
	h_i=p_{L,i}\ell_i=p_{R,i}r_i.
	\]
	For $b\in H_i$, the fixed graph vertices are the two disjoint
	copies of $L_{i,b}=p_{L,i}^{-1}(b)$ and $R_{i,b}=p_{R,i}^{-1}(b)$.
	Its edges are the pairs $(\ell_i(x),r_i(x))$ for $x\in X(E)$
	with $h_i(x)=b$.
	
	By \cref{thm:proper-cut-incidence}, the graphs are connected exactly
	when both endpoint ranges are attained and
	$\mathbb C[\ell_i]\cap\mathbb C[r_i]=\mathbb C[h_i]$ on $X(E)$.
	Indeed, connectivity on the prescribed nonempty vertex parts forces
	coverage. With coverage, the components within each fiber determine
	precisely the functions common to the two endpoints.

	\begin{definition}[Repair certificates and conflicting pairs]
		\label{def:joint-certificates}
		Work in the finite formation envelope $\Omega$ of
		\cref{def:safe-formation-envelope}.
		Fix a target closure $c_0$ on $D$ and a nonempty family of prescribed local realized states, each recorded with its position $v$.
		Let $\mathscr M(c_0)$ be the family of proper meet-irreducible $c_0$-closed sets.
		For $x,y\in\Omega$, define their exact descriptor-agreement set by
		$
		A(x,y)=\{d\in D:d(x)=d(y)\}.
		$		
		\begin{enumerate}[label=\textup{(\roman*)},leftmargin=*]
			\item \textbf{Repair certificates.}
			Let $\Gamma$ index the following three types of requirements,
			with distinct tags for the different types.
			For each requirement $\gamma$, let $\mathscr C_\gamma$
			be the family of all subsets of $\Omega$ satisfying the
			corresponding certificate condition:
			\begin{enumerate}[label=\textup{(\alph*)},leftmargin=*]
				\item \emph{Recovery.}
				For each $M\in\mathscr M(c_0)$, a certificate is a
				two-element set $\{x,y\}\subseteq\Omega$ such that
				$
				A(x,y)=M.
				$
				
				\item \emph{Local coverage.}
				For each prescribed local state $z\in Z_v$,
				a certificate is a singleton $\{x\}\subseteq\Omega$
				such that $x_v=z$.
				
				\item \emph{Separator cuts.}
				For each declared separator $(i,b)$ and each
				nontrivial cut of its fixed bipartite vertex set,
				a certificate is a singleton $\{x\}\subseteq\Omega$
				with $h_i(x)=b$ whose edge $(\ell_i(x),r_i(x))$ crosses that cut.
			\end{enumerate}
			The members of $\mathscr C_\gamma$ are the
			\textbf{repair certificates} for requirement $\gamma$.
			
			\item \textbf{Conflicting pairs.}
			A two-element set $\{x,y\}\subseteq\Omega$ is a
			\textbf{conflicting pair} for $c_0$ if there exist
			$S\subseteq D$ and $e\in c_0(S)$ such that
			\[
			d(x)=d(y)\quad\text{for every }d\in S,
			\qquad e(x)\ne e(y).
			\]
			Equivalently, $A(x,y)$ is not $c_0$-closed.
			Write
			\[
			\mathscr P^+
			=
			\bigl\{
			\{x,y\}\subseteq\Omega:
			x\ne y,\ 
			c_0(A(x,y))\ne A(x,y)
			\bigr\}
			\]
			for the family of conflicting pairs.
		\end{enumerate}
	\end{definition}

	To verify the equivalence in \textup{(ii)}, suppose that $S\subseteq A(x,y)$ and $e\in c_0(S)$ but $e\notin A(x,y)$.
	Monotonicity gives $e\in c_0(A(x,y))\setminus A(x,y)$.
	Conversely, if $A(x,y)$ is not closed, extensivity supplies
	$e\in c_0(A(x,y))\setminus A(x,y)$; taking $S=A(x,y)$ gives a violated dependency.

	\begin{definition}[Simultaneous faithful repair and joint edit compatibility]
		\label{def:joint-edit-compatibility}
		Fix the literal tuple-edit contract of
		Definitions~\ref{def:safe-formation-envelope} and~\ref{def:joint-certificates}.
		\begin{enumerate}[label=\textup{(\roman*)},leftmargin=*]
			\item \textbf{Faithful repair.}
			A set $E\subseteq\mathcal U$ is a \emph{simultaneous faithful repair} if
			$X(E)$ has determination closure $c_0$, realizes the prescribed local
			states, and meets all declared commuting-separator exactness and
			full-endpoint-range requirements.
			
			\item \textbf{Joint compatibility.}
			For one choice $C_\gamma\in\mathscr C_\gamma$ per requirement
			$\gamma\in\Gamma$, put
			\[
			E_{\mathbf C}=\bigcup_{\gamma\in\Gamma}\sigma(C_\gamma).
			\]
			The choice is \emph{jointly edit-compatible} if it activates no
			conflicting pair, namely $\sigma(P)\not\subseteq E_{\mathbf C}$ for all
			$P\in\mathscr P^+$.
			This test includes every pair activated by the selected edits.
		\end{enumerate}
	\end{definition}

	\begin{theorem}[Simultaneous faithful repair: criterion and minimum cost]
		\label{thm:joint-faithful-repair}
		Assume the prescribed endpoint alphabets and the separator factorizations on the allowed envelope specified above.
		\begin{enumerate}[label=\textup{(\roman*)},leftmargin=*]
			\item \emph{Repair criterion.}
			A set of tuple edits $E\subseteq\mathcal U$ is a simultaneous
			faithful repair exactly when both conditions below hold.
			
			\begin{equation}\label{eq:joint-conditions}
				\boxed{\begin{aligned}
						&\forall\gamma\in\Gamma\ \exists C\in\mathscr C_\gamma:
						\sigma(C)\subseteq E,\\
						&\forall P\in\mathscr P^+:\quad \sigma(P)\not\subseteq E.
				\end{aligned}}
			\end{equation}
			
			\item \emph{Minimum cost and native realization.}
			Thus the exact minimum native repair cost, without any safe-envelope assumption, is
			\begin{equation}\label{eq:unconditional-cost}
				\tau_{\rm FSN}=
				\min\{w(E):E\subseteq\mathcal U\text{ satisfies }
				\eqref{eq:joint-conditions}\}.
			\end{equation}
			An empty minimum has value $+\infty$.
			In the permitted commutant interpretation above, every feasible repair also has the displayed native output.
		\end{enumerate}
	\end{theorem}

	\begin{proof}
		\leavevmode
		\begin{enumerate}[label=\textup{(\roman*)},leftmargin=*]
			\item \emph{Repair criterion.}
			Because the coverage list is nonempty, satisfying its certificates
			ensures that the repaired source is nonempty, even if the initial
			source was empty. The activation law \eqref{eq:activation} shows that
			a candidate certificate $C$ is
			present in the repaired source exactly when $\sigma(C)\subseteq E$.
			The classical finite agreement-set criterion \cite{BeeriDFS1984} says that, once all
			positive dependencies of $c_0$ hold, exact recovery $c_{X(E)}=c_0$ is equivalent to
			occurrence of
			every proper meet-irreducible $M$ as an exact agreement set.
			Local state coverage is the second class of certificates.
			Connectivity of a graph on its prescribed vertices is equivalent to having an edge
			across every nontrivial cut, giving the third class and the stated commuting-native
			algebra equality.
			Finally, the repaired source violates some positive dependency precisely when it
			contains both members of a pair $P\in\mathscr P^+$, equivalently $\sigma(P)\subseteq
			E$.
			These equivalences give \eqref{eq:joint-conditions}.
			
			\item \emph{Minimum cost and native realization.} Minimizing $w(E)$ over the feasible edit sets gives \eqref{eq:unconditional-cost}, with value $+\infty$ if none exist. Every feasible $X(E)$ is nonempty, so the stated commutant construction supplies its native output.\qedhere
		\end{enumerate}
	\end{proof}

	Consider one native package whose admissible full states, with
	descriptors $(A,B)$, are
	\[
	x_0=(0,0),\qquad x_1=(1,1),\qquad x_2=(1,0),
	\]
	with $x_0$ initially present and $x_1,x_2$ independently insertable by unary native tuple edits.
	Let $c_0$ be generated by $A\to B$.
	Its proper meet-irreducibles are $\varnothing$ and $\{B\}$.
	The pair $(x_0,x_1)$ witnesses the first and $(x_0,x_2)$ the second.
	Each edit on its own preserves $A\to B$. Both edits together admit
	the pair $(x_1,x_2)$, which agrees on $A$ and differs on $B$.
	Every edit set supplying both witnesses therefore violates the required
	dependency.

	\subsection{Repair under safe envelopes}
	\label{subsec:positive-safe-repair}

	If the entire envelope satisfies the target dependencies, the remaining
	cost is the weight of the union of the chosen certificate supports.
	Shared tuple edits incur a single charge.

	\begin{definition}[Sound formation envelope for a target recovery closure]
		\label{def:positive-safe}
		An allowed envelope $\Omega$ is \textbf{$c_0$-sound} if every positive dependency of $c_0$ holds on $\Omega$, equivalently $\mathscr P^+=\varnothing$.
		Equivalently, every permitted set of tuple insertions preserves
		all required dependencies, including when all edits are used together.
	\end{definition}

	\begin{theorem}[Minimum simultaneous repair cost in a sound formation envelope]
		\label{thm:safe-joint-cost}
		Suppose the allowed full formation envelope is $c_0$-sound.
		\begin{enumerate}[label=\textup{(\roman*)},leftmargin=*]
			\item \emph{Existence and exact minimum.}
			A simultaneous faithful, local-state-exact and commuting-native-exact repair exists if and only if every $\mathscr C_\gamma$ is nonempty.
			When it exists, its exact minimum native tuple-edit cost is
			\begin{equation}\label{eq:safe-cost}
				\quad\tau_{\rm FSN}=
				\min_{\,C_\gamma\in\mathscr C_\gamma\ (\gamma\in\Gamma)}
				w\!\left(\bigcup_{\gamma\in\Gamma}\sigma(C_\gamma)\right).\quad
			\end{equation}
			
			\item \emph{Sharp bounds.}
			If $\lambda_\gamma=\min_{C\in\mathscr C_\gamma}w(\sigma(C))$,
			then
			\begin{equation}\label{eq:safe-sandwich}
				\quad\max_{\gamma\in\Gamma}\lambda_\gamma
				\ \le\ \tau_{\rm FSN}\ \le\
				\sum_{\gamma\in\Gamma}\lambda_\gamma.\quad
			\end{equation}
			Both inequalities are sharp as general finite certificate-support inequalities:
			maximal sharing attains the lower bound, while pairwise disjoint indispensable supports attain the upper bound.
		\end{enumerate}
	\end{theorem}

	\begin{proof}\leavevmode
		\begin{enumerate}[label=\textup{(\roman*)},leftmargin=*]
			\item \emph{Existence and exact minimum.}
			Every repaired source is a subset of the full envelope; hence positive soundness is automatic for all edits.
			Any successful repair contains a certificate for each requirement, so its edit set contains
			the union of their supports.
			Conversely, selecting one candidate certificate per requirement and inserting the
			union of its missing native tuples activates every selected certificate.
			Additional states cannot invalidate a negative recovery witness, local coverage, or
			connectivity of graphs on fixed vertices, while positive soundness
			is protected by the envelope hypothesis.
			Thus the union is a valid repair and its minimum cost is exact.
			Existence follows by taking all allowed edits.
			\item \emph{Sharp bounds.}
			The lower bound in \eqref{eq:safe-sandwich} holds for every feasible repair; choosing
			independently cheapest certificates and inserting their union proves the upper bound.
			The stated extreme support arrangements show sharpness at the certificate level.\qedhere
		\end{enumerate}
	\end{proof}
	
	In the next two corollaries, (S) denotes prescribed local-state coverage and (\Esep) the specified exact separator contract.
	
	\begin{corollary}[Simultaneous native cycle repair]
		\label{cor:cycle-all-witness}
		Suppose the native atlas is a connected finite graph with a protected spanning tree,
		$\beta$ unprotected chord relations, and a $c_0$-sound full envelope.
		Suppose the original global source is nonempty, covers all prescribed local states,
		and has exact commuting native separator graphs on fixed vertex sets.
		Let $q$ be the number of proper meet-irreducibles of $c_0$ not yet globally witnessed.
		Assume each such meet-irreducible has a paired witness over the protected tree, whose
		component tuples are all permitted in the allowed chord relations.
		Then native tuple additions of total unit cost at most
		\begin{equation}\label{eq:cycle-q}
			\quad 2\beta q\quad
		\end{equation}
		restore the entire intended recovery closure while preserving every
		positive dependency, local-state coverage, and native separator
		exactness.  In particular, for one closing edge, at most $2q$ native
		tuple additions suffice.  The examples in \cref{thm:multi-native-repair} attain the
		single-witness factor $2$.
	\end{corollary}
	
	\begin{proof}
		Choose a protected-tree pair witness for each missing
		meet-irreducible.  By the source-level target repair argument of
		\cref{thm:multi-native-repair}, at most two missing tuples per chord
		activate each chosen pair.  Insert the union of these additions; its
		unit cost is at most $2\beta q$.  The $c_0$-sound full envelope
		prevents any unwanted positive violation.  Existing nonimplication
		witnesses, local-state coverage, and connectedness of each separator
		graph persist under addition.  The complete meet-irreducible witness
		family therefore yields exact recovery.
	\end{proof}
	
	\begin{corollary}[Additive global-state completion bound]
		\label{cor:row-completion}
		Suppose instead that one may add whole global states from a fixed
		$c_0$-sound ambient source $\Omega$ already satisfying all intended
		local coverage and separator conditions. Assume moreover that
		$c_\Omega=c_0$, equivalently under soundness that $\Omega$ contains
		an exact-agreement witness for every proper meet-irreducible of $c_0$.
		Use a nonempty list of prescribed local states as in
		\cref{def:joint-certificates}, and fix $K\subseteq\Omega$.
		Let $q(K)$ count the proper meet-irreducibles still lacking witnesses.
		Let $s(K)$ count missing prescribed states over all sites.
		For each separator, use its graph on the fixed vertex set and let
		$r_i$ be its graphic rank. Put $\delta_i(K)=r_i(\Omega)-r_i(K)$.
		The minimum $\mu(K)$ counts states added from $\Omega\setminus K$
		to meet (F), (S), and (\Esep) together. The following bounds use
		zero for the maximum if there are no separator tests.
		\begin{equation}\label{eq:row-bound}
			\quad\max_i\delta_i(K)\le\mu(K)\le
			2q(K)+s(K)+\sum_i\delta_i(K).\quad
		\end{equation}
		The bounds count added global states. A tuple insertion can
		admit several of these states at once.
	\end{corollary}
	
	\begin{proof}
		For each missing meet-irreducible choose two witnessing global states
		in $\Omega$; for each missing local state choose an ambient extension.
		Then extend to full rank in each separator's graphic matroid.  The
		basis-extension axiom adds at most $\delta_i(K)$ states for separator
		$i$, and previous ranks can only increase.  All added states are in
		the sound envelope.  Conversely, one state adds at most one graphic
		rank unit to a fixed separator, so at least $\delta_i(K)$ states are
		needed for every $i$.
	\end{proof}
	
	\begin{example}[A sharp three-state simultaneous completion]
		\label{ex:three-row}
		Let the two endpoint observations be binary $A,B$, take the common
		separator to be constant, and use the $c_0$-closure with no
		nontrivial functional dependencies.  Starting with
		$K=\{(0,0)\}$ inside the sound four-state envelope
		$\Omega=\{0,1\}^2$, add precisely the states
		$(0,1)$ and $(1,0)$.  The resulting pairs witness the proper
		meet-irreducibles $\{A\}$ and $\{B\}$, cover both local binary
		state sets, and provide a connected three-edge bipartite separator
		graph on four prescribed vertices.  Two additions are necessary:
		any connected graph on those four vertices has at least three edges.
		Thus the lower graphic-rank bound in
		\eqref{eq:row-bound} is attained.  The same example can be read
		as a small relational join, as commuting diagonal observables, or as
		a finite two-package interface; none of these interpretations
		changes the proof.
	\end{example}

	\subsection{Conflict hypergraphs and compatible certificate selections}
	\label{sec:selected-safe-repairs}

	A tuple insertion can activate several states.
	Even when the full envelope is unsafe, its minimal forbidden edit sets determine which certificate supports can be chosen together
	\cite{AbramskyCaru,BeeriDFS1984,Schrijver2003}.
	Safe edit sets are independent sets of the resulting hypergraph, as in functional-dependency repair and hitting-set diagnosis
	\cite{LivshitsKimelfeldWijsen2021,LivshitsKimelfeldRoy2018,Reiter1987}.

	\begin{definition}[Edit conflicts, safe edit sets, and certificate supports]
		\label{def:edit-conflicts}
		Retain the finite native atlas, target closure $c_0$, edit ground set
		$\mathcal U$, activation supports $\sigma(x)$, and joint requirement
		families $\mathscr C_\gamma$ of
		Definitions~\ref{def:safe-formation-envelope} and~\ref{def:joint-certificates}.
		Assume $X(\varnothing)$ obeys all positive dependencies of $c_0$.
		\begin{enumerate}[label=\textup{(\roman*)},leftmargin=*]
			\item \textbf{Minimal conflicts.}
			Define the \textbf{native edit-conflict hypergraph}
			\[
			\mathcal H^+_{\!\sigma}
			=\min_{\subseteq}\bigl\{
			\sigma(x)\cup\sigma(y):\{x,y\}\in\mathscr P^+\bigr\},
			\]
			where $\min_{\subseteq}$ discards all proper supersets and duplicates.
			Its hyperedges are nonempty.
			
			\item \textbf{Safe edit sets.}
			The \textbf{safe edit complex} is
			\[
			\mathcal I_\sigma=\{E\subseteq\mathcal U:
			H\not\subseteq E\text{ for all }
			H\in\mathcal H^+_{\!\sigma}\}.
			\]
			The maximal elements of $\mathcal I_\sigma$ are the
			\textbf{maximal safe edit sets}.
			
			\item \textbf{Certificate supports.}
			For each requirement $\gamma$ retain the inclusion-minimal supports
			\[
			\mathcal S_\gamma
			=\min_{\subseteq}\{\sigma(C): C\in\mathscr C_\gamma\}.
			\]
		\end{enumerate}
	\end{definition}

	\begin{theorem}[Faithful repair from compatible supports: existence and minimum cost]
		\label{thm:selected-safe-repair}
		Retain Definition~\ref{def:edit-conflicts} and the finite readout and fixed-vertex
		hypotheses of \cref{thm:joint-faithful-repair}.
		Then the following statements are equivalent.
		\begin{enumerate}[label=(\roman*)]
			\item Some native edit set repairs all recovery, local-coverage,
			and commuting-separator requirements while preserving every positive
			dependency of $c_0$.
			
			\item There are supports $S_\gamma\in\mathcal S_\gamma$
			for all $\gamma\in\Gamma$ such that
			$\bigcup_{\gamma\in\Gamma}S_\gamma\in\mathcal I_\sigma$.
			
			\item Some maximal safe edit set $M\in\operatorname{Max}(\mathcal I_\sigma)$
			contains at least one support from $\mathcal S_\gamma$ for every $\gamma$.
		\end{enumerate}
		For nonnegative edit weights the exact minimum is
		\begin{equation}\label{eq:selected-cost}
			\ \tau_{\rm FSN}=
			\min_{\substack{S_\gamma\in\mathcal S_\gamma\ (\gamma\in\Gamma)\\
					\forall H\in\mathcal H^+_{\!\sigma},\
					H\not\subseteq\bigcup_\gamma S_\gamma}}
			w\!\left(\bigcup_\gamma S_\gamma\right).\
		\end{equation}
		Equivalently,
		\begin{equation}\label{eq:shell-decomposition}
			\ \tau_{\rm FSN}=
			\min_{M\in\operatorname{Max}(\mathcal I_\sigma)}
			\min_{\substack{S_\gamma\in\mathcal S_\gamma\\
					S_\gamma\subseteq M\ \forall\gamma}}
			w\!\left(\bigcup_\gamma S_\gamma\right).\
		\end{equation}
		An inner minimum over an empty candidate family has value $+\infty$.
		No $c_0$-soundness of the full envelope $\Omega$ is required.
	\end{theorem}
	
	\begin{proof}
		By the native activation law, a conflicting pair $\{x,y\}$
		is present in $X(E)$ exactly when $\sigma(x)\cup\sigma(y)\subseteq E$.
		Hence the repair preserves all positive dependencies exactly when $E\in\mathcal I_\sigma$.
		The desired negative witnesses, prescribed local states, and native
		separator cut edges are present precisely when each $\gamma$ has an
		activated certificate.
		\begin{itemize}
			\item \emph{\textup{(i)} $\Longleftrightarrow$ \textup{(ii)}.}
			Every activated certificate support in a successful repair $E$ contains a minimal candidate support for the same requirement. Choose one such $S_\gamma\subseteq E$ for each $\gamma$.
			Their union $E_*=\bigcup_\gamma S_\gamma$ is still safe, because
			$\mathcal I_\sigma$ is downward closed, and still activates all
			selected certificates.
			Thus it is a successful repair, of no larger cost than $E$.
			Conversely, any safe union of certificate supports is a successful repair.
			\item \emph{\textup{(ii)} $\Longleftrightarrow$ \textup{(iii)}.}
			Every safe set in a finite hypergraph lies in a maximal safe set. Conversely,
			the union of supports contained in a maximal safe set is safe by downward closure.
			\item \emph{Minimum cost.}
			The safe subunion constructed above has no greater cost than the given repair.
			Minimizing over all such unions gives \eqref{eq:selected-cost}. Grouping them by the maximal safe sets containing them gives \eqref{eq:shell-decomposition}.
		\end{itemize}
	\end{proof}
	
	\begin{corollary}[Sharp cost bounds within a maximal safe edit set]
		\label{cor:shell-cost-bounds}
		If a maximal safe edit set $M$ contains a candidate support for every $\gamma$,
		let $\lambda_\gamma(M)=
		\min\{w(S):S\in\mathcal S_\gamma,\ S\subseteq M\}$.
		Then a simultaneous repair exists inside $M$, and its minimum cost
		$\tau(M)$ satisfies
		\[
		\quad
		\max_\gamma\lambda_\gamma(M)
		\le\tau(M)\le\sum_\gamma\lambda_\gamma(M).\quad
		\]
		The extreme support-sharing examples of \cref{thm:safe-joint-cost} show that the two
		inequalities cannot be improved as universal certificate-support estimates.
		The same assertion can hold with $\Omega$ globally unsound.
	\end{corollary}
	
	\begin{proof}
		Every subset of a maximal safe edit set is safe.
		The least expensive supports for each requirement therefore have a safe union,
		yielding the upper bound; every valid union contains a support for each
		requirement, yielding the lower bound.
	\end{proof}
	
	\begin{example}[A globally unsafe envelope with an exact repair]
		\label{ex:unsafe-but-repairable}
		Consider one finite source with descriptor readouts $(A,B)$ and
		initial state $x_0=(0,0)$.  Three independently editable admissible
		states are
		\[
		x_1=(1,1),\qquad x_2=(1,0),\qquad x_3=(2,0),
		\]
		with unit-cost native edits $e_1,e_2,e_3$, respectively.  Let $c_0$
		be the closure generated by $A\to B$.  The full envelope is unsafe,
		since $x_1,x_2$ agree on $A$ but disagree on $B$.  Nevertheless
		$E=\{e_1,e_3\}$ is safe.  The pair $\{x_0,x_1\}$ has exact
		agreement $\varnothing$, and $\{x_0,x_3\}$ has exact agreement
		$\{B\}$.  These are precisely the proper meet-irreducibles of
		$c_0$, so the three-state repaired source is faithful.  Both edits
		are indispensable, and its exact repair cost is two.  
	\end{example}

	\subsection{Diversity criteria and edit-conflict complexity}
	\label{subsec:repair-compatibility-thresholds}

	\Cref{thm:diverse-compatibility} bounds the number of support choices
	needed for simultaneous repair when every conflict is detected by a pair.
	\Cref{thm:finite-conflict-clutter-realization} shows why this detection
	hypothesis matters: a star can realize any antichain of conflicts.

	Three editable unary factors on binary states exhibit a conflict invisible to
	pairwise tests. Initially only $(0,0,0)$ is present, and edit $e_i$
	admits value $1$ at factor $i$. A protected observer has the eight
	triples as states. A protected factor equates its state to the triple
	of unary values, and both descriptors are read at the observer.
	Using all three edits admits $(1,1,1)$. Assign descriptor readouts
	$(A,B)$ so that $(1,0,0)$ has values $(1,1)$, $(1,1,1)$ has values
	$(1,0)$, and every other permitted state has a distinct fresh
	$A$-value.  The required dependency $A\to B$ is violated only
	when all three edits are inserted.  The three edits form a size-three hyperedge of $\mathcal H^+_{\!\sigma}$.
	Each two-edit subset is safe.

	\begin{definition}[Binary-detectable support compatibility]
		\label{def:binary-compatibility}
		A family $\mathcal S_1,\ldots,\mathcal S_q$ of nonempty	inclusion-minimal candidate support families has \textbf{binary-detectable support compatibility} (or is \textbf{binary-detectable}) if every candidate support is individually safe and, for every finite selection of their options, an unsafe union contains an unsafe union of two selected options.
		An option $S$ in one family forbids an option $T$ in another if
		$S\cup T\notin\mathcal I_\sigma$.
	\end{definition}

	\begin{theorem}[Diversity criterion for jointly safe repair supports]
		\label{thm:diverse-compatibility}
		Suppose there are $q\ge1$ certificate requirements, each with at
		least $t$ distinct candidate supports, and these families are
		binary-detectable.  Suppose that, for every option $S$ in one
		family and every different requirement, $S$ forbids at most
		$\Delta$ options in that other requirement.  If
		\begin{equation}\label{eq:diversity-threshold}
			\quad t>(q-1)\Delta,\quad
		\end{equation}
		then a simultaneous positive-sound repair exists.
		If each candidate support has at most $b$ editable tuples of
		unit cost, a repair can be chosen with at most $bq$ native tuple
		insertions.  No positive-safety assumption on the full permitted
		formation envelope is needed.
	\end{theorem}

	\begin{proof}
		Choose options successively, one from each requirement.  After
		choosing $j<q$ options, at most $j\Delta$ options in the next
		family are forbidden.  The numerical hypothesis ensures an
		available option.  All selected options are individually safe and
		pairwise compatible.  Binary detectability makes their entire union
		safe.  By \cref{thm:selected-safe-repair}, this union is a
		valid simultaneous repair.  Its cardinality is at most the sum of
		the $q$ support sizes, proving the quantitative bound.
	\end{proof}

	\begin{corollary}[Simultaneous chord repair with jointly safe witness supports]
		\label{cor:unsafe-cycle-repair}
		Under the protected-spanning-tree conditions of
		Corollary~\ref{cor:cycle-all-witness}, drop the requirement that the
		full envelope be $c_0$-sound.
		Let $q$ count the missing proper meet-irreducible recovery requirements.
		\begin{enumerate}[label=\textup{(\roman*)},leftmargin=*]
			\item \emph{Cost of repair.}
			Suppose the $q$ requirements admit a jointly safe selection of
			protected-tree paired-witness lifts.
			Then all required witnesses can be restored simultaneously while
			preserving positive dependencies, pre-existing local-state coverage,
			and finite commuting-native separator exactness, with at most
			$
			2\beta q\quad
			$
			unit-cost native tuple insertions.
			\item \emph{Existence of a safe selection.}
			A sufficient condition for such a selection is that the lift-support families satisfy \cref{thm:diverse-compatibility}, with every lift supported on at most $2\beta$ editable chord tuples.
		\end{enumerate}
	\end{corollary}

	\begin{proof}
		\leavevmode
		\begin{enumerate}[label=\textup{(\roman*)},leftmargin=*]
			\item \emph{Cost of repair.}
			By the protected-tree lift hypothesis, each selected paired witness
			requires at most one native tuple per witness per editable chord,
			so its support has at most $2\beta$ elements.
			The joint safety hypothesis protects every positive dependency, regardless of other states activated by these edits.
			Existing coverage and native separator connectivity persist under additional global states and fixed separator readouts.
			The selected paired witnesses establish	all missing proper meet-irreducibles.
			Apply \cref{thm:selected-safe-repair} and sum the support sizes.
			
			\item \emph{Existence of a safe selection.}
			The diversity condition supplies a safe selection by
			\cref{thm:diverse-compatibility}.\qedhere
		\end{enumerate}
	\end{proof}

	\paragraph{Arbitrary edit conflicts.}
	The realization of an arbitrary edit-conflict hypergraph uses the
	hidden-variable encoding of a finite relation \cite[Section~2.2, Definition~7]{BacchusChenVanBeekWalsh2002}
	and a monotone Boolean description of forbidden supports
	\cite[Section~3]{Green2007}.

	\begin{theorem}[Realizing any family of minimal edit conflicts on a star]
		\label{thm:finite-conflict-clutter-realization}
		Let $\mathcal U$ be finite and let $\mathcal H$ be any antichain of
		nonempty subsets of $\mathcal U$, with $\mathcal H=\varnothing$
		permitted. There is a finite atlas of realized records as in
		Definition~\ref{def:multi-native-atlas}. It has one independently editable
		bit factor per $u\in\mathcal U$. Protected binary constraints form a star
		whose main observer has $1+|\mathcal H|$ records. The initial global
		state set is nonempty. For the single dependency $A\to B$, the
		hypergraph of minimal edit conflicts is exactly $\mathcal H$. Consequently the allowed
		native-record interpretation and this binary tree topology alone
		impose no uniform bound on edit-conflict rank.
	\end{theorem}
	
	\begin{proof}
		Take observer states
		$\mathcal T=\{\varnothing\}\cup\mathcal H$ and, for each
		$u\in\mathcal U$, a bit site with state set $\{0,1\}$.
		The protected observer--bit relation is
		\[
		R_u=\{(S,b)\in\mathcal T\times\{0,1\}:
		b=\mathbf 1_{\{u\in S\}}\}.
		\]
		The unary factor at bit site $u$ initially admits only $0$;
		its one allowed edit, identified with $u$, admits $1$ as well.
		The full envelope has one global state $x_S$ per $S\in\mathcal T$,
		since all bit values are forced by the observer. Its activation support
		is $\sigma(x_S)=S$, so after edits $E\subseteq\mathcal U$,
		\[
		X(E)=\{x_{\varnothing}\}\cup
		\{x_H:H\in\mathcal H,\ H\subseteq E\}.
		\]
		In particular the initial state $x_{\varnothing}$ survives every edit.
		Read $A=0$ on every observer state, and put
		$B(x_{\varnothing})=0$, $B(x_H)=1$ for $H\in\mathcal H$.
		The positive-conflict pairs are exactly
		$\{x_{\varnothing},x_H\}$, whose support unions are $H$.
		Every edge of the antichain $\mathcal H$ is inclusion-minimal.
		Thus $\mathcal H^+_{\!\sigma}=\mathcal H$.
		
		For the native-record interpretation, realize every nonempty local
		state set $Z$ as the characters of the multiplication algebra
		$\mathbb C^Z$ on $\ell^2(Z)$, tagged by the corresponding
		access--commutant datum as in
		Proposition~\ref{prop:spectral-observation-interpretation}.
		Keep the displayed readouts and protected relations as the specified
		observational formation data. This is precisely the permitted finite
		native-record class. Taking $\mathcal H=\{\mathcal U\}$ on
		arbitrarily large nonempty edit grounds proves the final assertion.
	\end{proof}
	
	\subsection{Recovery compatibility across interacting interfaces}
	\label{subsec:recovery-compatibility-complex}
	
	With positive dependencies, coverage, and separator requirements fixed,
	the witness families attainable by one repair form a downward-closed
	family. Using activation constraints from constraint satisfaction and
	database repair \cite{BacchusChenVanBeekWalsh2002,LivshitsKimelfeldWijsen2021},
	\cref{thm:cyclic-recovery-compatibility-universality} realizes each such
	family on triangles sharing an edge, with one editable tuple per chord.

	\begin{definition}[Recovery compatibility complex]
		\label{def:recovery-compatibility-complex}
		Fix the finite allowed-edit contract of
		\cref{def:safe-formation-envelope}, a target closure $c_0$,
		and a collection $\mathscr Q$ of proper meet-irreducibles of
		$c_0$ not witnessed by the initial source.
		Assume that the	initial source is nonempty, satisfies every positive dependency of $c_0$, and is realized under the stipulated native law.
		Any prescribed local-state coverage and commuting separator tests are also assumed to hold initially and are to remain
		satisfied.
		\begin{enumerate}[label=\textup{(\roman*)},leftmargin=*]
			\item \textbf{Joint recovery.}
			A subset $\mathscr S\subseteq\mathscr Q$ is \emph{jointly recoverable}
			if some edit set $E$ makes $X(E)$ nonempty, preserves all positive
			dependencies of $c_0$, meets these fixed coverage and separator
			clauses, and supplies an exact-agreement pair for every
			$M\in\mathscr S$.  Require in addition the native realization
			of $X(E)$ in the stipulated native law.
			
			\item \textbf{Compatibility complex.}
			The family $\mathcal K_{\rm rec}$ of jointly recoverable subsets
			is the \textbf{recovery compatibility complex}.
			It is downward closed: one repair witnessing $\mathscr S$ also witnesses every subset of $\mathscr S$.
		\end{enumerate}
	\end{definition}

	\begin{theorem}[Universality of interacting-cycle recovery compatibility]
		\label{thm:cyclic-recovery-compatibility-universality}
		Let $q\ge1$, and let $\mathcal K$ be a downward-closed family of
		subsets of $[q]=\{1,\ldots,q\}$ containing the empty set.
		Some singletons may be absent.
		There is a finite multivalued native formation atlas with the
		following properties.
		\begin{enumerate}[label=\textup{(\roman*)}]
			\item \emph{Cycle construction.} Its constraint graph has positions
			$o,h,v_1,\ldots,v_q$, a protected spanning tree with edges
			$oh,ov_1,\ldots,ov_q$, and exactly $q$ independently editable
			chords $hv_i$.  Each fundamental cycle is a triangle, and all
			of them share $oh$.  Every relation is binary.  Edit $e_i$
			inserts exactly one native tuple into the chord $hv_i$.
			\item \emph{Initial recovery witnesses.} For the descriptor set $D=\{A,B,d_1,\ldots,d_q\}$,
			let $c_0$ be the closure generated by $A\to B$.  The initial
			source is nonempty and satisfies $c_0$; it witnesses exactly
			the proper meet-irreducibles
			\[
			M_A=D\setminus\{A,B\},\qquad M_B=D\setminus\{A\},
			\]
			while the remaining proper meet-irreducibles are
			$M_i=D\setminus\{d_i\}$, $i\in[q]$.
			Each $M_i$ has an exact paired witness over the protected tree.
			\item \emph{Effect of editing the chords.} For every $E\subseteq[q]$, the source $X(E)$ is nonempty,
			and its native realization is explicit in the fixed
			access--commutant realization family.  Moreover,
			\begin{equation}\label{eq:cyclic-compatibility-equivalences}
				\begin{split}
					X(E)\text{ satisfies every implication of }c_0
					&\iff E\in\mathcal K,\\
					X(E)\text{ has an exact-}M_i\text{ witness}
					&\iff i\in E.
				\end{split}
			\end{equation}
			The initial prescribed local-state coverage and one nontrivial,
			fixed commuting-separator exactness test remain satisfied for
			all $E$.
			\item \emph{Compatibility and minimum cost.} Thus the recovery compatibility complex of the $q$
			missing meet-irreducibles is exactly $\mathcal K$.  For
			nonnegative edit weights $w_i$, the minimum cost of a repair
			preserving the positive, native, coverage and named separator
			clauses while witnessing every $M_i$ with $i\in S$ is
			\[
			\sum_{i\in S}w_i\quad\text{if }S\in\mathcal K,
			\]
			and no such repair exists if $S\notin\mathcal K$.
			In particular, a fully recovery-faithful repair exists precisely
			when $[q]\in\mathcal K$.
		\end{enumerate}
	\end{theorem}
	\begin{proof}\leavevmode
		\begin{enumerate}[label=\textup{(\roman*)},leftmargin=*]
			\item \emph{Cycle construction.}
			Let $\mathcal H$ be the antichain of inclusion-minimal nonfaces
			of $\mathcal K$. Its members are nonempty and may be singletons.  Write
			$\mathbf0$ for the zero vector in $\{0,1\}^q$ and
			$\mathbf e_i$ for the $i$th unit vector.
			
			Take the observer-record set at $o$ to be the disjoint union
			\[
			T=\{b_0,b_1,b_2,b_3\}\sqcup
			\{u_i:i\in[q]\}\sqcup\{t_H:H\in\mathcal H\}.
			\]
			Each $r\in T$ has a hub value $\eta(r)\in\{0,1\}$,
			a bit vector $\epsilon(r)\in\{0,1\}^q$, and the descriptor
			readouts $(A(r),B(r),d_1(r),\ldots,d_q(r))$ given by
			\[
			\begin{array}{c|c|c|c|c}
				r & \eta(r) & \epsilon(r) & (A(r),B(r)) & (d_j(r))_{j=1}^q\\ \hline
				b_0 & 0 & \mathbf0 & (0,0) & \mathbf0\\
				b_1 & 1 & \mathbf0 & (1,1) & \mathbf0\\
				b_2 & 0 & \mathbf0 & (2,0) & \mathbf0\\
				b_3 & 1 & \mathbf0 & (0,0) & \mathbf0\\
				u_i & 0 & \mathbf e_i & (0,0) & \mathbf e_i\\
				t_H & 0 & \mathbf1_H & (0,1) & \mathbf0
			\end{array}
			\]
			The hub site $h$ has records $\{0,1\}$, and each bit site
			$v_i$ has records $\{0,1\}$.  On the protected tree use
			\[
			F_{oh}=\{(r,\eta(r)):r\in T\},\qquad
			F_{ov_i}=\{(r,\epsilon_i(r)):r\in T\}.
			\]
			The initial and allowed chord relations are, respectively,
			\[
			F_{hv_i}=\{(0,0),(1,0)\},\qquad
			\widehat F_{hv_i}=F_{hv_i}\cup\{(0,1)\}.
			\]
			Here $e_i$ means insertion of $(0,1)$ into $F_{hv_i}$.
			A global state is uniquely determined by its observer record.
			Its exact activation rule is therefore
			\begin{equation}\label{eq:cyclic-compatibility-activation}
				X(E)=\{x_{b_0},x_{b_1},x_{b_2},x_{b_3}\}
				\cup\{x_{u_i}:i\in E\}
				\cup\{x_{t_H}:H\in\mathcal H,\ H\subseteq E\}.
			\end{equation}
			All descriptor readouts occur at $o$, so global consistency is
			automatic.
			\item \emph{Initial recovery witnesses.}
			The target closure generated by $A\to B$ has precisely the proper meet-irreducibles listed in the statement:
			its closed subsets are those containing $B$ whenever they
			contain $A$.  The pairs $(b_0,b_1)$ and $(b_0,b_2)$ have exact
			agreement sets $M_A$ and $M_B$, respectively.  The only
			observer records with $d_i=1$ are $u_i$; therefore an
			exact-$M_i$ pair exists in $X(E)$ if and only if $u_i$ is
			present, since $(b_0,u_i)$ realizes that agreement pattern.
			Every $u_i$ is available over the protected tree, proving the
			protected-tree witness assertion.
			
			\item \emph{Effect of editing the chords.}
			In the baseline and among all $u_i$, equal $A$-values have
			equal $B$-values: $A=0$ always gives $B=0$, while $b_1$
			and $b_2$ have distinct $A$-values.  Each record $t_H$ has
			$A=0$, $B=1$, and its presence violates $A\to B$ together
			with $b_0$.  Thus $X(E)$ is $c_0$-sound precisely when no
			$H\in\mathcal H$ is contained in $E$, equivalently when
			$E\in\mathcal K$.  This proves
			\eqref{eq:cyclic-compatibility-equivalences}.
			
			For a fixed nontrivial commuting separator test, take the
			left readout $\ell=A$ with prescribed range $\{0,1,2\}$,
			the right readout $r=\eta$ with range $\{0,1\}$, and the
			constant common separator.  The four baseline states supply
			the bipartite edges
			\[
			(0,0),\ (1,1),\ (2,0),\ (0,1),
			\]
			which form a spanning tree on its five prescribed vertices.
			The incidence graph stays connected under all insertions,
			so, in the permitted commutant interpretation,
			$\mathbb C[\ell]\cap\mathbb C[r]=\mathbb C[1]$ and both
			endpoint ranges remain exact.  Prescribe as local coverage
			the four baseline observer records, the two hub records, and
			the bit value $0$ at every $v_i$; they are covered already by
			$X(\varnothing)$ and remain covered.
			Each finite local record set is realized by the characters
			of its diagonal multiplication algebra.  For every $E$,
			the nonempty finite set $X(E)$ has the allowed native output
			$\mathbb C^{X(E)}$ acting on $\ell^2(X(E))$, with all the
			displayed readouts, as in
			\cref{def:safe-formation-envelope}.  This establishes the
			native realization for the stated law.
			
			\item \emph{Compatibility and minimum cost.}
			Finally, if a repair is sound and witnesses $M_i$ for every
			$i\in S$, the established equivalences force
			$S\subseteq E\in\mathcal K$, hence $S\in\mathcal K$.
			Conversely, $E=S$ meets all the requirements whenever
			$S\in\mathcal K$.  Nonnegative edit costs then give the
			stated exact minimum.  By
			\cref{thm:witness-globalization}, full recovery faithfulness
			requires and, under soundness, is implied by witnessing
			all $q$ missing meet-irreducibles; this is possible exactly
			when $[q]\in\mathcal K$.\qedhere
		\end{enumerate}
	\end{proof}
	
	\begin{corollary}[Unbounded obstruction size for simultaneous recovery repair]
		\label{cor:cyclic-recovery-no-bounded-local-test}
		For every integer $k\ge1$, there exists a finite native atlas
		with $k+1$ interacting fundamental cycles and $k+1$
		initially missing recovery witnesses such that every proper
		subfamily of those witnesses can be restored simultaneously,
		with positive dependencies, the specified coverage, one
		commuting separator and native realization preserved, but no
		repair can restore all $k+1$ witnesses under these clauses.
		Consequently there is no universal fixed $k$ for which the
		feasibility of every subfamily of at most $k$ missing
		exact-agreement requirements implies simultaneous faithful
		repair in this native-record class.
	\end{corollary}
	\begin{proof}
		Apply \cref{thm:cyclic-recovery-compatibility-universality}
		with $q=k+1$ and
		$\mathcal K=2^{[q]}\setminus\{[q]\}$.  Its unique minimal
		nonface is $[q]$.  Every proper subset is therefore
		jointly recoverable and the full set is not.
	\end{proof}
	
	\section{Fixed-law comparison and probe-relative non-collapse}
	\label{sec:fixed-law-noncollapse}
	
	The marked Level-II presentation has a faithful relational encoding
	(\cref{prop:rel-representation}). For heaps and access--commutant
	sources, \cref{thm:heap-commutant-comparison} compares chart and cut
	observations under the two native laws. For group sources, the allowed
	attachments determine which marked inputs can be distinguished.
	\Cref{thm:unbounded-native-resource-separation} compares sources with
	fixed representation resources.
	
	In \cref{thm:complete-compositional-invariants}, subgroup probes
	need not preserve every chart state; completeness concerns recovery
	or numerical outcomes. The subdirect tests in
	\cref{subsec:subdirect-attachment-tests} instead use a stronger fixed
	preservation contract and prescribed attachment class.
	
	\subsection{Comparison of heaps and access--commutant sources}
	\label{subsec:cross-native-comparison}
	
	A cut of a finite heap determines a canonical quotient of its two
	side images. \Cref{thm:heap-commutant-comparison} identifies the
	corresponding common algebra of diagonal observations and describes
	how it behaves under matching.
	
	\paragraph{The two fixed native laws.}
	A nonempty \emph{heap} can be viewed as a group with no chosen
	identity element. Formally, a heap is a set $A$ with a ternary operation
	$[x,y,z]$. The two identity laws are $[x,y,y]=x$ and $[x,x,y]=y$.
	The associativity law is $[[x,y,z],u,v]=[x,y,[z,u,v]]$.
	The native formation law $\mathsf J_{\mathrm{hp}}$ is defined as follows.
	Start with a finite diagram of nonempty finite heaps $A_i$.
	Each matching constraint $R_e$ is a subheap of the relevant chart
	product. A formation datum also specifies a nonempty finite tuple
	of generators $g_1,\ldots,g_s$ in the matching relation $W$.
	The law returns the generated subheap $U=\langle
	g_1,\ldots,g_s\rangle_{\mathrm{hp}}\subseteq W$, with its induced chart and descriptor
	maps.
	The native law can use the following packaged endpoint.
	\[
	\begin{gathered}
		\mathfrak P^{\mathrm{hp}}_U=
		(\{U\},U,U\times\{U\},\mathsf L_U,\mathsf R_U),\\
		\mathsf L_U(U)=\mathcal P(U),\quad
		\mathsf R_U(U)=\{V\subseteq U:V\ne\varnothing,\ V\text{ a subheap}\}.
	\end{gathered}
	\]
	The law is partial when the generator tuple is incompatible.
	All nonempty finite subheaps of $W$ occur by listing their elements as generators.
	
	The second law is precisely the preassigned $\Xi^{\mathrm{com}}_{\mathcal H}$ of
	\eqref{eq:standalone-commutant-law}:
	a concrete unital $*$-subalgebra $M\subseteq B(\mathcal H)$ is sent to its represented commutant endpoint $\mathfrak P_{M'}$.
	A native comparison retains the ambient Hilbert space and the
	concrete representation.
	For a finite nonempty source $U$, put
	\[
	\mathcal H_U=\ell^2(U),\qquad
	\mathcal A_U=\operatorname{diag}\mathbb C^U\subseteq B(\mathcal H_U),
	\qquad \mathcal I(U)=
	(\Xi^{\mathrm{com}}_{\mathcal H_U},
	(\mathrm{com},\mathcal A_U),\mathfrak P_{\mathcal A_U'}).
	\]
	The full diagonal algebra is maximal abelian, so
	$\mathcal A_U'=\mathcal A_U$.
	
	\begin{theorem}[Heap sources and their access--commutant observation realizations]
		\label{thm:heap-commutant-comparison}
		Let $U$ be a nonempty finite subheap of a product of finite heaps, with jointly
		injective chart maps $p_i:U\to X_i$ and retained finite-valued readouts $f_d$
		obtained from these charts.
		Fix a family $\mathcal T$ of full cuts, and at a cut $t:L\mid R$ write
		\[
		\alpha_t=\ker p_L,\quad \beta_t=\ker p_R,\quad
		C_t(U)=U/(\alpha_t\vee\beta_t),\quad
		U_L=p_L(U),\quad U_R=p_R(U).
		\]
		The following statements hold, without choosing an origin in any heap.
		\begin{enumerate}[label=(\roman*),nosep]
			\item \emph{Native strong cut reconstruction.}
			Every cut relation is a full subheap fiber product
			\begin{equation}
				U\ \cong\ U_L\times_{C_t(U)}U_R.
				\label{eq:heap-fiber}
			\end{equation}
			The quotient maps are heap homomorphisms.
			Every pair of side states with the same quotient value occurs.
			
			\item \emph{Exact commutant observation and recovery.}
			For a list of descriptors $S$ put $M_U(S)=\operatorname{alg}^*\{\mathbb
			C[f_d]:d\in S\}$ and $K_U(S)=M_U(S)'$ on $\mathcal H_U$.
			Then
			\begin{equation}
				d\in c_U(S)\iff \mathbb C[f_d]\subseteq M_U(S),\qquad
				K_U(S\cup T)=K_U(S)\cap K_U(T),
				\label{eq:recovery}
			\end{equation}
			while the two complete cut-side function algebras intersect in
			\begin{equation}
				\mathbb C[U_L]\cap\mathbb C[U_R]
				=\mathbb C[C_t(U)]\quad\text{inside }\mathbb C^U.
				\label{eq:cut-algebra}
			\end{equation}
			Chart-state surjectivity is equivalent to injectivity of the corresponding pullback $p_i^*: \mathbb C^{X_i}\to\mathbb C^U$.
			
			\item \emph{Unchanged separators.}
			If a named heap readout $y_t:U\twoheadrightarrow Y_t$ factors
			through both sides, it induces $\lambda_t:C_t(U)\twoheadrightarrow
			Y_t$.
			Strong reconstruction over the unchanged named interface and equality of its
			represented common-information algebra are equivalent:
			\begin{equation}
				U=U_L\times_{Y_t}U_R
				\iff\ \lambda_t\text{ is an isomorphism}
				\iff\ \mathbb C[U_L]\cap\mathbb C[U_R]
				=y_t^*\mathbb C^{Y_t}.
				\label{eq:named}
			\end{equation}
			
			\item \emph{Naturality and substitution.}
			For every chartwise compatible heap homomorphism $q:U\to V$ there is a canonical
			map $\bar q_t:C_t(U)\to C_t(V)$ (surjective if $q$ is),
			and the square
			\begin{equation}
				\begin{array}{ccc}
					\mathbb C^{C_t(V)}&\xrightarrow{\ \bar q_t^*\ }&
					\mathbb C^{C_t(U)}\\
					\big\downarrow\scriptstyle{c_{t,V}^*}&&
					\big\downarrow\scriptstyle{c_{t,U}^*}\\
					\mathbb C^V&\xrightarrow{\ q^*\ }&\mathbb C^U
				\end{array}
				\label{eq:naturality}
			\end{equation}
			commutes.
			For finite heap interface maps $p_i:H_i\to Y$, if their matching fiber product
			$H_1\times_YH_2$ is nonempty, the source algebra additionally obeys the canonical
			substitution isomorphism
			\begin{equation}
				\mathbb C^{H_1\times_YH_2}
				\cong\mathbb C^{H_1}\otimes_{\mathbb C^Y}\mathbb C^{H_2}.
				\label{eq:balanced}
			\end{equation}
			These maps compare source observations. An isomorphism of native
			endpoints would require a comparison of their contexts as well.
		\end{enumerate}
	\end{theorem}
	
	\begin{proof}\leavevmode
		\begin{enumerate}[label=\textup{(\roman*)},leftmargin=*]
			\item \emph{Native strong cut reconstruction.}
			The ternary heap operation is a Mal'cev term, so the congruences
			$\alpha_t,\beta_t$ permute.  Their join is their relational composite.
			Two partial chart tuples have equal images in $C_t(U)$ exactly when
			their lifts are connected by an $\alpha_t$-step and a $\beta_t$-step;
			the intermediate element supplies the missing rectangle.  The
			joint-coordinate embedding identifies $U$ with the full fiber
			product, proving (i).
			
			\item \emph{Exact commutant observation and recovery.}
			On $\ell^2(U)$ the algebra $M_U(S)$ is precisely the diagonal
			functions constant on the fibers of $f_S$: products of characteristic
			functions of individual readout fibers give all joint-fiber
			idempotents.  Thus its inclusion criterion is exactly the definition
			of determination.  Taking the commutant of the algebra generated by
			two families proves the second identity of (ii).  A function factors
			through both $p_L$ and $p_R$ precisely when it is constant on the
			equivalence relation generated by their two kernels, which is
			$\alpha_t\vee\beta_t$.  This gives \eqref{eq:cut-algebra}.  A pullback on finite
			function algebras is injective precisely when its underlying map is
			onto.
			
			\item \emph{Unchanged separators.}
			The quotient $y_t$ is constant on both projection congruences, so
			induces $\lambda_t$.  Since the cut relation is already the full
			fiber product over $C_t(U)$, the larger fiber product over $Y_t$
			adds no pairs if and only if $\lambda_t$ is injective.  Equality of
			the finite function subalgebras in \eqref{eq:named} is equivalent to the same
			condition.  This proves (iii).
			
			\item \emph{Naturality and substitution.}
			Chartwise compatibility carries each of $\alpha_t$ and
			$\beta_t$ into the corresponding kernel on $V$, hence induces
			$\bar q_t$.  It is onto when $q$ is onto; the displayed square
			commutes because the two underlying quotient maps do.  Finally the
			orthogonal idempotents $1_{h_1}\otimes 1_{h_2}$ in the balanced
			tensor product survive exactly when $p_1(h_1)=p_2(h_2)$.  These
			surviving idempotents form the full standard basis of the function
			algebra on the matching fiber product.  This proves \eqref{eq:balanced} and its
			associativity under further nonempty matching.\qedhere
		\end{enumerate}
	\end{proof}
	
	\subsection{Compatibility of concrete spectral representations}
	\label{subsec:fixed-ambient}
	
	In a fixed Hilbert space, products of spectral projections
	determine joint support. Given local multiplicities, a joint
	realization requires a positive integral refinement.
	
	\begin{theorem}[Fixed-ambient spectral-support and multiplicity criteria]
		\label{thm:fixed-ambient-criterion}
		Let $M_L,M_R\subseteq B(\mathcal H)$ be concrete finite-dimensional
		unital abelian $*$-subalgebras, with finite character sets $X_L,X_R$
		and minimal spectral projections $P_x\in M_L$, $Q_z\in M_R$.
		\begin{enumerate}[label=(\roman*),nosep]
			\item \emph{Joint spectral support.} They admit a simultaneous classical joint-spectrum relation
			inside the fixed ambient space exactly when the two algebras
			commute elementwise. In that case
			\begin{equation}
				R_{\mathcal H}=\{(x,z):P_xQ_z\ne0\}
				\subseteq X_L\times X_R,
				\quad
				M_L\vee M_R\cong\mathbb C^{R_{\mathcal H}}.
				\label{eq:support}
			\end{equation}
			The common algebra $M_L\cap M_R$ is the function algebra on the
			connected components of the bipartite support graph of
			$R_{\mathcal H}$.  For onto named separator maps
			$X_L\twoheadrightarrow Y\leftarrow X_R$ compatible with the joint
			support, exact unchanged-interface reconstruction holds if and only if
			\begin{equation}
				P_xQ_z\ne0\quad\text{for every }(x,z)\in X_L\times_Y X_R.
				\label{eq:atomic-probe}
			\end{equation}
			This atom-product criterion is stronger than equality of the common
			observable algebra with $\mathbb C^Y$.
			\item \emph{Prescribed multiplicities.} More generally, fix a nonempty finite joint source
			$U\subseteq\prod_i X_i$ and \emph{unitary-equivalence classes} of
			faithful diagonal representations of $\mathbb C^{X_i}$ with integer
			multiplicities $m_i(x)\ge1$.  There exists a common diagonal
			representation with \emph{exactly} the joint spectrum $U$ whose
			restrictions have those local multiplicities if and only if there
			are positive integers $n_u$ solving
			\begin{equation}
				\sum_{u\in U:\,p_i(u)=x}n_u=m_i(x)
				\quad\text{for all }i,x.
				\label{eq:multiplicity}
			\end{equation}
			Under common amplification $k$, replace the right side by
			$k m_i(x)$.  Existence of a nonnegative solution would instead
			permit loss of joint source support.  For already fixed operator
			embeddings in one Hilbert space, the actual products of spectral
			projections must be checked in addition to \eqref{eq:multiplicity}.
		\end{enumerate}
	\end{theorem}
	
	\begin{proof}\leavevmode
		\begin{enumerate}[label=\textup{(\roman*)},leftmargin=*]
			\item \emph{Joint spectral support.}
			If $M_L$ and $M_R$ commute, their minimal spectral projections
			commute and the products $P_xQ_z$ are mutually orthogonal
			projections summing to the identity.  The nonzero ones are the
			minimal spectral projections of $M_L\vee M_R$, giving \eqref{eq:support}.
			The two component algebras consist of functions of the first or
			second coordinate on that support.  A function common to both is
			constant along each edge and thus along every connected component
			of the incidence graph.  Consequently \eqref{eq:atomic-probe} is precisely the
			full-rectangle condition over the named separator.  Without
			commutation no such family of joint spectral projections exists.
			
			\item \emph{Prescribed multiplicities.}
			A joint diagonal representation decomposes the ambient
			Hilbert space into its nonzero joint eigenspaces
			$\bigoplus_{u\in U}\mathbb C^{n_u}$.  Restricting to the $i$-th
			observable sums these eigenspaces over the fibers of $p_i$, proving
			necessity of \eqref{eq:multiplicity}.  Conversely, positive integral solutions define
			that Hilbert direct sum and commuting diagonal coordinate
			representations; their local spectral multiplicities are exactly
			$m_i(x)$, so each local restriction is unitarily equivalent to its
			specified representation.  Positivity gives joint support exactly
			$U$.  A fixed concrete arrangement of the individual eigenspaces
			need not coincide with this newly built joint representation.\qedhere
		\end{enumerate}
	\end{proof}
	
	\begin{example}[An affine heap and its exact commutant cut]
		\label{ex:affine-heap}
		Take three two-point heaps $A=B=C=\mathbb F_2$ with ternary
		operation $[x,y,z]=x-y+z$ (there is \emph{no designated zero} in
		the native signature), and impose the affine subheap constraint
		\[
		U=\{(a,b,c)\in\mathbb F_2^3:a+b+c=1\}
		=\{001,010,100,111\}.
		\]
		The three generators $001,010,100$ produce $111$ by the ternary
		heap law, so $U$ is an output of $\mathsf J_{\mathrm{hp}}$.
		It is not a subgroup of the fixed coordinate groups because it
		omits $000$.  All three chart projections are surjective; any
		single coordinate varies while another is held fixed, and any
		two coordinates determine the third.  Across $a\mid bc$, the
		canonical interface is the two-point heap $C_t(U)\cong\mathbb F_2$,
		with named maps $a$ on the left and $1+b+c$ on the right.  The
		four-point source is their \emph{full} fiber product.  On
		$\mathcal H_U=\mathbb C^4$, the diagonal source
		$\mathcal A_U\cong\mathbb C^4$ has its required commutant-native
		endpoint, and the common cut algebra is precisely
		$\mathbb C[a]\cong\mathbb C^2$.
		
		The two-dimensional canonical local Hilbert representations
		($m_i(0)=m_i(1)=1$) cannot be retained on one common ambient space
		with joint support $U$: each chart value occurs in two distinct
		words, making \eqref{eq:multiplicity} impossible at amplification one.  At
		amplification $k=2$, taking $n_u=1$ for each of the four words
		solves every equation, so the exact threshold is $k=2$.  Source
		comparison and fixed local representation preservation are
		therefore different contracts.
	\end{example}
	
	\begin{example}[Incomplete joint support and noncommuting readouts]\label{ex:operator-nogo}
		Put $\mathcal H=\mathbb C^3$ with basis indexed by
		$R=\{(0,0),(0,1),(1,1)\}$.  The algebras $M_L$ and $M_R$ of
		functions of the first and second coordinate commute.  Their
		common algebra is only $\mathbb C1$ because the bipartite support
		graph is connected.  Their joint spectrum has only three
		points, since the pair $(1,0)$ is missing from the compatible rectangle.
		In the native context
		$\Xi^{\mathrm{com}}_{\mathbb C^2}$ take
		$M_Z=\operatorname{span}\{I,\sigma_z\}$ and
		$M_X=\operatorname{span}\{I,\sigma_x\}$ with the usual Pauli
		matrices.  Both are abelian unital native inputs with two-point
		spectra, but they do not commute and
		\[
		(M_Z\vee M_X)'=B(\mathbb C^2)'=\mathbb C I,
		\qquad \mathbb C^{X_Z\times X_X}\cong\mathbb C^4.
		\]
		Indeed an operator commuting with $\sigma_z$ is diagonal, and a
		diagonal operator commuting with $\sigma_x$ has equal diagonal
		entries. The anticommutation of the two Pauli matrices rules out a
		simultaneous spectral decomposition by \cref{thm:fixed-ambient-criterion}.
	\end{example}

	\subsection{Distinguishing sources by subdirect attachments}
	\label{subsec:subdirect-attachment-tests}

	The continuations here attach new native charts to marked subgroup relations.
	\Cref{thm:probe-complete-noncollapse} determines exactly when two
	marked inputs give the same answers and counts the subgroup tests needed
	to distinguish all of them.

	\paragraph{Native model and allowed attachments.}
	Let $S$ be a finite nonabelian simple group, set $A=B=S^n$ ($n\geq1$),
	and let
	$
	\mathscr R_n(S)=\{R\leq A\times B:\pi_A R=A,\ \pi_B R=B\}.
	$
	The input is the marked subgroup $R\in\mathscr R_n(S)$ with its
	specified projections. For each $T\in\mathscr R_n(S)$, the test
	attaches a chart with group $T$ using two binary relations.
	These are the graphs of its projections.
	\[
	G_A(T)=\{(t,a)\in T\times A:\pi_A(t)=a\},\qquad
	G_B(T)=\{(t,b)\in T\times B:\pi_B(t)=b\}.
	\]
	Both projection graphs are subgroups. Matching against them tests $T\leq R$.
	
	The continuation also carries five one-bit charts $P_1,P_2,P_3,Z_L,Z_R$.
	The only further formation edge is the unchanged equality edge
	$Z_L=Z_R$, and the two declared recovery descriptors are
	$d_1=P_1,d_2=P_2$ with the identity determination closure on
	$\{d_1,d_2\}$.  Designate the full cut
	\[
	\mathsf L=\{P_2,Z_L\},\qquad
	\mathsf R=\{A,B,T,P_1,P_3,Z_R\}
	\]
	with named separator $z=Z_L=Z_R\in C_2$. For a subgroup candidate $U$
	covering the chart states, let $p_{\mathsf L}$ and $p_{\mathsf R}$ be
	its complete side projections. Define the canonical cut quotient on
	the two actual side images by
	\[
	C_{\mathsf L\mid\mathsf R}(U)
	=U/(\ker p_{\mathsf L}\vee\ker p_{\mathsf R}).
	\]
	Since $z$ factors through both projections and is onto $C_2$, it
	induces a surjection
	$\lambda_U:C_{\mathsf L\mid\mathsf R}(U)\twoheadrightarrow C_2$.
	The cut budget measures the kernel order
	$\beta_{\mathsf L\mid\mathsf R}(U)=|\ker\lambda_U|$.
	The \emph{unchanged} cut requirement is $\beta_{\mathsf L\mid\mathsf R}(U)=1$.
	By \eqref{eq:named}, this gives exact common information and full
	fiber-product reconstruction over the named separator.
	
	Finally fix the same native endpoint law for every test, independently of
	$R,T$: put $Q=C_2$, $E=C_2\times C_2$,
	\[
	f_U:U\to Q,\quad f_U(u)=P_1(u),
	\qquad g:E\twoheadrightarrow Q,\quad g(e_1,e_2)=e_1.
	\]
	For each candidate $U$ the allowed endpoint pairs form the pullback
	$U\times_Q E$. A realization chooses a group homomorphism
	$s:U\to E$ with $gs=f_U$. Its graph is a subgroup of this pullback,
	projecting isomorphically onto $U$ and carrying the prescribed
	chart maps. The table presentation of the graph, with the
	projection to $E$ retained, gives the realized marked endpoint
	in the fixed group-presentation context. The projection $g$
	has a homomorphic section, so $s(u)=(f_U(u),0)$ supplies this clause
	for every candidate. The state image of the attached chart is
	the clause that detects subgroup inclusion in the theorem below.
	
	The test $\operatorname{Acc}(R;T)$ accepts if a proper subgroup $U$
	of the full matching source meets all four requirements. It must
	cover every original and attached chart, have the prescribed closure
	on the two descriptors, satisfy the cut test $\beta=1$, and admit a
	selector for the fixed native endpoint law.
	Within this declared continuation language, $R\equiv_{\mathrm{native}}R'$
	means that all admissible future tests give the same Boolean outcome
	on $R$ and $R'$.  An exact deterministic boundary summary may merge
	inputs only within these equivalence classes.
	
	\begin{theorem}[Subgroup probes: separation, summary size, and minimum probe families]
		\label{thm:probe-complete-noncollapse}
		For every finite nonabelian simple $S$ and $n\geq1$, in the explicitly marked probe language permitting a chart naming each full subgroup $T\in\mathscr R_n(S)$, the following hold. The counting and minimum-probe assertions in (iii)--(iv) concern this same family.
		\begin{enumerate}[label=\textup{(\roman*)},leftmargin=2.5em]
			\item \emph{Probe acceptance.} For all $R,T\in\mathscr R_n(S)$,
			\begin{equation}\label{eq:probe-inclusion}
				\operatorname{Acc}(R;T)=1
				\quad\Longleftrightarrow\quad T\leq R.
			\end{equation}
			Whenever the test accepts, a source can be chosen with index two
			in the full postmatching source, with two independent retained recovery readouts,
			\emph{exact} unchanged named-separator information $C_{\mathsf L\mid
				\mathsf R}(U)\cong C_2$, and a homomorphic native endpoint selector.
			\item \emph{Contextual equivalence.} For the continuation language containing these probes,
			\begin{equation}\label{eq:semantic-iff}
				R\equiv_{\mathrm{native}} R'
				\quad\Longleftrightarrow\quad R=R',
				\qquad R,R'\in\mathscr R_n(S).
			\end{equation}
			Thus the coarsest deterministic boundary summary for the marked
			inputs in $\mathscr R_n(S)$ under these tests has exactly
			$|\mathscr R_n(S)|$ states, and its minimum fixed-length binary
			encoding needs $\lceil\log_2|\mathscr R_n(S)|\rceil$ bits.
			\item \emph{Number of distinct sources.} Writing $a=|\operatorname{Aut}(S)|$, there are exactly
			\begin{equation}\label{eq:exact-count}
				|\mathscr R_n(S)|
				=N_n(S):=\sum_{k=0}^{n}\binom{n}{k}^{\!2} k!\,a^k.
			\end{equation}
			Every member of this family is a perfect group, has full original
			chart images, and is encoded by a partial matching of the $n$ left
			and $n$ right simple factors, labeled by automorphisms of $S$.
			The classification gives an exact finite record. It also provides
			a way to compute the outcome of each attachment directly. For fixed
			$S$, the minimum memory is $\Theta(n\log n)$ bits.
			\item \emph{Minimum number of probes.} Among fixed families of the subgroup probes
			$\operatorname{Acc}(-;T)$ in \eqref{eq:probe-inclusion}, the
			least number that distinguishes all inputs is $a^n n!-1$.
			One such family consists of all but one of the graphs of
			automorphisms of $S^n$.
		\end{enumerate}
	\end{theorem}
	\begin{proof}\leavevmode
		\begin{enumerate}[label=\textup{(\roman*)},leftmargin=*]
			\item \emph{Probe acceptance.}
			Write the simultaneous matching source, including the five one-bit
			charts, as
			\[
			W_{R,T}\cong(R\cap T)\times C_2^3\times C_2,
			\]
			where the final $C_2$ is the common $Z_L=Z_R$ bit.  The new $T$ chart
			has image exactly $R\cap T$.  Thus state-exactness of that chart
			forces $T\leq R$, regardless of any subsequently selected source.
			
			Conversely suppose $T\leq R$.  The source becomes
			$W_{R,T}\cong T\times C_2^3\times C_2$.  Choose
			\begin{equation}\label{eq:strict-source}
				U_T=\{(t,p_1=x,p_2=y,p_3=x,z):
				t\in T,\ x,y,z\in C_2\}.
			\end{equation}
			It is a native subgroup of index $2$, and since $T$ is subdirect,
			all three of the groups $A,B,T$ retain their full state images.  All
			five one-bit charts also remain surjective.  The variables $x,y$ vary independently. Thus each of $d_1,d_2$ varies
			while the other stays fixed, which gives the identity recovery closure.
			The two cut-side images are $C_2^{\{y,z\}}$ and
			$T\times C_2^{\{x,z\}}$, and their joint relation is their complete
			fiber product over the original separator $z\in C_2$.  Hence its
			canonical Goursat quotient is \emph{precisely} the nontrivial named
			separator, with $\beta=1$.  The endpoint relation has a homomorphic
			selector $s(u)=(x,0)$, since $g(s(u))=f_U(u)$, establishing (i).
			
			\item \emph{Contextual equivalence.}
			For distinct subdirect relations $R,R'$, choose one of them as $T$
			so that the other fails to contain it. By (i), this test accepts one
			and rejects the other. An exact deterministic summary must therefore
			keep separate states for them.
			Conversely, equal marked subgroups have the same outcome under
			every native continuation. The continuation classes are therefore the marked
			subgroups themselves \cite{ComonTATA}, giving the state and bit minima in (ii).
			
			\item \emph{Number of distinct sources.}
			By Goursat's lemma \cite[Theorem~2.1]{BauerSenZvengrowski2015},
			a subdirect relation $R\leq S^n\times S^n$ is
			determined by two normal kernels and an oriented isomorphism of their
			quotients. Since $S$ is nonabelian simple, every normal subgroup of $S^n$ is
			a product of coordinate factors \cite{Scott1980,LucchiniMoscatielloSpiga2020}. If the quotient has $k$
			surviving simple factors, choose the retained $k$ left coordinates
			and $k$ right coordinates in $\binom{n}{k}^2$ ways.  Every isomorphism
			of $S^k$ permutes its $k$ minimal normal factors and applies an
			automorphism of $S$ on each; there are $k!a^k$ choices.  These data
			produce distinct marked subgroups and cover every such subgroup,
			proving \eqref{eq:exact-count}.  In fact each is a product of
			$2n-k$ copies of $S$, so every $R$ is perfect and all its abelian
			quotients are trivial.
			
			The record (retained coordinates, their bijection, and the automorphism
			labels) has $O_S(n\log n)$ bits and reconstructs $R$ exactly.  The
			$k=n$ summand gives $N_n(S)\geq a^n n!$, while
			$N_n(S)\leq(n+1)a^n(n!)^2$ follows by bounding each summand;
			Stirling's formula completes the matching-order bit estimate.
			All matching operations can be computed by reconstructing the finite
			relation from its record and applying subgroup pullback, intersection
			and projection.
			
			\item \emph{Minimum number of probes.}
			View an automorphism graph as a perfect matching with
			automorphism labels. Its probe accepts $R$ exactly when that
			perfect matching extends the partial labeled matching for $R$.
			The labeled edges common to all such completions recover the
			partial matching. An edge inconsistent with a fixed match occurs
			in no completion. At free endpoints a proposed labeled edge can be
			excluded by changing the match or its label. Hence the full set
			of these probes distinguishes all inputs.
			If $j$ coordinates remain unmatched on each side, the number of
			accepting completions is $a^j j!$. The distinct possible counts
			differ by at least two, since $a\ge3$ (inner automorphisms alone
			form a copy of $S$). Two distinct completion sets of the same size
			also differ in at least two elements. Deleting any one probe
			therefore preserves separation. Conversely, a full graph input
			can accept only the probe for that graph: a subdirect $T$ contained
			in it projects onto the entire left group, so equals the graph.
			Omitting two full graphs from any fixed probe family makes the two
			corresponding inputs reject every probe. This gives the lower bound.\qedhere
			
		\end{enumerate}
	\end{proof}
	
	\begin{corollary}[Exact summary-state counts for marked $A_5$ subgroup probes]
		\label{cor:a5-count}
		For $S=A_5$, $|\operatorname{Aut}(S)|=120$, and the exact number of
		semantically necessary deterministic states is
		\[
		N_n(A_5)=\sum_{k=0}^{n}\binom{n}{k}^{\!2} k!\,120^k.
		\]
		In particular, $N_1=121$, $N_2=29\,281$, and
		$N_3=10\,628\,281$, requiring respectively $7,15,24$ bits.
		This count exceeds the $60^n n!$ inner-automorphism graphs.
		The permitted attachments distinguish all subdirect correspondences,
		including those defined by proper quotients and outer automorphisms.
		Abelianization and sign quotients lose this orientation information.
	\end{corollary}
	
	\begin{proof}
		The standard $A_5$ automorphism group is $S_5$ of order $120$
		\cite[$A_5$ entry]{ConwayEtAl1985Atlas}.
		Substitute into \eqref{eq:exact-count}; the $k=n$ summand alone is
		$120^n n!>60^n n!$, with further positive summands.  The original
		Goursat kernels and quotients are products of copies of the
		perfect group $A_5$, so their abelianizations are trivial.
	\end{proof}
	
	\begin{proposition}[One summary state without subdirect attachments]
		\label{prop:restricted-grammar-boundary}
		Consider the tests with no attached chart $T$. They require coverage
		of the original charts $A,B$, preserve the two private nondependencies,
		and keep the same separator and native law. Every $R\in\mathscr R_n(S)$
		passes. A single state suffices to summarize these outcomes.
	\end{proposition}
	\begin{proof}
		Without the test chart $T$, take the proper native subgroup of
		$R\times C_2^3\times C_2$ defined by $P_3=P_1$.
		Subdirectness supplies original chart coverage, $P_1,P_2$ remain
		independent, the protected cut remains the full fiber product over
		the common $Z$ bit, and $s(u)=(P_1(u),0)$ supplies the same fixed-law
		selector. Thus every input has the same outcome.
	\end{proof}

	\section{Discussion and conclusions}
	\label{sec:discussion}

	Structural packaging specifies how realized packages form assemblies
	and what those assemblies must preserve. Its formation relations
	determine the admitted combinations of states, its interfaces specify
	the observations used for matching, and its native realization laws
	determine how formations produce packaged endpoints. These data connect
	the existence of an assembly with the preservation of local states,
	recovery and non-recovery judgments, and prescribed common information.
	
	The results establish exact witness-lifting criteria for finite trees,
	faithful native realization criteria for the equational atlases studied
	here, and sharp source minima under specified preservation conditions.
	The group and linear constructions identify the interface information
	needed for composition and give explicit tests for cyclic obstructions.
	The comparison of native laws further shows that the resources required
	for faithful assembly depend on the permitted realizations, even when
	the local chart data are fixed.
	
	The remaining questions concern the design of sufficient interfaces,
	the effective selection of native realizations, and the preservation
	properties attainable on an individual atlas.
	
	\Needspace{5\baselineskip}
	\begin{problem}[Interface design and effective realization]
		\label{prob:native-cycle-obstruction}
		\leavevmode
		\begin{enumerate}[label=\textup{(\roman*)},leftmargin=*]
			\item \emph{Minimum interface signatures.}
			For heterogeneous interface groups with independently
			prescribed partial boundaries, determine the minimum cost
			of a local signature sufficient for the permitted cyclic
			attachments.
			
			\item \emph{Attainable spectra.}
			Determine how the attainable spectra of
			\cref{thm:contract-filtered-spectrum,thm:joint-fixed-kernel-spectrum,thm:coupled-extension-spectrum,thm:affine-coupled-spectrum}
			change when the coupling is unmarked, the reference family
			is restricted, the realizations are correlated, or a
			boundary map is noninjective.
			
			\item \emph{Effective profile selection.}
			In the oriented-side construction of
			\cref{thm:oriented-side-native-globalization},
			which succinct native presentations admit polynomial-time
			profile selection or short certificates of infeasibility?
		\end{enumerate}
	\end{problem}

	\begin{problem}[Attainable preservation properties and uniform implications]
		\label{prob:contract-logic-beyond-trees}
		\leavevmode
		\begin{enumerate}[label=\textup{(\roman*)},leftmargin=*]
			\item \emph{A prescribed atlas.}
			Fix a finite equational atlas and a family of protected
			cuts. Characterize, in terms of the chart maps and readouts,
			the truth assignments to state coverage, local recovery,
			global recovery, and exact cut information that its native
			sources can realize.
			
			\item \emph{Further uniform implications.}
			For which individual varieties that are not
			congruence-permutable do these preservation properties
			satisfy uniform implications beyond those in
			\cref{cor:unrestricted-contract-logic}?
		\end{enumerate}
	\end{problem}

	\section*{Use of AI tools}
	\addcontentsline{toc}{section}{Use of AI tools}
	The author used ChatGPT in mathematical proof drafting, finite checks,
	literature comparison, and manuscript preparation.
	Claude assisted only with review and organization.
	The author takes responsibility for the statements, proofs, attribution,
	and references.

\end{document}